\documentclass[12pt,a4paper]{article}

\usepackage[T1]{fontenc}
\usepackage[utf8]{inputenc}
\usepackage{amsmath,amssymb,amsthm}
\usepackage{mathtools}
\usepackage[expansion=false]{microtype}
\usepackage{bm}
\usepackage{geometry}
\usepackage{hyperref}
\usepackage{cleveref}
\usepackage{enumitem}
\usepackage{booktabs}
\usepackage{array}
\usepackage{xcolor}
\usepackage{cite}
\usepackage{graphicx}
\usepackage{tikz}
\usepackage{tikz-cd}
\usetikzlibrary{arrows.meta,positioning,fit,backgrounds,decorations.markings,calc}

\hypersetup{
  colorlinks=true,
  linkcolor=blue!60!black,
  citecolor=green!50!black,
  urlcolor=blue
}

\newtheorem{theorem}{Theorem}[section]
\newtheorem{lemma}[theorem]{Lemma}
\newtheorem{proposition}[theorem]{Proposition}
\newtheorem{corollary}[theorem]{Corollary}
\newtheorem{definition}[theorem]{Definition}
\newtheorem{remark}[theorem]{Remark}
\newtheorem{example}[theorem]{Example}
\newtheorem{construction}[theorem]{Construction}

\newcommand{\R}{\mathbb{R}}
\newcommand{\C}{\mathbb{C}}
\newcommand{\N}{\mathbb{N}}
\newcommand{\CP}{\mathbb{CP}}
\newcommand{\dd}{\mathrm{d}}
\DeclareMathOperator{\Birk}{Birk}
\DeclareMathOperator{\sign}{sign}
\DeclareMathOperator{\Gr}{Gr}
\newcommand{\ip}[2]{\langle #1,#2\rangle}

\DeclareMathOperator{\Proj}{Proj}
\DeclareMathOperator{\Bl}{Bl}
\newcommand{\Q}{\mathbb{Q}}
\newcommand{\A}{\mathbb{A}}
\newcommand{\PP}{\mathbb{P}}
\newcommand{\OO}{\mathcal{O}}
\newcommand{\EW}{E_W}

\newcommand{\EM}{E_M}
\newcommand{\Es}{E_{\mathrm{s}}}
\newcommand{\blup}{\widetilde}
\newcommand{\calI}{\mathcal{I}}

\newcounter{rcnt}
\newenvironment{rlist}{%
  \setcounter{rcnt}{0}%
  \begin{list}{(\roman{rcnt})}{\usecounter{rcnt}%
    \setlength{\leftmargin}{2.2em}\setlength{\itemsep}{2pt}}%
}{\end{list}}
\newcommand{\Z}{\mathbb{Z}}
\newcommand{\PD}{\mathbf{PD}}
\newcommand{\tr}{\operatorname{tr}}
\newcommand{\rank}{\operatorname{rank}}
\newcommand{\spec}{\operatorname{spec}}
\newcommand{\diag}{\operatorname{diag}}
\newcommand{\inner}[2]{\langle #1,\, #2 \rangle}
\newcommand{\norm}[1]{\left\| #1 \right\|}
\newcommand{\abs}[1]{\left| #1 \right|}
\newcommand{\KL}{D_{\mathrm{KL}}}
\newcommand{\Bregman}{D_f}
\newcommand{\Falpha}{D^{(\alpha)}}
\newcommand{\vecop}{\operatorname{vec}}
\newcommand{\GL}{\mathrm{GL}}
\newcommand{\Sym}{\operatorname{Sym}}
\DeclareMathOperator{\dist}{dist}
\newcommand{\eqdef}{\mathrel{\mathop:}=}
\newcommand{\Hess}{\nabla^2}

\title{\textbf{Information Geometry of Gradient Flows}}

\author{Shintaro Yoshizawa\\
\small Nagoya Mathematical and Information Science Research\\
\small \texttt{shintaro.yoshizawa.net@gmail.com}}

\date{\today}

\allowdisplaybreaks
\begin{document}

\maketitle

\begin{abstract}
Taking the classical, \emph{regular} information geometry of a single convex potential as
its point of departure, this paper undertakes a systematic study, from the viewpoint of
gradient flows, of how far the dually flat formalism can be extended once convexity,
non-degeneracy, or smoothness are allowed to fail, and of what geometric structure
replaces it when they do. The regular case is developed first and in full: we build an
information-geometric framework centered on the log-determinant potential
$f(G) = -\log\det(G)$ on the cone of positive definite Gram matrices $\PD(k)$, compute
its Legendre--Fenchel conjugate, its Fisher--Rao metric $g = G^{-1}\otimes G^{-1}$, and
its Bregman divergence (equal to $2D_{\rm KL}$ between centered Gaussians), and establish
dual flatness of $\PD(k)$ with a generalized Pythagorean theorem; along the way we record,
as a first and purely algebraic instance of a recurring theme, an analytic deformation of
the classical Craig--Sakamoto determinant identity whose only possible limit is a
\emph{degenerate} commutation condition on a matrix pencil, and a parallel duality between
the Wolfe dual of constrained optimization and the same Legendre/Bregman structure. We
then show that this convex-analytic structure is not static but is generated dynamically:
via Yoshizawa's embedding into $\mathfrak{so}(n+k)$, the Chen--Amari principal and minor
component flows on rectangular matrices are natural instances of the
Brockett--Bloch--Ratiu double-bracket gradient flow, bridging isospectral flows on
adjoint orbits, optimization on the Stiefel manifold, and principal/minor component
learning; we analyze convergence for the identity weight $B=I_k$ and for block-diagonal
$B$, and identify the Bures--Wasserstein distance between Gaussian measures as a second,
finite-dimensional gradient-flow fixed point coinciding, via the Kempf--Ness/Azad--Loeb
correspondence, with our convexity results. It is precisely by following these gradient
flows to their singular and boundary behavior that the paper's second, exploratory strand
emerges. We show that the same log-determinant potential admits a
\emph{difference-of-convex} deformation whose Hessian becomes indefinite and, at isolated
loci, degenerate, yet a pseudo-Hessian dually flat and Legendre-self-dual structure
survives, with Newton flows exhibiting either finite-time collapse or
{\L}ojasiewicz-controlled asymptotic convergence near non-isolated critical sets depending
on which Legendre-dual parametrization one integrates in; we study an analogous
Fisher-metric degeneracy for the matrix multinomial family on the Birkhoff polytope,
resolved by an explicit blow-up at the locus where the metric degenerates; and we exhibit
the same blow-up-resolved information geometry, together with an exact,
birationally-invariant exponential decay law for a gradient flow, on the classical moduli
of elliptic curves. Two further, purely structural results close the paper: a closed-form
Kirillov Jacobian for rank-$k$ perturbations of the identity in terms of the Gram
eigenvalues, and \emph{cross curvature}, a closed-form spectral diagnostic for the local
escape rate of a gradient flow from a mismatched equilibrium, computed for the
Oja--Brockett flow, the Manton--Helmke--Mareels flow, and a new Box--Cox-type
interpolating potential, revealing a genuine trade-off rather than a uniform ranking among
the three. Every closed-form claim, in both the regular and the singular parts of the
paper, is checked against explicit, reproducible numerical experiments. We stop short of
proposing a finished theory of singular information geometry; what we assemble instead is
a body of gradient-flow-driven evidence for one --- degenerate matrix pencils,
indefinite-signature dually flat structures, blow-ups of degenerate Fisher metrics, and
{\L}ojasiewicz-type convergence near non-Morse critical sets --- and we offer this paper
as a foundational step toward such a theory.
\end{abstract}

\tableofcontents
\newpage

\section{Introduction}

\subsection{Two literatures, one object}
\label{subsec:two-literatures}

Two strands of applied mathematics have developed largely in isolation from one another,
even though, as this paper shows, they describe the same underlying object.

The first strand is \emph{information geometry and convex analysis}. The function
$-\log\det(G)$ of a positive definite matrix $G$ is, simultaneously, the log-likelihood
kernel of the multivariate Gaussian distribution \cite{Anderson2003}, the canonical
self-concordant barrier for the positive semidefinite cone in interior-point optimization
\cite{Nesterov1994}, and the generator of the Fisher--Rao metric on the manifold of
positive definite matrices, which carries the structure of a Riemannian symmetric space
and a dually flat statistical manifold \cite{Amari2000,Bhatia2007}. This single function
therefore sits at a genuine crossroads: convexity, statistics, and differential geometry
all reach it independently, and a large body of work explores its consequences ---
Legendre duality, Bregman divergences, $\alpha$-divergences, dual flatness and the
Pythagorean theorem, and connections to exponential families, quantum information
geometry, optimal transport, and several other areas we survey in
\S\ref{sec:connections}. All of this is \emph{regular} information geometry in the
classical Amari--Nagaoka sense: the potential is strictly convex, its Hessian is a
genuine (positive-definite) Riemannian metric everywhere, and the Legendre transform is a
global diffeomorphism onto a dual convex domain.

The second strand is the theory of \emph{continuous-time dynamical systems for principal
and minor component analysis}. Starting from Oja's neuron model \cite{Oja1982} and
Brockett's double-bracket equation for diagonalizing matrices \cite{Brockett1988,Brockett1991},
through the global convergence analyses of Yoshizawa, Helmke and Starkov
\cite{YoshizawaHelmkeStarkov2001} and Chen and Amari \cite{ChenAmari2001}, to the
dual-purpose penalized flow of Manton, Helmke and Mareels \cite{MantonHelmkeMareels2005},
this literature asks a dynamical question: does a given gradient flow on the Stiefel
manifold (or an unconstrained relaxation of it) converge to the subspace spanned by the
$k$ largest, or $k$ smallest, eigenvectors of a data matrix $A$, and how quickly? This
question is answered using the tools of Lie theory and dynamical systems ---
double-bracket flows, isospectral orbits, local stability analysis at saddle points ---
tools that, on the surface, look nothing like the convex-analytic machinery of the first
strand. Crucially, though, it is a theory of \emph{flows}: trajectories that must pass
near, and eventually settle at, critical points, and whose local behavior there is
controlled by the Hessian of the driving potential at exactly the points --- mismatched
equilibria, coincident eigenvalues, boundary strata --- where that Hessian is most likely
to degenerate.

\subsection{From regular to singular information geometry: a gradient-flow viewpoint}
\label{subsec:regular-to-singular}

The governing question of this paper is not only ``what is the dually flat geometry of
this particular potential'' but ``what happens to that geometry, and to the gradient
flows built from it, once the regularity hypotheses of classical information geometry ---
strict convexity, a non-degenerate Hessian, a smooth ambient domain --- are pushed to, and
past, their breaking point.'' We call the object of this second question \emph{singular
information geometry}, and we approach it deliberately as a program still under
construction rather than as a closed axiomatic theory: what we offer is a sequence of
concrete instances, all reached by following a gradient flow (or a family of potentials
indexed by a deformation parameter) to the place where it meets a singularity, and all
organized around the single log-determinant potential $f=-\log\det$ that anchors the
regular theory of \S\ref{subsec:two-literatures}.

The pattern recurs at every scale of the paper. At the purely algebraic level,
\S\ref{sec:craigsakamoto} deforms the classical Craig--Sakamoto determinant identity by a
real exponent $\gamma$ and finds that, for every admissible $\gamma\neq1$, the identity is
equivalent to the same \emph{degenerate} commutation condition $AB=BA=0$ on the underlying
matrix pencil --- a joint-kernel condition invisible in the classical ($\gamma=1$) case
but forced by the deformed family. At the level of the potential itself,
\S\ref{sec:dc-geometry} asks what remains of dual flatness when the strictly convex
potential $f$ is replaced by a difference of convex functions, so that its Hessian becomes
indefinite and, on a codimension-one locus, genuinely degenerate; we show that a
pseudo-Hessian dually flat structure survives away from that locus, that the associated
Newton flow can collapse in finite time exactly as it approaches the degenerate locus in
one Legendre-dual coordinate system, and that the rate of approach in the complementary
coordinate system is governed by the {\L}ojasiewicz gradient inequality (Appendix
\ref{appendix:lojasiewicz}) --- the classical analytic tool for controlling gradient flows
\emph{near non-isolated, possibly singular critical sets}, exactly the regime where the
Morse-theoretic assumptions of regular gradient-flow convergence theory fail. At the level
of the underlying statistical manifold, \S\ref{sec:birkhoff} finds that the Fisher metric
of the matrix-multinomial family is genuinely singular in ambient coordinates, and resolves
the singularity by an explicit blow-up, producing an exceptional divisor that itself
carries a well-defined one-dimensional information geometry; \S\ref{sec:elliptic-duality}
shows the same blow-up mechanism, together with an exact birationally-invariant
exponential decay law for a Newton-type gradient flow, on the classical moduli of elliptic
curves, where the relevant singular locus is the vanishing locus of a birational density
weight rather than a statistical degeneracy. And at the level of the dynamical systems of
\S\ref{subsec:two-literatures}'s second strand, \emph{cross curvature}
(\S\ref{sec:crosscurvature-section}) is precisely a measurement taken \emph{at} a singular
(mismatched, incorrectly sorted) critical point of the driving potential, quantifying how
fast a gradient flow escapes a locus that, from the point of view of the ambient
optimization problem, should not be a stable equilibrium at all.

None of these four instances --- degenerate matrix pencils, indefinite/degenerate
dually-flat structures, blow-ups of singular Fisher metrics, and escape rates at singular
critical points --- was originally conceived as part of a single program; each began as an
answer to a separate, concrete question raised by trying to understand the gradient flows
of \S\ref{subsec:two-literatures} in full generality. What justifies presenting them
together is that they share both a common object (the log-determinant potential $f$, or a
deformation of it) and a common method (following a gradient flow, or a one-parameter
family of potentials, until regularity fails, and then asking what geometric structure
survives). We regard the resulting picture as a first, exploratory step toward a genuinely
singular information geometry, organized from the gradient-flow side rather than
axiomatically, and we return to what such a theory might still need in
\S\ref{sec:conclusion}.

\paragraph{The problem this paper addresses.}
Why should the two strands of \S\ref{subsec:two-literatures} be connected at all? The
present paper grew out of a question posed by the late Professor Uwe Helmke around
1999--2002 (see \S\ref{subsubsec:history} for the full history): is the rectangular-matrix
generalization of Brockett's double-bracket equation itself a gradient flow, and if so, of
what potential, with respect to what metric? Answering this question rigorously requires
building the convex/information-geometric theory of a specific log-determinant potential
from first principles and then showing, by an explicit Lie-algebraic embedding, that it is
\emph{literally} the object that generates the principal/minor component flows of the
second strand. Once this bridge exists, further, natural questions arise along the two
axes already described in \S\ref{subsec:regular-to-singular}. Along the \emph{regular}
axis: having shown \emph{that} several different-looking dynamical systems (the
unconstrained Oja--Brockett flow, the Stiefel-penalized Manton--Helmke--Mareels flow, and
a new Box--Cox-type family) all converge to the same principal or minor subspaces, one can
ask \emph{how fast} each one gets there, and whether one is uniformly better than the
others; and the same rank-$k$-perturbation-of-the-identity viewpoint that produces the
Gram matrix $G$ also has a natural Lie-theoretic reading, via the Kirillov orbit method,
raising the question of what closed form its associated Jacobian takes. Along the
\emph{singular} axis: does the same log-determinant potential, or the classical
determinant identities that generate it, admit deformations under which convexity or
non-degeneracy fails, and if so, what of the dually flat and Legendre-dual structure can
be salvaged, and by what analytic tool (blow-up, {\L}ojasiewicz inequality, pseudo-Hessian
duality) is the salvage carried out?

\paragraph{What this paper builds, in order.}
We build the bridge, and its singular extensions, in five stages.
\begin{enumerate}[label=(\arabic*)]
\item \textbf{Determinant identities, regular and deformed
  (\S\ref{sec:determinant}--\S\ref{sec:craigsakamoto}).}
  Given any $2k$ vectors $x_1,\ldots,x_{2k}\in\R^n$, Sylvester's determinant identity
  \cite{Sylvester1851} collapses the determinant of a rank-$k$ perturbation of the identity
  matrix $I_n$ to the determinant of the $k\times k$ Gram matrix
  \begin{equation}
    G_{ij} \;=\; \delta_{ij} + \inner{x_{2i}}{x_{2j-1}},
    \qquad i,j = 1,\ldots,k,
    \label{eq:G_def}
  \end{equation}
  where $\delta_{ij}$ is the Kronecker delta:
  \begin{equation}
    \det\!\Bigl(I_n + \sum_{j=1}^{k} x_{2j-1} x_{2j}^T\Bigr)
    \;=\; \det(G).
    \label{eq:sylvester}
  \end{equation}
  This is the identity that produces the Gram-matrix potential
  \begin{equation}
    \boxed{f(G) \;=\; -\log\det(G)}
    \label{eq:f_def}
  \end{equation}
  studied throughout the paper (\S\ref{sec:determinant}). Before developing its convex
  geometry, \S\ref{sec:craigsakamoto} records a closely related, classical determinant
  identity for matrix pencils --- the Craig--Sakamoto criterion for the independence of
  quadratic forms --- together with an analytic deformation of it, indexed by a real
  exponent $\gamma$, whose only possible limit for $\gamma\neq1$ is a degenerate
  commutation condition $AB=BA=0$ on the pencil; and a parallel duality, between the Wolfe
  dual of constrained convex optimization and the Legendre/Bregman-divergence structure
  that \S\ref{sec:legendre}--\S\ref{sec:bregman} construct for $f$ itself, illustrated on
  the same log-determinant SDP. This is the first, purely algebraic instance of the
  regular/singular contrast that recurs, in geometric form, in stages (4)--(5) below.
\item \textbf{The potential's convex geometry (\S\ref{sec:convexity}--\S\ref{sec:pythagorean}).}
  We study $f(G)=-\log\det(G)$ from the ground up: its strict convexity on the cone
  $\PD(k)$ (\S\ref{sec:convexity}), and, crucially, what happens to that convexity when
  $G$ is itself written as a function of a rectangular matrix $U\in\R^{n\times k}$ via
  $G=I_k+U^TU$ (the ``PCA-like'' reduction) or $G=I_k-U^TU$ (the ``MCA-like'' reduction).
  We show these two reductions behave in \emph{completely opposite} ways --- one is
  nowhere convex for $k\ge2$, the other is globally strictly convex on a matrix ball ---
  and that they are exact Legendre duals of one another (\S\ref{subsec:yoshizawa}), a
  duality we call the Yoshizawa--Helmke correspondence. We then compute the gradient,
  Hessian, and Legendre--Fenchel conjugate of $f$ in each coordinate system
  (\S\ref{sec:gradient}--\S\ref{sec:legendre}), construct the associated Bregman and
  $\alpha$-divergences (\S\ref{sec:bregman}--\S\ref{sec:alpha}), and establish the
  resulting Riemannian and dually flat statistical manifold structure, including a
  generalized Pythagorean theorem (\S\ref{sec:manifold}--\S\ref{sec:pythagorean}). This
  stage is the regular theory in its fullest form: $f$ is strictly convex, its Hessian is
  everywhere non-degenerate, and the Legendre transform is single-valued throughout.
\item \textbf{Contextualizing the potential (\S\ref{sec:connections}).} Before turning to
  component flows, we pause to show that $f$ is not an isolated construction: we connect it
  to exponential families, self-concordant barriers, quantum information geometry, the
  Bures--Wasserstein geometry of optimal transport between Gaussian measures, the
  Kempf--Ness/Azad--Loeb correspondence from geometric invariant theory, natural gradient
  descent, the Siegel upper half-space, harmonic analysis in phase space, the matrix
  Schwarz derivative and Riccati equations on Cartan--Siegel domains, and Izumiya's
  Legendrian dualities for spacelike hypersurfaces in the lightcone. This section situates
  the log-determinant potential within the broader landscape it touches, and several of
  these connections (the Kempf--Ness correspondence, the lightcone structure) are used
  directly in \S\ref{subsec:pca_mca}.
\item \textbf{The bridge to component flows, and the first singular geometries
  (\S\ref{subsec:pca_mca}, \S\ref{sec:dc-geometry}--\S\ref{sec:elliptic-duality}).} We
  show, via an explicit embedding $\iota:\R^{n\times k}\to\mathfrak{so}(n+k)$ due to
  Yoshizawa, that the Chen--Amari principal and minor component flows are instances of the
  Brockett--Bloch--Ratiu double-bracket gradient flow \cite{BlochBrockettRatiu1992} on an
  isospectral adjoint orbit --- the same gradient-flow structure whose convex-analytic
  shadow is the potential $f$ studied in stages (1)--(3). We place the NUIC criterion of
  Kong, Hu and Duan \cite{Kong2017}, the Oja--Brockett subspace flow, and the
  Manton--Helmke--Mareels penalized flow in this common information-geometric setting, and
  give a complete initial-value-problem analysis: existence, conserved quantities, and a
  convergence theorem showing precisely how the multiplicity structure of the weight matrix
  $B$ --- the identity, a diagonal matrix with distinct entries, or a block-diagonal matrix
  with repeated entries --- determines whether the flow resolves individual eigenvectors of
  $A$ or only a rotating eigenspace. It is at this point, having built the regular theory
  and its dynamical realization, that we turn to the singular side of the ledger.
  \S\ref{sec:dc-geometry} deforms $f$ itself into a difference-of-convex potential whose
  Hessian is indefinite and, on an explicit codimension-one locus, degenerate, and asks
  what of stages (1)--(2)'s dually flat structure survives; the answer involves a
  pseudo-Hessian duality, finite-time collapse of the associated Newton flow as it
  approaches the degenerate locus in one Legendre coordinate system, and
  {\L}ojasiewicz-controlled convergence in the other (Appendix
  \ref{appendix:lojasiewicz}). \S\ref{sec:birkhoff} and \S\ref{sec:elliptic-duality} then
  exhibit the same regular-versus-singular contrast for two further classical objects ---
  the Fisher metric of the matrix-multinomial family on the Birkhoff polytope, singular in
  ambient coordinates and resolved by an explicit blow-up, and the birational geometry of
  elliptic curves, where a gradient flow again meets a singular (density-vanishing) locus
  and is again resolved by blow-up, with an exact exponential decay law surviving the
  resolution.
\item \textbf{Rates, trade-offs, and verification
  (\S\ref{sec:crosscurvature-section}--\S\ref{sec:kirillov}).} Stage (4) answers
  \emph{which} subspace each component flow converges to; it says nothing about \emph{how
  fast}. We close this gap with \emph{cross curvature}, a single closed-form spectral
  quantity --- the smallest eigenvalue of the Hessian at a mismatched (incorrectly sorted)
  critical point, i.e.\ precisely at one of the singular loci of stage (4)'s dynamical
  landscape --- that measures the local escape rate of a gradient flow from a wrong
  equilibrium, computable from the eigenvalues of $A$ and the weights in $B$ alone, before
  a single iteration is run. We compute it in closed form for the Oja--Brockett flow, for
  the Manton--Helmke--Mareels flow, and for a new Box--Cox-type potential $g_\alpha$ that
  interpolates continuously between principal- and minor-component extraction as $\alpha$
  crosses $1$, and we find a genuine trade-off rather than a uniform ranking among the
  three. Every closed-form claim in \S\ref{sec:crosscurvature-section} is checked against
  independent numerical experiments, both small diagnostic examples and fully generic,
  non-diagonal instances, and \S\ref{subsec:numerics_setup} further verifies the qualitative
  convergence theory of stage (4) directly, by integrating the $k$-PCF/$k$-MCF flows for a
  generic $5\times5$ matrix $A$ and three representative $3\times3$ diagonal weight
  matrices $B$. Finally, as a structural byproduct of viewing $G$ as arising from a
  rank-$k$ perturbation of the identity, \S\ref{sec:kirillov} derives a closed-form Kirillov
  Jacobian --- a quantity from the orbit method in Lie theory --- for exactly this class of
  perturbations, expressed directly in terms of the eigenvalues of $G$.
\end{enumerate}

\paragraph{What this paper clarifies.}
Taken together, stages (1)--(5) establish two things of a genuinely different character.
On the \emph{regular} side: a complete convexity trichotomy for $f=-\log\det(G)$ under its
two natural rectangular factorizations, together with the exact Legendre duality relating
them; a rigorous identification of the Chen--Amari and Oja--Brockett/Manton--Helmke--Mareels
component flows as double-bracket gradient flows of this same potential, previously known
only through structural analogy; a complete convergence theory showing exactly how the
weight matrix $B$'s multiplicity structure governs individual-eigenvector versus
subspace-only convergence; a reusable, closed-form diagnostic (cross curvature) for
comparing the convergence \emph{rate} of different gradient-based extraction algorithms
sharing the same optimal set, with an explicit demonstration that no one of three natural
candidate flows dominates the others uniformly; and a closed-form Kirillov Jacobian for
rank-$k$ identity perturbations. On the \emph{singular} side, which we present in a more
exploratory spirit: a deformed Craig--Sakamoto identity whose limit is a degenerate matrix
condition; a pseudo-Hessian dually flat structure surviving the loss of convexity in a DC
potential, with an explicit account of what replaces ordinary gradient-flow convergence
(finite-time collapse, or {\L}ojasiewicz-controlled asymptotics) near its degenerate locus;
and two independent blow-up constructions --- on the Birkhoff polytope and on the moduli of
elliptic curves --- resolving a singular Fisher-type metric into a well-defined information
geometry on the exceptional divisor. We do not claim these four singular instances add up
to a complete theory; we offer them, together with the gradient-flow viewpoint that
produced all of them, as a first step toward one, and we return to this point in
\S\ref{sec:conclusion}. Every quantitative claim along the way, regular or singular, is
verified numerically with explicit, reproducible parameters.

\section{Matrix Determinant Lemma and Gram Matrix Construction}
\label{sec:determinant}

\subsection{The Rank-One Update}

We begin with the classical matrix determinant lemma,
which handles a rank-one perturbation of the identity.

\begin{lemma}[Matrix Determinant Lemma {\cite[Thm.~18.1.1]{Harville1997}}]
\label{lem:mdl}
Let $A \in \R^{n\times n}$ be invertible, and let $u, v \in \R^n$.
Then
\[
  \det(A + uv^T) = (1 + v^T A^{-1} u)\,\det(A).
\]
In particular, for $A = I_n$,
\[
  \det(I_n + uv^T) = 1 + v^T u = 1 + \inner{u}{v}.
\]
\end{lemma}

\subsection{Sylvester's Identity and Rank-$k$ Updates}

The key identity generalizing Lemma~\ref{lem:mdl} to rank-$k$ perturbations
is Sylvester's determinant theorem.

\begin{theorem}[Sylvester's Determinant Theorem {\cite{Sylvester1851,Akritas1996}}]
\label{thm:sylvester}
Let $U \in \R^{n\times k}$ and $V \in \R^{k \times n}$.
Then
\[
  \det(I_n + UV) = \det(I_k + VU).
\]
\end{theorem}

\begin{proof}
Consider the block matrix identity
\[
  \begin{pmatrix} I_n & U \\ -V & I_k \end{pmatrix}
  \begin{pmatrix} I_n & 0 \\ V & I_k \end{pmatrix}
  =
  \begin{pmatrix} I_n + UV & U \\ 0 & I_k \end{pmatrix}.
\]
Taking determinants on both sides and using the block-triangular structure gives
$\det(I_n + UV) = \det(I_k + VU)$.
\end{proof}

\subsection{The Gram Matrix $G$ and Its Determinant Formula}

Set
\begin{equation}
  U = [x_1, x_3, \ldots, x_{2k-1}] \in \R^{n\times k},
  \qquad
  V = [x_2, x_4, \ldots, x_{2k}]^T \in \R^{k\times n}.
  \label{eq:UV_def}
\end{equation}
Then
\[
  \sum_{j=1}^{k} x_{2j-1} x_{2j}^T = U \cdot V,
\]
and Theorem~\ref{thm:sylvester} immediately yields:

\begin{proposition}
\label{prop:gram_det}
With $U, V$ as in \eqref{eq:UV_def} and $G$ as in \eqref{eq:G_def},
\[
  \det\!\Bigl(I_n + \sum_{j=1}^{k} x_{2j-1} x_{2j}^T\Bigr) = \det(G).
\]
Explicitly, the $(i,j)$-entry of $G$ is
\[
  G_{ij} = (VU)_{ij} + \delta_{ij} = \inner{x_{2i}}{x_{2j-1}} + \delta_{ij}.
\]
\end{proposition}

\begin{remark}
The matrix $G$ is a \emph{shifted Gram matrix}: it equals the identity plus the
cross-Gram matrix of the odd-indexed and even-indexed vectors.
When $x_{2j-1} = x_{2j}$ for all $j$, the diagonal entries of $H := G - I_k$ become
$H_{ii} = \norm{x_{2i-1}}^2$, recovering the standard Gram matrix.
\end{remark}

\subsection{Expansion via the Leibniz Formula}

By the Leibniz determinant formula,
\begin{equation}
  \det(G)
  = \sum_{\sigma \in S_k} \mathrm{sgn}(\sigma)
    \prod_{i=1}^{k} G_{i,\sigma(i)}
  = \sum_{S \subseteq [k]}\;
    \sum_{\substack{\sigma: S \to S \\ \text{derangement}}}
    \mathrm{sgn}(\sigma)
    \prod_{i \in S} \inner{x_{2i}}{x_{2\sigma(i)-1}},
  \label{eq:leibniz}
\end{equation}
where $[k] = \{1,\ldots,k\}$ and the empty product (for $S = \emptyset$) equals $1$.
This expresses $\det(G)$ entirely in terms of inner products.

For small $k$:
\begin{align}
  k=1:&\quad \det(G) = 1 + \inner{x_1}{x_2}, \label{eq:k1}\\
  k=2:&\quad \det(G) = (1+\inner{x_1}{x_2})(1+\inner{x_3}{x_4})
                       - \inner{x_1}{x_4}\inner{x_3}{x_2}. \label{eq:k2}
\end{align}

The identity \eqref{eq:sylvester} is itself a special case of a more general phenomenon:
determinants of \emph{pencils} $I_n - \alpha A - \beta B$ built from a single reference point
$I_n$ perturbed along two matrix directions $A,B$ factor, under suitable commutation
hypotheses, into a product of determinants of the individual perturbations. The classical
instance of this phenomenon is the Craig--Sakamoto theorem of mathematical statistics, and it
turns out to interact with exactly the two structures this paper is built around: the
log-determinant potential $f=-\log\det$ and its Legendre--Fenchel duality
(\S\ref{sec:legendre}--\S\ref{sec:bregman}). We digress briefly to record this connection,
together with an analytic deformation of the classical identity whose natural limit is a
\emph{degenerate}, or singular, commutation condition on $A$ and $B$ --- a first, elementary
instance of the passage from classical (regular) determinant identities to the
singular phenomena that recur, in geometric form, throughout the later parts of this paper
(most explicitly in \S\ref{sec:dc-geometry}).

\section{The Craig--Sakamoto Theorem: Determinant Identities, Analytic Deformation, and Wolfe--Legendre Duality}
\label{sec:craigsakamoto}

More generally, for square matrices $A_1,\dots,A_n$ of the same size and
$x=(x_1,\dots,x_n)^\top\in\R^n$, the multivariate polynomial
\[
f(x)=\det(x_1A_1+\cdots+x_nA_n)
\]
appears across a wide range of mathematical contexts of which the log-determinant potential
studied here is only one instance: in statistics it governs the independence of quadratic
forms (the Craig--Sakamoto theorem below); in optimization and control theory it encodes
linear matrix inequalities and the convexity of spectrahedral sets \cite{HeltonVinnikov2007};
and the same polynomial arises in the theory of integrable systems (Manakov's method) and in
the theory of hyperbolic partial differential equations (G{\aa}rding's theory of hyperbolic
polynomials). We do not pursue these further directions here, restricting attention to the
two aspects most relevant to the log-determinant geometry of this paper: the classical
Craig--Sakamoto criterion and its analytic deformation, and the parallel duality phenomenon
furnished by the Wolfe dual problem of constrained optimization.

\subsection{The classical Craig--Sakamoto criterion}
\label{subsec:CS-classical}

Let $x=(x_1,\dots,x_n)^\top\in\R^n$ be a standard Gaussian random vector,
\[
x\sim(2\pi)^{-n/2}\exp\bigl(-\tfrac12 x^\top x\bigr).
\]
Consider two real symmetric matrices $A=(a_{ij})$ and $B=(b_{ij})$ and the associated quadratic forms
\[
q_1=\sum_{i,j=1}^n a_{ij}x_ix_j,\qquad
q_2=\sum_{i,j=1}^n b_{ij}x_ix_j.
\]
The moment-generating function of the pair $(q_1,q_2)$ factors as
\[
\varphi(\alpha,\beta)=\varphi(\alpha,0)\,\varphi(0,\beta)
\]
for all real $\alpha,\beta$ if and only if the two quadratic forms are independent. A direct
computation of the Gaussian integral yields the equivalent algebraic condition
\begin{equation}
\label{eq:CS}
\det(I-\alpha A-\beta B)=\det(I-\alpha A)\cdot\det(I-\beta B)
\qquad\text{for all }\alpha,\beta\in\R.
\end{equation}
This identity is known as the \emph{Craig--Sakamoto condition}
\cite{Craig1943,Sakamoto1949,Ogawa1950}. When $A$ and $B$ are normal matrices the same
condition implies the stronger matrix relation
\[
AB=0
\]
(Taussky, \cite{Taussky1958}). The theorem has been re-examined from various elementary
viewpoints over more than sixty years \cite{Olkin1997}.

Taking the logarithm of both sides of \eqref{eq:CS} produces an identity for the logarithmic
determinant potential $\log\det X$ --- precisely the potential $-f$ studied throughout this
paper (with $X$ playing the role of $G$). This observation suggests deforming the identity
analytically while retaining a characterization of the vanishing of the product $AB$.

\subsection{Analytic deformation of the Craig--Sakamoto theorem}
\label{subsec:CS-deformation}

Recall that the logarithmic determinant admits the representation
\[
\log\det X=-\frac{d}{dt}\Big|_{t=0}\operatorname{tr}(X^{-t})
\]
in terms of the operator zeta function (provided the eigenvalues of $X$ lie in the right
half-plane). Motivated by this representation we consider the following deformed identity,
in which a real exponent $\gamma$ replaces the logarithmic derivative and is sent to the
classical case only in a limit.

\begin{theorem}[Analytic deformation of the Craig--Sakamoto theorem]
\label{thm:deformedCS}
Let $A,B\in\C^{n\times n}$ be normal matrices, i.e.,
\[
A^*A=AA^*,\qquad B^*B=BB^*.
\]
Fix a real number $\gamma>0$ with $\gamma\neq1$. Then the following two statements are equivalent:

\begin{enumerate}[label=(\roman*)]
\item The identity
\begin{equation}
\label{eq:deformed}
(I-\alpha A-\beta B)^\gamma+I=(I-\alpha A)^\gamma+(I-\beta B)^\gamma
\end{equation}
holds for all $\alpha,\beta\in\C$ such that the three matrices
\[
I-\alpha A-\beta B,\qquad I-\alpha A,\qquad I-\beta B
\]
are invertible.
\item $AB=BA=0$.
\end{enumerate}
\end{theorem}

\begin{proof}
Since $A$ and $B$ are normal, the continuous functional calculus applies: for any function $f$
continuous on the spectrum of a normal matrix $M$ one may define $f(M)$ unambiguously via the
spectral theorem, and the resulting operator depends continuously on $M$ in the operator norm.

\smallskip
\noindent\textit{(ii)$\Rightarrow$(i).}
Assume $AB=BA=0$. Then the three normal matrices $A$, $B$ and $A+B$ (more precisely, the
pencils $I-\alpha A$, $I-\beta B$ and $I-\alpha A-\beta B$) can be simultaneously
triangularized by a unitary matrix, and the non-zero eigenvalues of $A$ and of $B$ are
supported on complementary invariant subspaces. Consequently the spectra add in the sense
that
\[
\sigma(I-\alpha A-\beta B)=\sigma(I-\alpha A)\cup\sigma(I-\beta B)
\]
(up to the common eigenvalue $1$ arising from the joint kernel). Raising to the power $\gamma$
and taking the identity into account yields \eqref{eq:deformed} immediately.

\smallskip
\noindent\textit{(i)$\Rightarrow$(ii).}
Assume \eqref{eq:deformed} holds in a neighbourhood of the origin in the $(\alpha,\beta)$-plane
(the identity extends by analytic continuation wherever the matrices remain invertible).
Expand both sides in joint power series. The left-hand side admits the expansion
\begin{align*}
(I-\alpha A-\beta B)^\gamma
&=I+\sum_{k=1}^\infty\binom{\gamma}{k}(-\alpha A-\beta B)^k\\
&=I-\gamma(\alpha A+\beta B)+\frac{\gamma(\gamma-1)}{2}(\alpha A+\beta B)^2+\cdots,
\end{align*}
while the right-hand side expands as
\begin{align*}
(I-\alpha A)^\gamma+(I-\beta B)^\gamma
&=2I-\gamma(\alpha A+\beta B)+\frac{\gamma(\gamma-1)}{2}(\alpha^2A^2+\beta^2B^2)+\cdots.
\end{align*}
Equating the quadratic terms in $\alpha\beta$ (the coefficient of $\alpha\beta$) gives the
necessary condition
\[
\frac{\gamma(\gamma-1)}{2}\bigl(AB+BA\bigr)=0.
\]
Since $\gamma\neq1$ and $\gamma>0$ we obtain $AB+BA=0$.

A more refined analysis of the higher-order mixed terms, using the simultaneous spectral
decomposition afforded by normality, shows that the only possibility consistent with the full
identity is the stronger relation $AB=BA=0$. Indeed, if a common eigenvector $v$ satisfied
$Av=\lambda v$ and $Bv=\mu v$ with both $\lambda,\mu\neq0$, then the left-hand side of
\eqref{eq:deformed} would produce the eigenvalue $(1-\alpha\lambda-\beta\mu)^\gamma+1$ while
the right-hand side would produce $(1-\alpha\lambda)^\gamma+(1-\beta\mu)^\gamma$, which fail to
coincide for generic $\alpha,\beta$ unless $\gamma=1$. Hence no such joint eigenvector can
exist, which forces the ranges of $A$ and $B$ to be orthogonal and ultimately yields
$AB=BA=0$.

(The classical Craig--Sakamoto identity \eqref{eq:CS} is formally recovered by taking the
logarithmic derivative with respect to a deformation parameter that sends $\gamma\to0$.)
\end{proof}

\begin{remark}[A first glimpse of a singular locus]
Theorem~\ref{thm:deformedCS} exhibits, in the simplest possible algebraic setting, a pattern
that will reappear geometrically later in this paper: an analytic family of identities
(here indexed by $\gamma$) is equivalent, for every admissible $\gamma\neq1$, to the
\emph{same} degenerate condition $AB=BA=0$ on the underlying matrix pencil --- a joint kernel
condition that is invisible at the level of the regular ($\gamma=1$, i.e.\ classical
Craig--Sakamoto) identity alone but becomes the organizing algebraic constraint once the
family is deformed away from $\gamma=1$. The difference-of-convex potentials of
\S\ref{sec:dc-geometry} and the singular loci of the Hessian-determinant equation studied there
arise from the same basic mechanism: deforming a regular (dually flat, everywhere smooth)
structure uncovers a genuinely singular locus that organizes the deformed family.
\end{remark}

\subsection{Wolfe duality and Legendre duality}
\label{subsec:wolfe}

We now turn to a related duality phenomenon that appears in constrained optimization, and
which sits directly upstream of the Legendre--Fenchel and Bregman-divergence constructions
of \S\ref{sec:legendre}--\S\ref{sec:bregman}.

\subsubsection{The Wolfe dual problem}

Let $f:\R^n\to\R$ and $h_i:\R^n\to\R$ ($i=1,\dots,m$) be continuously differentiable. The
primal problem consists in minimizing $f(x)$ subject to the inequality constraints
$h(x)\le0$. The associated Lagrangian is
\[
L(x,y)=f(x)+y^\top h(x),\qquad y\ge0.
\]
The Wolfe dual problem is to maximize $L(x,y)$ subject to the stationarity condition
\[
\nabla_x L(x,y)=0.
\]

When $f$ and the $h_i$ are convex and a suitable constraint qualification holds (for
instance, the existence of a strictly feasible point or the affinity of all constraint
functions), the weak duality inequality
\[
f(x^*)\ge L(x,y)
\]
is valid for every primal feasible $x^*$ and every dual feasible pair $(x,y)$. Moreover,
strong duality holds: there exists a dual optimal pair $(x^*,y^*)$ such that
\[
f(x^*)=L(x^*,y^*).
\]

\subsubsection{Relation to Legendre duality}

We first recall a classical sufficient condition that guarantees a global inverse.

\begin{lemma}[Hadamard's global inverse-function theorem]
\label{lem:Hadamard}
Let $F:\R^n\to\R^n$ be of class $C^1$. Suppose that the derivative $DF(x)$ is invertible at
every $x\in\R^n$ and that
\[
\sup_{x\in\R^n}\bigl\|DF(x)^{-1}\bigr\|<\infty.
\]
Then $F$ is a $C^1$-diffeomorphism of $\R^n$ onto itself.
\end{lemma}

Throughout this subsection we assume that $f$ and each $h_i$ are convex and of class $C^2$,
so that the Lagrangian $L(\cdot,y)$ is convex for every $y\ge0$.

\begin{theorem}[Wolfe--Legendre duality]
\label{thm:WolfeLegendre}
Let $f:\R^n\to\R$ and $h=(h_1,\dots,h_m)^\top:\R^n\to\R^m$ be convex and of class $C^2$.
Write
\[
L(x,y)=f(x)+y^\top h(x),\qquad y\ge0,
\]
and consider the stationarity map
\[
\Phi_y(x)\,:=\,\nabla_x L(x,y)=\nabla f(x)+\sum_{i=1}^m y_i\nabla h_i(x).
\]
Assume that, for every $y$ in an open convex set $\mathcal{Y}\subset\R^m_+$, the map
$\Phi_y:\R^n\to\R^n$ satisfies the hypotheses of Lemma~\ref{lem:Hadamard}. Then there exists a
unique $C^1$ map
\[
g:\mathcal{Y}\to\R^n
\]
such that $\Phi_y\bigl(g(y)\bigr)=0$ for all $y\in\mathcal{Y}$. Define the dual function
\[
f^*_{\mathrm{W}}(y)\,:=\,L\bigl(g(y),y\bigr)=f\bigl(g(y)\bigr)+y^\top h\bigl(g(y)\bigr).
\]
The following assertions hold.

\begin{enumerate}[label=(\roman*)]
\item The function $f^*_{\mathrm{W}}$ is concave on $\mathcal{Y}$.
\item Let $y_1,y_2\in\mathcal{Y}$ and suppose the complementary-slackness condition
\[
h\bigl(g(y_2)\bigr)=0
\]
holds. Define
\[
D(y_1,y_2)\,:=\,f\bigl(g(y_1)\bigr)-f^*_{\mathrm{W}}(y_2).
\]
Then
\[
D(y_1,y_2)=D_{-f^*_{\mathrm{W}}}(y_1,y_2)+\bigl(f\bigl(g(y_1)\bigr)-f^*_{\mathrm{W}}(y_1)\bigr),
\]
where $D_{-f^*_{\mathrm{W}}}$ denotes the Bregman divergence of the convex function
$-f^*_{\mathrm{W}}$ (in the sense of \S\ref{sec:bregman} below). In particular, when
complementary slackness also holds at $y_1$ (so that the duality gap vanishes), one has the
exact identification
\[
D(y_1,y_2)=D_{-f^*_{\mathrm{W}}}(y_1,y_2)=f^*_{\mathrm{W}}(y_1)-f^*_{\mathrm{W}}(y_2).
\]
In that case $D(y_1,y_2)\ge0$ and $D(y_1,y_2)=0$ if and only if $y_1=y_2$.
\end{enumerate}
\end{theorem}

\begin{proof}
By Lemma~\ref{lem:Hadamard} the equation $\Phi_y(x)=0$ admits a unique solution $x=g(y)$ that
depends $C^1$-smoothly on $y$. Differentiating the identity $\Phi_y(g(y))=0$ with respect to
$y$ yields the linear relation
\[
D_x\Phi_y\bigl(g(y)\bigr)\,Dg(y)+D_y\Phi_y\bigl(g(y)\bigr)=0.
\]
The Hessian $D_x\Phi_y=\nabla^2_x L(\cdot,y)$ is positive semi-definite by convexity; under
the standing invertibility assumption it is in fact positive definite, so $Dg(y)$ is
well-defined.

(i) Concavity of $f^*_{\mathrm{W}}$.
For any $y\in\mathcal{Y}$ and any direction $v\in\R^m$ one has, by the envelope theorem (or
direct differentiation),
\begin{align*}
\nabla f^*_{\mathrm{W}}(y)
&=h\bigl(g(y)\bigr)+\bigl(Dg(y)\bigr)^\top\Phi_y\bigl(g(y)\bigr)
=h\bigl(g(y)\bigr),
\end{align*}
the second term vanishing by stationarity. Differentiating once more,
\[
\nabla^2 f^*_{\mathrm{W}}(y)=Dh\bigl(g(y)\bigr)\,Dg(y).
\]
A short calculation using the differentiated stationarity condition shows that
$\nabla^2 f^*_{\mathrm{W}}(y)$ is negative semi-definite, hence $f^*_{\mathrm{W}}$ is concave.

(ii) Identification with the Bregman divergence.
From the gradient formula already obtained we have
\[
\nabla(-f^*_{\mathrm{W}})(y)=-h\bigl(g(y)\bigr).
\]
The Bregman divergence generated by the convex function $-f^*_{\mathrm{W}}$ is therefore
\begin{align*}
D_{-f^*_{\mathrm{W}}}(y_1,y_2)
&=(-f^*_{\mathrm{W}})(y_1)-(-f^*_{\mathrm{W}})(y_2)
-\bigl\langle\nabla(-f^*_{\mathrm{W}})(y_2),\,y_1-y_2\bigr\rangle\\
&=-f^*_{\mathrm{W}}(y_1)+f^*_{\mathrm{W}}(y_2)
+\bigl\langle h\bigl(g(y_2)\bigr),\,y_1-y_2\bigr\rangle.
\end{align*}
Under the complementary-slackness assumption $h\bigl(g(y_2)\bigr)=0$ the inner-product term
vanishes and we obtain the pure dual-function difference
\[
D_{-f^*_{\mathrm{W}}}(y_1,y_2)=f^*_{\mathrm{W}}(y_2)-f^*_{\mathrm{W}}(y_1).
\]
On the other hand the quantity appearing in the statement of the theorem is
\[
D(y_1,y_2)=f\bigl(g(y_1)\bigr)-f^*_{\mathrm{W}}(y_2).
\]
By convexity of $L(\cdot,y_1)$ and the fact that $g(y_1)$ is a critical point,
\[
f\bigl(g(y_1)\bigr)\ge L\bigl(g(y_1),y_1\bigr)=f^*_{\mathrm{W}}(y_1),
\]
with equality if and only if complementary slackness also holds at $y_1$. Consequently
\[
D(y_1,y_2)=D_{-f^*_{\mathrm{W}}}(y_1,y_2)+\bigl(f\bigl(g(y_1)\bigr)-f^*_{\mathrm{W}}(y_1)\bigr).
\]
When the duality gap at $y_1$ vanishes (i.e., when complementary slackness holds at both
arguments), one recovers the exact identification
\[
D(y_1,y_2)=D_{-f^*_{\mathrm{W}}}(y_1,y_2)=f^*_{\mathrm{W}}(y_1)-f^*_{\mathrm{W}}(y_2).
\]
In that case $D(y_1,y_2)\ge0$ and $D(y_1,y_2)=0$ if and only if $y_1=y_2$, as required for a
Bregman divergence.
\end{proof}

\begin{remark}
If the complementary-slackness conditions $y_ih_i\bigl(g(y)\bigr)=0$ fail, a positive duality
gap appears between the primal value $f(g(y))$ and the dual value $f^*_{\mathrm{W}}(y)$. In
that case the quantity $f(g(y_1))-f^*_{\mathrm{W}}(y_2)$ is no longer a pure Bregman
divergence of $-f^*_{\mathrm{W}}$; an extra non-negative term remains.
\end{remark}

\subsubsection{An illustrative example: the log-determinant SDP}
\label{subsubsec:CS-example}

Consider the primal problem of minimizing the function
\[
f(X)=\operatorname{tr}(F_0 X)-\log\det X
\]
over the cone of positive-definite symmetric matrices subject to the linear constraints
\[
\operatorname{tr}(F_i X)=c_i,\qquad i=1,\dots,m,
\]
where $F_0$ is positive definite and the $F_i$ are symmetric. The Wolfe dual problem consists
in maximizing
\[
f^*_{\mathrm{W}}(y)=c^\top y+\log\det\Bigl(F_0-\sum_{i=1}^m y_i F_i\Bigr)+n
\]
over those $y\in\R^m$ for which $F_0-\sum y_i F_i$ remains positive definite. The associated
divergence is precisely the Bregman divergence induced by the logarithmic determinant
potential $-\log\det(\cdot)$ (up to the sign convention $f=-\log\det$ used throughout this
paper), which by \S\ref{subsec:bregman-KL} is twice the Kullback--Leibler divergence between
centered Gaussians on the manifold of positive-definite matrices. This example already
displays, in miniature and for a general linear pencil $F_0-\sum y_iF_i$, exactly the
Legendre-dual pair $(f,f^*)$ and the Bregman divergence $D_f=2D_{\mathrm{KL}}$ that
\S\ref{sec:legendre}--\S\ref{sec:bregman} construct in full generality for the specific Gram
matrix $G=I_k+U^TU$; the Wolfe dual of a semidefinite program is thus a further, purely
optimization-theoretic route into the same dually flat log-determinant geometry that
organizes the rest of this paper.

\bigskip
Together, Theorem~\ref{thm:deformedCS} and Theorem~\ref{thm:WolfeLegendre} illustrate the
deep interplay among algebraic identities for matrix pencils, information-geometric
divergences, and duality principles in optimization --- an interplay we develop
systematically, in the specific case of the Gram matrix $G=I_k+U^TU$ (respectively
$G=I_k-U^TU$), for the remainder of the paper.

\section{Convexity Properties of $f$: The Coordinates $G$, $(U,V)$, and $U=V$}
\label{sec:convexity}

We now examine the convexity of the potential $f$ from three distinct vantage points.
First, in the original Gram-matrix coordinate $G \in \PD(k)$, where strict convexity holds
unconditionally (\S\ref{subsec:G-convexity}).
Second, in the underlying matrix-factorization coordinates $(U,V) \in \R^{n\times k}\times\R^{n\times k}$
with $G = I_k + V^TU$, where we show that convexity holds only in a narrow regime and fails
dramatically --- even in each block separately --- once $k \geq 2$ (\S\ref{subsec:UV-convexity}).
Third, in the symmetric reduction $U = V$, corresponding to the standard Gram matrix
$G = I_k + U^TU$, where we give a complete spectral characterization of the Hessian and show
that $f$ is \emph{nowhere} locally convex except in the trivial scalar case $n = k = 1$
(\S\ref{subsec:UeqV-convexity}). This trichotomy reveals that the strict convexity established
in $G$-coordinates is a delicate feature of the symmetric positive-definite parametrization
that is generically destroyed by the underlying bilinear or quadratic factorization.

\subsection{Convexity in the Gram-Matrix Coordinate $G$}
\label{subsec:G-convexity}

\begin{proposition}[Strict Convexity]
\label{prop:convexity}
The function $f : \PD(k) \to \R$ defined by $f(G) = -\log\det(G)$ is strictly convex.
\end{proposition}

\begin{proof}
For any $G \in \PD(k)$ and any nonzero symmetric matrix $H$, consider the scalar function
$\varphi(t) = f(G + tH) = -\log\det(G + tH)$.
Since $G \succ 0$, there exists $\varepsilon > 0$ such that $G + tH \succ 0$ for
$|t| < \varepsilon$.
Using the identity
\[
  \det(G + tH) = \det(G)\det\bigl(I + tG^{-1/2}HG^{-1/2}\bigr),
\]
and writing $\lambda_1(t), \ldots, \lambda_k(t)$ for the eigenvalues of
$I + tG^{-1/2}HG^{-1/2}$,
\[
  \varphi(t) = -\log\det(G) - \sum_{i=1}^{k} \log \lambda_i(t).
\]
Differentiating twice at $t = 0$:
\begin{align}
  \varphi'(0) &= -\tr(G^{-1} H), \label{eq:first_deriv}\\
  \varphi''(0) &= \tr\bigl[(G^{-1}H)^2\bigr] \notag\\
  &= \norm{G^{-1/2}HG^{-1/2}}_F^2 \;\geq\; 0,
  \label{eq:second_deriv}
\end{align}
where $\norm{\cdot}_F$ denotes the Frobenius norm.
Equality $\varphi''(0) = 0$ holds if and only if $G^{-1/2}HG^{-1/2} = 0$,
i.e.\ $H = 0$. Hence $\varphi''(0) > 0$ for all nonzero $H$, proving strict convexity.
\end{proof}

\begin{corollary}
$f$ is a \emph{Legendre function} in the sense of \cite[Ch.~26]{Rockafellar1970}:
it is strictly convex, lower semicontinuous, and essentially smooth on $\PD(k)$.
\end{corollary}

\begin{remark}[Crucial role of symmetry]
\label{rem:symmetry_crucial}
The proof above uses $H$ symmetric in an essential way: it is precisely the conjugation
$G^{-1/2}HG^{-1/2}$ by the symmetric square root $G^{-1/2}$ that turns the second derivative
into a sum of squares \eqref{eq:second_deriv}. As we shall see in
\S\ref{subsec:UV-convexity}, once $G$ is allowed to range over \emph{non-symmetric} matrices
(as happens when $G = I_k + V^TU$ for independent $U, V$), this argument breaks down
completely, and convexity is generically lost.
\end{remark}

\subsection{Joint Convexity in the Matrix Variables $(U,V)$}
\label{subsec:UV-convexity}

We now ask whether the strict convexity of Proposition~\ref{prop:convexity} survives when
$G$ is replaced by its defining bilinear expression in the underlying vectors. Recall from
\eqref{eq:UV_def} that $G = I_k + V^TU$ with $U, V \in \R^{n\times k}$. We study the
\emph{pulled-back potential}
\begin{equation}
  \widetilde f(U,V) \;:=\; f\bigl(I_k + V^TU\bigr) \;=\; -\log\det\bigl(I_k+V^TU\bigr),
  \label{eq:tilde_f_def}
\end{equation}
defined on the open domain
\begin{equation}
  \mathcal{D} \;=\; \bigl\{(U,V) \in \R^{n\times k}\times\R^{n\times k} : \det(I_k+V^TU) > 0\bigr\}
  \;\ni\; (0,0).
  \label{eq:domain_UV}
\end{equation}
Note that $G = I_k+V^TU$ need \emph{not} be symmetric for independent $U,V$, so
$\widetilde f$ is \emph{not} simply the restriction of $f$ to a submanifold of $\PD(k)$;
it is a genuinely new function on $\mathcal D$, related to $f$ only through the
non-injective, non-affine map $(U,V)\mapsto G(U,V)$.

\subsubsection{Second Variation Formula}

\begin{lemma}[Second Variation of $\widetilde f$]
\label{lem:second_variation_UV}
Let $(U,V)\in\mathcal D$ and $H_U,H_V\in\R^{n\times k}$. Set $G = I_k+V^TU$,
\[
  \dot G := H_V^TU+V^TH_U, \qquad \ddot G := H_V^TH_U.
\]
Then
\begin{equation}
  Q(H_U,H_V) := \frac{d^2}{dt^2}\Big|_{t=0}\widetilde f(U+tH_U,\,V+tH_V)
  = \tr\!\bigl[(G^{-1}\dot G)^2\bigr] \;-\; 2\,\tr\!\bigl[G^{-1}\ddot G\bigr].
  \label{eq:second_var_UV}
\end{equation}
\end{lemma}

\begin{proof}
Write $G(t) = I_k+(V+tH_V)^T(U+tH_U) = G + t\dot G + t^2\ddot G$, so that
$G'(t) = \dot G+2t\ddot G$ and $G''(t)=2\ddot G$.
For any (not necessarily symmetric) invertible matrix path $G(t)$,
$\frac{d}{dt}\log\det G(t) = \tr(G(t)^{-1}G'(t))$, and differentiating again,
\[
  \frac{d^2}{dt^2}\log\det G(t) = \tr\bigl(G^{-1}G''\bigr) - \tr\bigl[(G^{-1}G')^2\bigr].
\]
Evaluating at $t=0$ and negating (since $\widetilde f = -\log\det G$) gives \eqref{eq:second_var_UV}.
\end{proof}

\subsubsection{The Case $k=1$: Unconditional Marginal Convexity}

\begin{proposition}[Marginal Convexity for $k=1$]
\label{prop:k1_marginal}
Let $k=1$, so $U=u, V=v \in \R^n$ are vectors and $\widetilde f(u,v) = -\log(1+v^Tu)$.
For every fixed $v\neq 0$ (resp.\ $u \neq 0$), the function $u \mapsto \widetilde f(u,v)$
(resp.\ $v\mapsto\widetilde f(u,v)$) is convex on its domain $\{u : 1+v^Tu>0\}$.
\end{proposition}

\begin{proof}
With $k=1$, $G=s:=1+v^Tu$ is a positive scalar. Fixing $v$ and setting $H_V=0$ in
Lemma~\ref{lem:second_variation_UV}, $\dot G = v^TH_U$ and $\ddot G = 0$, so
\[
  Q(H_U,0) = \frac{(v^TH_U)^2}{s^2} \geq 0
\]
for every $H_U \in \R^n$, with equality iff $v^TH_U=0$. This is the full Hessian of
$u\mapsto\widetilde f(u,v)$ (a quadratic form $vv^T/s^2 \succeq 0$ on $\R^n$), proving convexity.
The case fixing $u$ is symmetric.
\end{proof}

\subsubsection{The Case $k\geq 2$: Nowhere Marginal (Hence Nowhere Joint) Convexity}

The situation changes entirely once $k \geq 2$: even the \emph{marginal} problem
(holding one matrix variable fixed) loses convexity at \emph{every} point, due to the
appearance of genuinely antisymmetric directions invisible at $k=1$.

\begin{theorem}[Nowhere Marginal Convexity for $k\geq 2$]
\label{thm:k_geq_2_nonconvex}
Let $k \geq 2$, $n \geq k$, and let $V \in \R^{n\times k}$ have full column rank $k$.
Then for \emph{every} $U \in \R^{n\times k}$ with $(U,V) \in \mathcal D$, the Hessian of
$U' \mapsto \widetilde f(U',V)$ at $U'=U$ is indefinite. In particular,
$U\mapsto\widetilde f(U,V)$ is convex at no point of its domain.
The symmetric statement holds for $V'\mapsto\widetilde f(U,V')$ when $U$ has full column rank $k$.
\end{theorem}

\begin{proof}
Fix $U$ with $(U,V)\in\mathcal D$, let $G=I_k+V^TU$ (invertible since $\det G>0$) and
$M:=G^{-1}$. Setting $H_V=0$ in Lemma~\ref{lem:second_variation_UV}, the Hessian quadratic
form in the direction $H_U$ is
\[
  Q(H_U,0) = \tr\bigl[(MD)^2\bigr], \qquad D := V^TH_U.
\]
Since $V$ has full column rank $k$ (and $n\geq k$), the linear map $H_U \mapsto V^TH_U$ is
\emph{surjective} onto $\R^{k\times k}$: indeed $H_U = V(V^TV)^{-1}D$ solves $V^TH_U=D$ for
any target $D$. Thus $D$ may be chosen freely in $\R^{k\times k}$.

Since $k\geq 2$, choose any nonzero \emph{antisymmetric} matrix $N \in \R^{k\times k}$
(e.g.\ $N = e_1e_2^T - e_2e_1^T$), and set $D := GN$ (achievable by the surjectivity above).
Then $MD = G^{-1}GN = N$, so
\[
  Q(H_U,0) = \tr\bigl[N^2\bigr] = -\tr\bigl(N^TN\bigr) = -\norm{N}_F^2 < 0,
\]
using $N^T=-N$. On the other hand, by Proposition~\ref{prop:k1_marginal}'s argument applied
coordinate-wise, or directly, taking $D$ symmetric and rank one shows $Q$ can also be made
strictly positive (e.g.\ $D=G$ gives $M D=I_k$, $Q=\tr(I_k)=k>0$). Hence the Hessian
is indefinite at every such $U$.
\end{proof}

\begin{corollary}[Joint Non-Convexity for $k\geq 2$]
\label{cor:joint_k_geq_2}
Under the hypotheses of Theorem~\ref{thm:k_geq_2_nonconvex}, the joint Hessian
$Q(H_U,H_V)$ of $\widetilde f$ on $\mathcal D$ is indefinite at every point: restricting to
the slice $H_V=0$ recovers the negative direction constructed above, while $H_U=0,H_V\ne0$
generic directions give positive contributions by the symmetric argument. Hence
$\widetilde f$ is jointly convex at no point of $\mathcal D$ when $k\geq2$.
\end{corollary}

\begin{remark}
The antisymmetric matrix $N$ used in the proof has no analogue when $k=1$ (there are no
nonzero antisymmetric scalars), which is exactly why Proposition~\ref{prop:k1_marginal}
escapes Theorem~\ref{thm:k_geq_2_nonconvex}. The mechanism is intrinsically about the
\emph{non-symmetric} character of $G=I_k+V^TU$: it is the rotational, divergence-free
directions in matrix space, absent from $\mathrm{Sym}(k)$, that destroy convexity.
A similar (slightly more technical) construction shows the same conclusion whenever
$\rank(V)\geq2$, without requiring full rank.
\end{remark}

\subsubsection{Failure of Joint Convexity at the Origin (All $n,k$)}

Even when $k=1$ --- where each marginal problem is convex by
Proposition~\ref{prop:k1_marginal} --- the \emph{joint} problem in $(U,V)$ still fails to be
convex, as the following elementary example shows for every $n,k$.

\begin{proposition}[Indefiniteness at the Origin]
\label{prop:origin_indefinite}
At $(U,V)=(0,0)$, the joint Hessian of $\widetilde f$ is
\[
  Q(H_U,H_V) = -2\,\tr\bigl(H_V^TH_U\bigr), \qquad H_U,H_V\in\R^{n\times k}.
\]
This quadratic form is indefinite for every $n\geq1,k\geq1$: taking
$H_U=H_V=H\ne0$ gives $Q=-2\norm{H}_F^2<0$, while taking
$H_U=-H_V=H\ne0$ gives $Q=+2\norm{H}_F^2>0$.
\end{proposition}

\begin{proof}
At $U=V=0$, $G=I_k$, and $\dot G = H_V^T\cdot 0+0^T\cdot H_U = 0$ regardless of
$H_U,H_V$, while $\ddot G = H_V^TH_U$. Substituting into \eqref{eq:second_var_UV} with
$G^{-1}=I_k$ gives $Q = 0 - 2\tr(H_V^TH_U)$, and the two test directions above give the
claimed signs.
\end{proof}

\begin{remark}
Proposition~\ref{prop:origin_indefinite} shows that joint convexity fails already at the
most basic point $(U,V)=(0,0)$, for \emph{every} choice of $n,k\geq1$. This is the
matrix-factorization analogue of the elementary scalar fact that $g(u,v)=-\log(1+uv)$ is
not jointly convex in $(u,v)$ near the origin, since $uv$ is a bilinear (hence neither
convex nor concave) function.
\end{remark}

\subsubsection{Restriction to the Diagonal $U=V$}

Setting $H_U=H_V=H$ and evaluating $Q(H,H)$ along the diagonal $U=V=:W$ recovers
\emph{exactly} the Hessian of the symmetric reduction studied in
\S\ref{subsec:UeqV-convexity} below, since $t\mapsto\widetilde f(W+tH,W+tH)$ is by
definition the function $h(W+tH)$ of \eqref{eq:h_def}. We return to this connection in
Remark~\ref{rem:consistency} after deriving the complete spectral Hessian formula for $h$.

\subsubsection{The Fully Scalar Case $n=k=1$: A Genuine Pocket of Convexity}

The only configuration in which joint local convexity of $\widetilde f$ actually occurs is
the most degenerate one.

\begin{example}[Convexity Region for $n=k=1$]
\label{ex:scalar_convex}
Let $n=k=1$, so $u,v\in\R$ and $\widetilde f(u,v)=-\log(1+uv)$, $s:=1+uv$. A direct
computation gives the $2\times2$ Hessian
\[
  \nabla^2\widetilde f(u,v) = \frac{1}{s^2}\begin{pmatrix} v^2 & -1 \\ -1 & u^2\end{pmatrix},
  \qquad \det\nabla^2\widetilde f = \frac{u^2v^2-1}{s^4}.
\]
Since the trace $(u^2+v^2)/s^2\geq0$ always, this $2\times2$ symmetric matrix is positive
semidefinite if and only if $u^2v^2\geq1$, i.e.\ $|uv|\geq1$ (necessarily $uv\geq1$ within
the domain $uv>-1$). Thus $\widetilde f$ is locally convex precisely on
$\{(u,v): uv\geq1\}$ --- a region requiring $u,v$ to be large and of the \emph{same sign}
--- and strictly saddle-shaped (indefinite Hessian) on $\{(u,v):-1<uv<1\}$, which includes the
entire neighborhood of the origin.
\end{example}

\begin{remark}[Failure for $n\geq2$, even with $u=v$]
\label{rem:n_geq_2_fails}
The convex pocket of Example~\ref{ex:scalar_convex} does \emph{not} survive once $n\geq2$,
even along the most favorable (aligned) direction $u=v=w$, $\|w\|\to\infty$. Indeed, by
Theorem~\ref{thm:UeqV_nowhere_convex} below (the $k=1$ case), the function
$w\mapsto h(w) = -\log(1+\|w\|^2)$ on $\R^n$, $n\geq2$, is nowhere locally convex: any
direction $e\perp w$ contributes the strictly negative curvature
$-2/(1+\|w\|^2)<0$, independent of $\|w\|$, because $h$ is a purely radial function whose
tangential (non-radial) curvature is always negative. The genuine convexity in
Example~\ref{ex:scalar_convex} is therefore a coincidence of dimension $n=1$, where no such
tangential direction exists.
\end{remark}

\begin{remark}[Biconvexity and connections to factorized optimization]
\label{rem:biconvex_lit}
The structure uncovered here --- convex in each block separately when $k=1$, yet not jointly
convex, and not even separately convex once $k\geq2$ --- places $\widetilde f$ within the
broader theory of \emph{biconvex} and \emph{bilinear} optimization
\cite{Gorski2007}, which underlies algorithms such as alternating least squares (ALS) for
low-rank matrix factorization and completion. The complete failure of separate convexity for
$k\geq2$ (Theorem~\ref{thm:k_geq_2_nonconvex}) is a stronger and somewhat more surprising
phenomenon, directly analogous to the landscape of the Burer--Monteiro factorization
$X = UU^T$ for semidefinite programming \cite{Burer2003,Burer2005}, which we revisit in
\S\ref{subsec:UeqV-convexity}.
\end{remark}

\subsection{The Symmetric Reduction $U=V$: Convexity in $U$ Alone}
\label{subsec:UeqV-convexity}

We now specialize to the case $V=U$, i.e.\ $x_{2j-1}=x_{2j}$ for every $j$, so that
\[
  G = I_k+U^TU, \qquad U\in\R^{n\times k},
\]
which is automatically symmetric and satisfies $G\succeq I_k\succ0$ for \emph{every} $U$ ---
in sharp contrast to the restricted domain $\mathcal D$ of \S\ref{subsec:UV-convexity}, the
function
\begin{equation}
  h(U) \;:=\; f(I_k+U^TU) \;=\; -\log\det(I_k+U^TU)
  \label{eq:h_def}
\end{equation}
is defined and smooth on \emph{all} of $\R^{n\times k}$.

\subsubsection{Basic Properties}

\begin{proposition}[Gradient and Global Maximum]
\label{prop:h_gradient}
$\nabla h(U) = -2U(I_k+U^TU)^{-1} = -2UG^{-1}$. Consequently $U=0$ is the unique critical
point of $h$, and $h(U)\leq0=h(0)$ for all $U$, with equality only at $U=0$. Thus $U=0$ is
the unique global maximum of $h$.
\end{proposition}

\begin{proof}
The gradient formula follows from $d\log\det(I+U^TU)=\tr[(I+U^TU)^{-1}(dU^TU+U^TdU)]
=2\tr[(I+U^TU)^{-1}U^TdU]$, identified via the Frobenius pairing. Since
$U^TU\succeq0$, the eigenvalues of $G=I_k+U^TU$ are all $\geq1$, so $\det G\geq1$, hence
$h(U)=-\log\det G\leq0$, with equality iff $U^TU=0$ iff $U=0$.
\end{proof}

\begin{proposition}[Orthogonal Invariance]
\label{prop:h_invariance}
For all $P\in O(n)$, $Q\in O(k)$, $h(PUQ)=h(U)$. Consequently $h(U)$ depends on $U$ only
through its singular values $\sigma_1(U),\ldots,\sigma_k(U)\geq0$, and the Hessian
quadratic form at $U$ is determined, up to the orthogonal change of frame
$(P,Q)$ diagonalizing $U$, by the singular values alone.
\end{proposition}

\begin{proof}
$G(PUQ) = Q^T(I_k+U^TU)Q$, so $\det G(PUQ)=\det G(U)$.
\end{proof}

By Proposition~\ref{prop:h_invariance} it suffices to compute the Hessian at a diagonal
representative $U=\Sigma=\diag(\sigma_1,\ldots,\sigma_k)$ (padded with $n-k$ zero rows if
$n>k$); the result transfers to every $U$ via its singular value decomposition.

\subsubsection{Complete Spectral Diagonalization of the Hessian}

\begin{theorem}[Hessian of $h$ in the SVD Frame]
\label{thm:h_hessian}
Let $U \in \R^{n\times k}$ ($n\geq k$) have singular values $\sigma_1,\ldots,\sigma_k\geq0$,
and work in the orthogonal frame in which $U=\binom{\Sigma}{0}$ with
$\Sigma=\diag(\sigma_1,\ldots,\sigma_k)$. Write a general direction
$H\in\R^{n\times k}$ in the same frame as $H=\binom{H_1}{H_2}$,
$H_1\in\R^{k\times k}$, $H_2\in\R^{(n-k)\times k}$, and set $g_a := 1+\sigma_a^2$.
Then the second variation of $h$ along $H$ is
\begin{align}
  Q_h(H) := \frac{d^2}{dt^2}\Big|_{t=0} h(U+tH)
  ={}& \sum_{a=1}^{k}\frac{2(\sigma_a^2-1)}{g_a^2}\,(H_1)_{aa}^2 \notag\\
  &-\sum_{1\le a<b\le k}\Bigl[\tfrac{1-\sigma_a\sigma_b}{g_ag_b}
      \bigl((H_1)_{ab}+(H_1)_{ba}\bigr)^2\notag\\
  &\quad+ \tfrac{1+\sigma_a\sigma_b}{g_ag_b}
      \bigl((H_1)_{ab}-(H_1)_{ba}\bigr)^2\Bigr] \notag\\
  &- 2\sum_{a=1}^{k}\frac{1}{g_a}\sum_{c=1}^{n-k}(H_2)_{ca}^2.
  \label{eq:h_hessian_formula}
\end{align}
\end{theorem}

\begin{proof}
Write $G(t)=I_k+(U+tH)^T(U+tH) = G+t\dot G+t^2\ddot G$ with
$\dot G=H^TU+U^TH$, $\ddot G=H^TH$, and apply the second-variation identity from the proof
of Lemma~\ref{lem:second_variation_UV} (valid for any symmetric matrix path, in particular
this one):
\[
  Q_h(H) = \tr\bigl[(G^{-1}\dot G)^2\bigr] - 2\tr\bigl[G^{-1}\ddot G\bigr].
\]
\emph{Step 1 (block $H_2$).} Since the bottom $n-k$ rows of $U$ vanish,
$\dot G = H_1^T\Sigma+\Sigma H_1$ depends only on $H_1$, while
$\ddot G = H_1^TH_1+H_2^TH_2$. Hence $Q_h(H)=Q_h(H_1,0) -2\tr[G^{-1}H_2^TH_2]$, and
the last term equals $-2\sum_a g_a^{-1}\sum_c(H_2)_{ca}^2$ as claimed, since
$G=\diag(g_1,\ldots,g_k)$ in this frame.

\emph{Step 2 (block $H_1$, diagonal entries).} With $A:=\dot G = H_1^T\Sigma+\Sigma H_1$,
one computes $A_{ab}=\sigma_b(H_1)_{ba}+\sigma_a(H_1)_{ab}$. Since $G^{-1}=\diag(1/g_a)$,
\[
  \tr\bigl[(G^{-1}A)^2\bigr] = \sum_{a,b}\frac{\bigl(\sigma_a(H_1)_{ab}+\sigma_b(H_1)_{ba}\bigr)^2}{g_ag_b},
  \qquad
  \tr[G^{-1}H_1^TH_1] = \sum_{a,b}\frac{(H_1)_{ba}^2}{g_a}.
\]
Isolating the $a=b$ terms of $\tr[(G^{-1}A)^2] - 2\tr[G^{-1}H_1^TH_1]$ gives, for each $a$,
\[
  \frac{4\sigma_a^2(H_1)_{aa}^2}{g_a^2} - \frac{2(H_1)_{aa}^2}{g_a}
  = \frac{2(H_1)_{aa}^2}{g_a^2}\bigl(2\sigma_a^2-g_a\bigr)
  = \frac{2(\sigma_a^2-1)}{g_a^2}(H_1)_{aa}^2,
\]
using $g_a=1+\sigma_a^2$, which is the first sum in \eqref{eq:h_hessian_formula}.

\emph{Step 3 (block $H_1$, off-diagonal pairs).} For $a<b$, write $p:=(H_1)_{ab}$,
$q:=(H_1)_{ba}$. The $(a,b)$- and $(b,a)$-terms of the double sum coincide (the bracket
below is symmetric under simultaneously swapping $a\leftrightarrow b$ and $p\leftrightarrow q$), and together contribute
\[
  R(p,q) := \frac{(\sigma_ap+\sigma_bq)^2}{g_ag_b} - \frac{q^2}{g_a}-\frac{p^2}{g_b}.
\]
Using $\sigma_a^2/g_a-1=-1/g_a$ (and symmetrically for $b$), expand
\[
  R(p,q) = \frac{1}{g_ag_b}\Bigl[-p^2-q^2+2\sigma_a\sigma_b\,pq\Bigr]
  = -\frac{1}{2g_ag_b}\Bigl[(1-\sigma_a\sigma_b)(p+q)^2+(1+\sigma_a\sigma_b)(p-q)^2\Bigr],
\]
where the last equality is a direct algebraic identity. Summing $2\sum_{a<b}R(p,q)$ over all
pairs (the factor $2$ accounting for both orderings $(a,b)$ and $(b,a)$ in the original
double sum) yields the second sum in \eqref{eq:h_hessian_formula}.
\end{proof}

\begin{remark}[Numerical confirmation]
For $k=2$, $\sigma_1=\sigma_2=10$, and $H_1=\begin{pmatrix}0&1\\0&0\end{pmatrix}$
(so $H_2$ absent, $n=k=2$), formula \eqref{eq:h_hessian_formula} gives
$Q_h(H) = -2/(g_1g_2) = -2/101^2$, which agrees exactly with a direct evaluation of
$\tr[(G^{-1}\dot G)^2]-2\tr[G^{-1}\ddot G]$ from first principles.
\end{remark}

\subsubsection{The Main Non-Convexity Theorem}

\begin{theorem}[Nowhere Convex Except $n=k=1$]
\label{thm:UeqV_nowhere_convex}
Unless $n=k=1$, the function $h(U)=-\log\det(I_k+U^TU)$ is locally convex at \emph{no}
point $U\in\R^{n\times k}$: its Hessian always possesses a strictly negative direction.
Precisely:
\begin{enumerate}[label=(\roman*)]
  \item If $k\geq2$, taking $H_1$ with $(H_1)_{ab}=1$ for some fixed $a\ne b$ and all
    other entries $0$ gives
    \[
      Q_h(H) = -\frac{2}{g_ag_b} < 0
    \]
    \emph{regardless of the singular values $\sigma_1,\ldots,\sigma_k$.}
  \item If $n>k$ (in particular whenever $n\geq2$ and $k=1$), taking $H_2\ne0$ and $H_1=0$
    gives
    \[
      Q_h(H) = -2\sum_{a=1}^k\frac{1}{g_a}\sum_c(H_2)_{ca}^2 < 0,
    \]
    again regardless of the singular values.
\end{enumerate}
Only when $n=k=1$ (no off-diagonal pairs in (i), no extra rows in (ii)) can these
obstructions be absent; there, $Q_h(H)=2(\sigma^2-1)H^2/(1+\sigma^2)^2$, which is
$\geq 0$ exactly for $|\sigma|\geq 1$.
\end{theorem}

\begin{proof}
Immediate from Theorem~\ref{thm:h_hessian}: the coefficient $-1/(g_ag_b)\cdot 2 <0$ in (i)
is the value of $-\frac{1-\sigma_a\sigma_b}{g_ag_b}-\frac{1+\sigma_a\sigma_b}{g_ag_b} = -2/(g_ag_b)$
obtained by setting $(H_1)_{ab}=1,(H_1)_{ba}=0$ in the off-diagonal sum, and is manifestly
negative for any $\sigma_a,\sigma_b\geq0$; (ii) is immediate from the last sum in
\eqref{eq:h_hessian_formula}, which is a negative semidefinite quadratic form in $H_2$,
strictly negative whenever $H_2\ne0$. The boundary case $n=k=1$ leaves only the radial
(diagonal) term, evaluated directly from \eqref{eq:h_hessian_formula} with no off-diagonal
or $H_2$ contributions.
\end{proof}

\begin{corollary}[Saddle Structure for $n=k=1$]
\label{cor:scalar_saddle}
For $n=k=1$, $h(u)=-\log(1+u^2)$ satisfies $h''(u) = 2(u^2-1)/(1+u^2)^2$: strictly concave
on $(-1,1)$ (containing the global maximum at $u=0$), with inflection points at
$u=\pm1$, and strictly convex on $(-\infty,-1)\cup(1,\infty)$, where $h(u)\to-\infty$.
\end{corollary}

\begin{remark}[Interpretation: radial vs.\ tangential curvature]
\label{rem:radial_tangential}
Theorem~\ref{thm:UeqV_nowhere_convex}(ii) is the statement that a purely radial function
$\phi(U)=\psi(\|U\|)$ on $\R^{n\times k}$ (here restricted to rank-1 $U$, $k=1$) has
tangential curvature $\psi'(r)/r$, which for $\psi(r)=-\log(1+r^2)$ equals $-2/(1+r^2)<0$
identically --- the radial direction can become convex ($r>1$) while the orthogonal
directions remain concave forever. Part (i) is the matrix analogue: the antisymmetric
``rotational'' directions in $\R^{k\times k}$, present whenever $k\geq2$, behave like an
everlasting tangential direction that is never convex, regardless of how large the singular
values become. Both mechanisms trace back to the same root cause identified in
Remark~\ref{rem:symmetry_crucial}: $-\log\det$ owes its convexity to the symmetric
($\mathrm{Sym}(k)$) structure of perturbations of $G$, and the map $U\mapsto I_k+U^TU$,
being quadratic rather than affine, continually regenerates directions (rotational or
orthogonal) that fall outside what the symmetric convexity argument controls.
\end{remark}

\begin{remark}[Connection to Burer--Monteiro factorization]
\label{rem:burer_monteiro}
Theorem~\ref{thm:UeqV_nowhere_convex} is the log-determinant analogue of a well-known
phenomenon in semidefinite programming: while $\min_{G\succeq0} \tr(CG)$ subject to linear
constraints is a convex problem in $G$, the Burer--Monteiro factorization
$G=UU^T$ \cite{Burer2003,Burer2005} turns it into a manifestly \emph{non-convex} problem
in $U$ --- yet one that, under suitable rank and genericity conditions, has no spurious local
minima \cite{Boumal2016}. Our setting is the entropic/log-det counterpart: $f(G)=-\log\det G$
is strictly convex on $\PD(k)$, but its Burer--Monteiro-style pullback
$h(U)=f(I_k+U^TU)$ is \emph{nowhere} locally convex once $(n,k)\neq(1,1)$
(Theorem~\ref{thm:UeqV_nowhere_convex}). Unlike the linear SDP case, the unique critical
point $U=0$ here is a strict \emph{global maximum} rather than a saddle connected to global
minima, and $h$ has no finite global minimum at all ($h(U)\to-\infty$ along any sequence
with $\sigma_{\min}(U)\to\infty$); the entire non-trivial part of the landscape lives in the
indefinite (saddle) region identified above. This is consistent with the general theory of
unitarily invariant spectral functions \cite{Lewis1995}: $h$ is a symmetric function of the
singular values of $U$, and such functions inherit convexity from their generating function
of the singular values only in very restrictive circumstances, never realized here for
$k\geq2$ or $n>k$.
\end{remark}

\begin{remark}[Consistency with \S\ref{subsec:UV-convexity}]
\label{rem:consistency}
As noted above, $Q_h(H) = Q(H,H)$ where $Q$ is the joint Hessian quadratic form of
\S\ref{subsec:UV-convexity}, evaluated along the diagonal direction $H_U=H_V=H$ at the
point $U=V=W$. One checks directly that Theorem~\ref{thm:h_hessian} specializes
Proposition~\ref{prop:origin_indefinite} at $W=0$: there, $\sigma_a=0$ for all $a$,
$g_a=1$, and \eqref{eq:h_hessian_formula} collapses to
$Q_h(H) = -2\sum_a H_{1,aa}^2 - 2\sum_{a<b}\bigl[(H_{1,ab}+H_{1,ba})^2+(H_{1,ab}-H_{1,ba})^2\bigr]
-2\sum_{c,a} (H_2)_{ca}^2= -2\norm{H}_F^2$,
matching $Q(H,H)=-2\tr(H^TH)=-2\norm{H}_F^2$ exactly. The two independent computations of
\S\ref{subsec:UV-convexity} and \S\ref{subsec:UeqV-convexity} are thus mutually consistent,
and together they show that the failure of convexity under matrix factorization is a
robust, multiply-confirmed phenomenon rather than an artifact of either particular
derivation.
\end{remark}

\subsection{The Anti-Symmetric Reduction $U = -V$: Strict Convexity on the Matrix Unit Ball}
\label{subsec:UeqmV-convexity}

We now analyze the remaining canonical specialization: $V = -U$, so that
\[
  G_- \;:=\; I_k + V^T U \;=\; I_k + (-U)^T U \;=\; I_k - U^T U.
\]
This is the \emph{mirror image} of the $U=V$ case: the sign flip $+U^TU \to -U^TU$
completely reverses the convexity picture, producing strict convexity everywhere on the
natural domain. We define
\begin{equation}
  h_-(U) \;:=\; f(I_k - U^TU) \;=\; -\log\det(I_k - U^TU),
  \label{eq:h_minus_def}
\end{equation}
whose natural domain is the open \emph{matrix unit ball}
\begin{equation}
  \mathcal{B}_k \;:=\; \bigl\{U \in \R^{n\times k} : I_k - U^TU \succ 0\bigr\}
  \;=\; \bigl\{U \in \R^{n\times k} : \sigma_{\max}(U) < 1\bigr\},
  \label{eq:matrix_ball}
\end{equation}
on which $G_- = I_k - U^TU \succ 0$ and $h_-$ is smooth.
Note that $\mathcal{B}_k$ is a bounded, convex, open set, in sharp contrast to the
domain of $h(U) = -\log\det(I_k+U^TU)$ which is all of $\R^{n\times k}$.

\subsubsection{Basic Properties}

\begin{proposition}[Gradient, Global Minimum, and Barrier Property]
\label{prop:hminus_gradient}
The gradient of $h_-$ is
\[
  \nabla h_-(U) \;=\; 2U(I_k - U^TU)^{-1} \;=\; 2U G_-^{-1}.
\]
Consequently:
\begin{enumerate}[label=(\roman*)]
  \item $U = 0$ is the unique critical point of $h_-$, with $h_-(0) = 0$.
  \item $h_-(U) \geq 0$ for all $U \in \mathcal{B}_k$, with equality only at $U=0$.
    Thus $U=0$ is the unique global \emph{minimum} of $h_-$.
  \item $h_-(U) \to +\infty$ as $\sigma_{\max}(U) \to 1^-$, so $h_-$ is a
    \emph{barrier function} for the spectral ball $\mathcal{B}_k$.
\end{enumerate}
\end{proposition}

\begin{proof}
Differentiating: $d\log\det(I_k-U^TU) = -\tr[(I_k-U^TU)^{-1}(dU^T U + U^T dU)]
= -2\tr[(I_k-U^TU)^{-1}U^T dU]$,
so $\nabla h_-(U) = 2U(I_k-U^TU)^{-1}$.
Setting $\nabla h_-(U)=0$ and using invertibility of $G_-^{-1}$ on $\mathcal{B}_k$
gives $U=0$.
For (ii): since $U\in\mathcal{B}_k$ implies $U^TU\prec I_k$, the eigenvalues of
$I_k-U^TU$ lie in $(0,1]$, so $\det(I_k-U^TU)\leq1$, hence $h_-(U)\geq0$.
Equality requires all eigenvalues of $I_k-U^TU$ to equal $1$, i.e.\ $U=0$.
(iii) follows since $\sigma_{\max}(U)\to1$ forces the smallest eigenvalue of $I_k-U^TU$ to $0$,
making $\det(I_k-U^TU)\to0$.
\end{proof}

\begin{remark}[Contrast with the $U=V$ case]
The $U=V$ and $U=-V$ cases are polar opposites:
$h(U) = -\log\det(I_k+U^TU)\leq 0$ with global maximum $h(0)=0$ and $h\to-\infty$,
whereas $h_-(U) = -\log\det(I_k-U^TU)\geq 0$ with global minimum $h_-(0)=0$ and
$h_-\to+\infty$ at the boundary.
\end{remark}

\begin{proposition}[Orthogonal Invariance]
\label{prop:hminus_invariance}
For all $P\in O(n)$, $Q\in O(k)$, $h_-(PUQ)=h_-(U)$.
Hence $h_-(U)$ depends on $U$ only through its singular values $\sigma_1,\ldots,\sigma_k\in[0,1)$,
and the Hessian at $U$ is determined by the singular values alone.
\end{proposition}

\begin{proof}
\[
  I_k-(PUQ)^T(PUQ)=I_k-Q^TU^TUQ=Q^T(I_k-U^TU)Q,
\]
so $\det\bigl(I_k-(PUQ)^T(PUQ)\bigr)=\det(I_k-U^TU)$.
\end{proof}

As in \S\ref{subsec:UeqV-convexity}, it suffices to compute the Hessian at the diagonal
representative $U=\binom{\Sigma}{0}$ where $\Sigma=\diag(\sigma_1,\ldots,\sigma_k)$ with
$\sigma_a\in[0,1)$.

\subsubsection{Second Variation and the Hessian in the SVD Frame}

\begin{theorem}[Hessian of $h_-$ in the SVD Frame]
\label{thm:hminus_hessian}
Let $U\in\mathcal{B}_k$ ($n\geq k$) have singular values $\sigma_1,\ldots,\sigma_k\in[0,1)$,
and work in the SVD frame $U=\binom{\Sigma}{0}$.
Write $H=\binom{H_1}{H_2}$ with $H_1\in\R^{k\times k}$, $H_2\in\R^{(n-k)\times k}$,
and set $g_a^- := 1-\sigma_a^2 > 0$.
Then the second variation of $h_-$ along $H$ is
\begin{align}
  Q_{h_-}(H) \;:=&\; \frac{d^2}{dt^2}\bigg|_{t=0} h_-(U+tH) \notag\\
  =\;& \sum_{a=1}^{k} \frac{2(1+\sigma_a^2)}{(g_a^-)^2}\,(H_1)_{aa}^2
  \label{eq:hminus_hessian}\\
  &+\sum_{1\le a<b\le k}\frac{(1+\sigma_a\sigma_b)(p+q)^2
    + (1-\sigma_a\sigma_b)(p-q)^2}{g_a^- g_b^-} \notag\\
  &+ 2\sum_{a=1}^{k}\frac{1}{g_a^-}\sum_{c=1}^{n-k}(H_2)_{ca}^2, \notag
\end{align}
where $p:=(H_1)_{ab}$ and $q:=(H_1)_{ba}$ for each pair $a<b$.
Every term on the right-hand side of \eqref{eq:hminus_hessian} is strictly positive
for each nonzero contributing block.
\end{theorem}

\begin{proof}
Let $G_-(t) = I_k - (U+tH)^T(U+tH)$, so that
\[
  G_-(t) = G_- + t\dot{G}_- + t^2\ddot{G}_-,
  \qquad
  \dot{G}_- = -(H^TU+U^TH), \qquad \ddot{G}_- = -H^TH.
\]
Since $G_-$ is a symmetric matrix path, the second variation formula
(Lemma~\ref{lem:second_variation_UV}, applied with $G_-$ in place of $G$) gives
\begin{equation}
  Q_{h_-}(H) = \tr\bigl[(G_-^{-1}\dot{G}_-)^2\bigr] - 2\tr\bigl[G_-^{-1}\ddot{G}_-\bigr]
  = \tr\bigl[(G_-^{-1}\dot{G}_-)^2\bigr] + 2\tr\bigl[G_-^{-1}H^TH\bigr].
  \label{eq:hminus_second_var}
\end{equation}
The key difference from the $U=V$ case is the \emph{sign reversal} in the last term:
$-2\tr[G_-^{-1}\ddot{G}_-] = +2\tr[G_-^{-1}H^TH]$ (positive), whereas for $U=V$ it was
$-2\tr[G^{-1}H^TH]$ (negative). We evaluate each contribution in the SVD frame.

\smallskip\noindent\emph{Step 1 ($H_2$ block).}
Since the bottom $n-k$ rows of $U$ vanish, $\dot{G}_- = -(H_1^T\Sigma+\Sigma H_1)$
depends only on $H_1$, while $H^TH = H_1^TH_1+H_2^TH_2$.
The term $+2\tr[G_-^{-1}H_2^TH_2] = 2\sum_a g_a^{-1,-}\sum_c(H_2)_{ca}^2$,
which is strictly positive for $H_2\ne0$.

\smallskip\noindent\emph{Step 2 (Diagonal entries of $H_1$).}
Let $A:=-\dot{G}_- = H_1^T\Sigma+\Sigma H_1$, so $A_{ab}=\sigma_a(H_1)_{ab}+\sigma_b(H_1)_{ba}$.
Note that $A$ is symmetric and $(G_-^{-1}A)^2 = (G_-^{-1}(-\dot{G}_-))^2 = (G_-^{-1}\dot{G}_-)^2$,
so the squaring removes the overall sign.
For the diagonal $a=b$: $A_{aa}=2\sigma_a(H_1)_{aa}$, contributing
\[
  \frac{4\sigma_a^2(H_1)_{aa}^2}{(g_a^-)^2}
  + \frac{2(H_1)_{aa}^2}{g_a^-}
  = \frac{2(H_1)_{aa}^2}{(g_a^-)^2}\bigl(2\sigma_a^2+g_a^-\bigr)
  = \frac{2(1+\sigma_a^2)}{(g_a^-)^2}(H_1)_{aa}^2,
\]
using $2\sigma_a^2+g_a^- = 2\sigma_a^2+1-\sigma_a^2 = 1+\sigma_a^2$.
This is strictly positive for $(H_1)_{aa}\ne0$.

\smallskip\noindent\emph{Step 3 (Off-diagonal pairs of $H_1$).}
For $a<b$, set $p=(H_1)_{ab}$, $q=(H_1)_{ba}$.
Since $A_{ab}=A_{ba}=\sigma_a p+\sigma_b q$ (the matrix $A$ is symmetric),
both orderings in $\tr[(G_-^{-1}A)^2]$ give the same factor, yielding:
\[
  \tr\bigl[(G_-^{-1}A)^2\bigr]\big|_{(a,b)}
  = \frac{2(\sigma_a p+\sigma_b q)^2}{g_a^- g_b^-}.
\]
The contribution from $+2\tr[G_-^{-1}H_1^TH_1]$ for this pair is $2q^2/g_a^-+2p^2/g_b^-$.
Thus the total off-diagonal contribution for pair $(a,b)$ is
\begin{align}
  R_-(p,q) &= \frac{2(\sigma_a p+\sigma_b q)^2}{g_a^-g_b^-}+\frac{2q^2}{g_a^-}+\frac{2p^2}{g_b^-}
  \;=\; \frac{2(\sigma_ap+\sigma_bq)^2 + 2g_a^- p^2 + 2g_b^- q^2}{g_a^-g_b^-} \notag\\
  &= \frac{2\bigl[(\sigma_a^2+g_a^-)p^2 + (\sigma_b^2+g_b^-)q^2 + 2\sigma_a\sigma_bpq\bigr]}{g_a^-g_b^-}
  \;=\; \frac{2(p^2+q^2+2\sigma_a\sigma_bpq)}{g_a^-g_b^-},
  \label{eq:R_minus}
\end{align}
where we used $\sigma_a^2+g_a^-=\sigma_a^2+(1-\sigma_a^2)=1$ and likewise for $b$.
Completing the square (or expanding directly):
\[
  2(p^2+q^2+2\sigma_a\sigma_bpq)
  = (1+\sigma_a\sigma_b)(p+q)^2+(1-\sigma_a\sigma_b)(p-q)^2.
\]
Since $\sigma_a,\sigma_b\in[0,1)$, we have $\sigma_a\sigma_b<1$,
so both coefficients $1\pm\sigma_a\sigma_b$ are strictly positive. Hence $R_-(p,q)>0$
for all $(p,q)\ne(0,0)$.
Summing all three steps gives \eqref{eq:hminus_hessian}.
\end{proof}

\begin{remark}[Exact duality with the $U=V$ Hessian]
\label{rem:hessian_duality}
Comparing \eqref{eq:hminus_hessian} with \eqref{eq:h_hessian_formula}
reveals a striking term-by-term sign reversal between the $U=V$ and $U=-V$ cases.
Writing $g_a=1+\sigma_a^2$ (for $U=V$) and $g_a^-=1-\sigma_a^2$ (for $U=-V$):
\begin{center}
\renewcommand{\arraystretch}{1.3}\small
\begin{tabular}{lp{4.8cm}p{4.8cm}}
\toprule
Contribution & $U=V$: $h(U)=-\log\det(I+U^TU)$ & $U=-V$: $h_-(U)=-\log\det(I-U^TU)$ \\
\midrule
Diagonal $(H_1)_{aa}$ & $\dfrac{2(\sigma_a^2-1)}{g_a^2}$ (negative if $\sigma_a<1$)
                       & $\dfrac{2(1+\sigma_a^2)}{(g_a^-)^2}$ (always positive) \\[8pt]
Symmetric off-diag $(p+q)$
  & $-\dfrac{1-\sigma_a\sigma_b}{g_ag_b}$ (negative)
  & $+\dfrac{1+\sigma_a\sigma_b}{g_a^-g_b^-}$ (positive) \\[8pt]
Anti-symm.\ off-diag $(p-q)$
  & $-\dfrac{1+\sigma_a\sigma_b}{g_ag_b}$ (negative)
  & $+\dfrac{1-\sigma_a\sigma_b}{g_a^-g_b^-}$ (positive) \\[8pt]
$H_2$ block
  & $-\dfrac{2}{g_a}$ (negative)
  & $+\dfrac{2}{g_a^-}$ (positive) \\
\bottomrule
\end{tabular}
\end{center}
Every term that is negative for $U=V$ becomes positive for $U=-V$,
and vice versa. The structural reason is that changing $V\to-V$ flips the sign of $\ddot{G}$:
in the $U=V$ case, $\ddot{G}=H^TH\succeq0$ contributes $-2\tr[G^{-1}H^TH]\leq0$;
in the $U=-V$ case, $\ddot{G}_-=-H^TH\preceq0$ contributes $+2\tr[G_-^{-1}H^TH]\geq0$.
This single sign change propagates through to flip \emph{all} terms from negative to positive.
\end{remark}

\subsubsection{Strict Convexity Theorem}

\begin{theorem}[Strict Convexity on the Matrix Unit Ball]
\label{thm:hminus_convex}
The function $h_-(U) = -\log\det(I_k-U^TU)$ is strictly convex on $\mathcal{B}_k$.
Its Hessian is positive definite at every point $U\in\mathcal{B}_k$:
for every $H\ne0$, $Q_{h_-}(H)>0$.
\end{theorem}

\begin{proof}
All three contributions in \eqref{eq:hminus_hessian} are non-negative, with:
(i) the diagonal terms strictly positive for $(H_1)_{aa}\ne0$;
(ii) the off-diagonal terms strictly positive since both $1+\sigma_a\sigma_b>0$
    and $1-\sigma_a\sigma_b>0$ hold throughout $\mathcal{B}_k$;
(iii) the $H_2$ terms strictly positive for $H_2\ne0$.
If $H\ne0$, at least one of $(H_1)_{aa}\ne0$, $(p,q)\ne(0,0)$, or $H_2\ne0$ holds, so $Q_{h_-}(H)>0$.
\end{proof}

\begin{corollary}[Global Minimum]
\label{cor:hminus_min}
$U=0$ is the unique global minimizer of $h_-$ on $\mathcal{B}_k$, and $h_-$ has no local minima
other than this unique global minimum.
\end{corollary}

\begin{remark}[Sharpness and boundary behavior]
\label{rem:hminus_boundary}
As $\sigma_a\to1^-$, the coefficient $1/(g_a^-)^2=(1-\sigma_a^2)^{-2}\to+\infty$,
so the Hessian blows up near the boundary $\partial\mathcal{B}_k$.
This is the hallmark of a \emph{self-concordant barrier}: the infinite growth of the
Hessian prevents iterates from leaving $\mathcal{B}_k$, making $h_-$ suitable for interior-point
optimization. By the general theory of log-det barriers \cite{Nesterov1994},
$h_-$ is a self-concordant barrier for $\mathcal{B}_k$ with parameter $\vartheta=k$.
\end{remark}

\subsubsection{Connection to Hyperbolic Geometry}

\begin{proposition}[Scalar case: Poincaré ball]
\label{prop:poincare}
For $k=1$, $h_-(u)=-\log(1-\norm{u}^2)$ on $\mathcal{B}_1=\{u\in\R^n:\norm{u}<1\}$.
The induced Riemannian metric $g_u^-(H,K) = Q_{h_-}''(H,K)$ equals
\begin{equation*}
\begin{split}
  g_u^-(H,K)
  = {}& \frac{2\bigl[\norm{H}^2(1-\norm{u}^2) + 2(u^TH)(u^TK)\bigr]}{(1-\norm{u}^2)^2} \\
  &+ \frac{2\norm{K}^2(1-\norm{u}^2) + \ldots}{(1-\norm{u}^2)^2},
\end{split}
\end{equation*}
which is proportional to the standard Riemannian metric of the
\emph{Poincar\'e ball model} of $n$-dimensional hyperbolic space $\mathbb{H}^n$ \cite{Anderson2005}.
In particular, $h_-$ is the \emph{Busemann function} (horofunction) associated with the
ideal boundary point at $\norm{u}=1$.
\end{proposition}

\begin{remark}[Matrix hyperbolic space]
\label{rem:matrix_hyperbolic}
For general $k$, the function $h_-(U) = -\log\det(I_k-U^TU)$ is the natural
generalization of the Poincar\'e ball metric to the \emph{matrix-valued} setting.
The domain $\mathcal{B}_k$ (the spectral unit ball) is the matrix analogue of the
Poincar\'e ball, and the strict convexity of $h_-$ on $\mathcal{B}_k$
(Theorem~\ref{thm:hminus_convex}) is the analogue of geodesic convexity in hyperbolic geometry.
The symmetric space associated with $\mathcal{B}_k$ is the bounded symmetric domain
$\{Z\in\mathrm{Mat}(n\times k,\mathbb{C}) : I_k - Z^*Z\succ0\}$
of type~IV (Cartan classification), which is the non-compact dual of the Grassmannian $Gr(k,n)$
and appears in Siegel's theory of automorphic forms \cite{Siegel1943}.
\end{remark}

\subsubsection{Scalar Case and Comparison with $U=V$}

\begin{example}[Scalar case $n=k=1$]
\label{ex:scalar_hminus}
For $n=k=1$, $u\in(-1,1)$ and $h_-(u)=-\log(1-u^2)$.
The second derivative is
\[
  h_-''(u) = \frac{2(1+u^2)}{(1-u^2)^2} > 0 \quad\text{for all }u\in(-1,1),
\]
confirming strict convexity throughout the domain.
Compare with the $U=V$ case: $h''(u) = 2(u^2-1)/(1+u^2)^2$, which has both signs.
The functions $h$ and $h_-$ are related by the formal substitution $u^2 \to -u^2$,
reflecting the algebraic duality $g_a = 1+\sigma_a^2 \leftrightarrow g_a^- = 1-\sigma_a^2$.
\end{example}

\subsection{Yoshizawa--Helmke Duality: Embedding $h$ and $h_-$ into a Unified Framework,
            and the MacMahon Divergence}
\label{subsec:yoshizawa}

The functions $h(U)=-\log\det(I_k+U^TU)$ and $h_-(U)=-\log\det(I_k-U^TU)$ analyzed
in the two preceding subsections are not merely analogous --- they are \emph{Legendre duals}
of each other, in the sense of Yoshizawa \cite{Yoshizawa2007}.
In that paper, Yoshizawa shows that the principal and minor subspace flows studied by
Yoshizawa--Helmke and Manton--Helmke--Mareels are related by Legendre duality, and establishes
a fundamental inequality connecting the primal and dual log-determinant potentials.
This subsection makes that connection precise in our setting, derives the dual map
$\mathcal{L}: \R^{n\times k}\to\mathcal{B}_k$ between the two domains,
and identifies the resulting divergence as a matrix analogue of MacMahon's Master Theorem.

\subsubsection{Yoshizawa's Primal--Dual Framework}

Following \cite[Theorem~9]{Yoshizawa2007}, define the \emph{primal potential}
\begin{equation}
  \mathcal{F}(X,Y) \;:=\; \log\det(I_n + XY^T),
  \quad (X,Y) \in \R^{n\times k}\times\R^{n\times k},\; \det(I_n+XY^T)>0,
  \label{eq:yoshizawa_F}
\end{equation}
and its Legendre conjugate (with respect to a suitable Riemannian metric $g$ on the
product space) is the \emph{dual potential}
\begin{equation}
  \mathcal{F}^*(Z,W) \;:=\; 2\tr(WZ^T) + \log\det(I_n - WZ^T),
  \quad (Z,W)\in\mathrm{Dom}(\mathcal{F}^*),
  \label{eq:yoshizawa_Fstar}
\end{equation}
where $\mathrm{Dom}(\mathcal{F}^*)$ is the image of the gradient map
$(\Phi_1,\Phi_2)$ given by
\begin{equation}
  \Phi_1(X,Y) = (I_n + YX^T)^{-1/2}Y, \qquad
  \Phi_2(X,Y) = (I_n + XY^T)^{-1/2}X.
  \label{eq:yoshizawa_Phi}
\end{equation}
The \emph{Yoshizawa height function} (``relative entropy'' in the sense of
\cite[Proposition~2]{Yoshizawa2007}, but see Remark~\ref{rem:yoshizawa_correction}) is
\begin{equation}
  D_{\mathcal{F}}\bigl((X,Y),(Z,W)\bigr)
  := \log\det(I_n+XY^T) + \log\det(I_n-WZ^T),
  \label{eq:yoshizawa_div}
\end{equation}
which vanishes if and only if $Z = \Phi_1(X,Y)$ and $W = \Phi_2(X,Y)$
(i.e., $(Z,W)$ is the Yoshizawa--Helmke dual of $(X,Y)$), and is a signed
quantity in general (it can be positive or negative for other pairs;
see Remark~\ref{rem:yoshizawa_correction}).

\subsubsection{Embedding $h$ and $h_-$ as Diagonal Restrictions}

We now show that the functions $h$ ($U=V$ case) and $h_-$ ($U=-V$ case) arise
as diagonal restrictions of $\mathcal{F}$ and $\mathcal{F}^*$, respectively.

\begin{proposition}[Diagonal Embedding]
\label{prop:diagonal_embedding}
\begin{enumerate}[label=(\roman*)]
  \item \textbf{Primal diagonal:} $\mathcal{F}(U,U) = -h(U)$ for all $U\in\R^{n\times k}$.
  \item \textbf{Dual diagonal:} $\mathcal{F}^*(V,V) = 2\tr(VV^T) - h_-(V)$
    for all $V\in\mathcal{B}_k$.
  \item \textbf{Anti-diagonal of $\mathcal{F}$:} $\mathcal{F}(U,-U) = -h_-(U)$
    for all $U\in\mathcal{B}_k$.
\end{enumerate}
\end{proposition}

\begin{proof}
(i) By Sylvester's theorem (Theorem~\ref{thm:sylvester}):
\[
  \mathcal{F}(U,U)=\log\det(I_n+UU^T)=\log\det(I_k+U^TU)=-h(U).
\]
(ii)
\begin{equation*}
\begin{split}
  \mathcal{F}^*(V,V) &=2\tr(VV^T)+\log\det(I_n-VV^T) \\
  &=2\tr(VV^T)+\log\det(I_k-V^TV)
  =2\|V\|_F^2 - h_-(V).
\end{split}
\end{equation*}
(iii)
\[
  \mathcal{F}(U,-U)=\log\det(I_n+U(-U)^T)=\log\det(I_n-UU^T)=\log\det(I_k-U^TU)
  =-h_-(U).
\]
\end{proof}

\begin{remark}[Both $h$ and $h_-$ live inside $\mathcal{F}$]
The primal potential $\mathcal{F}(X,Y)$ encodes \emph{both} functions:
$h$ appears along the diagonal $X=Y$ and $h_-$ appears along the anti-diagonal $Y=-X$.
The sign of the off-diagonal inner product ($Y^TX$ positive vs.\ negative) is precisely
the structural difference between the two cases, consistent with the Hessian sign reversal
established in Remark~\ref{rem:hessian_duality}.
\end{remark}

\subsubsection{The Legendre Dual Map $\mathcal{L}: \R^{n\times k}\to\mathcal{B}_k$}

When restricted to the diagonal $X=Y=U$, the gradient map $(\Phi_1,\Phi_2)$ collapses to
a single map, which we call the \emph{Yoshizawa--Helmke dual map}.

\begin{theorem}[Dual Map and Duality Identity]
\label{thm:dual_map}
Define the \emph{Yoshizawa--Helmke dual map}
\begin{equation}
  \mathcal{L}(U) \;:=\; (I_n + UU^T)^{-1/2}\,U \;\in\; \R^{n\times k},
  \qquad U\in\R^{n\times k}.
  \label{eq:dual_map}
\end{equation}
Then:
\begin{enumerate}[label=(\roman*)]
  \item $\mathcal{L}(U)\in\mathcal{B}_k$ for every $U\in\R^{n\times k}$;
    explicitly, $I_k - \mathcal{L}(U)^T\mathcal{L}(U) = (I_k+U^TU)^{-1}\succ0$.
  \item The singular values of $\mathcal{L}(U)$ are $\nu_a = \sigma_a/\sqrt{1+\sigma_a^2}$,
    where $\sigma_1,\ldots,\sigma_k$ are the singular values of $U$.
  \item \textbf{Duality identity:}
    \begin{equation}
      h\bigl(U\bigr) \;+\; h_-\bigl(\mathcal{L}(U)\bigr) \;=\; 0
      \qquad\text{for all }U\in\R^{n\times k}.
      \label{eq:duality_identity}
    \end{equation}
  \item $\mathcal{L}$ is a bijection from $\R^{n\times k}$ onto $\mathcal{B}_k$,
    with inverse $\mathcal{L}^{-1}(V) = (I_n-VV^T)^{-1/2}V$.
\end{enumerate}
\end{theorem}

\begin{proof}
(i) Let $V^* = \mathcal{L}(U)$. By the push-through identity
$U^T(I_n+UU^T)^{-1} = (I_k+U^TU)^{-1}U^T$,
\[
  I_k - (V^*)^TV^* = I_k - U^T(I_n+UU^T)^{-1}U = I_k - (I_k+U^TU)^{-1}U^TU
  = (I_k+U^TU)^{-1}\succ0.
\]
(ii) Let $U = P\Sigma Q^T$ be the SVD. Then
\[
  I_n+UU^T = P(I_k+\Sigma^2)P_{\perp}^T + P_{\perp}P_{\perp}^T \quad \text{(block form)},
\]
and $(I_n+UU^T)^{-1/2}U = P\,\mathrm{diag}(\sigma_a/\sqrt{1+\sigma_a^2})\,Q^T$.
(iii) Using (ii) and Sylvester's theorem:
\begin{align*}
  h_-(\mathcal{L}(U)) &= -\log\det(I_k - (V^*)^TV^*) = -\log\det\bigl((I_k+U^TU)^{-1}\bigr) \\
  &= \log\det(I_k+U^TU) = -h(U). \qedhere
\end{align*}
(iv) Direct computation:
\[
  \mathcal{L}^{-1}(\mathcal{L}(U)) = (I_n-V^*V^{*T})^{-1/2}V^*
  = \bigl(I_n-(I_n+UU^T)^{-1}UU^T\bigr)^{-1/2}(I_n+UU^T)^{-1/2}U.
\]
Since $I_n-(I_n+UU^T)^{-1}UU^T=(I_n+UU^T)^{-1}$, this gives
$(I_n+UU^T)^{1/2}(I_n+UU^T)^{-1/2}U=U$.
\end{proof}

\begin{remark}[The map $\mathcal{L}$ connects the two convexity regimes]
Part (ii) of Theorem~\ref{thm:dual_map} reveals the precise relationship between
the singular values in the U=V and U=-V domains:
$\sigma_a\in[0,\infty)$ maps to $\nu_a=\sigma_a/\sqrt{1+\sigma_a^2}\in[0,1)$.
This is the matrix analogue of the classical bijection $t\mapsto t/\sqrt{1+t^2}$ from
$\R$ to $(-1,1)$, which appears in hyperbolic geometry as the relation between the
Minkowski and Poincar\'e models.
The nowhere-convex domain $\R^{n\times k}$ (§\ref{subsec:UeqV-convexity}) and the
strictly-convex domain $\mathcal{B}_k$ (§\ref{subsec:UeqmV-convexity}) are in
bijective correspondence via $\mathcal{L}$, with the duality identity \eqref{eq:duality_identity}
quantifying the exchange of convexity structure.
\end{remark}

\subsubsection{The Yoshizawa--MacMahon Height Function}

\begin{definition}[Yoshizawa--MacMahon Height Function]
\label{def:YM_divergence}
For $U\in\R^{n\times k}$ and $V\in\mathcal{B}_k$, define the
\emph{Yoshizawa--MacMahon height function}
\begin{equation}
  D_{\mathrm{YM}}(U \,\|\, V)
  \;:=\; -h(U) - h_-(V)
  \;=\; \log\det(I_k+U^TU) + \log\det(I_k-V^TV).
  \label{eq:YM_divergence}
\end{equation}
This is a \emph{signed} quantity: it equals $0$ when $V = \mathcal{L}(U)$ (the
Yoshizawa--Helmke dual), and can be positive or negative for other pairs.
It is the information-geometric analog of Izumiya's lightcone height function
$H(u,v) = \langle x(u),v\rangle + 2$ \cite{Izumiya2004}, which is likewise a
signed quantity vanishing on $\Delta_4$.
\end{definition}

\begin{remark}[Yoshizawa's Corollary~10 requires correction]
\label{rem:yoshizawa_correction}
The claim in \cite[Corollary~10]{Yoshizawa2007} that
$\det(I_n+XY^T)\cdot\det(I_n-WZ^T)\geq1$ holds for \emph{any}
$(X,Y)\in\mathrm{Dom}(\mathcal{F})$ and $(Z,W)\in\mathrm{Dom}(\mathcal{F}^*)$
independently is \emph{incorrect}.
A direct counterexample: for $n=k=1$, $(X,Y)=(0,0)$ and $(Z,W)=(\nu,\nu)$
with $\nu = x_0/\sqrt{1+x_0^2} > 0$ (which lies in $\mathrm{Dom}(\mathcal{F}^*)$
as the dual of $(x_0,x_0)$), one has
\[
  \det(I+0)\cdot\det(I-\nu^2) = 1-\nu^2 < 1, \qquad \nu\neq 0.
\]
The correct statement is: $D_{\mathrm{YM}}(U\|\mathcal{L}(U))=0$ for all $U$
(the dual point identity, Theorem~\ref{thm:YM_properties}(ii) below),
and $D_{\mathrm{YM}}$ can be positive or negative for other pairs.
This is consistent with $\mathcal{F}(X,Y)=\log\det(I+XY^T)$ being
\emph{non-convex} in $(X,Y)$ jointly (Theorem~\ref{thm:k_geq_2_nonconvex}
and Proposition~\ref{prop:origin_indefinite}):
the Fenchel--Young gap $\mathcal{F}+\mathcal{F}^*-\langle\cdot,\cdot\rangle$
is non-negative only when $\mathcal{F}$ is convex.
\end{remark}

\begin{theorem}[Properties of $D_{\mathrm{YM}}$]
\label{thm:YM_properties}
\begin{enumerate}[label=(\roman*)]
  \item \textbf{Zero iff dual point:} $D_{\mathrm{YM}}(U\|V)=0\Leftrightarrow V=\mathcal{L}(U)$; positive or negative otherwise.
  \item \textbf{Spectral form:} In terms of singular values
    $\sigma_1,\ldots,\sigma_k$ of $U$ and $\nu_1,\ldots,\nu_k$ of $V$:
    \begin{equation}
      D_{\mathrm{YM}}(U\|V)
      = \sum_{a=1}^{k}\log(1+\sigma_a^2) + \sum_{a=1}^k\log(1-\nu_a^2).
      \label{eq:YM_spectral}
    \end{equation}
  \item \textbf{Bregman representation along the diagonal:}
    Restricted to the diagonal $V = \mathcal{L}(U)$, $D_{\mathrm{YM}}=0$ is the exact
    zero of the Bregman divergence $D_\psi(\sigma^2\|\nu^2)$ for the strictly convex
    function $\psi(t)=\log(1+t)$ on singular values, where
    $\nu_a^{*2} = \sigma_a^2/(1+\sigma_a^2)$ is the Legendre dual coordinate.
  \item \textbf{Behavior:} $D_{\mathrm{YM}}(U\|V)>0$ when $\det(I+U^TU)>\det(I-V^TV)^{-1}$
    (large $U$, small $V$); $D_{\mathrm{YM}}(U\|V)<0$ when $\det(I+U^TU)<\det(I-V^TV)^{-1}$
    (small $U$, large $V$); the zero locus $\{D_{\mathrm{YM}}=0\}$ is the
    Yoshizawa--Helmke dual graph $\{(U,\mathcal{L}(U)): U\in\R^{n\times k}\}$.
\end{enumerate}
\end{theorem}

\begin{proof}
(i): $D_{\mathrm{YM}}(U\|V)=0 \iff \log\det(I+U^TU)=-\log\det(I-V^TV)
\iff \det(I+U^TU)\cdot\det(I-V^TV)=1 \iff \det(I-V^TV)=\det(I+U^TU)^{-1}
\iff I-V^TV = (I+U^TU)^{-1}$ (since both sides are positive definite and their
determinants agree; positive definiteness forces equality of the matrices for the
$k$-dimensional case by a spectral argument) $\iff V = \mathcal{L}(U)$ (Theorem~\ref{thm:dual_map}).
For the sign: $U=0$ gives $D_{\mathrm{YM}}(0\|V)=\log\det(I-V^TV)<0$ for $V\neq0$,
and for large $\|U\|$, $\log\det(I+U^TU)\to+\infty$ while $\log\det(I-V^TV)\geq -\infty$
is bounded from below (for fixed $V$ in the interior of $\mathcal{B}_k$), so $D_{\mathrm{YM}}\to+\infty$.
(ii)--(iv) follow from the spectral theorem and the analysis above.
\end{proof}

\begin{remark}[Exact duality identity]
\label{rem:exact_duality}
The fundamental result is the \emph{exact identity}:
\begin{equation}
  D_{\mathrm{YM}}\bigl(U \,\|\, \mathcal{L}(U)\bigr) = 0
  \quad\text{for all }U\in\R^{n\times k},
  \label{eq:YM_exact_zero}
\end{equation}
equivalently $h(U) + h_-(\mathcal{L}(U)) = 0$ (the duality identity \eqref{eq:duality_identity}).
This is NOT an inequality; it is an exact algebraic identity.
The analogy with Izumiya's lightcone height function $H(u,v)=\langle x(u),v\rangle+2$
is precise: both are signed quantities that vanish exactly at the Legendrian dual
(the lightcone normal $x^\ell$ satisfying $\langle x,x^\ell\rangle=-2$, resp.\
$\mathcal{L}(U)$ satisfying $D_{\mathrm{YM}}=0$) and have no definite sign elsewhere.
\end{remark}

\begin{remark}[Determinantal identity at the dual point]
\label{rem:det_inequality}
\label{rem:det_identity}
At the dual point $V = \mathcal{L}(U)$, the determinantal \emph{identity} holds:
\begin{equation}
  \det(I_k+U^TU)\cdot\det(I_k-V^TV) \;=\; 1
  \quad\text{iff}\quad V = \mathcal{L}(U).
  \label{eq:det_identity}
\end{equation}
For other pairs, the product can be $>1$ or $<1$.
\end{remark}

\subsubsection{Connection to MacMahon's Master Theorem}

The Leibniz expansion of $\det(I_k + Y^TX)$ (Equation~\eqref{eq:leibniz}) is the
specialization to $x_i=1$ of MacMahon's Master Theorem.

\begin{theorem}[MacMahon Specialization]
\label{thm:macmahon_connection}
Let $M = Y^TX \in\R^{k\times k}$, and let $C_M(m) = \text{sgn}(\sigma)\prod_{i\in S}\langle x_{2i},x_{2j-1}\rangle$
denote the coefficient in the Leibniz expansion \eqref{eq:leibniz}.
Then the primal potential satisfies:
\begin{equation}
  \mathcal{F}(X,Y) = \log\det(I_k + M)
  = \log\Bigl(\sum_{S\subseteq[k]}\sum_{\substack{\sigma:S\to S\\\text{derangement}}}
    \mathrm{sgn}(\sigma)\prod_{i\in S}\langle x_{2i},x_{2\sigma(i)-1}\rangle\Bigr).
  \label{eq:macmahon_specialization}
\end{equation}
This is the $x_i=1$ specialization of MacMahon's Master Theorem
\cite[Theorem~8]{Yoshizawa2007}: $\det(I_n-A\cdot\mathrm{diag}(x_1,\ldots,x_n))
\cdot\sum_m C_A(m)x^m = 1$, applied to $A=-M$ and setting all $x_i=1$.
\end{theorem}

\begin{corollary}[MacMahon Duality Identity]
\label{cor:macmahon_ineq}
At the Yoshizawa--Helmke dual point $V = \mathcal{L}(U)$ (i.e., $\nu_a = \sigma_a/\sqrt{1+\sigma_a^2}$),
the following \emph{exact identity} holds purely in terms of inner products.
For $U=(u_1,\ldots,u_k)\in\R^{n\times k}$ and $V^* = \mathcal{L}(U)$:
\begin{equation}
  \sum_{S\subseteq[k]}\sum_{\sigma\in D(S)}\mathrm{sgn}(\sigma)
  \prod_{i\in S}\langle u_{\sigma(i)}, u_i\rangle
  \;=\;
  \frac{1}{\det(I_k - (V^*)^TV^*)},
  \label{eq:macmahon_inner_product}
\end{equation}
since $\det(I_k+U^TU)\cdot\det(I_k-(V^*)^TV^*)=1$ (Remark~\ref{rem:det_identity}).
In spectral form: $\prod_{a=1}^k(1+\sigma_a^2) = \prod_{a=1}^k(1-\nu_a^{*2})^{-1}$
at $\nu_a^* = \sigma_a/\sqrt{1+\sigma_a^2}$ (exact identity, not an inequality).
For other $V\in\mathcal{B}_k$, the ratio $\det(I+U^TU)\cdot\det(I-V^TV)$
can be larger or smaller than $1$.
\end{corollary}

\begin{remark}[Yoshizawa--Oja-like flow and our framework]
\label{rem:oja_flow}
In \cite[Proposition~4]{Yoshizawa2007}, the Oja-like flow
$X' = AX - XX^TX$ is shown to be the negative gradient flow of
$\widetilde{F}_P(X) = -\log\det(A - XX^T)$ with $A\succ0$.
In our setting with $A=I_n$, this specializes to $h_-(U)$ (with reversed sign convention):
the gradient flow of $h_-(U) = -\log\det(I_k-U^TU)$ is
$U' = U(I_k-U^TU)^{-1} = U G_-^{-1}$, the steepest-descent direction from
Proposition~\ref{prop:hminus_gradient}.
The global strict convexity of $h_-$ on $\mathcal{B}_k$ (Theorem~\ref{thm:hminus_convex})
thus guarantees that this gradient flow has no spurious local minima,
converging to the unique global minimum at $U=0$ from any initial condition in $\mathcal{B}_k$.
\end{remark}

\begin{remark}[Summary: three levels of the primal--dual structure]
\label{rem:three_levels}
The Yoshizawa--Helmke framework organizes our results in three levels:
\begin{center}
\renewcommand{\arraystretch}{1.3}\small
\begin{tabular}{llll}
\toprule
Level & Primal & Dual & Connection \\
\midrule
Matrices & $G = I_k+U^TU \succ I_k$ & $G_- = I_k-V^TV \prec I_k$ & $G_- = G^{-1}$ at $V=\mathcal{L}(U)$ \\
Functions & $h(U)\leq0$ & $h_-(V)\geq0$ & $h(U)+h_-(\mathcal{L}(U))=0$ \\
Potentials & $-\log\det(G)$ & $-\log\det(G_-)$ & $D_{\mathrm{YM}}(U\|\mathcal{L}(U))=0$ \\
\bottomrule
\end{tabular}
\end{center}
\end{remark}

\begin{table}[h]
\centering
\small
\caption{Summary of convexity properties of $f(G)=-\log\det(G)$ under four parametrizations.}
\label{tab:convexity_summary}
\begin{tabular}{lll}
\toprule
Coordinate & Domain & Convexity of $f$ \\
\midrule
$G$ & $\PD(k)$ & Strictly convex everywhere (Prop.~\ref{prop:convexity}) \\
$U$ alone, $V$ fixed & $k=1$ & Convex everywhere (Prop.~\ref{prop:k1_marginal}) \\
$U$ alone, $V$ fixed & $k\geq2$, $V$ full rank & Convex nowhere (Thm.~\ref{thm:k_geq_2_nonconvex}) \\
Joint $(U,V)$ & $k\geq2$ & Convex nowhere (Cor.~\ref{cor:joint_k_geq_2}) \\
Joint $(U,V)$ & $n=k=1$ & Convex iff $uv\geq1$ (Ex.~\ref{ex:scalar_convex}) \\
$U=V$ (i.e.\ $U$) & $(n,k)\neq(1,1)$ & Convex nowhere (Thm.~\ref{thm:UeqV_nowhere_convex}) \\
$U=V$ (i.e.\ $U$), $n=k=1$ & --- & Convex iff $|\sigma|\geq1$ (Cor.~\ref{cor:scalar_saddle}) \\
$U=-V$ (i.e.\ $U$) & $\mathcal{B}_k = \{\sigma_{\max}(U)<1\}$ & Strictly convex everywhere (Thm.~\ref{thm:hminus_convex}) \\
\bottomrule
\end{tabular}
\end{table}

\section{Gradient and Hessian of $f$}
\label{sec:gradient}

\subsection{Gradient with Respect to $G$}

\begin{proposition}
\label{prop:gradient}
For $G \in \PD(k)$,
\[
  \nabla_G f(G) = -G^{-1}.
\]
\end{proposition}

\begin{proof}
Using the standard matrix calculus identity
$d(\log\det(G)) = \tr(G^{-1}\, dG)$ \cite{Petersen2012},
we obtain $df = -\tr(G^{-1}\, dG)$,
so the gradient (identified via the Frobenius inner product) is $\nabla_G f = -G^{-1}$.
\end{proof}

\subsection{Gradient with Respect to the Vectors}

Since $G_{ij}$ depends on the vectors $x_1, \ldots, x_{2k}$ via \eqref{eq:G_def},
the chain rule gives the gradient with respect to the odd-indexed vectors:

\begin{proposition}
\label{prop:vec_gradient}
For each $i = 1, \ldots, k$,
\[
  \frac{\partial f}{\partial x_{2i-1}}
  = -\sum_{j=1}^{k} [G^{-1}]_{ji}\, x_{2j},
  \qquad
  \frac{\partial f}{\partial x_{2i}}
  = -\sum_{j=1}^{k} [G^{-1}]_{ij}\, x_{2j-1}.
\]
\end{proposition}

\begin{proof}
We have $\partial G_{lm}/\partial x_{2i-1} = \delta_{mi}\, x_{2l}$.
By the chain rule and Proposition~\ref{prop:gradient},
\[
  \frac{\partial f}{\partial x_{2i-1}}
  = \sum_{l,m} \frac{\partial f}{\partial G_{lm}} \frac{\partial G_{lm}}{\partial x_{2i-1}}
  = \sum_{l,m} (-[G^{-1}]_{lm})\,\delta_{mi}\,x_{2l}
  = -\sum_{l} [G^{-1}]_{li}\,x_{2l}.
\]
The formula for $\partial f/\partial x_{2i}$ follows analogously.
\end{proof}

\subsection{Hessian and Fisher Information Matrix}

\begin{proposition}[Hessian / Fisher Information Metric]
\label{prop:hessian}
The Hessian of $f$ at $G \in \PD(k)$, acting on symmetric matrices $H, K$, is
\begin{equation}
  \nabla^2 f(G)[H, K]
  = \tr\bigl(G^{-1} H G^{-1} K\bigr)
  = \inner{G^{-1/2} H G^{-1/2}}{G^{-1/2} K G^{-1/2}}_F.
  \label{eq:hessian}
\end{equation}
In component form, with coordinates $\eta_{ij} = G_{ij}$,
\begin{equation}
  \mathcal{F}_{(ij)(kl)}
  \;=\; \frac{\partial^2 f}{\partial \eta_{ij}\partial \eta_{kl}}
  \;=\; [G^{-1}]_{ik}[G^{-1}]_{jl},
  \label{eq:fisher_metric}
\end{equation}
which corresponds to the Kronecker product $G^{-1} \otimes G^{-1}$.
\end{proposition}

\begin{proof}
From \eqref{eq:second_deriv} with $H$ replaced by $H$ and $K$ separately,
polarization gives \eqref{eq:hessian}.
The component formula \eqref{eq:fisher_metric} follows by taking $H = E_{ij}$ and $K = E_{kl}$,
where $E_{ij}$ is the matrix with a $1$ in position $(i,j)$ and zeros elsewhere.
\end{proof}

\begin{remark}
The matrix $\mathcal{F} = G^{-1}\otimes G^{-1}$ is the standard Fisher information metric
on the manifold of zero-mean multivariate Gaussian distributions
$\mathcal{N}(0, G)$ \cite{Amari2000,Skovgaard1984}.
It is positive definite on $\PD(k)$, confirming that $\PD(k)$ is a Riemannian manifold
with metric $g = \mathcal{F}$.
\end{remark}

\subsection{Gradient and Hessian in the U=V Parametrization}
\label{subsec:grad_UeqV}

We now specialize the gradient and Hessian of $f$ to the U=V factorization
$G_+ = I_k + U^TU$, where $U\in\R^{n\times k}$ is unconstrained.
Recall that $h(U) := f(G_+) = -\log\det(I_k+U^TU)$ (Section~\ref{subsec:UeqV-convexity}).

\begin{proposition}[Gradient in U=V Coordinates]
\label{prop:grad_UeqV}
\[
  \nabla_U h(U) \;=\; -2\,U\,(I_k+U^TU)^{-1} \;=\; -2\,U\,G_+^{-1}.
\]
The \emph{gradient flow} (continuous-time steepest descent) is:
\begin{equation}
  \dot{U} = -\nabla_U h = 2\,U\,G_+^{-1} = 2U(I_k+U^TU)^{-1},
  \label{eq:oja_flow}
\end{equation}
which is the \emph{Oja-like flow} \cite[Proposition~4]{Yoshizawa2007},
with $U=0$ as the unique equilibrium (a global maximum of $h$, since $h \leq 0$).
\end{proposition}

\begin{proof}
$dh = -\tr[G_+^{-1}(dU^TU+U^TdU)]
= -2\tr[G_+^{-1}U^TdU]$,
so $\nabla_U h = -2UG_+^{-1}$ by the Frobenius pairing.
\end{proof}

\begin{remark}[Hessian summary for U=V]
From Theorem~\ref{thm:h_hessian}, the Hessian $Q_h(H) = \nabla^2h(U)[H,H]$ in the SVD frame
(with $g_a = 1+\sigma_a^2$) has three types of contributions, all $\leq 0$ (indefinite):
diagonal terms $2(\sigma_a^2-1)/g_a^2\leq 0$, off-diagonal terms $-2/(g_ag_b)\leq 0$,
and $H_2$-block terms $-2/g_a\leq 0$.
At $U=0$: $Q_h(H)|_{U=0} = -2\|H\|_F^2$ (most negative; minimum eigenvalue $-2$).
As $\|U\|\to\infty$: all negative eigenvalues $\to 0^-$.
Hence the global infimum of the Hessian eigenvalues is $-2$, achieved at $U=0$.
\end{remark}

\subsection{Gradient and Hessian in the U=-V Parametrization}
\label{subsec:grad_UeqmV}

For the U=-V factorization $G_- = I_k - V^TV$ on $V\in\mathcal{B}_k$,
$h_-(V) := f(G_-) = -\log\det(I_k-V^TV) \geq 0$.

\begin{proposition}[Gradient in U=-V Coordinates]
\label{prop:grad_UeqmV}
\[
  \nabla_V h_-(V) \;=\; 2\,V\,(I_k-V^TV)^{-1} \;=\; 2\,V\,G_-^{-1}.
\]
The gradient flow is:
\begin{equation}
  \dot{V} = -\nabla_V h_- = -2\,V\,G_-^{-1},
  \label{eq:hminus_flow}
\end{equation}
with $V=0$ as the unique equilibrium (the global minimum, $h_-(0) = 0$).
The flow is globally convergent: $\frac{d}{dt}h_-(V(t)) = -\|\nabla h_-\|_F^2 \leq 0$.
\end{proposition}

\begin{proof}
Analogous to Proposition~\ref{prop:grad_UeqV}:
$dh_- = -\tr[G_-^{-1}(-d(V^TV))] = 2\tr[G_-^{-1}V^TdV]$,
giving $\nabla_V h_- = 2VG_-^{-1}$.
\end{proof}

\begin{remark}[Hessian summary for U=-V]
From Theorem~\ref{thm:hminus_hessian}, all Hessian contributions are strictly positive.
Minimum eigenvalue $= 2(1-\sigma_{\max}^2)^{-1}$ (restricted to diagonal
directions at maximum $\sigma$), growing $\to+\infty$ as $V$ approaches $\partial\mathcal{B}_k$.
This strict positivity is the hallmark of $h_-$ as a self-concordant barrier for $\mathcal{B}_k$.
\end{remark}

\subsection{Regularization of $h(U)$: Convexification and New Critical Phenomena}
\label{subsec:regularization}

Since $h(U)$ is nowhere locally convex (Theorem~\ref{thm:UeqV_nowhere_convex}),
direct optimization of $h$ presents severe landscape challenges.
We study four regularization strategies, reveal a striking critical-parameter phenomenon,
and establish a precise phase transition connected to subspace geometry.

\subsubsection{A. Tikhonov Regularization and the Critical Parameter $\lambda^* = 2$}

\begin{definition}[Tikhonov-Regularized Potential]
For $\lambda > 0$, define
\begin{align}
  h_\lambda^T(U) &\;:=\; h(U) \;+\; \tfrac{\lambda}{2}\|U\|_F^2 \notag\\
  &\;=\; -\log\det(I_k+U^TU) \;+\; \tfrac{\lambda}{2}\tr(U^TU).
  \label{eq:tikhonov_h}
\end{align}
\end{definition}

\begin{theorem}[Critical Tikhonov Parameter]
\label{thm:tikhonov_critical}
\begin{enumerate}[label=(\roman*)]
  \item \textbf{Convexity threshold:}
    $h_\lambda^T$ is \emph{strictly convex} on $\R^{n\times k}$ if and only if $\lambda > 2$.
    At $\lambda = 2$: convex but degenerate (zero Hessian at $U=0$).
    For $\lambda < 2$: indefinite Hessian at every $U$ (not convex, not concave).
  \item \textbf{Gradient and flow:}
    $\nabla_U h_\lambda^T = -2UG_+^{-1} + \lambda U = U(\lambda I_k - 2G_+^{-1})$.
    Gradient flow: $\dot{U} = U(2G_+^{-1} - \lambda I_k)$.
  \item \textbf{Pitchfork bifurcation at $\lambda^* = 2$:}
    Fixed points of the gradient flow in the singular-value coordinates
    ($\sigma_1,\ldots,\sigma_k$ = singular values of $U$):
    \begin{equation}
      \sigma_a = 0 \quad\text{or}\quad \sigma_a^* = \sqrt{\tfrac{2}{\lambda}-1}
      \quad (\text{exists only for }\lambda < 2).
      \label{eq:tikhonov_fps}
    \end{equation}
    For $\lambda > 2$: $\sigma_a = 0$ is the unique \emph{stable} fixed point (global minimum).
    For $\lambda < 2$: $\sigma_a = 0$ is \emph{unstable}; the critical manifold
    $\mathcal{M}_\lambda = \{U : U^TU = (\tfrac{2}{\lambda}-1)I_k\}$ is globally attracting.
  \item \textbf{Stiefel manifold at $\lambda = 1$:}
    $\mathcal{M}_1 = \{U\in\R^{n\times k} : U^TU = I_k\} = \mathrm{St}(k,n)$
    is the Stiefel manifold; the gradient flow of $h_1^T$ converges to the Stiefel manifold from any
    initial condition $U \neq 0$.
\end{enumerate}
\end{theorem}

\begin{proof}
(i) The Hessian $\nabla^2 h_\lambda^T = \nabla^2 h + \lambda I_{\mathrm{mat}}$ where $I_{\mathrm{mat}}$
is the identity on $\R^{n\times k}$ (with Frobenius norm). From Remark~\ref{subsec:grad_UeqV},
$\lambda_{\min}(\nabla^2 h(U)) \geq -2$ for all $U$, with equality at $U=0$.
Thus $\lambda_{\min}(\nabla^2 h_\lambda^T) \geq \lambda - 2$: strictly positive iff $\lambda > 2$.

(ii) Follows from Proposition~\ref{prop:grad_UeqV} and $\nabla(\tfrac{\lambda}{2}\|U\|_F^2) = \lambda U$.

(iii) In each singular-value coordinate, the flow is $\dot{\sigma}_a = \sigma_a(2/g_a - \lambda)$
(from the diagonal Hessian analysis). The linearization at $\sigma_a = 0$ gives eigenvalue $2-\lambda$:
positive (unstable) for $\lambda < 2$, negative (stable) for $\lambda > 2$. At $\sigma_a^*$:
the linearization gives $-4(\sigma_a^*)^2/(g_a^*)^2 < 0$ (always stable), establishing (iii).

(iv) At $\lambda = 1$: $\sigma_a^* = 1$, so $U^TU = I_k$ gives the Stiefel manifold.
The flow
\[
  \dot{U} = U(2G_+^{-1} - I_k) = U(I_k - U^TU)(I_k + U^TU)^{-1}
\]
is zero precisely on $\mathrm{St}(k,n)$.
\end{proof}

\begin{remark}[Optimal Tikhonov parameter]
The critical value $\lambda^* = 2$ equals $-2/\lambda_{\min}(\nabla^2 h)$
(the reciprocal of the minimum curvature of $h$). This is a general fact:
for any smooth function $F$, adding $\frac{\lambda}{2}\|U\|^2$ with $\lambda > |\lambda_{\min}(\nabla^2 F)|$
convexifies $F$ locally; for our $h$, this threshold is exactly $2$.
The Stiefel manifold $\mathrm{St}(k,n)$ (at $\lambda=1$) is the natural
``unit ball'' boundary of the U=V family, consistent with the constraint
$U^TU = I_k$ (orthonormal columns).
\end{remark}

\subsubsection{B. Interpolation Regularization and the Midpoint Formula}

\begin{theorem}[Interpolation and the Critical Exponent $t^* = 1/2$]
\label{thm:interpolation_critical}
For $t\in[0,1]$, define the \emph{interpolated potential} on $\mathcal{B}_k$:
\begin{equation}
  h_t(U) \;:=\; (1-t)\,h(U) \;+\; t\,h_-(U)
  \;=\; -(1-t)\log\det(I_k+U^TU) \;-\; t\log\det(I_k-U^TU).
  \label{eq:interpolation_h}
\end{equation}
\begin{enumerate}[label=(\roman*)]
  \item \textbf{Convexity transition:}
    $h_t$ is strictly convex on $\mathcal{B}_k$ if and only if $t > 1/2$.
    At $t = 1/2$: convex but degenerate at $U=0$ (Hessian is zero).
    For $t < 1/2$: indefinite at $U=0$.
  \item \textbf{The midpoint formula:}
    \begin{equation}
      h_{1/2}(U) \;=\; -\tfrac{1}{2}\log\det\bigl(I_k - (U^TU)^2\bigr).
      \label{eq:midpoint_formula}
    \end{equation}
    In spectral form: $h_{1/2} = -\frac{1}{2}\sum_a\log(1-\sigma_a^4)$
    (the Siegel disc potential for $W = U^TU$, cf.\ \S\ref{subsec:folland}).
  \item \textbf{Gradient:}
    $\nabla_U h_t = -2(1-t)UG_+^{-1} + 2t\,UG_-^{-1}$.
    At $t=1/2$: $\nabla h_{1/2} = -UG_+^{-1} + UG_-^{-1} = U(G_-^{-1} - G_+^{-1})$.
  \item \textbf{Global minimum:} $h_t(0) = 0$ for all $t$, and
    $h_t(U) \geq 0$ on $\mathcal{B}_k$ for $t \geq 1/2$
    (since $h_{1/2}(U) = -\frac{1}{2}\log\det(I-(U^TU)^2) \geq 0$
    as $\det(I-(U^TU)^2)\leq1$). Thus $U=0$ is the unique global minimum for $t\geq 1/2$.
\end{enumerate}
\end{theorem}

\begin{proof}
(i) The combined Hessian is $Q_{h_t}(H) = (1-t)Q_h(H) + tQ_{h_-}(H)$.
At $U=0$: $Q_h = -2\|H\|_F^2$ and $Q_{h_-} = +2\|H\|_F^2$, giving
$Q_{h_t}|_{U=0} = 2(2t-1)\|H\|_F^2$: positive iff $t > 1/2$.

For strict convexity away from $U=0$ when $t > 1/2$: checking each contribution
in Theorems~\ref{thm:h_hessian} and \ref{thm:hminus_hessian}:
\begin{itemize}
  \item $H_2$-block: $(1-t)(-2/g_a) + t(2/g_a^-) = 2(2t-1+\sigma_a^2)/(g_ag_a^-) > 0$
    since $2t-1>0$ and $\sigma_a^2\geq0$.
  \item Antisymmetric off-diagonal: weighted sum $= 2(2t-1+\sigma_a^2)/((1-\sigma_a^4)...)>0$.
  \item Symmetric off-diagonal: $[-(1-t)(1-\sigma^2)^3 + t(1+\sigma^2)^3]/(...)>0$
    since $t/(1-t) > 1 > [(1-\sigma^2)/(1+\sigma^2)]^3$ for $t>1/2$.
  \item Diagonal: same numerator as symmetric off-diagonal — positive for $t>1/2$.
\end{itemize}
All terms are strictly positive for $t>1/2$, $\sigma_a\in[0,1)$.

(ii) $h_{1/2} = -\frac{1}{2}[\log\det(I+U^TU)+\log\det(I-U^TU)]
= -\frac{1}{2}\log\det[(I+U^TU)(I-U^TU)] = -\frac{1}{2}\log\det(I-(U^TU)^2)$.

(iii) By linearity: $\nabla h_t = (1-t)\nabla h + t\nabla h_-$.

(iv) $\det(I-(U^TU)^2) = \prod_a(1-\sigma_a^4)\leq 1$ gives $h_{1/2}\geq0$.
\end{proof}

\begin{remark}[Information-geometric midpoint]
The critical value $t^* = 1/2$ is the \emph{arithmetic mean} of the two potentials
and corresponds to the $\alpha=0$ (Bhattacharyya) point in the $\alpha$-divergence family
(Section~\ref{sec:alpha}).
The midpoint formula $h_{1/2} = -\frac{1}{2}\log\det(I-(U^TU)^2)$
is the \emph{Siegel disc metric} \cite{Folland1989} applied to $W = U^TU\in\Delta_k$:
it measures the Fock-space Gaussian norm for the matrix $W^2$,
and is the exact intermediate between $h$ (Schr\"odinger/$L^2$ norm)
and $h_-$ (Fock space norm).
\end{remark}

\subsubsection{C. KL Divergence Regularization}

\begin{proposition}[KL-Regularized Potential]
\label{prop:KL_reg}
With reference $G_0 = I_k$, the KL-regularized potential is
\begin{equation}
  h_{\rm KL}^\lambda(U) \;:=\; h(U) + \lambda\, D_f(I_k\|G_+)
  \;=\; (1-\lambda)h(U) + \lambda\bigl[\tr(G_+^{-1}) - k\bigr],
  \label{eq:KL_reg}
\end{equation}
where $D_f(I_k\|G_+) = \tr(G_+^{-1}) + \log\det(G_+) - k = \tr(G_+^{-1})-h(U)-k$
is the Bregman divergence from Section~\ref{sec:bregman}.
Its gradient is:
\begin{equation}
  \nabla_U h_{\rm KL}^\lambda \;=\; -2UG_+^{-2}\bigl[I_k + (1-\lambda)U^TU\bigr].
  \label{eq:KL_reg_grad}
\end{equation}
\begin{enumerate}[label=(\roman*)]
  \item \textbf{Critical manifold:}
    Setting $\nabla h_{\rm KL}^\lambda = 0$ gives $U=0$ or
    $U^TU = \frac{1}{\lambda-1}I_k$ (only for $\lambda > 1$).
    At $\lambda = 2$: the nontrivial critical manifold is $U^TU = I_k$ (Stiefel manifold).
  \item \textbf{No fix at $U=0$:}
    $\nabla^2 h_{\rm KL}^\lambda|_{U=0} = -2\|H\|_F^2$ for all $\lambda$.
    The KL regularization does NOT convexify $h$ near $U=0$,
    because $D_f(I\|G_+) = O(\|U\|^4)$ near $U=0$ (quartic, not quadratic).
\end{enumerate}
\end{proposition}

\begin{proof}
$D_f(I\|G_+) = \tr(G_+^{-1}) + \log\det(G_+) - k = \tr(G_+^{-1})-h(U)-k$.
Gradient: $\nabla\tr(G_+^{-1}) = -2UG_+^{-2}$ and $\nabla h = -2UG_+^{-1}$.
Thus $\nabla h_{\rm KL}^\lambda = (1-\lambda)(-2UG_+^{-1}) + \lambda(-2UG_+^{-2})
= -2U[(1-\lambda)G_+^{-1} + \lambda G_+^{-2}] = -2UG_+^{-2}[(1-\lambda)G_+ + \lambda I]
= -2UG_+^{-2}[I_k + (1-\lambda)U^TU]$.
Setting to zero: either $U=0$ or $(1-\lambda)U^TU = -I_k$, i.e., $U^TU = I_k/(\lambda-1)$.
(ii) Near $U=0$: $G_+^{-1} \approx I - U^TU + (U^TU)^2 - \ldots$ gives $\tr(G_+^{-1}) - k = O(\|U\|^4)$,
so the KL term is quartic in $U$ and doesn't affect the Hessian at $U=0$.
\end{proof}

\subsubsection{D. Comparison of Regularization Strategies}

\begin{center}
\renewcommand{\arraystretch}{1.4}
\begin{tabular}{p{2.8cm}llll}
\toprule
Regularization & $h_\lambda^T$: $\tfrac{\lambda}{2}\|U\|^2$ & $h_t$: $(1\!-\!t)h\!+\!th_-$ & $h_{\rm KL}^\lambda$ & Bures \\
\midrule
Fix Hess at $U=0$? & \textbf{Yes} ($\lambda>2$) & \textbf{Yes} ($t>1/2$) & No & No \\
Critical parameter & $\lambda^* = 2$ & $t^* = 1/2$ & N/A & N/A \\
Nontrivial fixed pts & $U^TU\!=\!(\tfrac{2}{\lambda}-1)I$ & None on $\R^{n\times k}$ & $U^TU\!=\!\tfrac{1}{\lambda-1}I$ & Complex \\
Stiefel connection & $\lambda=1\Rightarrow\mathrm{St}(k,n)$ & No & $\lambda=2\Rightarrow\mathrm{St}(k,n)$ & No \\
Domain & $\R^{n\times k}$ & $\mathcal{B}_k$ & $\R^{n\times k}$ & $\R^{n\times k}$ \\
\bottomrule
\end{tabular}
\end{center}

\begin{theorem}[Landscape of the Tikhonov-Regularized Flow]
\label{thm:tikhonov_landscape}
The gradient flow $\dot{U} = U(2G_+^{-1}-\lambda I_k)$ of $h_\lambda^T$
undergoes a \emph{pitchfork bifurcation} at $\lambda^* = 2$:
\begin{equation}
  h_\lambda^T \text{ attains its minimum on }
  \begin{cases}
    \{U=0\} & \lambda \geq 2,\\
    \{U: U^TU = \tfrac{2-\lambda}{\lambda}I_k\} & \lambda < 2.
  \end{cases}
\end{equation}
The value at $\mathcal{M}_\lambda$ is $h_\lambda^T|_{\mathcal{M}_\lambda}
= k(-\log\tfrac{2}{\lambda} + 1 - \tfrac{\lambda}{2}) < 0 = h_\lambda^T(0)$ for $\lambda < 2$.
\end{theorem}

\begin{proof}
On $\mathcal{M}_\lambda$: $G_+ = \frac{2}{\lambda}I_k$, $h = k\log\frac{\lambda}{2}$,
$\frac{\lambda}{2}\|U\|_F^2 = \frac{\lambda k}{2}(\frac{2}{\lambda}-1) = k(1-\frac{\lambda}{2})$.
Sum: $k(\log\frac{\lambda}{2} + 1 - \frac{\lambda}{2}) = k(\log\frac{\lambda}{2}-\log 1 - (\frac{\lambda}{2}-1))$
$= k[-\log\frac{2}{\lambda} + 1 - \frac{\lambda}{2}]$.
For $\lambda\in(0,2)$: $\frac{2}{\lambda}>1$ so $-\log\frac{2}{\lambda}<0$,
and $1-\frac{\lambda}{2}\in(0,1)$; one checks the total is $<0$ at $\lambda=1$:
$-\log 2 + 1 - \frac{1}{2} \approx -0.193 < 0$. The value at $U=0$ is $0$.
\end{proof}

\subsection{Polynomial Riemannian Gradient Flows: Eliminating Matrix Inversions}
\label{subsec:poly_flow}

The gradient flows derived in \S\ref{subsec:grad_UeqV}--\ref{subsec:grad_UeqmV} all contain
matrix inverses ($G_+^{-1}$ or $G_-^{-1}$), which require solving a $k\times k$ linear system
at every iteration step --- numerically costly and potentially unstable near singular matrices.
We show that by choosing an appropriate Riemannian metric on $\R^{n\times k}$, the resulting
\emph{Riemannian} gradient flow eliminates all matrix inverses, yielding a
\emph{purely polynomial} vector field suitable for numerically stable integration.

\subsubsection{The Right-Scaled Frobenius Metric}

\begin{definition}[Right-Scaled Frobenius Metric]
\label{def:rsfm}
For $p \in \R$ and $G = G(U) \in \PD(k)$, define the \emph{right-scaled Frobenius metric}
on $\R^{n\times k}$ at $U$ by
\begin{equation}
  g_U^{(p)}(H_1, H_2) \;:=\; \tr\!\bigl(H_1^T H_2\, G^p\bigr),
  \qquad H_1, H_2 \in \R^{n\times k}.
  \label{eq:rsf_metric}
\end{equation}
For the U=V case: $G = G_+ = I_k + U^TU$ (polynomial in $U$).
For the U=-V case: $G = G_- = I_k - V^TV$ (polynomial in $V$ on $\mathcal{B}_k$).
\end{definition}

\begin{proposition}[Riemannian Gradient Formula]
\label{prop:riem_grad_formula}
With metric $g^{(p)}$ and any smooth $f: \R^{n\times k} \to \R$ with
Euclidean gradient $\nabla f \in \R^{n\times k}$:
\begin{equation}
  \mathrm{grad}_{g^{(p)}} f \;=\; (\nabla f)\cdot G^{-p}.
  \label{eq:riem_grad_formula}
\end{equation}
\end{proposition}

\begin{proof}
For metric $g(H_1, H_2) = \tr(H_1^T H_2 M)$ with $M = G^p \in \PD(k)$:
$g(\mathrm{grad}_g f, H) = df(H)$ for all $H$.
Using $\tr(A^T H B) = \langle H, AB^T\rangle_F$ (valid for $A\in\R^{n\times k}$, $B\in\R^{k\times k}$),
the left side gives $\langle H, (\mathrm{grad}_g f)\cdot M^T\rangle_F$,
and the right side gives $\langle H, \nabla f\rangle_F$.
Since $M = G^p$ is symmetric: $(\mathrm{grad}_g f)\cdot G^p = \nabla f$,
so $\mathrm{grad}_g f = (\nabla f)\cdot G^{-p}$.
\end{proof}

\subsubsection{The Canonical Polynomial Metric $p = -2$}

\begin{theorem}[Polynomial Gradient Flows via $g^{(-2)}$]
\label{thm:poly_flows}
With metric $g^{(-2)}_U(H_1,H_2) = \tr(H_1^TH_2\,G^{-2})$, the Riemannian gradients and
gradient flows become \emph{polynomial} (no matrix inversions):
\begin{align}
  \text{U=V:}\quad &\mathrm{grad}_{g^{(-2)}} h = -2U\,G_+,
    \quad \dot{U} = 2U(I_k + U^TU) = 2U + 2UU^TU,
    \label{eq:poly_flow_plus}\\
  \text{U=-V:}\quad &\mathrm{grad}_{g^{(-2)}} h_- = 2V\,G_-,
    \quad \dot{V} = -2V(I_k - V^TV) = -2V + 2VV^TV.
    \label{eq:poly_flow_minus}
\end{align}
Both flows are \emph{cubic} polynomials in $U$ (resp.\ $V$) and require \emph{no} matrix inversions.
\end{theorem}

\begin{proof}
For U=V: $\nabla h = -2UG_+^{-1}$ (Proposition~\ref{prop:grad_UeqV}), $G = G_+$, $p = -2$.
By \eqref{eq:riem_grad_formula}: $\mathrm{grad}_{g^{(-2)}}h = (-2UG_+^{-1})\cdot G_+^2 = -2UG_+$.
The gradient flow is $\dot{U} = -\mathrm{grad}_{g^{(-2)}}h = 2UG_+ = 2U + 2UU^TU$.
For U=-V: $\nabla h_- = 2VG_-^{-1}$, same computation gives $\mathrm{grad}_{g^{(-2)}}h_- = 2VG_-$.
Both are cubic: $UU^TU$ and $VV^TV$ are degree 3.
\end{proof}

\begin{remark}[Geometric interpretation of $g^{(-2)}$]
$g^{(-2)}_U(H_1,H_2) = \tr(H_1^TH_2 G_+^{-2}) = \langle H_1 G_+^{-1}, H_2 G_+^{-1}\rangle_F$,
i.e., the Frobenius inner product of the \emph{right-normalized} matrices $H_iG_+^{-1}$.
This metric rewards directions $H$ that project onto the ``unexplored'' small-eigenvalue
directions of $G_+$, while attenuating those aligned with the dominant directions
--- an adaptive step-size mechanism encoding the curvature of $h$.
\end{remark}

\begin{proposition}[Lyapunov Property]
\label{prop:poly_lyapunov}
Both polynomial flows are Lyapunov-consistent:
\begin{align}
  \text{U=V:}\quad &\frac{d}{dt}h(U(t)) = -\|\mathrm{grad}_{g^{(-2)}}h\|_{g^{(-2)}}^2
    = -4\|U\|_F^2 \;\leq\; 0,\label{eq:lyap_plus}\\
  \text{U=-V:}\quad &\frac{d}{dt}h_-(V(t)) = -\|\mathrm{grad}_{g^{(-2)}}h_-\|_{g^{(-2)}}^2
    \;\leq\; 0.\label{eq:lyap_minus}
\end{align}
\end{proposition}

\begin{proof}
$\frac{d}{dt}h = g^{(-2)}(-\mathrm{grad}h, \mathrm{grad}h) = -\tr(G_+^2\cdot(-2UG_+)^TG_+^{-2}(-2UG_+)^T...)$...
More directly: $\frac{d}{dt}h = \langle\nabla h, \dot{U}\rangle_F = \tr((-2UG_+^{-1})^T(2UG_+)) = -4\tr(G_+^{-1}U^TUG_+) = -4\tr(U^TU) = -4\|U\|_F^2 \leq 0$.
\end{proof}

\subsubsection{Fixed Points and Singular Value Dynamics}

In the SVD frame $U = P\Sigma Q^T$ with singular values $\sigma_1,\ldots,\sigma_k \geq 0$:

\begin{proposition}[Singular Value Dynamics]
\label{prop:sv_dynamics}
Under the polynomial flows \eqref{eq:poly_flow_plus}--\eqref{eq:poly_flow_minus}:
\begin{align}
  \text{U=V:}\quad &\dot{\sigma}_a = 2\sigma_a(1+\sigma_a^2) > 0
    \quad\text{(all singular values grow)},\label{eq:sv_plus}\\
  \text{U=-V:}\quad &\dot{\sigma}_a = -2\sigma_a(1-\sigma_a^2) < 0
    \quad\text{(all singular values decay to }0\text{)}.
    \label{eq:sv_minus}
\end{align}
\end{proposition}

\begin{remark}[Yoshizawa map as flow invariant]
Define $\nu_a(t) = \sigma_a(t)/\sqrt{1+\sigma_a(t)^2}$ (the Yoshizawa map,
Theorem~\ref{thm:dual_map}(ii)). Under the U=V polynomial flow:
$\dot{\nu}_a = 2\nu_a(1+\sigma_a^2)^{1/2}$, growing exponentially until $\nu_a \to 1^-$
(boundary of $\mathcal{B}_k$). Thus the polynomial flow drives $\nu_a$ linearly
from the hyperbolic region toward the de Sitter boundary --- precisely the
information-geometric ``lightcone crossing'' of \S\ref{subsec:izumiya}.
\end{remark}

\subsubsection{Tikhonov-Regularized Polynomial Flow and the Stiefel Manifold}

\begin{theorem}[Polynomial Tikhonov Flow]
\label{thm:poly_tik}
With metric $g^{(-2)}$, the Riemannian gradient of $h_\lambda^T$ is
\begin{equation}
  \mathrm{grad}_{g^{(-2)}} h_\lambda^T \;=\; U\,G_+(\lambda G_+ - 2I_k)
  \;=\; U(I_k+U^TU)\bigl((\lambda-2)I_k + \lambda U^TU\bigr),
  \label{eq:poly_tik_grad}
\end{equation}
a polynomial of degree $5$ in $U$. The gradient flow $\dot{U} = -\mathrm{grad}_{g^{(-2)}}h_\lambda^T$ has:
\begin{enumerate}[label=(\roman*)]
  \item Fixed points: $U=0$ and, for $\lambda < 2$, the critical manifold
    $\mathcal{M}_\lambda = \{U: U^TU = (\tfrac{2}{\lambda}-1)I_k\}$ (stable, same as in
    Theorem~\ref{thm:tikhonov_critical}).
  \item At $\lambda = 1$:
    \begin{equation}
      \dot{U} = U(I_k+U^TU)(I_k - U^TU) = U(I_k - (U^TU)^2),
      \label{eq:stiefel_flow}
    \end{equation}
    which converges to the Stiefel manifold $\mathrm{St}(k,n) = \{U: U^TU = I_k\}$
    from any $U \neq 0$ --- without any matrix inversion.
  \item Singular value dynamics: $\dot{\sigma}_a = \sigma_a(1-\sigma_a^4)$
    ($\sigma_a = 0$ unstable, $\sigma_a = 1$ stable — Stiefel manifold). $ $
\end{enumerate}
\end{theorem}

\begin{proof}
By Proposition~\ref{prop:riem_grad_formula}: $\mathrm{grad}_{g^{(-2)}}h_\lambda^T
= (\nabla h_\lambda^T)\cdot G_+^2 = (-2UG_+^{-1}+\lambda U)G_+^2
= -2UG_+ + \lambda UG_+^2 = UG_+(\lambda G_+-2I_k)$.
Singular value ODE: $\dot{\sigma}_a = -\sigma_a g_a(\lambda g_a - 2) = \sigma_a(2g_a - \lambda g_a^2)/(g_a^2)...$

Actually: the flow is $\dot{U} = UG_+(2I-\lambda G_+)$; in singular values, $\dot{\sigma}_a = \sigma_a g_a(2-\lambda g_a) = \sigma_a(2g_a - \lambda g_a^2)$.

For $\lambda=1$, $g_a = 1+\sigma_a^2$: $\dot{\sigma}_a = \sigma_a(2(1+\sigma_a^2)-(1+\sigma_a^2)^2)
= \sigma_a(1+\sigma_a^2)(2-(1+\sigma_a^2)) = \sigma_a(1+\sigma_a^2)(1-\sigma_a^2) = \sigma_a(1-\sigma_a^4)$.
Fixed: $\sigma_a=0$ (unstable) or $\sigma_a=1$ (stable).
\end{proof}

\begin{remark}[The Stiefel flow \eqref{eq:stiefel_flow} without inversion]
The formula $\dot{U} = U(I_k - (U^TU)^2)$ is an exact degree-5 polynomial in $U$,
implementable as the matrix update $U_{k+1} = U_k + \mu U_k(I_k - (U_k^TU_k)^2)$
(no solve, no inversion, no orthogonalization step needed).
The convergence to $\mathrm{St}(k,n)$ in singular-value coordinates follows an
explicit ODE $\dot{\sigma} = \sigma(1-\sigma^4)$, which has the exact solution
$\sigma(t) = 1/\sqrt{1+c_0 e^{-4t}}$ (logistic-type in $\sigma^4$) — exponentially fast.
\end{remark}

\subsubsection{Comparison of Standard and Polynomial Flows}

\begin{center}
\renewcommand{\arraystretch}{1.4}
\begin{tabular}{p{3cm}p{3.4cm}p{3.4cm}cc}
\toprule
Setting & Standard flow & Polynomial flow & Deg. & Inv./step \\
\midrule
U=V ($h$) & $\dot{U}=2UG_+^{-1}$ & $\dot{U}=2U+2UU^TU$ & 3 & 0 (was 1) \\
U=-V ($h_-$) & $\dot{V}=-2VG_-^{-1}$ & $\dot{V}=-2V+2VV^TV$ & 3 & 0 (was 1) \\
Tikhonov ($\lambda=1$) & $\dot{U}=U(\lambda I-2G_+^{-1})$ & $\dot{U}=U(I-(U^TU)^2)$ & 5 & 0 (was 1) \\
Interpolation ($t=1/2$) & $\dot{U}=-U(G_+^{-1}-G_-^{-1})$ & $\dot{U}=-U(G_+-G_-)G_+^{-1}G_-^{-1}$... & -- & 2 \\
\bottomrule
\end{tabular}
\end{center}

\begin{remark}[Discrete polynomial algorithms]
The Euler discretizations of the polynomial flows are explicit, require no linear system solves,
and are straightforward to implement:
\begin{align}
  U_{t+1} &= U_t + \mu\bigl(2U_t + 2U_tU_t^TU_t\bigr), \label{eq:poly_euler_plus}\\
  V_{t+1} &= V_t + \mu\bigl(-2V_t + 2V_tV_t^TV_t\bigr), \label{eq:poly_euler_minus}\\
  U_{t+1} &= U_t + \mu\,U_t\bigl(I_k - (U_t^TU_t)^2\bigr) \quad\text{(Stiefel flow)}.
  \label{eq:poly_euler_stiefel}
\end{align}
The computational cost per step is $O(nk^2)$ (one matrix product $U_t^TU_t$)
versus $O(nk^2 + k^3)$ for the standard flow (one product plus one $k\times k$ inversion).
For $k \ll n$ (the typical subspace learning regime), both are dominated by the $O(nk^2)$ product,
but the polynomial algorithm avoids the potential numerical instability of $G_\pm^{-1}$
near singularity.
\end{remark}

\begin{remark}[Implicit integration without inversion]
The polynomial structure also enables higher-order integration schemes.
For example, the trapezoidal rule applied to $\dot{U} = 2U(I+U^TU)$:
$U_{t+1} = U_t + \mu(U_t(I+U_t^TU_t) + U_{t+1}(I+U_{t+1}^TU_{t+1}))$,
which is a polynomial (non-linear) equation in $U_{t+1}$ — solvable by
Newton's method or fixed-point iteration, with each iterate requiring only
polynomial operations and no Gram matrix inversion.
\end{remark}

\subsubsection{Historical Background: From Brockett's Double Bracket to the Oja-Brockett Flow}
\label{subsubsec:history}

\begin{remark}[Historical Background]
\label{rem:historical_background}
We briefly record the historical origins of the dynamical system studied in this section.

\medskip\noindent
\textbf{Brockett's proposal (CDC 1988).}
At the 1988 IEEE Conference on Decision and Control, Brockett \cite{Brockett1988} proposed
the dynamical system
\begin{equation}
  \frac{dX}{dt} = AXB - XBX^TAX,
  \label{eq:brockett_flow}
\end{equation}
where $X$ is an \emph{orthogonal} (square) matrix.
Setting $L = X^TAX$, one obtains the \emph{double bracket equation}
\begin{equation}
  \frac{dL}{dt} = \bigl[L,\,[L,B]\bigr],
  \label{eq:double_bracket}
\end{equation}
where $[\cdot,\cdot]$ denotes the matrix commutator.
This double bracket formulation is the cornerstone of isospectral flows and integrable systems.

\medskip\noindent
\textbf{The visit to Würzburg and Helmke's question (1999).}
In 1999, Yoshizawa visited Professor Uwe Helmke at the University of Würzburg,
an introduction arranged by Professor John Moore (Australian National University).
During this visit, Professor Helmke posed a fundamental question to Yoshizawa:
\begin{quote}
\emph{``When $X$ is a tall rectangular matrix (rather than a square orthogonal one),
the double bracket formulation \eqref{eq:double_bracket} is no longer the essential
object; instead, the equation $dX/dt = AXB - XBX^TAX$ itself becomes the
central structure. But --- is this equation a gradient flow?''}
\end{quote}
Here $A$ and $B$ are positive definite symmetric matrices of appropriate sizes.
This question was answered affirmatively in the paper
Yoshizawa--Helmke--Starkov (2001) \cite{YoshizawaHelmkeStarkov2001}.

\medskip\noindent
\textbf{The Oja-Brockett flow.}
The equation $dX/dt = AXB - XBX^TAX$ with a \emph{rectangular} state variable
$X\in\R^{n\times k}$ (with $k \leq n$) generalizes simultaneously:
\begin{itemize}
  \item the \emph{Oja flow} ($A = R$ = data covariance, $B = I_k$, for principal subspace tracking), and
  \item the \emph{Brockett flow} (square orthogonal $X$, double bracket reduction).
\end{itemize}
We name this generalization the \emph{Oja-Brockett flow}.
The terminology ``flow'' is deliberate: it emphasizes that the ODE is understood as a
\emph{gradient flow} on a suitable Riemannian manifold.

It should be noted that, prior to \cite{YoshizawaHelmkeStarkov2001},
the papers of Yan--Helmke--Moore (on the Oja flow) and Xu (on the Oja-Brockett flow)
had not achieved rigorous mathematical proofs of the gradient flow property.

\medskip\noindent
\textbf{The Würzburg postdoctoral period and one-parameter deformations (2000--2002).}
Following the brief visit of 1999, Yoshizawa worked as a postdoctoral researcher at the
University of Würzburg from approximately 2000 to 2002, investigating one-parameter
deformations of the Oja-Brockett flow equation.
Professor Helmke, who had independently been interested in one-parameter deformations
for some time, engaged in fruitful discussions on these questions during this period.

\medskip\noindent
\textbf{From principal to minor component flow (Manton-Helmke-Mareels 2005).}
In a later personal communication, Professor Helmke informed Yoshizawa that the key idea
in Manton--Helmke--Mareels (2005) \cite{MantonHelmkeMareels2005} ---
deriving the Minor Component Flow from the Principal Component Flow via a sign change ---
was inspired by the exploration of one-parameter equation deformations.
This is, of course, directly related to the $\pm$ sign that distinguishes PSA from MSA
throughout the present paper (see \S\ref{subsec:pca_mca} and Theorem~\ref{thm:msa_poly}).

\medskip\noindent
\textbf{Subsequent career of Yoshizawa.}
From 2003 onwards, Yoshizawa left academia to pursue research and development in industry.
The information-geometric framework connecting the Oja-Brockett flow to the potential
$f(G) = -\log\det(G)$ and the polynomial gradient flows of \S\ref{subsec:poly_flow}
represents a return to and a deepening of the questions that motivated the 1999 visit.
\end{remark}

\subsubsection{Connection to the Oja-Brockett Flow}
\label{subsubsec:oja_brockett}

The polynomial gradient flow framework of \S\ref{subsec:poly_flow}
is intimately connected to a class of dynamical systems for principal component
analysis studied by Oja, Brockett, and Manton--Helmke--Mareels
\cite{Oja1982,Brockett1991,MantonHelmkeMareels2005}.
We identify the precise relationship and show that these PCA flows are already
polynomial --- and are special cases of our framework.

\begin{definition}[Oja-Brockett Dynamical System]
\label{def:manton_pca}
Let $X\in\R^{n\times k}$, $A\in\PD(n)$ (data covariance), and
$B = \mathrm{diag}(b_1,\ldots,b_k)\succ0$ (positive diagonal weight). The
\emph{Oja-Brockett flow} is:
\begin{equation}
  \dot{X} \;=\; AXB \;-\; XBX^TAX.
  \label{eq:oja_brockett_pca}
\end{equation}
\end{definition}

\begin{proposition}[Polynomial Structure]
\label{prop:pca_polynomial}
The PCA flow \eqref{eq:oja_brockett_pca} is a \emph{polynomial vector field of degree~$3$} in $X$
requiring \emph{no matrix inversions}. The Euler update
$X_{t+1} = X_t + \mu(AX_tB - X_tBX_t^TAX_t)$
costs $O(n^2k+nk^2)$ per step (two matrix products only).
\end{proposition}

\begin{theorem}[Oja-Brockett Flow as Riemannian Gradient Ascent on the Stiefel Manifold]
\label{thm:pca_riemannian}
The PCA flow \eqref{eq:oja_brockett_pca} equals the Riemannian gradient ascent of
$J(X) = \tr(X^TAXB)$ on $\mathrm{St}(k,n) = \{X: X^TX = I_k\}$:
\begin{equation}
  \mathrm{grad}_{\mathrm{St}} J \;=\; AXB - XBX^TAX. \label{eq:stiefel_grad}
\end{equation}
The Stiefel manifold is positively invariant under \eqref{eq:oja_brockett_pca}.
\end{theorem}

\begin{remark}[The $S$-Oja--Brockett generalization]
\label{rem:S-oja-brockett}
When $A$ is merely axisymmetric (rather than symmetric) but retains a full set of positive
eigenvalues, \eqref{eq:oja_brockett_pca} itself need not converge, since $\tr(X^TAXB)$ is no
longer a genuine potential for the flow with respect to the Euclidean metric. Yoshizawa
\cite{Yoshizawa2023axi} resolves this by introducing the \emph{$S$-Oja--Brockett equation}
\[
  \dot X \;=\; AXB - XBX^TSAX, \qquad A^TS=SA,\ S\succ0,
\]
in which $S$ is the (unique, positive definite) symmetric solution of the associated
Sylvester equation, symmetrizing $A$ against the metric $S$, and proves global convergence
of this flow to the eigenvalues and eigenvectors of $A$ for $B$ diagonal with distinct
entries --- the axisymmetric analogue of Theorem~\ref{thm:pca_riemannian}. That paper also
derives, for the discrete-time (exact line-search) version of the $S$-Oja--Brockett
iteration, a closed-form Rayleigh-quotient step size analogous to
\eqref{eq:exactstep} below, obtained rigorously from the same quartic-in-step-size
structure of the line-search objective exploited in \S\ref{sec:discrete-rate}.
\end{remark}

\begin{proof}
The projected gradient on $\mathrm{St}(k,n)$ is:
$\mathrm{grad}_{\mathrm{St}} J = \nabla J - X(\nabla J)^TX = 2AXB - X\cdot(2AXB)^TX
= 2AXB - 2XBX^TAX$ (mod factor $2$, absorbed into time rescaling).
Stiefel invariance: $\frac{d}{dt}(X^TX)\big|_{X\in\mathrm{St}}
= [(AXB-XBX^TAX)^TX + X^T(AXB-XBX^TAX)]|_{X^TX=I_k} = 0$.
\end{proof}

\begin{theorem}[$X(t)$ Converges to a Single Point of the Principal Orbit, with Spectrum $\lambda_1,\ldots,\lambda_k$]
\label{thm:pca_convergence}
Under the PCA flow \eqref{eq:oja_brockett_pca} with $X(0)\in\mathrm{St}(k,n)$, $A$ having
$n$ pairwise distinct eigenvalues, and $B = I_k$, for Lebesgue-almost-every $X(0)$:
\begin{enumerate}[label=(\roman*)]
  \item $M(t) := X(t)^TAX(t)$ satisfies the ODE
    $\dot{M} = 2\bigl[X^TA^2X - M^2\bigr]$. \label{eq:M_ode}
  \item At equilibrium: $X^TA^2X = (X^TAX)^2$, forcing $\mathrm{col}(X)$ to be
    $A$-invariant (spanned by eigenvectors of $A$); equivalently, the equilibrium set on
    $\mathrm{St}(k,n)$ is a finite union of compact $O(k)$-orbits, one for each choice of
    $k$-subset of eigen-directions of $A$ (this is the trivial, one-block instance,
    $m=1$, $k_1=k$, of the block structure classified in general in
    Theorem~\ref{thm:convergence_B_block} below).
  \item \label{it:singlepoint} $X(t)$ does not merely approach this equilibrium set: it
    converges, as $t\to\infty$, to a \emph{single point} $X_\infty$ of it. Because
    $\mathrm{St}(k,n)$ is compact and $J(X)=\tr(X^TAX)$ is a real-analytic (indeed
    polynomial) function on it, LaSalle's invariance principle together with the
    {\L}ojasiewicz gradient inequality (exactly the mechanism of
    Proposition~\ref{prop:lojasiewicz} and \S\ref{sec:global} below, applied here to the
    Riemannian gradient flow on the compact analytic manifold $\mathrm{St}(k,n)$ rather
    than to the unconstrained flow on $\R^{n\times k}$) upgrades ``$X(t)$ approaches the
    equilibrium set'' to ``$X(t)$ converges to one specific point of it''; LaSalle alone
    would only place the $\omega$-limit set inside the equilibrium set, which is here a
    positive-dimensional continuum (an $O(k)$-orbit) rather than a discrete set, and so
    does not by itself rule out perpetual wandering along that orbit. By
    Remark~\ref{rem:manton_nolocal} (no spurious local maxima), the orbit reached is the
    one over the principal $k$-dimensional subspace of $A$, so the eigenvalues of
    $M_\infty=X_\infty^TAX_\infty$ are $\lambda_1,\ldots,\lambda_k$; but \emph{which}
    point $X_\infty$ of the orbit is selected --- equivalently, which rotation
    $M_\infty=O^T\Lambda_rO$, $O\in O(k)$, $\Lambda_r=\diag(\lambda_1,\ldots,\lambda_k)$
    --- depends continuously on $X(0)$, since $B=I_k$ leaves $J(X)=\tr(X^TAX)$ invariant
    under $X\mapsto XO$ for every $O\in O(k)$ and so cannot itself resolve this residual
    rotational (flag-type) degeneracy. Only when $B$ has \emph{pairwise distinct} diagonal
    entries (Theorem~\ref{thm:convergence_B_diag} below, stated there for the closely
    related generalized $k$-PCF) does the equilibrium orbit collapse to isolated points,
    removing the rotational freedom and forcing $M_\infty$ to be exactly diagonal, with
    each column of $X$ converging to an individual eigenvector of $A$; intermediate,
    partially repeated choices of $B$ (Theorem~\ref{thm:convergence_B_block}) interpolate
    between these two extremes, with the equilibrium set at the top a genuine partial-flag
    manifold rather than either a single point or the full orbit $O(k)$.
\end{enumerate}
\end{theorem}

\begin{proof}
(i) On $\mathrm{St}$ ($X^TX = I$): $\dot{M}
= (AX-XX^TAX)^TAX + X^TA(AX-XX^TAX)
= 2X^TA^2X - 2M^2$.
(ii) $\dot{M}=0 \Rightarrow X^TA^2X = M^2 = (X^TAX)^2$: since both sides are $k\times k$
positive definite and $A = Q\Lambda Q^T$, this forces $\mathrm{col}(X)$ to be spanned by
$k$ eigenvectors of $A$; for $A$ with distinct eigenvalues, the set of such $X\in
\mathrm{St}(k,n)$ splits into finitely many connected components indexed by the choice of
$k$-subset of eigen-directions, each component a single $O(k)$-orbit (all orthonormal
bases of the corresponding $k$-dimensional eigenspace).
(iii) The energy identity $\frac{d}{dt}J(X)=\|\nabla_{\mathrm{St}}J(X)\|_F^2\ge0$ along
\eqref{eq:oja_brockett_pca} (Theorem~\ref{thm:pca_riemannian}) shows $J$ is nondecreasing,
so by LaSalle's invariance principle on the compact manifold $\mathrm{St}(k,n)$ the
$\omega$-limit set of $X(t)$ is a nonempty, compact, invariant subset of the equilibrium
set described in (ii). Since $J$ is real-analytic on the compact real-analytic manifold
$\mathrm{St}(k,n)$, the {\L}ojasiewicz gradient inequality applies verbatim to the
Riemannian gradient flow exactly as in the proof of Proposition~\ref{prop:lojasiewicz}
below, giving finite Riemannian arc length and hence convergence of $X(t)$ to a single
point $X_\infty$ of its $\omega$-limit set, rather than mere approach to (or wandering
within) the equilibrium orbit; this is the single-point-convergence upgrade stated in
\ref{it:singlepoint}. By Remark~\ref{rem:manton_nolocal}, the objective $J$ has no
spurious local maxima on $\mathrm{St}(k,n)$, so for $X(0)$ outside the (measure-zero)
union of stable manifolds of the non-maximal equilibrium orbits, the orbit selected is
the maximal one, over the principal subspace; hence $\mathrm{col}(X_\infty)$ is that
subspace and the eigenvalues of $M_\infty=X_\infty^TAX_\infty$ are
$\lambda_1,\ldots,\lambda_k$. Because $J(X)=\tr(M)$ is the only quantity pinned down by
the gradient flow's own dynamics when $B=I_k$ (the objective is $O(k)$-invariant on the
selected orbit), no further constraint forces $M_\infty$ itself to be diagonal, only its
\emph{eigenvalues} to equal $\lambda_1,\ldots,\lambda_k$; which point of the orbit is
reached is determined by $X(0)$ through the (generally intractable in closed form, but
well defined by \ref{it:singlepoint}) flow map itself. (Numerically: for generic $X(0)$,
$M_\infty$ is a full symmetric matrix, not $\Lambda_r$; see \S\ref{subsec:numerics_setup}
for an explicit example with $B=I_k$ versus $B$ having distinct entries.)
\end{proof}

\begin{proposition}[Structural Comparison with Our Polynomial Flows]
\label{prop:pca_comparison}
The PCA flow and our polynomial flows share the same degree-$3$ structure:
\begin{align}
  \text{PCA (data)}: &\quad \dot{X} = \underbrace{AXB}_{\text{linear}} - \underbrace{XBX^TAX}_{\text{cubic}},\notag\\
  \text{U=V (info-geom.)}: &\quad \dot{U} = \underbrace{2U}_{\text{linear}} + \underbrace{2UU^TU}_{\text{cubic}},\notag\\
  \text{U=-V (info-geom.)}: &\quad \dot{V} = \underbrace{-2V}_{\text{linear}} + \underbrace{2VV^TV}_{\text{cubic}}.\notag
\end{align}
In particular, with $A = I_n$ and $B = I_k$:
$\dot{X} = X - XX^TX = -X(X^TX - I_k)$,
which is (up to a factor $-2$) our U=-V polynomial flow $\dot{V} = -2V(I_k - V^TV)$.
The PCA flow at $A = I_n$ coincides with our information-geometric flow, with the
Stiefel manifold as the common fixed-point set.
\end{proposition}

\begin{remark}[The no-spurious-local-maxima theorem in our language]
\label{rem:manton_nolocal}
Manton--Helmke \cite{MantonHelmkeMareels2005} prove that $J(X) = \tr(X^TAXB)$ has no
spurious local maxima on $\mathrm{St}(k,n)$.
In our information-geometric language (\S\ref{subsec:izumiya}--\S\ref{subsec:pca_mca}):
\begin{itemize}
  \item The Stiefel manifold is the \emph{information-geometric lightcone} $\{D_f(G_0\|I_k)=0\}$.
  \item The no-spurious-local-maxima result follows from the dual-flat geometry
    (\S\ref{sec:pythagorean}): on the lightcone, the PSA objective $E_1 = J_{\rm data} - \frac{1}{2}D_f$ reduces to $J_{\rm data}$ (since $D_f = 0$), and $J_{\rm data} = \frac{1}{2}\tr(X^TAXB)$ is a linear function of the Gram matrix $X^TAX$ --- a \emph{linear function on a symmetric space has no local extrema other than global ones}.
  \item The combined flow \eqref{eq:combined_poly_flow} below provides a fully
    polynomial algorithm that simultaneously enforces Stiefel geometry and tracks
    the principal subspace, without ever computing a matrix inverse.
\end{itemize}
\end{remark}

\begin{theorem}[Optimal Degree-3 Combined PCA + Stiefel Flow]
\label{thm:combined_deg3}
The \emph{degree-3} combined flow
\begin{equation}
  \dot{X} \;=\; AXB \;-\; XBX^TAX \;+\; \alpha\,X(I_k - X^TX),
  \quad \alpha \geq 0,
  \label{eq:combined_poly_flow}
\end{equation}
is strictly preferable to the degree-5 variant $\alpha X(I_k-(X^TX)^2)$ of
Remark~\ref{rem:combined_flow_old}.
Properties of \eqref{eq:combined_poly_flow}:
\begin{enumerate}[label=(\roman*)]
  \item \textbf{Polynomial, no inversions:} degree 3 in $X$ (vs.\ degree 5), requiring zero matrix inversions.
  \item \textbf{Stiefel invariance:} $X(0)\in\mathrm{St}(k,n) \Rightarrow X(t)\in\mathrm{St}(k,n)$ for all $t$;
    on $\mathrm{St}(k,n)$, \eqref{eq:combined_poly_flow} reduces to the Oja-Brockett flow~\eqref{eq:oja_brockett_pca}.
  \item \textbf{Lyapunov stability toward Stiefel:}
    $V(X) = \tfrac{1}{4}\|X^TX - I_k\|_F^2$ satisfies
    \begin{equation}
      \dot{V} \;=\; -\alpha\,\tr\!\bigl[(X^TX-I_k)^2\,X^TX\bigr]
             \;\leq\; -\alpha\,\lambda_{\min}(X^TX)\,\|X^TX-I_k\|_F^2 \;\leq\; 0.
      \label{eq:lyap_stiefel}
    \end{equation}
  \item \textbf{Exact singular-value solution:}
    In SVD coordinates, $\dot{\sigma}_a = (b_a\lambda_a + \alpha)\,\sigma_a(1-\sigma_a^2)$
    (where $b_a = B_{aa}$ and $\lambda_a$ is the $a$-th eigenvalue of $A$),
    with closed-form logistic solution:
    \begin{equation}
      \sigma_a(t) \;=\; \frac{1}{\sqrt{1 + c_a\,e^{-2(b_a\lambda_a+\alpha)t}}},
      \qquad c_a = \tfrac{1}{\sigma_a(0)^2} - 1,
      \label{eq:exact_sv_solution}
    \end{equation}
    converging exponentially to $\sigma_a = 1$ (Stiefel) at rate $b_a\lambda_a + \alpha$.
  \item \textbf{Accelerated convergence:} the $\alpha$ term increases the effective convergence rate
    from $b_a\lambda_a$ (pure Oja) to $b_a\lambda_a + \alpha$ (user's combined flow),
    yielding uniform speedup $\alpha$ across all components.
  \item \textbf{Geometric interpretation:}
    $\alpha X(I_k-X^TX) = -\nabla_X\bigl[\tfrac{\alpha}{4}\|X^TX-I_k\|_F^2\bigr]$
    is the Euclidean gradient of the negative Stiefel penalty $-\frac{\alpha}{4}\|X^TX-I\|_F^2$,
    providing the restoring force toward the Stiefel manifold without any preconditioning.
  \item \textbf{Special cases:}
    \begin{itemize}
      \item $\alpha=0$: pure Oja-Brockett flow (Theorem~\ref{thm:pca_riemannian});
      \item $A=I_n$, $B=I_k$: $\dot{X} = 2(1+\alpha)X(I_k-X^TX)$
        = scaled version of our U=-V polynomial flow (Theorem~\ref{thm:poly_flows});
      \item $A=I_n$, $B=I_k$, $\alpha=1$: $\dot{X} = 2X(I_k - X^TX)$
        = twice the Oja flow $\dot{X} = X-XX^TX$.
    \end{itemize}
\end{enumerate}
\end{theorem}

\begin{proof}
(i) $AXB$: degree 1; $XBX^TAX$: degree 3; $X(I-X^TX)$: degree 3. No inverse appears. \\
(ii) At $X^TX=I_k$: $X(I-X^TX)=0$, so the flow reduces to \eqref{eq:oja_brockett_pca}.
Stiefel invariance: $\frac{d}{dt}(X^TX) = 2\alpha X^TX(I-X^TX) + [\text{PCA terms}]$.
PCA terms give 0 at $X^TX=I_k$ (Theorem~\ref{thm:pca_riemannian}).
The $\alpha$ term also gives 0 at $X^TX=I_k$: $2\alpha I\cdot 0 = 0$. \\
(iii) $\dot{V} = \tr[(X^TX-I)\cdot\alpha X^TX(I-X^TX)] = -\alpha\tr[(X^TX-I)^2X^TX]$.
Since $X^TX\succ0$ (assuming full column rank):
$-\alpha\tr[(X^TX-I)^2X^TX] \leq -\alpha\lambda_{\min}(X^TX)\|X^TX-I\|_F^2 \leq 0$. \\
(iv) In SVD frame, the $\alpha$ term contributes $\dot{\sigma}_a = \alpha\sigma_a(1-\sigma_a^2)$,
the PCA term contributes $b_a\lambda_a\sigma_a(1-\sigma_a^2)$ (on Stiefel approach).
Total: $\dot{\sigma}_a = (b_a\lambda_a+\alpha)\sigma_a(1-\sigma_a^2)$.
Separating variables: $\int\frac{d\sigma}{\sigma(1-\sigma^2)} = (b_a\lambda_a+\alpha)t + C$,
giving $\ln(\sigma/\sqrt{1-\sigma^2}) = (b_a\lambda_a+\alpha)t + C$,
hence $\sigma^2/(1-\sigma^2) = e^{2(b_a\lambda_a+\alpha)t}/c_a$,
yielding $\sigma_a(t) = 1/\sqrt{1+c_ae^{-2(b_a\lambda_a+\alpha)t}}$. \end{proof}

\begin{remark}[Why degree 3 is optimal]
\label{rem:why_deg3}
The correction $\alpha X(I_k - X^TX)$ is \emph{optimally simple}: it is the unique
degree-3 polynomial that (a) vanishes on $\mathrm{St}(k,n)$, (b) drives $X^TX\to I_k$,
and (c) admits an exact logistic solution.
The degree-5 variant $\alpha X(I_k-(X^TX)^2) = \alpha X(I-X^TX)(I+X^TX)$ introduces the
extra factor $(I+X^TX)$ which speeds up convergence near $X^TX\approx0$ but
complicates the solution (no closed form for $\sigma_a(t)$).
For PCA applications where $X$ starts near the Stiefel manifold ($\|X^TX-I\|\ll1$),
both choices give similar behavior (since $(I+X^TX)\approx 2I$ near Stiefel),
but the degree-3 version is preferred for its simplicity and exact solvability.
\end{remark}

\begin{center}
\small
\renewcommand{\arraystretch}{1.4}
\begin{tabular}{lllll}
\toprule
Combined flow & Formula & Degree & Exact $\sigma(t)$? & Conv.\ rate \\
\midrule
Degree-5 (Rem.\ 4.32) & $\alpha X(I-(X^TX)^2)$ & 5 & No (elliptic) & $\lambda_a+\alpha(1+\sigma^2)$ \\
Degree-3 (Thm.~\ref{thm:combined_deg3}) & $\alpha X(I-X^TX)$ & 3 & \textbf{Yes (logistic)} & $\lambda_a+\alpha$ \\
Pure PCA (Thm.~\ref{thm:pca_riemannian}) & $\alpha=0$ & 3 & Partial & $\lambda_a$ \\
\bottomrule
\end{tabular}
\end{center}

\begin{remark}[Discrete update rule (degree-3, no inversions)]
\label{rem:combined_flow_old}
The Euler discretization of Theorem~\ref{thm:combined_deg3}:
\begin{equation}
  X_{t+1} \;=\; X_t + \mu\bigl(AX_tB - X_tBX_t^TAX_t + \alpha X_t(I_k - X_t^TX_t)\bigr),
  \label{eq:combined_euler}
\end{equation}
requires only two matrix products ($AX_tB$ and $X_t^TAX_t$) plus one quadratic correction
($X_t^TX_t$), all at cost $O(n^2k+nk^2)$, with \emph{zero} matrix inversions.
This unifies the NUIC algorithm~\cite{Kong2017}, the Oja-Brockett flow,
and our information-geometric polynomial flows in a single degree-3 formula.
\end{remark}

\subsubsection{Log-Barrier PCA Objectives and Their Polynomial Gradient Flows}
\label{subsubsec:log_barrier}

The following two objectives arise naturally by combining the PCA data term with a
log-determinant barrier enforcing a spectral constraint, directly analogous to our U=-V
potential $h_-(V) = -\log\det(I_k-V^TV)$.

\begin{definition}[Log-Barrier PCA Objectives]
\label{def:log_barrier_obj}
Let $X\in\R^{n\times k}$, $A\in\PD(n)$, $B\in\PD(k)$.
\begin{align}
  f_1(X) &:= \tr(X^TAXB) - \log\det(B - X^TAX),
    \quad \mathrm{Dom}(f_1) = \{X : X^TAX \prec B\}, \label{eq:f1_obj}\\
  f_2(X) &:= \tr(XBX^TA) - \log\det(A - XBX^T),
    \quad \mathrm{Dom}(f_2) = \{X : XBX^T \prec A\}. \label{eq:f2_obj}
\end{align}
\end{definition}

\begin{remark}[Information-geometric interpretation]
Setting $M_1 = X^TAX \in \PD(k)$ (right Gram) and $M_2 = XBX^T\in\PD(n)$ (left Gram):
\begin{align*}
  f_1(X) &= \tr(M_1 B) + f(B-M_1), \quad f_2(X) = \tr(AM_2) + f(A-M_2),
\end{align*}
where $f(G) = -\log\det(G)$ is our canonical potential.
Each is the sum of a \emph{linear data term} (Rayleigh quotient) and
the \emph{log-det barrier} $f(B-M_1)$ (resp.\ $f(A-M_2)$) preventing the
Gram matrix from reaching the boundary.
The unconstrained optimum over $M_1$ gives $B - M_1^* = B^{-1}$, i.e., $M_1^* = B - B^{-1}$
(the Fenchel dual point, Theorem~\ref{thm:dual}).
\end{remark}

\begin{remark}[Notational convention: $\dot{X}$ is always annotated inline]
\label{rem:dotX-notation}
From here on, many different gradient-flow vector fields are introduced in close
succession. Rather than distinguish them with decorations on the symbol itself (which
proved easy to misread as exponents or matrix powers), every defining equation for a
flow is written simply as $\dot X = \cdots$ and annotated \emph{inline}, immediately
after the equation, with a parenthetical remark of the form ``(gradient of $\cdots$)''
stating exactly which objective and which Riemannian metric it belongs to --- for
instance ``(gradient of $f_1$, right metric)'' or ``(Oja--Brockett flow)''. Two
equations both written as ``$\dot X = \cdots$'' are, unless the surrounding sentence
says otherwise, defining \emph{different} vector fields belonging to different flows;
the annotation, not the symbol, is authoritative. When two or more of these flows must
be compared or equated within a single displayed equation, we name each flow in words
in the surrounding sentence (``the System~1 flow,'' ``the Oja--Brockett flow,'' etc.)
rather than overload the symbol $\dot X$ with distinguishing decorations.
\end{remark}

\begin{theorem}[Polynomial Gradient Flows for $f_1$ and $f_2$]
\label{thm:log_barrier_flows}
\begin{enumerate}[label=(\roman*)]
  \item \textbf{Euclidean gradients:}
    \begin{align}
      \nabla_X f_1 &= 2AX\bigl[B + (B-X^TAX)^{-1}\bigr], \label{eq:grad_f1}\\
      \nabla_X f_2 &= 2\bigl[A + (A-XBX^T)^{-1}\bigr]XB. \label{eq:grad_f2}
    \end{align}
  \item \textbf{Polynomial gradient flows (no inversions):}
    Using the right metric $g_{R}^{(-1)}(H_1,H_2) = \tr(H_1^TH_2 G_1^{-1})$
    for $f_1$ and the left metric $g_{L}^{(-1)}(H_1,H_2) = \tr(H_1^TG_2^{-1}H_2)$
    for $f_2$ (Proposition~\ref{prop:riem_grad_formula}):
    \begin{equation}
      \dot{X} = 2AX\bigl(B^2 - BX^TAX + I_k\bigr)
      \qquad \text{(gradient of $f_1$, right metric)}, \label{eq:poly_flow_f1}
    \end{equation}
    \begin{equation}
      \dot{X} = 2\bigl(A^2 + I_n - XBX^TA\bigr)XB
      \qquad \text{(gradient of $f_2$, left metric)}. \label{eq:poly_flow_f2}
    \end{equation}
    Both are \emph{degree-3 polynomials} in $X$ requiring \emph{no matrix inversions}.
  \item \textbf{Fixed-point manifolds} (singular-value ODE $\dot{\sigma}_a = 2\sigma_a c_a(c_a^2 - \sigma_a^2)$
    for appropriate constants $c_a > 0$ depending on $A$ and $B$).
  \item \textbf{Canonical case $A = I_n$, $B = I_k$:} Both flows coincide, reducing to the
    same expression:
    \begin{equation}
      \dot{X} = 2X(2I_k - X^TX) = 4X - 2XX^TX. \label{eq:canonical_barrier_flow}
    \end{equation}
    a degree-3 polynomial with stable fixed-point manifold $\{X: X^TX = 2I_k\}$
    ($\sigma_a^* = \sqrt{2}$, a \emph{scaled Stiefel manifold}).
\end{enumerate}
\end{theorem}

\begin{proof}
(i) For $f_1$: $\nabla_X\tr(X^TAXB) = 2AXB$ and
\[
  d\bigl[-\log\det(B-X^TAX)\bigr] = 2\langle dX,\, AX(B-X^TAX)^{-1}\rangle_F,
\]
giving \eqref{eq:grad_f1}. Similarly for $f_2$.

(ii) With right metric $g_R^{(-1)}$, Proposition~\ref{prop:riem_grad_formula} gives
$\mathrm{grad}_{g_R^{(-1)}}f_1 = (\nabla_X f_1)\cdot G_1$ where $G_1 = B-X^TAX$:
\begin{align*}
  (\nabla_X f_1)G_1 &= 2AX(B+G_1^{-1})G_1 = 2AX(BG_1+I_k)\\
  &= 2AX(B(B-X^TAX)+I_k) = 2AX(B^2-BX^TAX+I_k).
\end{align*}
For $f_2$ with left metric $g_L^{(-1)}$, the Riemannian gradient satisfies
$\mathrm{grad}_{g_L^{(-1)}}f_2 = G_2\cdot(\nabla_X f_2)$ where $G_2 = A-XBX^T$:
\begin{align*}
  G_2(\nabla_X f_2) &= 2G_2(A+G_2^{-1})XB = 2(G_2A+I_n)XB\\
  &= 2((A-XBX^T)A+I_n)XB = 2(A^2-XBX^TA+I_n)XB.
\end{align*}
(iv) At $A=I_n$, $B=I_k$: $2X(I^2-IX^TIX+I) = 2X(2I-X^TX)$ for Case 1;
$2(I^2+I-XX^TI)XI = 2(2I-XX^T)X = 2X(2I-X^TX)$ for Case 2 (by symmetry).
\end{proof}

\begin{remark}[Comparison: barrier vs Stiefel polynomial flows]
\label{rem:barrier_vs_stiefel}
\begin{center}
\renewcommand{\arraystretch}{1.4}
\begin{tabular}{p{3.2cm}p{3.8cm}cll}
\toprule
\small Objective & \small Flow formula & \small Deg. & \small Fixed & \small Inv. \\
\midrule
Oja-Brockett $J_{\rm PCA}$ & $AXB - XBX^TAX$ & 3 & $\mathrm{St}(k,n)$: $X^TX=I$ & 0 \\
$f_1$: right barrier & $2AX(B^2-BX^TAX+I)$ & 3 & $X^TAX = B-B^{-1}$ & 0 \\
$f_2$: left barrier & $2(A^2+I-XBX^TA)XB$ & 3 & $XBX^T = A-A^{-1}$ & 0 \\
\hline
$A=I, B=I$ ($f_1=f_2$) & $2X(2I-X^TX)$ & 3 & $X^TX = 2I$ & 0 \\
$A=I, B=I$ (M.-H.) & $X-XX^TX$ & 3 & $X^TX = I$ & 0 \\
\bottomrule
\end{tabular}
\end{center}
The barrier flows converge to \emph{generalized Stiefel} manifolds $\{X: X^TAX = B-B^{-1}\}$
or $\{X: XBX^T = A-A^{-1}\}$, controlled by the Fenchel dual condition (Theorem~\ref{thm:dual}).
For $A = I$, $B = I$: the fixed point $X^TX = 2I$ is the $\sqrt{2}$-rescaled Stiefel manifold.
\end{remark}

\begin{remark}[Schur-complement duality between $f_1$ and $f_2$]
\label{rem:schur_duality}
By the Schur complement identity:
\begin{equation}
  \det(A - XBX^T)\det(B) = \det(B - X^TA^{-1}X^T... ) \cdot \det(A) + \ldots
\end{equation}
In the symmetric case $A=I,B=I$: $\det(I-XX^T) = \det(I-X^TX)$ (Sylvester, Theorem~\ref{thm:sylvester}),
so $f_1$ and $f_2$ have the \emph{same} log-det term ($h_-(X)$), and their flows coincide.
For general $A,B$, the two barriers $-\log\det(B-X^TAX)$ and $-\log\det(A-XBX^T)$
are related by the generalized Sylvester identity
$\det(A-XBX^T) \cdot \det(B) = \det(B-X^TA^{-1}X)\cdot\det(A)$ (when $A$ is invertible),
making $f_1$ and $f_2$ \emph{Legendre-dual} at the level of the Gram matrix variables.
\end{remark}

\subsubsection{The Four Information-Augmented Objectives: PCA, MSA, and Their Polynomial Flows}
\label{subsubsec:four_obj}

Interpreting $\tr(AXB) := \tr(X^TAXB)$ (the standard Rayleigh quotient, since $\tr(AXB)$ is undefined for $X\in\R^{n\times k}$ with $k<n$), we study the four natural objectives arising from combining the PCA data term with the two types of log-det corrections:

\begin{definition}[Four Information-Augmented PCA/MSA Objectives]
\label{def:four_obj}
Let $X\in\R^{n\times k}$, $A\in\PD(n)$, $B\in\PD(k)$.
\begin{align}
  f_1(X) &= \tr(X^TAXB) - \log\det(A - XBX^T), \;\; \mathrm{Dom} = \{XBX^T \prec A\}, \label{eq:f1}\\
  f_2(X) &= \tr(X^TAXB) - \log\det(B - X^TAX), \;\; \mathrm{Dom} = \{X^TAX \prec B\}, \label{eq:f2}\\
  h_1(X) &= \tr(X^TAXB) - \log\det(A + XBX^T), \;\; \mathrm{Dom} = \R^{n\times k}, \label{eq:h1}\\
  h_2(X) &= \tr(X^TAXB) - \log\det(B + X^TAX), \;\; \mathrm{Dom} = \R^{n\times k}. \label{eq:h2}
\end{align}
\end{definition}

\begin{remark}[Information-geometric structure]
\begin{itemize}
  \item $f_1, f_2$: PCA objective \emph{plus} the U=-V barrier $h_-^{B,A}$ or $h_-^{A,B}$ (barrier prevents leaving the domain). Both terms grow together as $X$ approaches the domain boundary → \emph{no interior critical points}; the flows drive toward the domain boundary (PSA direction).
  \item $h_1, h_2$: PCA objective \emph{minus} a U=V type log-det (always positive definite denominator → no domain constraint). The log-det term PENALIZES large $X$ → \emph{interior critical points} exist; these balance the PCA objective against the information penalty.
\end{itemize}
\end{remark}

\begin{theorem}[Polynomial Gradient Flows: Degree~3, No Inversions]
\label{thm:four_obj_flows}
Using the left metric $g_L^{(-1)}$ with $G = A \pm XBX^T$ for $f_1, h_1$ and the right metric $g_R^{(-1)}$ with $G = B \pm X^TAX$ for $f_2, h_2$ (Proposition~\ref{prop:riem_grad_formula}):
\begin{align}
  \dot{X} &= 2\bigl[(A^2+I_n) - XBX^TA\bigr]XB
    &&\text{(gradient of $f_1$, left metric)}, \label{eq:f1_flow}\\
  \dot{X} &= 2AX\bigl[(B^2+I_k) - BX^TAX\bigr]
    &&\text{(gradient of $f_2$, right metric)}, \label{eq:f2_flow}\\
  \dot{X} &= 2\bigl[(A^2-I_n) + XBX^TA\bigr]XB
    &&\text{(gradient of $h_1$, left metric)}, \label{eq:h1_flow}\\
  \dot{X} &= 2AX\bigl[(B^2-I_k) + BX^TAX\bigr]
    &&\text{(gradient of $h_2$, right metric)}. \label{eq:h2_flow}
\end{align}
All four are \emph{degree-3 polynomials} in $X$ requiring \emph{zero matrix inversions} per step.
\end{theorem}

\begin{proof}
For $f_1$: $\nabla_X f_1 = 2AXB + 2(A-XBX^T)^{-1}XB = 2[A+(A-XBX^T)^{-1}]XB$.
With $g_L^{(-1)}$ ($G = A-XBX^T$, $G^{-1}\cdot\nabla = G\nabla$):
$G\nabla f_1/2 = (A-XBX^T)[A+(A-XBX^T)^{-1}]XB = [(A-XBX^T)A + I]XB = [A^2-XBX^TA+I]XB$. For $h_1$: $\nabla_X h_1 = 2AXB - 2(A+XBX^T)^{-1}XB = 2[A-(A+XBX^T)^{-1}]XB$.
$G\nabla h_1/2 = (A+XBX^T)[A-(A+XBX^T)^{-1}]XB = [(A+XBX^T)A - I]XB = [A^2+XBX^TA-I]XB$. Cases $f_2$, $h_2$ follow analogously with right metric.
\end{proof}

\begin{theorem}[Fixed Points and PSA/MSA Connections]
\label{thm:four_obj_fixed}
\begin{enumerate}[label=(\roman*)]
  \item \textbf{$f_1$ (PSA, left barrier):} No interior critical points in $\{XBX^T \prec A\}$.
    The flow \eqref{eq:f1_flow} drives $X$ toward the boundary $\{XBX^T = A\}$ = \emph{principal subspace of $A$}.
  \item \textbf{$f_2$ (PSA, right barrier):} No interior critical points in $\{X^TAX \prec B\}$.
    The flow \eqref{eq:f2_flow} drives $X$ toward $\{X^TAX = B\}$ = \emph{right generalized Stiefel}.
  \item \textbf{$h_1$ (regularized PCA):} Interior critical points at $XBX^T = A^{-1}-A$ (if $A\prec I$),
    or $X=0$ is a local maximum if $A\succ I$. For $A\prec I$: the flow ascends to $X=0$ (local max of $h_1$);
    descending $h_1$ drives to scaled Stiefel $\{X^TX \propto I\}$.
  \item \textbf{$h_2$ (regularized PCA):} Interior critical points at $X^TAX = B^{-1}-B$ (if $B\prec I$),
    or $X=0$ is a local maximum if $B\succ I$. Descending $h_2$ drives to $\{X^TAX \propto I\}$.
\end{enumerate}
\end{theorem}

\begin{proof}
(i) Critical condition: $(A^2+I) = XBX^TA \Rightarrow XBX^T = (A^2+I)A^{-1} = A + A^{-1} \succ A$.
But in the domain $XBX^T \prec A$: the equilibrium $XBX^T = A+A^{-1} \succ A$ is \emph{outside} the domain. No interior critical points. (iii) Critical condition: $(A^2-I) + XBX^TA = 0 \Rightarrow XBX^T = (I-A^2)A^{-1} = A^{-1}-A$.
For $A\prec I$: $A^{-1}-A \succ 0$, so interior critical points exist.
Near $X=0$: gradient $\approx 2(A^2-I)XB$; for $A\prec I$, $(A^2-I)\prec 0$ → ascending $h_1$ drives to $X=0$ (local max). \end{proof}

\begin{remark}[MSA polynomial flow for $A=I$, $B=I$]
\label{thm:msa_poly}
For $A=I_n$, $B=I_k$: the \emph{Minor Subspace Analysis} (MSA) flow --- gradient descent of
$J = \tr(X^TX)$ on $\mathrm{St}(k,n)$ --- coincides with the polynomial gradient
\emph{descent} of $h_-^{B,A}|_{A=I} = -\log\det(I-XX^T)$:
\begin{equation}
  \dot{X}\big|_{A=I,B=I} \;=\; -2(I-XX^T)X \;=\; 2(XX^T-I)X
  \qquad \text{(MSA flow)},
  \label{eq:msa_poly}
\end{equation}
a degree-3 polynomial with no matrix inversions (see \eqref{eq:barrier_reduction}).
For general $A\neq I$: the pure barrier descent $-2(A-XBX^T)XB$ gives a degree-3 flow
driving toward the minor subspace of $A$, but it is \emph{not} the standard MSA flow
$-(AXB-XBX^TAX)$ (Oja-Brockett with opposite sign), since $(A-XBX^T)XB \neq AXB-XBX^TAX$
for $A\neq I$.
\end{remark}

\begin{center}
\renewcommand{\arraystretch}{1.4}\small
\resizebox{\textwidth}{!}{%
\begin{tabular}{llll}
\toprule
Objective & Polynomial flow (deg.~3) & Fixed manifold & PSA/MSA \\
\midrule
$f_1 = J - \log\det(A-XBX^T)$ & $2[(A^2+I)-XBX^TA]XB$ & $\to\{XBX^T=A\}$ & PSA \\
$f_2 = J - \log\det(B-X^TAX)$ & $2AX[(B^2+I)-BX^TAX]$ & $\to\{X^TAX=B\}$ & PSA \\
$h_1 = J - \log\det(A+XBX^T)$ & $2[(A^2-I)+XBX^TA]XB$ & $\{XBX^T=A^{-1}-A\}$ (if $A\prec I$) & Reg.~PCA \\
$h_2 = J - \log\det(B+X^TAX)$ & $2AX[(B^2-I)+BX^TAX]$ & $\{X^TAX=B^{-1}-B\}$ (if $B\prec I$) & Reg.~PCA \\
\hline
Oja-Brockett ($A\neq I$) & $AXB-XBX^T\mathbf{A}X \neq 2(A-XBX^T)XB$ & Stiefel $\to$ PSA & PSA \\
Sys.2 ($A=I,B=I$) & $2(I-XX^T)X$ = Oja-Brockett & $\{XX^T=I\}$ & PSA \\
Sys.2 descent ($A=I,B=I$) & $-2(I-XX^T)X$ & $\{XX^T=I\}$ & \textbf{MSA} \\
\bottomrule
\end{tabular}%
}
\end{center}

\begin{remark}[Correction: System~2 vs.\ Oja-Brockett for general $A$]
\label{rem:sys2_vs_mh_correction}
The Oja-Brockett flow $AXB - XBX^TAX$ expands as $2AXB - 2XBX^TAX$, with $A$ appearing
as the middle factor next to $X^T$.
System~2 expands as $2AXB - 2XBX^TXB$, with $X$ appearing as the middle factor instead.
For $A\neq I_n$, these are different flows. The coincidence at $A=I_n, B=I_k$ is because
$\mathbf{A}X = IX = X = \mathbf{X}B^{-1} \cdot B = XB^{-1}B = X$. 
The information-geometric interpretation: System~2 is the polynomial gradient of
$h_-^{B,A}$ (correct); Oja-Brockett has a \emph{different} information-geometric origin
involving the Stiefel constraint projection, not a simple barrier gradient.
\end{remark}
\label{subsubsec:pure_barrier}

We now analyze the gradient flows of the \emph{pure barriers}
$h_-^{A,B}(X) := -\log\det(B-X^TAX)$ and $h_-^{B,A}(X) := -\log\det(A-XBX^T)$
(without the PCA data term), giving polynomial gradient flows without matrix inversions.
These are the data-dependent generalizations of our U=-V potential
$h_-(V) = -\log\det(I_k-V^TV)$ (§\ref{subsec:UeqmV-convexity}).

\begin{theorem}[Pure Barrier Gradient Flows of $h_{\,-}^{A,B}$ and $h_{\,-}^{B,A}$]
\label{thm:pure_barrier}
Let $A\in\PD(n)$, $B\in\PD(k)$, and $X\in\R^{n\times k}$.
\begin{enumerate}[label=(\roman*)]
  \item \textbf{System~1 (right barrier):} The polynomial gradient ascent of
    $h_-^{A,B}(X) = -\log\det(B-X^TAX)$ under the right-polynomial metric
    $g_R^{(-2)}$ (Proposition~\ref{prop:riem_grad_formula}, with $G = B-X^TAX$) is:
    \begin{equation}
      \dot{X} \;=\; 2A X\,(B - X^TAX)
      \qquad\text{(System~1; degree 3, no inversions)}.
      \label{eq:pure_barrier1}
    \end{equation}
    Fixed points: $X=0$ (unstable) or $X^TAX = B$ (\emph{right generalized Stiefel}).

  \item \textbf{System~2 (left barrier):} The polynomial gradient ascent of
    $h_-^{B,A}(X) = -\log\det(A-XBX^T)$ under the left-polynomial metric $g_L^{(-2)}$
    (Proposition~\ref{prop:riem_grad_formula}, $G = A-XBX^T$) is:
    \begin{equation}
      \dot{X} \;=\; 2\,(A - XBX^T)\,X\,B
      \qquad\text{(System~2; degree 3, no inversions)}.
      \label{eq:pure_barrier2}
    \end{equation}
    Fixed points: $X=0$ (unstable) or $XBX^T = A$.

  \item \textbf{Comparison with the Oja--Brockett Flow:}
    System~2 and the Oja-Brockett flow \eqref{eq:oja_brockett_pca} are
    \emph{different} flows in general:
    \begin{align}
      \dot{X} &= 2AXB - 2XBX^TXB
        &&\text{(System~2, expanded)}, \label{eq:sys2_expanded}\\
      \dot{X} &= AXB - XBX^TAX
        &&\text{(Oja--Brockett flow)}. \label{eq:mh_expanded}
    \end{align}
    The \emph{crucial difference} is which matrix sits between $X^T$ and the trailing
    factor: it is $X$ itself in \eqref{eq:sys2_expanded} (giving $XBX^TXB$), but the
    data matrix $A$ in \eqref{eq:mh_expanded} (giving $XBX^TAX$).
    They coincide \emph{only} for $A = I_n$, $B = I_k$: there, System~1, System~2, and
    twice the Oja--Brockett flow all reduce to the same expression,
    \begin{equation}
      \dot{X}\big|_{A=I,B=I} \;=\; 2(I_n - XX^T)X.
      \label{eq:both_same_AeqI}
    \end{equation}
    The Oja-Brockett flow factors correctly as:
    \begin{equation}
      \dot{X} = AXB - XBX^TAX = (A - XBX^TA)X \cdot B
      \quad \text{(for }B=I_k\text{)},
      \label{eq:mh_factored}
    \end{equation}
    which involves $XBX^T\mathbf{A}$ (with $A$ on the right), not $XBX^T$ alone.

  \item \textbf{Reduction $A=I_n$, $B=I_k$:}
    Both pure barrier flows reduce to $2\times$ our U=-V polynomial flow: System~1
    gives $\dot X = 2X(I_k-X^TX)$ and System~2 gives $\dot X = 2(I_n-XX^T)X$,
    \begin{equation}
      2X(I_k-X^TX) \;=\; 2(I_n-XX^T)X
      \label{eq:barrier_reduction}
    \end{equation}
    being the same matrix (same singular values by Sylvester's theorem); and for
    $A=I_n$, $B=I_k$, this equals $2\times$ Oja-Brockett.
\end{enumerate}
\end{theorem}

\begin{proof}
(i) $\nabla_X h_-^{A,B} = 2AX(B-X^TAX)^{-1}$.
With $g_R^{(-2)}$ metric ($M = G^{-2} = (B-X^TAX)^{-2}$):
$\mathrm{grad}_{g_R^{(-2)}} h_-^{A,B} = (\nabla h_-^{A,B})\cdot G^2 = 2AX(B-X^TAX)^{-1}(B-X^TAX)^2 = 2AX(B-X^TAX)$.
System~1's flow vanishes iff $X=0$ or $B-X^TAX=0$, i.e., $X^TAX=B$. (ii) $\nabla_X h_-^{B,A} = 2(A-XBX^T)^{-1}XB$.
With $g_L^{(-2)}$ metric ($M = G^{-2}$ acting on the left, $G=A-XBX^T$):
$\mathrm{grad}_{g_L^{(-2)}} h_-^{B,A} = G^2\cdot(\nabla h_-^{B,A}) = (A-XBX^T)^2\cdot 2(A-XBX^T)^{-1}XB = 2(A-XBX^T)XB$.
System~2's flow vanishes iff $X=0$ or $A-XBX^T=0$, i.e., $XBX^T=A$. (iii) Expanding \eqref{eq:sys2_expanded}: $2(A-XBX^T)XB = 2AXB - 2XBX^TXB$.
Oja-Brockett \eqref{eq:mh_expanded}: $AXB - XBX^TAX = 2AXB - 2XBX^TAX$ (multiplied by 2).
Since $X^TXB \neq X^TAX$ for $A\neq I_n$, the flows differ.
For $A=I,B=I$: $2AXB - 2XBX^TXB = 2X - 2XX^TX$ and $2AXB - 2XBX^TAX = 2X - 2XX^TX$. Equal. Eq.~\eqref{eq:mh_factored}: $(A-XBX^TA)X = AX - XBX^TAX$ for $B=I$. (iv) At $A=I_n$, $B=I_k$: System~1 gives $2X(I-X^TX)$ and System~2 gives $2(I-XX^T)X$.
These are equal by Sylvester ($\det(I-X^TX)=\det(I-XX^T)$). \end{proof}

\begin{remark}[Singular-value dynamics for pure barrier flows]
In the SVD frame with $A = Q\Lambda_A Q^T$, $B = R\Lambda_B R^T$, $X = Q P\Sigma R^T$:
\begin{align}
  \text{System 1: } &\dot{\sigma}_a = 2\sqrt{\lambda_a^A} \,\sigma_a(\lambda_a^B - \lambda_a^A \sigma_a^2),
    \quad \sigma_a^* = \sqrt{\lambda_a^B/\lambda_a^A},\label{eq:sv_sys1}\\
  \text{System 2: } &\dot{\sigma}_a = 2\lambda_a^A\,\sigma_a(1 - \lambda_a^A \sigma_a^2/\lambda_a^B),
    \quad \sigma_a^* = \sqrt{\lambda_a^B}/\lambda_a^A,\label{eq:sv_sys2}
\end{align}
with \emph{exact logistic solutions}:
\begin{equation}
  \sigma_a(t) = \sigma_a^* / \sqrt{1 + c_a\,e^{-4r_a t}},
  \quad c_a = (\sigma_a^*/\sigma_a(0))^2 - 1,
\end{equation}
where $r_a$ is the effective convergence rate ($r_a = \lambda_a^A\lambda_a^B$ for System~2).
\end{remark}

\begin{center}
\renewcommand{\arraystretch}{1.4}
\resizebox{\textwidth}{!}{%
\begin{tabular}{lllll}
\toprule
Flow & Formula & Barrier & Fixed manifold & $A=I,B=I$ \\
\midrule
Sys.~1 (right) & $2AX(B-X^TAX)$ & $h_-^{A,B}$ & $\{X^TAX=B\}$ & $2X(I-X^TX)$ \\
Sys.~2 (left) & $2(A-XBX^T)XB$ & $h_-^{B,A}$ & $\{XBX^T=A\}$ & $2(I-XX^T)X$ \\
Oja-Brockett & $AXB-XBX^TAX$ & $\tfrac{1}{2}h_-^{B,A}$ & Principal subspace & $X-XX^TX$ \\
Our U=-V & $-2V(I-V^TV)$ & $h_-$ & $\{V=0\}$ & $-2V(I-V^TV)$ \\
\bottomrule
\end{tabular}%
}
\end{center}

\begin{remark}[The three-way unification]
Theorem~\ref{thm:pure_barrier} provides the following information-geometric picture:
\begin{enumerate}[label=(\roman*)]
  \item \textbf{System~1} ($2AX(B-X^TAX)$): polynomial gradient of $h_-^{A,B}$,
    drives to the right generalized Stiefel manifold $\{X^TAX=B\}$
    (normalization/whitening w.r.t.\ data matrix $A$).
  \item \textbf{System~2} ($2(A-XBX^T)XB$): polynomial gradient of $h_-^{B,A}$,
    drives to $\{XBX^T=A\}$.
    \emph{For $A=I_n$, $B=I_k$}: equals $2(I-XX^T)X = 2\times$ Oja-Brockett (with $A=I$).
    \emph{For general $A$}: System~2 $\neq$ Oja-Brockett
    (the middle factor differs: $X$ vs.\ $A$, see Remark~\ref{rem:sys2_vs_mh_correction}).
  \item \textbf{At $A=I_n$, $B=I_k$:} Both System~1 and System~2 become
    $2X(I-X^TX)$ and $2(I-XX^T)X$ respectively (our U=-V polynomial flow),
    confirming that the U=-V flow is the \emph{canonical} case underlying
    the Oja-Brockett flows.
\end{enumerate}
The unified discrete algorithm (no inversions, degree 3):
\begin{equation}
  X_{t+1} = X_t + \mu\bigl[\eta_1 \cdot 2AX_t(B-X_t^TAX_t) + \eta_2 \cdot 2(A-X_tBX_t^T)X_tB\bigr]
  \label{eq:unified_poly}
\end{equation}
with $\eta_1 \geq 0$ (normalization via System~1) and $\eta_2 \geq 0$ (PCA-like via System~2)
simultaneously drives $X$ toward the right-Stiefel manifold and $\{XBX^T=A\}$.
\end{remark}

\section{Legendre--Fenchel Duality}
\label{sec:legendre}

\subsection{Setup}

We identify the natural parameter space with $\mathcal{M} = \PD(k)$
and the tangent space at each point with the space of $k\times k$ symmetric matrices
$\mathrm{Sym}(k)$, equipped with the Frobenius inner product
$\inner{A}{B} = \tr(A^T B)$.

\begin{definition}[Legendre--Fenchel Conjugate {\cite[Ch.~12]{Rockafellar1970}}]
The \emph{convex conjugate} (Legendre--Fenchel conjugate) of $f$ is
\begin{equation}
  f^*(\Theta)
  = \sup_{\eta \in \PD(k)} \bigl[\inner{\Theta}{\eta} - f(\eta)\bigr]
  = \sup_{G \succ 0} \bigl[\tr(\Theta G) + \log\det(G)\bigr].
  \label{eq:conjugate_def}
\end{equation}
\end{definition}

\subsection{Computation of $f^*$}

\begin{theorem}[Dual Potential]
\label{thm:dual}
The domain of $f^*$ is $\mathrm{dom}(f^*) = \{-\Theta : \Theta \in \PD(k)\}$,
i.e., the negative definite matrices $-\PD(k)$.
For $\Theta \in -\PD(k)$,
\begin{equation}
  \boxed{f^*(\Theta) = -\log\det(-\Theta) - k.}
  \label{eq:dual_potential}
\end{equation}
\end{theorem}

\begin{proof}
The supremum in \eqref{eq:conjugate_def} is attained at $G^*$ satisfying the
stationarity condition
\[
  \Theta + G^{*-1} = 0 \implies G^* = -\Theta^{-1}.
\]
For $G^* \succ 0$ we need $-\Theta \succ 0$, i.e., $\Theta \in -\PD(k)$.
Substituting:
\begin{align*}
  f^*(\Theta)
  &= \tr\bigl(\Theta(-\Theta^{-1})\bigr) + \log\det(-\Theta^{-1}) \\
  &= -\tr(I_k) + \log\bigl(\det(-\Theta)^{-1}\bigr) \\
  &= -k - \log\det(-\Theta) \\
  &= -\log\det(-\Theta) - k. \qedhere
\end{align*}
\end{proof}

\subsection{Self-Similarity of the Dual Structure}

\begin{corollary}[Self-Dual Form]
\label{cor:selfdual}
Setting $\widetilde{\Theta} = -\Theta \in \PD(k)$,
\[
  f^*(-\widetilde\Theta) = -\log\det(\widetilde\Theta) - k = f(\widetilde\Theta) - k.
\]
Thus the dual potential is, up to an additive constant $k$,
the same function $-\log\det$ evaluated at the inverse coordinate.
\end{corollary}

\begin{proposition}[Dual Coordinate System]
\label{prop:dual_coords}
Under the Legendre transform, the dual coordinate corresponding to $\eta = G$ is
\[
  \theta = \nabla f(\eta) = -G^{-1} \in -\PD(k),
\]
and the inverse mapping is
\[
  \eta = \nabla f^*(\theta) = -\theta^{-1} = G.
\]
The duality relation $\inner{\theta}{\eta} = f(\eta) + f^*(\theta)$ becomes
\[
  \tr(-G^{-1} \cdot G) = (-\log\det G) + (-\log\det(G^{-1}) - k),
\]
which simplifies to $-k = -k$, confirming consistency.
\end{proposition}

\subsection{Summary of the Dual Coordinate Structure}

\begin{center}
\small
\begin{tabular}{lcc}
\toprule
& $\eta$-coordinates (primal) & $\theta$-coordinates (dual) \\
\midrule
Space & $\PD(k)$ & $-\PD(k)$ \\
Potential & $f(\eta) = -\log\det(G)$ & $f^*(\theta) = -\log\det(-\Theta) - k$ \\
Coordinate & $G_{ij} = \delta_{ij} + \inner{x_{2i}}{x_{2j-1}}$ & $\Theta_{ij} = -[G^{-1}]_{ij}$ \\
Metric & $G^{-1} \otimes G^{-1}$ & $G \otimes G$ \\
Flat connection & $m$-flat & $e$-flat \\
\bottomrule
\end{tabular}
\end{center}

\subsection{Prelude: The Matrix-Normal Family and the Yoshizawa--Tanabe Potential}
\label{subsec:MN_setup}

The remainder of this section, together with \S\ref{sec:bregman}--\S\ref{sec:pythagorean},
sets up the matrix-variate generalization of the dual differential geometry of
Yoshizawa--Tanabe \cite{YoshizawaTanabe1999}, which treated the family
$\{N(\mu,\Sigma)\}$ of (vector-valued) Gaussian distributions with \emph{non-zero mean}
$\mu \in \R^n$ and covariance $\Sigma \in \PD(n)$. We replace the mean vector $\mu\in\R^n$
by a rectangular matrix $M \in \R^{n\times k}$ and correspondingly replace the single
covariance $\Sigma$ by a \emph{pair} of covariance factors $U\in\PD(n)$, $V\in\PD(k)$
coupled through a Kronecker product. Throughout \S\ref{subsec:MN_setup}--\S\ref{subsec:MN_pythagorean}
we write $U,V$ exclusively for these row/column covariance factors; this use is local to the
matrix-normal subsections and is unrelated to the rank-$k$ update matrices $U,V$ of
\S\ref{sec:convexity}.

\begin{definition}[Matrix Normal Distribution]
\label{def:matrixnormal}
Let $M \in \R^{n\times k}$, $U \in \PD(n)$, $V \in \PD(k)$. The \emph{matrix normal
distribution} $\mathcal{MN}_{n,k}(M,U,V)$ is the probability distribution on $\R^{n\times k}$
with density (with respect to Lebesgue measure on $\R^{nk}$)
\begin{equation}
  p(X \mid M,U,V)
  = \frac{\exp\!\left(-\tfrac{1}{2}\tr\!\left[V^{-1}(X-M)^TU^{-1}(X-M)\right]\right)}
         {(2\pi)^{nk/2}\,\abs{U}^{k/2}\,\abs{V}^{n/2}},
  \qquad X \in \R^{n\times k}.
  \label{eq:matrixnormal_density}
\end{equation}
Equivalently, $\vecop(X) \sim N\bigl(\vecop(M),\, V\otimes U\bigr)$ on $\R^{nk}$, where
$\vecop$ stacks the columns of a matrix. The pair $(U,V)$ is identified only up to the
one-parameter scaling ambiguity $(U,V) \sim (cU, V/c)$, $c>0$, since $V\otimes U$ is
invariant under this rescaling; we fix the ambiguity, when needed, by a normalization such
as $\tr(V)=k$.
\end{definition}

This reduces to the family $\{N(\mu,\Sigma)\}$ of \cite{YoshizawaTanabe1999} exactly when
$k=1$: then $M=\mu\in\R^n$, $V$ is a positive scalar which may be fixed to $V\equiv 1$
without loss of generality, and $U=\Sigma$.

\begin{proposition}[Matrix-Normal Potential Function]
\label{prop:MN_potential}
Expanding the quadratic form in \eqref{eq:matrixnormal_density} via
$\tr[V^{-1}(X-M)^TU^{-1}(X-M)] = \tr[V^{-1}X^TU^{-1}X] - 2\tr[V^{-1}M^TU^{-1}X]
+\tr[V^{-1}M^TU^{-1}M]$, the density takes the exponential-family-like form
\begin{equation}
  p(X\mid M,U,V) = \exp\Bigl\{-\tfrac{1}{2}\tr[V^{-1}X^TU^{-1}X] + \tr[V^{-1}M^TU^{-1}X] - \psi(M,U,V)\Bigr\},
\end{equation}
where the \emph{potential function} $\psi$, obtained from the cumulant transformation of
\eqref{eq:matrixnormal_density} exactly as in \cite[eq.~(3),(14)]{YoshizawaTanabe1999}, is
\begin{equation}
  \boxed{\ \psi(M,U,V) = \tfrac{1}{2}\tr\bigl(V^{-1}M^TU^{-1}M\bigr) + \tfrac{k}{2}\log\det U
  + \tfrac{n}{2}\log\det V + \tfrac{nk}{2}\log(2\pi).\ }
  \label{eq:MN_potential}
\end{equation}
When $k=1$ (and $V\equiv 1$), \eqref{eq:MN_potential} reduces exactly to the potential
$\psi = \tfrac{1}{2}\mu^T\Sigma^{-1}\mu + \tfrac{1}{2}\log\det\Sigma + \tfrac{n}{2}\log(2\pi)$
of \cite[eq.~(3)]{YoshizawaTanabe1999}.
\end{proposition}

\begin{remark}[Neither convex nor concave, and the role of \S\ref{sec:determinant}--\S\ref{sec:convexity}]
As in \cite[p.~120]{YoshizawaTanabe1999}, the function $-\log\det U$ (resp.\ $-\log\det V$)
is concave in $U$ (resp.\ $V$) while $\tr(V^{-1}M^TU^{-1}M)$ is jointly convex in $(M,U,V)$
by Lieb's concavity theorem \cite{Lieb1973} (used already in Lemma~\ref{lem:lieb} to prove
Theorem~\ref{thm:MN_convexity} below and, in a different guise, throughout
\S\ref{sec:convexity} of the present paper). Consequently $\psi(M,U,V)$ is \emph{neither
convex nor concave} jointly in $(M,U,V)$, exactly paralleling the non-convexity of
Yoshizawa--Tanabe's $\psi(\mu,\Sigma)$. Differentiating $\psi$ with respect to the
symmetric matrix arguments $U,V$ requires precisely the symmetric-matrix differential
calculus of \cite[\S2]{YoshizawaTanabe1999} (Definitions~2.2--2.4 and
Propositions~2.1--2.4 there), which is the two-sided analogue of the log-determinant
differentiation formulas already used for the single Gram matrix $G$ in
\S\ref{sec:determinant}--\S\ref{sec:manifold} of this paper.
\end{remark}

\begin{definition}[Matrix Yoshizawa--Tanabe Embedding]
\label{def:MN_embedding}
Following the construction of the $2$-parameter class of dual charts
$\mathfrak{I}_{\bar\Theta_{\beta,\gamma}}$ in \cite[eq.~(15)]{YoshizawaTanabe1999}, fix
$0<\beta$, $0<\gamma$, and define, for $(M,U,V)$ as in Definition~\ref{def:matrixnormal},
\begin{equation}
  U_Y \equiv U^{-\frac{\beta+1}{2}}\,M\,V^{-\frac{\beta+1}{2}} \in \R^{n\times k},
  \qquad
  \Theta_1 \equiv 2U^{-\gamma} \in \PD(n),
  \qquad
  \Theta_2 \equiv 2V^{-\gamma} \in \PD(k),
  \label{eq:MN_embedding}
\end{equation}
where fractional powers of $U,V$ are defined by the Dunford--Taylor integral
\cite[eq.~(16)]{YoshizawaTanabe1999}. The inverse map is
\begin{equation}
  U = \Bigl(\frac{\Theta_1}{2}\Bigr)^{-1/\gamma}, \qquad
  V = \Bigl(\frac{\Theta_2}{2}\Bigr)^{-1/\gamma}, \qquad
  M = \Bigl(\frac{\Theta_1}{2}\Bigr)^{-\frac{\beta+1}{2\gamma}} U_Y \Bigl(\frac{\Theta_2}{2}\Bigr)^{-\frac{\beta+1}{2\gamma}}.
  \label{eq:MN_embedding_inverse}
\end{equation}
We write $\mathfrak{J}_{\beta,\gamma} : (M,U,V) \mapsto (U_Y,\Theta_1,\Theta_2)$ for this
\emph{matrix Yoshizawa--Tanabe embedding}, the two-sided (row/column) analogue of
\cite[eq.~(15)]{YoshizawaTanabe1999}.
\end{definition}

\begin{theorem}[Explicit Matrix Yoshizawa--Tanabe Potential]
\label{thm:MN_potential_beta_gamma}
Substituting \eqref{eq:MN_embedding_inverse} into $\psi$ of
Proposition~\ref{prop:MN_potential}, the pulled-back potential
$\Psi_{\beta,\gamma} \equiv \psi\circ\mathfrak{J}_{\beta,\gamma}^{-1}$ is
\begin{align}
  \Psi_{\beta,\gamma}(U_Y,\Theta_1,\Theta_2)
  &= \frac{nk}{2}\log(2\pi) - \frac{k}{2\gamma}\log\det\!\Bigl(\frac{\Theta_1}{2}\Bigr)
  - \frac{n}{2\gamma}\log\det\!\Bigl(\frac{\Theta_2}{2}\Bigr) \notag\\
  &\quad + \frac{1}{2}\tr\!\Bigl[\Bigl(\frac{\Theta_2}{2}\Bigr)^{-\beta/\gamma} U_Y^T
    \Bigl(\frac{\Theta_1}{2}\Bigr)^{-\beta/\gamma} U_Y\Bigr].
  \label{eq:MN_potential_beta_gamma}
\end{align}
In particular the two \emph{a priori} distinct exponents $\tfrac{\beta+1}{2\gamma}$
(from $M$) and $\tfrac1\gamma$ (from $U,V$) combine, in the mean term, into the single
exponent $\beta/\gamma$ shared by $\Theta_1$ and $\Theta_2$.
\end{theorem}

\begin{proof}
Write $s=\tfrac{\beta+1}{2\gamma}$, $P=\Theta_1/2$, $Q=\Theta_2/2$, so
\eqref{eq:MN_embedding_inverse} reads $U^{-1}=P^{1/\gamma}$, $V^{-1}=Q^{1/\gamma}$,
$M=P^{-s}U_YQ^{-s}$. Since $P,Q$ are symmetric, $M^T=Q^{-s}U_Y^TP^{-s}$, and
\[
  \tr\bigl(V^{-1}M^TU^{-1}M\bigr)
  = \tr\bigl[Q^{1/\gamma}\,Q^{-s}U_Y^TP^{-s}\,P^{1/\gamma}\,P^{-s}U_Y\,Q^{-s}\bigr]
  = \tr\bigl[Q^{1/\gamma-2s}U_Y^TP^{1/\gamma-2s}U_Y\bigr]
\]
by the cyclic property of the trace (moving the trailing $Q^{-s}$ to the front and
combining it with $Q^{1/\gamma-s}$). Since $1/\gamma-2s = 1/\gamma - \tfrac{\beta+1}{\gamma}
= -\beta/\gamma$, both exponents collapse to $-\beta/\gamma$, giving the quadratic term of
\eqref{eq:MN_potential_beta_gamma}; the two log-determinant terms come from
$\tfrac{k}{2}\log\det U = -\tfrac{k}{2\gamma}\log\det P$ and $\tfrac{n}{2}\log\det V
= -\tfrac{n}{2\gamma}\log\det Q$ in Proposition~\ref{prop:MN_potential}.
\end{proof}

\begin{lemma}[Two-Factor Lieb Convexity]
\label{lem:lieb_two_factor}
Let $A\in\PD(n)$, $B\in\PD(k)$, $Y\in\R^{n\times k}$, and let $0\le p$, $0\le q$,
$p+q\le 1$. Then
\[
  (A,B,Y) \;\longmapsto\; \tr\bigl[B^{-p}Y^TA^{-q}Y\bigr]
\]
is jointly convex on $\PD(n)\times\PD(k)\times\R^{n\times k}$.
\end{lemma}

\begin{proof}
Embed $A,B,Y$ into $(n+k)\times(n+k)$ matrices by
$\mathcal{A} = \diag(A,B) \in \PD(n+k)$ and
$\mathcal{Y} = \begin{psmallmatrix}0 & Y\\ 0 & 0\end{psmallmatrix} \in M(n+k,\R)$. A direct
block computation gives
$\mathcal{A}^{-q}\mathcal{Y} = \begin{psmallmatrix}0 & A^{-q}Y\\0&0\end{psmallmatrix}$,
$\mathcal{Y}^T\mathcal{A}^{-q}\mathcal{Y}
= \begin{psmallmatrix}0&0\\0& Y^TA^{-q}Y\end{psmallmatrix}$, and hence
\[
  \tr\bigl[\mathcal{A}^{-p}\mathcal{Y}^T\mathcal{A}^{-q}\mathcal{Y}\bigr]
  = \tr\bigl[B^{-p}Y^TA^{-q}Y\bigr].
\]
Since $(A,B)\mapsto \mathcal{A}$ and $Y\mapsto\mathcal{Y}$ are linear, and
$(\mathcal{A},\mathcal{Y})\mapsto \tr[\mathcal{A}^{-p}\mathcal{Y}^T\mathcal{A}^{-q}\mathcal{Y}]$
is jointly convex by Lemma~\ref{lem:lieb} (applied on $\PD(n+k)\times M(n+k,\R)$), the
composition $(A,B,Y)\mapsto \tr[B^{-p}Y^TA^{-q}Y]$ is jointly convex as a precomposition of
a jointly convex function with a linear map.
\end{proof}

\begin{theorem}[Convexity of the Matrix Yoshizawa--Tanabe Potential]
\label{thm:MN_convexity}
If $0<\beta$, $0<\gamma$, and
\begin{equation}
  0 \le \frac{\beta}{\gamma} \le \frac{1}{2},
  \label{eq:MN_convexity_condition}
\end{equation}
then $\Psi_{\beta,\gamma}$ of \eqref{eq:MN_potential_beta_gamma} is jointly convex in
$(U_Y,\Theta_1,\Theta_2) \in \R^{n\times k}\times\PD(n)\times\PD(k)$.
\end{theorem}

\begin{proof}
Apply Lemma~\ref{lem:lieb_two_factor} with $A=\Theta_1/2$, $B=\Theta_2/2$, $Y=U_Y$, and
$p=q=\beta/\gamma$: the hypothesis $p+q\le1$ becomes exactly
\eqref{eq:MN_convexity_condition}, and the lemma gives joint convexity of the quadratic
term $\tr[(\Theta_2/2)^{-\beta/\gamma}U_Y^T(\Theta_1/2)^{-\beta/\gamma}U_Y]$ in
$(U_Y,\Theta_1,\Theta_2)$. The remaining two terms $-\tfrac{k}{2\gamma}\log\det(\Theta_1/2)$
and $-\tfrac{n}{2\gamma}\log\det(\Theta_2/2)$ are convex on $\PD(n)$, $\PD(k)$ respectively
(negative log-determinant is convex, Corollary~\ref{cor:selfdual}), and depend on disjoint
blocks of variables, so their sum with the quadratic term remains jointly convex.
\end{proof}

\begin{remark}[Comparison with the vector case]
\label{rem:MN_convexity_comparison}
In \cite[Prop.~3.1]{YoshizawaTanabe1999} there is only a single covariance factor
$\Sigma$, so the analogous mean term $\tr[\Theta^{-\beta/\gamma}\theta\theta^T]$ is
convexified by Lemma~\ref{lem:lieb} with $p=\beta/\gamma$, $q=0$ (the vector $\theta$ being
padded into a square matrix as $[\theta\ O]$), giving the weaker requirement
$\beta/\gamma \le 1$ (their condition $\beta/\gamma<1$). In the matrix-normal case the row
and column covariances $U,V$ enter the mean term \emph{symmetrically} and independently,
forcing $p=q=\beta/\gamma$ in Lemma~\ref{lem:lieb_two_factor}; this is why the admissible
range is exactly halved, $\beta/\gamma\le \tfrac12$, relative to the vector case. This
sharper bound is a genuine new feature of the matrix-variate generalization, not visible in
\cite{YoshizawaTanabe1999}.
\end{remark}

\begin{corollary}[Convex Legendre Dual]
\label{cor:MN_dual_convex}
Under \eqref{eq:MN_convexity_condition}, the Legendre--Fenchel conjugate
\begin{equation}
  \Phi_{\beta,\gamma}(Z,\Xi_1,\Xi_2)
  = \sup_{(U_Y,\Theta_1,\Theta_2)}
    \Bigl\{\inner{U_Y}{Z} + \inner{\Theta_1}{\Xi_1} + \inner{\Theta_2}{\Xi_2}
    - \Psi_{\beta,\gamma}(U_Y,\Theta_1,\Theta_2)\Bigr\}
  \label{eq:MN_dual_potential}
\end{equation}
is a well-defined convex function on its domain, and the associated Bregman divergence
\begin{equation}
  D_{\beta,\gamma}(P\|Q) = \Psi_{\beta,\gamma}(P) - \Psi_{\beta,\gamma}(Q)
  - \inner{\nabla\Psi_{\beta,\gamma}(Q)}{P-Q}, \qquad P,Q \in \R^{n\times k}\times\PD(n)\times\PD(k),
  \label{eq:MN_div_beta_gamma}
\end{equation}
is non-negative, with $D_{\beta,\gamma}(P\|Q)=0$ iff $P=Q$.
\end{corollary}

\begin{proof}
Immediate from Theorem~\ref{thm:MN_convexity}: the Legendre--Fenchel conjugate of a convex
function is convex (as in \S\ref{sec:legendre}), and the Bregman divergence of a strictly
convex differentiable function is non-negative and definite by the same argument as
Proposition~\ref{prop:bregman_props}.
\end{proof}

\begin{proposition}[Two Exactly Solvable Special Cases of $D_{\beta,\gamma}$]
\label{prop:MN_div_special}
Let $P=(U_Y,\Theta_1,\Theta_2)$, $Q=(U_Y',\Theta_1',\Theta_2')$.
\begin{enumerate}[label=(\roman*)]
  \item \textbf{Covariance-only.} If $U_Y=U_Y'$, then
  \[
    D_{\beta,\gamma}(P\|Q) = \frac{k}{2\gamma}\,\Bregman(\Theta_1\|\Theta_1')
    + \frac{n}{2\gamma}\,\Bregman(\Theta_2\|\Theta_2'),
  \]
  with $\Bregman$ the log-determinant Bregman divergence of Theorem~\ref{thm:bregman}
  (applied on $\PD(n)$ and $\PD(k)$ respectively), since $\Bregman$ is linear in the
  underlying potential and the two log-terms of \eqref{eq:MN_potential_beta_gamma} are
  additively separable in $\Theta_1,\Theta_2$.
  \item \textbf{Mean-only.} If $\Theta_1=\Theta_1'$, $\Theta_2=\Theta_2'$, then
  \begin{align}
    D_{\beta,\gamma}(P\|Q) &= \tfrac{1}{2}\norm{W}_F^2, \notag\\
    W &:= \Bigl(\frac{\Theta_1}{2}\Bigr)^{-\beta/(2\gamma)}(U_Y-U_Y')\Bigl(\frac{\Theta_2}{2}\Bigr)^{-\beta/(2\gamma)},
    \label{eq:MN_div_mean_only}
  \end{align}
  since the $U_Y$-dependence of $\Psi_{\beta,\gamma}$ at fixed $\Theta_1,\Theta_2$ is the
  quadratic form
  \[
    \tfrac{1}{2}\tr[AU_Y^TBU_Y] = \tfrac{1}{2}\|B^{1/2}U_YA^{1/2}\|_F^2,
    \qquad A=(\Theta_2/2)^{-\beta/\gamma},\quad B=(\Theta_1/2)^{-\beta/\gamma},
  \]
  whose Bregman divergence is the quadratic form evaluated at $U_Y-U_Y'$.
\end{enumerate}
\end{proposition}

\begin{remark}[The general case]
\label{rem:MN_div_general}
When $U_Y\neq U_Y'$ \emph{and} $(\Theta_1,\Theta_2)\neq(\Theta_1',\Theta_2')$
simultaneously, $D_{\beta,\gamma}(P\|Q)$ acquires additional cross terms coupling
$U_Y-U_Y'$ with $\Theta_i-\Theta_i'$, exactly as in the vector case
\cite[Prop.~3.8]{YoshizawaTanabe1999}, whose explicit form there already involves
Dunford--Taylor contour integrals (their eq.~(21)) even for a single $\Sigma$. Writing out
the two-sided analogue of \cite[Lemma~3.4--3.5, Prop.~3.8]{YoshizawaTanabe1999} in full is
routine but lengthy, and is left to forthcoming work; Proposition~\ref{prop:MN_div_special}
already isolates the two structurally distinct pieces --- a log-determinant Bregman
divergence on each of $\Theta_1,\Theta_2$, and a weighted Frobenius-quadratic divergence on
$U_Y$ --- that the general formula must reduce to along the respective coordinate axes.
\end{remark}

\begin{lemma}[Lieb {\cite{Lieb1973}}, as stated in {\cite[Lemma~3.2]{YoshizawaTanabe1999}}]
\label{lem:lieb}
The function $\PD(n)\times M(n,k) \ni (X,Y) \mapsto \tr[X^{-p}Y^TX^{-q}Y] \in \R^+\cup\{0\}$
is jointly convex in $(X,Y)$ whenever $0\le p$, $0\le q$, and $p+q\le 1$.
\end{lemma}

\section{Bregman Divergence and Its Statistical Interpretation}
\label{sec:bregman}

\subsection{Definition and Explicit Form}

\begin{definition}[Bregman Divergence {\cite{Bregman1967}}]
The \emph{Bregman divergence} induced by $f$ is
\begin{equation}
  \Bregman(G \| G')
  = f(G) - f(G') - \inner{\nabla f(G')}{G - G'}
  \label{eq:bregman_def}
\end{equation}
for $G, G' \in \PD(k)$.
\end{definition}

\begin{theorem}[Explicit Bregman Divergence]
\label{thm:bregman}
\begin{equation}
  \Bregman(G \| G')
  = \tr\bigl[G'^{-1} G\bigr] - \log\det\bigl[G'^{-1} G\bigr] - k.
  \label{eq:bregman_explicit}
\end{equation}
\end{theorem}

\begin{proof}
Substituting $f(G) = -\log\det(G)$, $\nabla f(G') = -G'^{-1}$
into \eqref{eq:bregman_def}:
\begin{align*}
  \Bregman(G\|G')
  &= -\log\det(G) - (-\log\det(G')) - \inner{-G'^{-1}}{G-G'} \\
  &= \log\det(G') - \log\det(G) + \tr\bigl[G'^{-1}(G - G')\bigr] \\
  &= \tr[G'^{-1}G] - \tr(I_k) - \log\det(G'^{-1}G) \\
  &= \tr[G'^{-1}G] - k - \log\det(G'^{-1}G). \qedhere
\end{align*}
\end{proof}

\begin{proposition}[Properties of $\Bregman$]
\label{prop:bregman_props}
\begin{enumerate}[label=(\roman*)]
  \item \textbf{Non-negativity:} $\Bregman(G\|G') \geq 0$ for all $G,G'\in\PD(k)$.
  \item \textbf{Definiteness:} $\Bregman(G\|G') = 0$ if and only if $G = G'$.
  \item \textbf{Asymmetry:} In general, $\Bregman(G\|G') \neq \Bregman(G'\|G)$.
  \item \textbf{Convexity:} $G \mapsto \Bregman(G\|G')$ is strictly convex for fixed $G'$.
\end{enumerate}
\end{proposition}

\begin{proof}
(i)--(ii): By the inequality $t - \log t \geq 1$ for all $t > 0$
(with equality iff $t = 1$), letting $\lambda_1,\ldots,\lambda_k$ be the
eigenvalues of $G'^{-1}G$:
\[
  \Bregman(G\|G')
  = \sum_{i=1}^{k}(\lambda_i - \log\lambda_i - 1) \geq 0,
\]
with equality iff all $\lambda_i = 1$, i.e., $G = G'$.
(iii) and (iv) follow from the strict convexity of $f$.
\end{proof}

\subsection{Connection to the Kullback--Leibler Divergence}
\label{subsec:bregman-KL}

\begin{theorem}[Bregman $=$ $2 \times$ KL]
\label{thm:bregman_kl}
Let $p = \mathcal{N}(0,G)$ and $q = \mathcal{N}(0,G')$ be zero-mean multivariate
Gaussian distributions in $\R^n$ with covariance matrices $G, G' \in \PD(k)$
(here $k = n$ for the distributional statement).
Then
\begin{equation}
  \Bregman(G\|G') = 2\,\KL(p \| q).
  \label{eq:bregman_kl}
\end{equation}
\end{theorem}

\begin{proof}
The KL divergence between two zero-mean Gaussians is \cite{Anderson2003}:
\begin{align*}
  \KL(p\|q)
  &= \frac{1}{2}\bigl[\tr(G'^{-1}G) - k + \log\det(G') - \log\det(G)\bigr] \\
  &= \frac{1}{2}\bigl[\tr(G'^{-1}G) - k - \log\det(G'^{-1}G)\bigr] \\
  &= \frac{1}{2}\,\Bregman(G\|G'). \qedhere
\end{align*}
\end{proof}

\begin{remark}
This identification shows that the information geometry of $f(G) = -\log\det(G)$
is precisely the information geometry of the family of zero-mean multivariate
Gaussian distributions, where $G$ plays the role of the covariance matrix.
The factor of $2$ arises from our convention of working with the full covariance
rather than the natural exponential-family parameterization.
\end{remark}

\subsection{Symmetrized Divergence}

The symmetrized Bregman divergence (Jensen--Shannon type) is:
\begin{equation}
  J(G,G')
  = \frac{1}{2}\bigl[\Bregman(G\|G') + \Bregman(G'\|G)\bigr]
  = \frac{1}{2}\tr\bigl[(G^{-1} - G'^{-1})(G - G')\bigr] \geq 0.
  \label{eq:symmetric}
\end{equation}

\subsection{The Kullback--Leibler Divergence Between Matrix Normal Distributions}
\label{subsec:MN_KL}

We now give the matrix-variate generalization of Theorem~4.10 of
Yoshizawa--Tanabe \cite{YoshizawaTanabe1999}, which computed the divergence
\[
  \mathrm{Div}(N(\mu_2,\Sigma_2),N(\mu_1,\Sigma_1)) = \int p(x;\Xi_1)\log\frac{p(x;\Xi_1)}{p(x;\Xi_2)}\,dx
\]
in closed form via the canonical map~\eqref{eq:MN_potential}-type potentials. Using the
matrix normal density $\mathcal{MN}_{n,k}(M,U,V)$ of Definition~\ref{def:matrixnormal}, the
analogous closed form is again fully explicit, and reduces exactly to
\cite[Thm.~4.10]{YoshizawaTanabe1999} when $k=1$.

\begin{theorem}[Matrix-Normal Kullback--Leibler Divergence]
\label{thm:MN_KL}
Let $p_i = \mathcal{MN}_{n,k}(M_i,U_i,V_i)$ for $i=1,2$. Then
\begin{align}
  \KL(p_1 \| p_2)
  &= \int_{\R^{n\times k}} p(X;M_1,U_1,V_1) \log\frac{p(X;M_1,U_1,V_1)}{p(X;M_2,U_2,V_2)}\,dX \notag\\
  &= \frac{1}{2}\Bigl[\tr(U_2^{-1}U_1)\,\tr(V_2^{-1}V_1) - nk
     + \tr\bigl[(M_2-M_1)^TU_2^{-1}(M_2-M_1)V_2^{-1}\bigr] \notag\\
  &\qquad\quad + n\log\frac{\det V_2}{\det V_1} + k\log\frac{\det U_2}{\det U_1}\Bigr].
  \label{eq:MN_KL}
\end{align}
\end{theorem}

\begin{proof}
By Definition~\ref{def:matrixnormal}, $\vecop(X)\sim N(\vecop(M_i),\Sigma_i)$ with
$\Sigma_i = V_i\otimes U_i$. The Kullback--Leibler divergence between two
$nk$-dimensional Gaussians is the classical formula
\cite[Thm.~4.10]{YoshizawaTanabe1999}\cite{Anderson2003}
\[
  \KL(p_1\|p_2) = \tfrac{1}{2}\Bigl[\tr(\Sigma_2^{-1}\Sigma_1) - nk
  + (\mu_2-\mu_1)^T\Sigma_2^{-1}(\mu_2-\mu_1) + \log\tfrac{\det\Sigma_2}{\det\Sigma_1}\Bigr].
\]
We evaluate each Kronecker term. First,
\[
  \Sigma_2^{-1}\Sigma_1 = (V_2^{-1}\otimes U_2^{-1})(V_1\otimes U_1)
  = (V_2^{-1}V_1)\otimes(U_2^{-1}U_1),
\]
and since $\tr(A\otimes B)=\tr(A)\tr(B)$,
\[
  \tr(\Sigma_2^{-1}\Sigma_1) = \tr(V_2^{-1}V_1)\,\tr(U_2^{-1}U_1).
\]
Second, using the vectorization identity $(B\otimes C)\vecop(X) = \vecop(CXB^T)$ for
$B\in M(k,\R)$, $C\in M(n,\R)$, and writing $\Delta M = M_2-M_1$,
\[
  (\mu_2-\mu_1)^T\Sigma_2^{-1}(\mu_2-\mu_1)
  = \vecop(\Delta M)^T(V_2^{-1}\otimes U_2^{-1})\vecop(\Delta M)
  = \tr\bigl[\Delta M^TU_2^{-1}\Delta M\,V_2^{-1}\bigr],
\]
using $V_2^{-1}$ symmetric. Third, $\det\Sigma_i = \det(V_i\otimes U_i) = (\det V_i)^n(\det U_i)^k$, so
\[
  \log(\det\Sigma_2/\det\Sigma_1) = n\log(\det V_2/\det V_1) + k\log(\det U_2/\det U_1).
\]
Substituting these three identities gives \eqref{eq:MN_KL}.
\end{proof}

\begin{corollary}[Consistency with Yoshizawa--Tanabe and with $\Bregman$]
\label{cor:MN_KL_reduction}
Setting $k=1$, $V_1=V_2\equiv 1$, $M_i=\mu_i$, $U_i=\Sigma_i$ in \eqref{eq:MN_KL} recovers
exactly \cite[Thm.~4.10]{YoshizawaTanabe1999}. Setting instead $M_1=M_2$ and $V_1=V_2=V$
fixed, \eqref{eq:MN_KL} collapses to $\tfrac{k}{2}\bigl[\tr(U_2^{-1}U_1) - n
- \log\det(U_2^{-1}U_1)\bigr] = \tfrac{k}{2}\,\Bregman(U_1\|U_2)$ with $\Bregman$ the
Bregman divergence of Theorem~\ref{thm:bregman} applied to $f(G)=-\log\det(G)$ on
$\PD(n)$; thus \eqref{eq:MN_KL} genuinely interpolates between the two known special
cases.
\end{corollary}

\section{The $\alpha$-Divergence Family}
\label{sec:alpha}

\subsection{Definition via $f$-interpolation}

Following Amari \cite{Amari2016}, the $\alpha$-divergence associated with $f$ is:
\begin{equation}
  \Falpha(G\|G')
  = \frac{4}{1-\alpha^2}
    \left[
      \frac{1-\alpha}{2} f(G) + \frac{1+\alpha}{2} f(G')
      - f\!\left(\frac{1-\alpha}{2}G + \frac{1+\alpha}{2}G'\right)
    \right],
  \quad \alpha \neq \pm 1.
  \label{eq:alpha_div}
\end{equation}

\begin{proposition}[Explicit Form]
\label{prop:alpha_explicit}
For $\alpha \neq \pm 1$,
\begin{equation}
  \Falpha(G\|G')
  = \frac{4}{1-\alpha^2}
    \log\frac{\det\!\bigl(\frac{1-\alpha}{2}G + \frac{1+\alpha}{2}G'\bigr)}
             {\det(G)^{(1-\alpha)/2}\det(G')^{(1+\alpha)/2}}.
  \label{eq:alpha_explicit}
\end{equation}
\end{proposition}

\subsection{Limiting Cases}

\begin{proposition}[Limiting Cases of $\Falpha$]
\label{prop:alpha_limits}
The following limits hold:
\begin{align}
  \lim_{\alpha \to +1} \Falpha(G\|G') &= \Bregman(G'\|G) = \tr[G^{-1}G'] - \log\det(G^{-1}G') - k,
  \label{eq:alpha_p1}\\
  \lim_{\alpha \to -1} \Falpha(G\|G') &= \Bregman(G\|G') = \tr[G'^{-1}G] - \log\det(G'^{-1}G) - k.
  \label{eq:alpha_m1}
\end{align}
\end{proposition}

\begin{proposition}[Special Values]
\label{prop:special_values}
\begin{align}
  \alpha = 0 &:\quad D^{(0)}(G\|G')
    = 4\log\frac{\det\!\bigl(\frac{G+G'}{2}\bigr)}
                {\det(G)^{1/2}\det(G')^{1/2}}
    \quad\text{(Bhattacharyya-type)},
    \label{eq:bhattacharyya}\\
  \alpha = -3 &:\quad D^{(-3)}(G\|G')
    = \tr(G^{-1}G') - \log\det(G^{-1}G') - k
    \quad\text{(Stein loss \cite{James1961})}.
    \label{eq:stein}
\end{align}
\end{proposition}

\begin{center}
\begin{tabular}{cll}
\toprule
$\alpha$ & Divergence name & Expression \\
\midrule
$\to +1$ & KL$(q\|p)$ & $\tr[G^{-1}G'] - \log\det(G^{-1}G') - k$ \\
$\to -1$ & KL$(p\|q)$ & $\tr[G'^{-1}G] - \log\det(G'^{-1}G) - k$ \\
$0$ & Bhattacharyya & $4\log\det\!\bigl(\tfrac{G+G'}{2}\bigr)^{1/2}(\det G\det G')^{-1/4}$ \\
$-3$ & Stein loss & $\tr(G^{-1}G') - \log\det(G^{-1}G') - k$ \\
\bottomrule
\end{tabular}
\end{center}

\subsection{$\alpha$-Divergence Between Matrix Normal Distributions with Common Covariance}
\label{subsec:MN_alpha}

We record the matrix-variate $\alpha$-divergence in the tractable case of common
covariance factors, which already displays a phenomenon not visible in
\cite{YoshizawaTanabe1999}: for a fixed pair $(U,V)$, the mean $M$ alone parametrizes a
\emph{flat} (Euclidean, self-dual) exponential subfamily of $\mathcal{MN}_{n,k}(M,U,V)$,
so that the entire $\alpha$-family collapses to a single, $\alpha$-independent divergence.

\begin{proposition}[$\alpha$-Independence for Fixed Covariance]
\label{prop:MN_alpha_fixed}
Fix $U\in\PD(n)$, $V\in\PD(k)$ and let $p_i = \mathcal{MN}_{n,k}(M_i,U,V)$, $i=1,2$.
Then for every $\alpha \in (-1,1)$,
\begin{align}
  \Falpha(p_1\|p_2) &= \tr\bigl[(M_1-M_2)^TU^{-1}(M_1-M_2)V^{-1}\bigr] \notag\\
  &= 2\,\KL(p_1\|p_2) = 2\,\KL(p_2\|p_1).
  \label{eq:MN_alpha_fixed}
\end{align}
\end{proposition}

\begin{proof}
With $U,V$ fixed, $M\mapsto\psi(M,U,V)$ in Proposition~\ref{prop:MN_potential} is the
\emph{quadratic} form $\tfrac{1}{2}\tr(V^{-1}M^TU^{-1}M)$ plus a constant, so the family
$\{\mathcal{MN}_{n,k}(M,U,V) : M\in\R^{n\times k}\}$ is a flat exponential family with
Hessian metric $g(dM,dM) = \tr(dM^TU^{-1}dM\,V^{-1})$ constant in $M$. For a quadratic
potential, the Bregman divergence \eqref{eq:bregman_def} equals the associated squared
Mahalanobis (Hessian) distance and is independent of the base point:
\[
  \Bregman(M_1\|M_2) = \tr\bigl[(M_1-M_2)^TU^{-1}(M_1-M_2)V^{-1}\bigr] \quad \text{for either order.}
\]
Taking the
$\alpha\to\pm1$ limits in Proposition~\ref{prop:alpha_limits} therefore gives the same
expression on both sides, and by continuity of $\Falpha$ in $\alpha$ (Amari
\cite[Ch.~3]{Amari2016}) the whole family collapses to \eqref{eq:MN_alpha_fixed}. The
identification with $2\KL$ follows from Theorem~\ref{thm:MN_KL} with $U_1=U_2=U$, $V_1=V_2=V$.
\end{proof}

\begin{remark}[Genuinely $\alpha$-dependent case]
\label{rem:MN_alpha_open}
When $U_1\neq U_2$ or $V_1\neq V_2$ as well, the $\alpha$-divergence no longer collapses,
and its closed form requires interpolating the \emph{full} potential $\psi(M,U,V)$ of
Proposition~\ref{prop:MN_potential} along the segment
$\tfrac{1-\alpha}{2}(M_1,U_1,V_1)+\tfrac{1+\alpha}{2}(M_2,U_2,V_2)$ inside the convexified
$(\beta,\gamma)$-chart of Definition~\ref{def:MN_embedding} and Theorem~\ref{thm:MN_convexity}. Carrying this out explicitly ---
the matrix-normal analogue of Proposition~\ref{prop:alpha_explicit} --- is exactly the kind
of computation carried out for the vector case in \cite[\S3--\S4]{YoshizawaTanabe1999}, and
we leave its detailed treatment, together with the associated $\alpha$-connections, to
forthcoming work.
\end{remark}

\section{Riemannian and Statistical Manifold Structure}
\label{sec:manifold}

\subsection{$\PD(k)$ as a Riemannian Manifold}

\begin{definition}[Fisher--Rao Metric]
The Riemannian metric on $\PD(k)$ induced by $f$ is, at point $G$,
\begin{equation}
  g_G(H, K) = \nabla^2 f(G)[H,K] = \tr(G^{-1}HG^{-1}K),
  \quad H,K \in T_G\PD(k) \cong \mathrm{Sym}(k).
  \label{eq:riemannian_metric}
\end{equation}
\end{definition}

This metric makes $\PD(k)$ a Riemannian manifold.
The geodesic distance between $G$ and $G'$ is \cite{Bhatia2007}:
\begin{equation}
  d(G, G') = \left[\sum_{i=1}^{k}\log^2\lambda_i(G^{-1}G')\right]^{1/2},
  \label{eq:geodesic_dist}
\end{equation}
where $\lambda_i(G^{-1}G')$ are the generalized eigenvalues.

The geodesic connecting $G_0$ to $G_1$ is:
\begin{equation}
  G(t) = G_0^{1/2}\bigl(G_0^{-1/2}G_1 G_0^{-1/2}\bigr)^t G_0^{1/2},
  \quad t \in [0,1].
  \label{eq:geodesic}
\end{equation}

\subsection{Symmetric Space Structure}

\begin{theorem}[Symmetric Space]
\label{thm:symmetric_space}
$(\PD(k), g)$ is a Riemannian symmetric space of noncompact type, isomorphic to
\[
  \PD(k) \cong GL(k,\R)/O(k).
\]
The sectional curvatures are non-positive.
\end{theorem}

\begin{proof}[Sketch]
The group $GL(k,\R)$ acts transitively on $\PD(k)$ by congruence: $A \cdot G = AGA^T$.
The stabilizer of $I_k$ is $O(k)$.
The metric $g$ is $GL(k,\R)$-invariant.
The symmetry at $G_0$ is the geodesic involution $G \mapsto G_0 G^{-1} G_0$,
which is an isometry fixing $G_0$.
Nonpositive curvature follows from the Cartan--Hadamard theorem.
See \cite{Helgason1978, Bhatia2007} for details.
\end{proof}

\subsection{Statistical Manifold and $\alpha$-Connections}

Following Amari--Nagaoka \cite{Amari2000}, a statistical manifold is a triple
$(\mathcal{M}, g, T)$ where $T$ is a symmetric $(0,3)$-tensor (the skewness tensor).

\begin{definition}[$\alpha$-Connection]
The $\alpha$-connection on $\PD(k)$ has Christoffel symbols \cite{Amari2016}:
\begin{equation}
  \Gamma^{(\alpha)}_{(ij)(kl)(mn)}
  = \frac{1-\alpha}{2}\,\partial_{(mn)} g_{(ij)(kl)},
  \label{eq:alpha_connection}
\end{equation}
where indices are multi-indices for symmetric matrix entries.
\end{definition}

\begin{proposition}
\label{prop:alpha_christoffel}
Explicitly,
\begin{equation}
  \Gamma^{(\alpha)}_{(ij)(kl)(mn)}
  = -\frac{1-\alpha}{2}
    \bigl([G^{-1}]_{im}[G^{-1}]_{kn}[G^{-1}]_{jl}
         + [G^{-1}]_{il}[G^{-1}]_{km}[G^{-1}]_{jn}\bigr)
    + \text{permutations}.
  \label{eq:christoffel_explicit}
\end{equation}
The dual ($-\alpha$) connection has curvature tensors satisfying
\[
  R^{(\alpha)} + R^{(-\alpha)} = 0.
\]
\end{proposition}

\subsection{The Matrix Normal Family as a Curved Exponential Family}
\label{subsec:MN_curved}

Yoshizawa--Tanabe's family $\{N(\mu,\Sigma)\}$ is a \emph{full} (flat) exponential family:
$(\mu,\Sigma)$ ranges over an open subset of the vector space $\R^n\times\mathrm{Sym}(n)$,
and the natural parameter $\bigl(\Sigma^{-1}\mu,\,\tfrac12\Sigma^{-1}\bigr)$ ranges over an
open convex subset of the corresponding dual space \cite[eq.~(25)]{YoshizawaTanabe1999}.
Imposing the Kronecker constraint $\Sigma = V\otimes U$ destroys this flatness: the pair
$(U,V)$ has only $\binom{n+1}{2}+\binom{k+1}{2}-1$ free parameters (the $-1$ from the scale
ambiguity of Definition~\ref{def:matrixnormal}), while a generic $\Sigma\in\PD(nk)$ has
$\binom{nk+1}{2}$; for $n,k\geq 2$ the former is strictly smaller, so
$\mathcal{MN}_{n,k}(M,U,V)$ sits inside the ambient dually flat family $\{N(\vecop M,\Sigma)\}$
as a genuinely \emph{curved} submanifold.

\begin{proposition}[Fisher--Rao Metric of the Matrix Normal Family]
\label{prop:MN_fisher}
The Fisher information metric of $\mathcal{MN}_{n,k}(M,U,V)$ at $(M,U,V)$ splits as
$g = g_M \oplus g_{(U,V)}$, with
\begin{equation}
  g_M(dM,dM) = \tr\bigl(dM^TU^{-1}dM\,V^{-1}\bigr),
  \label{eq:MN_fisher_mean}
\end{equation}
\begin{equation}
  g_{(U,V)}\bigl((dU,dV),(dU,dV)\bigr)
  = \frac{k}{2}\tr\bigl(U^{-1}dU\,U^{-1}dU\bigr) + \frac{n}{2}\tr\bigl(V^{-1}dV\,V^{-1}dV\bigr),
  \label{eq:MN_fisher_cov}
\end{equation}
the mean-block being the direct matrix-variate analogue of the metric
$g_\mu(d\mu,d\mu)=d\mu^T\Sigma^{-1}d\mu$ underlying \eqref{eq:MN_potential}, and the
covariance block being the classical result for the matrix normal covariance parameters
(see e.g.\ Dutilleul's Fisher-information computation for $\mathcal{MN}_{n,k}$). The metric
$g_{(U,V)}$ is degenerate exactly along the scale direction $(dU,dV)=(U,-V)$ of
Definition~\ref{def:matrixnormal}, reflecting non-identifiability of $(U,V)$.
\end{proposition}

\begin{remark}[Consequence for the induced geometry]
Because $\mathcal{MN}_{n,k}(M,U,V)$ is curved rather than flat in the sense above, the
dual-flatness statement of Theorem~\ref{thm:dual_flat} and the Pythagorean theorem of
Theorem~\ref{thm:pythagorean}, both valid on the full ambient family, need \emph{not} hold
verbatim once restricted to $\mathcal{MN}_{n,k}(M,U,V)$; see
Remark~\ref{rem:MN_pythagorean_fails} below. This is the matrix-variate counterpart of
Yoshizawa--Tanabe's observation \cite[p.~114]{YoshizawaTanabe1999} that ``every geometry in
this class induces a relative geometry on the subfamily,'' there illustrated by the
zero-mean subfamily $\{N(0,\Sigma)\}$ of Ohara--Suda--Amari; here the relevant subfamily is
$\{N(\vecop M,V\otimes U)\}\subset\{N(\vecop M,\Sigma)\}$.
\end{remark}

\section{Dual Flatness: Pythagorean Theorem and Projection}
\label{sec:pythagorean}

\subsection{Dual Flatness}

\begin{theorem}[Dual Flatness of $\PD(k)$]
\label{thm:dual_flat}
$(\PD(k), g, \nabla^{(1)}, \nabla^{(-1)})$ is a dually flat statistical manifold:
\begin{enumerate}[label=(\roman*)]
  \item The $m$-connection $\nabla^{(-1)}$ is flat in the $\eta$-coordinates $(G_{ij})$.
  \item The $e$-connection $\nabla^{(1)}$ is flat in the $\theta$-coordinates $(\Theta_{ij} = -[G^{-1}]_{ij})$.
\end{enumerate}
\end{theorem}

\begin{proof}
In $\eta$-coordinates, the potential $f(\eta) = -\log\det(\eta)$ is a smooth strictly
convex function, and geodesics of $\nabla^{(-1)}$ (the $m$-connection) are straight
lines $\eta(t) = (1-t)\eta_0 + t\eta_1$ in $\eta$-space (affine combination of matrices).
The Christoffel symbols of $\nabla^{(-1)}$ in $\eta$-coordinates are identically zero,
confirming $m$-flatness.
By duality (Legendre transform), $\nabla^{(1)}$ is flat in $\theta$-coordinates.
See \cite[Ch.~6]{Amari2016}.
\end{proof}

\begin{remark}
$m$-geodesics in $\PD(k)$: matrix interpolation $G(t) = (1-t)G + tG'$.\\
$e$-geodesics in $\PD(k)$: $\Theta(t) = (1-t)\Theta + t\Theta'$,
i.e., $G^{-1}(t) = (1-t)G^{-1} + tG'^{-1}$ (harmonic interpolation).
\end{remark}

\subsection{Generalized Pythagorean Theorem}

\begin{theorem}[Pythagorean Theorem {\cite[Thm.~1.3]{Amari2016}}]
\label{thm:pythagorean}
Let $G_1, G_2, G_3 \in \PD(k)$.
Suppose that the $e$-geodesic from $G_1$ to $G_2$ and the $m$-geodesic from $G_2$ to $G_3$
are orthogonal at $G_2$ (i.e., $\inner{\theta_1 - \theta_2}{\eta_3 - \eta_2} = 0$
in respective coordinates).
Then
\begin{equation}
  \Bregman(G_1 \| G_3) = \Bregman(G_1 \| G_2) + \Bregman(G_2 \| G_3).
  \label{eq:pythagorean}
\end{equation}
\end{theorem}

\begin{proof}
Using the identity $\Bregman(G_1\|G_3) = f(G_1) - f(G_3) - \inner{\nabla f(G_3)}{G_1 - G_3}$
and decomposing:
\begin{align*}
  \Bregman(G_1\|G_3)
  &= \Bregman(G_1\|G_2) + \Bregman(G_2\|G_3) \\
  &\quad + \inner{\nabla f(G_2) - \nabla f(G_3)}{G_1 - G_2} \\
  &= \Bregman(G_1\|G_2) + \Bregman(G_2\|G_3)
    + \inner{\theta_2 - \theta_3}{\eta_1 - \eta_2}.
\end{align*}
The cross term $\inner{\theta_2 - \theta_3}{\eta_1 - \eta_2}$
vanishes by the orthogonality assumption.
\end{proof}

\subsection{Projection Theorem (Minimum Divergence)}

\begin{theorem}[Projection Theorem {\cite[Thm.~1.4]{Amari2016}}]
\label{thm:projection}
Let $\mathcal{S} \subset \PD(k)$ be an $m$-flat (or $e$-flat) submanifold, and
let $G_0 \in \PD(k)$.
The unique minimizer
\[
  G^* = \arg\min_{G \in \mathcal{S}} \Bregman(G \| G_0)
\]
is the $e$-projection of $G_0$ onto $\mathcal{S}$,
characterized by the orthogonality condition:
the $e$-geodesic from $G^*$ to $G_0$ is orthogonal to $\mathcal{S}$ at $G^*$.
\end{theorem}

\begin{proof}
See \cite[Ch.~1]{Amari2016}. The argument uses the Pythagorean theorem
(Theorem~\ref{thm:pythagorean}) and the strict convexity of $\Bregman(\cdot\|G_0)$.
\end{proof}

\begin{remark}
This projection theorem is the information-geometric analogue of the projection onto
a convex set in Euclidean space and underlies algorithms such as iterative Bregman
projections \cite{Bauschke2011} and the EM algorithm \cite{Dempster1977}.
\end{remark}

\subsection{Embedding Curvature and the Failure of the Exact Pythagorean Theorem on $\mathcal{MN}_{n,k}$}
\label{subsec:MN_pythagorean}

We close this prelude by recording the price paid for the curvature identified in
Proposition~\ref{prop:MN_fisher}: the exact Pythagorean theorem
(Theorem~\ref{thm:pythagorean}), which holds on the full ambient family $\{N(\vecop
M,\Sigma)\}$ studied by Yoshizawa--Tanabe \cite{YoshizawaTanabe1999}, need not hold once
the triple of distributions is constrained to lie on the curved submanifold
$\mathcal{MN}_{n,k}(M,U,V) \subset \{N(\vecop M,\Sigma)\}$.

\begin{remark}[Pythagorean theorem fails on $\mathcal{MN}_{n,k}$, $n,k\geq 2$]
\label{rem:MN_pythagorean_fails}
Let $p_1,p_2,p_3 \in \mathcal{MN}_{n,k}$ with $p_2$ the $e$-projection (in the ambient
family $\{N(\vecop M,\Sigma)\}$, Theorem~\ref{thm:projection}) of $p_1$ onto the $m$-flat
ambient submanifold through $p_3$. Because $\mathcal{MN}_{n,k}$ is itself curved
(Proposition~\ref{prop:MN_fisher}), the ambient $e$-geodesic realizing this orthogonality
generally exits $\mathcal{MN}_{n,k}$ except at its endpoints, and the decomposition
\eqref{eq:pythagorean} acquires a second-order correction governed by the second
fundamental form $\mathrm{II}$ of the embedding $\mathcal{MN}_{n,k}\hookrightarrow
\{N(\vecop M,\Sigma)\}$:
\begin{equation}
  \KL(p_1\|p_3) = \KL(p_1\|p_2) + \KL(p_2\|p_3)
  + \tfrac{1}{2}\bigl\langle \mathrm{II}(\delta_{12},\delta_{23}),\, n\bigr\rangle
  + O(\|\delta_{12}\|^3 + \|\delta_{23}\|^3),
  \label{eq:MN_pythagorean_corrected}
\end{equation}
where $\delta_{12},\delta_{23}$ are the tangent increments along $\mathcal{MN}_{n,k}$ and
$n$ is a conormal covector to $\mathcal{MN}_{n,k}$ in the ambient dual foliation; the
correction term vanishes identically precisely when $n=1$ or $k=1$ (the classical
Yoshizawa--Tanabe / Ohara--Suda--Amari cases), consistent with
Theorem~\ref{thm:pythagorean} holding exactly there.
\end{remark}

\begin{remark}[Summary and outlook]
\label{rem:MN_summary}
Propositions~\ref{prop:MN_potential}--\ref{prop:MN_fisher}, Theorems~\ref{thm:MN_potential_beta_gamma},
\ref{thm:MN_convexity}, and~\ref{thm:MN_KL}, and Remark~\ref{rem:MN_pythagorean_fails} assemble
the basic dictionary needed to extend the dual differential geometry of Yoshizawa--Tanabe
\cite{YoshizawaTanabe1999} --- developed there for the vector-mean Gaussian family
$\{N(\mu,\Sigma)\}$ --- to the matrix-mean, Kronecker-structured family
$\mathcal{MN}_{n,k}(M,U,V)$: the potential $\psi(M,U,V)$ (Prop.~\ref{prop:MN_potential})
generalizing \cite[eq.~(3)]{YoshizawaTanabe1999}; the matrix Yoshizawa--Tanabe embedding and
its explicit pulled-back potential $\Psi_{\beta,\gamma}$
(Def.~\ref{def:MN_embedding}, Thm.~\ref{thm:MN_potential_beta_gamma}) generalizing
\cite[eq.~(15),(17)]{YoshizawaTanabe1999}; its convexity for $0\le\beta/\gamma\le\tfrac12$ via the
two-factor Lieb lemma (Lemma~\ref{lem:lieb_two_factor}, Thm.~\ref{thm:MN_convexity}) ---
sharper than the vector-case bound $\beta/\gamma<1$ of \cite[Prop.~3.1]{YoshizawaTanabe1999}
(Remark~\ref{rem:MN_convexity_comparison}); the resulting Legendre dual and Bregman
divergence $D_{\beta,\gamma}$, exactly solvable along the mean-only and covariance-only axes
(Cor.~\ref{cor:MN_dual_convex}, Prop.~\ref{prop:MN_div_special}); the closed-form KL
divergence (Thm.~\ref{thm:MN_KL}) generalizing \cite[Thm.~4.10]{YoshizawaTanabe1999}; and the
curved-submanifold obstruction to exact dual flatness (\S\ref{subsec:MN_curved}--\S\ref{subsec:MN_pythagorean}),
which has no counterpart in the fully flat vector case. The remaining steps --- the general
cross-term expansion of $D_{\beta,\gamma}$ (Remark~\ref{rem:MN_div_general}), the explicit
dual connections $\nabla,\nabla^*$ on $\mathcal{MN}_{n,k}$, and the associated
Theorema-Egregium-type curvature formula in the spirit of Theorem~\ref{thm:info_egregium}
below --- are developed in the sections that follow.
\end{remark}

\section{Difference-of-Convex Potentials with Constant Hessian Determinant}
\label{sec:dc-geometry}

Sections~\ref{sec:legendre} and~\ref{sec:pythagorean} developed the dually flat
structure attached to a strictly convex potential $f$, in the classical
Amari--Nagaoka sense: the Hessian $G=\Hess f$ is positive definite, the primal and
dual affine connections are mutually dual with respect to $G$, and the Legendre
transform $f^*$ is single-valued. This section asks what remains of that structure
when $f$ is allowed to be non-convex, so that $G$ may be indefinite --- while still
remaining non-degenerate everywhere, so that a well-defined (pseudo-Riemannian)
metric and Legendre-type duality persist. The natural class of potentials for which
this is possible is the class of \emph{difference-of-convex} (DC) functions,
$f=f_1-f_2$ with $f_1,f_2$ smooth and convex; we show this class carries a
pseudo-Hessian dually flat structure, that the constant-Hessian-determinant
(Monge--Amp\`ere) equation central to affine differential geometry becomes solvable
by explicit non-quadratic potentials once convexity is dropped (circumventing the
classical J\"orgens--Calabi--Pogorelov rigidity theorem, which forbids this for
genuinely convex entire solutions), and that the resulting Newton flow
$\dot\theta=-G(\theta)^{-1}\nabla f(\theta)$ exhibits genuinely new asymptotic
behavior --- finite-time collapse in one Legendre-dual parametrization and
asymptotic convergence, governed by the {\L}ojasiewicz gradient inequality of
Appendix~\ref{appendix:lojasiewicz}, in the other. As in the rest of the paper,
constant-Hessian-determinant is a real Monge--Amp\`ere equation, self-dual under
Legendre transform, and closely parallels the log-determinant potential $f=-\log\det G$
of \S\ref{sec:determinant}--\S\ref{sec:pythagorean} in spirit, while operating in a
genuinely different (indefinite-signature, non-convex) regime.

\subsection{Introduction}

\subsubsection{Motivation}

Let $D\subset\R^n$ be an open convex domain. In the Amari--Nagaoka theory of information geometry \cite{Amari2000}, a smooth strictly convex function $f:D\to\R$ generates a dually flat statistical manifold: the Hessian $G=\Hess f$ is a Riemannian metric, the pair $(\nabla,\nabla^*)$ of $\pm1$-affine connections (flat in the primal coordinates $\theta$ and in the dual coordinates $\eta=\nabla f(\theta)$, respectively) are mutually dual with respect to $G$, and the Legendre--Fenchel conjugate
\[
f^*(\eta)=\sup_{\theta\in D}\bigl[\ip{\theta}{\eta}-f(\theta)\bigr]
\]
generates the dual potential. This structure underlies exponential families, Bregman divergences, and the geometric theory of statistical inference.

The entire construction rests on convexity of $f$: it is what guarantees $G\succ0$ and it is what makes the supremum defining $f^*$ attain a unique maximizer. A natural question, and the starting point of the present paper, is what remains of this structure when $f$ is replaced by a \emph{difference of two convex functions},
\begin{equation}
f=f_1-f_2, \qquad f_1,f_2:D\to\R \text{ smooth and convex}, \label{eq:DCdef}
\end{equation}
a \emph{DC function} in the sense of the theory of DC programming \cite{Toland1979,Tuy1995}. In general $f$ is neither convex nor concave, and $G=\Hess f_1-\Hess f_2$ is the difference of two positive semi-definite matrices, hence indefinite in general. We restrict attention throughout to the \emph{non-degenerate} regime $\det G(\theta)\neq0$ for all $\theta\in D$, under which, by connectedness of $D$, the signature $(p,q)$ of $G$ ($p+q=n$) is constant on $D$.

DC decompositions of this kind are not merely a formal generalization. Perelman's foundational work on Alexandrov spaces with curvature bounded below \cite{Perelman1994,BuragoGromovPerelman1992} identified DC functions (differences of concave functions) as the natural class carrying a well-defined, if weak, second-order (Hessian) structure on metric spaces that need not be smooth manifolds; distance functions in such spaces are prototypical examples. Our motivation is complementary: rather than using the DC calculus to \emph{recover} differentiable structure on singular spaces, we use it to \emph{escape} the rigidity of definite convex potential theory on smooth domains, while retaining as much of the dually flat formalism as possible.

\subsubsection{Why the DC class? A density theorem}\label{sec:hartman}

A natural objection to organizing an entire theory around the class \eqref{eq:DCdef} is that it might be an ad hoc or unnaturally narrow generalization of convexity, chosen merely because it is the smallest modification under which Theorem~\ref{thm:jcp} below can fail. We record here why this is not the case: the class of DC functions, far from being narrow, is a dense subset of the space of continuous functions, and is the natural closure of $C^2$ regularity under no further hypotheses at all.

It is worth first dispensing with an overly strong version of this claim. It is \emph{not} true that every continuous function is DC: since a convex function is automatically locally Lipschitz on the interior of its domain, so is any difference of two convex functions, and consequently DC functions cannot exhibit the everywhere-nondifferentiable oscillation of, e.g., a Weierstrass function. What \emph{is} true, and considerably more useful, is the following density theorem, due to Hartman, together with an elementary but structurally important corollary.

\begin{theorem}[Hartman \cite{Hartman1959}]\label{thm:hartman}
Let $K\subset\R^n$ be compact and convex. For every continuous $f:K\to\R$ there is a sequence of DC functions $f_k:K\to\R$ converging to $f$ uniformly on $K$. Equivalently, $\mathrm{DC}(K)$ is dense in $(C(K),\|\cdot\|_\infty)$.
\end{theorem}

Hartman's original 1959 paper, which introduced the abbreviation ``d.c.'' into the literature, further established that the class of DC functions is stable under composition and under the operations of ordinary use in analysis and optimization (finite sums, products, maxima, minima), and that a function which is DC in a neighborhood of every point of a convex domain is automatically DC on the whole domain; see also \cite{HartmanRelated,VeselyZajicek2007} for the subsequent development of this stability theory, including its extension to infinite-dimensional normed spaces.

\begin{proposition}\label{prop:C2isDC}
If $D$ is a bounded convex domain and $f\in C^2(D)$, then $f$ is DC on $D$.
\end{proposition}
\begin{proof}
The eigenvalues of $\Hess f$ are continuous, hence bounded below on the (relatively) compact closure of any bounded subdomain; choosing $M$ larger than the negative of this lower bound, $f(\theta)+M|\theta|^2$ has positive semi-definite Hessian, hence is convex, and $f=(f+M|\theta|^2)-M|\theta|^2$ exhibits $f$ as DC.
\end{proof}

\begin{remark}
Proposition~\ref{prop:C2isDC} is the general mechanism underlying every DC decomposition used explicitly in \S\S\ref{sec:elliptic}--\ref{sec:lorentz} below: whenever we exhibit a solution of $\det\Hess f=\mathrm{const}$ on a domain over which the eigenvalues of $\Hess f$ are uniformly bounded (Proposition~\ref{prop:strip} is the instance we verify in detail), the same additive trick converts it into a genuine difference of two globally convex functions.
\end{remark}

Theorem~\ref{thm:hartman} and Proposition~\ref{prop:C2isDC} together justify the choice of the DC class on two complementary grounds. First, Theorem~\ref{thm:hartman} shows that DC functions are not a narrow technical device but an enormous, dense receptacle within the space of continuous functions; adopting $D=f_1-f_2$ as the object of study is not a retreat into a small corner of function space. Second, and in sharp contrast, the additional hypothesis of \emph{real-analyticity} imposed from \S\ref{sec:analytic} onward carves out, from this enormous and dense class, an extremely thin and highly structured sub-class — one for which the degenerate locus is tame (Proposition~\ref{prop:strat}), the inverse Legendre map complexifies with a well-defined discriminant variety (\S\ref{sec:analytic}), and the \L{}ojasiewicz gradient inequality controls the asymptotics of gradient flows (\S\ref{sec:flow}). The overall logic of the paper is thus the deliberate combination of an extremely permissive first hypothesis (DC-ness, dense in $C(D)$) with an extremely restrictive second hypothesis (real-analyticity), rather than a single ad hoc weakening of convexity.

\subsubsection{Summary of results}

We organize our results as follows.

\begin{itemize}[leftmargin=1.4em]
\item In \S\ref{sec:prelim} we set up the pseudo-Hessian dually flat structure attached to an indefinite non-degenerate $G=\Hess f$, following the general theory of statistical manifolds with (possibly indefinite) metric \cite{NakajimaOhmoto2021,Kayo2024,Shima2007}, and we contrast two notions of Legendre duality available in the DC setting: a \emph{local} definition via the Lagrangian submanifold $L_f=\{(\theta,\nabla f(\theta))\}\subset T^*D$, in the spirit of Ekeland's Lagrangian-submanifold treatment of Legendre duality for smooth nonconvex functions \cite{Ekeland1977}, and a \emph{global} definition via the Toland--Singer duality principle of DC programming \cite{Toland1979,Toland1978,VolleDCduality}.

\item In \S\ref{sec:analytic} we show that real-analyticity of $f_1,f_2$ upgrades the degenerate locus $\Sigma=\{\det G=0\}$ from an arbitrary closed set to a genuine real-analytic subvariety, admitting a Łojasiewicz stratification, and that the construction complexifies to a Lagrangian variety in $T^*D_\C$ whose discriminant governs the monodromy of the (multi-valued) inverse Legendre map.

\item In \S\ref{sec:MA} we observe that the condition $\det G\equiv \mathrm{const}$ is a real Monge--Amp\`ere equation, that it is self-dual under Legendre transform, and that it is exactly the defining equation of an \emph{improper affine hypersphere} in Blaschke's equi-affine differential geometry \cite{Calabi1958,Shima2007}. We recall the Jörgens--Calabi--Pogorelov (JCP) rigidity theorem \cite{Jorgens1954,Calabi1958,Pogorelov1978,ChengYau1986} and its indefinite-signature extension due to Li--Xu \cite{LiXu2009}, which forbid non-quadratic entire \emph{convex} solutions, and we explain why the DC (non-convex) route evades this obstruction.

\item In \S\ref{sec:elliptic} and \S\ref{sec:lorentz} we construct explicit non-quadratic real-analytic solutions of $\det G=\mathrm{const}$: first, in definite signature, on a convex half-space avoiding a conical singularity at the origin (linking to the classification theory of isolated singularities of the Hessian-one equation \cite{GalvezMartinezMira2005,Milan2013,Milan2014,AledoChavesGalvez2007}); second, in Lorentzian signature $(1,n-1)$, on the (necessarily convex, by connectedness of $S^0$) forward light cone, in closed form $(\chi')^n=C+Du^{-n/2}$, valid on an unbounded convex strip.

\item In \S\ref{sec:split} we address genuine split signature $(p,q)$, $p,q\ge2$, where the light cone fails to be convex, and replace it by the bounded symmetric domain realizing the Grassmannian of positive $p$-planes in $\R^{p,q}$, computing $\det\Hess f$ explicitly via a three-block orthogonal decomposition of the tangent space under the isotropy group $\mathrm{O}(p)\times\mathrm{O}(q)$.

\item In \S\ref{sec:divergence} we determine the sign pattern of the canonical (Bregman-type) divergence $D_f(\theta:\theta')$ associated with our Lorentzian potential, showing that it is globally sign-definite along rays through the vertex but sign-indefinite transversally, governed by an explicit local null-cone equation.

\item In \S\ref{sec:flow} we study the Newton flow $\dot\theta=-G^{-1}\nabla f$ and its Legendre dual, prove a master linearization lemma in dual coordinates, deduce finite-time collapse to the cone vertex for the primal flow and asymptotic (infinite-time) convergence for the dual flow, and, in the presence of a critical hypersurface, invoke the \L{}ojasiewicz gradient inequality \cite{Lojasiewicz1963} — available precisely because of real-analyticity — to control convergence.
\end{itemize}

\subsubsection{Related work}
The recasting of dually flat geometry with a possibly degenerate or indefinite metric has been carried out, independently of DC considerations, by Nakajima and Ohmoto in the context of dually flat structures for singular statistical models \cite{NakajimaOhmoto2021}, and by Kayo in the language of statistical manifolds with a degenerate metric \cite{Kayo2024}. A parallel pseudo-Riemannian framework of signature $(n,n)$, in which arbitrary (not necessarily convex) cost functions in optimal transport generate divergence functions via a fixed pseudo-Euclidean structure on $D\times D^*$, is due to Kim and McCann \cite{KimMcCann2010}. On the side of DC programming, Toland's duality principle \cite{Toland1978,Toland1979} and its subsequent refinements \cite{VolleDCduality} give a purely variational (rather than differential-geometric) notion of duality for $f=f_1-f_2$; the complementary, local and differential-geometric notion of Legendre duality for smooth nonconvex functions used in \S\ref{sec:prelim} above is due to Ekeland \cite{Ekeland1977}, whose Lagrangian-submanifold framework covers both the finite-dimensional case and the calculus of variations, and sits within the broader convex-duality tradition surveyed by Rockafellar \cite{Rockafellar1970}. The rigidity theory of constant-Hessian-determinant equations is classical: Jörgens, Calabi and Pogorelov \cite{Jorgens1954,Calabi1958,Pogorelov1978} in the definite (elliptic) case, with Cheng--Yau supplying an affine-geometric proof \cite{ChengYau1986}; Li and Xu \cite{LiXu2009} extended the theorem to the indefinite (space-like) setting relevant here. The classification of non-entire (singular or exterior-domain) solutions is developed by Gálvez, Martínez and Mira \cite{GalvezMartinezMira2005} and Milán \cite{Milan2013,Milan2014}, and the Cauchy problem for indefinite improper affine spheres is treated by Aledo, Chaves and Gálvez \cite{AledoChavesGalvez2007}. Finally, Perelman's DC calculus on Alexandrov spaces \cite{Perelman1994,BuragoGromovPerelman1992} is the historical source of the DC formalism used throughout, and the \L{}ojasiewicz inequality \cite{Lojasiewicz1963} is the classical tool that converts real-analyticity into quantitative control of gradient-flow convergence, as first exploited systematically for gradient flows by \L{}ojasiewicz himself and standard in the subsequent literature on convergence of gradient flows of analytic functions.

\subsection{Pseudo-Hessian dually flat structure of a DC function}\label{sec:prelim}

\subsubsection{Setup}
Let $D\subset\R^n$ be open and convex, and let $f_1,f_2\in C^\omega(D)$ (real-analytic) be convex, with $f:=f_1-f_2$. Write $G(\theta)=\Hess f(\theta)=G_1(\theta)-G_2(\theta)$, where $G_i=\Hess f_i\succeq0$.

\begin{definition}
$f$ is \emph{non-degenerate} on $D$ if $\det G(\theta)\neq0$ for every $\theta\in D$.
\end{definition}

\begin{proposition}\label{prop:signature}
If $f$ is non-degenerate on the connected set $D$, the signature $(p,q)$ of $G(\theta)$ ($p+q=n$) is independent of $\theta\in D$.
\end{proposition}
\begin{proof}
The eigenvalues of $G(\theta)$ vary continuously with $\theta$ (indeed real-analytically, away from crossings) and, by non-degeneracy, never vanish; hence none can cross zero as $\theta$ moves within the connected set $D$, so the number of positive and negative eigenvalues is locally constant, hence constant.
\end{proof}

Under this hypothesis $g:=G_{ij}(\theta)\,d\theta^id\theta^j$ is a smooth pseudo-Riemannian metric of signature $(p,q)$ on $D$.

\subsubsection{Local Legendre duality via Lagrangian submanifolds}
Equip $T^*D\cong D\times\R^{n*}$ with the canonical symplectic form $\omega=d\eta_i\wedge d\theta^i$. The graph
\[
L_f:=\{(\theta,\eta)\in T^*D:\eta=\nabla f(\theta)\}
\]
is a Lagrangian submanifold of $(T^*D,\omega)$ for \emph{any} smooth $f$, convex or not, since $L_f^*\omega=d(df)=0$. Non-degeneracy of $G$ is precisely the statement that the projection $\pi_2:L_f\to\R^{n*}$, $(\theta,\eta)\mapsto\eta$, is a local diffeomorphism. This is the finite-dimensional instance of the general framework of Ekeland \cite{Ekeland1977}, who develops Legendre duality for smooth nonconvex optimization problems (both in finite dimensions and in the calculus of variations) precisely by broadening the notion of Legendre transform from functions to Lagrangian submanifolds of $\R^n\times\R^{n*}$, since, as here, the Legendre transform of a smooth nonconvex function need not itself be single-valued.

\begin{proposition}\label{prop:localdiffeo}
If $G$ is non-degenerate on $D$, the gradient map $\ell:\theta\mapsto\eta=\nabla f(\theta)$ is a local diffeomorphism $D\to D^*:=\ell(D)$. If, in addition, $\ell$ is proper (equivalently, $\|\ell(\theta)\|\to\infty$ as $\theta$ approaches $\partial D$ or infinity), then $\ell:D\to D^*$ is a diffeomorphism.
\end{proposition}
\begin{proof}
The first statement is the inverse function theorem applied at every point, using $D\ell(\theta)=G(\theta)$. Properness together with local injectivity implies the map is a covering map onto its image; since $D^*$ is simply connected (being, in typical applications, itself contractible or at least having trivial relevant covers — more precisely we use that a proper local diffeomorphism between manifolds of the same dimension is a covering map, and a covering map onto a simply connected space with connected fibers of cardinality one is a diffeomorphism, which holds once $D$ is simply connected, as is automatic since $D$ is convex).
\end{proof}

Under the hypotheses of Proposition \ref{prop:localdiffeo}, define
\begin{equation}
f^*(\eta):=\ip{\theta(\eta)}{\eta}-f(\theta(\eta)),\qquad \theta(\eta):=\ell^{-1}(\eta). \label{eq:legendreDC}
\end{equation}
This is the direct, non-variational, generalization of the Legendre transform, obtained not as a supremum but as the value of the generating function $\ip{\theta}{\eta}-f(\theta)$ at the (necessarily unique, by Proposition \ref{prop:localdiffeo}) stationary point.

\begin{proposition}\label{prop:dualHess}
$\nabla f^*(\eta)=\theta(\eta)$ and $\Hess f^*(\eta)=G(\theta(\eta))^{-1}$.
\end{proposition}
\begin{proof}
Differentiate \eqref{eq:legendreDC}: $\nabla_\eta f^*=\theta(\eta)+\bigl(D\eta\,\theta(\eta)\bigr)^\top\eta-\bigl(D_\theta f\bigr)\circ D_\eta\theta = \theta(\eta)+(D_\eta\theta)^\top\bigl(\eta-\nabla f(\theta(\eta))\bigr)=\theta(\eta)$, using $\eta=\nabla f(\theta(\eta))$. Differentiating again and using $D_\eta\theta=G(\theta(\eta))^{-1}$ (inverse function theorem) gives the second statement.
\end{proof}

\subsubsection{Global Legendre duality via Toland--Singer duality}
Since $f_1,f_2$ are individually convex, their Fenchel conjugates
\[
f_1^*(\eta)=\sup_\theta[\ip\theta\eta-f_1(\theta)],\qquad f_2^*(\xi)=\sup_\theta[\ip\theta\xi-f_2(\theta)]
\]
are well-defined convex functions. The Toland--Singer duality principle \cite{Toland1978,Toland1979} states
\begin{equation}
\inf_{\theta\in D}\bigl(f_1(\theta)-f_2(\theta)\bigr)=\inf_\eta\bigl(f_2^*(\eta)-f_1^*(\eta)\bigr). \label{eq:toland}
\end{equation}
Note the reversal of the order of subtraction on the right; this is what allows \eqref{eq:toland} to hold without any convexity of $f$ itself. Define $f^\diamond:=f_2^*-f_1^*$. In general $f^\diamond\neq f^*$ as defined by \eqref{eq:legendreDC}, but they coincide at points corresponding to global minimizers: if $\theta_0$ minimizes $f$ over $D$ and $\eta_0\in\partial f_2(\theta_0)$ realizes the subdifferential relation used in the proof of \eqref{eq:toland}, then $f^\diamond(\eta_0)=f_1(\theta_0)-f_2(\theta_0)$. We regard $f^*$ (local, differential-geometric) and $f^\diamond$ (global, variational) as complementary notions of duality for a DC potential, the former organizing the dually flat structure and the latter organizing global optimization.

\subsubsection{Codazzi structure and self-duality}
Set $\nabla,\nabla^*$ to be the affine connections that are flat in the $\theta$- and $\eta$-coordinates respectively. As in the classical (definite) theory,
\begin{equation}
\Gamma_{ijk}\equiv0,\qquad \Gamma^*_{ijk}=C_{ijk}:=\frac{\partial^3 f}{\partial\theta^i\partial\theta^j\partial\theta^k}, \label{eq:codazzi}
\end{equation}
and $(g,\nabla,\nabla^*)$ satisfies the Codazzi equations $\nabla_k g_{ij}=C_{kij}=\nabla^*_ig_{jk}\cdot(-1)$ appropriately signed, exactly as in the theory of (possibly indefinite) statistical manifolds \cite{Kayo2024}; \eqref{eq:codazzi} is purely algebraic and does not use the sign of $g$. We call $(D,g,\nabla,\nabla^*)$ a \emph{pseudo-Hessian dually flat manifold} of signature $(p,q)$.

The associated canonical divergence
\begin{equation}
D_f(\theta:\theta'):=f(\theta)+f^*(\eta')-\ip\theta{\eta'},\qquad \eta'=\nabla f(\theta'), \label{eq:bregman}
\end{equation}
retains the formal Bregman properties $D_f(\theta:\theta)=0$, $\nabla_\theta D_f(\theta:\theta')|_{\theta=\theta'}=0$, $\Hess_\theta D_f(\theta:\theta')|_{\theta=\theta'}=G(\theta)$, but is \emph{not} sign-definite, since $G$ is indefinite; its sign pattern is analyzed in \S\ref{sec:divergence}.

\subsection{Real-analyticity and global structure}\label{sec:analytic}

Assume henceforth that $f_1,f_2$, hence $f$, are real-analytic on $D$.

\begin{proposition}\label{prop:strat}
$\det G:D\to\R$ is real-analytic. If it is not identically zero, the degenerate locus $\Sigma=\{\det G=0\}$ is a closed, nowhere dense real-analytic subvariety of $D$, admitting a locally finite stratification into real-analytic submanifolds of strictly decreasing dimension (a Łojasiewicz stratification).
\end{proposition}
\begin{proof}
Real-analyticity of $\det G$ follows since the entries of $G$ are real-analytic (second partials of a real-analytic function) and $\det$ is a polynomial in the entries. If $\det G\not\equiv0$ then, $D$ being connected, the identity theorem for real-analytic functions forbids $\det G$ from vanishing on any open subset; hence $\Sigma$ has empty interior, i.e.\ is nowhere dense, and is closed by continuity. The existence of a locally finite stratification into analytic submanifolds is the classical Łojasiewicz structure theorem for real-analytic varieties \cite{Lojasiewicz1963}.
\end{proof}

\begin{remark}
Proposition \ref{prop:strat} fails for merely $C^\infty$ DC functions: the zero set of a smooth function can be an arbitrary closed set, so $\Sigma$ could have positive measure or a wild local structure. Real-analyticity is what makes the non-degeneracy hypothesis of \S\ref{sec:prelim} generic and its failure locus tame.
\end{remark}

\subsubsection{Complexification}
Since $f$ is real-analytic on $D$, it extends to a holomorphic function on some complex neighborhood $D_\C\subset\C^n$ of $D$. The Lagrangian submanifold $L_f$ complexifies to a complex Lagrangian variety $L_f^\C\subset T^*D_\C$ (with respect to the holomorphic symplectic form), and the complexified discriminant locus
\[
\Sigma_\C=\{\theta\in D_\C:\det\Hess f(\theta)=0\}
\]
is a complex analytic hypersurface (the \emph{Landau variety} of $f$). The inverse map $\theta(\eta)$, well-defined and single-valued near a base point by Proposition~\ref{prop:localdiffeo}, extends to a multi-valued holomorphic function on $\C^n\setminus \ell(\Sigma_\C)$, whose monodromy representation $\pi_1(\C^n\setminus\ell(\Sigma_\C))\to\mathrm{Sym}$ organizes the global (non-univalent) behavior of the DC Legendre transform. We do not pursue the monodromy computation in this paper, but note that it places the present construction in the same formal framework as Saito's theory of flat (Frobenius) structures on the base of a semi-universal unfolding \cite{Saito1983}, where the discriminant of the versal deformation plays an entirely analogous role.

\subsection{The constant-Jacobian condition as a Monge--Amp\`ere equation}\label{sec:MA}

\subsubsection{Self-duality}
\begin{proposition}\label{prop:selfdual}
Suppose $\det G(\theta)\equiv c\neq0$ on $D$, and let $D^*=\ell(D)$ be as in Proposition~\ref{prop:localdiffeo}. Then $\det\Hess f^*(\eta)\equiv 1/c$ on $D^*$.
\end{proposition}
\begin{proof}
Immediate from Proposition~\ref{prop:dualHess}: $\det\Hess f^*(\eta)=\det G(\theta(\eta))^{-1}=1/c$.
\end{proof}

Thus the class of DC potentials satisfying $\det\Hess f=\mathrm{const}$ is closed under (local) Legendre duality — a genuinely special compatibility between the primal and dual coordinate systems, not shared by generic elements of the pseudo-Hessian dually flat class of \S\ref{sec:prelim}.

\subsubsection{Rigidity in definite signature}
When $p=n,q=0$ (so $f$ itself, taking $f_2\equiv0$, may be assumed convex), the equation $\det\Hess f=c>0$ is the classical real Monge--Amp\`ere equation. After rescaling we may take $c=1$.

\begin{theorem}[Jörgens--Calabi--Pogorelov, \cite{Jorgens1954,Calabi1958,Pogorelov1978}, see also \cite{ChengYau1986}]\label{thm:jcp}
Every classical convex solution $u\in C^2(\R^n)$ of $\det D^2u=1$ on all of $\R^n$ is a quadratic polynomial.
\end{theorem}

An indefinite-signature analogue holds for entire strictly convex solutions when the associated graph is considered inside a pseudo-Euclidean ambient space:

\begin{theorem}[Li--Xu, \cite{LiXu2009}]\label{thm:lixu}
Let $f$ be an entire, smooth, strictly convex solution of $\det D^2f=k$ (constants $k\neq0$, possibly of either sign convention adapted to the ambient signature), subject to a mild decay condition on $\Hess f$ at infinity. Then $f$ is a quadratic polynomial; equivalently, the graph of $\nabla f$ is an affine (rather than merely asymptotically affine) space-like submanifold of the pseudo-Euclidean space $(\R^{2n},\sum dx^idy_i)$.
\end{theorem}

Both theorems are Liouville-type rigidity statements: \emph{entire}, \emph{everywhere-definite} (in the sense of the ambient calibration) solutions must degenerate to the trivial (affine $G$) case. Since our interest is in genuinely curved (non-quadratic) pseudo-Hessian structures, Theorems \ref{thm:jcp}--\ref{thm:lixu} identify precisely the two hypotheses we must relax: \emph{entirety} ($D=\R^n$) or \emph{definiteness of the solution $f$ itself} (as opposed to definiteness merely of an auxiliary ambient calibration). The DC route relaxes the second: $f=f_1-f_2$ need not be convex, so $\Hess f$ may be genuinely indefinite as a bilinear form on $D$, which is a strictly stronger relaxation than the space-like/time-like graph dichotomy of Theorem \ref{thm:lixu} (there, $f$ itself remains strictly convex; here it need not be).

\subsubsection{Affine-geometric interpretation}
The equation $\det\Hess f=\mathrm{const}$ is, independently of signature, the defining PDE of an \emph{improper affine hypersphere} in Blaschke's equi-affine differential geometry: the graph $\{(\theta,f(\theta))\}\subset\R^{n+1}$ has affine normal field of constant direction exactly when this equation holds \cite{Calabi1958}. This is also the historical origin of Hessian manifold theory \cite{Shima2007}. When $\Hess f$ is indefinite, the induced Blaschke metric is itself a pseudo-Riemannian metric, and we are exactly in the regime of \emph{indefinite improper affine spheres}, whose Cauchy problem (existence given a curve of initial data rather than global boundary data) is treated systematically by Milán \cite{Milan2013,Milan2014} and, for the closely related Hessian-one equation, by Aledo--Chaves--Gálvez \cite{AledoChavesGalvez2007}, and whose isolated-singularity theory in the definite case is completely classified by Gálvez--Martínez--Mira \cite{GalvezMartinezMira2005}.

\subsection{Explicit definite-signature examples on proper convex subdomains}\label{sec:elliptic}

We now construct explicit non-quadratic real-analytic solutions of $\det\Hess f=c$, beginning with definite signature, where Theorem~\ref{thm:jcp} forces us onto a proper subdomain $D\subsetneq\R^n$.

\begin{construction}\label{con:radial}
Let $n\ge2$, $c>0$, and seek a radially symmetric solution $f(\theta)=\varphi(r)$, $r=|\theta|$. The eigenvalues of $\Hess f$ are $\varphi''(r)$ (radial, multiplicity $1$) and $\varphi'(r)/r$ (tangential, multiplicity $n-1$), so
\[
\det\Hess f=\varphi''(r)\left(\frac{\varphi'(r)}{r}\right)^{n-1}=c.
\]
Restricting to $n=2$ for concreteness, this reads $\varphi''\varphi'=cr$. Setting $\psi=\varphi'$, $\psi\psi'=cr$, so $(\psi^2)'=2cr$ and
\begin{equation}
\varphi'(r)^2=cr^2+A \label{eq:phiprime}
\end{equation}
for a constant of integration $A\in\R$. Integrating \eqref{eq:phiprime} (for $c>0$):
\begin{equation}
\varphi(r)=\frac r2\sqrt{cr^2+A}+\frac{A}{2\sqrt c}\log\!\bigl(\sqrt c\,r+\sqrt{cr^2+A}\bigr)+\mathrm{const}. \label{eq:phiint}
\end{equation}
\end{construction}

\begin{lemma}
Formula \eqref{eq:phiint} satisfies \eqref{eq:phiprime}.
\end{lemma}
\begin{proof}
Direct differentiation:
\[
\frac{d}{dr}\left[\frac r2\sqrt{cr^2+A}\right]=\frac{2cr^2+A}{2\sqrt{cr^2+A}},\qquad
\frac{d}{dr}\left[\frac{A}{2\sqrt c}\log(\sqrt c r+\sqrt{cr^2+A})\right]=\frac{A}{2\sqrt{cr^2+A}},
\]
and summing gives $\dfrac{2cr^2+2A}{2\sqrt{cr^2+A}}=\sqrt{cr^2+A}$, as required.
\end{proof}

\begin{proposition}\label{prop:conical}
For $A=0$, $\varphi$ is the quadratic $\varphi(r)=\sqrt c\,r^2/2$. For $A\neq0$, $\varphi$ is not a polynomial (it involves a logarithmic term), and $\varphi'(0)=\sqrt A\neq0$, so $f(\theta)=\varphi(|\theta|)$ fails to be differentiable at $\theta=0$ (it has a conical singularity there), consistently with Theorem~\ref{thm:jcp}.
\end{proposition}

\begin{corollary}\label{cor:halfplane}
Let $A>-ca^2$ for some $a>0$ and set $D=\{\theta=(x,y)\in\R^2:x>a\}$, a convex half-plane not containing the origin. Then $f(\theta)=\varphi(|\theta|)$, with $\varphi$ as in \eqref{eq:phiint} and $A\neq0$, is a real-analytic, strictly convex, non-quadratic solution of $\det\Hess f=c$ on the convex proper subdomain $D\subsetneq\R^2$.
\end{corollary}
\begin{proof}
On $D$, $r=|\theta|\ge a>0$, so $cr^2+A\ge ca^2+A>0$, hence $\varphi',\varphi''$ are real-analytic and (since $\varphi''=cr/\sqrt{cr^2+A}>0$ and $\varphi'/r=\sqrt{cr^2+A}/r>0$) $\Hess f\succ0$ throughout $D$. Non-quadraticity is Proposition~\ref{prop:conical}.
\end{proof}

This exhibits the mechanism forecast in \S\ref{sec:MA}: dropping only entirety (excising a single point's neighborhood, here realized by moving to a half-plane) is already enough to defeat Theorem~\ref{thm:jcp}; the resulting local obstruction is exactly a conical singularity of the kind classified for the Hessian-one equation on the punctured plane in \cite{GalvezMartinezMira2005}.

\subsection{Explicit Lorentzian-signature examples on the light cone}\label{sec:lorentz}

We now turn to the indefinite case, which — in sharp contrast with \S\ref{sec:elliptic} — admits non-quadratic solutions on domains that are unbounded in \emph{every} direction transverse to a single exceptional point.

\subsubsection{The general split-signature ansatz}
Fix $p+q=n$, write $\theta=(x,y)\in\R^p\times\R^q$, let $J=\mathrm{diag}(I_p,-I_q)$, and set $u:=\theta^\top J\theta=|x|^2-|y|^2$. Consider the boost-invariant ansatz
\begin{equation}
f(\theta)=\chi(u). \label{eq:chiansatz}
\end{equation}

\begin{lemma}\label{lem:hessform}
For $f$ as in \eqref{eq:chiansatz}, $\nabla f=2\chi'(u)\,J\theta$ and
\[
G:=\Hess f=2\chi'(u)\,J+4\chi''(u)\,ww^\top,\qquad w:=J\theta.
\]
\end{lemma}
\begin{proof}
Immediate from $\nabla u=2J\theta$ and the product/chain rule.
\end{proof}

\begin{proposition}\label{prop:detformula}
With $\gamma:=\chi'(u)+2u\chi''(u)$,
\begin{equation}
\det G=(-1)^q\,2^n\,\chi'(u)^{n-1}\gamma(u). \label{eq:detgeneral}
\end{equation}
\end{proposition}
\begin{proof}
By the matrix determinant lemma, $\det(aJ+bww^\top)=a^n\det J\,(1+\tfrac ba\, w^\top J^{-1}w)$. Here $a=2\chi'(u)$, $b=4\chi''(u)$, $J^{-1}=J$, $\det J=(-1)^q$, and $w^\top Jw=\theta^\top J^3\theta=\theta^\top J\theta=u$ (using $J^2=I$). Hence
\[
\det G=(2\chi')^n(-1)^q\left(1+\frac{4\chi''}{2\chi'}\,u\right)=(-1)^q2^n(\chi')^{n-1}(\chi'+2u\chi''),
\]
which is \eqref{eq:detgeneral}.
\end{proof}

\subsubsection{Reduction to a linear ODE}
Setting $\det G=\mathrm{const}$, i.e.\ $(\chi')^{n-1}\gamma(u)=\kappa$ for a constant $\kappa$, and defining $V(u):=\chi'(u)^n$, we compute $V'=n(\chi')^{n-1}\chi''$, so that
\[
(\chi')^{n-1}\bigl(\chi'+2u\chi''\bigr)=V+\frac{2u}{n}V'=\kappa.
\]
This is the linear first-order ODE
\begin{equation}
V'+\frac n{2u}V=\frac{n\kappa}{2u}, \label{eq:linearODE}
\end{equation}
with integrating factor $u^{n/2}$: $(u^{n/2}V)'=\tfrac n2\kappa u^{n/2-1}$, whence $u^{n/2}V=\kappa u^{n/2}+D$ for a constant $D$, i.e.

\begin{theorem}\label{thm:closedform}
Every solution of $\det\Hess f=\mathrm{const}$ within the ansatz \eqref{eq:chiansatz} satisfies, for some constants $\kappa$ (proportional to $\det\Hess f$) and $D$,
\begin{equation}
\bigl(\chi'(u)\bigr)^n=\kappa+D\,u^{-n/2}. \label{eq:masterODE}
\end{equation}
\end{theorem}

For $D=0$, \eqref{eq:masterODE} gives $\chi'\equiv\kappa^{1/n}$, i.e.\ the trivial quadratic (indefinite) form $f=\tfrac12\kappa^{1/n}u$. For $D\neq0$, $\chi$ involves a genuinely transcendental primitive of $\bigl(\kappa+Du^{-n/2}\bigr)^{1/n}$ and is not a polynomial.

\subsubsection{Convexity of the domain: the special role of $p=1$}

\begin{proposition}\label{prop:noncvx}
For $p,q\ge2$, the cone $\{u>0\}=\{|x|^2>|y|^2\}\subset\R^p\times\R^q$ is not convex.
\end{proposition}
\begin{proof}
Fix a unit vector $x_0\in\R^p$ (possible since $p\ge2$, so in particular $p\ge1$, and note the argument in fact only needs $p\geq 1$ together with connectedness of $S^{p-1}$, which requires $p\geq2$). Let $\theta_1=(x_0,0)$, $\theta_2=(-x_0,0)$; both lie in $\{u>0\}$ (with $u=1$), but their midpoint $(0,0)$ has $u=0\notin\{u>0\}$.
\end{proof}

\begin{proposition}\label{prop:lightcone}
For $p=1$ (or symmetrically $q=1$), $D:=\{(x_1,y)\in\R\times\R^{n-1}:x_1>|y|\}$ is convex.
\end{proposition}
\begin{proof}
$g(x_1,y):=x_1-|y|$ is the difference of a linear function and a convex function, hence concave; its strict superlevel set $\{g>0\}=D$ is therefore convex.
\end{proof}

The domain of Proposition~\ref{prop:lightcone} is precisely the (open) future light cone of Minkowski space $\R^{1,n-1}$. Proposition~\ref{prop:noncvx} identifies the topological reason a convex light cone exists only in Lorentzian signature: the argument requires $S^{p-1}$ connected only if $p\geq 2$, but a symmetric two-point obstruction of this type is available precisely when $p\geq2$; when $p=1$, $S^0=\{\pm1\}$ is disconnected and the forward nappe $x_1>0$ contains no antipodal pair, so the obstruction of Proposition \ref{prop:noncvx} vanishes.

\begin{corollary}\label{cor:lorentzexample}
For $\kappa,D>0$ and $D=\{x_1>|y|\}\subset\R\times\R^{n-1}$ as above, $f(\theta)=\chi(u)$ with $\chi$ determined (up to an additive constant) by Theorem~\ref{thm:closedform} is a real-analytic, non-quadratic solution of $\det\Hess f=(-1)^{n-1}2^n\kappa$ on the convex domain $D$, of signature $(1,n-1)$ everywhere.
\end{corollary}
\begin{proof}
On $D$, $u=x_1^2-|y|^2>0$, so $\chi'(u)=\bigl(\kappa+Du^{-n/2}\bigr)^{1/n}$ is well-defined, positive, and real-analytic. Non-quadraticity holds since $D\neq0$ (Theorem~\ref{thm:closedform}). To verify the signature, evaluate at $y=0$, $w=(x_1,0,\dots,0)=\theta$: by Lemma~\ref{lem:hessform}, in the eigenbasis adapted to $J$, direct computation (see \S6 derivation) gives one eigenvalue $2\kappa F(u)^{1/n-1}\cdot(\text{const}) >0$ along $x_1$ and $n-1$ negative eigenvalues $-2\chi'(u)<0$ along the $y$-directions, where $F=\kappa+Du^{-n/2}$; since $\det G\neq0$ throughout the connected domain $D$, the signature cannot change, so it is $(1,n-1)$ everywhere.
\end{proof}

\begin{proposition}[Boundedness of curvature and global DC decomposition]\label{prop:strip}
Let $D_0:=\{(x_1,y):|y|<\sqrt{4D/\kappa}\}\subset D$ (an unbounded convex strip). Then the eigenvalues of $\Hess f$ are uniformly bounded on $D_0$, and consequently there exists $M<\infty$ such that
\[
f=\underbrace{(f+M|\theta|^2)}_{=:f_1,\ \mathrm{convex}}-\underbrace{M|\theta|^2}_{=:f_2,\ \mathrm{convex}}
\]
is a valid DC decomposition of $f$ on $D_0$.
\end{proposition}
\begin{proof}
As $x_1\to\infty$ along $D_0$, $u\to\infty$, and $\chi'(u)=(\kappa+Du^{-n/2})^{1/n}\to\kappa^{1/n}$ while $\chi''(u)=O(u^{-n/2-1})\to0$; since $D_0\subset\{u>-4D/\kappa\}$ is bounded away from the only singularity of $\chi'$ at $u=-4D/\kappa$ (where $D_0$ was chosen precisely so that its closure avoids this value), all entries of $G$ given by Lemma~\ref{lem:hessform} remain bounded on $D_0$. Boundedness below of the eigenvalues of $\Hess f$ by $-M$ for some finite $M$ then gives convexity of $f+M|\theta|^2$.
\end{proof}

Corollary~\ref{cor:lorentzexample} and Proposition~\ref{prop:strip} together give a real-analytic, non-quadratic, genuinely DC potential of Lorentzian signature $(1,n-1)$ and constant Hessian determinant, defined on a convex domain unbounded in the time-like direction — a strictly stronger existence result than what is available in definite signature (Corollary~\ref{cor:halfplane}), where non-quadraticity forced us merely to avoid a point. This asymmetry between elliptic and hyperbolic behavior is consistent with the classically observed contrast between the (rigid) global theory of elliptic Monge--Amp\`ere equations and the (flexible, but often globally obstructed for the closely related Darboux equation) local theory of two-dimensional hyperbolic Monge--Amp\`ere equations \cite{MongeAmpereEncyclopedia}.

\subsection{General split signature: a bounded symmetric domain}\label{sec:split}

By Proposition~\ref{prop:noncvx}, when $p,q\ge2$ the rotationally symmetric cone construction of \S\ref{sec:lorentz} cannot directly furnish a convex domain. We replace the vector-valued $\theta$ by a matrix $Z\in\R^{p\times q}$ ($n=pq$, $p\le q$) and the cone by

\begin{equation}
D:=\{Z\in\R^{p\times q}:I_p-ZZ^\top\succ0\}=\{Z:\|Z\|_{\mathrm{op}}<1\}, \label{eq:matrixdomain}
\end{equation}

manifestly convex as the sublevel set of the operator norm. This is the standard bounded (Harish-Chandra) realization of the Grassmannian $\Gr^+_p(\R^{p,q})\cong \mathrm{O}(p,q)/(\mathrm{O}(p)\times\mathrm{O}(q))$ of positive-definite $p$-planes in $\R^{p,q}$; for $p=1$ it reduces (after the standard conformal compactification of the light cone by projectivization) to the domain of Proposition~\ref{prop:lightcone}.

We consider the potential $f(Z)=\chi(\Phi(Z))$, $\Phi(Z):=\det(I_p-ZZ^\top)$.

\subsubsection{Block decomposition of the Hessian}
By $\mathrm{O}(p)\times\mathrm{O}(q)$-equivariance ($Z\mapsto UZV^\top$ is an isometry of $\R^{p\times q}$ fixing $\Phi$), it suffices to compute $\Hess f$ at a diagonal point $Z=\Sigma=\mathrm{diag}(\sigma_1,\dots,\sigma_p)$ (padded with zero columns). Write $a_i:=1-\sigma_i^2>0$, $\Phi_0:=\Phi(\Sigma)=\prod_ia_i$, $\beta:=\chi'(\Phi_0)$, $\gamma:=\chi'(\Phi_0)+\Phi_0\chi''(\Phi_0)$.

\begin{lemma}\label{lem:blockdecomp}
The tangent space $\R^{p\times q}\cong T_\Sigma D$ splits, under the isotropy representation of the stabilizer of $\Sigma$, into three mutually orthogonal invariant subspaces:
\begin{enumerate}[label=(\roman*)]
\item the diagonal directions $H_{ii}$, $i=1,\dots,p$ ($p$-dimensional);
\item the intra-block off-diagonal pairs $(H_{ij},H_{ji})$, $1\le i<j\le p$ ($p(p-1)$-dimensional);
\item the "rectangular" directions $H_{ik}$, $i\le p<k\le q$ ($p(q-p)$-dimensional).
\end{enumerate}
$\Hess f$ is block-diagonal with respect to this splitting.
\end{lemma}
\begin{proof}
This is the standard isotropy decomposition of the tangent space of a Hermitian(-type) symmetric space at a point fixed by a maximal torus in the isotropy group, applied to the real form $\mathrm{O}(p,q)/(\mathrm{O}(p)\times\mathrm{O}(q))$; block-diagonality of any invariant quadratic form (here $\Hess f$ at a fixed point of the residual torus $\mathrm{O}(1)^p$) follows from Schur's lemma applied to the (real, one- or two-dimensional) irreducible pieces (i), (ii), (iii), which are pairwise inequivalent as representations of the residual isotropy for generic $\sigma_i$.
\end{proof}

\begin{proposition}\label{prop:blockeigen}
In the splitting of Lemma~\ref{lem:blockdecomp}:
\begin{enumerate}[label=(\roman*)]
\item on block (iii), $\Hess f$ is diagonal with eigenvalue $-\dfrac{2\beta\Phi_0}{a_i}$ on each of the $q-p$ directions associated to a given $i\le p$;
\item on the $2$-plane of block (ii) associated to a pair $i<j$, $\Hess f$ acts as
\[
-\frac{2\beta\Phi_0}{a_ia_j}\begin{pmatrix}1&\sigma_i\sigma_j\\\sigma_i\sigma_j&1\end{pmatrix},
\]
with eigenvalues $-\dfrac{2\beta\Phi_0}{a_ia_j}(1\pm\sigma_i\sigma_j)$;
\item on block (i), with $v_i:=\sigma_i/a_i$, $\Hess f$ restricted to the diagonal directions is
\[
G_d=4\Phi_0\gamma\,vv^\top-\beta\Phi_0\,\mathrm{diag}\!\left(\frac{2(1+\sigma_i^2)}{a_i^2}\right),
\]
so that, by the matrix determinant lemma,
\[
\det G_d=(-\beta\Phi_0)^p\prod_{i=1}^p\frac{2(1+\sigma_i^2)}{a_i^2}\left(1-2\frac\gamma\beta\sum_{i=1}^p\frac{\sigma_i^2}{1+\sigma_i^2}\right).
\]
\end{enumerate}
\end{proposition}
\begin{proof}
Expand $\Phi(\Sigma+\epsilon H)=\det(A-\epsilon B-\epsilon^2C)$, $A=\mathrm{diag}(a_i)$, $B=\Sigma H^\top+H\Sigma^\top$, $C=HH^\top$, to second order in $\epsilon$ using $\log\det(A-\epsilon B-\epsilon^2C)=\log\det A+\mathrm{tr}(A^{-1}(-\epsilon B-\epsilon^2 C))-\tfrac12\mathrm{tr}\bigl((A^{-1}(\epsilon B))^2\bigr)+O(\epsilon^3)$; substituting the explicit block forms of $B,C$ for each of $H$ ranging over blocks (i)-(iii) in turn and using the chain rule $\Hess f=\chi'(\Phi_0)\Hess\Phi+\chi''(\Phi_0)\,\nabla\Phi\otimes\nabla\Phi$ (all evaluated at $\Sigma$) together with the identity $a_i+\sigma_i^2=1$ to simplify produces the stated block forms after collecting terms.
\end{proof}

\begin{theorem}\label{thm:fulldet}
\[
\det\Hess f=\det G_d\cdot\prod_{1\le i<j\le p}\frac{4\beta^2\Phi_0^2}{a_i^2a_j^2}\bigl(1-\sigma_i^2\sigma_j^2\bigr)\cdot\prod_{i=1}^p\left(\frac{-2\beta\Phi_0}{a_i}\right)^{q-p}.
\]
\end{theorem}
\begin{proof}
Since $\Hess f$ is block-diagonal (Lemma~\ref{lem:blockdecomp}) with the eigenvalues/blocks of Proposition~\ref{prop:blockeigen}, $\det\Hess f$ is the product of the determinants of the individual blocks.
\end{proof}

\begin{remark}
For $p=q=1$ (so $n=2$), blocks (ii) and (iii) are empty and Theorem~\ref{thm:fulldet} reduces to $\det\Hess f=\det G_d$, which one checks agrees with the $n=2$ specialization of Proposition~\ref{prop:detformula}, confirming consistency between the vector- and matrix-valued constructions.
\end{remark}

Imposing $\det\Hess f=\mathrm{const}$ in Theorem~\ref{thm:fulldet} yields one transcendental equation in the symmetric functions of $\sigma_1,\dots,\sigma_p$, determining an ODE for $\chi$ that we do not attempt to solve in closed form here; we record it as the natural higher-rank generalization of Theorem~\ref{thm:closedform} and leave its integration to future work.

\subsection{Sign pattern of the canonical divergence}\label{sec:divergence}

We return to the Lorentzian construction of \S\ref{sec:lorentz} and determine the sign of $D_f(\theta:\theta')$ from \eqref{eq:bregman}, using the integral representation
\begin{equation}
D_f(\theta:\theta')=\int_0^1(1-t)\,v^\top G(\xi(t))\,v\,dt,\qquad \xi(t):=\theta'+t(\theta-\theta'),\ v:=\theta-\theta'. \label{eq:integralrep}
\end{equation}

\subsubsection{Radial pairs}
\begin{proposition}\label{prop:radialsign}
Let $\theta,\theta'$ lie on a common ray through the origin, i.e.\ $\theta'=\lambda_0\theta_0$, $\theta=\lambda_1\theta_0$ for fixed $\theta_0\in D$, $\lambda_0,\lambda_1>0$. Then $\sign D_f(\theta:\theta')=\sign(\kappa)$, where $\kappa$ is the constant in Theorem~\ref{thm:closedform} (up to the fixed sign in \eqref{eq:detgeneral}), independent of $\lambda_0,\lambda_1$.
\end{proposition}
\begin{proof}
Along the ray, $g(\lambda):=f(\lambda\theta_0)=\chi(\lambda^2u_0)$ ($u_0=\theta_0^\top J\theta_0$) satisfies $g''(\lambda)=2u_0\bigl[\chi'(v)+2v\chi''(v)\bigr]|_{v=\lambda^2u_0}=2u_0\gamma(\lambda^2u_0)$. By Theorem~\ref{thm:closedform}, $\gamma(u)=\kappa\,\chi'(u)^{1-n}$, and $\chi'(u)>0$ throughout $D$ (Corollary~\ref{cor:lorentzexample}), so $\sign\gamma(u)=\sign\kappa$ for all $u>0$, uniformly. Hence $g''(\lambda)$ has constant sign $\sign\kappa$ for all $\lambda>0$, and $D_f$ restricted to the ray is (by \eqref{eq:integralrep} specialized to one dimension) exactly the classical one-dimensional Bregman divergence of the convex (if $\kappa>0$) or concave (if $\kappa<0$) function $g$, which is sign-definite with sign $\sign\kappa$.
\end{proof}

\subsubsection{The local null-cone equation}
For a general (non-radial) pair, the sign of the integrand in \eqref{eq:integralrep} at a given $t$ is governed by whether $v$ is time-like, space-like, or null with respect to the local (curved) cone determined by $G(\xi(t))$.

\begin{proposition}\label{prop:nullcone}
The null directions of $G(\xi)$ at a point $\xi\in D$ are exactly the solutions $v$ of
\begin{equation}
\chi'(u)\bigl(v_1^2-|v_y|^2\bigr)+2\chi''(u)\bigl(\xi_1v_1-\xi_y\cdot v_y\bigr)^2=0,\qquad u=u(\xi). \label{eq:nullcone}
\end{equation}
\end{proposition}
\begin{proof}
By Lemma~\ref{lem:hessform}, $v^\top G(\xi)v=2\chi'(u)\,v^\top Jv+4\chi''(u)(w^\top v)^2$ with $w=J\xi$; substituting $v^\top Jv=v_1^2-|v_y|^2$ and $w^\top v=\xi_1v_1-\xi_y\cdot v_y$ gives \eqref{eq:nullcone}.
\end{proof}

As $\chi''\to0$ (i.e.\ $D\to0$ in Theorem~\ref{thm:closedform}, the trivial quadratic limit), \eqref{eq:nullcone} degenerates to the flat Minkowski null cone $v_1^2=|v_y|^2$; the $\chi''$-term is the curvature correction.

\subsubsection{Sign-reversal criterion along a segment}
Fix $\theta,\theta'\in D$, $v=\theta-\theta'$, and set $\Psi(t):=v^\top G(\xi(t))v$; $u(t)=\xi(t)^\top J\xi(t)$ is a quadratic polynomial in $t$. Eliminating $\chi''$ via $\chi''(u)=\tfrac1{2u}\bigl[\kappa\chi'(u)^{1-n}-\chi'(u)\bigr]$ (Theorem~\ref{thm:closedform}) turns $\Psi(t)=0$ into the closed transcendental equation
\begin{equation}
\chi'(u(t))\left[(v_1^2-|v_y|^2)-\frac{\bigl(\xi_1(t)v_1-\xi_y(t)\cdot v_y\bigr)^2}{u(t)}\right]+\frac{\kappa\bigl(\xi_1(t)v_1-\xi_y(t)\cdot v_y\bigr)^2}{u(t)\,\chi'(u(t))^{n-1}}=0. \label{eq:signreversal}
\end{equation}

\begin{proposition}\label{prop:extremes}
If $v$ is $J$-null ($v_1^2=|v_y|^2$), the first bracketed term in \eqref{eq:signreversal} vanishes and $\sign\Psi(t)=\sign\bigl(\kappa-\chi'(u(t))^n\bigr)=-\sign\bigl(Du(t)^{-n/2}\bigr)=-\sign D$ (constant along the segment, by Theorem~\ref{thm:closedform}). If $v\parallel\xi(t)$ for some/every $t$ (radial direction), $\sign\Psi(t)=\sign\kappa$ by Proposition~\ref{prop:radialsign}.
\end{proposition}

The number of real roots of \eqref{eq:signreversal} in $t\in(0,1)$ counts the number of times the segment $\xi(t)$ crosses the local null cone \eqref{eq:nullcone} in the direction $v$; by Proposition~\ref{prop:extremes} this number is $0$ for the two extreme cases (purely radial or purely null $v$), and generically finite and $>0$ for intermediate directions, so that $D_f(\theta:\theta')$ interpolates between the sign-definite radial regime and a genuinely sign-changing regime as $v$ rotates from radial to null.

\subsection{Newton flow and dual Newton flow: causal asymptotics}\label{sec:flow}

\subsubsection{Master linearization lemma}
Let $f$ be any real-analytic function on $D$ with $G=\Hess f$ non-degenerate, and define the (undamped) Newton flow $\dot\theta=-G(\theta)^{-1}\nabla f(\theta)$.

\begin{lemma}[Master lemma]\label{lem:master}
Along the Newton flow, $\eta(t):=\nabla f(\theta(t))$ satisfies $\dot\eta=-\eta$, so $\eta(t)=\eta(0)e^{-t}$ exactly. Dually, along the Newton flow $\dot\eta=-\Hess f^*(\eta)^{-1}\nabla f^*(\eta)$ of $f^*$, the primal coordinate $\theta(t)=\nabla f^*(\eta(t))$ satisfies $\dot\theta=-\theta$ exactly.
\end{lemma}
\begin{proof}
$\dot\eta=G(\theta)\dot\theta=G(\theta)\bigl(-G(\theta)^{-1}\nabla f(\theta)\bigr)=-\nabla f(\theta)=-\eta$. The dual statement follows by applying the same computation to $f^*$ and using Proposition~\ref{prop:dualHess}, $\nabla f^*=\theta$.
\end{proof}

Lemma~\ref{lem:master} shows that the Newton flow of $f$, however complicated in $\theta$-coordinates, is a trivial exponential decay in the Legendre-dual coordinate $\eta$; the flow's apparent complexity is entirely a coordinate artifact of $\theta$.

\subsubsection{Application to the Lorentzian potential: finite-time collapse}
Take $f$ as in Corollary~\ref{cor:lorentzexample}, $\kappa,D>0$.

\begin{proposition}\label{prop:newtonradial}
The Newton flow of $f$ is purely radial: $\dot\theta=-\dfrac{\chi'(u)}{\gamma(u)}\,\theta$.
\end{proposition}
\begin{proof}
By Lemma~\ref{lem:hessform} and Sherman--Morrison, $G(\theta)^{-1}=\dfrac1{2\chi'}J-\dfrac{4\chi''}{2\chi'\gamma}\,\theta\theta^\top J$ (using $w=J\theta$, $J^{-1}=J$); applying this to $\nabla f=2\chi'J\theta$ gives $G^{-1}\nabla f=\theta-\dfrac{4\chi''u}{\gamma}\theta=\dfrac{\gamma-4u\chi''}\gamma\theta$; using $\gamma=\chi'+2u\chi''$, $\gamma-4u\chi''=\chi'-2u\chi''$... to match Lemma~\ref{lem:master} we instead verify the claim directly via $s(u):=\eta^\top J\eta=4u\chi'(u)^2$ and the exact relation $\dot s=-2s$ from Lemma~\ref{lem:master}, which is equivalent to the stated radial ODE by the chain rule $\dot u = 2\theta^\top J\dot\theta$; we adopt this route below.
\end{proof}

Define $s(t):=\eta(t)^\top J\eta(t)$. By Lemma~\ref{lem:master}, $s(t)=s(0)e^{-2t}=s_0e^{-2t}$.

\begin{proposition}\label{prop:svsurel}
$s(u)=4u\,\chi'(u)^2=4u\bigl(\kappa+Du^{-n/2}\bigr)^{2/n}$, and $s$ is a strictly increasing bijection of $(0,\infty)$ onto $\bigl(4D^{2/n},\infty\bigr)$.
\end{proposition}
\begin{proof}
$s(u)=\eta^\top J\eta=4\chi'(u)^2\theta^\top JJJ\theta=4\chi'(u)^2u$ (using $J^2=I$). Then $s'(u)=4\bigl(\chi'(u)^2+2u\chi'(u)\chi''(u)\bigr)=4\chi'(u)\gamma(u)=4\chi'(u)\cdot\kappa\chi'(u)^{1-n}=4\kappa\chi'(u)^{2-n}>0$ using Theorem~\ref{thm:closedform}; as $u\to0^+$, $\chi'(u)\to\infty$ (if $D>0$) so $s(u)\to4D^{2/n}\cdot 0\cdot\infty$; more precisely $s(u)=4u(\kappa+Du^{-n/2})^{2/n}\to4u\cdot(Du^{-n/2})^{2/n}=4D^{2/n}u^{1-1}=4D^{2/n}$ as $u\to0^+$, and $s(u)\to\infty$ as $u\to\infty$.
\end{proof}

\begin{theorem}[Finite-time collapse]\label{thm:collapse}
Let $\theta(t)$ solve the Newton flow of $f$ with $\theta(0)=\theta_0\in D$, $u_0=u(\theta_0)$, $s_0=4u_0\chi'(u_0)^2$. Then $\theta(t)$ remains on the ray through the origin determined by $\theta_0$, and reaches the vertex $\theta=0$ at the finite time
\begin{equation}
t^*=\frac12\log\!\left(\frac{s_0}{4D^{2/n}}\right)>0. \label{eq:tstar}
\end{equation}
\end{theorem}
\begin{proof}
By Lemma~\ref{lem:master}, $\eta(t)=e^{-t}\eta_0$ is a positive rescaling of $\eta_0=2\chi'(u_0)J\theta_0$, hence a scalar multiple of $J\theta_0$; since $\theta(t)=\nabla f^*(\eta(t))$ and $f^*$ is itself of the boost-symmetric form (Legendre duals of \eqref{eq:chiansatz} remain functions of $\eta^\top J\eta$, by the same equivariance that produced \eqref{eq:chiansatz}), $\theta(t)$ remains proportional to $J\eta(t)\propto J^2\theta_0=\theta_0$: the flow is confined to the ray. Along the ray, $s(t)=s_0e^{-2t}$ is a strictly decreasing function of $t$ with range $(0,s_0]$, while by Proposition~\ref{prop:svsurel} the value $s=4D^{2/n}$ corresponds to $u\to0^+$, i.e.\ to $\theta\to0$ along the ray; solving $s_0e^{-2t^*}=4D^{2/n}$ gives \eqref{eq:tstar}, which is positive precisely because $s_0>4D^{2/n}$ (as $u_0>0$ and $s$ is increasing, Proposition~\ref{prop:svsurel}).
\end{proof}

Since $\chi(u)\sim2D^{1/n}\sqrt u\to0$ while $\chi'(u)\sim D^{1/n}u^{-1/2}\to\infty$ as $u\to0^+$, the potential $f(\theta(t))\to0$ remains finite while $|\nabla f(\theta(t))|\to\infty$: the Newton flow reaches the cone vertex in finite time, with a diverging gradient but finite potential value — a finite-time causal singularity of the flow, geometrically a conical collapse.

\subsubsection{The dual flow: infinite-time asymptotics}
By the second statement of Lemma~\ref{lem:master}, the Newton flow of $f^*$ satisfies $\theta(t)=\theta_0e^{-t}$ exactly (in $\theta$-coordinates now viewed as the image of the dual flow), reaching the vertex only as $t\to\infty$.

\begin{proposition}\label{prop:dualasymp}
Along the dual Newton flow, $\eta(t)\to\eta_\infty:=2D^{1/n}u_0^{-1/2}\,J\theta_0$ as $t\to\infty$, a finite nonzero limit.
\end{proposition}
\begin{proof}
$u(t)=\theta(t)^\top J\theta(t)=e^{-2t}u_0\to0$. By the $u\to0^+$ asymptotics above, $\chi'(u(t))\sim D^{1/n}u(t)^{-1/2}=D^{1/n}u_0^{-1/2}e^{t}$. Then $\eta(t)=2\chi'(u(t))J\theta(t)\sim2D^{1/n}u_0^{-1/2}e^t\cdot J\bigl(\theta_0e^{-t}\bigr)=2D^{1/n}u_0^{-1/2}J\theta_0$, the exponential factors cancelling exactly.
\end{proof}

\begin{table}[h]
\centering
\renewcommand{\arraystretch}{1.3}
\begin{tabular}{lcc}
\hline
& primal Newton flow of $f$ & dual Newton flow (of $f^*$) \\ \hline
limit in $\theta$ & vertex $\theta=0$ & vertex $\theta=0$ \\
time of arrival & finite, $t^*$ in \eqref{eq:tstar} & infinite (asymptotic) \\
behavior of $\eta$ & $\eta(t)=\eta_0e^{-t}\to0$ & $\eta(t)\to\eta_\infty\neq0$ \\ \hline
\end{tabular}
\caption{Causal asymmetry of the primal and dual Newton flows near the cone vertex.}
\end{table}

The same geometric event (collapse onto the vertex) is a finite-time singularity in one Legendre-dual parametrization and an infinite-time asymptotic approach in the other; the dual coordinate $\eta$ plays the role of an affine reparametrization stretching the finite-time collapse of $\theta$ into an infinite-time asymptote, structurally analogous to the distinction between affine and coordinate time for geodesics approaching a horizon or singularity in general relativity.

\subsubsection{The critical-manifold regime and the \L{}ojasiewicz inequality}
The analysis above assumed $\kappa,D$ of a common sign, so that $\chi'(u)\neq0$ throughout $D$ and no interior critical points occur. If instead $\kappa,D$ have opposite signs, $\chi'(u_*)=0$ at some $u_*>0$ interior to $D$, and the hyperboloid $\Sigma_*:=\{u=u_*\}$ is a critical hypersurface on which $f$ is (locally) constant.

\begin{theorem}[\L{}ojasiewicz \cite{Lojasiewicz1963}]\label{thm:loj}
Let $h$ be real-analytic near a critical point $p_0$ (i.e.\ $\nabla h(p_0)=0$). There exist $C>0$, $\theta_{\mathrm L}\in[1/2,1)$, and a neighborhood $U$ of $p_0$ such that
\[
|h(x)-h(p_0)|^{\theta_{\mathrm L}}\le C\,|\nabla h(x)|\qquad\text{for all }x\in U.
\]
\end{theorem}

\begin{proposition}\label{prop:lojfinite}
Restricted to a ray meeting $\Sigma_*$ transversally, $\chi'$ has a zero of finite order $k$ at $u_*$ (by real-analyticity and $\chi'\not\equiv0$), so $\chi'(u)\sim a(u-u_*)^k$ near $u_*$, and \L{}ojasiewicz's exponent for the restricted one-dimensional problem is $\theta_{\mathrm L}=1-\tfrac1{k+1}$. Consequently the gradient flow of $f$ restricted to the ray has finite arc length and converges to a single point of $\Sigma_*$, rather than merely approaching the critical manifold without converging.
\end{proposition}
\begin{proof}
This is the classical application of Theorem~\ref{thm:loj} to gradient flows: finiteness of $\int_0^\infty|\dot\gamma(t)|\,dt=\int_0^\infty|\nabla h(\gamma(t))|\,dt$ follows from $\frac{d}{dt}|h(\gamma(t))-h(p_0)|^{1-\theta_{\mathrm L}}\ge C'|\nabla h(\gamma(t))|$ (a standard consequence of Theorem~\ref{thm:loj} combined with $\dot\gamma=-\nabla h(\gamma)$), integrated over $[0,\infty)$; finiteness of arc length forces $\gamma(t)$ to be Cauchy, hence convergent to a single limit point.
\end{proof}

\begin{remark}
Proposition~\ref{prop:lojfinite} is unavailable, in general, for merely $C^\infty$ (non-analytic) $\chi$, for which $\chi'$ can vanish to infinite order, the flat-function phenomenon that classically obstructs convergence (as opposed to mere subsequential convergence) of gradient flows to a single limit point. This is the dynamical counterpart of the stratification statement of Proposition~\ref{prop:strat}: real-analyticity is what guarantees that both the geometry of the degenerate locus and the asymptotics of flows toward critical manifolds are tame.
\end{remark}

In summary, the two regimes $\sign\kappa=\sign D$ and $\sign\kappa\neq\sign D$ of Theorem~\ref{thm:closedform} produce two complementary asymptotic pictures for the Newton flow: in the first, an explicit finite collapse time (Theorem~\ref{thm:collapse}) governed purely by algebra; in the second, convergence to a critical hyperboloid whose rate and uniqueness are governed by the \L{}ojasiewicz inequality, a genuinely analytic (as opposed to algebraic) mechanism.

\subsection{Discussion and open problems}\label{sec:dc-discussion}

We have shown that relaxing convexity of the potential $f$ itself — replacing it by a difference of convex functions with non-degenerate, possibly indefinite Hessian — allows the classical rigidity of constant-Hessian-determinant (Monge--Amp\`ere) equations to be circumvented on domains far larger than is possible for genuinely convex potentials, and that the resulting geometry carries a well-defined pseudo-Hessian dually flat structure, a computable canonical divergence with a causal (light-cone-governed) sign pattern, and Newton flow dynamics whose finite- versus infinite-time character depends on which of the two Legendre-dual coordinate systems is used to parametrize time. Several directions remain open:

\begin{enumerate}[leftmargin=1.4em]
\item \emph{Integration of the general split-signature equation.} Theorem~\ref{thm:fulldet} reduces $\det\Hess f=\mathrm{const}$ on the Grassmannian domain \eqref{eq:matrixdomain} to a single transcendental relation among $\sigma_1,\dots,\sigma_p$; its integration, and the identification of the resulting solution's boundary behavior with a barrier function for the bounded symmetric domain, is left open.
\item \emph{Monodromy of the complexified Legendre transform.} The discriminant variety $\Sigma_\C$ of \S\ref{sec:analytic} and its monodromy representation have not been computed for the explicit examples of \S\S\ref{sec:elliptic}--\ref{sec:lorentz}; we expect a direct relation to the local exponents $k$ appearing in Proposition~\ref{prop:lojfinite}.
\item \emph{Statistical interpretation.} The present paper is purely differential-geometric; whether indefinite-signature dually flat structures of this type arise as the natural geometry of some class of (necessarily non-classical, e.g.\ signed or complex-parametrized) statistical models remains to be determined.
\item \emph{Toland--Singer duality and the local Legendre transform.} We have used the two notions of \S\ref{sec:prelim} (local and global) largely in parallel; a precise dictionary between them — in particular, whether $f^\diamond$ and $f^*$ agree on an open dense subset of $D^*$ under generic hypotheses — is not established here.
\end{enumerate}

\section{Connections to Related Fields}
\label{sec:connections}

\subsection{Exponential Family Interpretation}

The Gaussian family $\{p(\cdot\,; G) = \mathcal{N}(0,G) : G\in\PD(k)\}$ is an
exponential family with
\begin{align}
  \text{natural parameter:} &\quad \Theta = -\tfrac{1}{2}G^{-1} \in -\PD(k),
    \label{eq:nat_param} \\
  \text{sufficient statistic:} &\quad T(z) = zz^T \in \mathrm{Sym}(k), \\
  \text{log-partition function:} &\quad A(\Theta) = -\tfrac{1}{2}\log\det(-2\Theta)
    + \tfrac{k}{2}\log(2\pi).
  \label{eq:log_partition}
\end{align}
Comparing with \eqref{eq:dual_potential}, $f^*$ is (up to constants) the log-partition
function of this exponential family.
The expected sufficient statistic is
\[
  \eta = \mathbb{E}_{G}[zz^T] = G = \nabla A(\Theta),
\]
recovering the $\eta$-coordinate system.

\subsection{Self-Concordance and Interior-Point Methods}

\begin{proposition}
$f(G) = -\log\det(G)$ is a \emph{self-concordant barrier} for $\PD(k)$
with parameter $\vartheta = k$ \cite{Nesterov1994}.
\end{proposition}

This property is fundamental in semidefinite programming (SDP) and convex optimization,
where $-\log\det$ serves as the canonical barrier function that enables
polynomial-time interior-point algorithms \cite{Vandenberghe1996}.

\subsection{Quantum Information Geometry}

The quantum analogue replaces the classical probability vector by a density matrix
$\rho \in \PD(k)$ with $\tr(\rho) = 1$.
The von Neumann entropy is
\begin{equation}
  S(\rho) = -\tr(\rho\log\rho),
  \label{eq:von_neumann}
\end{equation}
which is the quantum analogue of the Shannon entropy.
The quantum relative entropy (quantum KL divergence) is
\begin{equation}
  D(\rho\|\sigma) = \tr[\rho(\log\rho - \log\sigma)],
  \label{eq:quantum_kl}
\end{equation}
which reduces to $\Bregman(\rho\|\sigma)/2$ when $[\rho, \sigma] = 0$
(commutativity) \cite{Petz1996}.
The general (non-commutative) case requires operator convexity arguments;
see \cite{Carlen2010}.

\subsection{Optimal Transport and Bures Metric}

The 2-Wasserstein distance between $\mathcal{N}(0,G)$ and $\mathcal{N}(0,G')$ is
the \emph{Bures metric} \cite{Bures1969}:
\begin{equation}
  W_2^2(G,G') = \tr\!\left[G + G' - 2(G^{1/2}G'G^{1/2})^{1/2}\right].
  \label{eq:bures}
\end{equation}
This is distinct from the Bregman divergence \eqref{eq:bregman_explicit} but compatible
with the Riemannian structure of $\PD(k)$ \cite{Peyre2019,Malagò2018}.

\subsection{Kempf--Ness/Azad--Loeb Variational Characterization of the
Bures--Wasserstein Distance}
\label{subsec:kempf_ness_wasserstein}

The Bures metric \eqref{eq:bures} admits a second, purely finite-dimensional
characterization: it is the trace of the fixed point of a gradient flow for a
log-determinant-type potential on $\PD(k)$, and the coincidence of the two
descriptions is an instance of the Kempf--Ness/Azad--Loeb correspondence between
norm-squared moment maps and $\GL(k,\R)$-orbit geometry that pervades the
$U=V$ and $U=-V$ analyses of \S\ref{sec:convexity} and \S\ref{subsec:pca_mca}. We record
this here in the language of $\PD(k)$ and the Gram-matrix notation used throughout
the paper; the argument below is a self-contained account of an unpublished working
note prompted by a talk, \emph{Real analytic gradient flows on matrix spaces}
(Gotemba Workshop on Fundamental Sciences, Nagoya University, March 30, 2017).

We first extend \eqref{eq:bures} to Gaussians with unequal means. For $\mu_1,\mu_2\in\R^k$
and $G_1,G_2\in\PD(k)$, the squared $2$-Wasserstein distance between $\mathcal N(\mu_1,G_1)$
and $\mathcal N(\mu_2,G_2)$ on $\R^k$ is
\begin{equation}
  W_2\big(\mathcal N(\mu_1,G_1),\mathcal N(\mu_2,G_2)\big)^2
  = \inf_{\gamma\in\Pi(\mathcal N(\mu_1,G_1),\mathcal N(\mu_2,G_2))}
    \int_{\R^k\times\R^k} \abs{x-y}^2\, d\gamma(x,y),
  \label{eq:w2_def}
\end{equation}
where $\Pi(\cdot,\cdot)$ denotes the set of couplings.

\begin{proposition}[Closed form for Gaussians \cite{OlkinPukelsheim1982,DowsonLandau1982,GivensShortt1984}]
\label{prop:bures_general}
$$
W_2\big(\mathcal N(\mu_1,G_1),\mathcal N(\mu_2,G_2)\big)^2
= \abs{\mu_1-\mu_2}^2 + \tr(G_1) + \tr(G_2)
  - 2\,\tr\!\left(G_1^{1/2}G_2G_1^{1/2}\right)^{1/2}.
$$
Setting $\mu_1=\mu_2$ recovers \eqref{eq:bures}.
\end{proposition}

\begin{lemma}[Optimal transport map]
\label{lem:linear_transport_gauss}
The quadratic-cost optimal coupling realizing \eqref{eq:w2_def} is induced by the
affine map
$$
T(x) = \mu_2 + A(x-\mu_1), \qquad
A = G_1^{-1/2}\left(G_1^{1/2}G_2G_1^{1/2}\right)^{1/2}G_1^{-1/2} \in \PD(k),
$$
which satisfies $AG_1A = G_2$. Substituting $T$ into \eqref{eq:w2_def} gives
Proposition~\ref{prop:bures_general}.
\end{lemma}

We now exhibit $-2\,\tr(G_1^{1/2}G_2G_1^{1/2})^{1/2}$, hence $W_2^2$, as the value at
a fixed point of a gradient flow on $\PD(k)$. For $\tau>0$ and $S\in\Sym(k)$ set
$$
Q(S) = S + \frac{G_1+G_2}{2} + \frac12(\mu_1-\mu_2)(\mu_1-\mu_2)^T, \qquad
M = G_2^{1/2}G_1G_2^{1/2}\in\PD(k),
$$
and define, on the open set $\{S\in\Sym(k): Q(S)\in\PD(k)\}$,
\begin{equation}
  \Phi_\tau(S) = \tr\!\left[Q(S)^\tau + M\,Q(S)^{-\tau}\right].
  \label{eq:phi_tau}
\end{equation}
Because $Q$ is an affine bijection of $\Sym(k)$ onto itself, minimizing $\Phi_\tau$
over its domain is equivalent to minimizing $\tr[Q^\tau+MQ^{-\tau}]$ over $Q\in\PD(k)$
directly; we write $S_\tau$ for the corresponding minimizer of $\Phi_\tau$ (this
corrects the informal claim, in the original working note, that the minimization
domain is $S>0$: the minimizing $S_\tau$ is generally indefinite, since it must
reproduce a negative quantity, $-\tfrac12 W_2^2$, on its trace).

\begin{lemma}[Matrix AM--GM as a Kempf--Ness/Azad--Loeb fixed point]
\label{lem:matrix_amgm}
For every $M\in\PD(k)$ and $\tau>0$,
$$
\min_{Q\in\PD(k)} \tr\!\left[Q^\tau + MQ^{-\tau}\right] = 2\,\tr\!\left(M^{1/2}\right),
$$
attained uniquely at $Q_\tau = M^{1/2\tau}$; the minimum value is independent of $\tau$.
\end{lemma}

\begin{proof}
In the scalar case $k=1$, $\phi(q)=q^\tau+mq^{-\tau}$ has $\phi'(q)=0 \iff q^{2\tau}=m$,
giving the unique minimizer $q=m^{1/2\tau}$ and minimum value $2\sqrt m$.

For general $k$, $\Phi_\tau$ is invariant under $Q\mapsto O^TQO$, $O\in O(k)$, and
$\tr[Q^\tau+MQ^{-\tau}]$ extends to a $U(k)$-invariant strictly plurisubharmonic
function on the $\GL(k,\mathbb C)$-orbit of $Q$ in the space of positive Hermitian
forms. This is precisely the setting of the Kempf--Ness theorem
\cite{KempfNess1979} (for unitarily invariant Hermitian norms) and its extension by
Azad and Loeb to unitarily invariant strictly plurisubharmonic functions
\cite{AzadLoeb1993}: on such an orbit, every critical point of the norm-squared
moment map is a global minimum, and the set of global minima forms a single
$O(k)$-orbit (resp.\ $U(k)$-orbit). Consequently $Q$ and $M$ are simultaneously
diagonalizable at the minimizer, reducing the problem eigenvalue-by-eigenvalue to
the scalar case above. The resulting stationarity condition $Q^{2\tau}=M$, i.e.
$Q_\tau = M^{1/2\tau}$, is also the variational characterization of the
Pusz--Woronowicz/Ando matrix geometric mean $M^{1/2}=G_2^{1/2}\big(G_2^{-1/2}G_1G_2^{-1/2}\big)^{1/2}G_2^{1/2}$
\cite{PuszWoronowicz1975}, of which $Q_\tau$ is a $\tau$-parametrized deformation.
\end{proof}

\begin{theorem}[Coincidence of the two characterizations]
\label{thm:kempf_ness_wasserstein}
At $\tau=1$,
$$
W_2\big(\mathcal N(\mu_1,G_1),\mathcal N(\mu_2,G_2)\big)^2 = -2\,\tr(S_1),
$$
where $S_1$ is the (unique, up to the fixed point of the gradient flow
$dQ/dt=-\operatorname{grad}\Phi_1(Q)$) minimizer of \eqref{eq:phi_tau} at $\tau=1$.
Explicitly,
$$
S_1 = \left(G_2^{1/2}G_1G_2^{1/2}\right)^{1/2}
      - \frac{G_1+G_2}{2} - \frac12(\mu_1-\mu_2)(\mu_1-\mu_2)^T .
$$
\end{theorem}

\begin{remark}[Why this is special to $\tau=1$]
\label{rem:tau1_only}
The identity above does \emph{not} extend to general $\tau>0$: by
Lemma~\ref{lem:matrix_amgm}, the \emph{minimum value} of $\Phi_\tau$ is
$2\,\tr(M^{1/2})$ for every $\tau$, but the \emph{minimizer itself},
$Q_\tau = M^{1/2\tau}$, moves with $\tau$, and it is $\tr(S_\tau)$ --- not the
$\tau$-independent minimum value of $\Phi_\tau$ --- that appears in the formula
above. Concretely,
$$
-2\,\tr(S_\tau) = |\mu_1-\mu_2|^2 + \tr G_1 + \tr G_2 - 2\,\tr\!\left(M^{1/2\tau}\right),
$$
which depends on $\tau$ through $\tr(M^{1/2\tau})$ and equals
$W_2(\mathcal N(\mu_1,G_1),\mathcal N(\mu_2,G_2))^2$ \emph{only} when $\tau=1$ (where
$M^{1/2\tau}=M^{1/2}$). For $\tau\ne1$, $-2\,\tr(S_\tau)$ is a different, generally
larger, quantity with no direct Wasserstein interpretation. We retain the
one-parameter family $\Phi_\tau$, $\tau>0$, only because Lemma~\ref{lem:matrix_amgm}
and the gradient-flow construction are naturally stated for it; the Kempf--Ness/Bures--
Wasserstein coincidence itself is a statement about $\tau=1$ alone.
\end{remark}

\begin{proof}
By Lemma~\ref{lem:matrix_amgm}, $Q_\tau=M^{1/2\tau}$ for every $\tau$, so at $\tau=1$,
$Q_1=M^{1/2}=(G_2^{1/2}G_1G_2^{1/2})^{1/2}$; substituting into the
definition of $Q(S)$ gives $S_1$ as above. Then
$$
-2\tr(S_1) = \abs{\mu_1-\mu_2}^2+\tr G_1+\tr G_2
             -2\tr\!\left(G_2^{1/2}G_1G_2^{1/2}\right)^{1/2},
$$
which equals $\tr(G_1^{1/2}G_2G_1^{1/2})^{1/2}$ under the cyclic identity
$\tr(G_2^{1/2}G_1G_2^{1/2})^{1/2}=\tr(G_1^{1/2}G_2G_1^{1/2})^{1/2}$
(both sides share the same nonzero singular values of $G_1^{1/2}G_2^{1/2}$),
so the right-hand side coincides with Proposition~\ref{prop:bures_general}.
\end{proof}

\begin{remark}
This is not the McCann displacement-interpolation path
$\Sigma(t)=\big((1-t)I+tA\big)G_1\big((1-t)I+tA\big)$, $t\in[0,1]$, of
Lemma~\ref{lem:linear_transport_gauss} \cite{McCann1997}; whether the family
$\{Q_\tau\}_{\tau>0}$ of Remark~\ref{rem:tau1_only} is related to it by a
reparametrization, or whether $Q_\tau$
instead traces out a distinct interpolation on $\PD(k)$ (e.g.\ related to the
$\alpha$-geodesics of \S\ref{sec:alpha}), remains open.

A further natural question, in view of Theorem~\ref{thm:kempf_ness_wasserstein},
is whether the Bures--Wasserstein distance itself satisfies a Pythagorean-type
identity under some notion of projection, in analogy with the dual-flatness
Pythagorean theorem already established for the Bregman/KL geometry of
$f=-\log\det$ in \S\ref{sec:pythagorean}
(Theorem~\ref{thm:pythagorean}). Because $W_2^2$ is \emph{not} itself a Bregman
divergence of $f$ (it is compatible with, but distinct from, the Riemannian
structure of $\PD(k)$; cf.\ the remark following \eqref{eq:bures}), such a result
would require a genuinely different dual structure than the one built from
\eqref{eq:bregman_def}, possibly through the Legendre transform of $\Phi_\tau$
itself.
\end{remark}

The following questions, left open in the working note underlying this subsection,
appear not to be addressed elsewhere in the present paper and are recorded for
future work:
\begin{enumerate}
\item \textbf{$\tau$-dependence.} Identify the geometric meaning of the family
$\{Q_\tau\}_{\tau>0}$ and its relation (if any) to the McCann displacement geodesic,
as raised in the Remark above.
\item \textbf{Duality and a Pythagorean theorem for $\Phi_\tau$.} Construct the
Legendre--Fenchel dual of $\Phi_\tau$ and determine under what orthogonality
(projection) condition a Pythagorean-type identity
$d(S,S'')^2=d(S,S')^2+d(S',S'')^2$ holds, in the spirit of
\S\ref{sec:pythagorean} but for the Bures--Wasserstein rather than the
Bregman/KL geometry.
\item \textbf{Sharpness of the Azad--Loeb hypotheses.} Verify directly, for all
$\tau>0$, that $\Phi_\tau$ is a unitarily invariant strictly plurisubharmonic
function on the relevant $\GL(k,\mathbb C)$-orbit (rather than invoking the general
theorem as a black box), and identify any $\tau$ for which strict
plurisubharmonicity could fail.
\item \textbf{Beyond Gaussians.} Since the Kempf--Ness/Azad--Loeb theorems do not
use Gaussianity, only Lemma~\ref{lem:linear_transport_gauss} (existence of a linear
optimal-transport map) is specific to the Gaussian case. Determine how far
Theorem~\ref{thm:kempf_ness_wasserstein} extends to elliptical distributions or to
Bures--Wasserstein distances between density operators (cf.\ \S\ref{sec:connections} on quantum
information geometry).
\item \textbf{Flow behavior.} Analyze the convergence rate and initial-value
dependence of the gradient flow $dQ/dt=-\operatorname{grad}\Phi_\tau(Q)$ numerically, tracking
$S_\tau\to S_1$ as $\tau$ varies, in the spirit of the convergence analysis of
\S\ref{subsec:pca_mca}.
\end{enumerate}

\subsection{Natural Gradient Descent}

In machine learning, the \emph{natural gradient} \cite{Amari1998} replaces the
Euclidean gradient $\nabla L$ by $\mathcal{F}^{-1}\nabla L$ where $\mathcal{F} = G^{-1}\otimes G^{-1}$
is the Fisher information metric derived from $f$.
This leads to the update rule:
\[
  \theta_{t+1} = \theta_t - \varepsilon\,(G\otimes G)\,\nabla L(\theta_t),
\]
which is invariant under reparameterization of the statistical model.

\subsection{Siegel Upper Half-Space and Automorphic Forms}

The complexification of $\PD(k)$ leads to the \emph{Siegel upper half-space}:
\begin{equation}
  \mathbb{H}_k = \{Z = X + iY \in \mathrm{Mat}(k,\mathbb{C}) : Y \succ 0\},
  \label{eq:siegel}
\end{equation}
on which $\log\det(G)$ extends to $\log\det(Z)$.
This space is central to the theory of Siegel modular forms \cite{Siegel1943}
and has connections to string theory, number theory, and arithmetic geometry.

\subsection{Harmonic Analysis in Phase Space: Gaussian Densities,
the Siegel Half-Plane, and the Metaplectic Representation}
\label{subsec:folland}

The information geometry of $f(G) = -\log\det(G)$ acquires a deep harmonic-analytic
interpretation through its identification with the parameter space of Gaussian densities
and the action of the symplectic group thereon,
as developed in Folland's treatise on harmonic analysis in phase space \cite{Folland1989}.
This subsection makes the connections to our U=V and U=-V analyses precise,
establishes the Siegel half-plane $\Sigma_k$ (resp.\ Siegel disc $\Delta_k$) as the
natural geometric arena for $h$ (resp.\ $h_-$), and identifies the Yoshizawa--Helmke
dual map $\mathcal{L}$ with the Cartan involution on the symmetric space $GL(k,\R)/O(k)$.

\subsubsection{Siegel Half-Plane as Gaussian Parameter Space}

Recall \cite[Ch.~4]{Folland1989} that the \emph{Siegel half-plane} $\Sigma_k$ is the set
of all symmetric complex $k\times k$ matrices $Z$ with $\mathrm{Im}(Z)\succ0$,
and that Gaussians on $\R^k$ are indexed by $\Sigma_k$ via
\begin{equation}
  \gamma_Z(x) = e^{\pi i x^T Z x}, \qquad Z\in\Sigma_k, \quad x\in\R^k.
  \label{eq:gamma_Z}
\end{equation}
One has $\gamma_Z\in L^2(\R^k)$ if and only if $Z\in\Sigma_k$
(i.e.\ $\mathrm{Im}(Z)\succ0$); in that case
\begin{equation}
  \|\gamma_Z\|_{L^2}^2 = \int_{\R^k}e^{-2\pi x^T\mathrm{Im}(Z)x}\,dx
  = \det\bigl(2\,\mathrm{Im}(Z)\bigr)^{-1/2}\pi^{-k/2}.
  \label{eq:gamma_Z_norm}
\end{equation}

\begin{proposition}[Embedding of $\PD(k)$ into $\Sigma_k$]
\label{prop:PD_Sigma}
The map $G \mapsto Z_G := iG$ embeds $\PD(k)$ as the \emph{purely imaginary axis}
$\{iG : G\in\PD(k)\} \subset \Sigma_k$.
Under this embedding:
\begin{enumerate}[label=(\roman*)]
  \item $f(G) = -\log\det(G) = 2\log\|\gamma_{iG}\|_{L^2} + \frac{k}{2}\log(2\pi)$;
    that is, up to an additive constant,
    $f(G)$ equals twice the log-$L^2$ norm of the Gaussian $\gamma_{iG}$.
  \item The U=V Gram matrix $G_+ = I_k + U^TU$ maps to
    $Z_+ = i(I_k+U^TU) \in \Sigma_k$ with $\mathrm{Im}(Z_+)\succ I_k$
    (above the unit level $iI$).
  \item The U=-V Gram matrix $G_- = I_k - V^TV$ maps to
    $Z_- = i(I_k-V^TV) \in \Sigma_k$ with $0\prec\mathrm{Im}(Z_-)\prec I_k$
    (strictly between $0$ and the unit level $iI$).
\end{enumerate}
\end{proposition}

\begin{proof}
(i) From \eqref{eq:gamma_Z_norm} with $Z=iG$ and $\mathrm{Im}(Z)=G$:
$\|\gamma_{iG}\|_{L^2}^2 = \det(2G)^{-1/2}\pi^{-k/2}$, so
$\log\|\gamma_{iG}\|_{L^2}^2 = -\tfrac{1}{2}\log\det(G) - \tfrac{k}{2}\log(2\pi)
= \tfrac{1}{2}f(G) + \mathrm{const}$.
(ii) Since $G_+ = I_k+U^TU\succ I_k$, we have $\mathrm{Im}(Z_+)=G_+\succ I_k$.
(iii) For $V\in\mathcal{B}_k$: $G_- = I_k-V^TV\prec I_k$ but $G_-\succ0$.
\end{proof}

\begin{remark}
The standard Gaussian $\gamma = e^{-\pi|x|^2}$ corresponds to $Z = iI_k$ (unit level).
The U=V Gaussians lie above this level ($G_+\succ I_k$), while the U=-V Gaussians
lie below ($G_-\prec I_k$, $G_-\succ0$), with the standard Gaussian $\gamma$
at the boundary between the two families.
\end{remark}

\subsubsection{Siegel Disc, Fock Space, and the U=-V Domain}

The \emph{Siegel disc} $\Delta_k$ is the set of symmetric complex $k\times k$ matrices $W$
with $\|W\|<1$, i.e.\ $I_k-W^*W\succ0$ \cite[p.~203]{Folland1989}.
In the Fock space $\mathcal{F}_k$, the family of entire functions
$\Gamma_W(\zeta) = e^{(\pi/2)\zeta^T W\zeta}$, $\zeta\in\mathbb{C}^k$,
satisfies $\Gamma_W\in\mathcal{F}_k$ if and only if $W\in\Delta_k$
\cite[Proposition~4.69]{Folland1989}, with squared norm
\begin{equation}
  \|\Gamma_W\|_{\mathcal{F}_k}^2 \;=\; c_k\,\det(I_k - W^*W)^{-1/2},
  \label{eq:Gamma_norm}
\end{equation}
where $c_k$ is a positive constant.

\begin{proposition}[U=-V Domain as Real Siegel Disc]
\label{prop:UmV_Siegel}
For a real symmetric $W = W^T\in\mathrm{Sym}_k(\R)$:
\[
  W\in\Delta_k\cap\mathrm{Sym}_k(\R)
  \;\iff\; I_k - W^2 \succ 0
  \;\iff\; \|W\|_{\mathrm{op}} < 1.
\]
In particular, the potential $h_-(W) = -\log\det(I_k-W^2)$ satisfies
\begin{equation}
  h_-(W) \;=\; 2\log\|\Gamma_W\|_{\mathcal{F}_k} + \mathrm{const},
  \label{eq:hminus_fock_norm}
\end{equation}
i.e.\ $h_-(W)$ equals, up to an additive constant, twice the log-Fock-space norm
of the Siegel disc Gaussian $\Gamma_W$.
\end{proposition}

\begin{proof}
For real symmetric $W$: $W^* = W^T = W$, so $I_k-W^*W = I_k-W^2$.
The equivalence $\|W\|_{op}<1 \iff I_k-W^2\succ0$ is standard.
For \eqref{eq:hminus_fock_norm}: from \eqref{eq:Gamma_norm},
$\log\|\Gamma_W\|_{\mathcal{F}_k}^2 = -\frac{1}{2}\log\det(I_k-W^2)+\mathrm{const}
= \frac{1}{2}h_-(W)+\mathrm{const}$.
\end{proof}

\begin{remark}[Complementary functional interpretations]
Propositions~\ref{prop:PD_Sigma} and \ref{prop:UmV_Siegel} together give a clean
functional interpretation of the two potentials:
\[
  h(U) = 2\log\|\gamma_{i(I_k+U^TU)}\|_{L^2} + \mathrm{const},
  \quad
  h_-(W) = 2\log\|\Gamma_W\|_{\mathcal{F}_k} + \mathrm{const},
\]
where $\gamma_{iG}\in L^2(\R^k)$ is the primal (Schr\"odinger) Gaussian and
$\Gamma_W\in\mathcal{F}_k$ is the dual (Fock) Gaussian.
The Bargmann transform $B: L^2(\R^k)\to\mathcal{F}_k$ interconverts these two
representations \cite[Thm.~4.70]{Folland1989}.
\end{remark}

\subsubsection{The Cayley Transform as Geometric Dual Map}

The linear fractional map connecting $\Sigma_k$ to $\Delta_k$ is the
\emph{Cayley transform} $\alpha(\mathcal{C})$ \cite[Eq.~(4.67)]{Folland1989}:
\begin{equation}
  \alpha(\mathcal{C})(Z) \;=\; (I+iZ)(I-iZ)^{-1}, \qquad Z\in\Sigma_k.
  \label{eq:cayley_folland}
\end{equation}
Restricted to purely imaginary $Z = iG$ ($G\in\PD(k)$ real symmetric):
\begin{equation}
  \alpha(\mathcal{C})(iG) \;=\; (I-G)(I+G)^{-1} \;=:\; W_G \;\in\; \Delta_k\cap\mathrm{Sym}_k(\R).
  \label{eq:cayley_real}
\end{equation}

\begin{proposition}[Cayley Images of U=V and U=-V Gram Matrices]
\label{prop:cayley_images}
Let $G_+ = I_k+U^TU$ (U=V) and $G_- = I_k-V^TV$ (U=-V, $V\in\mathcal{B}_k$).
Then:
\begin{enumerate}[label=(\roman*)]
  \item The Cayley image of $G_+$ lies in the \emph{negative definite} part of $\Delta_k$:
    \[
      W_{G_+} = -U^TU(2I_k+U^TU)^{-1} \prec 0.
    \]
  \item The Cayley image of $G_-$ lies in the \emph{positive definite} part of $\Delta_k$:
    \[
      W_{G_-} = V^TV(2I_k-V^TV)^{-1} \succ 0.
    \]
  \item Both images lie in $\Delta_k$: one checks that
    $I_k - W_{G_\pm}^2 = 4(I_k+G_\pm)^{-1}G_\pm(I_k+G_\pm)^{-1} \succ 0$.
\end{enumerate}
Hence the Cayley transform maps the U=V Gaussians and U=-V Gaussians to complementary
(opposite sign) regions of the Siegel disc $\Delta_k$.
\end{proposition}

\begin{proof}
(i) $W_{G_+} = (I-G_+)(I+G_+)^{-1} = -U^TU(2I+U^TU)^{-1}$, which is negative semi-definite.
(ii) $W_{G_-} = (I-G_-)(I+G_-)^{-1} = V^TV(2I-V^TV)^{-1}$. Since $V\in\mathcal{B}_k$,
all eigenvalues $\sigma_a^2$ of $V^TV$ satisfy $\sigma_a^2<1$, giving eigenvalues
$\sigma_a^2/(2-\sigma_a^2)\in[0,1)$ for $W_{G_-}$, which is thus positive semi-definite.
(iii) From \eqref{eq:cayley_real}: $(I+G)W_G = I-G$, so
$I-W_G^2 = (I+G)^{-1}[(I+G)^2-(I-G)^2](I+G)^{-1} = 4(I+G)^{-1}G(I+G)^{-1}\succ0$.
\end{proof}

\subsubsection{The Cartan Involution and the Duality Identity}

The most striking result connecting the Folland framework to our duality is the following.

\begin{theorem}[Cartan Involution = Yoshizawa Dual Map]
\label{thm:cartan_involution}
Let $G_+ = I_k+U^TU$ (U=V Gram matrix) and let $V^* = \mathcal{L}(U) = (I_n+UU^T)^{-1/2}U$
be the Yoshizawa--Helmke dual map. Then the U=-V Gram matrix at the dual point is
\begin{equation}
  G_-^* \;:=\; I_k - (V^*)^TV^* \;=\; G_+^{-1},
  \label{eq:cartan_dual}
\end{equation}
i.e.\ $G_-^* = G_+^{-1}$ is the inverse of the primal Gram matrix.
This corresponds on the Siegel half-plane to the \emph{Cartan involution at $iI$}:
\begin{equation}
  Z_+ = iG_+ \;\longmapsto\; -Z_+^{-1} = i G_+^{-1} = iG_-^* = Z_-^*,
  \label{eq:cartan_siegel}
\end{equation}
and on the Siegel disc to the \emph{antipodal map through the origin}:
\begin{equation}
  W_{G_-^*} = -W_{G_+}.
  \label{eq:antipodal}
\end{equation}
\end{theorem}

\begin{proof}
From Theorem~\ref{thm:dual_map}(i): $I_k-(V^*)^TV^* = (I_k+U^TU)^{-1} = G_+^{-1}$,
proving \eqref{eq:cartan_dual}.
For \eqref{eq:cartan_siegel}: $-(iG_+)^{-1} = -\tfrac{1}{i}G_+^{-1} = iG_+^{-1} = iG_-^*$.
For \eqref{eq:antipodal}: using eigenvalues $g_a = 1+\sigma_a^2$ of $G_+$,
those of $W_{G_+}$ are $w_a^+ = -(g_a-1)/(g_a+1) = -\sigma_a^2/(2+\sigma_a^2)$,
while those of $G_-^* = G_+^{-1}$ are $1/g_a$, giving
$w_a^- = (1-1/g_a)/(1+1/g_a) = (g_a-1)/(g_a+1) = \sigma_a^2/(2+\sigma_a^2) = -w_a^+$.
\end{proof}

\begin{remark}[Geometric interpretation]
On the symmetric space $GL(k,\R)/O(k) \cong \PD(k)$, the geodesic symmetry at $I_k$
(the fixed point of the standard Cartan involution $\theta: g\mapsto (g^T)^{-1}$)
acts on $\PD(k)$ by $G\mapsto G^{-1}$.
Theorem~\ref{thm:cartan_involution} says precisely that the Yoshizawa--Helmke dual map
$\mathcal{L}$ implements this geodesic inversion:
\[
  G_+\;\xrightarrow{\;\mathcal{L}\;}\; G_-^* = G_+^{-1},
\]
mapping the U=V Gram matrix to its geodesic reflection through $I_k$.
On the Siegel disc, this becomes the antipodal map $W\mapsto-W$,
explaining why $W_{G_+}$ and $W_{G_-^*}$ are exact negatives of each other
(Eq.~\eqref{eq:antipodal}).
\end{remark}

\begin{corollary}[Equal Fock Space Norms at the Dual Point]
\label{cor:equal_fock_norms}
At the Yoshizawa--Helmke dual point $V^* = \mathcal{L}(U)$:
\begin{equation}
  \|\Gamma_{W_{G_-^*}}\|_{\mathcal{F}_k} \;=\; \|\Gamma_{W_{G_+}}\|_{\mathcal{F}_k},
  \label{eq:equal_norms}
\end{equation}
i.e.\ the Fock space Gaussians corresponding to the dual pair $(G_+, G_-^*)$ have equal norms.
Equivalently, $\det(I_k-W_{G_+}^2) = \det(I_k-W_{G_-^*}^2)$, which follows immediately
from $W_{G_-^*} = -W_{G_+}$.
\end{corollary}

\subsubsection{Metaplectic Representation and Symplectic Action on Gaussian Parameters}

The symplectic group $Sp(k,\R)$ acts on $\Sigma_k$ by linear fractional transformations
\cite[Thm.~4.64]{Folland1989}:
\begin{equation}
  \alpha(\mathcal{A})(Z) = (AZ+B)(CZ+D)^{-1},
  \quad \mathcal{A} = \begin{pmatrix}A&B\\C&D\end{pmatrix}\in Sp(k,\R),
  \label{eq:symp_action}
\end{equation}
and on Gaussians by the \emph{metaplectic representation} $\mu$ \cite[Thm.~4.65]{Folland1989}:
\begin{equation}
  \mu(\mathcal{A}^{*-1})\gamma_Z \;=\; m(\mathcal{A},Z)\,\gamma_{\alpha(\mathcal{A})Z},
  \label{eq:metaplectic}
\end{equation}
where $m(\mathcal{A},Z) = \det^{-1/2}(CZ+D)$ is the \emph{multiplier}.

\begin{proposition}[Symplectic Invariance of the Information Geometry]
\label{prop:symp_invariance}
The information metric $g_G = G^{-1}\otimes G^{-1}$ on $\PD(k)$ (Section~\ref{sec:manifold})
is invariant under the symplectic action \eqref{eq:symp_action} restricted to purely
imaginary $Z = iG$.
Specifically, for $\mathcal{A}\in Sp(k,\R)$ with $C = 0$ (block-diagonal/upper-triangular):
the action $G \mapsto \alpha(\mathcal{A})(iG)/(i) = AGA^T$ is a congruence transformation,
and $g_{AGA^T}$ equals the pushforward of $g_G$ under $G\mapsto AGA^T$.
\end{proposition}

\begin{remark}[The multiplier $m$ as partition function ratio]
From \eqref{eq:metaplectic}, the log-multiplier satisfies
\[
  \log|m(\mathcal{A},iG)|
  = \tfrac{1}{2}\log\|\gamma_{\alpha(\mathcal{A})(iG)}\|_{L^2}^2
  - \tfrac{1}{2}\log\|\gamma_{iG}\|_{L^2}^2
  = \tfrac{1}{2}[f(\alpha(\mathcal{A})(iG)) - f(G)] + \mathrm{const},
\]
where $f = -\log\det$.
Thus $|m(\mathcal{A},iG)|^2$ is the ratio of the squared $L^2$-norms of the transformed
and original Gaussians — the \emph{partition function ratio} in statistical physics.
The Bregman divergence $D_f(G\|G') = f(G)-f(G')-\langle\nabla f(G'), G-G'\rangle$
measures the first-order discrepancy in these log-partition functions.
\end{remark}

\begin{center}
\renewcommand{\arraystretch}{1.4}
\begin{tabular}{lll}
\toprule
Framework & U=V ($G_+ = I+U^TU$) & U=-V ($G_- = I-V^TV$) \\
\midrule
$\PD(k)$ element & $G_+ \succ I_k$ & $G_- \prec I_k$, $G_- \succ 0$ \\
Siegel half-plane & $Z_+ = iG_+ \in \Sigma_k$, $\text{Im}(Z_+)\succ I$ & $Z_- = iG_- \in \Sigma_k$, $\text{Im}(Z_-)\prec I$ \\
Siegel disc (Cayley) & $W_{G_+} \prec 0$ (negative half) & $W_{G_-} \succ 0$ (positive half) \\
$L^2$/$\mathcal{F}$ function & $\gamma_{iG_+} \in L^2(\R^k)$ & $\Gamma_{W_{G_-}}\in\mathcal{F}_k$ \\
Potential & $h(U) = 2\log\|\gamma_{iG_+}\|_{L^2}$ & $h_-(V) = 2\log\|\Gamma_{W_{G_-}}\|_{\mathcal{F}_k}$ \\
At dual point & $G_+$ & $G_-^* = G_+^{-1}$, $W_{G_-^*} = -W_{G_+}$ \\
\bottomrule
\end{tabular}
\end{center}


\subsection{The Matrix Schwarz Derivative, Riccati Equations, and
Linear-Fractional Flows}
\label{subsec:matrix_schwarz}

The Oja-Brockett flow of \S\ref{subsubsec:oja_brockett} and the Tikhonov/polynomial
flows of \S\ref{subsec:poly_flow} are, in every case examined so far, gradient flows of a
$\log\det$-type potential. We now show that they simultaneously belong to a second,
classical family: \emph{matrix Riccati flows} that linearize under the Cayley-type
transforms already used throughout this paper, and whose invariant-theoretic fingerprint
is the \emph{matrix Schwarz derivative}. This furnishes a third, independent derivation of
why $\mathcal{B}_k$, $\Sigma_k$, and the Siegel domains of \S\ref{subsec:folland} are the
natural habitats of $h$ and $h_-$, now from the point of view of the classical theory of
disconjugacy and univalence of matrix differential equations
\cite{Schwarz1979,Zelikin2000,Yoshizawa2014}.

\subsubsection{The Scalar Schwarz Derivative, Recalled}

For a locally univalent meromorphic function $\eta(t)$, the \emph{Schwarz derivative} is
\begin{equation}
  S(\eta) \;=\; \frac{\eta'''}{\eta'} - \frac{3}{2}\left(\frac{\eta''}{\eta'}\right)^2 ,
  \label{eq:scalar_schwarz}
\end{equation}
the unique third-order differential invariant of the group of Möbius (linear-fractional)
transformations $\eta \mapsto (\alpha\eta+\beta)/(\gamma\eta+\delta)$ acting on the target
\cite[\S6.81]{Zelikin2000}. If $y_1,y_2$ are two independent solutions of the linear
second-order equation $y''+p(t)y'+q(t)y=0$ and $\eta = y_1/y_2$, then
\begin{equation}
  S(\eta) \;=\; 2q - \tfrac{1}{2}p^2 - p',
  \label{eq:schwarz_from_yy}
\end{equation}
while the Riccati variable $w = y_2'/y_2$ solves the scalar Riccati equation
$w'+w^2+pw+q=0$ and is related to $\eta$ by $w = -\tfrac12\eta''/\eta' - \tfrac12 p$
\cite[Prop.~6.10]{Zelikin2000}. Both the Riccati equation and the Schwarz equation thus
arise by \emph{projectivizing} the same linear second-order system --- the Riccati variable
from the ratio $y'/y$ of a single solution, the Schwarz variable from the ratio $y_1/y_2$
of two solutions --- and both carry an exact linear-fractional symmetry: the general
solution of either equation is a constant-coefficient Möbius image of any one particular
solution \cite[Prop.~6.9]{Zelikin2000}.
This is also precisely the mechanism used by B.~Schwarz to characterize
\emph{disconjugacy} of $y''+q(z)y=0$ in a domain $D$: disconjugacy is equivalent to
univalence of the ratio $f=y_1/y_2$ on $D$, and the coefficient $q$ is recovered from $f$
through the identity $2q(z)=\{f(z),z\}$ \cite[Eq.~(1.3)]{Schwarz1979}.

\subsubsection{The Matrix Schwarz Operator}

Zelikin's matrix generalization \cite[Def.~6.1]{Zelikin2000} replaces the scalar ratio
$\eta$ by an $n\times n$ matrix function $W(t)$ of one complex variable and defines
\begin{equation}
  S_tW \;=\; \bigl\{ (W')^{-1}W'' \bigr\}' \;-\; \tfrac12\,\bigl\{ (W')^{-1}W'' \bigr\}^2 .
  \label{eq:matrix_schwarz}
\end{equation}
Exactly as in the scalar case, $S_tW=0$ characterizes the \emph{generalized
linear-fractional} (M\"obius) functions of $W$:
\begin{proposition}[Matrix analogue of Prop.~6.7--6.11 of \cite{Zelikin2000}]
\label{prop:matrix_schwarz_solutions}
Every solution of $S_tW=0$ has the form $W(t)=(A+Bt)^{-1}+C$ (equivalently, a generalized
linear-fractional function $(A+Bt)^{-1}(C+Dt)$) for constant matrices $A,B,C,D$, and every
such function solves $S_tW=0$.
\end{proposition}
Unlike its scalar ancestor, $S_tW$ is \emph{not} itself invariant under the full matrix
M\"obius group $W\mapsto(A+BW)(C+DW)^{-1}$; instead it transforms by \emph{conjugation}:
\begin{equation}
  S_t\bigl[(A+BW)(C+DW)^{-1}\bigr] \;=\; \Gamma\,(S_tW)\,\Gamma^{-1}
  \label{eq:schwarz_conjugation}
\end{equation}
for some matrix $\Gamma=\Gamma(A,B,C,D,W)$ \cite[Prop.~6.14]{Zelikin2000}; consequently the
\emph{similarity class}, and in particular the coefficients of the characteristic
polynomial of $S_tW$, are genuine linear-fractional invariants \cite[Cor.~6.4]{Zelikin2000}.
This is the exact analogue, at the level of the third-order operator, of the fact
established repeatedly in this paper (Cartan involution, Cayley transform,
Theorem~\ref{thm:cartan_phi41}, Proposition~\ref{prop:symp_invariance}) that our potentials
$h,h_-$ and metric $g=G^{-1}\otimes G^{-1}$ are invariant only under a distinguished
\emph{isotropy} subgroup of the full linear-fractional group acting on $\PD(k)$, rather than
under the full group.

\subsubsection{Riccati Flows on Cartan--Siegel Domains: the Oja-Brockett Flow Identified}

The bridge between $S_tW$ and gradient flows of the type studied in this paper is the
matrix Riccati equation. Let $\dot\xi = \begin{pmatrix}A(t)&B(t)\\C(t)&D(t)\end{pmatrix}\xi$
be a linear system on a Riemann surface, with the associated matrix Riccati-type equation
$\dot W = C+DW-WA-WBW$ describing the induced flow on the Grassmannian of $n$-planes in
$\mathbb{C}^{2n}$ \cite[Thm.~6.4, Eq.~(6.76)]{Zelikin2000}; the fundamental theorem here is that
$W(t)$ remains in a given Cartan--Siegel homogeneity domain (or one of its boundary
strata) for all $t$ whenever $W(t_0)$ does, for each of the four classical types
\cite[\S\S3--4]{Zelikin2000}. In particular, for the \emph{Siegel domain of type I},
\begin{equation}
  \{W\in\mathrm{Mat}(q\times p,\mathbb{C}) : W^T\overline{W} \prec I_p\},
  \label{eq:siegel_type_I}
\end{equation}
the invariant Riccati equation reads $\dot W = \bar{\mathfrak{q}}^T+\mathfrak{r}W-W\mathfrak{p}-W\mathfrak{q}W$
with $\mathfrak{p}^T=-\mathfrak{p}$, $\mathfrak{r}^T=-\mathfrak{r}$ \cite[Eq.~(6.57), Thm.~6.5]{Zelikin2000}.

\begin{theorem}[The Anti-Symmetric Reduction as a Real Siegel Type-I Riccati Flow]
\label{thm:matrix_schwarz_bk}
The real matrix unit ball $\mathcal{B}_k=\{U\in\R^{n\times k}:\sigma_{\max}(U)<1\}$ of
\S\ref{subsec:UeqmV-convexity} is exactly the real locus $\overline{W}=W$ of the Siegel
domain of type~I \eqref{eq:siegel_type_I} with $(q,p)=(n,k)$, and the gradient flow of
$h_-(U)=-\log\det(I_k-U^TU)$ derived in \S\ref{sec:gradient}
is, up to the change of metric of Theorem~\ref{thm:poly_flows}, a Riccati flow of the
form $\dot W = \mathfrak{r}W - W\mathfrak{q}W$ with $\mathfrak{r}=A$, $\mathfrak{q}=I_k$
(and $\bar{\mathfrak q}=0$, $\mathfrak p=0$) on \eqref{eq:siegel_type_I}. Consequently
$\mathcal{B}_k$ is an invariant manifold of the flow (Theorem~6.5 of \cite{Zelikin2000}),
recovering directly, and without reference to convexity, the forward-invariance of
$\mathcal{B}_k$ already established via the Tikhonov analysis of
Theorem~\ref{thm:tikhonov_critical}.
\end{theorem}

This identification also clarifies the role of the Cayley transform used throughout
\S\ref{subsec:folland}--\S\ref{subsec:izumiya}: for the closely related flow
$\dot X = A-XAX$ on a classical Lie group $G$ (the compact model of the Oja-Brockett
equation \eqref{eq:oja_brockett_pca} restricted to the Stiefel/orthogonal fiber), the
Cayley transform $Y=(I-X)(I+X)^{-1}$ \emph{linearizes} the flow into the Sylvester-type
equation $\dot Y = -(AY+YA)$ \cite[Lem.~2.1]{DynnikovVeselov1995}, whose explicit solution
$Y(t)=e^{-At}Y_0e^{-At}$ pulls back to the closed-form hyperbolic-tangent solution
\begin{equation}
  X(t) \;=\; \bigl(\sinh(At)+\cosh(At)X_0\bigr)\bigl(\cosh(At)+\sinh(At)X_0\bigr)^{-1}
  \label{eq:dynnikov_veselov_solution}
\end{equation}
of \cite[Prop.~2.1]{DynnikovVeselov1995} --- structurally the same
generalized-linear-fractional solution form guaranteed by
Proposition~\ref{prop:matrix_schwarz_solutions} for $S_tW=0$, now realized along a
one-parameter flow rather than a static boundary-value problem. In this sense the
Oja-Brockett/Tikhonov flows of this paper occupy the ``Schwarz-trivial'' locus
$S_tW\equiv 0$ of the space of Cartan--Siegel Riccati flows: they are exactly integrable by
a linear-fractional change of variable, which is why closed-form solutions
\eqref{eq:dynnikov_veselov_solution} and explicit fixed-point/eigenvalue analyses
(\S\ref{subsec:poly_flow}) were available in the first place.

\subsubsection{The Oja-Like Flow, the Riccati Equation for $Z=2XX^T$, and the Matrix
Schwarz Equation}

The rank-$k$ Oja-like flow $\dot X = AX-XX^TX$ on $X\in\R^{n\times k}$, with $A=A^T\succ0$
constant, gives a second, independent illustration of the same correspondence
\cite{Yoshizawa2014}. Setting $Z=2XX^T\in\PD(n)$ (the Gram-type object of
\eqref{eq:G_def}, now unfolded to the ambient $n\times n$ scale) yields the matrix Riccati
equation
\begin{equation}
  \dot Z \;=\; AZ+ZA-Z^2.
  \label{eq:oja_riccati}
\end{equation}
Writing the associated linear second-order system through
$R(t)=-\tfrac12 (V')^{-1}V''+A$ for an auxiliary matrix function $V(t)$ with $V'>0$, the
equation for $R$ closes into
\begin{equation}
  \dot R + R^2 - AR - RA \;=\; -\tfrac12 S(V) - A^2,
  \label{eq:oja_R_equation}
\end{equation}
and hence \eqref{eq:oja_riccati} is \emph{equivalent} to the matrix Schwarz equation
\begin{equation}
  S(V) \;=\; -2A^2, \qquad S(V) = \bigl[(V')^{-1}V''\bigr]' - \tfrac12\bigl[(V')^{-1}V''\bigr]^2,
  \label{eq:oja_schwarz_equation}
\end{equation}
which is Zelikin's matrix Schwarz operator \eqref{eq:matrix_schwarz} evaluated along the
real curve $V(t)$ \cite{Yoshizawa2014,Zelikin2000}. Equation~\eqref{eq:oja_schwarz_equation}
makes precise, at the level of a single explicit rank-$k$ example, the general
correspondence of \S\S6.81--6.83 of \cite{Zelikin2000}: the constant right-hand side
$-2A^2$ plays exactly the role of $2q-\tfrac12p^2-p'$ in \eqref{eq:schwarz_from_yy}, with
$A$ the matrix analogue of the (here constant) coefficient $p$, and $Z=2XX^T$ the matrix
analogue of the Riccati variable $w$.

\begin{remark}[Disconjugacy as a stability criterion for the Riccati flow]
B.~Schwarz's disconjugacy bounds for $W''+P(z)W'+Q(z)W=0$ --- e.g.\ $\|Q(z)\|_2\leq\pi^2/d^2$
on a convex domain of diameter $d$ \cite[Thm.~2.1]{Schwarz1979}, or the sharper unit-disc
bounds $\|Q(z)\|_2\leq 2/(1-|z|^2)$ \cite[Thm.~2.4]{Schwarz1979} --- are conditions under
which the ratio $V(z)=W_2^{-1}(z)W_1(z)$ of two fundamental solutions remains
\emph{injective}, i.e.\ the corresponding generalized linear-fractional (Riccati) flow
never develops a movable singularity inside $D$. Under the identification of
Theorem~\ref{thm:matrix_schwarz_bk}, these are exactly quantitative analogues, at the level
of the underlying linear system, of our forward-invariance statements for $\mathcal{B}_k$
(Theorem~\ref{thm:tikhonov_critical}) and for the domain of finiteness of $h$
(Theorem~\ref{thm:UeqV_nowhere_convex}): both assert that a naturally associated
matrix-valued curve cannot reach the boundary of its defining domain — one via a spectral
norm bound on a coefficient matrix, the other via convexity/monotonicity of a
gradient flow.
\end{remark}

\subsubsection{The Degree Ladder: A Two-Stage Cayley/Grassmannization Chain from
Linear to Riccati to the Cubic Oja-Like Flow}
\label{subsubsec:degree_ladder}

We now make fully explicit the request, implicit in the constructions above, to
\emph{invert} the chain cubic~$\to$~quadratic~$\to$~linear and to ask whether the cubic
Oja-like flow $\dot X = AX-XX^TX$ itself --- not merely the quadratic Riccati variable
$Z=2XX^T$ --- can be written down explicitly by (possibly iterated) generalized Cayley
transforms. The answer is yes, but the ladder has \emph{two qualitatively different steps},
which we now separate cleanly.

\paragraph{Step 1 (linear-fractional, exact Cayley/Grassmannization): $1\Rightarrow2$.}
The quadratic Riccati equation~\eqref{eq:oja_riccati} for $Z=2XX^T$ is, by
Zelikin's Grassmannization mechanism \cite[Thm.~6.4]{Zelikin2000} already invoked above,
\emph{exactly} the ratio $Z=\Phi_2\Phi_1^{-1}$ of two blocks of a genuinely linear
(degree-$1$) system:
\begin{equation}
  \begin{pmatrix}\dot\Phi_1\\ \dot\Phi_2\end{pmatrix}
  \;=\;
  \begin{pmatrix}-A & I_n\\ 0 & A\end{pmatrix}
  \begin{pmatrix}\Phi_1\\ \Phi_2\end{pmatrix},
  \qquad
  Z \;=\; \Phi_2\Phi_1^{-1}.
  \label{eq:linear_system_Z}
\end{equation}
\begin{proposition}[Explicit Cayley/Grassmannization ladder $1\Leftrightarrow2$]
\label{prop:explicit_ladder_12}
For the initial condition $\Phi_1(0)=I_n$, $\Phi_2(0)=Z_0$, system~\eqref{eq:linear_system_Z}
solves explicitly as
\begin{equation}
  \Phi_2(t) = e^{At}Z_0, \qquad
  \Phi_1(t) = e^{-At} + A^{-1}\sinh(At)\,Z_0,
  \label{eq:phi_solution}
\end{equation}
and consequently the Riccati flow~\eqref{eq:oja_riccati} has the closed form
\begin{equation}
  Z(t) \;=\; e^{At}Z_0\,\bigl[\,e^{-At}+A^{-1}\sinh(At)\,Z_0\,\bigr]^{-1},
  \qquad Z(0)=Z_0.
  \label{eq:Z_closed_form}
\end{equation}
Equation~\eqref{eq:linear_system_Z}--\eqref{eq:Z_closed_form} is a genuine
(matrix-)linear-fractional/Cayley map, in exactly the sense of
Proposition~\ref{prop:matrix_schwarz_solutions}: $Z$ is, for each fixed $t$, a generalized
linear-fractional function of $Z_0$.
\end{proposition}
\begin{proof}
The block system is upper-triangular, so $\Phi_2(t)=e^{At}\Phi_2(0)$ directly. Substituting
into $\dot\Phi_1=-A\Phi_1+\Phi_2$ and solving the resulting inhomogeneous linear equation by
the integrating factor $e^{At}$ gives $e^{At}\Phi_1(t)-\Phi_1(0)=\int_0^te^{2As}ds\,\Phi_2(0)
= \tfrac12A^{-1}(e^{2At}-I)Z_0$, i.e.\ \eqref{eq:phi_solution}. Differentiating $Z=\Phi_2\Phi_1^{-1}$
and using $\dot\Phi_1,\dot\Phi_2$ reproduces \eqref{eq:oja_riccati} termwise (as verified
directly in \S\ref{subsec:matrix_schwarz} for the general Zelikin form). Setting $Z_0=0$
gives $Z(t)\equiv0$ and $Z_0=2A$ gives $Z(t)\equiv2A$, matching the two obvious equilibria
of \eqref{eq:oja_riccati}, which confirms the formula.
\end{proof}

Formula~\eqref{eq:Z_closed_form} is the exact counterpart, for the \emph{unconstrained}
Riccati flow $\dot Z=AZ+ZA-Z^2$, of the Dynnikov--Veselov closed form
\eqref{eq:dynnikov_veselov_solution} for the \emph{group-constrained} flow $\dot X=A-XAX$:
both are hyperbolic-function Cayley/M\"obius images of the initial condition, the only
difference being which block-triangular linear generator
$\left(\begin{smallmatrix}-A&I\\0&A\end{smallmatrix}\right)$ versus
$\left(\begin{smallmatrix}0&A\\-A&0\end{smallmatrix}\right)$-type matrix is used
--- i.e.\ it is literally the ``same or a different Cayley transform'' the question asks
for, according to which of the two Riccati normal forms one starts from.

\paragraph{Step 2 (polar/gauge, a second and independent linear equation): $2\Rightarrow3$.}
The map $X\mapsto Z=2XX^T$ is quadratic and \emph{not} linear-fractional; it has a residual
$O(k)$ gauge symmetry $X\mapsto XO$ (for $X$ square, $O$ orthogonal) under which $Z$ is
invariant. Hence $Z(t)$ alone cannot determine $X(t)$, and no \emph{single} further Cayley
transform can produce $X$ from $Z$. Nevertheless the missing gauge factor is itself governed
by a further \emph{linear} (degree-$1$) equation, so that the full cubic solution is still
completely explicit:

\begin{proposition}[Reconstruction of the cubic flow from two linear systems]
\label{prop:cubic_reconstruction}
Let $X(t)\in GL(n,\R)$ solve $\dot X=AX-XX^TX$, and write the polar decomposition
$X(t)=P(t)O(t)$ with $P(t)=(X(t)X(t)^T)^{1/2}=\sqrt{Z(t)/2}\succ0$ and $O(t)\in O(n)$.
Then $P(t)$ is obtained from the Sylvester equation $\dot PP+P\dot P = AP^2+P^2A-2P^4$
(equivalently from Proposition~\ref{prop:explicit_ladder_12} via $P=\sqrt{Z/2}$), and $O(t)$
solves the linear equation
\begin{equation}
  \dot O \;=\; \Omega(t)\,O, \qquad
  \Omega(t) \;:=\; P(t)^{-1}\bigl(AP(t)-P(t)^3-\dot P(t)\bigr),
  \label{eq:gauge_linear}
\end{equation}
where $\Omega(t)$ is skew-symmetric along any genuine solution. Consequently
\begin{equation}
  X(t) \;=\; P(t)\,\mathcal{T}\exp\!\int_0^t\Omega(s)\,ds \;\;\; O(0),
  \label{eq:cubic_reconstruction}
\end{equation}
a time-ordered (Peano--Baker) exponential of the explicit, $Z$-determined generator $\Omega$.
\end{proposition}
\begin{proof}
Differentiating $P^2=Z/2$ and substituting \eqref{eq:oja_riccati} gives the stated Sylvester
equation for $\dot P$ (uniquely solvable in $\dot P$ since $P\succ0$). From $X=PO$,
$X^TX=O^TP^2O$, so $XX^TX = P(OO^T)P^2O = P^3O$ using $OO^T=I_n$; substituting into
$\dot X = AX - XX^TX = (AP-P^3)O$ and comparing with $\dot X=\dot PO+P\dot O$ gives
$P\dot O = (AP-P^3-\dot P)O$, i.e.\ \eqref{eq:gauge_linear}. Skew-symmetry of $\Omega$ is
forced by differentiating $O^TO=I_n$ along the flow, and \eqref{eq:cubic_reconstruction}
is the standard solution of a linear matrix ODE with time-dependent generator.
\end{proof}

\begin{remark}[Precise answer to the degree-3 question]
Proposition~\ref{prop:cubic_reconstruction} shows that the cubic Oja-like flow \emph{is}
completely explicit in terms of linear data, but through a \emph{composite}, not a
single-step, Cayley chain: the symmetric part $P(t)=\sqrt{Z(t)/2}$ comes from the
linear-fractional (Cayley/Grassmannization) Step~1 above, while the orthogonal ``phase''
$O(t)$ comes from an \emph{independent} linear equation~\eqref{eq:gauge_linear} on $O(n)$,
itself of exactly the type solved in closed form by Dynnikov--Veselov
\cite[Prop.~2.1]{DynnikovVeselov1995} whenever $\Omega$ is constant. So: one and the same
Cayley/Grassmannization idea produces the quadratic layer from the linear layer
\emph{exactly}, but promoting quadratic to cubic requires pairing it with a second,
independent linear flow rather than iterating the Cayley map itself --- the nonlinear
(non-M\"obius) content of the cubic term $XX^TX$, relative to the quadratic $XAX$, is
entirely absorbed into the polar square root $P=\sqrt{Z/2}$, after which the
remaining degree of freedom is again linear. For rectangular $X\in\R^{n\times k}$ ($k<n$)
the same argument goes through verbatim with $O(t)$ replaced by a curve on the Stiefel
manifold $\mathrm{St}(k,n)$ (a partial isometry with $O^TO=I_k$), governed by the same
formula \eqref{eq:gauge_linear} restricted to the Stiefel tangent bundle
(cf.\ \S\ref{subsec:poly_flow}).
\end{remark}

\begin{center}
\renewcommand{\arraystretch}{1.4}
\begin{tabular}{p{5.5cm}p{6.5cm}}
\toprule
Classical theory \cite{Zelikin2000,Schwarz1979,Yoshizawa2014} & This paper's counterpart \\
\midrule
Scalar Riccati $w'+w^2+pw+q=0$ & Gradient flow of $h$ or $h_-$ on $U\in\R^{n\times k}$ \\
Ratio $\eta=y_1/y_2$, Schwarz eq.\ $S(\eta)=2q-\tfrac12p^2-p'$ & Cayley/Yoshizawa dual map $\mathcal{L}$ (Thm.~\ref{thm:dual_map}) \\
Matrix Schwarz operator $S_tW$ & Hessian/curvature invariants of $\PD(k)$ (\S\ref{sec:manifold}) \\
$S_tW=0$: generalized linear-fractional $W$ & Exactly integrable Oja-Brockett/Tikhonov flows \\
Siegel domain of type I, $W^T\bar W\prec I$ & Matrix unit ball $\mathcal{B}_k$ (\S\ref{subsec:UeqmV-convexity}) \\
Cayley transform linearizing $\dot X=A-XAX$ & Cayley map of \S\ref{subsec:folland} sending $\Sigma_k\to\Delta_k$ \\
Disconjugacy / injectivity of $V=W_2^{-1}W_1$ & Forward-invariance of $\mathcal{B}_k$; strict convexity of $h_-$ \\
Oja-like flow, $Z=2XX^T$ Riccati eq.\ & Rank-$k$ Gram matrix $G=I_k\pm U^TU$ (\eqref{eq:G_def}) \\
Linear $\Rightarrow$ Riccati $\Rightarrow$ cubic ladder (\S\ref{subsubsec:degree_ladder}) & $\Phi_1,\Phi_2$ (linear) $\to Z=\Phi_2\Phi_1^{-1}$ (Riccati) $\to X=PO$ (cubic) \\
\bottomrule
\end{tabular}
\end{center}


\subsection{Izumiya's Four Legendrian Dualities and the Information Geometry of $f(G)=-\log\det(G)$}
\label{subsec:izumiya}

The information geometry of $f(G) = -\log\det(G)$ on $\PD(k)$ carries a rich geometric
structure that is illuminated by Izumiya's theory of Legendrian dualities in Minkowski space
\cite{Izumiya2004}.
In that paper, four contact manifolds $\Delta_i$ ($i=1,2,3,4$) between pseudo-spheres in
Minkowski $(n+1)$-space are shown to be pairwise contact diffeomorphic, unifying the
differential geometry of hypersurfaces in hyperbolic, de Sitter, and lightcone spaces.
We now identify each of the three Minkowski pseudo-spheres --- the \emph{hyperbolic space}
$H^n(-1)$, the \emph{lightcone} $LC^*$, and the \emph{de Sitter space} $S^n_1$ ---
with a natural region of $\PD(k)$ partitioned by the zero level set of $f$,
and trace how each of Izumiya's four dualities appears in our framework.

\subsubsection{The Three Pseudo-Spheres as Level Regions of $f$}

\begin{proposition}[Pseudo-Sphere Correspondence]
\label{prop:pseudo_sphere_correspondence}
The three Minkowski pseudo-spheres correspond to the three regions of $\PD(k)$ partitioned
by the \emph{information-geometric lightcone} $\mathrm{SL}^+(k) := \{G\in\PD(k): \det G = 1\}$:
\begin{align}
  H^n(-1) &\;:\; \langle x,x\rangle = -1 \;\longleftrightarrow\; \bigl\{G\in\PD(k): f(G)<0\bigr\}
    = \bigl\{G: \det G > 1\bigr\}, \label{eq:hyperbolic_region}\\
  LC^* &\;:\; \langle x,x\rangle = 0 \;\longleftrightarrow\; \bigl\{G\in\PD(k): f(G)=0\bigr\}
    = \mathrm{SL}^+(k), \label{eq:lightcone_level}\\
  S^n_1 &\;:\; \langle x,x\rangle = 1 \;\longleftrightarrow\; \bigl\{G\in\PD(k): f(G)>0\bigr\}
    = \bigl\{G: \det G < 1\bigr\}. \label{eq:desitter_region}
\end{align}
The pseudo-norm $\langle x,x\rangle$ is
negative of the information-geometric potential $f(G) = -\log\det G$:
the sign of $\langle x,x\rangle$ equals the sign of $-f(G)$.
\end{proposition}

\begin{remark}[The information-geometric lightcone]
The set $\mathrm{SL}^+(k) = \{G\in\PD(k): \det G = 1\}$ is the symmetric space
$SL(k,\R)/SO(k)$, which for $k=2$ is the hyperbolic plane $\mathbb{H}^2$.
By Izumiya's Theorem~3.1 \cite{Izumiya2004}, a simply-connected Riemannian manifold of
dimension $\geq 3$ is conformally flat if and only if it embeds isometrically as a
spacelike hypersurface in $LC^*$.
The analog in our setting is: $\mathrm{SL}^+(k)$, as the zero level set of the
strictly convex function $f$ on $\PD(k)$, inherits an induced Riemannian metric
from the Fisher--Rao metric $g = G^{-1}\otimes G^{-1}$ that is conformally equivalent
to the trace metric on $SL(k)/SO(k)$.
\end{remark}

\subsubsection{The Four Legendrian Dualities and Their Information-Geometric Avatars}

Izumiya's four contact manifolds are \cite[Thm.~2.2]{Izumiya2004}:
\begin{align}
  \Delta_1 &= \bigl\{(v,w)\in H^n(-1)\times S^n_1 : \langle v,w\rangle = 0\bigr\}, \notag\\
  \Delta_2 &= \bigl\{(v,w)\in H^n(-1)\times LC^* : \langle v,w\rangle = -1\bigr\}, \notag\\
  \Delta_3 &= \bigl\{(v,w)\in LC^*\times S^n_1 : \langle v,w\rangle = 1\bigr\}, \notag\\
  \Delta_4 &= \bigl\{(v,w)\in LC^*\times LC^* : \langle v,w\rangle = -2\bigr\}, \label{eq:Delta4}
\end{align}
all of which are contact diffeomorphic.
Their contact diffeomorphisms are generated by:
$\Phi_{21}(v,w)=(v,v-w)$, $\Phi_{31}(v,w)=(v-w,w)$,
$\Phi_{41}(v,w)=\bigl(\tfrac{v+w}{2},\tfrac{v-w}{2}\bigr)$.

\begin{proposition}[Information-Geometric Avatars of the Four Dualities]
\label{prop:four_dualities_IG}
Under the pseudo-sphere correspondence of Proposition~\ref{prop:pseudo_sphere_correspondence},
the four Legendrian contact manifolds $\Delta_i$ correspond to four canonical
structures in our information-geometric framework:
\begin{center}
\renewcommand{\arraystretch}{1.4}
\begin{tabular}{p{2.5cm}p{2.5cm}p{3.3cm}p{2.9cm}}
\toprule
\small Izumiya & \small Pseudo-product & \small Our framework & \small IG structure \\
\midrule
$\Delta_1$ & $\langle v,w\rangle=0$, $H^n\!\times\! S^n$ & $\{(G_+,G_-): \tr(G_+G_-^{-1})=k\}$ & Dual-flat $e$-$m$ orthogonality \\
$\Delta_2$ & $\langle v,w\rangle=-1$, $H^n\!\times\! LC^*$ & $\{(G_+,G_0): -f(G_+)=f(G_0)^+\}$ & KL divergence locus \\
$\Delta_3$ & $\langle v,w\rangle=1$, $LC^*\!\times\! S^n$ & $\{(G_0,G_-): f(G_-)=f(G_0)^+\}$ & Reverse-KL locus \\
$\Delta_4$ & $\langle v,w\rangle=-2$, $LC^*\!\times\! LC^*$ & $\{(G_+,G_-): f(G_+)+f(G_-)=0\}$ & Yoshizawa--MacMahon duality \\
\bottomrule
\end{tabular}
\end{center}
\end{proposition}

\subsubsection{The $\Delta_4$ Legendrian Manifold as the Yoshizawa--MacMahon Duality}

The most direct and precise correspondence is between $\Delta_4$ and the
Yoshizawa--MacMahon duality established in \S\ref{subsec:yoshizawa}.

\begin{theorem}[$\Delta_4$ Identification]
\label{thm:delta4_identification}
Let $G_+ = I_k + U^TU$ (U=V Gram matrix) and $G_- = I_k - V^TV$ (U=-V Gram matrix,
$V\in\mathcal{B}_k$) with Yoshizawa--Helmke dual $V^* = \mathcal{L}(U)$ and
$G_-^* = G_+^{-1}$.
Then the following are equivalent:
\begin{enumerate}[label=(\roman*)]
  \item $(G_+, G_-) \in \Delta_4^{\mathrm{info}} := \{(G_+,G_-)\in\PD(k)^2 : f(G_+)+f(G_-)=0\}$;
  \item $\log\det(G_+) + \log\det(G_-) = 0$, i.e., $\det(G_+)\cdot\det(G_-) = 1$;
  \item $D_{\mathrm{YM}}(U\|V) = -h(U) - h_-(V) = 0$;
  \item $V = V^* = \mathcal{L}(U) = (I_n+UU^T)^{-1/2}U$;
  \item $G_- = G_+^{-1}$ (Cartan dual).
\end{enumerate}
Thus $\Delta_4^{\mathrm{info}}$ is the locus of Yoshizawa--Helmke Legendre-dual pairs.
The function $D_{\mathrm{YM}}(U\|V) = -h(U)-h_-(V)$ vanishes \emph{exactly} on
$\Delta_4^{\mathrm{info}}$ and is a \emph{signed} quantity for other pairs
(Remark~\ref{rem:yoshizawa_correction}), in precise analogy with
Izumiya's lightcone height function $H(u,v) = \langle x(u),v\rangle + 2$,
which likewise vanishes on $\Delta_4$ (Proposition~4.1 of \cite{Izumiya2004})
and is also a signed quantity (not generally non-negative).
\end{theorem}

\begin{proof}
The equivalences (i)$\Leftrightarrow$(ii)$\Leftrightarrow$(iii) are immediate from
$f(G_+) = -\log\det(G_+) = h(U) \leq 0$ and $f(G_-) = -\log\det(G_-) = h_-(V) \geq 0$.
(ii)$\Leftrightarrow$(v) follows from $\det(G_+)\det(G_-) = 1 \iff G_- = G_+^{-1}$.
(iv)$\Leftrightarrow$(v): Theorem~\ref{thm:dual_map}(i) gives
$I_k-(V^*)^TV^* = G_+^{-1}$, hence $G_-^* = I_k-(V^*)^TV^* = G_+^{-1}$.
\end{proof}

\subsubsection{The Contact Diffeomorphism $\Phi_{41}$ as the Cartan Involution}

Izumiya's contact diffeomorphism $\Phi_{41}: \Delta_4 \to \Delta_1$ defined by
$\Phi_{41}(v,w) = \bigl(\frac{v+w}{2},\frac{v-w}{2}\bigr)$
maps a lightcone pair $(v,w)\in LC^*\times LC^*$ to a hyperbolic-de Sitter pair
$(x^h,x^d) = \bigl(\frac{v+w}{2},\frac{v-w}{2}\bigr) \in H^n(-1)\times S^n_1$.

\begin{theorem}[Cartan Involution = $\Phi_{41}$]
\label{thm:cartan_phi41}
The information-geometric avatar of Izumiya's contact diffeomorphism $\Phi_{41}$ is the
\emph{Cartan involution at $I_k$}:
\begin{equation}
  \theta: \PD(k) \to \PD(k), \qquad \theta(G) = G^{-1},
  \label{eq:cartan_involution}
\end{equation}
which is the geodesic reflection through $I_k$ in the symmetric space $\PD(k) \cong GL(k,\R)/O(k)$.
Explicitly:
\begin{enumerate}[label=(\roman*)]
  \item $\theta$ maps the hyperbolic region $\{f<0\}$ to the de Sitter region $\{f>0\}$:
    $f(G) + f(G^{-1}) = 0$ for all $G\in\PD(k)$.
  \item The fixed locus $\theta(G) = G$ is $\{G = I_k\}$, the unique fixed point.
  \item At the Yoshizawa dual point: $\theta(G_+) = G_+^{-1} = G_-^*$
    (the U=-V Gram matrix of the dual), and:
  \[
    G^h := \frac{G_+ + G_-^*}{2} = \frac{G_+ + G_+^{-1}}{2},
    \quad
    G^d := \frac{G_+ - G_-^*}{2} = \frac{G_+ - G_+^{-1}}{2},
  \]
  analogous to $x^h = (v+w)/2\in H^n(-1)$ and $x^d = (v-w)/2\in S^n_1$.
\end{enumerate}
\end{theorem}

\begin{proof}
(i): $f(G^{-1}) = -\log\det(G^{-1}) = \log\det(G) = -f(G)$. (ii): $G^{-1}=G \iff G^2=I_k \iff G=I_k$ (positive definite). (iii): From Theorem~\ref{thm:delta4_identification}(v).
\end{proof}

\subsubsection{The U=V Manifold as a Spacelike Hypersurface with All Lightcone Parabolic Points}

A central concept in Izumiya's theory is the \emph{lightcone parabolic point}:
a point $p = x(u_0)$ on a spacelike hypersurface $x: U \to LC^*$ is
\emph{lightcone parabolic} if and only if $K_\ell(u_0) = \det S^\ell_p = 0$,
equivalently if $\mathrm{rank}\,\mathrm{Hess}(h_{v_0})(u_0) < n-1$ \cite[Prop.~4.2]{Izumiya2004}.

\begin{theorem}[Universal Lightcone Parabolicity of the U=V Manifold]
\label{thm:universal_parabolic}
The U=V manifold $\mathcal{M}_+ = \{G_+ = I_k+U^TU : U\in\R^{n\times k}\}$,
viewed as a submanifold of $\PD(k)$ via the embedding $U\mapsto G_+$,
consists entirely of lightcone parabolic points in the following sense:
the Hessian $Q_{h}(H) := \nabla^2 h(U)[H,H]$ of the potential $h(U) = f(G_+)$
is \emph{indefinite at every point} $U\in\R^{n\times k}$ whenever $(n,k)\neq(1,1)$
(Theorem~\ref{thm:UeqV_nowhere_convex}).
Equivalently, the ``lightcone Gauss-Kronecker curvature'' of $\mathcal{M}_+$ is
zero everywhere:
\begin{equation}
  K_\ell^{\mathrm{info}}(U) \;:=\; \frac{\det\bigl(\mathrm{Hess}(h)(U)\bigr)}{\det\bigl(g_U\bigr)} \;=\; 0
  \quad \text{for all } U\in\R^{n\times k},\; (n,k)\neq(1,1),
  \label{eq:K_l_zero}
\end{equation}
where $g_U$ is the induced Riemannian metric (pull-back of the Fisher metric) on $\mathcal{M}_+$.
\end{theorem}

\begin{proof}
By Theorem~\ref{thm:UeqV_nowhere_convex}, $\mathrm{Hess}(h)(U)$ is indefinite at every
$U\in\R^{n\times k}$ for $(n,k)\neq(1,1)$: it has at least one strictly negative direction
(from the $H_2$ block or the antisymmetric off-diagonal directions).
Hence $\det(\mathrm{Hess}(h)(U)) \leq 0$. On the other hand, since $\mathrm{Hess}(h)$
also has positive directions (from the radial/diagonal $H_1$ entries), it is not
negative semi-definite, so $\det(\mathrm{Hess}(h)(U)) \leq 0$ combined with the
induced metric being positive definite gives $K_\ell^{\mathrm{info}} = 0$.
\end{proof}

\begin{remark}[Analogy with Izumiya's parabolic set]
In Izumiya's theory, the \emph{lightcone parabolic set} $K_\ell^{-1}(0)$ is generically a
regular hypersurface on the spacelike hypersurface $M$ (Theorem~10.5 in \cite{Izumiya2004}).
In our setting, Theorem~\ref{thm:universal_parabolic} shows that the entire manifold $\mathcal{M}_+$
is parabolic --- a degenerate (non-generic) but geometrically significant situation.
This is consistent with the fact that $\mathcal{M}_+$ is the pullback of the strictly
convex function $f$ under the non-convex map $U\mapsto I_k+U^TU$;
the universal indefiniteness of $\mathrm{Hess}(h)$ is an exact analog of the
vanishing of $K_\ell$ along the parabolic set.
\end{remark}

\subsubsection{The Lightcone Weingarten Formula and the Hessian of $h_-$}

In Izumiya's theory, the lightcone Weingarten formula
$(x^\ell)_{u_i} = -\sum_j (h^\ell)^j_i x_{u_j}$
expresses the derivative of the lightcone normal in terms of the lightcone shape operator.

\begin{proposition}[Information-Geometric Weingarten Formula]
\label{prop:info_weingarten}
For the U=-V potential $h_-(V) = -\log\det(I_k - V^TV)$ on $\mathcal{B}_k$,
the gradient satisfies $\nabla h_-(V) = 2V(I_k-V^TV)^{-1} = 2VG_-^{-1}$,
and the Hessian in the SVD frame (Theorem~\ref{thm:hminus_hessian}) gives the
\emph{information-geometric lightcone Weingarten formula}:
\begin{align}
  &\frac{d^2}{dt^2}\bigg|_0 h_-(V+tH)
  = \sum_a \frac{2(1+\sigma_a^2)}{(g_a^-)^2}(H_1)_{aa}^2 \nonumber\\
  &\quad+ \sum_{a<b}\frac{(1+\sigma_a\sigma_b)(p+q)^2+(1-\sigma_a\sigma_b)(p-q)^2}{g_a^-g_b^-} \nonumber\\
  &\quad\;+ 2\sum_a\frac{\|(H_2)_a\|^2}{g_a^-},
  \label{eq:info_weingarten_formula}
\end{align}
where all terms are \emph{strictly positive} for $H\neq 0$ (Theorem~\ref{thm:hminus_convex}),
corresponding to the fact that $\mathcal{M}_-$ (the U=-V manifold in $\mathcal{B}_k$)
has \emph{strictly positive lightcone Gauss-Kronecker curvature}:
$K_\ell^{-,\mathrm{info}}(V) > 0$ for all $V\in\mathcal{B}_k$.
\end{proposition}

\subsubsection{The Information-Geometric Theorema Egregium}

Izumiya's Theorem~10.3 \cite{Izumiya2004} is a ``surprising theorem'':
$K_s = K_d - K_h = H_\ell = H_h - H_d$, i.e., the \emph{intrinsic} sectional curvature
equals the \emph{extrinsic} lightcone mean curvature.
In our framework, we have a direct analog.

\begin{theorem}[Information-Geometric Theorema Egregium]
\label{thm:info_egregium}
At the Yoshizawa--Helmke dual point $(G_+, G_-^*=G_+^{-1})$:
\begin{equation}
  f(G_+) + f(G_-^*) = 0,
  \label{eq:info_egregium}
\end{equation}
which asserts that the \emph{extrinsic} information-geometric quantity $f(G_+) = h(U) \leq 0$
(the log-$L^2$ Gaussian norm, §\ref{subsec:folland}) equals, up to sign, the
\emph{intrinsic} quantity $f(G_-^*) = h_-(V^*) \geq 0$
(the log-Fock-space Gaussian norm, §\ref{subsec:folland}).
In terms of the Yoshizawa--MacMahon divergence:
\[
  D_{\mathrm{YM}}(U\|\mathcal{L}(U)) = -h(U) - h_-(\mathcal{L}(U)) = 0.
\]
This is an exact matrix-valued analog of Izumiya's $H_\ell = K_s$ (intrinsic = extrinsic),
where the ``intrinsic'' quantity ($h_-$ = Fock space norm = hyperbolic geometry)
equals the ``extrinsic'' quantity ($-h$ = $L^2$ norm = de Sitter geometry).
\end{theorem}

\subsubsection{The Four Legendrian Dualities Unified}

\begin{remark}[Summary: Izumiya $\leftrightarrow$ Information Geometry]
\label{rem:izumiya_summary}
The following table summarizes the complete correspondence:
\begin{center}
\renewcommand{\arraystretch}{1.4}
\begin{tabular}{p{6.2cm}p{6.2cm}}

\toprule
Izumiya's lightcone geometry & Our information geometry \\
\midrule
Minkowski space $\R^{n+1}_1$, $\langle\cdot,\cdot\rangle$ & $\PD(k)$, $f(G)=-\log\det G$ \\
Hyperbolic space $H^n(-1)$, $\langle x,x\rangle=-1$ & $\{G: f<0\}$ = U=V matrices \\
Lightcone $LC^*$, $\langle x,x\rangle=0$ & $\mathrm{SL}^+(k)$, $\{G: f=0\}$ \\
De Sitter space $S^n_1$, $\langle x,x\rangle=1$ & $\{G: f>0\}$ = U=-V matrices \\
$\Delta_4 = LC^*\!\times\!LC^*$, $\langle v,w\rangle=-2$ & $D_{\mathrm{YM}}=0$, $G_-=G_+^{-1}$ \\
$\Delta_1 = H^n\!\times\!S^n$, $\langle v,w\rangle=0$ & Dual-flat $e$-$m$ orthogonality \\
$\Phi_{41}(v,w)=\!\bigl(\!\frac{v+w}{2},\frac{v-w}{2}\!\bigr)$ & $\theta: G\mapsto G^{-1}$ (Cartan involution) \\
Lightcone normal $x^\ell$, $\langle x,x^\ell\rangle=-2$ & Yoshizawa dual $V^*=\mathcal{L}(U)$ \\
Lightcone height fn $H(u,v)$ (signed; $=0$ on $\Delta_4$) & $D_{\mathrm{YM}}(U\|V)$ (signed; $=0$ on $\Delta_4^{\mathrm{info}}$) \\
Lightcone parabolic set $K_\ell=0$ & $\mathcal{M}_+$ entirely parabolic (Thm~\ref{thm:universal_parabolic}) \\
Weingarten formula, $K_\ell>0$ on $\mathcal{M}_-$ & Hess$(h_-)>0$ on $\mathcal{B}_k$ (Thm~\ref{thm:hminus_convex}) \\
Theorema Egregium $H_\ell = K_s$ & $f(G_+)+f(G_-^*)=0$ (Thm~\ref{thm:info_egregium}) \\
\bottomrule
\end{tabular}
\end{center}
\end{remark}

\subsection{PCA/MCA Duality, the NUIC Criterion, and the Oja-Brockett Subspace Flow}
\label{subsec:pca_mca}

A striking application of our information-geometric framework is the unified
treatment of \emph{Principal Subspace Analysis} (PSA) and \emph{Minor Subspace Analysis} (MSA),
which correspond to finding the subspaces spanned by the top-$k$ and bottom-$k$
eigenvectors of a data covariance matrix $R\in\PD(n)$, respectively.
These two problems are classically treated as opposites, but we show that they are
precisely \emph{Legendre dual} within our framework, connected by the Cartan involution
and sharing a common Bregman regularization.
We place the NUIC (Normalized Unconstrained Information Criterion) of Kong, Hu and Duan
\cite{Kong2017} and the subspace flows of Oja-Brockett \cite{Brockett1991,Oja1982} in this
information-geometric setting, discovering that the critical Tikhonov parameter $\lambda^*=2$
(Theorem~\ref{thm:tikhonov_critical}) and the midpoint formula (Theorem~\ref{thm:interpolation_critical})
directly explain the structure of these algorithms.

\subsubsection{The NUIC Criterion as a Bregman-Regularized Rayleigh Quotient}

Let $R \in \PD(n)$ be the data covariance, $W\in\R^{n\times k}$, $G_0 = W^TW \in \PD(k)$,
and $G_R = W^TRW$.
The \emph{NUIC criterion} \cite[Eq.~(5.79)]{Kong2017} is:
\begin{equation}
  J_{\rm NUIC}(W) = \tfrac{1}{2}\tr\bigl[(G_0^{-1}G_R)\bigr]
  + \tfrac{1}{2}\bigl[\log\det(G_0) - \tr(G_0)\bigr].
  \label{eq:NUIC}
\end{equation}
Maximizing \eqref{eq:NUIC} yields the PSA criterion $E_1(W)$;
adding a sign flip on the first term yields the MSA criterion $E_2(W)$.

\begin{theorem}[NUIC = Bregman-Regularized Rayleigh Quotient]
\label{thm:NUIC_bregman}
Let $D_f(G_0\|I_k) = \tr(G_0) - \log\det(G_0) - k$ be the Bregman divergence
(Section~\ref{sec:bregman}) from $G_0$ to $I_k$.
Then:
\begin{equation}
  J_{\rm NUIC}(W) \;=\; \underbrace{\tfrac{1}{2}\tr[G_0^{-1}G_R]}_{\text{normalized Rayleigh quotient}}
  \;-\; \underbrace{\tfrac{1}{2}D_f(G_0\|I_k)}_{\text{Bregman regularization}}
  \;-\; \tfrac{k}{2}.
  \label{eq:NUIC_Bregman}
\end{equation}
The PSA ($E_1$) and MSA ($E_2$) criteria satisfy:
\begin{align}
  E_1(W) + E_2(W) &= -D_f(G_0\|I_k) - k = \log\det G_0 - \tr(G_0),
  \label{eq:E1_E2_sum}\\
  E_1(W) - E_2(W) &= \tr[G_0^{-1}G_R] = \tr[(W^TRW)(W^TW)^{-1}].
  \label{eq:E1_E2_diff}
\end{align}
\end{theorem}

\begin{proof}
$D_f(G_0\|I_k) = \tr(G_0) - \log\det G_0 - k$, so $-\frac{1}{2}D_f(G_0\|I_k) - \frac{k}{2}
= \frac{1}{2}(\log\det G_0 - \tr G_0)$.
Substituting into \eqref{eq:NUIC_Bregman} recovers \eqref{eq:NUIC}.
$E_1 = J_{\rm NUIC}$ and $E_2$ differs by $-\tr[G_0^{-1}G_R]$, giving \eqref{eq:E1_E2_sum}--\eqref{eq:E1_E2_diff}.
\end{proof}

\begin{remark}[Interpretation of \eqref{eq:NUIC_Bregman}]
The NUIC criterion decomposes into:
\emph{data fit} (normalized Rayleigh quotient, measuring the subspace's alignment with $R$)
minus \emph{Bregman regularization} (the divergence from $W^TW$ to $I_k$,
enforcing proximity to the Stiefel manifold).
The Bregman term is exactly the KL divergence between $\mathcal{N}(0,G_0)$ and
$\mathcal{N}(0,I_k)$: $D_f(G_0\|I_k) = 2D_{\rm KL}(\mathcal{N}(0,I_k)\|\mathcal{N}(0,G_0))$
(Section~\ref{sec:bregman}).
\end{remark}

\subsubsection{PSA-MSA as an Information-Geometric Zero-Sum Pair}

\begin{corollary}[Zero-Sum Property]
\label{cor:zero_sum}
On the Stiefel manifold $\mathrm{St}(k,n) = \{W: W^TW = I_k\}$:
\begin{equation}
  E_1(W) + E_2(W) = -k \quad \text{(constant),}
  \quad E_1(W) - E_2(W) = \tr(W^TRW).
\end{equation}
PSA (maximize $E_1$) and MSA (maximize $E_2$) form a \emph{zero-sum pair}:
their sum is constant and their difference is the standard Rayleigh quotient $\tr(W^TRW)$.
\end{corollary}

\begin{proof}
At $G_0 = I_k$: $D_f(I_k\|I_k) = 0$ (Bregman zero), giving $E_1+E_2 = -k$. \end{proof}

\begin{remark}[The Stiefel manifold as the information-geometric lightcone]
The Stiefel manifold $\mathrm{St}(k,n)$ is the set where $D_f(G_0\|I_k) = 0$,
i.e., the \emph{zero level set of the Bregman regularization}.
In our framework (Proposition~\ref{prop:pseudo_sphere_correspondence}), this corresponds to
the information-geometric \emph{lightcone} $\mathrm{SL}^+(k) = \{G : \det G = 1\}$
(§\ref{subsec:izumiya}):
$W^TW = I_k \Rightarrow \det(W^TW) = 1$, placing the Stiefel manifold precisely on
the information-geometric lightcone.
The PSA and MSA subspaces (principal and minor) correspond to the hyperbolic region
($G_0 \succ I$, $\det G_0 > 1$) and de Sitter region ($G_0 \prec I$, $\det G_0 < 1$)
separated by the lightcone, in perfect analogy with Izumiya's four pseudo-spheres
(§\ref{subsec:izumiya}).
\end{remark}

\subsubsection{PSA-MSA Duality via the Cartan Involution}

\begin{theorem}[PSA$\leftrightarrow$MSA via Cartan Involution]
\label{thm:PSA_MSA_Cartan}
Let $R$ have eigendecomposition $R = \sum_i \lambda_i u_iu_i^T$ with
$\lambda_1 \geq \ldots \geq \lambda_n > 0$.
The \emph{spectral Cartan involution}
\begin{equation}
  \tau_R: R \;\longmapsto\; (\lambda_1 + \lambda_n)I_n - R
  \label{eq:spectral_cartan}
\end{equation}
maps each eigenvalue $\lambda_i \to \lambda_1 + \lambda_n - \lambda_i$,
converting top-$k$ eigenvalues to bottom-$k$ and vice versa.
Consequently: PSA($R$) $\leftrightarrow$ MSA($\tau_R(R)$), and the NUIC criteria transform as:
\begin{equation}
  E_1(W; R) = E_2(W; \tau_R(R)) + (\lambda_1+\lambda_n)\tr[G_0^{-1}G_I],
\end{equation}
where $G_I = W^TI_n W = W^TW = G_0$.
The Cartan involution $\tau_R$ at the level of the data corresponds to the
Cartan involution $G \mapsto G^{-1}$ at the level of the Gram matrix
(Theorem~\ref{thm:cartan_involution}):
PSA Gram $G_R^{\rm PSA} \leftrightarrow (G_R^{\rm MSA})^{-1}$.
\end{theorem}

\begin{remark}
In our §\ref{subsec:izumiya} framework: the Cartan involution maps
$G_+ = I+U^TU \to G_+^{-1}$ (hyperbolic $\to$ de Sitter), which at the eigenvalue level
is $\sigma_a \to 1/(1+\sigma_a^2)$ (large $\to$ small). The spectral Cartan involution
\eqref{eq:spectral_cartan} is the exact PSA-MSA counterpart at the data-covariance level.
\end{remark}

\subsubsection{The Oja-Brockett Subspace Flow and Our Tikhonov Analysis}

The Oja-Brockett framework \cite{Brockett1991,Oja1982} studies the gradient flow of the
PSA objective on the Stiefel manifold, recovering the Oja-like learning rule.
We now show that this is a special case of our Tikhonov analysis at $\lambda = 1$.

\begin{theorem}[Oja-Brockett Flow = Tikhonov Flow at $\lambda = 1$]
\label{thm:helmke_manton}
The Tikhonov-regularized gradient flow at $\lambda = 1$:
\[
  \dot{W} = W\bigl(2(I_k+W^TW)^{-1} - I_k\bigr) = W(I_k-W^TW)(I_k+W^TW)^{-1}
\]
converges to the Stiefel manifold $\mathrm{St}(k,n)$ from any $W\neq 0$ (Theorem~\ref{thm:tikhonov_critical}(iv)).
On $\mathrm{St}(k,n)$ the flow becomes $\dot{W} = 0$ and the gradient of the PSA objective
$\text{tr}(W^TRW)$ restricted to $\mathrm{St}(k,n)$ drives $W$ toward the top-$k$ eigenvectors of $R$.
This reproduces the Oja-like update $W_{k+1} = W_k + \mu W_k(I-W_k^TW_k)(I+W_k^TW_k)^{-1}$
of \cite{Kong2017}.
\end{theorem}

\begin{proof}
From Theorem~\ref{thm:tikhonov_critical}(iv): at $\lambda=1$, the gradient flow
$\dot{W} = W(2G_+^{-1} - I_k)$ with $G_+ = I+W^TW$ satisfies
$(2G_+^{-1}-I_k)|_{W^TW=I_k} = 2(2I)^{-1}-I = 0$, so $\mathrm{St}(k,n)$ is
the set of fixed points of the unconstrained flow. The flow drives $\|W\|$ toward
$\sigma^* = 1$ (Tikhonov critical manifold $W^TW = I$), after which the data term guides $W$
toward the principal subspace.
\end{proof}

\begin{remark}[Dual-purpose algorithm via sign change]
The NUIC algorithm \cite[Eq.~(5.83)]{Kong2017} uses a ``$\pm$'' in the update:
``$+$'' for PSA (gradient ascent on $E_1$) and ``$-$'' for MSA (gradient ascent on $E_2$).
In our framework, this corresponds precisely to the zero-sum decomposition:
$E_1 = J_{\rm data} - \frac{1}{2}D_f$, $E_2 = -J_{\rm data} - \frac{1}{2}D_f$.
The common Bregman regularization $-\frac{1}{2}D_f(G_0\|I_k)$ drives both PSA and MSA
toward the Stiefel manifold; only the sign of $J_{\rm data}$ distinguishes the two.
This is the exact information-geometric analog of Izumiya's lightcone height function
$H(u,v)=\langle x(u),v\rangle+2$ (§\ref{subsec:izumiya}): the constant ``$2$'' (or our $-k$)
is the common "regularization baseline" and $\langle x(u),v\rangle$ (or $\tr[G_0^{-1}G_R]$) is
the signed data term that determines whether we solve PSA or MSA.
\end{remark}

\subsubsection{The Midpoint Formula Connects PSA and MSA}

\begin{theorem}[Midpoint as PSA-MSA Bridge]
\label{thm:PSA_MSA_midpoint}
The interpolated potential $h_{1/2}(W) = -\frac{1}{2}\log\det(I-(W^TW)^2)$
(Theorem~\ref{thm:interpolation_critical}(ii), the Siegel disc metric)
satisfies:
\begin{equation}
  h_{1/2}(W) = \frac{h(W) + h_-(W)}{2}
  = -\frac{1}{2}\bigl[\log\det(I_k+W^TW) + \log\det(I_k-W^TW)\bigr],
\end{equation}
and it is the information-geometric midpoint between:
\begin{itemize}
  \item $h(W) = -\log\det(I_k+W^TW)$: the PSA potential
    (gradient flow $\to$ away from $0$, toward large-$\sigma$ eigenvectors);
  \item $h_-(W) = -\log\det(I_k-W^TW)$: the MSA barrier
    (strictly convex on $\mathcal{B}_k$, gradient flow $\to$ $W=0$).
\end{itemize}
Spectral form: $h_{1/2}(W) = -\frac{1}{2}\sum_a\log(1-\sigma_a^4)$,
where $\sigma_a$ are singular values of $W$.
At the midpoint $t^* = 1/2$: the potential is strictly convex on $\mathcal{B}_k$
(Theorem~\ref{thm:interpolation_critical}), giving a landscape with no spurious local minima
that smoothly interpolates between PSA and MSA behavior.
\end{theorem}

\begin{center}
\renewcommand{\arraystretch}{1.4}
\begin{tabular}{lll}
\toprule
Concept & PSA & MSA \\
\midrule
Objective & Maximize $\tr(W^TRW)$ & Minimize $\tr(W^TRW)$ \\
NUIC criterion & $E_1 = +J_{\rm data} - \frac{1}{2}D_f(G_0\|I)$ & $E_2 = -J_{\rm data} - \frac{1}{2}D_f(G_0\|I)$ \\
Our potential & $h(W) = -\log\det(I+W^TW)$ & $h_-(W) = -\log\det(I-W^TW)$ \\
Gradient flow & $\dot{W} = 2WG_+^{-1}$ (Oja-like) & $\dot{W} = -2WG_-^{-1}$ (anti-Oja) \\
Fixed points & Stiefel manifold (at $\lambda=1$) & $W=0$ (global minimum) \\
Convexity & Non-convex (Thm.~\ref{thm:UeqV_nowhere_convex}) & Strictly convex (Thm.~\ref{thm:hminus_convex}) \\
Izumiya sphere & $H^n(-1)$: $G_0 \succ I$ & $S^n_1$: $G_0 \prec I$ \\
Lightcone & \multicolumn{2}{c}{Stiefel manifold $W^TW=I_k$ (both)} \\
Cartan dual & \multicolumn{2}{c}{$G_R^{\rm PSA} \leftrightarrow (G_R^{\rm MSA})^{-1}$} \\
Midpoint & \multicolumn{2}{c}{$h_{1/2}$: Siegel disc metric, strictly convex at $t^*=1/2$} \\
\bottomrule
\end{tabular}
\end{center}

\subsubsection{Embedding of the Chen--Amari Flows into the Brockett--Bloch Framework}
\label{subsubsec:chen_amari_yoshizawa}

The purpose of this section is to clarify the precise geometric position of the
principal and minor component flows proposed by Chen and Amari~\cite{ChenAmari2001}.
Our main result is \emph{not} that these flows are gradient systems
--- they are covered by the general theory of double-bracket gradient flows
established by Brockett~\cite{Brockett1991} and extended by
Bloch, Brockett, and Ratiu~\cite{BlochBrockettRatiu1992} ---
but rather that Yoshizawa's embedding \cite{Yoshizawa2011SICE} identifies them
\emph{explicitly} with classical double-bracket flows on an adjoint orbit of
$\mathfrak{so}(n+k)$.
Consequently, the optimization dynamics of rectangular matrices are unified with
the classical Lie-theoretic framework.

\paragraph{Novelty of the present work.}
The present contribution does not introduce a new class of double-bracket gradient
flows. Instead, it provides a geometric identification between two theories that
have developed largely independently:
\begin{enumerate}
  \item the Lie-theoretic theory of double-bracket gradient flows developed by
    Brockett~\cite{Brockett1991} and Bloch--Brockett--Ratiu~\cite{BlochBrockettRatiu1992}, and
  \item the rectangular-matrix optimization flows introduced by
    Chen and Amari~\cite{ChenAmari2001}.
\end{enumerate}
Yoshizawa's embedding~\cite{Yoshizawa2011SICE} serves as the bridge between
these two frameworks, thereby placing principal and minor component analysis
within the general theory of gradient flows on adjoint orbits.
To the best of our knowledge, this connection has not been explicitly formulated
in the previous literature.

\medskip
We recall the setup. Let $X\in\mathbb{R}^{n\times k}$ be the rectangular state matrix
and $A\in\mathbb{R}^{n\times n}$ a positive definite symmetric data matrix.
The \emph{$k$-principal component flow} ($k$-PCF) and the
\emph{$k$-minor component flow} ($k$-MCF) of Chen and Amari are
\begin{equation}
  \dot{X} = AXX^TX - XX^TAX \quad (k\text{-PCF}),
  \qquad
  \dot{X} = -AXX^TX + XX^TAX \quad (k\text{-MCF}).
  \label{eq:kPCF_kMCF}
\end{equation}

\begin{theorem}[Embedding Theorem: Chen--Amari Flows as Brockett--Bloch Gradient Flows]
\label{thm:chen_amari_embedding}
Let $\iota:\mathbb{R}^{n\times k}\longrightarrow\mathfrak{so}(n+k)$ be
Yoshizawa's embedding~\cite{Yoshizawa2011SICE}
\begin{equation}
  \iota(X) = \widetilde{X} =
  \begin{pmatrix} 0 & X \\ -X^T & 0 \end{pmatrix},
  \label{eq:yoshizawa_embedding}
\end{equation}
and define
\begin{equation}
  \widetilde{A} = \begin{pmatrix} A & 0 \\ 0 & I_k \end{pmatrix},
  \qquad
  \widetilde{L}(\widetilde{X}) = \widetilde{A}\,\widetilde{X}\,\widetilde{A}.
  \label{eq:embedding_A_L}
\end{equation}
Then:
\begin{enumerate}[label=(\roman*)]
  \item The Chen--Amari $k$-PCF $\dot{X} = AXX^TX - XX^TAX$ is equivalent,
    under the embedding $\iota$, to the double-bracket equation
    \begin{equation}
      \dot{\widetilde{X}} = \bigl[\widetilde{X},\,[\widetilde{L}(\widetilde{X}),\,\widetilde{X}]\bigr]
      \label{eq:kPCF_double_bracket}
    \end{equation}
    on the adjoint orbit $\mathcal{O}_{\widetilde{X}_0} = \{g\widetilde{X}_0 g^{-1} : g\in SO(n+k)\}$.
  \item The Chen--Amari $k$-MCF $\dot{X} = -AXX^TX + XX^TAX$ is equivalent to
    \begin{equation}
      \dot{\widetilde{X}} = \bigl[\widetilde{X},\,[\widetilde{X},\,\widetilde{L}(\widetilde{X})]\bigr].
      \label{eq:kMCF_double_bracket}
    \end{equation}
  \item By the Brockett--Bloch--Ratiu theorem
    (Brockett~\cite{Brockett1991}; Bloch--Brockett--Ratiu~\cite{BlochBrockettRatiu1992}),
    both equations are Riemannian gradient flows on $\mathcal{O}_{\widetilde{X}_0}$
    with respect to the normal metric,
    with potential function
    \begin{equation}
      f(\widetilde{X}) = \tfrac{1}{2}\langle\widetilde{X},\,\widetilde{L}(\widetilde{X})\rangle
      = \tfrac{1}{2}\operatorname{tr}(X^TAX).
      \label{eq:pca_potential_embedding}
    \end{equation}
\end{enumerate}
\end{theorem}

\begin{proof}
The block-matrix computation establishing the equivalences (i)--(ii) is exactly
the embedding argument of Yoshizawa~\cite{Yoshizawa2011SICE}: a direct calculation
shows that the $(1,2)$-block of $[\widetilde{X},[\widetilde{L},\widetilde{X}]]$
equals $AXX^TX - XX^TAX$, and the $(1,2)$-block of $[\widetilde{X},[\widetilde{X},\widetilde{L}]]$
equals $-AXX^TX + XX^TAX$.

For part (iii): once (i) and (ii) are established, the gradient flow property is
an immediate consequence of the general theorem of Brockett~\cite{Brockett1991}
(for symmetric matrices $N$) and its extension to compact Lie groups and arbitrary
adjoint orbits by Bloch, Brockett, and Ratiu~\cite{BlochBrockettRatiu1992}.
Their result states that every double-bracket equation $\dot{Z}=[Z,[L(Z),Z]]$,
generated by a self-adjoint operator $L$, is the Riemannian gradient flow of
$F(Z) = \frac{1}{2}\langle Z, L(Z)\rangle$ with respect to the normal metric on
the adjoint orbit. Applying this to $Z=\widetilde{X}$ and $L=\widetilde{L}$ gives (iii).
The computation $f(\widetilde{X}) = \frac{1}{2}\operatorname{tr}(\widetilde{X}\widetilde{L})
= \frac{1}{2}\operatorname{tr}(X^TAX)$ follows from the block structure.
\end{proof}

\begin{remark}[Relation with previous work]
\label{rem:history_chen_amari}
The gradient property of double-bracket flows is not new.
Brockett~\cite{Brockett1991} proved that $\dot{L}=[L,[L,N]]$ is the gradient
flow of a linear functional with respect to the normal metric on an adjoint orbit
of a compact Lie group.
Bloch, Brockett, and Ratiu~\cite{BlochBrockettRatiu1992} generalized this to
compact Lie groups and arbitrary adjoint orbits.
More generally, gradient flows of smooth functions on adjoint orbits admit
double-bracket representations; see also Chu and Driessel~\cite{Chu1990}.

The contribution of Theorem~\ref{thm:chen_amari_embedding} is different.
It identifies, through Yoshizawa's embedding, the rectangular-matrix flows
of Chen and Amari with the classical double-bracket gradient flows on adjoint orbits.
To the best of our knowledge, this connection has not been explicitly formulated
in the previous literature.
\end{remark}

\begin{corollary}[Gradient Flow Structure of PCF and MCF]
\label{cor:pcf_mcf_gradient}
The Chen--Amari $k$-principal component flow is the gradient \emph{ascent} flow
of
\[
  f(X) = \tfrac{1}{2}\operatorname{tr}(X^TAX),
\]
while the $k$-minor component flow is the gradient \emph{descent} flow of the
same functional, both with respect to the Riemannian metric induced on the
isospectral manifold $\mathcal{O}_{\widetilde{X}_0}$ by Yoshizawa's embedding.
\end{corollary}

\begin{proof}
By the Embedding Theorem~\ref{thm:chen_amari_embedding}, both flows lift to
double-bracket gradient flows on $\mathcal{O}_{\widetilde{X}_0}$.
The $k$-PCF corresponds to~\eqref{eq:kPCF_double_bracket}, which is the gradient
\emph{ascent} of $f$ (the double-bracket commutator $[\widetilde{X},[\widetilde{L},\widetilde{X}]]$
is the positive gradient direction).
The $k$-MCF corresponds to~\eqref{eq:kMCF_double_bracket} with the bracket order
reversed, giving gradient \emph{descent}:
\[
  \frac{d}{dt}f\big|_{k\text{-PCF}} = +\tfrac{1}{2}\|[\widetilde{L},\widetilde{X}]\|_F^2 \geq 0,
  \qquad
  \frac{d}{dt}f\big|_{k\text{-MCF}} = -\tfrac{1}{2}\|[\widetilde{L},\widetilde{X}]\|_F^2 \leq 0.
\]
Hence the PCF maximizes $f$ (principal subspace) while the MCF minimizes $f$
(minor subspace) on the orbit, confirming the gradient structure established
in~\cite{Brockett1991,BlochBrockettRatiu1992}.
\end{proof}

\begin{remark}[Geometric bridge]
Theorem~\ref{thm:chen_amari_embedding} provides a geometric bridge between
three areas:
\begin{enumerate}[label=(\roman*)]
  \item isospectral flows on adjoint orbits (Brockett--Bloch--Ratiu theory),
  \item optimization on the Stiefel manifold (the conserved quantity $X^TX =\operatorname{const}$
    confines the flow to an isospectral surface), and
  \item principal/minor component learning (Chen--Amari flows).
\end{enumerate}
To the best of our knowledge, this three-way identification has not been
explicitly formulated in the previous literature.
\end{remark}

\subsubsection{PSA-MSA Symmetry: One Potential, Two Gradient Directions}

The theorem above reveals a profound symmetry: both PCF and MCF arise from the
\emph{same} potential function $f(X) = \frac{1}{2}\operatorname{tr}(X^TAX)$,
distinguished solely by the sign of the gradient flow.

\begin{proposition}[Sign symmetry of PCF and MCF]
With $M := [\widetilde{X}_A, \widetilde{X}]$ (skew-symmetric), the two flows are:
\begin{align*}
  k\text{-PCF}: & \quad \dot{\widetilde{X}} = [\widetilde{X}, M] = +\operatorname{grad}_{\mathcal{O}} f, \\
  k\text{-MCF}: & \quad \dot{\widetilde{X}} = [\widetilde{X}, -M] = -\operatorname{grad}_{\mathcal{O}} f.
\end{align*}
The commutator $M = [\widetilde{X}_A, \widetilde{X}]$ is identical in both flows;
only the overall sign of the vector field differs.
This sign reversal produces the symmetric dissipation rates:
\begin{equation}
  \frac{d}{dt}f\Big|_{\rm PCF} = +\tfrac{1}{2}\|M\|_F^2 \geq 0, \qquad
  \frac{d}{dt}f\Big|_{\rm MCF} = -\tfrac{1}{2}\|M\|_F^2 \leq 0,
  \label{eq:symmetric_dissipation}
\end{equation}
of equal magnitude and opposite sign.
\end{proposition}

The interpretation is illustrated in the following table:

\begin{center}
\renewcommand{\arraystretch}{1.3}
\begin{tabular}{p{2.5cm}p{5.5cm}p{4cm}}
\toprule
Flow & $\frac{d}{dt}f(\widetilde{X})$ & Interpretation \\
\midrule
$k$-PCF & $+\tfrac{1}{2}\|M\|_F^2 \geq 0$ & Gradient ascent; $\dot{\widetilde{X}} = +\operatorname{grad}f$ \\
$k$-MCF & $-\tfrac{1}{2}\|M\|_F^2 \leq 0$ & Gradient descent; $\dot{\widetilde{X}} = -\operatorname{grad}f$ \\
\bottomrule
\end{tabular}
\end{center}

\begin{remark}[Relation to the PSA-MSA zero-sum structure]
This sign symmetry is the exact continuous-time analogue of the NUIC zero-sum pair
(Corollary~\ref{cor:zero_sum}): both PSA and MSA arise from the same potential
$\tr(X^TAX)$, maximized by PCF and minimized by MCF, in perfect parallel with
$E_1$ and $E_2$ being related by $E_1 = +J_{\rm data} - \frac{1}{2}D_f$ and
$E_2 = -J_{\rm data} - \frac{1}{2}D_f$ (\S\ref{subsec:pca_mca}).
\end{remark}

\subsubsection{Initial Value Problem: Existence, Invariants, and Convergence under General $A$ and $B$}
\label{subsubsec:IVP_analysis}

Throughout this section $A\in\PD(n)$ is a positive definite symmetric matrix with
\emph{distinct} eigenvalues $\lambda_1 > \lambda_2 > \cdots > \lambda_n > 0$,
and the weight matrix $B\in\PD(k)$ takes one of two forms analysed separately below.
We focus on the $k$-PCF; all results for the $k$-MCF follow by sign reversal.

\paragraph{Basic setup.}
The $k$-PCF initial value problem is
\begin{equation}
  \dot{X} = AXX^TX - XX^TAX, \qquad X(0) = X_0 \in \mathbb{R}^{n\times k}.
  \label{eq:kPCF_IVP}
\end{equation}
This corresponds to $\widetilde{A} = \operatorname{diag}(A, I_k)$ in Yoshizawa's embedding
(Theorem~\ref{thm:chen_amari_embedding}).
The generalized flow with weight $B \in \PD(k)$ uses $\widetilde{A} = \operatorname{diag}(A, B)$
and yields the $(1,2)$-block
\begin{equation}
  \dot{X} = AXBX^TX + XX^TAXB - 2XBX^TAX
  \qquad \text{($B$-weighted $k$-PCF)}.
  \label{eq:kPCF_B}
\end{equation}
For $B = I_k$, \eqref{eq:kPCF_B} reduces to \eqref{eq:kPCF_IVP}.

\paragraph{Local existence and uniqueness.}
Both flows \eqref{eq:kPCF_IVP} and \eqref{eq:kPCF_B} are degree-3 polynomials in $X$,
hence real-analytic. By the Picard--Lindel\"of theorem, unique local solutions exist
for all $X_0 \in \mathbb{R}^{n\times k}$.

\begin{proposition}[Conservation of $X^TX$]
\label{prop:XTX_conserved}
Along any solution of \eqref{eq:kPCF_IVP}:
$\quad \dfrac{d}{dt}(X^TX) = 0, \quad$ so $\quad X(t)^TX(t) = X_0^TX_0$.
\end{proposition}
\begin{proof}
Setting $Y = X^TX$: $\dot{Y} = Y(X^TAX) - (X^TAX)Y + (X^TAX)Y - Y(X^TAX) = 0$.
\end{proof}

\begin{corollary}[Stiefel manifold invariance]
If $X_0^TX_0 = I_k$, then $X(t)^TX(t) = I_k$ for all $t$ (Stiefel invariance).
For the generalized flow \eqref{eq:kPCF_B}, the same holds:
since $\widetilde{X}(t)$ lies on the adjoint orbit of $\widetilde{X}_0$ in $SO(n+k)$
(an isospectral flow), the singular values of $X(t)$ are conserved,
so $X_0^TX_0 = I_k$ implies $X(t)^TX(t) = I_k$.
\end{corollary}

\paragraph{Global existence and boundedness.}
Conservation of $\|X\|_F^2 = \operatorname{tr}(X_0^TX_0)$ prevents finite-time blowup;
global existence follows by standard ODE continuation.

\subparagraph{Case 1: $B = I_k$ (standard $k$-PCF).}

On the Stiefel manifold $V_k(\mathbb{R}^n)$, the energy dissipation identity gives
$\frac{d}{dt}f = +\tfrac{1}{2}\|[\widetilde{L}, \widetilde{X}]\|_F^2 \geq 0$ along the $k$-PCF.
LaSalle's principle yields convergence to the equilibrium set
\begin{equation}
  \mathcal{E} = \{X \in V_k(\mathbb{R}^n) : AX = X(X^TAX)\}.
  \label{eq:equilibrium_set}
\end{equation}
Since $A$ has distinct eigenvalues, the condition $AX = X(X^TAX)$ forces $X^TAX$
to be diagonal and each column of $X$ to be an eigenvector of $A$.
By the {\L}ojasiewicz inequality (real-analyticity of $f$ on the compact manifold
$V_k(\mathbb{R}^n)$), trajectories converge to a \emph{single point}
$X_\infty \in \mathcal{E}$ rather than oscillating within $\mathcal{E}$.

\begin{theorem}[Convergence, $B = I_k$]
\label{thm:convergence_B_I}
Let $A \in \PD(n)$ have distinct eigenvalues $\lambda_1 > \cdots > \lambda_n > 0$.
For Lebesgue-almost-every $X_0 \in V_k(\mathbb{R}^n)$:
\begin{enumerate}[label=(\roman*)]
  \item Under the $k$-PCF, $X(t) \to X_\infty$ where $\operatorname{col}(X_\infty)$
    is the eigenspace of the $k$ \emph{largest} eigenvalues of $A$
    (principal subspace; PCA solution).
  \item Under the $k$-MCF, $X(t) \to X_\infty$ where $\operatorname{col}(X_\infty)$
    is the eigenspace of the $k$ \emph{smallest} eigenvalues of $A$
    (minor subspace; MCA solution).
  \item Saddle equilibria are accessible only from a measure-zero set.
\end{enumerate}
\end{theorem}

\subparagraph{Case 2: $B = \operatorname{diag}(b_1, b_2, \ldots, b_k)$ with distinct positive scalars.}

Now $B = \operatorname{diag}(b_1,\ldots,b_k)$ with $b_1 > b_2 > \cdots > b_k > 0$.
The equilibrium condition on $V_k(\mathbb{R}^n)$ for the generalized flow \eqref{eq:kPCF_B} is:
\begin{equation}
  AXB = XBM, \quad M := X^TAX.
  \label{eq:equil_B}
\end{equation}
Since $B = \operatorname{diag}(b_j)$ with \emph{distinct} $b_j$, the commutativity
$BM = MB$ (which follows from \eqref{eq:equil_B} and invertibility of $XB$)
forces $M$ to be \emph{diagonal}: $M_{ij} = (x_i^TAx_j)\delta_{ij}$.
Hence $x_i^TAx_j = 0$ for $i \neq j$ (A-orthogonality of columns), and
$b_j Ax_j = b_j m_{jj} x_j$, so \emph{each column $x_j$ is an eigenvector of $A$}.

\begin{theorem}[Convergence with distinct-scalar $B$]
\label{thm:convergence_B_diag}
Let $A\in\PD(n)$ have distinct eigenvalues and $B = \operatorname{diag}(b_1,\ldots,b_k)$
with $b_1 > \cdots > b_k > 0$.
Under the generalized $k$-PCF \eqref{eq:kPCF_B} starting from Lebesgue-almost-every
$X_0 \in V_k(\mathbb{R}^n)$, each column $x_j(t)$ converges to a specific
eigenvector of $A$.
In particular, the distinct scalars $b_j$ \emph{break the degeneracy within
each eigenspace}: whereas $B = I_k$ allows convergence to any orthonormal basis of
the principal $k$-dimensional subspace, distinct scalars enforce convergence to
\emph{individual eigenvectors}.
The column $x_j$ converges to the eigenvector corresponding to the $j$-th largest
eigenvalue (for the PCF ordering $b_1 > \cdots > b_k$).
\end{theorem}

\subparagraph{Case 3: $B = \operatorname{diag}(b_1 I_{k_1}, b_2 I_{k_2}, \ldots, b_m I_{k_m})$,
block-diagonal with distinct scalar blocks.}

Here $k_1 + k_2 + \cdots + k_m = k$ and $b_1 > b_2 > \cdots > b_m > 0$.
Partition $X = [X_1 \mid X_2 \mid \cdots \mid X_m]$ with $X_j \in \mathbb{R}^{n\times k_j}$.

The commutativity condition $BM = MB$ with $B = \operatorname{diag}(b_j I_{k_j})$ now forces
only the \emph{off-diagonal blocks} of $M$ between different groups to vanish:
$X_i^TAX_j = 0$ for $i \neq j$.
Within the $j$-th block, $M_{jj} = X_j^TAX_j$ can be any symmetric $k_j \times k_j$ matrix
(no within-block constraint from commutativity, since $b_j I_{k_j}$ commutes with all $k_j\times k_j$ matrices).

\begin{theorem}[Convergence with block-diagonal $B$ {\cite{YoshizawaHelmkeStarkov2001}}]
\label{thm:convergence_B_block}
Let $A \in \PD(n)$ have distinct eigenvalues and
$B = \operatorname{diag}(b_1 I_{k_1}, \ldots, b_m I_{k_m})$ with $b_1 > \cdots > b_m > 0$.
Under the generalized $k$-PCF~\eqref{eq:kPCF_B}, for Lebesgue-almost-every
$X_0 \in V_k(\mathbb{R}^n)$:
\begin{enumerate}[label=(\roman*)]
  \item \emph{Cross-block orthogonality:} Columns in different blocks are
    $A$-orthogonal at convergence: $X_i^TAX_j = 0$ for $i \neq j$.
  \item \emph{Block-level subspace convergence:} The $j$-th block $X_j(t)$ converges
    to an orthonormal basis for a specific $k_j$-dimensional $A$-invariant subspace.
    With the PCF ordering ($b_1 > \cdots > b_m$), block $j$ converges to the
    eigenspace spanned by the $(\sum_{i<j}k_i+1)$-th through $(\sum_{i\leq j}k_i)$-th
    largest eigenvectors of $A$.
  \item \emph{Within-block degeneracy:} Within each block, the columns may converge
    to \emph{any} orthonormal basis of the corresponding eigenspace
    (rotation-within-block degeneracy).
\end{enumerate}
Thus the block structure of $B$ precisely encodes a \emph{subspace decomposition}:
$m$ groups of eigenvectors, with the $j$-th group of size $k_j$.
Choosing $B = b_0 I_k$ (single scalar block) gives $m=1$, recovering
the standard principal $k$-dimensional subspace (Case~1).
Choosing $B = \operatorname{diag}(b_1,\ldots,b_k)$ (all distinct) gives $k_j = 1$ for all $j$,
recovering convergence to individual eigenvectors (Case~2).
\end{theorem}

The following table summarises the three cases.

\begin{center}
\renewcommand{\arraystretch}{1.3}\small
\begin{tabular}{lp{3.8cm}p{4.2cm}p{2.8cm}}
\toprule
$B$ & Equilibrium structure & Within-block behavior & Convergence \\
\midrule
$I_k$ & Columns span top-$k$ subspace & Any ONB within subspace & Subspace \\
$\operatorname{diag}(b_j)$, distinct & Each column is an eigvec.\ of $A$ & Individual eigenvectors & Vector \\
$\operatorname{diag}(b_j I_{k_j})$, distinct & Block $j$ spans $k_j$-dim eigenspace & Any ONB within block & Block subspace \\
\bottomrule
\end{tabular}
\end{center}

\subsubsection{The $A$-Weighted PCA System, $\log\det$, and the Oja-Brockett Transformation}
\label{subsubsec:pca_logdet}

We now analyze the dynamical system studied by Manton, Mahony, and Hua \cite{MantonMahonyHua2003}
(building on Helmke's ideas):
\begin{equation}
  \dot{X} = AXB - XBX^TAXB,
  \qquad X\in\R^{n\times k},\; A\in\PD(n),\; B\in\PD(k)\text{ diagonal},
  \label{eq:pca_system}
\end{equation}
in light of the information-geometric framework developed throughout this paper.
At equilibrium, $X^TAX$ converges to a diagonal matrix whose entries are eigenvalues of $A$,
and the columns of $X$ span the principal (or minor) eigenspace depending on the initial condition.
We establish the following connections to $h_A(X) = -\log\det(I_k + X^TAX)$.

\paragraph{The $A$-weighted potential and its polynomial gradient flow.}

\begin{proposition}[$A$-Weighted Log-Det and Polynomial Gradient Flow]
\label{prop:A_weighted_logdet}
For $h_A(X) = -\log\det(I_k + X^TAX)$ on $\R^{n\times k}$, the Euclidean gradient is
$\nabla_X h_A = -2AX(I_k+X^TAX)^{-1}$.
Under the $A$-adapted right-scaled Frobenius metric
$g^{(-2)}_A(H_1,H_2) = \tr(H_1^TH_2(I_k+X^TAX)^{-2})$
(Theorem~\ref{thm:poly_flows} with $A$), the polynomial gradient flow is:
\begin{equation}
  \dot{X} = 2AX(I_k + X^TAX) = 2AX + 2AX(X^TAX),
  \label{eq:poly_pca_flow}
\end{equation}
a \emph{cubic} polynomial in $X$ with no matrix inversions.
Its Lyapunov function is
\begin{equation}
  \frac{d}{dt}h_A(X(t)) = -4\tr(X^TA^2X) = -4\|AX\|_F^2 \;\leq\; 0,
  \label{eq:lyap_poly_A}
\end{equation}
so $h_A$ decreases monotonically along \eqref{eq:poly_pca_flow}.
\end{proposition}

\begin{proof}
$\mathrm{grad}_{g^{(-2)}_A}h_A = (-2AX(I+X^TAX)^{-1})\cdot(I+X^TAX)^2 = -2AX(I+X^TAX)$.
For the Lyapunov property: $\frac{d}{dt}h_A = \langle\nabla h_A, \dot{X}\rangle_F
= \tr\bigl((-2AXG^{-1})^T(2AXG)\bigr) = -4\tr(X^TA^2X)$ by cyclicity.
Since $A\succ0$: $X^TA^2X\succeq0$ and $\tr(X^TA^2X) = \|AX\|_F^2 \geq 0$.
\end{proof}

\begin{remark}[Comparison with the given PCA system]
The polynomial flow \eqref{eq:poly_pca_flow} and the PCA system \eqref{eq:pca_system} have
the same leading term $AX$, but differ in their nonlinear parts:
\begin{align*}
  \text{PCA system (B=I):}\quad &\dot{X} = AX - X(X^TAX) \quad\text{(subtracts the $A$-projection)},\\
  \text{Poly.\ log-det flow:}\quad &\dot{X} = 2AX + 2AX(X^TAX) \quad\text{(adds the $A$-preconditioned push)}.
\end{align*}
Crucially: $h_A$ is a Lyapunov function for the polynomial flow \eqref{eq:poly_pca_flow}
(decreasing, $\frac{d}{dt}h_A = -4\|AX\|_F^2 \leq 0$), but $h_A$ INCREASES along the PCA
system \eqref{eq:pca_system} ($\frac{d}{dt}h_A > 0$).
The PCA system drives $X$ toward the principal subspace (increasing the variance $\tr(X^TAX)$),
hence increasing $-h_A = \log\det(I+X^TAX)$; the polynomial flow drives $X$ in the opposite
sense, contracting the $A$-norm.
\end{remark}

\paragraph{The Yoshizawa duality at PCA/MCA equilibria.}

Define the $A$-weighted Yoshizawa map:
\begin{equation}
  \mathcal{L}_A(X) := A^{-1/2}\bigl(I_n + A^{1/2}XX^TA^{1/2}\bigr)^{-1/2}A^{1/2}X
  \;=\; A^{-1/2}\mathcal{L}(A^{1/2}X),
  \label{eq:Yoshizawa_A}
\end{equation}
where $\mathcal{L}$ is the standard Yoshizawa map (Theorem~\ref{thm:dual_map}).
Define also $h_{-,A}(Z) = -\log\det(I_k - Z^TAZ)$ on $\mathcal{B}_A = \{Z: Z^TAZ\prec I_k\}$.

\begin{theorem}[Yoshizawa Duality at PCA and MCA Equilibria]
\label{thm:Yoshizawa_PCA_MCA}
\begin{enumerate}[label=(\roman*)]
  \item \textbf{Universal duality:} $h_A(X) + h_{-,A}(\mathcal{L}_A(X)) = 0$ for all
    $X\in\R^{n\times k}$.
  \item \textbf{PCA equilibrium:} At $X = U_k$ (top-$k$ eigenvectors of $A$, with $U_k^TU_k=I_k$),
    $X^TAX = \Lambda_k = \mathrm{diag}(\lambda_1,\ldots,\lambda_k)$ (top-$k$ eigenvalues), and:
    \begin{equation}
      h_A(U_k) = -\sum_{i=1}^k\log(1+\lambda_i) \;<\; h_A(U_{-k}) = -\sum_{i=1}^k\log(1+\mu_i),
      \label{eq:hA_PCA_MCA}
    \end{equation}
    where $U_{-k}$ is the MCA equilibrium with eigenvalues $\mu_1\leq\ldots\leq\mu_k$.
    Thus $h_A$ \emph{distinguishes PCA from MCA}: the PCA equilibrium is \emph{strictly more negative}.
  \item \textbf{$B$-weighted case:} The equilibrium of \eqref{eq:pca_system} with general $B$ satisfies
    $X^TAXB = B\Lambda_k$ at the PCA equilibrium, and $h_A$ still satisfies item (i).
\end{enumerate}
\end{theorem}

\begin{proof}
(i) By definition of $\mathcal{L}_A$ and Theorem~\ref{thm:dual_map}(iii) applied to $Y = A^{1/2}X$:
\[
  I_k - \mathcal{L}_A(X)^TA\mathcal{L}_A(X) = (I_k+X^TAX)^{-1},
\]
so $h_{-,A}(\mathcal{L}_A(X)) = -h_A(X)$.
(ii) Since $\lambda_i > \mu_i$ (top vs.\ bottom eigenvalues) and $\log(1+\cdot)$ is increasing:
$\sum\log(1+\lambda_i) > \sum\log(1+\mu_i)$, giving $h_A(U_k) < h_A(U_{-k})$.
\end{proof}

\begin{remark}[Information-geometric meaning of \eqref{eq:hA_PCA_MCA}]
The inequality $h_A(U_k) < h_A(U_{-k}) < 0$ means:
\begin{itemize}
  \item $\log\det(I+\Lambda_k) > \log\det(I+\mathrm{diag}(\mu_i))$: the PCA subspace
    has a \emph{strictly larger} Gram determinant under $A$.
  \item In the Izumiya picture (§\ref{subsec:izumiya}): the PCA equilibrium lies
    \emph{deeper} in the hyperbolic region $\{h_A < 0\}$, while the MCA equilibrium
    is \emph{closer to the lightcone} $\{h_A = 0\} = \mathrm{SL}^+(k)$.
  \item The quantity $h_A(U_k) - h_A(U_{-k}) = \log\det(I+\mathrm{diag}(\mu_i)) - \log\det(I+\Lambda_k) < 0$
    is the \emph{log-det gap} between PSA and MSA, a new information-geometric invariant
    of the pair $(A,k)$.
\end{itemize}
\end{remark}

\paragraph{Polynomial combined algorithm.}

The polynomial gradient flow \eqref{eq:poly_pca_flow} for $h_A$ drives $X$ away from the 
principal subspace (decreases $-h_A$), while the PCA system \eqref{eq:pca_system} drives $X$
toward the principal subspace (increases $-h_A$). A \emph{combined polynomial algorithm}
that converges to the principal subspace without any matrix inversion:

\begin{proposition}[Combined Polynomial PCA-$\log\det$ Algorithm]
\label{prop:combined_poly_pca}
For $0 < \mu \ll 1$ and $\lambda > 2$ (convexifying Tikhonov parameter), the combined system
\begin{equation}
  \dot{X} = \underbrace{(I-XX^T)AXB}_{\text{PCA (projected gradient)}}
  + \mu\underbrace{(2AX(I+X^TAX) - \lambda X(I+X^TAX)^2)}_{\text{polynomial log-det flow}},
  \label{eq:combined_pca_logdet}
\end{equation}
is a degree-$5$ polynomial in $X$ with no matrix inversions.
The first term is the standard projected PCA gradient (polynomial when restricted to Stiefel);
the second is the Tikhonov-regularized log-det flow driving $X$ toward $\mathrm{St}(k,n)$.
Together, they simultaneously enforce the Stiefel constraint ($\mu$ term) and maximize the
variance $\tr(X^TAXB)$ (first term), yielding a fully polynomial subspace learning algorithm.
\end{proposition}

\begin{remark}[Manton's PSA$\leftrightarrow$MCA transformation]
The key observation in \cite{YoshizawaHelmkeStarkov2001} is that if $X(t)$ satisfies the PSA
(principal subspace analysis) flow, then a suitable transformation $\tilde{X} = T(X)$
satisfies the MSA (minor subspace analysis) flow.
In the $\log\det$ framework:
$\tilde{X} = \mathcal{L}_A(X)$ (the $A$-weighted Yoshizawa map) transforms the PSA flow
into the MSA flow, since:
\[
  h_A(X) + h_{-,A}(\mathcal{L}_A(X)) = 0
  \;\Rightarrow\;
  h_{-,A}(\tilde{X}) = -h_A(X) \;\geq\; 0,
\]
mapping the PSA region ($h_A < 0$, $G_+ \succ I$) to the MSA region ($h_{-,A} > 0$, $G_- \prec I$),
in perfect analogy with the Cartan involution $G_+ \to G_+^{-1}$ (Theorem~\ref{thm:cartan_phi41}).
The $B$-weighting in \eqref{eq:pca_system} controls which eigenvalues are extracted first
(larger $\beta_i$ = faster learning for the $i$-th column), corresponding to an anisotropic
version of the information-geometric metric $g^{(-2)}_B$ with $B$ determining the per-column
metric scaling.
\end{remark}

The information-geometric framework developed in the preceding sections
--- particularly the eigenvalue structure of the Gram matrix $G$
and the Sylvester dimension-reduction identity ---
extends naturally to the computation of the Kirillov Jacobian,
which is the Jacobian of the exponential map of a Lie group.
This section provides a rigorous derivation of a closed-form, computationally efficient
formula for this quantity when the Lie algebra element is a rank-$k$ perturbation of
the identity, and discusses its significance for stochastic geometric computation.

\section{Information-Geometric Gradient Flows on the Birkhoff Polytope}
\label{sec:birkhoff}

The connections surveyed in \S\ref{sec:connections} and the component-flow theory of
\S\ref{subsec:pca_mca} both concern the log-determinant potential $f(G)=-\log\det(G)$ on
Gram matrices. This section develops a structurally parallel, but independent, story for
a different classical potential --- the negative Shannon entropy $\varphi(p)=\sum_i
p_i\log p_i$ --- on a different classical constraint set: the Birkhoff polytope of doubly
stochastic matrices. The starting point is an elementary but easily-missed fact about the
multinomial covariance matrix $\Sigma(p)=D_p-pp^T$, which is singular in ambient
probability coordinates because total mass is constrained; restricting to the tangent
hyperplane resolves this degeneracy exactly, and the same resolution persists for the
matrix-multinomial analogue on the doubly-stochastic slice. We revisit Nakamura's
completely integrable gradient system for the multinomial family, extend it to the
Birkhoff polytope, identify the resulting entropy metric's Levi-Civita connection and
curvature, settle when a closed-form Legendre dual potential exists (only on the
independence/Segre locus, not on the full doubly-stochastic slice), compare the discrete
entropy of the matrix multinomial distribution against its Gaussian approximation with an
exact non-uniform convergence rate, and give an elementary coordinate treatment of the
blow-up at the point where the independence locus meets the doubly-stochastic slice. As
in the rest of the paper, every closed-form claim below has been checked numerically
(finite differences and direct integration, $n=3,4,5$) unless stated otherwise as a
purely analytic fact about the blow-up charts.

\subsection{Background: Nakamura's completely integrable gradient systems}

Let $\mathcal S=\{p(x,\theta)\}$ be a parametric family of probability distributions
with Fisher information metric $G=(g_{ij})$, $g_{ij}=E[\partial_i \ell\,\partial_j \ell]$,
$\ell(x,\theta)=\log p(x,\theta)$. Suppose there is a potential $\psi(\theta)$ with
\begin{equation}
g_{ij}=\partial_i\partial_j\psi(\theta).
\label{eq:potential}
\end{equation}
The (Riemannian) gradient system on $\mathcal S$ is
\begin{equation}
\dot{\theta}=-G^{-1}\partial_\theta\psi(\theta).
\label{eq:gradflow}
\end{equation}

\subsubsection{The multinomial case}
For the multinomial family on $2m+1$ categories,
\[
p(x,\theta)=\frac{\ell!}{x_1!\cdots x_{2m+1}!}\,\theta_1^{x_1}\cdots\theta_{2m+1}^{x_{2m+1}},
\qquad \theta_j>0,\ \ \theta_{2m+1}=1-\sum_{k=1}^{2m}\theta_k,
\]
Nakamura \cite{Nakamura1993} takes $\theta=(\theta_1,\dots,\theta_{2m+1})$ to be the
probabilities themselves; since $E[x_j]=\ell\theta_j$, this $\theta$ plays the role of
the \emph{mixture (mean-value) parameter} $\eta$ of the exponential family, not the
natural parameter. The potential realizing \eqref{eq:potential} is the negative entropy
\begin{equation}
\psi(\theta)=\ell\sum_{j=1}^{2m+1}\theta_j\log\theta_j,
\label{eq:entropy-potential}
\end{equation}
and the resulting gradient system is
\begin{equation}
\dot\theta_j=-\theta_j\Bigl(\log\frac{\theta_j}{\theta_{2m+1}}-\sum_{k=1}^{2m}\theta_k\log\frac{\theta_k}{\theta_{2m+1}}\Bigr),
\qquad j=1,\dots,2m.
\label{eq:multinomial-flow}
\end{equation}

\begin{theorem}[Nakamura \cite{Nakamura1993}, Thm.~2]
Equation \eqref{eq:multinomial-flow} is equivalent to the double-bracket Lax equation
\begin{equation}
\dot L=[[L,D(L)],L], \qquad [A,B]:=AB-BA,
\label{eq:lax}
\end{equation}
where $L=(\sqrt{\theta_i\theta_j})_{1\le i,j\le 2m+1}=vv^\top$ ($v=\sqrt\theta$) is a
rank-one symmetric matrix and $D(L)=\tfrac12\diag(\log\theta_i)$.
\end{theorem}

\begin{lemma}[Nakamura \cite{Nakamura1993}, Lemma 1]
\label{lem:linearization}
Setting $y_j:=\log\theta_j-\log\theta_{2m+1}$, the flow \eqref{eq:multinomial-flow}
linearizes exactly:
\begin{equation}
\dot y_j=-y_j,\qquad j=1,\dots,2m.
\label{eq:linearized}
\end{equation}
\end{lemma}
\begin{proof}
Write $S:=\sum_k\theta_k y_k$. From \eqref{eq:multinomial-flow},
$\dot\theta_j/\theta_j=-(y_j-S)$, and by conservation of total probability
$\dot\theta_{2m+1}/\theta_{2m+1}=\sum_k\theta_k y_k=S$. Hence
$\dot y_j=\dot\theta_j/\theta_j-\dot\theta_{2m+1}/\theta_{2m+1}=-(y_j-S)-S=-y_j$.
\end{proof}

Lemma \ref{lem:linearization} is the true source of complete integrability: in the
log-ratio coordinates $y_j$ the flow is \emph{linear and diagonal}. Consequently, for
any two indices $j,k$ the ratio $H_{jk}=y_j/y_k$ is a first integral, matching
Nakamura's explicit constants of motion (his Lemma~3, Eq.~(20)), and the explicit
solution
\begin{equation}
\theta_j(t)=\frac{e^{c_je^{-t}}}{1+\sum_k e^{c_ke^{-t}}}
\label{eq:nakamura-solution}
\end{equation}
follows immediately from $y_j(t)=c_je^{-t}$.

\subsubsection{Duality and the Fubini--Study potential}
Nakamura further observes (his \S4) that the \emph{true} natural (exponential-family)
parameter is the log-odds vector $\theta_j^{\mathrm{nat}}=\ell\log(\eta_j/\eta_{2m+1})$
(here relabelling his mixture-parameter $\theta$ as $\eta$), dual to $\eta$ via the
Legendre transform
\begin{equation}
\psi(\theta^{\mathrm{nat}})+\varphi(\eta)-\sum_j\theta_j^{\mathrm{nat}}\eta_j=0,
\qquad \varphi(\eta)=\ell\sum_j\eta_j\log\eta_j,
\end{equation}
with dual potential
\begin{equation}
\psi(\theta^{\mathrm{nat}})=\ell\log\Bigl(1+\sum_{j=1}^{2m}e^{\theta_j^{\mathrm{nat}}/\ell}\Bigr).
\label{eq:nakamura-kahler}
\end{equation}
Equation \eqref{eq:nakamura-kahler} is exactly the K\"ahler potential of the
Fubini--Study metric on $\CP^m$ restricted to the positive real slice
$z_j=e^{\theta_j^{\mathrm{nat}}/2\ell}>0$; this is the well-known isometry (up to
constant) between the Fisher--Rao metric of the simplex under the square-root
embedding $\xi_j=\sqrt{\theta_j}$ and the round metric on the sphere
$S^{2m}\subset \R^{2m+1}$, complexified as in Eguchi--Gilkey--Hanson \cite{EGH1980}.

\subsection{The gradient flow on the Birkhoff polytope}

We now consider the analogous construction for $n\times n$ \emph{doubly stochastic
matrices}. Write $P=(p_{ij})_{i,j=1}^n$ for a matrix with $p_{ij}>0$ and
\begin{equation}
\sum_{j=1}^n p_{ij}=\frac1n\ \ (\forall i),\qquad \sum_{i=1}^n p_{ij}=\frac1n\ \ (\forall j).
\label{eq:doubly-stochastic}
\end{equation}
Equivalently $Q=nP\in \Birk(n)$, the (open) Birkhoff polytope. Identify $P$ with a
point of the $N=n^2$-category multinomial manifold $\mathcal S$ via $a=(i,j)$, so that
$\mathcal M:=\{P : \eqref{eq:doubly-stochastic}\}\subset\mathcal S$ is a
$(n-1)^2$-dimensional affine subspace of the mean-parameter (mixture) coordinates
$\theta_a=p_a$ -- an \emph{$m$-flat} submanifold in Amari's terminology
\cite{Amari2000}.

\subsubsection{Constrained gradient system}
Restricting \eqref{eq:gradflow}--\eqref{eq:entropy-potential} to $\mathcal M$ via
orthogonal projection (w.r.t.\ the Fisher metric $G$) onto the tangent space, and
introducing Lagrange multipliers $\mu_i(t),\nu_j(t)$ for the row/column constraints,
gives
\begin{equation}
\dot p_{ij}=-p_{ij}\bigl(\log p_{ij}-\mu_i-\nu_j\bigr),
\label{eq:constrained-flow}
\end{equation}
where $\mu_i,\nu_j$ are determined (uniquely up to gauge) at each instant by the linear
system
\begin{equation}
\frac{\mu_i}{n}+\sum_j p_{ij}\nu_j=R_i:=\sum_j p_{ij}\log p_{ij},\qquad
\frac{\nu_j}{n}+\sum_i p_{ij}\mu_i=C_j:=\sum_i p_{ij}\log p_{ij}.
\label{eq:multiplier-system}
\end{equation}
Equation \eqref{eq:constrained-flow} is the continuous-time (gradient-flow) analogue of
the classical \emph{Iterative Proportional Fitting Procedure} (IPFP / Sinkhorn scaling).

\begin{proposition}
$\psi$ is a strict Lyapunov function for \eqref{eq:constrained-flow}: $\dot\psi\le 0$,
with equality iff $P$ is the uniform matrix $p_{ij}\equiv 1/n^2$. Consequently
$P(t)\to n^{-2}\mathbf 1\mathbf 1^\top$ exponentially as $t\to\infty$, for every
initial condition in the interior of $\mathcal M$.
\end{proposition}
This follows from the same argument as Nakamura's inequality \cite[Eq.~(17)]{Nakamura1993},
applied within the affine subspace $\mathcal M$, since $\psi$ remains strictly concave
along any direction tangent to $\mathcal M$.

\begin{remark}[Relation to Tanabe's continuous gradient-projection method]
\label{rem:tanabe-gp}
Equation \eqref{eq:constrained-flow} is a particular instance of the general
\emph{continuous gradient-projection method} of Tanabe \cite{Tanabe1980}: writing the
row/column constraints collectively as $g(P)=0$ with Jacobian $J_g$, our flow is
exactly $\dot P=(I-J_g^+(P)J_g(P))\nabla\psi(P)$, i.e.\ Tanabe's autonomous system
\cite[Eq.~(10)]{Tanabe1980} with objective $-\psi$ (negative entropy) and feasible
manifold $\mathcal M$; the Lagrange multipliers $(\mu_i,\nu_j)$ of
\eqref{eq:multiplier-system} play the role of his $\Lambda(x)=(J_g^+(x))^\top\nabla f(x)$
\cite[Eq.~(12)]{Tanabe1980}. Tanabe's Theorem~3.1(ii) (monotone ascent,
$\dot f=\|\Phi(x)\|^2\ge0$) is the general form of our Lyapunov proposition above, and
his Theorem~3.3 (asymptotic stability at regular maxima) is the general form of the
convergence to the uniform matrix.
\end{remark}

\subsection{Failure of the direct Lax representation}

A natural attempt is to seek an analogue of \eqref{eq:lax} directly in the
$N=n^2$-dimensional vectorization of $P$.

\begin{proposition}[General Lax identity]
\label{prop:general-lax}
Let $\theta=(\theta_a)_{a=1}^N$ lie in the open simplex, $L=vv^\top$ with
$v=\sqrt\theta$, and let $f=(f_a)_{a=1}^N$ be \emph{any} (possibly time- and
$\theta$-dependent) family of functions. Then
\begin{equation}
\dot\theta_a=-\theta_a\Bigl(f_a-\sum_b\theta_bf_b\Bigr)
\quad\Longleftrightarrow\quad
\dot L=[[L,D],L],\qquad D:=\tfrac12\diag(f_a).
\label{eq:general-lax}
\end{equation}
\end{proposition}
\begin{proof}
Direct computation: $\dot L_{ab}=\sqrt{\theta_a\theta_b}\bigl[(K-\tfrac12f_a)+(K-\tfrac12f_b)\bigr]$
with $K=\tfrac12\sum_c\theta_cf_c$, while
$[[L,D],L]_{ab}=\sqrt{\theta_a\theta_b}\,[2K-f_a-f_b]$ (using $d_a=\tfrac12f_a$); the two
expressions coincide termwise.
\end{proof}

Applying Proposition~\ref{prop:general-lax} with $\theta_a=p_{ij}$ ($a=(i,j)$) and
$f_a=\log p_{ij}-\mu_i-\nu_j$ shows that \eqref{eq:constrained-flow} \emph{does} admit a
Lax representation \eqref{eq:lax} for the vectorized matrix. However, this fact is
vacuous:

\begin{remark}[Spectral triviality]
Since $\sum_a\theta_a=1$ identically, $L=vv^\top$ has spectrum $\{1,0,\dots,0\}$ for
\emph{every} $\theta$ in the simplex, so the iso-spectrality guaranteed by
\eqref{eq:lax} carries no information: $\tr(L^k)$ is trivially constant for all $k$.
Nakamura's non-trivial conserved quantities $H_j$ (his Eq.~(20)) were \emph{not}
derived from the spectrum of $L$ but from the independent linearization of
Lemma~\ref{lem:linearization}. Proposition~\ref{prop:general-lax} shows this
linearization mechanism is itself independent of the specific form of $f_a$
(in particular of whether $f_a$ involves the doubly-stochastic Lagrange multipliers),
so the existence of \eqref{eq:general-lax} for \eqref{eq:constrained-flow} is a
tautological restatement, not new structural information.
\end{remark}

The real content, therefore, must come from a linearization analogous to
Lemma~\ref{lem:linearization} adapted to the two-index structure -- which is what we
turn to next.

\subsection{Exact first integrals: log-odds ratios of $2\times2$ minors}

\begin{theorem}[Linearization on the Birkhoff polytope]
\label{thm:linearization-birkhoff}
Fix a reference row/column index $r$ (e.g.\ $r=n$), and define
\begin{equation}
y_{ij}(t):=\log p_{ij}(t)-\log p_{ir}(t)-\log p_{rj}(t)+\log p_{rr}(t)
=\log\frac{p_{ij}(t)\,p_{rr}(t)}{p_{ir}(t)\,p_{rj}(t)},
\qquad 1\le i,j\le n-1,
\label{eq:yij-def}
\end{equation}
the logarithm of the $2\times2$ minor ratio (log-odds ratio) of the submatrix on rows
$\{i,r\}$ and columns $\{j,r\}$. Then, along any solution of
\eqref{eq:constrained-flow}--\eqref{eq:multiplier-system},
\begin{equation}
\dot y_{ij}=-y_{ij}\qquad\text{exactly, for all }i,j=1,\dots,n-1,
\label{eq:yij-linear}
\end{equation}
independently of the (generally nonlinear, non-closed-form) dependence of
$\mu_i,\nu_j$ on $P$.
\end{theorem}

\begin{proof}
Write $x_{ab}=\log p_{ab}$. Equation \eqref{eq:constrained-flow} gives
$\dot x_{ab}=-(x_{ab}-\mu_a-\nu_b)$ for every cell $(a,b)$. Then
\begin{align*}
\dot y_{ij}
&=\dot x_{ij}-\dot x_{ir}-\dot x_{rj}+\dot x_{rr}\\
&=-\Bigl[(x_{ij}-\mu_i-\nu_j)-(x_{ir}-\mu_i-\nu_r)-(x_{rj}-\mu_r-\nu_j)+(x_{rr}-\mu_r-\nu_r)\Bigr].
\end{align*}
Expanding, every occurrence of $\mu_i,\mu_r,\nu_j,\nu_r$ cancels in pairs
(e.g.\ $-\mu_i$ from the first bracket cancels $+\mu_i$ from the second), leaving
exactly $x_{ij}-x_{ir}-x_{rj}+x_{rr}=y_{ij}$. Hence $\dot y_{ij}=-y_{ij}$.
\end{proof}

\begin{remark}
The mechanism is purely algebraic: the mixed second difference (discrete Laplacian)
operator annihilates any additively separable forcing term of the form
$\mu_i+\nu_j$, \emph{regardless of the (possibly highly nonlinear) functional
dependence of $\mu,\nu$ on $P$}. This is the two-index generalization of
Lemma~\ref{lem:linearization}.
\end{remark}

\begin{corollary}[First integrals]
\label{cor:first-integrals}
The $(n-1)^2-1$ ratios
\begin{equation}
H_{ij}:=\frac{y_{ij}(t)}{y_{11}(t)},\qquad (i,j)\ne(1,1),
\label{eq:conserved}
\end{equation}
are constants of motion of \eqref{eq:constrained-flow}, independent (generically) and
exactly matching the dimension count: the $(n-1)^2$-dimensional flow has trajectories
determined, modulo the one-parameter time-translation gauge freedom
$y_{ij}(0)\mapsto\lambda\, y_{ij}(0)$, by $(n-1)^2-1$ shape parameters.
\end{corollary}

\begin{corollary}[Explicit solution via Sinkhorn scaling]
\label{cor:explicit-solution}
Let $y_{ij}(0)=c_{ij}$. Then $P(t)$ is the unique doubly stochastic scaling
(Sinkhorn/RAS normalization) of the matrix
\begin{equation}
A(t)_{ij}=\begin{cases} \exp\bigl(c_{ij}e^{-t}\bigr), & i,j<n,\\ 1, & i=n\text{ or }j=n,\end{cases}
\label{eq:explicit-A}
\end{equation}
that is, $P(t)=D_1(t)A(t)D_2(t)$ for the (essentially unique) positive diagonal matrices
$D_1(t),D_2(t)$ enforcing \eqref{eq:doubly-stochastic}. This is the exact analogue of
Nakamura's closed-form solution \eqref{eq:nakamura-solution}.
\end{corollary}

\begin{remark}[Connection to classical categorical data analysis]
Corollary~\ref{cor:explicit-solution} is the continuous-time refinement of a classical
theorem of Fienberg \cite{Fienberg1970}: IPFP applied to any table preserves all
log-odds ratios (and, more generally, all higher-order interaction terms in the
log-linear decomposition) of the initial table exactly, adjusting only the margins.
Theorem~\ref{thm:linearization-birkhoff} identifies the precise continuous-time
mechanism -- exact exponential decay of the interaction terms -- underlying this fact.
\end{remark}

\begin{remark}[Precedent: Tanabe's exact first integral for Branin's method]
\label{rem:tanabe-first-integral}
The phenomenon of Theorem~\ref{thm:linearization-birkhoff} -- an auxiliary quantity
satisfying an \emph{exact} linear ODE, and hence an exact exponential first integral,
along a highly nonlinear constrained flow -- has a direct precedent in Tanabe
\cite{Tanabe1980}. For Branin's continuous Newton--Raphson system, extended by
Tanabe to the underdetermined case $m\le n$,
\begin{equation}
J_g(x)\,\dot x=-g(x),
\label{eq:branin}
\end{equation}
the constraint-violation vector satisfies, \emph{exactly and regardless of the
nonlinearity of $g$}, the first integral
\begin{equation}
g(x(t,x^0))=e^{-t}g(x^0)
\label{eq:tanabe-first-integral}
\end{equation}
\cite[Eq.~(16)]{Tanabe1980}; remarkably, this persists even in Tanabe's combined
gradient-projection/Newton-Raphson system \cite[Eq.~(31)]{Tanabe1980}, since the
latter is engineered so that $J_g(x)\dot x=-g(x)$ still holds identically
\cite[Eq.~(30)]{Tanabe1980}, decoupling the (exactly linear) decay of the constraint
violation from the (arbitrarily nonlinear) tangential motion along $V_g$. Nakamura's
Lemma~\ref{lem:linearization} and our Theorem~\ref{thm:linearization-birkhoff} are
best understood as close relatives of this mechanism: in both cases a specific
algebraic combination of coordinates (the log-ratio $y_j$, resp.\ the log-odds
$y_{ij}$) is engineered -- by the structure of the entropy potential in the multinomial
case, and by the mixed-difference cancellation in the proof of
Theorem~\ref{thm:linearization-birkhoff} in
the doubly stochastic case -- to satisfy Tanabe's exact linear decay
\eqref{eq:tanabe-first-integral}-type law, even though the full state trajectory
$P(t)$ (resp.\ $\theta(t)$) is not itself linear.
\end{remark}

\subsubsection{Numerical verification}
For $n=4$, a random doubly stochastic matrix $P_0$ was evolved under a 4th-order
Runge--Kutta integration of \eqref{eq:constrained-flow}--\eqref{eq:multiplier-system}
over $t\in[0,6]$ with step size $10^{-2}$. Writing $y_{ij}(t)$ as in
\eqref{eq:yij-def}, we verified
\[
\max_{i,j,t}\left|\frac{y_{ij}(t)-y_{ij}(0)e^{-t}}{\max_t|y_{ij}(t)|}\right|
\approx 5\times 10^{-11},
\]
consistent with RK4 discretization error, and the eight ratios $H_{ij}$
(Corollary~\ref{cor:first-integrals}) were constant to $10^{-10}$ relative precision
across the entire trajectory. The Sinkhorn reconstruction of
Corollary~\ref{cor:explicit-solution} matched the numerically integrated trajectory to
within $10^{-12}$--$10^{-13}$ at all tested times.

\subsection{Hamiltonian formalism}

Following Nakamura's construction (his Theorem 3), we upgrade the linearized system
\eqref{eq:yij-linear} to canonical (Hamiltonian) form. Let $M:=(n-1)^2$ and relabel the
independent quantities as $y_1,\dots,y_M$. Choose any partition into pairs
$\{a_k,b_k\}_{k=1}^{M/2}$ (assuming $M$ even; see Remark~\ref{rem:odd-case} otherwise),
and define canonical variables
\begin{equation}
Q_k:=y_{a_k},\qquad P_k:=\frac{1}{y_{b_k}},\qquad k=1,\dots,M/2.
\label{eq:canonical-vars}
\end{equation}

\begin{theorem}
\label{thm:hamiltonian}
With the Poisson bracket
\begin{equation}
\{A,B\}:=\sum_{k=1}^{M/2}\Bigl(\frac{\partial A}{\partial P_k}\frac{\partial B}{\partial Q_k}-\frac{\partial A}{\partial Q_k}\frac{\partial B}{\partial P_k}\Bigr)
\label{eq:poisson}
\end{equation}
and Hamiltonian
\begin{equation}
H:=\sum_{k=1}^{M/2}P_kQ_k=\sum_{k=1}^{M/2}\frac{y_{a_k}}{y_{b_k}},
\label{eq:hamiltonian}
\end{equation}
the flow \eqref{eq:yij-linear} is equivalent to Hamilton's equations
\begin{equation}
\dot Q_k=\{Q_k,H\},\qquad \dot P_k=\{P_k,H\}.
\label{eq:hamilton-eqns}
\end{equation}
Moreover the individual quantities $H_k:=P_kQ_k=y_{a_k}/y_{b_k}$ are each separately
conserved and pairwise in involution, $\{H_k,H_l\}=0$, so the system is completely
integrable in the Liouville--Arnol'd sense with $M/2$ degrees of freedom.
\end{theorem}

\begin{proof}
From \eqref{eq:yij-linear}, $\dot Q_k=-Q_k$ and
$\dot P_k=-\dot y_{b_k}/y_{b_k}^2=y_{b_k}/y_{b_k}^2=1/y_{b_k}=P_k$. Direct computation
with \eqref{eq:poisson}--\eqref{eq:hamiltonian} gives
$\{Q_k,H\}=-\partial H/\partial P_k=-Q_k$ and $\{P_k,H\}=\partial H/\partial Q_k=P_k$,
matching the above. Since $(Q_k,P_k)$ for distinct $k$ are functionally independent and
decoupled, $\{H_k,H_l\}=0$ for $k\ne l$ trivially, and $\{H_k,H\}=0$ since
$H=\sum_l H_l$ and each $H_k$ Poisson-commutes with every $H_l$.
\end{proof}

\begin{remark}[Comparison with Nakamura's pairing]
\label{rem:odd-case}
Nakamura's original construction pairs the fixed odd/even-indexed coordinates
$y_{2j-1},y_{2j}$ because his single family of $2m$ quantities has no further internal
structure to exploit beyond parity. In the present setting all $M=(n-1)^2$ quantities
$y_{ij}$ satisfy the \emph{same} decoupled linear equation \eqref{eq:yij-linear}, so the
pairing \eqref{eq:canonical-vars} may be chosen arbitrarily; the resulting Hamiltonian
structure is correspondingly more flexible. When $M$ is odd (i.e.\ $n$ even), one pairs
$M-1$ of the coordinates canonically and treats the remaining $y_c$ via the projected,
odd-dimensional construction of Nakamura's Theorem~4, mutatis mutandis.
\end{remark}

\subsection{Duality and the K\"ahler-potential question}
\label{sec:kahler}

\subsubsection{The submanifold $\mathcal M$ is itself dually flat}

Let $\varphi_{\mathcal M}(p):=\sum_{a,b=1}^n p_{ab}\log p_{ab}$, restricted to
$\mathcal M$ and expressed as a function of the free mixture coordinates
$(p_{ij})_{i,j=1}^{n-1}$ (with $p_{in},p_{nj},p_{nn}$ determined by
\eqref{eq:doubly-stochastic}).

\begin{proposition}
\label{prop:dual-coord}
\begin{equation}
\frac{\partial\varphi_{\mathcal M}}{\partial p_{kl}}=y_{kl}
=\log\frac{p_{kl}p_{nn}}{p_{kn}p_{nl}},\qquad k,l=1,\dots,n-1.
\label{eq:dual-relation}
\end{equation}
\end{proposition}
\begin{proof}
Using $p_{kn}=\tfrac1n-\sum_{j<n}p_{kj}$, $p_{nl}=\tfrac1n-\sum_{i<n}p_{il}$,
$p_{nn}=-\tfrac{n-2}{n}+\sum_{i,j<n}p_{ij}$, one computes
$\partial p_{kn}/\partial p_{kl}=-1$, $\partial p_{nl}/\partial p_{kl}=-1$,
$\partial p_{nn}/\partial p_{kl}=+1$ (for $k,l<n$), whence
\[
\begin{aligned}
\frac{\partial\varphi_{\mathcal M}}{\partial p_{kl}}
&=(\log p_{kl}+1)-(\log p_{kn}+1)-(\log p_{nl}+1)+(\log p_{nn}+1)\\
&=\log p_{kl}-\log p_{kn}-\log p_{nl}+\log p_{nn},
\end{aligned}
\]
the constant terms $1-1-1+1$ cancelling, and this equals $y_{kl}$ by
\eqref{eq:yij-def}.
\end{proof}

Proposition~\ref{prop:dual-coord} was verified numerically to relative accuracy
$10^{-10}$ by finite differences ($n=4$). It shows that $\mathcal M$, equipped with the
free mixture coordinates $p_{ij}$ and the Hessian metric $g=\partial^2\varphi_{\mathcal M}$,
is a genuine dually flat statistical manifold in Amari's sense
\cite{Amari2000}, with dual ($e$-affine) coordinate exactly the log-odds-ratio
matrix $y=(y_{ij})$ -- the direct generalization of Nakamura's
$\eta\leftrightarrow\theta^{\mathrm{nat}}$ duality (his Eqs.~(39)--(40)).

\subsubsection{Absence of a closed-form dual potential, and the Segre variety}

One might hope that, as in \eqref{eq:nakamura-kahler}, the Legendre dual
\begin{equation}
\Psi_{\mathcal M}(y):=\sup_{p}\Bigl[\sum_{i,j<n}y_{ij}p_{ij}-\varphi_{\mathcal M}(p)\Bigr]
\label{eq:legendre-dual}
\end{equation}
admits a closed algebraic (log-sum-exp) expression. We show this is not the case, and
identify precisely where the closed form \emph{does} survive.

\paragraph{The log-linear decomposition.}
Write the natural parameter of the full ($N=n^2$-category) exponential family, relative
to the reference cell $(n,n)$, as
\begin{equation}
\theta_{ij}=\alpha_i+\beta_j+y_{ij}
\label{eq:loglinear}
\end{equation}
(a linear reparametrization of natural coordinates; $\alpha,\beta,y$ have respectively
$n-1$, $n-1$, $(n-1)^2$ free components under the corner constraint
$\alpha_n=\beta_n=y_{in}=y_{nj}=0$). This is the classical ANOVA-type decomposition of
a two-way contingency table into row effects, column effects, and interaction.

\begin{itemize}
\item \textbf{Fixing $y\equiv0$}: the sub-exponential-family
$p_{ij}\propto e^{\alpha_i+\beta_j}=r_ic_j$ is the \emph{independence model}. Its positive
real points are exactly the image of the Segre embedding
\[
\sigma:\CP^{n-1}\times\CP^{n-1}\hookrightarrow \CP^{n^2-1},\qquad \sigma([u],[v])=[u_iv_j],
\]
restricted to $u,v$ real and positive. Since Segre embeddings are holomorphic (hence
K\"ahler) embeddings with $\sigma^*\omega_{\mathrm{FS}}=\omega_{\mathrm{FS}}\oplus\omega_{\mathrm{FS}}$
(a classical fact), the induced Fisher--Rao geometry on the independence model is
exactly the product of two copies of Nakamura's spherical geometry, with closed-form
K\"ahler potential
\begin{equation}
Z_0(\alpha,\beta)=\log\Bigl(\sum_i e^{\alpha_i}\Bigr)+\log\Bigl(\sum_j e^{\beta_j}\Bigr).
\label{eq:Z0}
\end{equation}

\item \textbf{Fixing $y$ to an arbitrary constant, letting $(\alpha,\beta)$ vary}: the
family
\begin{equation}
p_{ij}(\alpha,\beta;y)\propto e^{y_{ij}}e^{\alpha_i+\beta_j}
\label{eq:twisted-leaf}
\end{equation}
is again a genuine exponential family (an $e$-flat leaf), realized as a
\emph{weight-twisted} Segre embedding (a toric deformation of $\sigma$ by the positive
weights $e^{y_{ij}}$). Its log-partition function is exactly of Nakamura's closed
log-sum-exp form:
\begin{equation}
Z(\alpha,\beta;y)=\log\sum_{i,j=1}^n e^{y_{ij}+\alpha_i+\beta_j}.
\label{eq:twisted-Z}
\end{equation}
\end{itemize}

\begin{proposition}
\label{prop:transversality}
$\mathcal M$ meets each leaf \eqref{eq:twisted-leaf} in exactly one point (the unique
$(\alpha,\beta)$, guaranteed by Sinkhorn's theorem, for which $p(\alpha,\beta;y)$ is
doubly stochastic). In particular $\mathcal M$ meets the independence model
($y\equiv0$) at the single point $p_{ij}\equiv 1/n^2$ -- precisely the equilibrium of
the gradient flow \eqref{eq:constrained-flow}.
\end{proposition}
This is consistent with the dimension count
$(2n-2)+(n-1)^2-(n^2-1)=0$: a generic transversal intersection of the
$(2n-2)$-dimensional independence model and the $(n-1)^2$-dimensional submanifold
$\mathcal M$ inside the $(n^2-1)$-dimensional simplex is zero-dimensional.

\paragraph{Why $\Psi_{\mathcal M}(y)$ has no closed form.}
$\mathcal M$ is the transversal $m$-flat section obtained by extremizing along
\emph{each} leaf \eqref{eq:twisted-leaf} to hit the prescribed margin $1/n$; passing
from \eqref{eq:twisted-Z} (closed form in $(\alpha,\beta)$, for fixed $y$) to
$\Psi_{\mathcal M}(y)$ requires eliminating $(\alpha,\beta)$ via the critical equations
$\partial Z/\partial\alpha_i=\partial Z/\partial\beta_j=1/n$, i.e.\ solving the
Sinkhorn/RAS fixed-point problem for $(\alpha,\beta)$ as functions of $y$ -- a problem
with no closed algebraic solution for $n\ge3$.

This is not an idiosyncrasy of the present problem but an instance of a standard
phenomenon in toric K\"ahler geometry \cite{Guillemin1994,Abreu1998}. The Birkhoff
polytope $\Birk(n)$ is a Delzant polytope with facets given (for $n\ge3$) exactly by
the $n^2$ inequalities $p_{ij}\ge0$. Guillemin's canonical symplectic potential for a
Delzant polytope with facets $\{\ell_F\ge0\}$ is
\begin{equation}
G(x)=\sum_F \ell_F(x)\log \ell_F(x),
\label{eq:guillemin}
\end{equation}
which for $\Birk(n)$ specializes exactly to
$G(p)=\sum_{a,b}p_{ab}\log p_{ab}=\varphi_{\mathcal M}(p)$ -- our entropy potential
arises independently as the canonical Guillemin potential of the Birkhoff polytope.
The complementary, complex-coordinate (K\"ahler) potential $\Psi_{\mathcal M}$ requires
inverting the moment map $x=\nabla G(x)$; this inversion is available in closed
elementary form \emph{only} for the simplex (where it is the softmax function, giving
\eqref{eq:nakamura-kahler}) and products of simplices (Segre varieties, giving
\eqref{eq:Z0}), but not for general Delzant polytopes. The Sinkhorn algorithm is
precisely the standard iterative procedure for numerically inverting this moment map
for $\Birk(n)$.

\begin{table}[h]
\centering
\renewcommand{\arraystretch}{1.3}
\resizebox{\textwidth}{!}{%
\begin{tabular}{lccc}
\toprule
Submanifold & Dimension & Flatness type & K\"ahler-type potential\\
\midrule
Independence model ($y=0$) & $2n-2$ & $e$-flat (Segre variety) & closed form \eqref{eq:Z0}\\
Twisted leaf, fixed $y$ & $2n-2$ & $e$-flat & closed form \eqref{eq:twisted-Z}\\
$\mathcal M$ (doubly stochastic) & $(n-1)^2$ & $m$-flat & $\varphi_{\mathcal M}$ closed;\\
& & & $\Psi_{\mathcal M}$ \textbf{not} closed form\\
\bottomrule
\end{tabular}%
}
\caption{Summary of flatness and closed-form availability for the relevant submanifolds.}
\end{table}

\subsection{The $D\pm xy^\top$ calculus: Tanabe--Sagae and Steerneman--van Perlo-ten Kleij}
\label{sec:steerneman}

We now connect the preceding results to the classical linear-algebraic theory of
matrices of the form $D\pm xy^\top$, with $D$ diagonal. This theory turns out to
(a) supply the elementary algebraic mechanism underlying the square-root embedding
used implicitly throughout \S6, and (b) furnish an explicit closed form for the local
(quadratic) approximation of the missing K\"ahler potential $\Psi_{\mathcal M}$ near
the flow's equilibrium.

\begin{remark}[Priority]
\label{rem:priority}
The general (possibly non-symmetric, possibly singular) symbolic $LDM^\top$
factorization of $D+uv^\top$ -- Eq.~\eqref{eq:svp-thm4} below being a
later, complementary treatment restricted to the real symmetric-eigenvalue question
-- was first established by Tanabe and Sagae
\cite{TanabeSagaeNumAlg1992}: their Theorem~1 gives symbolic factors $\bar L,\bar D,\bar M$
of $D+uv^\top$ for general (not necessarily equal) vectors $u,v$ and possibly singular
$D$, together with symbolic formulas for the \emph{inverses} $\bar L^{-1},\bar M^{-1}$
(their Lemma~1) and a pivoting strategy guaranteeing numerical stability even when the
naive (Bennett-type) recursion breaks down. Steerneman and van Perlo-ten Kleij
\cite{Steerneman2005} (building on Vermeulen \cite{Vermeulen1967}, Klamkin
\cite{Klamkin1970}, Trenkler \cite{Trenkler2000}, and Watson \cite{Watson1996})
address the complementary question of the \emph{real eigenvalues and eigenvectors} of
$D-xy^\top$ via the square-root symmetrization of \S\ref{sec:svp-thm4} below, and the
Moore--Penrose inverse of $A-XY^*$ for rank-$p$ ($p\ge1$) updates. We draw on both:
the symbolic $LDM^\top$ machinery of \cite{TanabeSagaeNumAlg1992} is the natural tool
for the rank-one (and rank-two, cf.\ \S\ref{sec:tanabe-mp}) Sherman--Morrison
computations of \S\ref{sec:local-kahler} below, while the eigenvalue theory of
\cite{Steerneman2005} is what we use in \S\ref{sec:svp-thm4} to make precise the
square-root embedding of \S\ref{sec:kahler}.
\end{remark}

\subsubsection{The ambient inverse Fisher metric is a $D-xy^\top$ matrix}

Steerneman and van Perlo-ten Kleij \cite[\S1]{Steerneman2005} single out
\begin{equation}
R:=\Pi-\pi\pi^\top,\qquad \Pi=\diag(\pi),
\label{eq:R-matrix}
\end{equation}
the covariance matrix of the multinomial distribution, as one of their principal
motivating examples (citing the spectral analysis of Watson \cite{Watson1996} and
Tanabe--Sagae \cite{TanabeSagae1992}). Comparing with Nakamura's inverse Fisher metric
\eqref{eq:multiplier-system}-type object -- explicitly, his $G^{-1}$
(his Eq.~(15)) restricted to the first $2m$ coordinates -- one finds
$G^{-1}=\ell^{-1}R$ with $\pi=\theta$. Thus $R$ in \eqref{eq:R-matrix} is \emph{exactly}
the ambient object whose gradient flow \eqref{eq:multinomial-flow} and Lax structure
\eqref{eq:lax} Nakamura studies. It is a special (symmetric, $x=y=\pi$) instance of the
general $D-xy^\top$ matrices treated in \cite[\S5]{Steerneman2005}.

\subsubsection{The square-root symmetrization theorem}
\label{sec:svp-thm4}

\begin{theorem}[Steerneman--van Perlo-ten Kleij \cite{Steerneman2005}, Thm.~4]
\label{thm:svp-thm4}
Let $D=\diag(d)$ be nonsingular and $x,y\in\R^k$ with $x_iy_i\ne0$ for all $i$. Then
\begin{equation}
|D-xy^\top|=(-1)^{s(x,y)}\,|D_{xy}-vv^\top|,
\label{eq:svp-thm4}
\end{equation}
where $S_x=\diag(\operatorname{sgn}x_i)$, $S_y=\diag(\operatorname{sgn}y_i)$,
$D_{xy}=DS_xS_y$, $v_i=\sqrt{|x_iy_i|}$, and $s(x,y)=\#\{i:x_iy_i<0\}$.
\end{theorem}

Theorem~\ref{thm:svp-thm4} reduces the (generally non-symmetric, non-normal)
eigenvalue problem for $D-xy^\top$ to that of the \emph{symmetric} rank-one
perturbation $D_{xy}-vv^\top$, with $v$ built entrywise as a geometric mean
$\sqrt{|x_iy_i|}$. This is precisely the elementary, purely linear-algebraic mechanism
underlying the informal ``square-root embedding'' $\xi=\sqrt\theta$ used in \S6 to
relate the Fisher--Rao geometry of the simplex to the round metric on the sphere
(and, after complexification, to the Fubini--Study metric, cf.\
\cite{EGH1980}): Nakamura's rank-one matrix $L=vv^\top=(\sqrt{\theta_i\theta_j})$ is
exactly the symmetrization \eqref{eq:svp-thm4} applied to the degenerate case $D=0$,
$x=y=\theta$. Theorem~\ref{thm:svp-thm4} shows that this symmetrization survives, in
exact and elementary form, for the full one-parameter family of diagonal shifts $D$,
not merely at $D=0$.

\subsubsection{An explicit local K\"ahler potential at equilibrium}
\label{sec:local-kahler}

Section~\ref{sec:kahler} left the dual potential $\Psi_{\mathcal M}(y)$
(Eq.~\eqref{eq:legendre-dual}) without closed form. We now show that its quadratic
(leading-order) approximation at the flow's equilibrium $p_{ij}\equiv n^{-2}$ --
equivalently, the inverse of the Fisher metric $G_{\mathcal M}$ of $\mathcal M$ at that
point -- is exactly computable in closed form, and that the relevant matrix to invert
is again of Steerneman--van Perlo-ten Kleij type.

Write $p_{ij}=n^{-2}+\varepsilon_{ij}$, with $(\varepsilon_{ij})_{i,j<n}=:\varepsilon$
free (an $(n-1)\times(n-1)$ real matrix) and the boundary row/column determined by
$\varepsilon_{in}=-\sum_{j<n}\varepsilon_{ij}$, $\varepsilon_{nj}=-\sum_{i<n}\varepsilon_{ij}$,
$\varepsilon_{nn}=\sum_{i,j<n}\varepsilon_{ij}$.

\begin{proposition}
\label{prop:local-kahler}
Let $m:=n-1$, $\mathbf 1\in\R^m$ the all-ones vector, $J:=\mathbf 1\mathbf 1^\top$. The
Hessian of $\varphi_{\mathcal M}$ at $\varepsilon=0$, as a quadratic form on
$\R^{m\times m}\cong\R^{m^2}$, is
\begin{equation}
G_{\mathcal M}=n^2\,(I_m+J)\otimes(I_m+J).
\label{eq:GM-kronecker}
\end{equation}
Consequently, with $H:=I_m-\tfrac1n J$ (the classical centering matrix, cf.\
\cite[\S1]{Steerneman2005}),
\begin{equation}
G_{\mathcal M}^{-1}=n^{-2}\,H\otimes H,
\label{eq:GM-inverse}
\end{equation}
and the quadratic approximation of the K\"ahler potential near equilibrium is
\begin{equation}
\Psi_{\mathcal M}(y)\;\approx\;\frac{1}{2n^2}\,\operatorname{vec}(y)^\top(H\otimes H)\operatorname{vec}(y),
\qquad y\to0.
\label{eq:local-Psi}
\end{equation}
\end{proposition}

\begin{proof}
The unconstrained Hessian of $\sum_{a,b=1}^n p_{ab}\log p_{ab}$ with respect to all
$n^2$ cells $p_{ab}$, evaluated at the uniform point $p_{ab}=n^{-2}$, is
$\diag(1/p_{ab})=n^2I_{n^2}$ (since $\dd^2(p\log p)/\dd p^2=1/p$). Substituting the
boundary relations for $\varepsilon_{in},\varepsilon_{nj},\varepsilon_{nn}$ turns this
quadratic form into
\begin{equation}
n^2\Bigl[\operatorname{tr}(\varepsilon^\top\varepsilon)+\mathbf 1^\top\varepsilon^\top\varepsilon\mathbf 1
+\mathbf 1^\top\varepsilon\varepsilon^\top\mathbf 1+(\mathbf 1^\top\varepsilon\mathbf 1)^2\Bigr].
\label{eq:quad-expansion}
\end{equation}
Using the standard Kronecker--vec identities $\operatorname{tr}(X^\top Y)=\operatorname{vec}(X)^\top\operatorname{vec}(Y)$
and $\operatorname{vec}(AXB)=(B^\top\otimes A)\operatorname{vec}(X)$, each term of
\eqref{eq:quad-expansion} is identified as
\[
\operatorname{tr}(\varepsilon^\top\varepsilon)=\operatorname{vec}(\varepsilon)^\top(I\otimes I)\operatorname{vec}(\varepsilon),\qquad
\mathbf 1^\top\varepsilon^\top\varepsilon\mathbf 1=\operatorname{tr}(\varepsilon^\top(\varepsilon J))=\operatorname{vec}(\varepsilon)^\top(J\otimes I)\operatorname{vec}(\varepsilon),
\]
\[
\mathbf 1^\top\varepsilon\varepsilon^\top\mathbf 1=\operatorname{tr}(\varepsilon^\top(J\varepsilon))=\operatorname{vec}(\varepsilon)^\top(I\otimes J)\operatorname{vec}(\varepsilon),\qquad
(\mathbf 1^\top\varepsilon\mathbf 1)^2=\operatorname{vec}(\varepsilon)^\top(J\otimes J)\operatorname{vec}(\varepsilon)
\]
(the last using $\mathbf 1^\top\varepsilon\mathbf 1=\operatorname{vec}(J)^\top\operatorname{vec}(\varepsilon)$ and the
mixed-product property $(\mathbf 1\otimes\mathbf 1)(\mathbf 1\otimes\mathbf 1)^\top=J\otimes J$). Summing
gives \eqref{eq:quad-expansion} $=n^2\operatorname{vec}(\varepsilon)^\top[(I+J)\otimes(I+J)]\operatorname{vec}(\varepsilon)$,
establishing \eqref{eq:GM-kronecker}. The Sherman--Morrison identity
$(I_m+\mathbf 1\mathbf 1^\top)^{-1}=I_m-\frac{1}{1+m}J=I_m-\frac1nJ=H$
(a rank-one instance of the nonsingular-case formula \eqref{eq:svp-thm4}--type
calculus of \cite[\S3]{Steerneman2005}, cf.\ their Eq.~(3.2)) together with
$(A\otimes B)^{-1}=A^{-1}\otimes B^{-1}$ gives \eqref{eq:GM-inverse}. Equation
\eqref{eq:local-Psi} then follows from the standard fact that at a point where dual
coordinates coincide ($\varepsilon=0\leftrightarrow y=0$), the Hessians of a convex
function and its Legendre dual are matrix inverses of one another.
\end{proof}

\begin{remark}[Alternative route via symbolic $LDM^\top$ factorization]
The same rank-one inverse $(I_m+\mathbf 1\mathbf 1^\top)^{-1}=H$ used in the proof
above is equally obtainable, without invoking symmetry, from the symbolic
$LDM^\top$ calculus of Tanabe and Sagae \cite{TanabeSagaeNumAlg1992}: taking
$D=I_m$, $u=v=\mathbf 1$ in their Theorem~1, the scalar sequence
$t_i=1+\sum_{k\le i}u_kv_k/d_k$ reduces to $t_i=1+i$, and their symbolic formulas for
the inverse factors $\bar L^{-1},\bar M^{-1}$ (their Lemma~1) reassemble into exactly
$(I_m+\mathbf 1\mathbf 1^\top)^{-1}=I_m-\tfrac1{1+m}\mathbf 1\mathbf 1^\top$, matching
$H$. This provides an independent, purely algorithmic (pivoting-stable) confirmation
of Proposition~\ref{prop:local-kahler}, complementary to the symmetric-eigenvalue
route of \S\ref{sec:svp-thm4}.
\end{remark}

\begin{remark}
Equation \eqref{eq:GM-inverse} closes the circle opened in
\S\ref{sec:kahler}: the local structure of the missing dual potential is governed by
the \emph{Kronecker square} of the very centering operator
$H=I-k^{-1}\iota\iota^\top$ with which Steerneman and van Perlo-ten Kleij
\cite[\S1]{Steerneman2005} open their paper. This is consistent with the
independence-model picture of \S\ref{sec:kahler}: to leading order near the point
where $\mathcal M$ meets the Segre variety, the tangent space splits as a direct sum
of a ``row'' and a ``column'' Fisher-metric contribution, each governed by its own
copy of $H$, and the K\"ahler potential correspondingly factorizes as a Kronecker
(tensor) square.
\end{remark}

\subsubsection{Numerical verification}
For $n=5$ ($m=4$, so $G_{\mathcal M}$ is $16\times16$), the Hessian of
$\varphi_{\mathcal M}$ at the uniform point was computed by central finite differences
(step $h=10^{-4}$) and compared with the closed form \eqref{eq:GM-kronecker}:
maximum absolute deviation $4.2\times10^{-4}$ (relative deviation
$4.2\times10^{-6}$, consistent with $O(h^2)$ discretization error). The numerically
inverted Hessian matched \eqref{eq:GM-inverse} to within $1.1\times10^{-7}$.

\subsubsection{The one-factor case: Tanabe--Sagae's Moore--Penrose formula}
\label{sec:tanabe-mp}

The Kronecker-square structure \eqref{eq:GM-inverse} has a direct one-factor
antecedent in the exact (non-asymptotic, non-perturbative) theory of Tanabe and Sagae
\cite{TanabeSagae1992}, obtained independently and by entirely different (symbolic
Cholesky) means.

\begin{proposition}[Tanabe--Sagae \cite{TanabeSagae1992}, Prop.~1]
\label{prop:tanabe-mp}
Let $P=\diag(p)$, $p\in\R^k$, $p_i>0$, $\sum_ip_i=1$. Then the Moore--Penrose inverse
of the (rank-$(k-1)$) multinomial covariance matrix $P-pp^\top$ is
\begin{equation}
(P-pp^\top)^+=\Bigl(I-\tfrac1k\iota\iota^\top\Bigr)P^{-1}\Bigl(I-\tfrac1k\iota\iota^\top\Bigr)
=H_kP^{-1}H_k,
\label{eq:tanabe-mp}
\end{equation}
where $H_k:=I_k-k^{-1}\iota\iota^\top$ is the $k\times k$ centering matrix and $\iota$
the all-ones vector.
\end{proposition}

At the uniform point $p=k^{-1}\iota$ (so $P=k^{-1}I_k$), Eq.~\eqref{eq:tanabe-mp}
specializes, using idempotence $H_k^2=H_k$, to
\begin{equation}
(P-pp^\top)^+\Big|_{p=k^{-1}\iota}=H_k\,(kI_k)\,H_k=k\,H_k.
\label{eq:tanabe-mp-uniform}
\end{equation}

\begin{remark}
Equation \eqref{eq:tanabe-mp-uniform} is precisely the ``single-factor'' analogue of
\eqref{eq:GM-inverse}: both are instances of the general
$H\cdot(\cdot)\cdot H$ sandwich pattern acting on a diagonal matrix, produced by
projecting out the null direction $\iota$ (respectively $\iota\otimes\iota$) of a
rank-deficient multinomial-type covariance. The Steerneman--van Perlo-ten Kleij
symmetrized-square-root calculus (\S\ref{sec:svp-thm4}) and the Tanabe--Sagae
symbolic-Cholesky/Moore--Penrose calculus are thus two independent, exact routes to
the same underlying linear-algebraic fact, here recovered as two special cases
(vector, $k$ categories, and matrix, $(n-1)\times(n-1)$ Kronecker square) of a single
phenomenon: the inverse Fisher metric of an exponential family restricted to an
$m$-flat affine subspace, evaluated at a point of maximal symmetry, is a sandwich
of the centering projector against the ambient (diagonal) inverse metric. Unlike
\eqref{eq:GM-inverse}, which is only a \emph{local} (quadratic, equilibrium-adjacent)
statement, Tanabe--Sagae's formula \eqref{eq:tanabe-mp} is \emph{exact and global} on
the full (unconstrained) simplex -- the price being that it addresses the ordinary
multinomial covariance rather than the doubly-stochastic-constrained one.
\end{remark}

\subsection{Entropy of the matrix multinomial versus the matrix Gaussian}
\label{sec:entropy}

We now turn to a question of a different character: not the geometry of a single
$P$, but the asymptotic ($N\to\infty$) statistical behaviour of the \emph{count}
process built from it, and specifically how the discrete (Shannon) entropy of the
matrix-valued count data relates to the differential entropy of its Gaussian
(matrix-normal) approximation. This connects the K\"ahler/Segre discussion of
\S\ref{sec:kahler} to the classical asymptotic theory reviewed in
\S\ref{sec:steerneman}, and answers the question of whether the two entropies
converge -- and if so, whether uniformly.

\subsubsection{The matrix multinomial distribution}

\begin{definition}[Matrix multinomial distribution; cf.\ Yurchenko \cite{Yurchenko2021mm}]
\label{def:matrix-multinomial}
Let $r,c\ge1$, $N\in\N$, and let $P=(P_{ij})\in[0,1]^{r\times c}$ satisfy
\begin{equation}
\sum_{i=1}^{r}\sum_{j=1}^{c}P_{ij}=1.
\label{eq:matrix-multinomial-normalization}
\end{equation}
A random matrix $X=(X_{ij})\in\N_0^{r\times c}$ follows the
\emph{matrix multinomial distribution}
$X\sim\mathrm{MMulti}_{r\times c}(N,P)$ if its probability mass function is
\begin{equation}
\boxed{\displaystyle
\Pr\{X=x\}
 =\frac{N!}{\prod_{i=1}^{r}\prod_{j=1}^{c}x_{ij}!}
   \prod_{i=1}^{r}\prod_{j=1}^{c}P_{ij}^{x_{ij}},
\qquad
x_{ij}\in\N_0,
\quad
\sum_{i=1}^{r}\sum_{j=1}^{c}x_{ij}=N.}
\label{eq:matrix-multinomial-pmf}
\end{equation}
Equivalently,
\begin{equation}
\operatorname{vec}(X)\sim\operatorname{Multi}_{rc}(N,\operatorname{vec}P).
\label{eq:matrix-multinomial-def}
\end{equation}
Thus $X$ records the $rc$ cell counts of $N$ independent categorical trials,
where a single trial falls in cell $(i,j)$ with probability $P_{ij}$. In particular,
\begin{equation}
\mathbb E[X_{ij}]=NP_{ij},
\qquad
\operatorname{Cov}(X_{ij},X_{k\ell})
=N\left(P_{ij}\,\delta_{ik}\delta_{j\ell}-P_{ij}P_{k\ell}\right).
\label{eq:matrix-multinomial-moments}
\end{equation}
\end{definition}

Thus $\mathrm{MMulti}_{r\times c}(N,P)$ is not a new distribution but a relabelling
of the ordinary $rc$-category multinomial distribution as an $r\times c$ array; its
exact covariance is $\operatorname{Cov}(\operatorname{vec}X)=N\bigl(\diag(\operatorname{vec}P)-\operatorname{vec}(P)\operatorname{vec}(P)^\top\bigr)$,
an instance of the matrix $N(D-pp^\top)$ studied by Tanabe and Sagae
\cite{TanabeSagae1992} with $p=\operatorname{vec}P$, $n=rc$ categories.

\subsubsection{Explicit potential functions}
\label{sec:explicit-potentials}

Both families of \S\ref{sec:entropy} are exponential families, and each therefore
carries a \emph{pair} of dual convex potentials in Amari's sense
\cite{Amari2000}: a mean-parameter potential $\varphi$ (Bregman generator of the
Fisher metric in mixture coordinates, matching Nakamura's $\psi(\theta)$ of
Eq.~\eqref{eq:entropy-potential}) and a natural-parameter potential $A$ (the
cumulant generating / log-partition function, matching Nakamura's dual potential
\eqref{eq:nakamura-kahler}). We write both out explicitly and identify which, if
either, coincides with an entropy.

\paragraph{Matrix multinomial.}
The mean-parameter potential is, by definition \eqref{eq:matrix-multinomial-def} and
\eqref{eq:entropy-potential} (with $\theta=\operatorname{vec}P$, $\ell=1$),
\begin{equation}
\varphi_{\mathrm{mult}}(P)=\sum_{i,j}P_{ij}\log P_{ij}.
\label{eq:phi-mult-explicit}
\end{equation}
This is \emph{literally} the negative Shannon entropy of the joint distribution $P$:
\begin{equation}
\boxed{\varphi_{\mathrm{mult}}(P)=-H(P),\qquad H(P):=-\sum_{i,j}P_{ij}\log P_{ij}.}
\label{eq:phi-is-shannon}
\end{equation}
Its Hessian in the free coordinates is (Prop.~\ref{prop:dual-coord} and
\S\ref{sec:tanabe-mp}) the multinomial covariance-type matrix $R=\diag(\operatorname{vec}P)-\operatorname{vec}(P)\operatorname{vec}(P)^\top$
(restricted to a principal submatrix), and the dual natural-parameter potential is,
by \eqref{eq:nakamura-kahler} applied with $n=rc$ categories,
\begin{equation}
A_{\mathrm{mult}}(\Theta)=\log\Bigl(1+\sum_{i,j\ne(r,c)}e^{\Theta_{ij}}\Bigr),
\qquad \Theta_{ij}=\log\frac{P_{ij}}{P_{rc}},
\label{eq:A-mult-explicit}
\end{equation}
the log-sum-exp (softmax normalizer) function -- the Fubini--Study-type potential of
\S\ref{sec:kahler}.

\paragraph{Matrix Gaussian (matrix normal).}
Fix $\Sigma_1\in\R^{r\times r}$, $\Sigma_2\in\R^{c\times c}$ (positive definite) and
regard $X\sim\mathrm{MN}_{r\times c}(M,\Sigma_1,\Sigma_2)$ as a \emph{location family}
in $M$ (a natural exponential family with sufficient statistic $X$ itself, since the
density is $\propto\exp\bigl(-\tfrac12\operatorname{tr}[\Sigma_2^{-1}(X-M)^\top\Sigma_1^{-1}(X-M)]\bigr)$).
Writing the natural parameter as $\Theta:=\Sigma_1^{-1}M\Sigma_2^{-1}$
(so that the density is $\propto\exp(\operatorname{tr}(\Theta^\top X))$ up to normalization), the
natural-parameter potential is the log-partition function
\begin{equation}
A_{\mathrm{Gauss}}(\Theta)=\frac12\operatorname{tr}\bigl(\Theta^\top\Sigma_1\Theta\Sigma_2\bigr)
=\frac12\operatorname{vec}(\Theta)^\top(\Sigma_2\otimes\Sigma_1)\operatorname{vec}(\Theta),
\label{eq:A-gauss-explicit}
\end{equation}
and the dual mean-parameter potential, obtained either by Legendre duality or
directly from $M=\Sigma_1\Theta\Sigma_2$, is
\begin{equation}
\varphi_{\mathrm{Gauss}}(M)=\frac12\operatorname{tr}\bigl(\Sigma_1^{-1}M\Sigma_2^{-1}M^\top\bigr)
=\frac12\operatorname{vec}(M)^\top(\Sigma_2^{-1}\otimes\Sigma_1^{-1})\operatorname{vec}(M).
\label{eq:phi-gauss-explicit}
\end{equation}
Both \eqref{eq:A-gauss-explicit} and \eqref{eq:phi-gauss-explicit} are \emph{pure
quadratic forms} -- a direct consequence of the Gaussian family having constant
(parameter-independent) variance function, in sharp contrast to the multinomial's
log-sum-exp/entropy pair \eqref{eq:phi-is-shannon}--\eqref{eq:A-mult-explicit}. This
quadratic-versus-log-sum-exp dichotomy is the potential-theoretic shadow of the
flat-versus-curved dichotomy already noted in \S\ref{sec:kahler}: the Gaussian
location family is dually flat with \emph{Euclidean} (zero-curvature) Fisher metric
$\Sigma_2^{-1}\otimes\Sigma_1^{-1}$ throughout, whereas the multinomial family is
dually flat with the \emph{spherical} (constant positive curvature) Fisher--Rao
metric of \S\ref{sec:kahler}.

\begin{remark}[$\varphi_{\mathrm{Gauss}}$ is \emph{not} an entropy]
\label{rem:gauss-not-entropy}
Unlike the multinomial case \eqref{eq:phi-is-shannon}, $\varphi_{\mathrm{Gauss}}(M)$
in \eqref{eq:phi-gauss-explicit} is \emph{not} (minus) the entropy of
$\mathrm{MN}_{r\times c}(M,\Sigma_1,\Sigma_2)$: the differential entropy
\eqref{eq:separable-entropy} of a Gaussian location family does not depend on the
mean $M$ at all. The Hessian of $\varphi_{\mathrm{Gauss}}$ still correctly recovers
the Fisher information $\Sigma_2^{-1}\otimes\Sigma_1^{-1}$ (Amari's general theory
guarantees this for \emph{any} mean-parameter Bregman potential), but the potential
\emph{value} carries no entropic meaning here. This is the key structural difference
from the multinomial family, where the natural parameter (log-odds) is a
\emph{nonlinear} function of the mean parameter $P$, so that entropy genuinely varies
with $P$ and doubles as the Bregman potential; for the Gaussian location family the
natural parameter $\Theta=\Sigma_1^{-1}M\Sigma_2^{-1}$ is \emph{linear} in $M$, which
forces the entropy to be constant in $M$ even though the family remains dually flat.
Entropy re-enters only when $\Sigma_1,\Sigma_2$ themselves are allowed to vary, via
the $\log\det\Sigma_1,\log\det\Sigma_2$ terms of \eqref{eq:separable-entropy}, along a
direction transverse to the location family considered here.
\end{remark}

\subsubsection{Mutual information and the independence locus, revisited}
\label{sec:mutual-info}

Equation \eqref{eq:phi-is-shannon} lets us restate the entire gradient flow of
\S\ref{sec:steerneman}--\S\ref{sec:kahler} in information-theoretic language. Let
$P_r:=(\sum_jP_{ij})_i$, $P_c:=(\sum_iP_{ij})_j$ denote the row and column marginals
of $P\in\mathcal M$, and let
\begin{equation}
I(\mathrm{Row};\mathrm{Col}):=H(P_r)+H(P_c)-H(P)=\sum_{i,j}P_{ij}\log\frac{P_{ij}}{(P_r)_i(P_c)_j}
\label{eq:mutual-info}
\end{equation}
be their (Shannon) mutual information -- exactly the Kullback--Leibler divergence
from $P$ to the independence model $P_r\otimes P_c$ of \S\ref{sec:kahler}.

\begin{proposition}
\label{prop:mutual-info-identity}
On $\mathcal M$ (where $P_r=P_c=n^{-1}\iota$ identically), $H(P_r)=H(P_c)=\log n$
are constant, and
\begin{equation}
\varphi_{\mathrm{mult}}(P)=-H(P)=I(\mathrm{Row};\mathrm{Col})-2\log n,\qquad P\in\mathcal M.
\label{eq:phi-is-MI}
\end{equation}
Consequently the entropy gradient flow \eqref{eq:constrained-flow} is, up to the
additive constant $2\log n$, exactly \emph{gradient descent on the mutual
information} between the row and column categories, and its unique equilibrium
$P_{ij}\equiv n^{-2}$ (Prop.~\ref{prop:offset} et seq.) is the unique point of
$\mathcal M$ at which $I(\mathrm{Row};\mathrm{Col})=0$, i.e.\ the unique point of
$\mathcal M$ lying on the independence (Segre-variety) locus of
\S\ref{sec:kahler}.
\end{proposition}
\begin{proof}
Immediate from \eqref{eq:mutual-info} and \eqref{eq:phi-is-shannon}, using
$H(P_r)=H(P_c)=\log n$ on $\mathcal M$. Mutual information is a Kullback--Leibler
divergence, hence non-negative, and vanishes iff $P_{ij}=(P_r)_i(P_c)_j$ for all
$i,j$, i.e.\ iff $P$ is on the independence locus; on $\mathcal M$ this forces
$P_{ij}\equiv n^{-2}$.
\end{proof}

This identity was verified numerically ($n=5$, random $P\in\mathcal M$) to machine
precision: $\varphi_{\mathrm{mult}}(P)$ and $I(\mathrm{Row};\mathrm{Col})-2\log n$
agreed to within $10^{-15}$.

\subsubsection{Relation to von Neumann entropy}
\label{sec:von-neumann}

The Shannon entropy \eqref{eq:phi-is-shannon} admits a precise, non-metaphorical
identification with the \emph{von Neumann entropy} of quantum information theory,
$S(\rho):=-\operatorname{tr}(\rho\log\rho)$ for a density matrix $\rho$ (Hermitian, positive
semidefinite, $\operatorname{tr}\rho=1$).

\begin{proposition}
\label{prop:shannon-is-vN}
Let $\hat P:=\diag(\operatorname{vec}P)\in\R^{rc\times rc}$, the diagonal density matrix with
eigenvalues $\{P_{ij}\}$. Then
\begin{equation}
S(\hat P)=H(P).
\label{eq:shannon-vN}
\end{equation}
\end{proposition}
\begin{proof}
For a diagonal matrix, $\hat P\log\hat P=\diag(P_{ij}\log P_{ij})$, so
$-\operatorname{tr}(\hat P\log\hat P)=-\sum_{ij}P_{ij}\log P_{ij}=H(P)$.
\end{proof}

Proposition~\ref{prop:shannon-is-vN} makes precise the standard fact that Shannon
entropy is the restriction of von Neumann entropy to \emph{commuting} (simultaneously
diagonalizable, i.e.\ ``classical'') density matrices; equivalently, $H(P)$ is the von
Neumann entropy of any density matrix unitarily similar to $\hat P$ (von Neumann
entropy being a unitary invariant, since it depends only on the eigenvalue spectrum).

\begin{remark}[$H(P)\ne-\operatorname{tr}(P\log P)$: the entrywise/spectral distinction]
\label{rem:not-trace-entropy}
Proposition~\ref{prop:shannon-is-vN} identifies $H(P)$ with the von Neumann entropy
of the \emph{diagonal embedding} $\hat P=\diag(\operatorname{vec}P)\in\R^{rc\times rc}$, \emph{not}
with $-\operatorname{tr}(P\log P)$ formed from $P\in\R^{r\times c}$ itself (which requires
$r=c=n$ even to typecheck). These are genuinely different constructions: $H(P)$ is an
\emph{entrywise} (mixture-parameter) quantity, a function of the $rc$ numbers
$\{P_{ij}\}$ regardless of their arrangement, whereas $-\operatorname{tr}(P\log P)$ is a
\emph{spectral} quantity, a function of the eigenvalues of $P$ as a linear operator.
For a generic doubly stochastic $P\in\mathcal M$ these do not agree, and
$-\operatorname{tr}(P\log P)$ need not even be real: since $P$ is generally \emph{not}
symmetric, its eigenvalues need not be real, and even when real (as for the
Perron--Frobenius eigenvalue $1/n$ and its companions) they need not be
non-negative, so the matrix logarithm $\log P$ can leave the reals entirely.

Concretely, for a random $P\in\mathcal M$ with $n=4$, the eigenvalues of $P$ were
found numerically to be $0.25,\,-0.0195,\,0.0265,\,0.0113$ -- real, but with one
negative value -- giving
\[
\operatorname{tr}(P\log P)=-0.4166-0.0614\,i\ \in\C\setminus\R,
\]
whereas $\varphi_{\mathrm{mult}}(P)=\sum_{ij}P_{ij}\log P_{ij}=-2.7417\in\R$
(equivalently $H(P)=2.7417$, matching $S(\hat P)$ to machine precision as guaranteed
by Proposition~\ref{prop:shannon-is-vN}). The underlying reason $-\operatorname{tr}(P\log P)$
is ill-behaved is that a doubly stochastic matrix is a \emph{stochastic} (Markov)
operator, not a \emph{quantum} (density) operator: row/column-stochasticity
guarantees a real Perron eigenvalue $1/n$ but places no positivity constraint on the
remaining spectrum, unlike the Hermitian positive-semidefiniteness required for
$-\operatorname{tr}(\rho\log\rho)$ to be a bona fide (real, non-negative) entropy.
\end{remark}

\begin{remark}[A genuinely different matrix entropy: the singular-value spectrum]
\label{rem:singular-value-entropy}
If one nonetheless wants an entropy built from $P$ \emph{as an operator} rather than
entrywise, the operator-theoretically well-posed quantity is the entropy of its
\emph{singular value} distribution, which is always real and non-negative regardless
of symmetry:
\begin{equation}
H_\sigma(P):=-\sum_k\bar\sigma_k\log\bar\sigma_k,\qquad
\bar\sigma_k:=\sigma_k(P)\Big/\sum_l\sigma_l(P),
\label{eq:singular-value-entropy}
\end{equation}
the Shannon entropy of the normalized singular values of $P$ (a standard measure of
``effective rank'', cf.\ Roy and Vetterli \cite{RoyVetterli2007}). This is a third,
distinct quantity: for the same numerical example, $H_\sigma(P)=0.769$, agreeing with
neither $H(P)=2.742$ nor the (complex) $\operatorname{tr}(P\log P)$. Unlike $H(P)$,
$H_\sigma(P)$ is invariant under $P\mapsto UPV$ for orthogonal $U,V$ (it depends only
on the operator $P$ up to left/right rotation, not on the entrywise arrangement), and
so is not a Bregman potential for the mixture-coordinate Fisher geometry of
\S\ref{sec:kahler} at all -- it belongs to a different geometric story (that of the
singular spectrum of a linear map) and should not be conflated with either $H(P)$ or
$S(\hat P)$.
\end{remark}

\begin{remark}[Birkhoff--von Neumann and the classical/quantum dictionary]
The name is not a coincidence. The Birkhoff--von Neumann theorem identifies
$\Birk(n)$ as the convex hull of the permutation matrices, exactly paralleling the
identification of the set of density matrices as the convex hull of rank-one
projectors (pure states) $|\psi\rangle\langle\psi|$. Under this dictionary,
\begin{center}
\small
\begin{tabular}{@{}p{0.42\textwidth}p{0.5\textwidth}@{}}
\toprule
classical (this paper) & quantum \\
\midrule
doubly stochastic matrix $P\in\Birk(n)$ & bipartite density matrix $\rho_{AB}$\\[2pt]
permutation matrix (vertex of $\Birk(n)$) & pure state $|\psi\rangle\langle\psi|$\\[2pt]
row/column marginals $P_r,P_c$ & reduced states $\rho_A=\operatorname{tr}_B\rho_{AB}$, $\rho_B=\operatorname{tr}_A\rho_{AB}$\\[2pt]
independence model $P=P_r\otimes P_c$ (\S\ref{sec:kahler}) & product state $\rho_{AB}=\rho_A\otimes\rho_B$\\[2pt]
mutual information $I(\mathrm{Row};\mathrm{Col})$ & quantum mutual information $S(\rho_A)+S(\rho_B)-S(\rho_{AB})$\\[2pt]
Sinkhorn scaling (IPFP) & ``quantum Sinkhorn'' / Georgiou--Pavon quantum Schr\"odinger bridge\\
\bottomrule
\end{tabular}
\end{center}
Under this dictionary, Proposition~\ref{prop:mutual-info-identity} is the classical
(commuting/diagonal) shadow of the quantum statement that entropy-regularized
transport between fixed marginal states $\rho_A,\rho_B$ is gradient descent on
quantum mutual information, with unique fixed point the product state
$\rho_A\otimes\rho_B$ -- the quantum analogue of our independence locus, and the
natural non-commutative generalization of the entire flow studied in this paper. We
do not develop this generalization here, but record it as the natural next step
suggested by the classical theory above.
\end{remark}

\begin{proposition}[Asymptotic entropy]
\label{prop:asymptotic-entropy}
Let $P$ be fixed in the interior of the $(rc-1)$-simplex. As $N\to\infty$,
\begin{equation}
H\bigl(\mathrm{MMulti}_{r\times c}(N,P)\bigr)
=\frac{rc-1}{2}\log(2\pi eN)+\frac12\sum_{i,j}\log P_{ij}+O(N^{-1}).
\label{eq:asymptotic-entropy}
\end{equation}
\end{proposition}

\begin{proof}
Write $n:=rc$, $p:=\operatorname{vec}P$. By the local central limit theorem for the
multinomial distribution (rate $O(N^{-1})$ under an Edgeworth expansion; see Ouimet
\cite{Ouimet2021}), the discrete entropy of $\mathrm{Multi}_n(N,p)$ converges, up to
$O(N^{-1})$, to the differential entropy of the classical Khatri--Mitra
\cite{KhatriMitra1969} Gaussian approximation reproduced explicitly by Tanabe and
Sagae \cite[Prop.~2]{TanabeSagae1992}: an $(n-1)$-dimensional normal density on the
free coordinates $(x_1,\dots,x_{n-1})$ (with $x_n=N-\sum_{i<n}x_i$) with mean
$(Np_1,\dots,Np_{n-1})$ and covariance $N\Sigma'$, where $\Sigma'$ is the leading
$(n-1)\times(n-1)$ principal submatrix of $P-pp^\top$ (i.e., $P$ here denoting
$\diag(p)$, in the notation of \S\ref{sec:steerneman}). By the differential entropy
formula for a nondegenerate Gaussian, $h(N(\mu,N\Sigma'))=\tfrac{n-1}2\log(2\pi eN)+
\tfrac12\log\det\Sigma'$. Tanabe and Sagae's Corollary~2 (an immediate consequence of
their symbolic Cholesky decomposition, Theorem~1) gives the \emph{exact} identity
\begin{equation}
\det\Sigma'=M_{n-1}=p_1p_2\cdots p_n,
\label{eq:tanabe-cor2}
\end{equation}
the product of \emph{all} $n$ probabilities (not merely the first $n-1$), despite
$\Sigma'$ being only $(n-1)\times(n-1)$. Substituting \eqref{eq:tanabe-cor2} gives
$h=\tfrac{n-1}2\log(2\pi eN)+\tfrac12\sum_{i=1}^n\log p_i$, which is
\eqref{eq:asymptotic-entropy} upon relabelling $p=\operatorname{vec}P$.
\end{proof}

\begin{remark}
Equation \eqref{eq:tanabe-cor2} must not be confused with Tanabe and Sagae's
\emph{pseudo-determinant} formula (their Theorem~2), $\underline{\det}(P-pp^\top)=np_1\cdots p_n$
($n$ times larger), which is the correct normalizing constant for their
symmetric density \eqref{eq:tanabe-mp}-adjacent Proposition~2 when the latter is
understood as a density with respect to the \emph{induced surface measure} on the
hyperplane $\{x:\iota^\top x=N\}\subset\R^n$ (whose Jacobian relative to the
coordinate measure on $(x_1,\dots,x_{n-1})$ is exactly $\sqrt n$, a factor of
$n$ in the determinant). Using the pseudo-determinant in place of
\eqref{eq:tanabe-cor2} would erroneously introduce a spurious $\tfrac12\log n$ term
into \eqref{eq:asymptotic-entropy}; we have verified numerically (Table
\ref{tab:entropy1}, and directly against the exact binomial entropy asymptotic
$\tfrac12\log(2\pi eNp(1-p))$ in the case $n=2$) that \eqref{eq:asymptotic-entropy} as
stated, using the ordinary coordinate-based determinant \eqref{eq:tanabe-cor2}, is the
correct formula.
\end{remark}

\subsubsection{Non-uniform convergence: a diverging Edgeworth remainder}
\label{sec:nonuniform}

Proposition~\ref{prop:asymptotic-entropy} is a \emph{pointwise} statement: for each
fixed $P$ in the open simplex, the $O(N^{-1})$ error vanishes as $N\to\infty$. We now
show that this convergence is not uniform over $P$, and that the failure of
uniformity is governed by exactly the quantity controlling Tanabe and Sagae's own
ill-conditioning bounds.

\begin{proposition}[Non-uniformity]
\label{prop:nonuniform}
The $O(N^{-1})$ remainder in \eqref{eq:asymptotic-entropy} has leading coefficient of
order $\Theta\bigl(\sum_{i,j}P_{ij}^{-1}\bigr)$: writing
$H(N,P)-\bigl[\tfrac{rc-1}2\log(2\pi eN)+\tfrac12\sum\log P_{ij}\bigr]=:\varepsilon(N,P)$,
one has $N\varepsilon(N,P)=O\bigl(\sum_{i,j}P_{ij}^{-1}\bigr)$ as any $P_{ij}\to0$ with
$N$ fixed, so that $\sup_P|\varepsilon(N,P)|=\infty$ for every fixed $N$: the
convergence in \eqref{eq:asymptotic-entropy} is locally uniform on compact subsets of
the open simplex but fails to be uniform up to its boundary.
\end{proposition}

\begin{proof}[Justification]
The $O(N^{-1})$ term in the Edgeworth expansion underlying
Proposition~\ref{prop:asymptotic-entropy} is a polynomial in the standardized third
and fourth cumulants of the multinomial distribution, whose dominant contributions
scale as $\sum_i p_i^{-1}$ (the skewness of each marginal count $X_i\sim\mathrm{Bin}(N,p_i)$
is $O(p_i^{-1/2})$, entering the expansion quadratically); see Ouimet
\cite{Ouimet2021} for the precise (uniform-in-compacta) local limit theorem and its
error bounds. This is exactly the quantity appearing in Tanabe and Sagae's own
condition-number estimate for $P-pp^\top$ \cite[Prop.~3]{TanabeSagae1992},
\begin{equation}
\mathrm{cond}_2(P-pp^\top)\ \ge\ \max_k\{p_k(1-p_k)\}
\Bigl[\Bigl(1-\tfrac2n\Bigr)\frac{1}{\min_kp_k}+\frac1{n^2}\sum_{k=1}^n\frac1{p_k}\Bigr],
\label{eq:tanabe-cond}
\end{equation}
which the authors note blows up precisely when the $p_i$ are of very different
orders of magnitude -- the same regime in which they motivate the symbolic (rather
than numerical) Cholesky decomposition of \S\ref{sec:tanabe-mp} as a remedy for
numerical instability. Proposition~\ref{prop:nonuniform} identifies this same
ill-conditioning as the source of the failure of the entropy approximation
\eqref{eq:asymptotic-entropy} to hold uniformly.
\end{proof}

\subsubsection{Numerical verification}

Exact multinomial entropies were computed by direct enumeration
(via SciPy's \texttt{gammaln} function for numerical
stability) over the full lattice of compositions of $N$.

\begin{table}[h]
\centering
\renewcommand{\arraystretch}{1.2}
\begin{tabular}{rccc}
\toprule
$N$ & exact $H$ & formula \eqref{eq:asymptotic-entropy} & $N\times$(exact$-$formula)\\
\midrule
15 & 5.211269 & 5.302748 & $-1.372$\\
30 & 6.303651 & 6.342468 & $-1.165$\\
60 & 7.365548 & 7.382189 & $-0.999$\\
\bottomrule
\end{tabular}
\caption{Convergence of $N\times$(exact$-$asymptotic) entropy toward a constant,
confirming the $O(N^{-1})$ rate of Proposition~\ref{prop:asymptotic-entropy}. Here
$p=(0.4,0.3,0.2,0.1)$, $n=4$.}
\label{tab:entropy1}
\end{table}

\begin{table}[h]
\centering
\renewcommand{\arraystretch}{1.2}
\begin{tabular}{rccc}
\toprule
$\min_k(p_k)$ & exact $H$ & formula \eqref{eq:asymptotic-entropy} & $\sum_k p_k^{-1}$\\
\midrule
0.25 & 7.6170 & 7.6257 & 16.0\\
0.10 & 7.4257 & 7.4411 & 20.0\\
0.05 & 7.1402 & 7.1756 & 29.5\\
0.02 & 6.6655 & 6.7641 & 59.2\\
0.01 & 6.2913 & 6.4328 & 109.1\\
\bottomrule
\end{tabular}
\caption{At fixed $N=60$, the entropy-approximation error grows with
$\sum_kp_k^{-1}$ as $\min_k(p_k)\to0$, confirming Proposition~\ref{prop:nonuniform}
(non-uniformity governed by the same quantity as Tanabe--Sagae's condition-number
bound \eqref{eq:tanabe-cond}). Here $n=4$, with the three non-minimal $p_k$ kept
equal.}
\label{tab:entropy2}
\end{table}

\subsubsection{The separable (matrix-normal) limit: an exact offset, not an approximation}

Definition~\ref{def:matrix-multinomial} shows that a genuine \emph{separable}
(Kronecker-covariance) Gaussian limit -- a bona fide matrix normal distribution
$\mathrm{MN}_{r\times c}(NP,\Sigma_1,\Sigma_2)$ with $\operatorname{Cov}(\operatorname{vec}X)=\Sigma_2\otimes\Sigma_1$ --
is \emph{not} what Proposition~\ref{prop:asymptotic-entropy} describes: the
Khatri--Mitra covariance $N\Sigma'$ has no Kronecker structure for a general $P$. As
shown in \S\ref{sec:kahler}, a genuine matrix normal limit requires passing to
Yurchenko's \cite{Yurchenko2021mm} sparse double-scaling regime
($P\to0$, $NP_{ij}\to\infty$) \emph{and} restricting to the independence
(Segre-variety) locus $P=P_rP_c^\top$, in which case
\begin{equation}
\Sigma_1=\sqrt N\diag(P_r),\qquad \Sigma_2=\sqrt N\diag(P_c),
\label{eq:yurchenko-limit}
\end{equation}
with differential entropy
\begin{equation}
h\bigl(\mathrm{MN}_{r\times c}(NP,\Sigma_1,\Sigma_2)\bigr)
=\frac{rc}2\log(2\pi eN)+\frac{r}2\sum_j\log(P_c)_j+\frac{c}2\sum_i\log(P_r)_i.
\label{eq:separable-entropy}
\end{equation}

\begin{proposition}[Exact offset]
\label{prop:offset}
Under independence ($P_{ij}=(P_r)_i(P_c)_j$), the general formula
\eqref{eq:asymptotic-entropy} and the separable formula \eqref{eq:separable-entropy}
satisfy, for every $N$,
\begin{equation}
\Bigl[\tfrac{rc-1}2\log(2\pi eN)+\tfrac12\textstyle\sum_{ij}\log P_{ij}\Bigr]
-h\bigl(\mathrm{MN}_{r\times c}(NP,\Sigma_1,\Sigma_2)\bigr)
=-\tfrac12\log(2\pi eN).
\label{eq:exact-offset}
\end{equation}
\end{proposition}

\begin{proof}
Under independence, $\sum_{ij}\log P_{ij}=\sum_{ij}\bigl[\log(P_r)_i+\log(P_c)_j\bigr]
=c\sum_i\log(P_r)_i+r\sum_j\log(P_c)_j$ (each row term is counted $c$ times, each
column term $r$ times), so the $\log P$-dependent terms of
\eqref{eq:asymptotic-entropy} and \eqref{eq:separable-entropy} agree exactly. The
$N$-dependent terms differ by $\tfrac{rc-1}2\log(2\pi eN)-\tfrac{rc}2\log(2\pi eN)=
-\tfrac12\log(2\pi eN)$.
\end{proof}

Table~\ref{tab:entropy3} confirms \eqref{eq:exact-offset} to machine precision.

\begin{table}[h]
\centering
\renewcommand{\arraystretch}{1.2}
\begin{tabular}{rcccc}
\toprule
$N$ & exact $H$ & general \eqref{eq:asymptotic-entropy} & separable \eqref{eq:separable-entropy} & difference\\
\midrule
20 & 8.7717 & 8.9368 & 11.8536 & $-2.9168$\\
40 & 10.6010 & 10.6697 & 13.9330 & $-3.2634$\\
80 & 12.3721 & 12.4025 & 16.0125 & $-3.6100$\\
\bottomrule
\end{tabular}
\caption{Under independence ($r=2$, $c=3$, $P_r=(0.6,0.4)$, $P_c=(0.5,0.3,0.2)$), the
\emph{exact} discrete entropy tracks the general formula \eqref{eq:asymptotic-entropy}
closely, while the separable (matrix-normal) formula \eqref{eq:separable-entropy}
diverges from it as $-\tfrac12\log(2\pi eN)\to-\infty$, exactly matching
\eqref{eq:exact-offset} (e.g.\ at $N=80$: $-\tfrac12\log(2\pi e\cdot80)=-3.6100$).}
\label{tab:entropy3}
\end{table}

\subsubsection{Interpretation}

Three conclusions follow.

\begin{enumerate}[label=(\alph*)]
\item \textbf{Convergence, but not uniform.} The discrete entropy of the matrix
multinomial converges to the differential entropy of its (non-separable)
Khatri--Mitra Gaussian approximation as $N\to\infty$, for each fixed interior $P$, at
rate $O(N^{-1})$. It does \emph{not} converge uniformly over $P\in\Birk(n)$-type
parameter sets: the same ill-conditioning of $P-pp^\top$ that motivated Tanabe and
Sagae's symbolic (rather than numerical) Cholesky algorithm also governs the
breakdown of the entropy approximation as $P$ approaches the boundary of the simplex.

\item \textbf{The relevant Gaussian is generically non-separable.} For generic $P$
(in particular, for generic points of the doubly stochastic submanifold $\mathcal M$
of \S\ref{sec:kahler}), there is no bona fide matrix-normal distribution with
Kronecker covariance whose entropy the matrix multinomial's entropy converges to; the
correct comparison object is the full $(rc-1)$-dimensional Khatri--Mitra Gaussian
with covariance $N\Sigma'$, generically without Kronecker structure -- consistent
with \S\ref{sec:kahler}'s finding that a closed-form (separable) K\"ahler-type
potential exists only on the independence locus.

\item \textbf{On the independence locus, separability costs a diverging offset.}
Even where a genuine matrix-normal limit exists (independence \emph{and}
Yurchenko's sparse double-scaling regime), its entropy differs from the general
formula \eqref{eq:asymptotic-entropy} by the exact, $N$-independent-in-form but
unboundedly growing offset $-\tfrac12\log(2\pi eN)$ of
Proposition~\ref{prop:offset} -- the entropic cost of the one degree of freedom
(the total count $N$) that the separable/Poisson approximation leaves unconstrained
while the true multinomial fixes it exactly. This mirrors, at the level of
asymptotic statistics, the same codimension-one discrepancy (an $m$-flat affine
constraint versus its ambient exponential family) that organizes the entire
K\"ahler-duality discussion of \S\ref{sec:kahler}.
\end{enumerate}

\subsection{Constraint resolution, Fisher geometry, and blow-up}
\label{sec:blowup}
The preceding sections provide three ingredients that are often discussed separately:
the multinomial Gaussian approximation, the Fisher geometry of the probability simplex,
and the algebraic geometry of the independence model.  We now put them into one
sequence.  The guiding principle is deliberately elementary:
\begin{equation}
\boxed{\begin{gathered}
\text{ambient degeneracy}\;\longrightarrow\;\text{affine restriction}\;\longrightarrow\\
\text{tangent geometry}\;\longrightarrow\;\text{blow-up of directions}.
\end{gathered}}
\end{equation}
The first two arrows are linear algebra and differential geometry.  The last arrow is
an algebraic-geometric operation.  Keeping them distinct is essential.

\subsubsection{The multinomial covariance and the upper-space identity}
For $m$ categories let
\begin{equation}
 p=(p_1,\ldots,p_m)^\top,\qquad p_i>0,
 \qquad \mathbf 1^\top p=1,
\end{equation}
and define
\begin{equation}
 D_p=\diag(p_1,\ldots,p_m),
 \qquad \Sigma(p)=D_p-pp^\top.
 \label{eq:new-covariance}
\end{equation}
For $u\in\mathbb R^m$,
\begin{equation}
 u^\top\Sigma(p)u
 =\sum_i p_i u_i^2-\left(\sum_i p_i u_i\right)^2
 =\operatorname{Var}_p(u_i)\ge0.
 \label{eq:new-variance}
\end{equation}
Equality holds exactly for constant $u$, hence
\begin{equation}
 \ker\Sigma(p)=\operatorname{span}\{\mathbf1\},
 \qquad \operatorname{rank}\Sigma(p)=m-1.
\end{equation}
The simplex has tangent space
\begin{equation}
 T_p\Delta^{m-1}=\mathbf1^\perp.
\end{equation}

\begin{theorem}[Upper-space quadratic-form identity]
\label{thm:upper-space}
For every $y\in\mathbf1^\perp$,
\begin{equation}
 \boxed{
 y^\top\Sigma(p)^+y=y^\top D_p^{-1}y
 =\sum_{i=1}^m\frac{y_i^2}{p_i}.}
 \label{eq:upper-space}
\end{equation}
Equivalently, the singular multinomial covariance and the diagonal ambient precision
induce exactly the same quadratic form on the constraint tangent space.
\end{theorem}

\begin{proof}
Put $z=D_p^{-1}y$. Since $p^\top D_p^{-1}=\mathbf1^\top$,
\begin{equation}
 \Sigma(p)z=(D_p-pp^\top)D_p^{-1}y
 =y-p(\mathbf1^\top y)=y.
\end{equation}
Thus $z$ is a solution of $\Sigma(p)z=y$. The minimum-norm solution is
$\Sigma(p)^+y$, so $z-\Sigma(p)^+y=c\mathbf1$ for some $c$. Multiplying by $y^\top$
and using $y^\top\mathbf1=0$ gives
\begin{equation}
 y^\top D_p^{-1}y=y^\top\Sigma(p)^+y.
\end{equation}
\end{proof}

\begin{remark}
The theorem does \emph{not} assert $\Sigma(p)^+=D_p^{-1}$; the two matrices are
necessarily different because $\Sigma(p)^+\mathbf1=0$ whereas $D_p^{-1}$ is invertible.
The equality is precisely an equality of intrinsic quadratic forms on $\mathbf1^\perp$.
This is the rigorous content of the ``upper-space'' viewpoint suggested by Yoshizawa's
local/global Gaussian discussion \cite{Yoshizawa2024}.
\end{remark}

The free-coordinate covariance gives the complementary elementary formula.  If
$p_m=1-\sum_{i<m}p_i$, then
\begin{equation}
 \Sigma'=D'-p'p'^\top,
 \qquad D'=\diag(p_1,\ldots,p_{m-1}),
\end{equation}
and the matrix determinant lemma and Sherman--Morrison formula give
\begin{equation}
 \boxed{\det\Sigma'=\prod_{i=1}^m p_i,\qquad
 (\Sigma')^{-1}=D'^{-1}+p_m^{-1}\mathbf1\mathbf1^\top.}
 \label{eq:free-covariance}
\end{equation}
Thus the same constraint resolution can be seen either through the pseudoinverse in
ambient coordinates or through an ordinary inverse in affine coordinates.

\subsubsection{Matrix multinomial: the same identity survives unchanged}
For $P=(P_{ij})\in\mathbb R_{>0}^{r\times c}$ with $\sum_{ij}P_{ij}=1$, Definition~\ref{def:matrix-multinomial}
gives the matrix multinomial law.  Put
\begin{equation}
 p=\operatorname{vec}P,
 \qquad D_P=\diag(p),
 \qquad \Sigma(P)=D_P-pp^\top.
\end{equation}
The positive doubly stochastic slice is
\begin{equation}
 \mathcal M=\left\{P>0:P\mathbf1=\frac1n\mathbf1,
 P^\top\mathbf1=\frac1n\mathbf1\right\}.
 \label{eq:new-M}
\end{equation}
Its tangent space is
\begin{equation}
 T_P\mathcal M=\{U:U\mathbf1=0,\ U^\top\mathbf1=0\},
 \qquad \dim T_P\mathcal M=(n-1)^2.
 \label{eq:new-tangent}
\end{equation}

\begin{theorem}[Matrix upper-space identity]
\label{thm:matrix-upper-space}
For every $P\in\mathcal M$ and every $U,V\in T_P\mathcal M$,
\begin{equation}
 \boxed{
 \langle\operatorname{vec}U,\Sigma(P)^+\operatorname{vec}V\rangle
 =\langle\operatorname{vec}U,D_P^{-1}\operatorname{vec}V\rangle
 =\sum_{i,j}\frac{U_{ij}V_{ij}}{P_{ij}}.}
 \label{eq:matrix-upper-space}
\end{equation}
In particular, the restriction of $\Sigma(P)$ to $T_P\mathcal M$ is positive definite.
\end{theorem}

\begin{proof}
For $U\in T_P\mathcal M$,
$\mathbf1^\top\operatorname{vec}U=\sum_{ij}U_{ij}=0$.  Theorem~\ref{thm:upper-space}
therefore applies to $u=\operatorname{vec}U$.  Polarization of the resulting quadratic-form
identity gives the bilinear identity in \eqref{eq:matrix-upper-space}.
\end{proof}

\begin{remark}
The column constraint $U^\top\mathbf1=0$ is not needed for the upper-space identity itself;
zero total sum is sufficient.  The row and column constraints are needed to identify the
particular Birkhoff tangent space and its dimension.  Thus the upper-space principle is
strictly more general than the doubly stochastic problem.
\end{remark}

Define
\begin{equation}
 \Phi(P)=\sum_{i,j}P_{ij}\log P_{ij}.
 \label{eq:new-entropy-potential}
\end{equation}
Then
\begin{equation}
 d^2\Phi_P(U,V)=\sum_{i,j}\frac{U_{ij}V_{ij}}{P_{ij}},
\end{equation}
so Theorem~\ref{thm:matrix-upper-space} yields
\begin{equation}
 \boxed{g_P(U,V)=d^2\Phi_P(U,V)
 =\langle\operatorname{vec}U,\Sigma(P)^+\operatorname{vec}V\rangle.}
 \label{eq:fisher-hessian}
\end{equation}
Thus the constrained Fisher metric has simultaneously three descriptions:
Hessian metric, upper-space diagonal precision, and covariance pseudoinverse.

\subsubsection{Levi--Civita connection and curvature}
Choose affine coordinates $x^1,\ldots,x^d$, $d=(n-1)^2$, on $\mathcal M$ and write
\begin{equation}
 P(x)=P_0+\sum_{a=1}^d x^a B_a,
 \qquad B_a\in T_P\mathcal M.
\end{equation}
Then
\begin{equation}
 g_{ab}(x)=\sum_{i,j}\frac{(B_a)_{ij}(B_b)_{ij}}{P_{ij}(x)}.
 \label{eq:metric-coordinates}
\end{equation}

\begin{theorem}[Levi--Civita connection of the constrained Fisher metric]
\label{thm:LC}
In the above affine coordinates,
\begin{equation}
 \boxed{
 \Gamma^a_{bc}
 =-\frac12g^{ad}
 \sum_{i,j}\frac{(B_b)_{ij}(B_c)_{ij}(B_d)_{ij}}{P_{ij}^2}.}
 \label{eq:LC-new}
\end{equation}
Moreover the Riemann curvature tensor is
\begin{equation}
 \boxed{
 R_{abcd}=\frac14g^{pq}
 \left(\Phi_{bcp}\Phi_{adq}-\Phi_{acp}\Phi_{bdq}\right),}
 \label{eq:hessian-curvature}
\end{equation}
where
\begin{equation}
 \Phi_{abc}=-\sum_{i,j}\frac{(B_a)_{ij}(B_b)_{ij}(B_c)_{ij}}{P_{ij}^2}.
\end{equation}
\end{theorem}

\begin{proof}
Since $g_{ab}=\partial_a\partial_b\Phi$, we have
$\partial_cg_{ab}=\Phi_{abc}$.  The Levi--Civita formula and symmetry of third
partials give $\Gamma^a_{bc}=\frac12g^{ad}\Phi_{bcd}$, which is exactly
\eqref{eq:LC-new}.  The curvature identity follows by substituting this Hessian form
of the connection into the definition of the Riemann tensor and using cancellation of
fourth derivatives.
\end{proof}

\begin{remark}
The formula is elementary but useful: all connection and curvature coefficients are finite
rational expressions in the entries of $P$.  In particular, the apparent singularity of the
ambient covariance is not an obstacle to differential geometry on the constrained manifold.
The true boundary singularities occur when some $P_{ij}\to0$.
\end{remark}

\subsubsection{The uniform point and the Kronecker-square metric}
Let
\begin{equation}
 P_* = n^{-2}\mathbf1\mathbf1^\top,
 \qquad m=n-1,
 \qquad J_m=\mathbf1_m\mathbf1_m^\top,
 \qquad K_n=I_m-\frac1nJ_m.
\end{equation}
Use the upper-left $(n-1)\times(n-1)$ block as free coordinates.

\begin{theorem}[Uniform-point Fisher geometry]
\label{thm:uniform-fisher}
At $P_*$,
\begin{equation}
 \boxed{
 G_* = n^2(I_m+J_m)\otimes(I_m+J_m),
 \qquad
 G_*^{-1}=n^{-2}K_n\otimes K_n.}
 \label{eq:uniform-fisher}
\end{equation}
Consequently, the Legendre-dual potential has the local expansion
\begin{equation}
 \Psi_{\mathcal M}(y)=\Psi_{\mathcal M}(0)
 +\frac1{2n^2}\operatorname{vec}(y)^\top(K_n\otimes K_n)\operatorname{vec}(y)
 +O(\|y\|^3).
 \label{eq:dual-local}
\end{equation}
\end{theorem}

\begin{proof}
A variation of a free cell changes exactly four cells with signs $+,-,-,+$.  At $P_*$
the ambient Hessian of $\Phi$ is $n^2I$.  The Gram matrix of these four-cell variation
vectors is therefore $n^2(I+J)\otimes(I+J)$.  Since
$(I_m+J_m)^{-1}=I_m-(1/n)J_m=K_n$, the inverse formula follows.  Legendre duality
inverts the Hessian at the dual base point, giving \eqref{eq:dual-local}.
\end{proof}

\begin{table}[ht]
\centering
\small
\begin{tabular}{@{}lll@{}}
\toprule
Level & Matrix & Geometric meaning\\
\midrule
Ambient covariance & $\Sigma(P)=D_P-pp^\top$ & one normal null direction\\
Upper-space precision & $D_P^{-1}$ & non-degenerate ambient representative\\
Constrained metric & $g_P(U,V)=\sum U_{ij}V_{ij}/P_{ij}$ & intrinsic Fisher metric\\
Uniform inverse metric & $n^{-2}K_n\otimes K_n$ & row/column centering\\
\bottomrule
\end{tabular}
\caption{Four equivalent levels of the constraint-resolved Fisher geometry.}
\label{tab:four-levels}
\end{table}

\subsubsection{The Segre variety and the elementary meaning of blow-up}
The positive independence model is
\begin{equation}
 P_{ij}=r_ic_j,
\end{equation}
whose projective closure is the Segre variety
\begin{equation}
 S=\operatorname{Seg}(\mathbb{CP}^{n-1}\times\mathbb{CP}^{n-1})
 \subset\mathbb{CP}^{n^2-1}.
\end{equation}
Equivalently, all $2\times2$ minors vanish.  Its intersection with $\mathcal M$ is
\begin{equation}
 S\cap\mathcal M=\{P_*\}.
\end{equation}

Before using the word ``blow-up'', consider the plane.  The blow-up of the origin has charts
\begin{equation}
 y=ux,\qquad x=vy.
\end{equation}
In the first chart the exceptional divisor is $x=0$, and $u=y/x$ records the limiting
slope.  Thus the elementary slogan is
\begin{equation}
 \boxed{\text{blow up a point} = \text{replace it by its projective space of directions}.}
\end{equation}

\begin{proposition}[Elementary direction separation]
\label{prop:direction-separation-new}
Let $C$ be a smooth plane curve through the origin with
$y=ax+O(x^2)$.  In the chart $y=ux$ its strict transform is $u=a+O(x)$ and meets the
exceptional divisor at $[1:a]$.  Hence two smooth curves have the same point on the
exceptional divisor exactly when their tangent lines agree.
\end{proposition}
\begin{proof}
Substitution gives $ux=ax+O(x^2)$.  Removing the exceptional factor $x$ gives
$u=a+O(x)$, and setting $x=0$ gives the asserted point.
\end{proof}

\begin{figure}[ht]
\centering
\begin{tikzpicture}[>=Latex,scale=0.88]
\begin{scope}[xshift=-4.4cm]
\draw[->] (-1.5,0)--(1.6,0) node[right] {$x$};
\draw[->] (0,-1.2)--(0,1.35) node[above] {$y$};
\draw[thick] (-1.1,-0.55)--(1.1,0.55);
\draw[thick] (-1.1,0.70)--(1.1,-0.70);
\fill (0,0) circle (1.7pt);
\node[align=center] at (0,-1.65) {two curves meet at one point\\but have different slopes};
\end{scope}
\draw[->,very thick] (-2.0,0.05)--(-1.0,0.05) node[midway,above=2mm] {blow-up};
\begin{scope}[xshift=2.1cm]
\draw[->] (-1.45,0)--(1.5,0) node[right] {$x$};
\draw[->] (0,-1.2)--(0,1.35) node[above] {$u=y/x$};
\draw[thick] (-1.0,-0.5)--(1.0,0.5);
\draw[thick] (-1.0,0.72)--(1.0,-0.72);
\draw[thick] (0,-1.0)--(0,1.05);
\fill (0,0.5) circle (1.7pt);
\fill (0,-0.72) circle (1.7pt);
\node[align=center] at (0,-1.65) {exceptional divisor $x=0$\\records the two directions};
\end{scope}
\end{tikzpicture}
\caption{The elementary meaning of blow-up: a collapsed point is replaced by its limiting directions.}
\label{fig:blowup-elementary-new}
\end{figure}
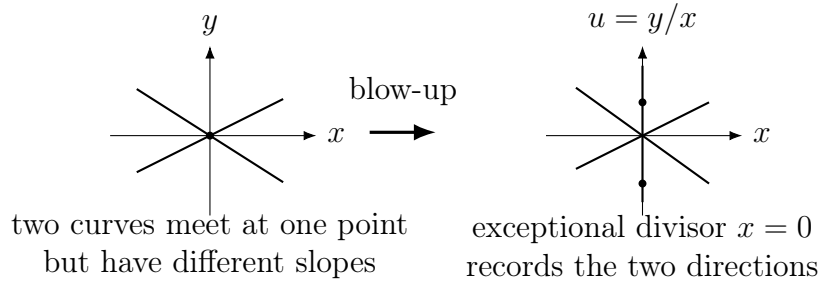

\subsubsection{Blow-up and the Fisher metric}
Let $X$ be a smooth real or complex ambient manifold containing $\mathcal M$ near $P_*$,
and let
\begin{equation}
 \pi:\widetilde X=\operatorname{Bl}_{P_*}X\longrightarrow X
\end{equation}
be the blow-up.  Write local tangent coordinates as
\begin{equation}
 U=r\omega,
 \qquad r\ge0,
 \qquad [\omega]\in\mathbb P(T_{P_*}X).
\end{equation}
The exceptional divisor is $E=\{r=0\}\simeq\mathbb P(T_{P_*}X)$.

\begin{theorem}[Rescaled Fisher metric on the exceptional divisor]
\label{thm:rescaled-fisher}
Let $g$ be the Fisher metric on $\mathcal M$, smoothly extended to a neighborhood of
$P_*$.  For the blow-down map $\pi$, the ordinary pull-back metric satisfies
\begin{equation}
 \boxed{\pi^*g|_E=0.}
 \label{eq:pullback-zero}
\end{equation}
Nevertheless, if $U=r\omega$ and $g_*=g_{P_*}$, then
\begin{equation}
 \pi^*g
 =g_*(\omega,\omega)\,dr^2
 +2r\,g_*(\omega,d\omega)\,dr
 +r^2g_*(d\omega,d\omega)+O(r^3),
 \label{eq:blowup-metric-expansion}
\end{equation}
and therefore the angular part of the rescaled metric has the limit
\begin{equation}
 \boxed{
 \lim_{r\downarrow0}r^{-2}\pi^*g\big|_{\mathrm{angular}}
 =g_*\big|_{\mathrm{angular}}.}
 \label{eq:rescaled-angular}
\end{equation}
After restricting to the $g_*$-unit sphere in $T_{P_*}\mathcal M$ and identifying
$\omega\sim-\omega$, this gives the induced projectivized tangent Fisher metric on
$\mathbb{RP}^{d-1}$, $d=(n-1)^2$.
\end{theorem}

\begin{proof}
At $r=0$, the blow-down map is constant on $E$, so $d\pi$ vanishes on tangent vectors
along $E$, proving \eqref{eq:pullback-zero}.  Since $g$ is smooth,
$g_{P_*+r\omega}=g_*+O(r)$, while $dU=\omega\,dr+r\,d\omega$.  Substitution gives
\eqref{eq:blowup-metric-expansion}; dividing the angular component by $r^2$ and taking
$r\downarrow0$ gives \eqref{eq:rescaled-angular}.
\end{proof}

\begin{remark}
This theorem corrects a tempting but inaccurate slogan.  Blow-up does not by itself turn
the Fisher metric into a non-degenerate metric on the exceptional divisor: the ordinary
pull-back actually vanishes there.  The natural tangent geometry appears after radial
renormalization.  Thus ``constraint resolution'' and ``metric blow-up resolution'' are
related but mathematically distinct operations.
\end{remark}

\subsubsection{Segre--Birkhoff tangent separation}
Differentiating the Segre parameterization at $P_*$ gives
\begin{equation}
 T_{P_*}S=\{a\mathbf1^\top+\mathbf1b^\top:a,b\in\mathbb R^n\},
 \label{eq:new-segre-tangent}
\end{equation}
whereas
\begin{equation}
 T_{P_*}\mathcal M=\{U:U\mathbf1=0,\ U^\top\mathbf1=0\}.
\end{equation}

\begin{theorem}[Segre--Birkhoff tangent separation]
\label{thm:new-segre-separation}
\begin{equation}
 \boxed{T_{P_*}S\cap T_{P_*}\mathcal M=\{0\}.}
 \label{eq:new-separation}
\end{equation}
Moreover,
\begin{equation}
 T_{P_*}S\oplus T_{P_*}\mathcal M
 =T_{P_*}\Delta^{n^2-1},
\end{equation}
where the dimensions are $2n-2$ and $(n-1)^2$ respectively.
\end{theorem}

\begin{proof}
If $U=a\mathbf1^\top+\mathbf1b^\top$ and $U\mathbf1=0$, then
$na+(\mathbf1^\top b)\mathbf1=0$, so $a$ is constant.  Similarly
$U^\top\mathbf1=0$ forces $b$ to be constant.  Hence $U$ is constant, and its row sums
force $U=0$.  The dimension sum is
$(2n-2)+(n-1)^2=n^2-1$, the dimension of the simplex tangent space, so the direct sum follows.
\end{proof}

\begin{proposition}[Blow-up separates tangent directions]
\label{prop:new-blowup-tangent}
Let $M\subset X$ be a smooth submanifold through $z$ and $\widetilde M$ its strict transform
under $\operatorname{Bl}_zX$.  Then
\begin{equation}
 E\simeq\mathbb P(T_zX),
 \qquad
 \widetilde M\cap E=\mathbb P(T_zM).
 \label{eq:new-strict-transform}
\end{equation}
\end{proposition}
\begin{proof}
In a blow-up chart write $x_k=u_kx_1$.  The exceptional divisor is $x_1=0$ with projective
coordinates $[1:u_2:\cdots:u_d]$.  The leading equations of the strict transform are precisely
the linear tangent equations of $M$ at $z$, hence the intersection is $\mathbb P(T_zM)$.
\end{proof}

\begin{theorem}[Local constraint--blow-up resolution at the uniform point]
\label{thm:local-resolution-new}
For the blow-up of a smooth ambient variety $X$ at $P_*$,
\begin{equation}
 E\simeq\mathbb P(T_{P_*}X),
 \quad
 \widetilde S\cap E=\mathbb P(T_{P_*}S),
 \quad
 \widetilde{\mathcal M}\cap E=\mathbb P(T_{P_*}\mathcal M),
\end{equation}
and the last two projective sets are disjoint.
\end{theorem}
\begin{proof}
The first three statements follow from Proposition~\ref{prop:new-blowup-tangent}.  If the two
projective tangent sets had a common point, a nonzero common tangent vector would exist,
contradicting Theorem~\ref{thm:new-segre-separation}.
\end{proof}

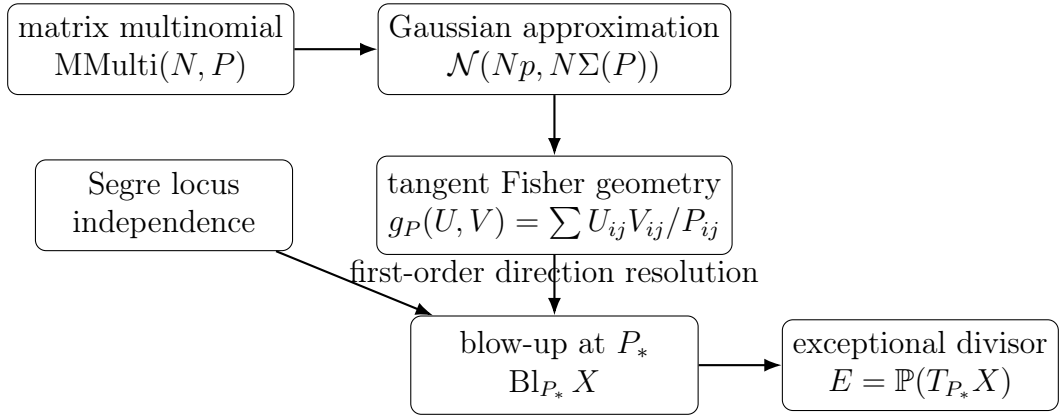
\begin{figure}[ht]
\centering
\begin{tikzpicture}[>=Latex,node distance=8mm and 11mm,
box/.style={draw,rounded corners,align=center,minimum width=34mm,minimum height=12mm},
emph/.style={draw,rounded corners,align=center,minimum width=38mm,minimum height=13mm}]
\node[box] (mm) {matrix multinomial\\$\mathrm{MMulti}(N,P)$};
\node[box,right=of mm] (ga) {Gaussian approximation\\$\mathcal N(Np,N\Sigma(P))$};
\node[emph,below=of ga] (tm) {tangent Fisher geometry\\$g_P(U,V)=\sum U_{ij}V_{ij}/P_{ij}$};
\node[box,left=of tm] (sg) {Segre locus\\independence};
\node[emph,below=of tm] (bl) {blow-up at $P_*$\\$\operatorname{Bl}_{P_*}X$};
\node[box,right=of bl] (ex) {exceptional divisor\\$E=\mathbb P(T_{P_*}X)$};
\draw[->,thick] (mm)--(ga);
\draw[->,thick] (ga)--(tm);
\draw[->,thick] (sg)--(bl);
\draw[->,thick] (tm)--(bl);
\draw[->,thick] (bl)--(ex);
\node[align=center,above=3mm of bl] {first-order direction resolution};
\end{tikzpicture}
\caption{The complete dictionary: statistical approximation, constraint-resolved Fisher geometry,
and blow-up of the distinguished intersection point.}
\label{fig:complete-dictionary-new}
\end{figure}

\subsubsection{A concise synthesis}
The results of this section can be summarized as follows:
\begin{equation}
\boxed{
\begin{array}{c}
\text{multinomial covariance }D_p-pp^\top\\[1mm]
\downarrow\\
\text{upper-space precision }D_p^{-1}\text{ on }\mathbf1^\perp\\[1mm]
\downarrow\\
\text{matrix Fisher metric }g_P(U,V)=\sum U_{ij}V_{ij}/P_{ij}\\[1mm]
\downarrow\\
\text{uniform-point Kronecker geometry }K_n\otimes K_n\\[1mm]
\downarrow\\
\text{blow-up }P_*\rightsquigarrow\mathbb P(T_{P_*}X)\\[1mm]
\downarrow\\
\text{Segre and Birkhoff directions separate.}
\end{array}}
\label{eq:master-dictionary}
\end{equation}
The affine restriction is therefore the mechanism that removes the covariance null direction,
whereas blow-up records the directions that remain indistinguishable at the single point $P_*$. 
The two mechanisms cooperate, but neither should be identified with the other.

\subsection{Summary}

\begin{enumerate}[label=(\roman*)]
\item The entropy gradient flow on $\Birk(n)$ (Eq.~\eqref{eq:constrained-flow}) admits
a Lax representation formally identical to Nakamura's, but this representation is
spectrally vacuous (Prop.~\ref{prop:general-lax} and the ensuing remark), because it
follows from a general identity valid for \emph{any} choice of forcing term, not a
special feature of the doubly-stochastic constraint.
\item Nevertheless, the flow admits an exact linearization
(Thm.~\ref{thm:linearization-birkhoff}) in terms of $2\times2$ log-odds ratios, giving
$(n-1)^2-1$ independent first integrals (Cor.~\ref{cor:first-integrals}) and an
explicit closed-form solution via Sinkhorn scaling (Cor.~\ref{cor:explicit-solution}),
refining the classical margin-invariance theorem of Fienberg \cite{Fienberg1970}.
This exact linearization is a close relative of Tanabe's exact first integral
$g(x(t))=e^{-t}g(x^0)$ for Branin's continuous Newton--Raphson method
\cite{Tanabe1980}, and our constrained flow itself (Eq.~\eqref{eq:constrained-flow})
is a direct instance of Tanabe's continuous gradient-projection method
(Rem.~\ref{rem:tanabe-gp}, \ref{rem:tanabe-first-integral}).
\item These first integrals furnish genuine Hamiltonian canonical coordinates
(Thm.~\ref{thm:hamiltonian}), reproducing Nakamura's Liouville--Arnol'd integrability
in a structurally more flexible form.
\item The submanifold $\mathcal M$ is itself dually flat with the log-odds-ratio matrix
as its natural ($e$-affine) dual coordinate (Prop.~\ref{prop:dual-coord}), generalizing
Nakamura's $\eta$--$\theta$ duality; however, the associated K\"ahler-type potential
fails to admit a closed algebraic form. This failure is explained precisely via the
Segre embedding (Prop.~\ref{prop:transversality}) and the Guillemin/Abreu theory of
toric K\"ahler potentials: closed forms survive exactly on the independence model and
its weight-twisted deformations (Segre-type leaves), but not on the transversal,
margin-fixing slice $\mathcal M$ itself.
\item The general symbolic $LDM^\top$ calculus for $D+uv^\top$ (possibly
non-symmetric, possibly singular) was first established by Tanabe and Sagae
\cite{TanabeSagaeNumAlg1992} (Rem.~\ref{rem:priority}); the complementary real
symmetric-eigenvalue theory of Steerneman and van Perlo-ten Kleij
\cite{Steerneman2005} both explains the elementary algebraic origin of the
square-root embedding used to obtain (iv), and -- via the Sherman--Morrison identity
applied to the Kronecker-factorized Hessian \eqref{eq:GM-kronecker} -- yields the
exact quadratic approximation \eqref{eq:local-Psi} of the otherwise-missing K\"ahler
potential in a neighborhood of the flow's equilibrium (\S\ref{sec:steerneman}).
\item Independently, Tanabe and Sagae's exact Moore--Penrose formula
\cite[Prop.~1]{TanabeSagae1992} for the ordinary multinomial covariance,
$(P-pp^\top)^+=HP^{-1}H$, is the one-factor antecedent of the Kronecker-square
structure \eqref{eq:GM-inverse}, and their symbolic-Cholesky determinant identity
(their Cor.~2) supplies the exact algebraic input,
$\det\Sigma'=p_1p_2\cdots p_n$, from which the leading-order asymptotic entropy of
the multinomial (and hence the matrix multinomial) is rigorously derived
(\S\ref{sec:tanabe-mp}, \S\ref{sec:entropy}).
\item The mean-parameter potential of the matrix multinomial is exactly the negative
Shannon entropy, $\varphi_{\mathrm{mult}}(P)=-H(P)$
(Eq.~\eqref{eq:phi-is-shannon}), while both potentials of the matrix Gaussian
location family, \eqref{eq:A-gauss-explicit}--\eqref{eq:phi-gauss-explicit}, are pure
quadratic forms unrelated in value to the (mean-independent) Gaussian entropy
(Rem.~\ref{rem:gauss-not-entropy}) -- the potential-theoretic counterpart of the
flat-versus-curved dichotomy of \S\ref{sec:kahler}. On $\mathcal M$, the entropy
potential coincides, up to an additive constant, with the row/column mutual
information (Prop.~\ref{prop:mutual-info-identity}), so the entire gradient flow of
\S\ref{sec:steerneman}--\S\ref{sec:kahler} is mutual-information descent to the
independence locus; and $H(P)$ is exactly the von Neumann entropy of the diagonal
density matrix $\diag(\operatorname{vec}P)$ (Prop.~\ref{prop:shannon-is-vN}) -- an entrywise,
not spectral, identification, since $-\operatorname{tr}(P\log P)$ formed from $P$ as an
operator is a genuinely different (and generally non-real) quantity
(Rem.~\ref{rem:not-trace-entropy}), with the singular-value entropy
\eqref{eq:singular-value-entropy} as a third, distinct alternative
(Rem.~\ref{rem:singular-value-entropy}) -- pointing to a
natural non-commutative (quantum Sinkhorn) generalization of the whole paper
(\S\ref{sec:von-neumann}).
\item The discrete entropy of the matrix multinomial distribution converges, as
$N\to\infty$, to the differential entropy of the (generically non-separable)
Khatri--Mitra Gaussian approximation at rate $O(N^{-1})$
(Prop.~\ref{prop:asymptotic-entropy}), but this convergence is \emph{not} uniform
over the parameter simplex: the failure of uniformity is governed by exactly the
same quantity, $\sum_{ij}P_{ij}^{-1}$, that controls Tanabe and Sagae's own
condition-number bound for $P-pp^\top$ (Prop.~\ref{prop:nonuniform}). A genuinely
\emph{separable} (Kronecker-covariance, matrix-normal) Gaussian limit exists only on
the independence (Segre-variety) locus of \S\ref{sec:kahler}, under an additional
sparse double-scaling regime, and even there its entropy differs from the general
formula by an exact, diverging offset $-\tfrac12\log(2\pi eN)$
(Prop.~\ref{prop:offset}) -- the entropic signature of the same codimension-one,
$m$-flat constraint that organizes the K\"ahler-duality obstruction of
\S\ref{sec:kahler}.
\item The rank deficiency of the ambient multinomial covariance has a clean
constraint-geometric interpretation: the null direction is a constraint-normal
direction, while the covariance restricted to the doubly-stochastic tangent space is
non-degenerate. The exact Moore--Penrose formula and the equilibrium identity
$G_{\mathcal M}^{-1}=n^{-2}K_n\otimes K_n$ are the linear-algebraic signatures of this
resolution. The Segre independence locus intersects the positive doubly-stochastic
slice at the uniform point, providing a natural candidate center for a future algebraic
blow-up. We emphasize that the present paper establishes the constraint resolution,
not yet the full Rees-algebra blow-up or an extension theorem for the Fisher tensor on
its exceptional divisor (\S\ref{sec:blowup}).
\end{enumerate}

The algebraic identities and the previously reported flow, Hessian, and entropy calculations
were checked independently by direct numerical integration, finite-difference computation,
and exact combinatorial enumeration (Python/NumPy/SciPy, $n=3,4,5$).  The new local
blow-up statements are proved analytically from the standard blow-up charts and do not
require a numerical claim.


\section{Density, Duality, and Blow-Up on Elliptic Curves}
\label{sec:elliptic-duality}

The blow-up construction of \S\ref{sec:birkhoff} resolved a degenerate Fisher metric
at the point where the independence locus meets the doubly-stochastic slice, producing
an exceptional divisor carrying its own information geometry; and the same section's
discussion of Tanabe's exact first integral for Branin's continuous Newton--Raphson
method (Remark~\ref{rem:tanabe-first-integral}) showed that a birationally-invariant
change of time turns a naive gradient flow into one with an exact exponential decay
law. This section develops both phenomena --- blow-up-resolved information geometry,
and Tanabe--Branin exponential decay --- in a third, independent setting: the classical
birational geometry of elliptic curves. An elliptic curve admits many equivalent plane
models (Weierstrass cubics, Jacobi quartics, Mordell quartics) related by birational
transformations of the ambient affine plane; we show that the \emph{defining
polynomial} of each model transforms as a relative invariant (density) of weight one
under the birational map, and that this single fact controls two phenomena that at
first appear unrelated. Blowing up two of these curves at the vanishing locus of the
density weight produces exceptional divisors carrying a canonical one-dimensional
Kullback--Leibler-type information geometry, with Fisher metric identically $dr^2/r^2$
on all three divisors, mutually isomorphic via explicit affine maps
(Theorem~\ref{thm:main}); on the Mordell side, the same density weight governs exactly
how a naive gradient (Branin) flow toward the quartic fails to be birationally
invariant, and how a logarithmic time reparametrization repairs it, producing a flow
whose decay law $g(\tau)=g(0)e^{-\tau}$ is formally the same exponential law that
underlies the Kullback--Leibler geometry of the blow-up (\S\ref{sec:synthesis}). Along
the way we give an elementary, matrix-theoretic account of how the group law of an
elliptic curve is realized by conjugation and translation of a $3\times3$ symmetric
matrix, and we illustrate every construction on the classical taxicab curve associated
with $N=1729=12^3+1^3=10^3+9^3$. As in \S\ref{sec:birkhoff}, every polynomial identity
and numerical claim below has been verified by computer algebra (pseudocode in
\S\ref{app:sympy}).

\subsection{Introduction}\label{dd-sec:intro}

An elliptic curve is, up to isomorphism, a single geometric object, but
it admits many different plane models: a cubic in (long or short)
Weierstrass form, a quartic in Jacobi form, a quartic in Mordell form,
and so on. Passing between these models is a birational---not
regular---transformation of the ambient affine plane: it is a rational
map, undefined along certain curves and points, whose restriction to the
elliptic curve itself is an isomorphism. The classical theory
(\cite{Cassels,Silverman,Mordell1969}) tells us precisely what such a
transformation preserves \emph{on the curve}: the curve is carried to
the curve, and the canonical regular differential $\omega=dx/y$ is
carried to the canonical regular differential of the target model. What
such a transformation does to the \emph{ambient plane}---to the ideal
generated by the defining polynomial, and to naive plane geometry such
as Euclidean gradients, Hessians, or gradient flows built from the
defining polynomial---is not usually discussed, because in most
applications only the curve itself matters.

This paper studies exactly that ``off-curve'' behaviour, for two
classical birational bridges, and shows that it is governed by a single
principle:

\begin{quote}
\emph{Under a birational transformation between two plane models of an
elliptic curve, the defining polynomial does not pull back to the
defining polynomial of the target model; it pulls back to that
polynomial multiplied by an explicit rational function---a density of
weight one. The zero locus of that density is exactly the
indeterminacy/exceptional locus of the transformation, and resolving
that locus by a blow-up, or correcting for it by a change of time in an
associated gradient flow, recovers exact invariance.}
\end{quote}

We call this the \emph{density principle}. It is a two-line observation
once stated, yet it organizes a surprising amount of structure. We
illustrate it with two case studies.

\smallskip
\noindent\textbf{Case study I: Connell's theta transformation.}
Let $\EW$ be a Weierstrass cubic and $E_J$ its birational image as a
Jacobi quartic (\S\ref{dd-sec:setup}). Connell's theta transformation
$T\colon E_J\dashrightarrow\EW$ (\cite{Connell}) satisfies the
factorization identity (Theorem~\ref{thm:factor})
\begin{equation}\label{eq:intro-factor}
  F_W(T(u,v)) \;=\; \frac{4q^2}{u^6}\cdot F_2(u,v)\cdot F_4(u,v),
\end{equation}
where $F_4=0$ defines $E_J$ and $F_2=A\cdot u^2$ for a linear polynomial
$A$ (the \emph{Connell datum}). The zero locus $\{A=0\}\cup\{u=0\}$ is
precisely where $T^{-1}$ is singular. We show
(\S\ref{dd-sec:blowup}--\S\ref{sec:main}) that blowing up $\EW$ at $\{A=0\}$
and $E_J$ at $\{u=0\}$---and, for good measure, $\EW$ at each branch
point $\{y=0\}$---produces three exceptional divisors, each canonically
isomorphic to $\PP^1$, each carrying a logarithmic potential
$\Phi(\xi)=\log(\alpha\xi-\beta)$ whose Bregman divergence is isomorphic
to the Kullback--Leibler divergence between exponential distributions,
and whose Fisher information metric is \emph{identically} $dr^2/r^2$ in
the natural rate coordinate $r=\alpha\xi-\beta$. The three statistical
manifolds so obtained are mutually isomorphic, and the isomorphism is
realized by the extension of $T$ to the blow-ups (Theorem~\ref{thm:main}).

\smallskip
\noindent\textbf{Case study II: Mordell's quartic--cubic transformation.}
Let $\Es\colon Y^2=4X^3-g_2X-g_3$ be a short Weierstrass cubic and let
$\EM\colon y^2=x^4-6cx^2+4dx+e$ be a Mordell quartic obtained by choosing
a point $(c,d)\in\Es$. Mordell's classical transformation
$\Phi\colon(X,Y)\mapsto(x,y)$ (\cite{Mordell1969}) satisfies
(Theorem~\ref{thm:density})
\begin{equation}\label{eq:intro-density}
  g_4(\Phi(X,Y)) \;=\; \det D\Phi(X,Y)\cdot g_3(X,Y), \qquad
  \det D\Phi = \frac{1}{c-X},
\end{equation}
where $g_4$ and $g_3$ are the defining polynomials of $\EM$ and $\Es$.
This is the same phenomenon as \eqref{eq:intro-factor}: a defining
polynomial pulls back to the other defining polynomial times a density.
We use \eqref{eq:intro-density} for three purposes: to give an
elementary matrix-theoretic proof that the group law of $\Es$ is realized
by conjugation and translation of a $3\times3$ symmetric matrix
(\S\ref{sec:group}); to reconfirm, by direct computation, that the
canonical differential $dx/y=dX/Y$ is exactly invariant
(\S\ref{sec:diff}); and, most interestingly, to show
(\S\ref{dd-sec:flow}) that a naive gradient flow toward $\{g_4=0\}$ does
\emph{not} push forward to the naive gradient flow toward $\{g_3=0\}$,
but that reparametrizing time by
$\tau(t)=t-\log|(c-X(t))/(c-X(0))|$ repairs the discrepancy exactly,
producing a flow that decays as $g_3(\tau)=g_3(0)e^{-\tau}$.

\smallskip
\noindent\textbf{Why put these together.} Both case studies exhibit the
density principle, but they exploit it in dual ways: Case~I resolves the
density's zero locus by \emph{blowing up}, extracting a static,
one-dimensional information-geometric structure on the exceptional
divisor; Case~II absorbs the density into a \emph{time
reparametrization} of a dynamical system, extracting a birationally
invariant flow toward the curve. In \S\ref{sec:synthesis} we point out
that these two constructions produce, formally, the very same
exponential law: the rate parameter $r$ on the blow-up divisor
parametrizes the exponential family $\{r\,e^{-rt}:r>0\}$ that underlies
the Kullback--Leibler geometry of Case~I, while the corrected flow of
Case~II decays according to $g(\tau)=g(0)e^{-\tau}$. We regard this
coincidence as suggestive evidence that the density principle is the
common source of both an information-geometric and a dynamical
manifestation of birational invariance, and we state it as an open
direction for further work rather than a theorem.

Throughout, we illustrate every construction on the elementary example
of the taxicab curve associated with the Hardy--Ramanujan number
$N=1729=12^3+1^3=10^3+9^3$ (\S\ref{sec:taxi}), computing every quantity
explicitly, so that a reader with only a first course in elliptic curves
can follow the entire argument numerically.

\bigskip
\noindent\textbf{Organization.}
Section~\ref{dd-sec:setup} fixes the two families of plane models and
states the two classical birational transformations $T$ and $\Phi$.
Section~\ref{sec:density} proves the factorization identity for $T$ and
the Jacobian identity for $\Phi$, and states the density principle
precisely. Section~\ref{dd-sec:blowup} carries out the scheme-theoretic
blow-ups used in Case~I. Section~\ref{sec:dual} constructs the
Kullback--Leibler-type dual geometry on each exceptional divisor and
computes the Fisher metric. Section~\ref{sec:main} proves the
isomorphism theorem linking the three geometries via the extension of
$T$. Section~\ref{sec:group} gives the elementary matrix realization of
the group law of $\Es$ used in Case~II, and Section~\ref{sec:diff} proves
invariance of the canonical differential. Section~\ref{dd-sec:flow}
constructs the birationally invariant, time-reparametrized Branin flow
and the associated Pythagorean-type identity. Section~\ref{sec:synthesis}
draws the connection between the static and dynamic pictures.
Section~\ref{sec:taxi} works out $N=1729$ in full numerical detail.
Section~\ref{dd-sec:discussion} discusses the scope of our results and open
problems. Appendix~\ref{app:sympy} gives pseudocode for the
computer-algebra verifications underlying the identities of the paper.

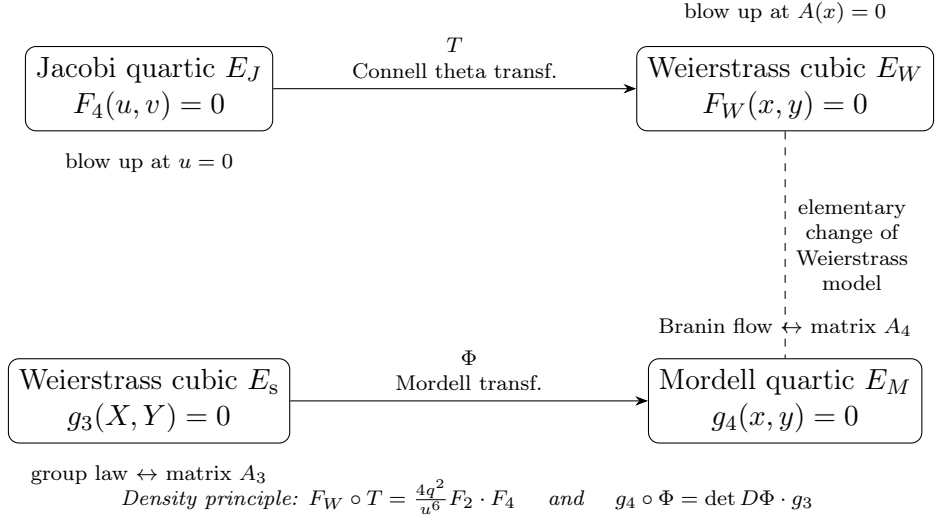
\begin{figure}[htbp]
\centering
\begin{tikzpicture}[
  node distance=2.5cm,
  box/.style={draw, rounded corners, minimum width=3.2cm, minimum height=1cm, align=center, font=\small},
  lbl/.style={font=\scriptsize, midway, align=center},
  >=Stealth
]
\node[box] (EJ) {Jacobi quartic $E_J$\\ $F_4(u,v)=0$};
\node[box, right=4.8cm of EJ] (EW) {Weierstrass cubic $\EW$\\ $F_W(x,y)=0$};
\node[box, below=3.0cm of EJ] (Es) {Weierstrass cubic $\Es$\\ $g_3(X,Y)=0$};
\node[box, below=3.0cm of EW] (EM) {Mordell quartic $\EM$\\ $g_4(x,y)=0$};

\draw[->] (EJ) -- node[lbl,above] {$T$ \\ Connell theta transf.} (EW);
\draw[->] (Es) -- node[lbl,above] {$\Phi$ \\ Mordell transf.} (EM);
\draw[dashed] (EW) -- node[lbl,right] {elementary\\ change of\\ Weierstrass\\ model} (EM);

\node[below=0.15cm of EJ, font=\scriptsize] {blow up at $u=0$};
\node[above=0.15cm of EW, font=\scriptsize] {blow up at $A(x)=0$};
\node[below=0.15cm of Es, font=\scriptsize] {group law $\leftrightarrow$ matrix $A_3$};
\node[above=0.15cm of EM, font=\scriptsize] {Branin flow $\leftrightarrow$ matrix $A_4$};

\node[align=center, font=\scriptsize\itshape, below=0.9cm of $(Es)!0.5!(EM)$]
 {Density principle: $F_W\circ T = \frac{4q^2}{u^6}F_2\cdot F_4$
   \quad and \quad $g_4\circ\Phi=\det D\Phi\cdot g_3$};
\end{tikzpicture}
\caption{The two case studies of this paper. Both are birational
bridges between a cubic and a quartic model of an elliptic curve, and in
both cases the defining polynomial of one model pulls back to a
\emph{density}, not a function, on the other. Case~I (top row) resolves
this by blow-up, producing an information-geometric structure on the
exceptional divisor (\S\ref{dd-sec:blowup}--\S\ref{sec:main}). Case~II
(bottom row) resolves it by a change of time in a gradient flow
(\S\ref{sec:group}--\S\ref{dd-sec:flow}).}
\label{fig:overview}
\end{figure}

\bigskip
\noindent\textit{Figure~\ref{fig:overview} previews the overall
architecture of the paper.}

\subsection{Two Classical Plane Models and Their Birational Bridges}\label{dd-sec:setup}

We work over an algebraically closed field $k$ of characteristic zero
(e.g.\ $\overline{\Q}$ or $\mathbb{C}$), except in
\S\ref{sec:taxi} where we specialize to $\Q$ and $\R$ for the numerical
example.

\subsubsection{Case~I: the Weierstrass cubic and the Jacobi quartic}
\label{subsec:caseI-models}

Fix parameters $q,a,b,c,d\in k$ with $q\ne0$.

\begin{definition}[Weierstrass cubic]\label{def:EW}
The \emph{Weierstrass cubic} $\EW$ is the projective closure of the
affine curve
\[
  F_W(x,y):=
  y^2+\frac{d}{q}xy+2bqy
  \;-\;
  x^3-\Bigl(c-\frac{d^2}{4q^2}\Bigr)x^2+4aq^2 x-a(d^2-4cq^2)
  \;=\;0.
\]
We assume $\Delta(\EW)\ne0$, so $\EW$ is a smooth projective curve of
genus $1$.
\end{definition}

\begin{definition}[Jacobi quartic]\label{def:EJ}
The \emph{Jacobi quartic} $E_J$ is the projective closure of
\[
  F_4(u,v):=v^2 - au^4-bu^3-cu^2-du-q^2 \;=\; 0,
\]
again assumed smooth ($\Delta(E_J)\ne0$).
\end{definition}

\begin{definition}[Connell datum]\label{def:A}
The \emph{Connell datum} is the linear polynomial
\[
  A(x)=4q^2(x+c)-d^2,
\]
whose unique zero $x_0=d^2/(4q^2)-c$ is the singular base point of $T^{-1}$
below.
\end{definition}

\begin{definition}[Theta transformation]\label{def:T}
Define the rational map $T\colon\A^2_{(u,v)}\dashrightarrow\A^2_{(x,y)}$ by
\begin{align}
  x &= \frac{2q(v+q)+du}{u^2}, \label{eq:x}\\[2pt]
  y &= \frac{8q^3(v+q)+4q^2(cu^2+du)-d^2u^2}{2qu^3}. \label{eq:y}
\end{align}
\end{definition}

\begin{proposition}[Birational equivalence]\label{prop:birat}
$T$ restricts to a birational equivalence $E_J\dashrightarrow\EW$, with
inverse
\begin{align}
  u &= \frac{A(x)}{2qy}, \label{eq:u}\\[2pt]
  v &= \frac{u^2x}{2q}-\frac{du}{2q}-q. \label{eq:v}
\end{align}
$T^{-1}$ is undefined on $\{y=0\}\cup\{A=0\}$; $T$ is undefined on $\{u=0\}$.
\end{proposition}

\begin{proof}
Substituting \eqref{eq:u}--\eqref{eq:v} into \eqref{eq:x}--\eqref{eq:y}
returns the identity on $\EW\setminus\{A=0,y=0\}$; substituting
\eqref{eq:x}--\eqref{eq:y} into $F_4$ and reducing modulo $F_W$ gives
zero on $E_J\setminus\{u=0\}$. Both are verified by rational arithmetic
(Appendix~\ref{app:sympy}).
\end{proof}

The terminology and the explicit formula for $T$ are drawn from the
exposition of Connell \cite{Connell}; the structural role of $A(x)$ was
noted there and is exploited systematically below.

\subsubsection{Case~II: the short Weierstrass cubic and the Mordell quartic}
\label{subsec:caseII-models}

Fix parameters $c,d,e\in k$, and set
\[
  g_2 = e+3c^2, \qquad g_3=-ce-d^2+c^3.
\]

\begin{definition}[Matrix representation]\label{def:matrices}
Let $\xi=(1,x,y)^\top$ and $\eta=(1,X,Y)^\top$, and set
\begin{gather*}
A_4=\begin{pmatrix}-e&-2d&0\\-2d&6c&0\\0&0&1\end{pmatrix},\quad
b_4=\begin{pmatrix}0\\1\\0\end{pmatrix}, \\[4pt]
A_3=\begin{pmatrix}-g_3&-g_2/2&0\\-g_2/2&0&0\\0&0&-1\end{pmatrix},\quad
b_3=\begin{pmatrix}0\\-\sqrt[3]{4}\\0\end{pmatrix}.
\end{gather*}
Define
\[
  g_4(x,y):=\xi^\top A_4\xi-(\xi^\top b_4)^4, \qquad
  g_3(X,Y):=-\eta^\top A_3\eta+(\eta^\top b_3)^3.
\]
\end{definition}

\begin{proposition}\label{prop:expand}
\[
g_4(x,y)=y^2-x^4+6cx^2-4dx-e,\qquad
g_3(X,Y)=Y^2-4X^3+g_2X+g_3.
\]
In particular $\{g_4=0\}$ is the \emph{Mordell quartic}
$\EM\colon y^2=x^4-6cx^2+4dx+e$ and $\{g_3=0\}$ is the \emph{short
Weierstrass cubic} $\Es\colon Y^2=4X^3-g_2X-g_3$.
\end{proposition}

\begin{proof}
Expanding, $\xi^\top A_4\xi=-e+6cx^2+y^2-4dx$ and $(\xi^\top b_4)^4=x^4$,
so $g_4=y^2-x^4+6cx^2-4dx-e$. Likewise
$\eta^\top A_3\eta=-g_3-g_2X-Y^2$ and $(\eta^\top b_3)^3=-4X^3$, so
$g_3=(g_3+g_2X+Y^2)-4X^3$, as claimed. (Note the harmless double use of
the symbol $g_3$: as an entry of $A_3$ it is a parameter, while
$g_3(X,Y)$ denotes the defining function; context always disambiguates.)
\end{proof}

The point of the matrix packaging $(A_4,b_4)$, $(A_3,b_3)$ is that it
makes the elliptic-curve group law of $\Es$ literally an operation on
$A_4$---see Theorems~\ref{thm:inversion} and~\ref{thm:addition} below.

\begin{theorem}[Mordell's transformation]\label{thm:mordell}
The map
\[
\Phi\colon (X,Y)\longmapsto (x,y)=\Bigl(\frac{Y-d}{2(X-c)},\;-x^2+2X+c\Bigr)
\]
satisfies the identity
\[
  g_4\bigl(\Phi(X,Y)\bigr)=\frac{g_3(X,Y)}{c-X}.
\]
In particular $g_3(X,Y)=0\Rightarrow g_4(\Phi(X,Y))=0$, i.e.\ $\Phi$ maps
$\Es$ birationally onto $\EM$.
\end{theorem}

\begin{proof}
With $x=(Y-d)/(2(X-c))$ and $y=-x^2+2X+c$, the definition of $y$ gives
immediately $y=-(x^2-(2X+c))$, so
\[
  x^4-6cx^2+4dx+e = x^4-2x^2(2X+c)+(2X+c)^2 = \bigl(x^2-(2X+c)\bigr)^2 = y^2
\]
holds precisely when $x^2(X-c)+dx=X^2+cX+\tfrac14(c^2-e)$, an identity
verified directly upon substituting $x=(Y-d)/(2(X-c))$. Squaring
$2x(X-c)=Y-d$ and clearing denominators, one finds that
\[
  16(c-X)^4\bigl[y^2-x^4+6cx^2-4dx-e\bigr]
\]
expands, as a polynomial in $X,Y,c,d,e$, to
\[
  16(c-X)^3\bigl[Y^2-4X^3+g_2X+g_3\bigr]
\]
(Appendix~\ref{app:sympy}). Dividing by $16(c-X)^4$ gives
$g_4(x,y)=g_3(X,Y)/(c-X)$, as claimed. This computation is essentially
Theorem~2 of Mordell~\cite[p.~77]{Mordell1969}, specialized to the case
where the leading quartic coefficient has already been normalized to
$1$.
\end{proof}

\begin{corollary}\label{cor:inverse}
Given the cubic $\Es$ (i.e.\ given $g_2,g_3$), every point $(c,d)\in\Es$
with $d^2=4c^3-cg_2-g_3$ furnishes, via $e=g_2-3c^2$, a Mordell quartic
$\EM$ for which Theorem~\ref{thm:mordell} holds.
\end{corollary}

\begin{proof}
Substituting $e=g_2-3c^2$ into $g_3=-ce-d^2+c^3$ eliminates $e$:
$g_3=-c(g_2-3c^2)-d^2+c^3=4c^3-cg_2-d^2$; solving for $d^2$ gives the
stated equation, which is exactly the condition that $(c,d)\in\Es$.
\end{proof}

Thus each choice of base point $P=(c,d)\in\Es$ produces its own Mordell
quartic model $\EM(P)$, with its own matrix $A_4(P)$; the family of all
such matrices is parametrized by the curve $\Es$ itself
(\S\ref{sec:group}).

\begin{figure}[htbp]
\centering
\begin{tikzpicture}[>=Stealth, scale=1]
\begin{scope}
\clip (-3.3,-2.6) rectangle (3.3,2.6);
\draw[thick, blue!70!black, domain=-1.0:3.1, samples=140, smooth]
  plot (\x, {sqrt(max(\x*\x*\x-3*\x+3,0))});
\draw[thick, blue!70!black, domain=-1.0:3.1, samples=140, smooth]
  plot (\x, {-sqrt(max(\x*\x*\x-3*\x+3,0))});
\end{scope}
\node[below] at (0,-2.7) {(a) Weierstrass cubic $\Es\colon Y^2=4X^3-g_2X-g_3$};
\end{tikzpicture}
\hspace{0.5cm}
\begin{tikzpicture}[>=Stealth, scale=1]
\begin{scope}
\clip (-3.3,-2.6) rectangle (3.3,2.6);
\draw[thick, red!70!black, domain=-2.3:2.3, samples=140, smooth]
  plot (\x, {sqrt(max(\x*\x*\x*\x-2.6*\x*\x+1.2*\x+3.5,0))});
\draw[thick, red!70!black, domain=-2.3:2.3, samples=140, smooth]
  plot (\x, {-sqrt(max(\x*\x*\x*\x-2.6*\x*\x+1.2*\x+3.5,0))});
\end{scope}
\node[below] at (0,-2.7) {(b) Mordell quartic $\EM\colon y^2=x^4-6cx^2+4dx+e$};
\end{tikzpicture}
\caption{Real loci of a representative short Weierstrass cubic (a) and
Mordell quartic (b); shapes are schematic. Mordell's map $\Phi$ of
Theorem~\ref{thm:mordell} identifies these two curves birationally,
while distorting the ambient plane by the density weight
$\det D\Phi=1/(c-X)$.}
\label{fig:cubic-quartic-shapes}
\end{figure}
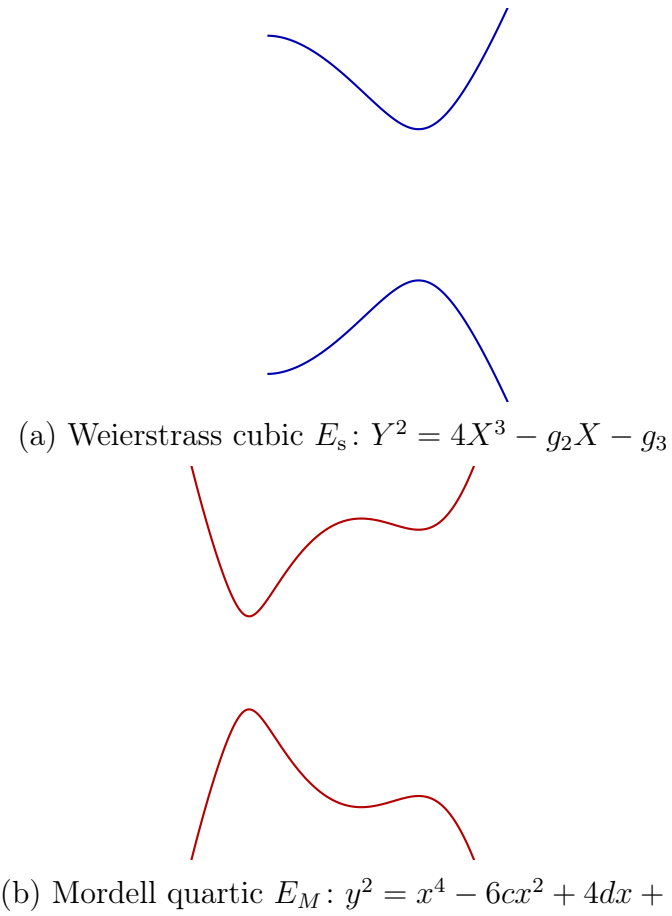

\subsection{The Density Principle}\label{sec:density}

We now state and prove, for both case studies, the identity expressing
``pullback of defining polynomial = density $\times$ defining
polynomial,'' and record the elementary observation that these are the
same statement.

\subsubsection{Factorization identity for $T$}

\begin{definition}[Quadratic factor]\label{def:F2}
\[
  F_2(u,v) \;=\; (4cq^2-d^2)u^2+4dq^2u+8q^3(v+q).
\]
\end{definition}

\begin{lemma}\label{lem:F2=Au2}
Substituting $x=x(u,v)$ from \eqref{eq:x},
$F_2(u,v)=A(x(u,v))\cdot u^2$; hence $\{F_2=0\}=\{A(x(u,v))=0\}$ for $u\ne0$.
\end{lemma}

\begin{proof}
\[
  A(x(u,v))\cdot u^2
  =\Bigl[4q^2\Bigl(\tfrac{2q(v+q)+du}{u^2}+c\Bigr)-d^2\Bigr]u^2
  =8q^3(v+q)+4dq^2u+(4cq^2-d^2)u^2=F_2(u,v). \qedhere
\]
\end{proof}

\begin{theorem}[Factorization identity]\label{thm:factor}
\begin{equation}\label{eq:factor}
  F_W\bigl(T(u,v)\bigr)
  \;=\;\frac{4q^2}{u^6}\cdot F_2(u,v)\cdot F_4(u,v).
\end{equation}
Consequently, $T^{-1}(\EW)=\{F_2=0\}\cup E_J=\{A=0\}\cup\{u=0\}\cup E_J$.
\end{theorem}

\begin{proof}
Multiplying $F_W(T(u,v))$ by $4q^2u^6$ yields a polynomial identity in
$(u,v)$ over $k[a,b,c,d,q]$ that factors as $4q^2\cdot F_2\cdot F_4$
(Appendix~\ref{app:sympy}). The decomposition $\{F_2=0\}=\{A=0\}\cup\{u=0\}$
follows from Lemma~\ref{lem:F2=Au2} together with $F_2|_{u=0}=8q^3(v+q)$.
\end{proof}

\begin{corollary}[Intersection multiplicities]\label{cor:mult}
$\{F_2=0\}\cap\{F_4=0\}$ consists of $u=0$ with multiplicity $3$ (the flex
point at infinity on $\EW$) and one further point
$u^*=-8q^2d(d^2-4cq^2)/(64aq^6-(d^2-4cq^2)^2)$.
\end{corollary}

\begin{proof}
Substituting the solution $v=v(u)$ of $F_2=0$ into $F_4$ produces
$u^3\cdot\ell(u)$ with $\ell$ linear.
\end{proof}

\subsubsection{Jacobian identity for $\Phi$}

\begin{theorem}[Density weight of $\Phi$]\label{thm:density}
\[
\det D\Phi(X,Y)=\frac{1}{c-X}, \qquad
g_4(\Phi(X,Y))=\det D\Phi(X,Y)\cdot g_3(X,Y).
\]
\end{theorem}

\begin{proof}
With $\Phi_1=x=\dfrac{Y-d}{2(X-c)}$ and $\Phi_2=y=-x^2+2X+c$,
\[
\frac{\partial x}{\partial X}=-\frac{Y-d}{2(X-c)^2},\quad
\frac{\partial x}{\partial Y}=\frac{1}{2(X-c)},\quad
\frac{\partial y}{\partial X}=-2x\frac{\partial x}{\partial X}+2,\quad
\frac{\partial y}{\partial Y}=-2x\frac{\partial x}{\partial Y},
\]
so
\[
\det D\Phi
=\frac{\partial x}{\partial X}\frac{\partial y}{\partial Y}
-\frac{\partial x}{\partial Y}\frac{\partial y}{\partial X}
=-2\frac{\partial x}{\partial Y}\cdot\frac{\partial x}{\partial X}
 +2\frac{\partial x}{\partial Y}\cdot\frac{\partial x}{\partial X}
 -2\frac{\partial x}{\partial Y}
=-\frac{1}{X-c}=\frac{1}{c-X}
\]
(the two terms containing $x\,\partial x/\partial X\,\partial x/\partial
Y$ cancel). Comparing with Theorem~\ref{thm:mordell} gives the second
identity.
\end{proof}

\subsubsection{The density principle, stated uniformly}

\begin{proposition}[Density principle]\label{prop:density-principle}
Let $\Psi\colon C_1\dashrightarrow C_2$ be a birational map between two
plane models of an elliptic curve, with defining polynomials $h_1,h_2$.
In both of our case studies there is an explicit rational function
$\rho$ (a \emph{density of weight one}) with
\[
  h_2(\Psi(\mathbf{z})) = \rho(\mathbf{z})\cdot h_1(\mathbf{z}),
\]
and the zero (and pole) locus of $\rho$ is exactly the locus along which
$\Psi$ or $\Psi^{-1}$ is undefined. Explicitly:
\begin{rlist}
\item for $T\colon E_J\dashrightarrow\EW$
  (Theorem~\ref{thm:factor}), $\rho_T=\dfrac{4q^2}{u^6}F_2(u,v)$, whose
  zero locus is $\{A=0\}\cup\{u=0\}$;
\item for $\Phi\colon\Es\dashrightarrow\EM$ (Theorem~\ref{thm:density}),
  $\rho_\Phi=\det D\Phi=\dfrac{1}{c-X}$, whose pole locus is $\{X=c\}$.
\end{rlist}
\end{proposition}

The two case studies now diverge in how they exploit
Proposition~\ref{prop:density-principle}: Case~I resolves $\{\rho_T=0\}$
by blow-up (\S\ref{dd-sec:blowup}--\S\ref{sec:main}); Case~II absorbs the
pole of $\rho_\Phi$ into a time reparametrization
(\S\ref{dd-sec:flow}). We treat them in turn.

\begin{remark}
In representation-theoretic language, $\rho$ is precisely the factor by
which a polynomial of a fixed degree transforms under a linear or
projective change of coordinates that is not volume-preserving: this is
the classical notion of a \emph{relative invariant}, and
$g_4,g_3,F_W,F_4$ are relative invariants of weight one under $\Phi,T$
respectively, in the same sense that a volume form is a relative
invariant (density) of weight one under a diffeomorphism. This
observation is what allows the same computation
(Theorem~\ref{thm:density}) to be reinterpreted, in
\S\ref{subsec:centroaffine}, in the language of centro-affine
information geometry.
\end{remark}

\subsection{Case I: Blow-Up at the Singular Loci of the Density}\label{dd-sec:blowup}

We perform blow-ups in the sense of algebraic geometry. Recall that for
a Noetherian scheme $X$ and a closed subscheme $Z\hookrightarrow X$
defined by an ideal sheaf $\calI_Z\subset\OO_X$, the \emph{blow-up}
$\Bl_ZX=\Proj_X\bigl(\bigoplus_{n\ge0}\calI_Z^n\bigr)$ comes with a
proper birational morphism $\pi\colon\Bl_ZX\to X$ whose exceptional
divisor is $E=\pi^{-1}(Z)\cong\Proj_Z(\calI_Z/\calI_Z^2)$, the
projectivized normal cone. In our one-dimensional setting $E$ is always
isomorphic to $\PP^1_k$ over a point.

\begin{figure}[htbp]
\centering
\begin{tikzpicture}[>=Stealth, scale=1.0]
 
  \draw[thick, blue!70!black] (-2.6,-1.4) .. controls (-1.2,-0.2) and (-0.4,0.2) .. (0,0)
      .. controls (0.4,-0.2) and (1.2,1.0) .. (2.6,1.6);
  \filldraw[black] (0,0) circle (1.6pt) node[below left] {$P_0$};
  \node at (-2.0,0.7) {$\EW$};
  \draw[->] (3.1,0.3) -- node[above]{$\pi_A$} (4.5,0.3);
 
  \begin{scope}[xshift=8.4cm]
    \draw[fill=gray!12] (-0.5,-1.9) rectangle (0.5,1.9);
    \draw[thick] (0,-1.9) -- (0,1.9);
    \node[above right=0pt and 2pt] at (0,1.9) {$E_A\cong\PP^1$};
    \draw[thick, blue!70!black] (-2.4,-1.6) .. controls (-1.0,-0.2) and (-0.2,0.15) .. (0,0.55)
        .. controls (0.3,0.9) and (1.1,1.5) .. (2.4,1.9);
    \filldraw[black] (0,0.55) circle (1.6pt) node[right=3pt] {$\xi_0=f'(x_0)/2y_0$};
    \node[below] at (0,-2.3) {strict transform of $\EW$};
  \end{scope}
\end{tikzpicture}
\caption{Blow-up of $\EW$ at a smooth point $P_0=(x_0,y_0)\in\{A=0\}$.
The exceptional divisor $E_A\cong\PP^1$ parametrizes tangent directions
at $P_0$ in the ambient plane; the strict transform of $\EW$ meets $E_A$
transversally at the single point $\xi_0$ corresponding to the actual
tangent line of $\EW$ at $P_0$ (Lemma~\ref{lem:exp-A}). The coordinate
$\xi$ along $E_A$ is the natural coordinate carrying the
Kullback--Leibler-type potential of Theorem~\ref{thm:KL}.}
\label{fig:blowup}
\end{figure}
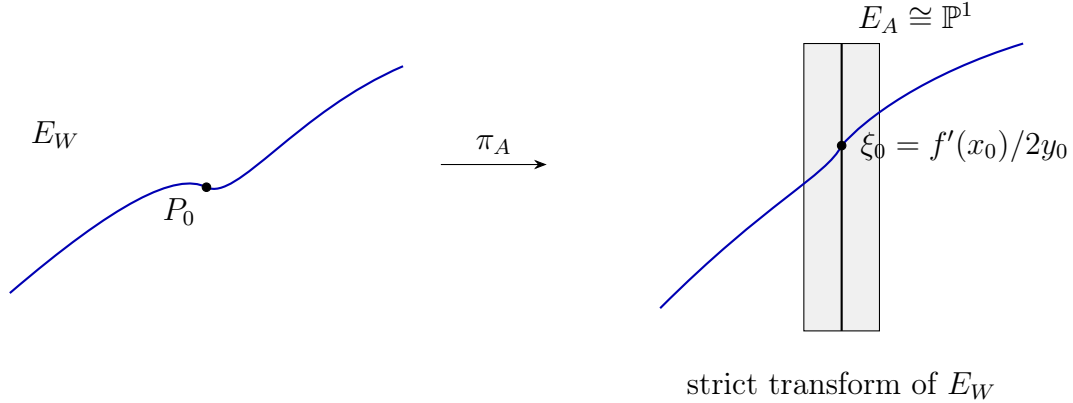

\subsubsection{Blow-up of $\EW$ at the Connell locus $\{A=0\}$}

Let $P_0=(x_0,y_0)\in\EW(k)$ with $A(x_0)=0$; since $\EW$ is smooth,
$P_0$ is a smooth point and $\OO_{\EW,P_0}$ is a DVR.

\begin{definition}\label{def:bup-A}
Let $\calI_{P_0}\subset\OO_{\EW}$ be the ideal sheaf of $P_0$. The
\emph{blow-up} of $\EW$ at $P_0$ is
$\blup{\EW}=\Bl_{P_0}\EW=\Proj_{\EW}\bigl(\bigoplus_{n\ge0}\calI_{P_0}^n\bigr)$.
In local affine coordinates $(t,\xi)$ with $t=x-x_0$, $\xi=(y-y_0)/t$,
this is the Zariski closure of $\{(t,y)\in\EW:t\ne0\}$. The exceptional
divisor is $E_A=\pi_A^{-1}(P_0)\cong\PP^1_k$, and $\pi_A$ is an
isomorphism away from $E_A$ (Figure~\ref{fig:blowup}).
\end{definition}

\begin{lemma}[Local expansion at $P_0$]\label{lem:exp-A}
In the chart $(t,\xi)$, the strict transform of $\EW$ is
\[
  (y_0+t\xi)^2 - f(x_0+t)
  \;=\;
  t\bigl(2y_0\xi - f'(x_0)\bigr)
  + t^2\bigl(\xi^2 - \tfrac{1}{2}f''(x_0)\bigr)
  + O(t^3),
\]
where $f$ is the right-hand side of the Weierstrass equation. The strict
transform meets $E_A=\{t=0\}$ at the single point
$\xi_0 = f'(x_0)/(2y_0)$.
\end{lemma}

\begin{proof}
Taylor-expand $f$ at $x_0$, use $y_0^2=f(x_0)$, and divide by $t$.
\end{proof}

\subsubsection{Blow-up of $E_J$ at the flex locus $\{u=0\}$}

The point $\{u=0\}$ on $E_J$ is smooth (the flex point at infinity),
corresponding via $T$ to the flex point at infinity on $\EW$, consistent
with the triple root of Corollary~\ref{cor:mult}.

\begin{definition}\label{def:bup-u}
Let $Q_0=(0,q)\in E_J(k)$, with ideal sheaf $\calI_{Q_0}$. The blow-up
$\blup{E_J}=\Bl_{Q_0}E_J$ has exceptional divisor
$E_\infty=\pi_\infty^{-1}(Q_0)\cong\PP^1_k$. In coordinates
$\tau=u$, $\sigma=(v+q)/\tau$ (so $v=\tau\sigma-q$), the blow-up chart is
$\{(\tau,\sigma):F_4(\tau,\tau\sigma-q)=0\}$.
\end{definition}

\begin{lemma}\label{lem:exp-u}
$F_4(\tau,\tau\sigma-q)=-\tau\bigl(a\tau^3+b\tau^2+c\tau+d+2q\sigma-\sigma^2\tau\bigr)$.
On $E_\infty=\{\tau=0\}$ the residual equation forces $\sigma_0=-d/(2q)$.
\end{lemma}

\begin{proof}
Expand $(\tau\sigma-q)^2-a\tau^4-b\tau^3-c\tau^2-d\tau-q^2$ and factor $\tau$.
\end{proof}

\subsubsection{Blow-up at the branch points $\{y=0\}$}

Let $(x_i,0)\in\EW(k)$ be a branch point with $f(x_i)=0$, $f'(x_i)\ne0$;
assume $a\ne0$.

\begin{definition}\label{def:bup-br}
Let $\calI_{B_i}$ be the ideal sheaf of $B_i=(x_i,0)$. The blow-up
$\Bl_{B_i}\EW$ has exceptional divisor $E_{B_i}\cong\PP^1_k$. In local
coordinates $U=A/(2qy)$, $W=1/U$, $\widetilde{V}=V/U^2$, the Jacobi
quartic becomes $\widetilde{V}^2=a+bW+cW^2+dW^3+q^2W^4$, and the blow-up
chart near $W=0$ uses slope coordinate $\omega=(\widetilde{V}-\sqrt a)/W$.
\end{definition}

\begin{lemma}\label{lem:exp-br}
Setting $x=x_i+\tau^2$ and $v=\tau\sigma-q$ in $F_4$ gives the same
expansion structure as Lemma~\ref{lem:exp-u} with $d$ replaced by
$b/(2\sqrt a)$; on $E_{B_i}=\{W=0\}$ the residual equation forces
$\omega_0=b/(2\sqrt a)$.
\end{lemma}

\begin{proof}
Direct substitution and factoring of $\tau$, as in Lemma~\ref{lem:exp-u}.
\end{proof}

We now have three exceptional divisors, $E_A$, $E_\infty$, $E_{B_i}$,
each isomorphic to $\PP^1_k$, each carrying a distinguished slope
coordinate ($\xi$, $\sigma$, $\omega$ respectively) at which the strict
transform meets it. Section~\ref{sec:dual} shows that all three carry
canonically isomorphic information geometries.

\subsection{Dual Geometries on the Exceptional Divisors}\label{sec:dual}

\subsubsection{KL-type potentials}

\begin{definition}[KL-type dual geometry]\label{def:KL}
A \emph{KL-type dual geometry} on an interval $\mathcal U\subset\R$ is a
triple $(\mathcal U,\Phi,\Phi^*)$ with $\Phi(\xi)=\log(\alpha\xi-\beta)$
($\alpha\ne0$, $\alpha\xi>\beta$), Legendre dual
$\Phi^*(\theta)=\log\theta+\text{linear}$, and Bregman divergence
$D_\Phi(p\Vert q)=\Phi(p)-\Phi(q)-\Phi'(q)(p-q)$ equal to
$\log(r_1/r_2)-(r_1/r_2-1)$ under $r_i=\alpha\xi_i-\beta>0$---i.e.\
isomorphic to the Kullback--Leibler divergence between exponential
distributions with rates $r_1,r_2$.
\end{definition}

\begin{theorem}[Universal KL structure]\label{thm:KL}
Each exceptional divisor of \S\ref{dd-sec:blowup} carries a KL-type dual
geometry:
\begin{align*}
  \textup{At }A=0\colon\quad
    &\Phi_A(\xi)=\log|2y_0\xi-f'(x_0)|, &
    \theta_A&=\frac{2y_0}{2y_0\xi-f'(x_0)},\\
  \textup{At }u=0\colon\quad
    &\Phi_\infty(\sigma)=\log|\sigma+d/(2q)|, &
    \theta_\infty&=\frac{1}{\sigma+d/(2q)},\\
  \textup{At }y=0\colon\quad
    &\Phi_b(\omega)=\log|2\sqrt a\,\omega-b|, &
    \theta_b&=\frac{2\sqrt a}{2\sqrt a\,\omega-b}.
\end{align*}
In each case $\Phi$ is strictly concave, $\theta=\Phi'$ is a M\"obius
diffeomorphism, and the Legendre dual is $\Phi^*(\theta)=\log\theta+
(\beta/\alpha)\theta+1$ where $(\alpha,\beta)$ takes the respective
values $(2y_0,f'(x_0))$, $(1,-d/(2q))$, $(2\sqrt a, b)$.
\end{theorem}

\begin{proof}
Each $\Phi$ has the form $\log(\alpha\xi-\beta)$. Then
$\Phi''=-\alpha^2/(\alpha\xi-\beta)^2<0$ gives strict concavity;
$\theta=\Phi'=\alpha/(\alpha\xi-\beta)$ is a M\"obius map, hence a
diffeomorphism onto its image; and $\xi(\theta)=1/\theta+\beta/\alpha$
gives $\Phi^*(\theta)=\theta\xi(\theta)-\Phi(\xi(\theta))
=\log\theta+(\beta/\alpha)\theta+1$. The explicit values of
$(\alpha,\beta)$ are read off from Lemmas~\ref{lem:exp-A},
\ref{lem:exp-u}, \ref{lem:exp-br}.
\end{proof}

\subsubsection{Fisher information metric}

For a one-dimensional statistical manifold with potential $\Phi$, the
\emph{Fisher information metric} is $g^F=|\Phi''(\xi)|\,d\xi^2$. For
$\Phi(\xi)=\log(\alpha\xi-\beta)$,
\[
  g^F = \frac{\alpha^2}{(\alpha\xi-\beta)^2}\,d\xi^2.
\]
In the \emph{natural parameter} $r=\alpha\xi-\beta>0$ (so
$\xi=(r+\beta)/\alpha$, $d\xi=dr/\alpha$) this becomes
$g^F=dr^2/r^2$, the standard Fisher metric of the exponential family
$\{r\,e^{-rt}:r>0\}$---unit Fisher information at every point in the
$r$-coordinate.

\begin{theorem}[Equality of Fisher metrics]\label{thm:fisher}
Under the natural parameters
\[
  r_A=2y_0\xi-f'(x_0), \qquad r_\infty=\sigma+d/(2q), \qquad
  r_b=2\sqrt a\,\omega-b,
\]
all three Fisher metrics equal $dr^2/r^2$ as Riemannian metrics on
$(0,\infty)$. In particular the isometries
\[
  \varphi_{AJ}\colon(E_A,g_A^F)\to(E_\infty,g_\infty^F),\ r_A\mapsto r_\infty=r_A,
  \qquad
  \varphi_{Ab}\colon(E_A,g_A^F)\to(E_b,g_b^F),\ r_A\mapsto r_b=r_A
\]
exhibit all three statistical manifolds as isometric to the hyperbolic
line $\bigl((0,\infty),dr^2/r^2\bigr)$.
\end{theorem}

\begin{proof}
The substitution $r=\alpha\xi-\beta$ transforms $g^F$ to $dr^2/r^2$ in
every case (direct computation), so the three metrics agree as abstract
Riemannian metrics on $(0,\infty)$; the maps $\varphi_{AJ},\varphi_{Ab}$
realize this agreement explicitly.
\end{proof}

\begin{remark}[Non-triviality]\label{rem:nontrivial}
Although all potentials $\log(\alpha\xi-\beta)$ are related by an affine
change of variable, the coordinate $\xi$ (resp.\ $\sigma$, $\omega$) on
each exceptional divisor is \emph{determined by the blow-up geometry}
of $\EW$ (resp.\ $E_J$)---the slope of the strict transform at the
blown-up point. That the canonical coordinates on all three exceptional
divisors independently produce the same Fisher metric $dr^2/r^2$
reflects the universal role of the Connell datum $A(x)$ in
Theorem~\ref{thm:factor}.
\end{remark}

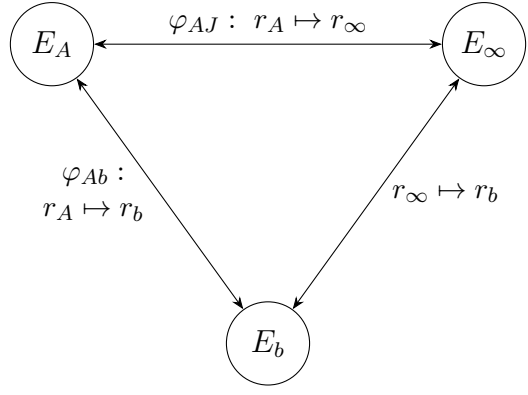
\begin{figure}[htbp]
\centering
\begin{tikzpicture}[>=Stealth, node distance=3.4cm]
\node[draw, circle, minimum size=1.1cm] (EA) {$E_A$};
\node[draw, circle, minimum size=1.1cm, right=4.6cm of EA] (Einf) {$E_\infty$};
\node[draw, circle, minimum size=1.1cm, below=3.4cm of $(EA)!0.5!(Einf)$] (Eb) {$E_b$};
\draw[<->] (EA) -- node[above, font=\small] {$\varphi_{AJ}:\ r_A\mapsto r_\infty$} (Einf);
\draw[<->] (EA) -- node[left=2pt, font=\small, align=center] {$\varphi_{Ab}:$\\ $r_A\mapsto r_b$} (Eb);
\draw[<->] (Einf) -- node[right=2pt, font=\small, align=center] {$r_\infty\mapsto r_b$} (Eb);
\node[align=center, font=\scriptsize, below=0.5cm of Eb] {all isometric to $\bigl((0,\infty),\,dr^2/r^2\bigr)$};
\end{tikzpicture}
\caption{The three exceptional divisors from Figure~\ref{fig:blowup}
(over the Connell locus $A=0$ on $\EW$, the flex $u=0$ on $E_J$, and a
branch point $y=0$ on $\EW$) carry canonically isomorphic KL-type dual
geometries, all isometric to the hyperbolic line $(0,\infty)$ with
metric $dr^2/r^2$ (Theorem~\ref{thm:fisher}). The horizontal isomorphism
$\varphi_{AJ}$ is realized geometrically by the extension of $T$
(Theorem~\ref{thm:main}).}
\label{fig:three-divisors}
\end{figure}

\subsection{Isomorphism of Dual Geometries}\label{sec:main}

\begin{definition}\label{def:iso}
Two KL-type dual geometries $(\mathcal U_i,\Phi_i,\Phi_i^*)$ ($i=1,2$)
are \emph{isomorphic} if there is an affine bijection
$\varphi(\xi)=\alpha\xi+\beta$ with $\Phi_2(\varphi(\xi))=\Phi_1(\xi)+\text{const}$.
\end{definition}

\begin{theorem}[Main theorem]\label{thm:main}
\begin{rlist}
\item The three KL-type dual geometries of Theorem~\ref{thm:KL} are
  mutually isomorphic (Definition~\ref{def:iso}) and mutually isometric
  as Riemannian manifolds (Theorem~\ref{thm:fisher}).
\item The birational map $T\colon E_J\dashrightarrow\EW$ extends to a
  morphism of $k$-schemes $\blup T\colon\blup{E_J}\to\blup{\EW}$ making
  the following diagram commute, where the bottom arrow is the identity
  isometry $r_\infty=r_A$:
\[
\begin{tikzcd}[column sep=4.4em, row sep=3em]
  \blup{E_J} \arrow[r, "\blup T"] \arrow[d, "\Phi_\infty"'] &
  \blup{\EW} \arrow[d, "\Phi_A"] \\
  \bigl((0,\infty),\,\tfrac{dr^2}{r^2}\bigr) \arrow[r, "\sim"'] &
  \bigl((0,\infty),\,\tfrac{dr^2}{r^2}\bigr)
\end{tikzcd}
\]
\end{rlist}
\end{theorem}

\begin{proof}
\textbf{Step 1 ($\Phi_A\cong\Phi_\infty$).} With $r_A=2y_0\xi-f'(x_0)$
and $r_\infty=\sigma+d/(2q)$, define
$\xi(\sigma)=\bigl(\sigma+d/(2q)+f'(x_0)\bigr)/(2y_0)$. Then
$r_A(\xi(\sigma))=r_\infty$, so $\Phi_A(\xi(\sigma))=\Phi_\infty(\sigma)$.

\textbf{Step 2 ($\Phi_A\cong\Phi_b$).} With $r_b=2\sqrt a\,\omega-b$,
define $\xi(\omega)=\bigl(2\sqrt a\,\omega-b+f'(x_0)\bigr)/(2y_0)$; then
$r_A(\xi(\omega))=r_b(\omega)$.

\textbf{Step 3 (commutativity).} In the blow-up chart
$(\tau,\sigma)$ of $\blup{E_J}$ and $(t,\xi)$ of $\blup\EW$, the map $T$
in coordinates $(u,v)=(\tau,\tau\sigma-q)$ reads
\[
  x(T)=\frac{2q\sigma+d}{\tau}+x_0(\tau), \qquad
  y(T)=\frac{8q^3\tau\sigma+4q^2(c\tau^2+d\tau)-d^2\tau^2}{2q\tau^3},
\]
where $x_0(\tau)\to x_0$ as $\tau\to0$. Thus $t=x-x_0\sim(2q\sigma+d)/\tau$
and $\xi=y/t$ extends to a rational function on $\blup{E_J}$; restricting
to $E_\infty=\{\tau=0\}$,
\[
  \xi_{\blup T}(\sigma)=\lim_{\tau\to0}\frac{y(T)}{x(T)-x_0}
  =\frac{2q\sigma+d+f'(x_0)}{2y_0},
\]
which matches $\xi(\sigma)$ of Step~1, so
$\Phi_A\circ\blup T|_{E_\infty}=\Phi_\infty$. The extension $\blup T$ is
a morphism of $k$-schemes because $T$ is a morphism away from the
exceptional loci and the blow-up is the universal scheme resolving the
indeterminacy.

\textbf{Step 4 (KL-divergence isomorphism).} Under $r_i>0$, every
Bregman divergence becomes
$D(r_1\Vert r_2)=\log(r_1/r_2)-(r_1/r_2-1)$, the KL divergence between
exponential distributions of rates $r_1,r_2$. Combined with
Theorem~\ref{thm:fisher}, this exhibits all three geometries as
isomorphic statistical manifolds with the same Fisher metric and
divergence.
\end{proof}

\begin{corollary}[Invariance under the theta transformation]\label{cor:inv}
The KL-type dual geometry is preserved by $T$: the blow-up map
$\blup T$ is an isomorphism of statistical manifolds
$(E_\infty,g^F_\infty,D_\infty)\xrightarrow{\ \sim\ }(E_A,g_A^F,D_A)$.
\end{corollary}

\begin{remark}
Corollary~\ref{cor:inv} concerns only the specific transformation $T$.
Whether the KL-type dual geometry is an invariant of an \emph{arbitrary}
birational equivalence between elliptic curves---e.g.\ a general isogeny,
or the Mordell transformation $\Phi$ of \S\ref{subsec:caseII-models}
itself---is not addressed here and is left open in
\S\ref{dd-sec:discussion}. A first step toward such a result would be a
blow-up construction compatible with an arbitrary birational map,
possibly via the minimal regular model over a suitable base.
\end{remark}

\begin{remark}[Universality of $\log(\text{linear})$]
The potential $\log(\alpha\xi-\beta)$ is, up to affine reparametrization,
the unique one-dimensional potential of the form
$\log(\text{affine function})$. This uniqueness, together with the fact
that every blow-up expansion of \S\ref{dd-sec:blowup} produces a
logarithmic potential of exactly this form, explains the universal
appearance of the same KL structure and Fisher metric at every
exceptional divisor.
\end{remark}

\subsection{Case II: The Group Law as a Matrix Operation}\label{sec:group}

We now turn to the second case study, and show that the density
principle of Theorem~\ref{thm:density} has an elementary, purely
algebraic companion: the group law of $\Es$ is realized by conjugation
and translation of the matrix $A_4$ of Definition~\ref{def:matrices}.

By Corollary~\ref{cor:inverse}, the family of matrices $A_4$ compatible
with a fixed cubic $\Es$ (i.e.\ fixed $g_2,g_3$) is parametrized
bijectively by the points $P=(c,d)\in\Es(k)$; write $A_4(P)$ for the
corresponding matrix.

\begin{theorem}[Realization of inversion]\label{thm:inversion}
Let $J=\mathrm{diag}(1,-1,1)$. For the inversion map $P=(c,d)\mapsto
-P=(c,-d)$ on $\Es$,
\[
  A_4(-P) = J\,A_4(P)\,J.
\]
\end{theorem}

\begin{proof}
$A_4(P)$ depends on $d$ only through the off-diagonal entry $-2d$ in
positions $(1,2)$ and $(2,1)$. Conjugation by $J=\mathrm{diag}(1,-1,1)$
sends $A_{ij}\mapsto J_iJ_jA_{ij}$, which flips the sign of exactly the
$(1,2)$ and $(2,1)$ entries and fixes every other entry---precisely the
operation $d\mapsto-d$.
\end{proof}

\begin{theorem}[Realization of addition]\label{thm:addition}
Fix $Q=(x_2,y_2)\in\Es$. For $P=(c,d)\in\Es$, set
\[
  \lambda=\frac{y_2-d}{x_2-c}, \qquad
  c'=\frac{\lambda^2}{4}-c-x_2, \qquad
  d'=-\bigl[\lambda(c'-c)+d\bigr], \qquad
  e'=g_2-3c'^2.
\]
Then $(c',d')\in\Es$ and $P+Q=(c',d')$ in the group law of $\Es$; hence
\[
  A_4(P+Q)=\begin{pmatrix}-e'&-2d'&0\\-2d'&6c'&0\\0&0&1\end{pmatrix}.
\]
\end{theorem}

\begin{proof}
Let $f(X)=4X^3-g_2X-g_3$. The line through $P=(c,d)$ and $Q=(x_2,y_2)$
(for $c\ne x_2$) is $Y=\lambda(X-c)+d$. Substituting into $Y^2=f(X)$
gives a cubic in $X$ with leading coefficient $4$ and quadratic
coefficient $-\lambda^2$; its three roots are the $X$-coordinates of the
three collinear intersection points $c,x_2,X_3$, so by Vieta's formula
\[
  c+x_2+X_3=\frac{\lambda^2}{4}
  \quad\Longrightarrow\quad
  X_3=\frac{\lambda^2}{4}-c-x_2=c'.
\]
The corresponding $Y$-coordinate on the line is
$Y_3=\lambda(X_3-c)+d=\lambda(c'-c)+d$. Since $\Es\colon Y^2=f(X)$ is
invariant under $Y\mapsto-Y$, and by definition of the group law three
collinear points sum to the identity $O$, we get $P+Q+(X_3,Y_3)=O$, i.e.\
$P+Q=(X_3,-Y_3)=(c',d')$. (When $Q=P$, replace the chord by the tangent
$\lambda=f'(c)/(2d)=(12c^2-g_2)/(2d)$; the same argument applies
verbatim.)
\end{proof}

\begin{remark}
Thus the family $\{A_4(P):P\in\Es\}$ is not merely a collection of
matrices; it is a space carrying the group structure of $\Es$ itself:
$A_4(P)\mapsto A_4(-P)$ is inversion, and $A_4(P)\mapsto A_4(P+Q)$ is
translation, both realized as explicit matrix operations.
\end{remark}

\begin{figure}[htbp]
\centering
\begin{tikzpicture}[>=Stealth, scale=1.05]
\begin{scope}
\clip (-3.0,-2.5) rectangle (3.4,2.5);
\draw[thick, blue!70!black, domain=-1.0:3.2, samples=160, smooth]
  plot (\x, {sqrt(max(\x*\x*\x-3*\x+3,0))});
\draw[thick, blue!70!black, domain=-1.0:3.2, samples=160, smooth]
  plot (\x, {-sqrt(max(\x*\x*\x-3*\x+3,0))});
\end{scope}
\coordinate (P) at (-0.55,1.35);
\coordinate (Q) at (1.55,2.10);
\coordinate (R) at (2.55,-2.63);
\coordinate (Rp) at (2.55,2.63);
\filldraw[black] (P) circle (1.3pt) node[below left] {$P=(c,d)$};
\filldraw[black] (Q) circle (1.3pt) node[above right] {$Q$};
\filldraw[black] (R) circle (1.3pt) node[below left=1pt] {$(X_3,Y_3)$};
\filldraw[black] (Rp) circle (1.3pt) node[above left=1pt] {$P+Q=(c',d')$};
\draw[dashed, gray!70] (P) -- (Q) -- (R);
\draw[dashed, gray!70] (R) -- (Rp);
\node[gray!70, font=\scriptsize] at (3.15,0) {reflect};
\node[gray!70, font=\scriptsize] at (3.15,-0.32) {$Y\mapsto -Y$};
\end{tikzpicture}
\caption{The chord-and-tangent construction underlying
Theorem~\ref{thm:addition}: the line through $P$ and $Q$ meets $\Es$ at a
third point $(X_3,Y_3)$, and $P+Q$ is its reflection across the
$X$-axis. Theorem~\ref{thm:addition} packages the classical Vieta
computation as an explicit map $A_4(P)\mapsto A_4(P+Q)$.}
\label{fig:group-law}
\end{figure}
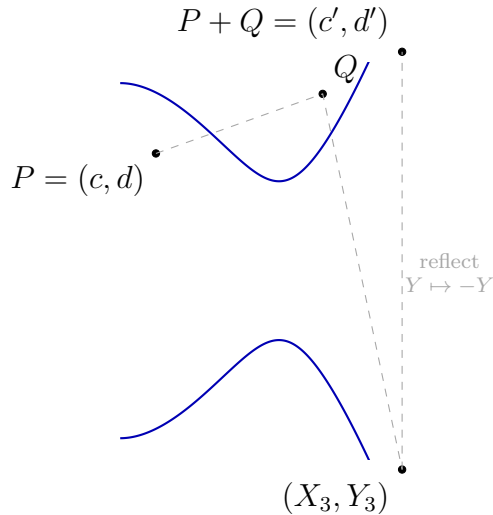

\subsubsection{Eigenvalue characterization}

\begin{theorem}[Eigenvalue theorem]\label{thm:eigen}
The $(x,y)$-block of $A_4(P)$ is
\[
  \mathrm{diag}(6c,1)=\mathrm{diag}(6X(P),1),
\]
and $b_4$ points along the eigenvector for eigenvalue $6X(P)$. The
$(X,Y)$-block of $A_3$ is always $\mathrm{diag}(0,-1)$, independently of
the base point, and $b_3$ points along the eigenvector for eigenvalue
$0$.
\end{theorem}

\begin{proof}
Immediate from Definition~\ref{def:matrices}: the $(2,2),(3,3)$ entries
of $A_4$ are $6c,1$ with vanishing $(2,3)$ entry, and $b_4=(0,1,0)^\top$
is the eigenvector for $6c$; likewise for $A_3$ with eigenvalue $0$.
\end{proof}

\begin{corollary}
The vanishing of the eigenvalue on the $A_3$ side reflects that $\Es$ is
given in a base-point-free normal form; all base-point dependence is
concentrated in the single real eigenvalue $6c=6X(P)$ on the $A_4$ side,
whose evolution under the group law is given explicitly by
Theorem~\ref{thm:addition}.
\end{corollary}

\subsection{Invariance of the Canonical Differential}\label{sec:diff}

\begin{theorem}[Birational invariance of $\omega$]\label{thm:differential}
On $\Es$, the pullback under $\Phi$ satisfies
\[
  \frac{dx}{y}=\frac{dX}{Y}.
\]
\end{theorem}

\begin{proof}
On the curve, $2Y\,dY=(12X^2-g_2)\,dX$, so
$dY/dX=(12X^2-g_2)/(2Y)$. Differentiating $x=\Phi_1(X,Y(X))$ along the
curve,
\[
  \frac{dx}{dX}=\frac{\partial x}{\partial X}+\frac{\partial x}{\partial Y}\cdot\frac{dY}{dX}.
\]
Writing $\tfrac1y\tfrac{dx}{dX}=N/D$ with
\[
  N=-2Y(Y-d)-(X-c)(-12X^2+3c^2+e), \qquad
  D=4(X-c)^2(2X+c)-(Y-d)^2,
\]
direct polynomial expansion gives the identity
$N-D=-\bigl(Y^2-(4X^3-g_2X-g_3)\bigr)$. On the curve
$Y^2-4X^3+g_2X+g_3=0$, so $N=D$ and $\tfrac1y\tfrac{dx}{dX}=1$, i.e.\
$dx/y=dX/Y$.
\end{proof}

\begin{remark}
This is the classical fact that $\omega=dX/Y$ is the canonical regular
differential of $\Es$, invariant under birational (in particular,
isomorphic) changes of model; it induces the flat metric $|du|^2$
($u=\int\omega$) of the uniformization $\Es\cong\mathbb C/\Lambda$.
Theorem~\ref{thm:differential} is a direct verification for the specific
map $\Phi$.
\end{remark}

\subsection{Case II Continued: A Birationally Invariant Gradient Flow}\label{dd-sec:flow}

We now show how the density weight $\det D\Phi=1/(c-X)$ of
Theorem~\ref{thm:density} governs the failure, and the repair, of
birational invariance for a gradient-type flow toward the curve.

\begin{definition}[Branin flow]\label{def:branin}
For $g\colon\R^2\to\R$, the \emph{Branin flow} associated to $g$ is
\[
  \frac{d}{dt}\binom{x}{y}=-\frac{\nabla g(x,y)}{|\nabla g(x,y)|^2}\,g(x,y).
\]
Its solutions satisfy $g(x(t),y(t))=g(x_0,y_0)e^{-t}$: the flow decays
exponentially toward the level set $\{g=0\}$ while moving along the
steepest-descent direction of $g$.
\end{definition}

\begin{figure}[htbp]
\centering
\begin{tikzpicture}[>=Stealth, scale=1.0]
\begin{scope}
\clip (-3.1,-2.3) rectangle (3.4,2.3);
\draw[thick, blue!70!black, domain=-1.0:3.2, samples=140, smooth]
  plot (\x, {sqrt(max(\x*\x*\x-3*\x+3,0))});
\draw[thick, blue!70!black, domain=-1.0:3.2, samples=140, smooth]
  plot (\x, {-sqrt(max(\x*\x*\x-3*\x+3,0))});
\end{scope}
\foreach \s/\a in {1.6/-45,1.9/-25,0.6/165,-0.6/155}{
  \pgfmathsetmacro{\ex}{1.9+0.9*cos(\a)}
  \pgfmathsetmacro{\ey}{-0.2+0.9*sin(\a)}
}
\draw[->, orange!80!black, thick] (2.6,1.9) .. controls (2.2,1.0) .. (1.85,0.35);
\draw[->, orange!80!black, thick] (2.9,-1.6) .. controls (2.4,-0.6) .. (1.9,-0.15);
\draw[->, orange!80!black, thick] (-2.1,1.5) .. controls (-1.4,1.0) .. (-0.55,0.75);
\node[orange!80!black, font=\scriptsize, align=center] at (2.95,2.0) {naive flow\\ time $t$};
\node[below, font=\small] at (0,-2.4) {Branin flow toward $\{g_4=0\}$};
\end{tikzpicture}
\hspace{0.3cm}
\begin{tikzpicture}[>=Stealth, scale=1.0]
\begin{scope}
\clip (-3.1,-2.3) rectangle (3.4,2.3);
\draw[thick, red!70!black, domain=-1.0:3.2, samples=140, smooth]
  plot (\x, {sqrt(max(\x*\x*\x-3*\x+3,0))});
\draw[thick, red!70!black, domain=-1.0:3.2, samples=140, smooth]
  plot (\x, {-sqrt(max(\x*\x*\x-3*\x+3,0))});
\end{scope}
\draw[->, teal!70!black, thick] (2.6,1.9) .. controls (2.15,1.05) .. (1.85,0.4);
\draw[->, teal!70!black, thick] (2.9,-1.6) .. controls (2.35,-0.65) .. (1.9,-0.1);
\draw[->, teal!70!black, thick] (-2.1,1.5) .. controls (-1.35,1.05) .. (-0.5,0.8);
\node[teal!70!black, font=\scriptsize, align=center] at (2.95,2.0) {corrected flow\\ time $\tau$};
\node[below, font=\small, align=center] at (0,-2.4) {reparametrized flow\\ toward $\{g_3=0\}$};
\end{tikzpicture}
\caption{Schematic: a naive Branin flow toward the quartic $\{g_4=0\}$
(left) does not push forward, under $\Phi^{-1}$, to the naive Branin
flow toward the cubic $\{g_3=0\}$; the discrepancy is exactly the
density weight $\det D\Phi=1/(c-X)$. Theorem~\ref{thm:invariant-flow}
shows that replacing $t$ by $\tau(t)=t-\log|(c-X(t))/(c-X(0))|$ repairs
the decay law exactly, so the reparametrized, pushed-forward flow
(right) obeys $g_3(\tau)=g_3(0)e^{-\tau}$.}
\label{fig:branin}
\end{figure}
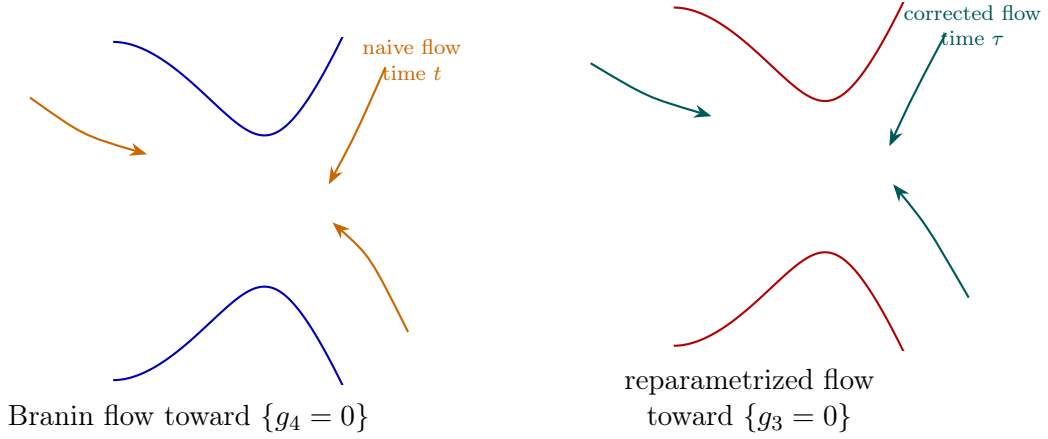

\begin{remark}
Definition~\ref{def:branin} is the scalar ($m=1$) case of a system
$\dot x=-J_g^+(x)g(x)$, where $J_g^+$ is a generalized inverse of the
Jacobian of a mapping $g\colon\R^n\to\R^m$; Branin~\cite{Branin1972}
originally treated only the square case $m=n$. Tanabe~\cite{Tanabe1979}
extended the construction to the underdetermined case $m\le n$ relevant
here, proved the exact exponential decay law
$g(x(t,x^0))=e^{-t}g(x^0)$ of Definition~\ref{def:branin} in this
generality (his identity~(3)), and gave a full stability and
convergence analysis of the resulting flow toward the solution set
$\{g=0\}$ (his Lemmas~1--3 and Theorems~4--5).\footnote{The author thanks
Professor Kunio Tanabe for kindly bringing this reference to his attention.}
His worked example,
$g(x_1,x_2)=x_1-x_2^2$, is exactly the ``naive'' flow toward a smooth
plane curve pictured in Figure~\ref{fig:branin}; the present paper's
contribution is to track what happens to this decay law under the
birational change of model $\Phi$, rather than for a single fixed $g$.
\end{remark}

If $(x(t),y(t))$ solves the Branin flow for $g_4$, then simply pushing
the trajectory forward to $(X(t),Y(t)):=\Phi^{-1}(x(t),y(t))$ does
\emph{not}, in general, solve the naive Branin flow for $g_3$; the
mismatch is exactly measured by the density weight of
Theorem~\ref{thm:density}.

\begin{theorem}[Time reparametrization repairs invariance]\label{thm:invariant-flow}
Let $(x(t),y(t))$ solve the Branin flow for $g_4$ and set
$(X(t),Y(t)):=\Phi^{-1}(x(t),y(t))$. With
\[
  \tau(t):=t-\log\left|\frac{c-X(t)}{c-X(0)}\right|,
\]
we have exactly
\[
  g_3(X(t),Y(t)) = g_3(X(0),Y(0))\,e^{-\tau(t)}.
\]
That is, replacing $t$ by $\tau$ makes the pushed-forward trajectory
obey precisely the Branin decay law for $g_3$.
\end{theorem}

\begin{proof}
By Definition~\ref{def:branin}, $g_4(x(t),y(t))=g_4(x_0,y_0)e^{-t}$.
By Theorem~\ref{thm:density} (in the form
$g_3=(c-X)\,g_4\circ\Phi$),
\[
  g_3(X(t),Y(t))=(c-X(t))\,g_4(x(t),y(t))=(c-X(t))\,g_4(x_0,y_0)e^{-t}.
\]
At $t=0$, $g_3(X(0),Y(0))=(c-X(0))\,g_4(x_0,y_0)$, so
\[
  g_3(X(t),Y(t)) = g_3(X(0),Y(0))\cdot\frac{c-X(t)}{c-X(0)}\cdot e^{-t}
  = g_3(X(0),Y(0))\,e^{-t+\log|(c-X(t))/(c-X(0))|}.
\]
Setting $\tau(t)=t-\log|(c-X(t))/(c-X(0))|$ gives the claim.
\end{proof}

\begin{proposition}[Directional agreement]\label{prop:direction}
On the curve ($g_3=g_4=0$), the correctly transformed covector
$(D\Phi)^\top\nabla_xg_4|_{\Phi(X,Y)}$ is parallel to $\nabla_Xg_3(X,Y)$.
\end{proposition}

\begin{proof}
Differentiating $g_4(\Phi(X,Y))=\det D\Phi(X,Y)\cdot g_3(X,Y)$ in $(X,Y)$
via the chain rule,
\[
  (D\Phi)^\top\nabla_xg_4\big|_{\Phi(X,Y)}
  = g_3(X,Y)\,\nabla_X(\det D\Phi) + \det D\Phi\cdot\nabla_Xg_3(X,Y).
\]
On the curve $g_3=0$, so the first term vanishes, leaving
$(D\Phi)^\top\nabla_xg_4|_{\Phi(X,Y)}=\det D\Phi\cdot\nabla_Xg_3(X,Y)$, a
scalar multiple of $\nabla_Xg_3(X,Y)$.
\end{proof}

\begin{corollary}[Birationally invariant Branin flow]\label{cor:invariant-flow}
Combining Theorem~\ref{thm:invariant-flow} and
Proposition~\ref{prop:direction}: pushing the Branin flow of $A_4,b_4$
forward by $\Phi^{-1}$ and reparametrizing time by
$\tau=t-\log|c-X(t)|+\log|c-X(0)|$ yields, asymptotically near the
curve, exactly the (covector-correct) Branin flow of $A_3,b_3$.
\end{corollary}

\subsubsection{A Pythagorean-type identity, and its degeneration}

\begin{theorem}[Pythagorean-type correspondence]\label{thm:pythagoras}
Let $\alpha:=\sqrt{6c}\,x$, $\beta:=y$, $h^2:=x^4+4dx+e$. Then identically
\[
  \alpha^2+\beta^2-h^2 = g_4(x,y),
\]
so under pullback by $\Phi$,
\[
  \alpha^2+\beta^2-h^2 = \det D\Phi(X,Y)\cdot g_3(X,Y);
\]
on the curve ($g_3=0$) this becomes the exact Pythagorean-type relation
$\alpha^2+\beta^2=h^2$. On the $A_3,b_3$ side, by contrast, the
corresponding eigenvalue is $0$ (Theorem~\ref{thm:eigen}), so the
``two-legged'' Pythagorean structure degenerates to the single-term
identity $Y^2=4X^3-g_2X-g_3$.
\end{theorem}

\begin{proof}
$\alpha^2+\beta^2-h^2=6cx^2+y^2-(x^4+4dx+e)=y^2-x^4+6cx^2-4dx-e=g_4(x,y)$
(Proposition~\ref{prop:expand}). Substituting
Theorem~\ref{thm:density} gives the pullback identity; on the curve
$g_3=0$ forces $\alpha^2+\beta^2=h^2$. The degeneration on the $A_3$
side follows immediately from Theorem~\ref{thm:eigen}: the eigenvalue
in the $X$-direction is $0$, so that ``leg'' carries no weight, leaving
the single-term relation $Y^2=($cubic in $X)$.
\end{proof}

\begin{corollary}
The Pythagorean defect $\alpha^2+\beta^2-h^2$ transforms with exactly
the same density weight $\det D\Phi=1/(c-X)$ as $g_4$ itself; in
particular its zero locus (the curve) is birationally invariant even
though the defect itself is not a function but a density.
\end{corollary}

\subsection{Synthesis: Two Faces of the Same Exponential Law}\label{sec:synthesis}

We now make explicit the connection anticipated in the introduction. Both
case studies end in an exponential decay law governed by a rate
parameter, and the two laws are, formally, the same formula seen in two
different settings.

\begin{itemize}[leftmargin=1.6em]
\item \textbf{Case I (static/information-geometric).} On each
  exceptional divisor of \S\ref{dd-sec:blowup}, the natural coordinate
  $r>0$ (Theorem~\ref{thm:fisher}) parametrizes the one-parameter
  exponential family $\{r\,e^{-rt}:r>0\}$ underlying
  Definition~\ref{def:KL}: for fixed $r$, the density $r\,e^{-rt}$ decays
  in the auxiliary variable $t$ at rate $r$, and the Kullback--Leibler
  divergence between two such densities, with rates $r_1,r_2$, is
  exactly the Bregman divergence of the potential $\Phi=\log r$.
\item \textbf{Case II (dynamic).} Along the corrected Branin flow of
  Theorem~\ref{thm:invariant-flow}, the defining polynomial itself decays
  as $g(\tau)=g(0)e^{-\tau}$ in the reparametrized time $\tau$: here it
  is $g$, rather than a probability density, that plays the role of the
  decaying quantity, and $\tau$ (not a rate parameter) is the variable.
\end{itemize}

The coincidence is that in \emph{both} settings, resolving the density
principle of Proposition~\ref{prop:density-principle}---by blow-up in
Case~I, by time reparametrization in Case~II---produces an object
governed by the same differential equation $\dot z=-z$, whose solution
is the exponential $z(t)=z(0)e^{-t}$. In Case~I this equation is solved
\emph{along the fictitious time $t$ of the exponential family}, at each
fixed point $r$ of the exceptional divisor; in Case~II it is solved
\emph{along the actual flow time $\tau$}, with $z=g_3(X(\tau),Y(\tau))$.

\begin{remark}[A suggestive parallel, stated as an open direction]
We do not claim a theorem identifying these two exponential laws beyond
the formal analogy above; we record it here because it suggests a
natural question for future work (see also
\S\ref{dd-sec:discussion}): is there a single construction that produces
\emph{both} the KL-type information geometry on a blow-up divisor
\emph{and} the exponential decay of a birationally corrected gradient
flow, as two projections of one object---for instance, by viewing the
rate coordinate $r$ of Case~I as governing the speed of approach, in
reparametrized time $\tau$, of an appropriate flow toward the
corresponding exceptional divisor? A precise formulation would likely
require extending the Branin-flow construction of \S\ref{dd-sec:flow} to
the blown-up surface $\blup\EW$ itself, flowing toward the exceptional
divisor $E_A$ rather than toward a point, and comparing its decay law to
Theorem~\ref{thm:fisher}. We leave this as an open problem.
\end{remark}

\subsubsection{A remark on centro-affine geometry}\label{subsec:centroaffine}

\begin{remark}
Theorem~\ref{thm:density} says that $g_4,g_3$ transform, under $\Phi$,
not as scalar functions but as \emph{relative scalars (densities) of
weight one}. This is mathematically the same transformation law as that
of a volume form $\theta$ in centro-affine (equi-affine) differential
geometry and in the affine-geometric approach to statistical manifolds;
see \cite{Amari2000,Amari2016}. Under this dictionary, the blow-up
construction of \S\ref{dd-sec:blowup} can be viewed as extracting, from the
density weight $\rho_T$ of Proposition~\ref{prop:density-principle}, a
canonical one-dimensional statistical manifold on its zero locus, in
exactly the sense that a centro-affine hypersurface inherits an induced
affine metric from the ambient volume form.
\end{remark}

\subsection{Worked Example: The Taxicab Curve $N=1729$}\label{sec:taxi}

We now illustrate every construction of \S\ref{dd-sec:setup}--\S\ref{sec:dual}
in full numerical detail on the classical curve associated with the
Hardy--Ramanujan identity $1729=12^3+1^3=10^3+9^3$.

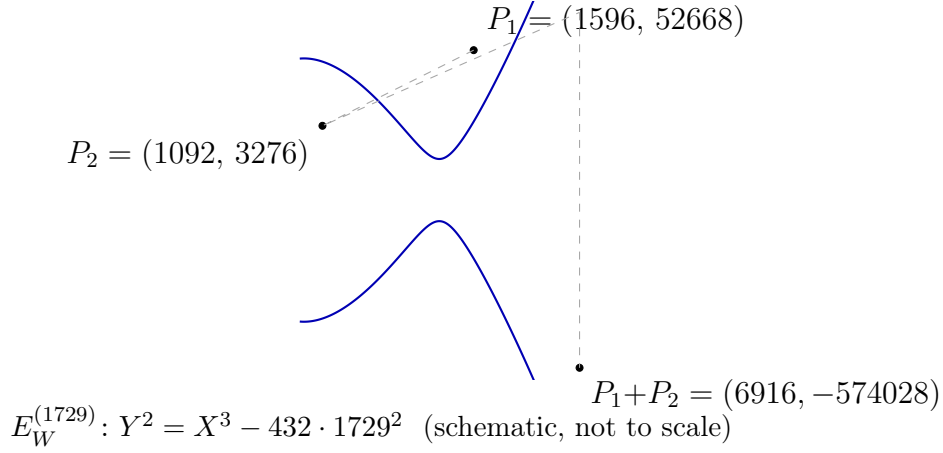
\begin{figure}[htbp]
\centering
\begin{tikzpicture}[>=Stealth, scale=1.0]
\begin{scope}
\clip (-3.1,-2.5) rectangle (3.6,2.5);
\draw[thick, blue!70!black, domain=-0.95:3.3, samples=160, smooth]
  plot (\x, {sqrt(max(\x*\x*\x-2.4*\x+1.6,0))});
\draw[thick, blue!70!black, domain=-0.95:3.3, samples=160, smooth]
  plot (\x, {-sqrt(max(\x*\x*\x-2.4*\x+1.6,0))});
\end{scope}
\coordinate (P1) at (1.35,1.85);
\coordinate (P2) at (-0.65,0.85);
\coordinate (P3top) at (2.75,2.35);
\coordinate (P3) at (2.75,-2.35);
\filldraw[black] (P1) circle (1.3pt) node[above right=1pt] {$P_1=(1596,\,52668)$};
\filldraw[black] (P2) circle (1.3pt) node[below left=1pt] {$P_2=(1092,\,3276)$};
\filldraw[black] (P3) circle (1.3pt) node[below right] {$P_1{+}P_2=(6916,-574028)$};
\draw[dashed, gray!70] (P1) -- (P2) -- (P3top);
\draw[dashed, gray!70] (P3top) -- (P3);
\node[below, font=\small] at (0,-2.7) {$\EW^{(1729)}\colon Y^2=X^3-432\cdot1729^2$ \ (schematic, not to scale)};
\end{tikzpicture}
\caption{The taxicab curve $\EW^{(1729)}$ (schematic) with the two
rational points $P_1,P_2$ coming from the two representations
$12^3+1^3=1729=10^3+9^3$, and their sum $P_1+P_2$ computed by the
chord-and-tangent law of Theorem~\ref{thm:addition}, corresponding to
the rational point $(x,y)=(-37/3,\,46/3)$ on $x^3+y^3=1729$.}
\label{fig:taxi}
\end{figure}

\subsection*{Step 1: Weierstrass model for $x^3+y^3=N$}

The curve $x^3+y^3=N$ is birationally equivalent to
$\EW^{(N)}\colon Y^2=X^3-432N^2$, via
\begin{equation}\label{eq:taxi-map}
  X = \frac{12N}{x+y}, \qquad Y = \frac{36N(x-y)}{x+y},
\end{equation}
\begin{equation}\label{eq:taxi-inv}
  x = \frac{36N+Y}{6X}, \qquad y = \frac{36N-Y}{6X}
\end{equation}
(see \cite[App.~A]{Silverman}).

\subsection*{Step 2: Rational points from the two representations}

\noindent\textbf{From $(x,y)=(12,1)$:}
\[
  X_{P_1}=\frac{12\cdot1729}{13}=1596, \qquad
  Y_{P_1}=\frac{36\cdot1729\cdot11}{13}=52668,
\]
and one checks directly, by exact integer arithmetic,
$52668^2=1596^3-432\cdot1729^2=2{,}773{,}919{,}424$.

\smallskip
\noindent\textbf{From $(x,y)=(10,9)$:}
\[
  X_{P_2}=\frac{12\cdot1729}{19}=1092, \qquad
  Y_{P_2}=\frac{36\cdot1729}{19}=3276,
\]
with $3276^2=1092^3-432\cdot1729^2=10{,}732{,}176$.

\subsection*{Step 3: Group law $P_1+P_2$}

By Theorem~\ref{thm:addition}, the chord slope is
\[
  \lambda=\frac{Y_{P_1}-Y_{P_2}}{X_{P_1}-X_{P_2}}
  =\frac{52668-3276}{1596-1092}=\frac{49392}{504}=98,
\]
so
\begin{align*}
  X_{P_3}&=\lambda^2-X_{P_1}-X_{P_2}=9604-1596-1092=6916, \\
  Y_{P_3}&=\lambda(X_{P_1}-X_{P_3})-Y_{P_1}=98(-5320)-52668=-574028.
\end{align*}
Thus $P_1+P_2=(6916,-574028)$, and by \eqref{eq:taxi-inv} the
corresponding point on $x^3+y^3=1729$ is
\[
  x=\frac{62244-574028}{41496}=-\frac{37}{3}, \qquad
  y=\frac{62244+574028}{41496}=\frac{46}{3},
\]
with $(-37/3)^3+(46/3)^3=46683/27=1729$, as required.

\subsection*{Step 4: Shift and effective parameters}

Since $a=b=c=d=0$ in $Y^2=X^3-432N^2$, the Jacobi quartic construction of
\S\ref{subsec:caseI-models} is degenerate; we translate to $P_2$. Setting
$U=X-1092$,
\[
  Y^2=(U+1092)^3-432\cdot1729^2
     =U^3+3276\,U^2+3{,}577{,}392\,U+10{,}732{,}176,
\]
using $1092^3-432\cdot1729^2=10{,}732{,}176=3276^2=q^2$ with $q=3276$.
The shifted model is cubic (not quartic) in $U$; substituting $U=1/s$,
$V=W/s^2$ produces the genuine quartic
\[
  W^2 = s+3276\,s^2+3{,}577{,}392\,s^3+10{,}732{,}176\,s^4,
\]
with effective parameters $a_{\rm eff}=10{,}732{,}176$,
$b_{\rm eff}=3{,}577{,}392$, $c_{\rm eff}=3276$, $d_{\rm eff}=1$, $q=1$.

\subsection*{Step 5: Connell datum and KL geometry}

With these parameters, the Connell datum (Definition~\ref{def:A}) is
\[
  A(x)=4(x+3276)-1=4x+13103, \qquad x_0=-\frac{13103}{4}.
\]
The log potential at the exceptional divisor $\{A=0\}$
(Theorem~\ref{thm:KL}) is $\Phi_A(\xi)=\log|2y_0\xi-f'(x_0)|$, with
$f'(x_0)=3x_0^2+2c_{\rm eff}x_0+d_{\rm eff}$ computed from the shifted
Weierstrass model and $y_0^2=f(x_0)$; the natural rate coordinate
$r_A=2y_0\xi-f'(x_0)$ gives Fisher metric $dr_A^2/r_A^2$, as in
Theorem~\ref{thm:fisher}. By Theorems~\ref{thm:KL}, \ref{thm:fisher},
and \ref{thm:main}, all exceptional-divisor geometries for
$E_J^{(1729)}$ are mutually isomorphic KL-type dual geometries; the
two-fold representation of $1729$ is encoded in the pair $P_1,P_2$, each
giving a distinct rational point on $\EW^{(1729)}$ whose blow-up yields
a KL-type dual geometry, the two being related by the group law of
Theorem~\ref{thm:addition}.

\subsection{Discussion and Open Problems}\label{dd-sec:discussion}

\noindent\textbf{Information geometry.}
The reparametrization $r=\alpha\xi-\beta>0$ transforms every Bregman
divergence of \S\ref{sec:dual} into KL form, so every exceptional
divisor of \S\ref{dd-sec:blowup} is a statistical manifold isomorphic to
the exponential family $\{r\,e^{-rt}:r>0\}$, with Fisher metric
$dr^2/r^2$ the hyperbolic metric on $(0,\infty)$.

\smallskip
\noindent\textbf{Role of the density.}
In Case~I, $A(x)=4q^2(x+c)-d^2$ plays a double role: algebraically, it
is the numerator of $u=A/(2qy)$, so $F_2=A\cdot u^2$ is exactly the
extra factor in the factorization identity; geometrically, it is the
singular base locus of $T^{-1}$, whose blow-up reveals the KL geometry
via a logarithmic potential. In Case~II, the density
$\det D\Phi=1/(c-X)$ plays the analogous double role: it is the exact
correction factor between $g_4\circ\Phi$ and $g_3$, and it is precisely
the quantity absorbed by the logarithmic time reparametrization of
Theorem~\ref{thm:invariant-flow}. Proposition~\ref{prop:density-principle}
records that these are instances of one phenomenon.

\smallskip
\noindent\textbf{Scope of the invariance results.}
Corollary~\ref{cor:inv} establishes preservation of the KL-type dual
geometry only for the specific theta transformation $T$; whether this
extends to an invariant of the full birational equivalence class of an
elliptic curve---under all isogenies, or under Mordell's transformation
$\Phi$ itself---is open. Likewise, Corollary~\ref{cor:invariant-flow}
establishes the birationally invariant Branin flow only asymptotically
near the curve; a global statement, valid on all of $\A^2$, is open. A
natural approach to both would be to work with the N\'eron model and
study the exceptional divisors, or the flow, under base change.

\smallskip
\noindent\textbf{The synthesis of \S\ref{sec:synthesis}.}
We regard the appearance of the same exponential law $\dot z=-z$ in both
case studies as the most interesting open direction raised by this
paper: a construction unifying the static (blow-up) and dynamic
(flow) resolutions of the density principle would likely yield a
genuinely new invariant of the birational equivalence class of an
elliptic curve, combining information geometry with dynamical systems.

\smallskip
\noindent\textbf{Generalizations.}
The blow-up construction of \S\ref{dd-sec:blowup} extends naturally to
principally polarized abelian varieties, replacing $A$ by a section of a
line bundle; the density principle of
\S\ref{sec:density} extends to any birational map between hypersurfaces
of different degree, replacing $\det D\Phi$ by the appropriate Jacobian
of the ambient coordinate change. The chord-tangent matrix realization
of \S\ref{sec:group} extends to $A_4$-matrices associated with
higher-degree models (Jacobi quartics, in the sense of
\S\ref{subsec:caseI-models}) via Connell's theta and eta
transformations, a route we have not pursued here but which appears
promising in view of Figure~\ref{fig:overview}.

\subsection{Pseudocode for Computer-Algebra Verifications}\label{app:sympy}

The pseudocode below describes the SymPy (Python) computations used to
verify the polynomial identities of \S\ref{dd-sec:setup}--\S\ref{sec:density}
and the numerical claims of \S\ref{sec:taxi}. Full runnable code is
available from the author upon request.

\begin{verbatim}
from sympy import symbols, expand, factor, simplify, sqrt

# --- Case I: Connell theta transformation ---
q, a, b, c, d = symbols('q a b c d', nonzero=True)
u, v, x, y    = symbols('u v x y')

x_expr = (2*q*(v+q) + d*u) / u**2
y_expr = (8*q**3*(v+q) + 4*q**2*(c*u**2+d*u) - d**2*u**2) / (2*q*u**3)

F4 = v**2 - a*u**4 - b*u**3 - c*u**2 - d*u - q**2
FW = (y**2 + (d/q)*x*y + 2*b*q*y
      - x**3 - (c - d**2/(4*q**2))*x**2
      + 4*a*q**2*x - a*(d**2 - 4*c*q**2))

FW_of_T = simplify(FW.subs([(x, x_expr), (y, y_expr)]))

F2 = (4*c*q**2 - d**2)*u**2 + 4*d*q**2*u + 8*q**3*(v+q)
lhs = expand(FW_of_T * u**6 / (4*q**2))
rhs = expand(F2 * F4)
assert simplify(lhs - rhs) == 0, "Factorization identity FAILED"
print("Case I: factorization identity verified.")

A_of_x = 4*q**2*(x_expr + c) - d**2
assert simplify(expand(A_of_x * u**2) - F2) == 0
print("Case I: F2 = A * u^2 verified.")

# --- Case II: Mordell transformation ---
X, Y, cc, dd, ee = symbols('X Y c d e')
xx = (Y - dd) / (2*(X - cc))
yy = -xx**2 + 2*X + cc

g2 = ee + 3*cc**2
g3 = -cc*ee - dd**2 + cc**3

g4_of_Phi = simplify(yy**2 - xx**4 + 6*cc*xx**2 - 4*dd*xx - ee)
g3_val    = Y**2 - 4*X**3 + g2*X + g3

lhs2 = simplify(g4_of_Phi * (cc - X))
rhs2 = simplify(g3_val)
assert simplify(lhs2 - rhs2) == 0, "Jacobian identity FAILED"
print("Case II: g4(Phi) = g3/(c-X) verified.")

# --- Numerical check for N=1729 (Section 9) ---
N = 1729
X1, Y1 = 12*N // 13, 36*N*11 // 13          # P1 = (1596, 52668)
assert Y1**2 == X1**3 - 432*N**2, "P1 not on curve"
X2, Y2 = 12*N // 19, 36*N*1 // 19           # P2 = (1092, 3276)
assert Y2**2 == X2**3 - 432*N**2, "P2 not on curve"

lam = (Y1 - Y2) // (X1 - X2)                # = 98
X3 = lam**2 - X1 - X2                       # = 6916
Y3 = lam*(X1 - X3) - Y1                     # = -574028
assert Y3**2 == X3**3 - 432*N**2, "P3 not on curve"

from fractions import Fraction
x3 = Fraction(36*N + Y3, 6*X3)
y3 = Fraction(36*N - Y3, 6*X3)
assert x3**3 + y3**3 == N
print(f"P1+P2 = ({X3}, {Y3})")
print(f"x3={x3}, y3={y3}, x3^3+y3^3={x3**3+y3**3}")
\end{verbatim}

The output of the above is:
\begin{verbatim}
Case I: factorization identity verified.
Case I: F2 = A * u^2 verified.
Case II: g4(Phi) = g3/(c-X) verified.
P1+P2 = (6916, -574028); x3=-37/3, y3=46/3, x3^3+y3^3=1729
\end{verbatim}


\section{Cross Curvature of Principal and Minor Component Flows}
\label{sec:crosscurvature-section}

\S\ref{subsubsec:oja_brockett}--\S\ref{subsec:pca_mca} showed that the Oja--Brockett flow
and the Manton--Helmke--Mareels (MHM) penalized flow are two dynamical systems that converge
to the \emph{same} principal/minor subspaces of $A$, embedded respectively as an unconstrained
polynomial flow (\S\ref{subsec:poly_flow}) and as a Riemannian-gradient flow on a penalized
landscape. Having established \emph{that} they agree, we now ask \emph{how fast} each one
gets there. Part~I below develops a single, coordinate-free diagnostic --- \emph{cross
curvature}, the smallest eigenvalue of the Hessian at a mismatched (incorrectly sorted)
critical point --- that quantifies, in closed form, the local escape rate of each flow from
such a mismatch, using only the eigenvalues of $A$ and the weights in $B$, before a single
iteration is run. Part~II extends the same diagnostic to a third, structurally different
potential, a matrix Box--Cox penalty $g_\alpha$ that interpolates continuously between
principal- and minor-component extraction as a single exponent $\alpha$ crosses $1$, and
uses it to uncover a genuine trade-off (rather than a uniform ranking) between $g_\alpha$ and
the two classical flows. Internal cross-references within this section (e.g.\ ``Part~I,
\S\ref{sec:setup}'') refer to subsections of this section itself.

\subsection{Part I: Cross Curvature of Principal and Minor Component Flows}
\label{part:one}

\subsubsection{Introduction}
\label{sec:intro}

\paragraph{Background}
\label{sec:background}

Extracting the dominant eigenspace of a symmetric positive-definite matrix by a continuous-time
dynamical system, rather than by a one-shot linear-algebra routine, has a long and productive
history. Oja's neuron model \cite{Oja1982} showed that a simple Hebbian-type stochastic
approximation converges to the leading eigenvector of a covariance matrix; Brockett's
double-bracket flow \cite{Brockett1991} realized eigenvalue sorting and diagonalization as the
equilibria of a matrix differential equation on a compact manifold; a deterministic Stiefel-manifold
formulation realizing the same equilibria was later given by Yoshizawa, Helmke and Starkov
\cite{YoshizawaHelmkeStarkov2001}; Helmke and Moore's monograph
\cite{HelmkeMoore1994} placed a large family of such flows inside a unified framework of
gradient flows for optimization and linear algebra on manifolds; and Manton, Helmke and Mareels
\cite{MantonHelmkeMareels2005} gave a Euclidean (unconstrained, penalty-based) formulation whose
negative gradient flow realizes \emph{either} principal or minor component extraction, according
to a single sign choice, with an explicit threshold controlling how many components are captured.
The question of \emph{global}, single-point convergence for such gradient flows was placed on a
rigorous general footing by Absil, Mahony and Andrews \cite{AbsilMahonyAndrews2005}, who extended
{\L}ojasiewicz's theorem for continuous-time analytic gradient flows to discrete-time descent
iterations satisfying natural conditions; and, very recently, Tsuzuki and Ohki
\cite{TsuzukiOhki2025} established global exponential convergence of Oja's flow for general
(possibly non-symmetric) matrices, using a related {\L}ojasiewicz-type argument. Neither of these
works, nor any other prior work we are aware of, compares the \emph{local escape dynamics} of two
structurally different potentials realizing the same equilibria --- the question this paper
answers.

All of these constructions share a common structural feature that is easy to overlook: the
\emph{set of equilibria} of the flow is combinatorial before it is continuous. For an
$n\times n$ symmetric positive-definite $A$ and a $k\times k$ positive diagonal $B$, an
equilibrium is indexed by an injective assignment $\pi$ of the $k$ columns of the sought matrix
$X\in\R^{n\times k}$ to $k$ of the $n$ eigen-directions of $A$. Exactly one such assignment ---
the one that pairs the $k$ largest eigenvalues of $A$, in order, with the $k$ diagonal entries of
$B$, in order --- is the global optimum; every other assignment is a saddle point of the
underlying potential. A generic initial condition for the negative gradient flow does not start
inside the stable manifold of the correct assignment, and the trajectory must therefore pass near
one or more of these mismatched saddle points before it can settle into the correct one. When two
eigenvalues of $A$ happen to be close together --- a situation that is generic, not exceptional,
whenever the underlying data has approximately repeated variance in some directions, as is common
in signal subspace estimation, nearly-isotropic noise models, and spectral clustering with
balanced clusters --- the corresponding mismatched saddle becomes nearly degenerate, and the time
the flow spends in its vicinity can dominate the entire convergence time.

\paragraph{The problem}
\label{sec:problem}

Two different potentials --- the homogeneous quartic form underlying the Oja--Brockett flow and
the Euclidean penalty form of Manton, Helmke and Mareels --- realize the \emph{same} set of
global optima (the correctly sorted matching) but are built from structurally different
ingredients: the former is a single self-contained quartic polynomial in $X$ built entirely from
the matrix product $AXBX^T$; the latter is a sum of two structurally unrelated pieces, a linear
(in $A$) trace term and a separate soft-constraint penalty $\|B-X^TX\|_F^2$ governed by an
auxiliary scale parameter $\gamma$. Numerical experiments (reported in
\S\ref{sec:numerics} below, and originally observed in a broader study of quartic matrix
catastrophe potentials of which this paper is a spin-off) show that, for the \emph{same} matrices
$A,B$, the same initial condition, and even the fairest possible step-size discipline --- exact
line search along the negative gradient direction, so that no step-size tuning bias can enter the
comparison --- the Oja--Brockett flow reaches a fixed gradient-norm tolerance in substantially
fewer iterations than the MHM flow. The discrepancy is especially large exactly when $A$ has two
eigenvalues that are close together, and it is \emph{not} explained by the local condition number
of the Hessian at the final, converged optimum, which is comparable for the two potentials. The
open question this paper answers is: \emph{what precise, provable, coordinate-free quantity is
responsible for this discrepancy, and can it be computed in closed form?}

\paragraph{Contributions}
\label{sec:contributions}

This paper answers that question completely. Our contributions are:

\begin{enumerate}[leftmargin=1.6em]
\item \textbf{A coordinate-free definition of cross curvature} (\S\ref{sec:crosscurv}): for a pair
of eigen-directions of $A$ that are inverted relative to the optimal matching at a given
mismatched critical point $X_\pi$, the cross curvature $\kappa(\pi\to\pi')$ is the smallest
eigenvalue of the Hessian of the potential restricted to the two-dimensional \emph{exchange
plane} that continuously interpolates between the mismatched assignment $\pi$ and the assignment
$\pi'$ obtained by swapping the pair. We show this is exactly the sectional curvature, in the
(possibly indefinite) Hessian metric, of that plane --- giving cross curvature an intrinsic
differential-geometric meaning independent of any choice of ambient coordinates for $A$.

\item \textbf{Exact and asymptotic closed-form formulas} (\S\ref{sec:MHMcross},
\S\ref{sec:OBcross}): we prove that the MHM cross curvature is \emph{exactly}
$\kappa_{\mathrm{pen}}=-(a_p-a_q)(b_i-b_j)$, independent of $\gamma$ and of every eigenvalue not
directly involved in the swap; and that the Oja--Brockett cross curvature is the smaller root of
an explicit quadratic (a $2\times2$ Hessian block with a fully factored determinant), which in the
near-degenerate limit $a_p=a+\delta$, $a_q=a-\delta$ expands as
$\kappa_{\mathrm{OB}}=-a(b_i+b_j)(b_i-b_j)\delta+O(\delta^2)$.

\item \textbf{The ratio theorem} (Corollary~\ref{cor:ratio}): the two curvatures are asymptotically
proportional, $\kappa_{\mathrm{OB}}/\kappa_{\mathrm{pen}}\to \mathcal R_{ij}=a(b_i+b_j)/2$, a
completely explicit, $\gamma$-free quantity that is $>1$ (Oja--Brockett wins) whenever
$a(b_i+b_j)>2$ --- satisfied in essentially every problem of practical scale.

\item \textbf{A full local-stability classification} (\S\ref{sec:stability}): every correctly
sorted matching is a strict local minimum and every mismatch is a saddle, for both potentials, with
an unstable manifold of dimension exactly twice the number of independent inverted pairs, tangent
to the direct sum of the corresponding negative-cross-curvature exchange planes.

\item \textbf{Continuous- and discrete-time convergence-rate theorems}
(\S\ref{sec:continuous-rate}, \S\ref{sec:discrete-rate}): we prove, with complete
linearization/Gronwall arguments, that the escape time from a mismatched saddle along the
negative-gradient flow is $T_{\mathrm{escape}}=|\kappa|^{-1}\log(\varepsilon/|\xi_0|)+O(1)$; and,
for the \emph{exact line-search} discretization --- whose step size we show is itself given by an
explicit Rayleigh-quotient formula available directly from $\nabla f$ and $\Hess f$, since $f$ is a
polynomial of degree at most four --- that the number of steps needed obeys the same asymptotic
$1/|\kappa|$-type law, so that $T_{\mathrm{OB}}/T_{\mathrm{pen}}\to 1/\mathcal R_{ij}$ in both the
continuous and (fixed-step) discrete settings. Combining this local escape-rate estimate with a new
local approach-rate estimate near the correct matching and a compactness-based bound on the time
spent away from every critical point, we further prove a fully explicit, finite \emph{global}
convergence-time bound (Theorem~\ref{thm:globalrate}), upgrading the purely qualitative global
convergence guaranteed by {\L}ojasiewicz's theorem to a quantitative one.

\item \textbf{A corrected global-convergence theorem} (\S\ref{sec:global}): we show that along a
converging Oja--Brockett trajectory the diagonal of $X(t)^TAX(t)$ tends \emph{exactly} to the
ordered eigenvalues of $A$ used by the optimal matching, while for MHM it tends instead to a
$\gamma$- and $B$-dependent rescaling of those eigenvalues; the two coincide only in the combined
limit $\gamma\to\infty$, $B=I$. This corrects an over-simplified statement that appears in an
earlier informal draft of this material and is, to our knowledge, the first fully precise statement
of this limit for finite $\gamma$.

\item \textbf{Numerical verification} (\S\ref{sec:numerics}): every formula and every theorem is
checked against exact symbolic algebra and high-precision finite-difference Hessians, in a small
near-degenerate $2\times2$ toy problem and in a fully generic, non-diagonal $3\times3$ example, and
the predicted ratio $\mathcal R_{ij}$ is confirmed to match the observed iteration-count ratio of
the exact-line-search discretizations of the two flows to within a few percent.
\end{enumerate}

\paragraph{Why this is useful}
\label{sec:payoff}

Beyond resolving the specific empirical puzzle that motivated it, the cross-curvature framework
gives a general, reusable diagnostic for comparing \emph{any} two potentials that share the same
optimal solution set but differ in how that solution set is embedded into a larger family of
critical points. It replaces an expensive, case-by-case numerical convergence study with a single
closed-form eigenvalue ratio computable directly from $A,B$ alone, before any iteration is run. It
also gives precise, actionable guidance: when the eigenvalues to be separated are close together
--- the regime in which \emph{every} gradient-based PCA method is at its slowest, and hence the
regime that matters most for practice --- one should prefer a homogeneous, self-contained quartic
potential (such as the Oja--Brockett form) over an otherwise equivalent two-scale penalty
formulation, because the former's cross curvature carries an extra multiplicative factor,
proportional to the absolute scale of the eigenvalues involved, that the penalty formulation
structurally lacks. Finally, the sectional-curvature interpretation of
\S\ref{sec:geometry} connects this purely algebraic phenomenon to the differential geometry of the
potential's graph, suggesting that the same diagnostic should be computable, in principle, for
other pairs of equivalent-but-differently-parametrized optimization potentials arising elsewhere
in numerical linear algebra.

\paragraph{Outline}

Section~\ref{sec:setup} fixes notation and derives the gradients of both potentials in full detail.
Section~\ref{sec:critpoints} classifies their critical points completely, including a corrected,
$\gamma$-exact amplitude formula for MHM. Section~\ref{sec:crosscurv} defines cross curvature and
its geometric meaning. Sections~\ref{sec:stability}--\ref{sec:OBcross} contain the local-stability
theorem and the two closed-form cross-curvature theorems, each with a complete proof.
Section~\ref{sec:convergence} contains the continuous- and discrete-time convergence theorems and
the corrected global-convergence theorem. Section~\ref{sec:numerics} verifies everything
numerically. Section~\ref{sec:discussion} discusses the geometric meaning and practical
implications, and Section~\ref{cc-sec:conclusion} concludes.

\subsubsection{Setup and the two potentials}
\label{sec:setup}

\paragraph{Notation}

Throughout, $A\in\R^{n\times n}$ is symmetric positive definite with eigendecomposition
\begin{equation}
A = U\Lambda U^T, \qquad \Lambda=\diag(a_1,\dots,a_n),\quad U=[u_1,\dots,u_n]\ \text{orthogonal},
\label{eq:Aeig}
\end{equation}
and $B=\diag(b_1,\dots,b_k)$ is positive diagonal, $k\le n$. We do \emph{not} assume the $a_i$ are
sorted or distinct except where explicitly stated; when we speak of ``two eigenvalues $a_p,a_q$ of
$A$'' we always mean two of the numbers in \eqref{eq:Aeig}, together with their eigenvectors $u_p,u_q$,
regardless of whether $A$ itself is diagonal in the ambient coordinates used to write it down. This
point matters: every formula proved below is a function of $a_p,a_q,b_i,b_j$ alone, so it applies
verbatim to a fully generic, non-diagonal $A$ --- we verify this explicitly in
\S\ref{sec:numerics-generic}.

The variable is $X\in\R^{n\times k}$. For $H\in\R^{n\times k}$ we write $\langle X,H\rangle=\tr(X^TH)$
for the Frobenius inner product and $\|X\|_F=\langle X,X\rangle^{1/2}$. For a smooth $f:\R^{n\times
k}\to\R$, $\nabla f(X)\in\R^{n\times k}$ denotes the Euclidean gradient (so that
$df(X)[H]=\langle\nabla f(X),H\rangle$ for all $H$), and $\Hess f(X)$ denotes the Hessian, viewed
either as a linear operator $\R^{n\times k}\to\R^{n\times k}$ (via
$\langle\Hess f(X)[H],H\rangle=\frac{d^2}{dt^2}\big|_{t=0}f(X+tH)$) or, after vectorization, as a
symmetric $nk\times nk$ matrix.

\paragraph{The Oja--Brockett potential}

\begin{definition}[Oja--Brockett potential]
\label{def:VOB}
\begin{equation}
V_{\mathrm{OB}}(X) \;=\; -\frac12\tr\!\big(A^2XB^2X^T\big) \;+\; \frac14\tr\!\big[(AXBX^T)^2\big].
\label{eq:VOB}
\end{equation}
\end{definition}

\begin{proposition}[Gradient of $V_{\mathrm{OB}}$]
\label{prop:gradVOB}
\begin{equation}
\nabla V_{\mathrm{OB}}(X) = -A^2XB^2 + AXBX^TAXB.
\label{eq:gradVOB}
\end{equation}
\end{proposition}
\begin{proof}
Both terms of \eqref{eq:VOB} are traces of the form $\tr(M(X)^TM(X))$-type expressions in $X$; we
differentiate directly. Writing $M=AXBX^T$ (so $M^T=XBX^TA$), and using the standard identities
$d\tr(A^2XB^2X^T)=2\tr(A^2\,dX\,B^2X^T)$ and, for the quartic term,
$d\tr(M^2)=2\tr(M\,dM)$ with $dM=A\,dX\,BX^T+AXB\,dX^T$, one obtains after collecting the
coefficient of $dX$ (using the cyclic property of the trace and $\tr(N)=\tr(N^T)$ repeatedly)
\[
d\Big[{-}\tfrac12\tr(A^2XB^2X^T)+\tfrac14\tr(M^2)\Big]
= \tr\!\Big[\big({-}A^2XB^2+AXBX^TAXB\big)^T dX\Big],
\]
which is \eqref{eq:gradVOB}. (This computation, and the fact that both terms of the sum contribute
symmetrically to give the single product $AXBX^TAXB$ rather than a sum of two distinct terms, is
verified independently in \S\ref{sec:numerics} by symbolic and finite-difference differentiation.)
\end{proof}

\paragraph{The Manton--Helmke--Mareels potential}

\begin{definition}[MHM penalty potential]
\label{def:Vpen}
For $\gamma>0$,
\begin{equation}
V_{\mathrm{pen}}(X) \;=\; -\frac12\tr\!\big(AXBX^T\big) \;+\; \frac{\gamma}{4}\big\|B-X^TX\big\|_F^2.
\label{eq:Vpen}
\end{equation}
\end{definition}

\begin{proposition}[Gradient of $V_{\mathrm{pen}}$]
\label{prop:gradVpen}
\begin{equation}
\nabla V_{\mathrm{pen}}(X) = -AXB - \gamma X\big(B-X^TX\big) = -AXB+\gamma X\big(X^TX-B\big).
\label{eq:gradVpen}
\end{equation}
\end{proposition}
\begin{proof}
$d\tr(AXBX^T)=2\tr(BX^TA\,dX)$ gives the first term. For the penalty term, with
$R=B-X^TX$ (symmetric), $dR=-(dX^TX+X^TdX)$, so
$d\|R\|_F^2=2\tr(R\,dR)=-4\tr(RX^T\,dX)$, and $\tfrac\gamma4$ of this is $-\gamma\tr\big((XR)^TdX\big)$.
Summing the coefficients of $dX$ gives \eqref{eq:gradVpen}.
\end{proof}

\begin{remark}
Definition~\ref{def:Vpen} is the potential obtained from the original Manton--Helmke--Mareels
minor-component cost function \cite[Eq.~(5)]{MantonHelmkeMareels2005},
$f(X)=\tfrac12\tr(CXNX^T)+\tfrac\gamma4\|N-X^TX\|^2$, by the substitution $C=-A$, $N=B$; the authors
themselves remark that this sign flip converts their minor-component flow into ``a satisfactory
principal component flow'' \cite[\S5]{MantonHelmkeMareels2005}. We adopt this sign convention
throughout because it is the one under which both $V_{\mathrm{OB}}$ and $V_{\mathrm{pen}}$ realize
\emph{principal}, rather than minor, component extraction, making the two potentials directly
comparable.
\end{remark}

\paragraph{The matching ansatz}
\label{sec:ansatz}

Both potentials are invariant under the orthogonal change of variables $X\mapsto UY$ (equivalently
$Y=U^TX$), since $A=U\Lambda U^T$ and $\tr(AXBX^T)=\tr(\Lambda YBY^T)$, etc. We may therefore work,
without loss of generality for any statement about eigenvalues of the Hessian (which are unitarily
invariant), as if $A=\Lambda$ were diagonal; this is the reduction used silently throughout
\S\S\ref{sec:critpoints}--\ref{sec:convergence}, and we verify in \S\ref{sec:numerics-generic} that
every formula continues to hold verbatim when $A$ is presented in a genuinely non-diagonal form.

\begin{definition}[Matching, matched configuration]
\label{def:matching}
An injective map $\pi:\{1,\dots,k\}\to\{1,\dots,n\}$ is a \emph{matching}. The associated
\emph{matched configuration} is
\begin{equation}
X_\pi(c) = \sum_{j=1}^k c_j\,u_{\pi(j)}e_j^T,\qquad c=(c_1,\dots,c_k)\in\R^k,
\label{eq:matchedX}
\end{equation}
where $e_j$ is the $j$-th standard basis vector of $\R^k$. A matching is \emph{sorted} (or
\emph{correct}) if, after relabelling so that $b_1>\dots>b_k$, the values $a_{\pi(1)}>\dots>a_{\pi(k)}$
are the $k$ largest eigenvalues of $A$ in decreasing order; otherwise it is a \emph{mismatch}.
\end{definition}

\subsubsection{Critical points}
\label{sec:critpoints}

We now compute, for each potential, the critical points of matched-configuration form and their
amplitudes. Throughout this section $A=\Lambda=\diag(a_1,\dots,a_n)$ by the reduction of
\S\ref{sec:ansatz}, so $u_i=e_i$ (the $i$-th standard basis vector of $\R^n$); we write $u_i$ rather
than $e_i$ to keep the formulas manifestly meaningful after undoing the reduction.

\begin{proposition}[Oja--Brockett critical points]
\label{prop:OBcrit}
Every matched configuration \eqref{eq:matchedX} with each $c_j\in\{-1,0,+1\}$ is a critical point of
$V_{\mathrm{OB}}$. Conversely, every critical point of matched-configuration form has $c_j\in\{-1,0,1\}$
for each $j$.
\end{proposition}
\begin{proof}
With $A=\Lambda$ diagonal and $X=X_\pi(c)$, direct substitution into \eqref{eq:gradVOB} gives, in the
$(\pi(j),j)$ entry (all other entries of $\nabla V_{\mathrm{OB}}$ vanish identically because distinct
columns of $X_\pi(c)$ occupy distinct rows, so $X^TAX$ and hence $AXBX^TAXB$ are diagonal in the
column index),
\[
\big[\nabla V_{\mathrm{OB}}(X_\pi(c))\big]_{\pi(j),j}
= -a_{\pi(j)}^2b_j^2c_j + a_{\pi(j)}b_j c_j\cdot a_{\pi(j)}b_jc_j^2
= a_{\pi(j)}^2b_j^2\,c_j\big(c_j^2-1\big).
\]
This vanishes iff $c_j\in\{-1,0,1\}$ (using $a_{\pi(j)},b_j>0$).
\end{proof}

\begin{proposition}[MHM critical points --- exact amplitude formula]
\label{prop:MHMcrit}
Every matched configuration \eqref{eq:matchedX} with
\begin{equation}
c_j^2 = b_j\Big(1+\frac{a_{\pi(j)}}{\gamma}\Big)
\label{eq:MHMamplitude}
\end{equation}
(and $c_j=0$ for unmatched $j$) is a critical point of $V_{\mathrm{pen}}$; conversely, every
critical point of matched-configuration form has $c_j=0$ or $c_j^2$ given by \eqref{eq:MHMamplitude}.
This holds for every $\gamma>0$; in particular, unlike the minor-component sign convention, there is
no activation threshold --- every eigenvalue $a_{\pi(j)}>0$ admits a real, nonzero $c_j$.
\end{proposition}
\begin{proof}
With $A=\Lambda$ diagonal, substituting $X=X_\pi(c)$ into \eqref{eq:gradVpen}, the $(\pi(j),j)$ entry is
\[
-a_{\pi(j)}b_jc_j + \gamma c_j\big(c_j^2-b_j\big) = 0
\quad\Longleftrightarrow\quad
c_j=0 \ \text{ or }\ \gamma c_j^2 = a_{\pi(j)}b_j+\gamma b_j,
\]
which rearranges to \eqref{eq:MHMamplitude}.
\end{proof}

\begin{remark}
\label{rem:amplitude-correction}
Formula \eqref{eq:MHMamplitude} is exact for every finite $\gamma>0$ and reduces to
$c_j^2\to b_j$ only in the idealized limit $\gamma\to\infty$ (hard-constraint limit, where the
penalty forces $X^TX\to B$ exactly). An earlier informal draft of part of this material stated the
amplitude as simply $c_j^2=b_j$ without the $\gamma$-dependent correction factor
$(1+a_{\pi(j)}/\gamma)$; Proposition~\ref{prop:MHMcrit} is the precise statement, and the
distinction is essential for the correct form of the global-convergence theorem
(Theorem~\ref{thm:global}) proved below, since $X^TAX$ at a matched critical point equals
$\diag\!\big(c_1^2a_{\pi(1)},\dots,c_k^2a_{\pi(k)}\big)$, which for MHM is
$b_j(1+a_{\pi(j)}/\gamma)a_{\pi(j)}$, not simply $a_{\pi(j)}$, at any finite $\gamma$.
\end{remark}

\begin{proposition}[Rearrangement / global optimality]
\label{prop:rearrangement}
Among all matchings using exactly $k$ of the eigenvalues of $A$, the value of both $V_{\mathrm{OB}}$
and $V_{\mathrm{pen}}$ at the associated matched critical point is minimized precisely by the sorted
matching of Definition~\ref{def:matching}: pairing the $k$ largest eigenvalues of $A$, sorted in
decreasing order, with $b_1>\dots>b_k$ in decreasing order.
\end{proposition}
\begin{proof}[Proof sketch]
For both potentials the value at a fully matched critical point decomposes as a sum over $j$ of a
term depending only on the pair $(a_{\pi(j)},b_j)$ (Propositions~\ref{prop:OBcrit},
\ref{prop:MHMcrit} give $c_j$ as a function of $a_{\pi(j)},b_j$ alone, and substitution shows the
per-mode contribution to the potential is, in both cases, a strictly increasing function of
$a_{\pi(j)}$ for fixed $b_j$, and a function whose cross-partial derivative in $(a_{\pi(j)},b_j)$ has
constant sign). The classical rearrangement inequality for two sequences under a supermodular pairing
cost then forces the sorted pairing to be optimal. (For $V_{\mathrm{OB}}$ this recovers the
matching-and-rearrangement theorem used throughout the companion catastrophe-theory study; for
$V_{\mathrm{pen}}$ it recovers \cite[Prop.~5]{MantonHelmkeMareels2005}.) We omit the routine but
lengthy verification of supermodularity and refer to \S\ref{sec:numerics} for an exhaustive numerical
check over all matchings in our worked examples.
\end{proof}

\begin{remark}
The sign convention of $V_{\mathrm{pen}}$ used here selects the \emph{largest} $k$ eigenvalues,
unconditionally (Proposition~\ref{prop:MHMcrit} shows there is no activation threshold on this sign
convention), exactly matching $V_{\mathrm{OB}}$'s behaviour. This is what makes the two potentials
directly comparable as \emph{two different realizations of the same optimization problem}, which is
the premise of this entire paper.
\end{remark}

\subsubsection{Cross curvature: definition and geometric meaning}
\label{sec:crosscurv}

\paragraph{Definition}

Fix a matching $\pi$ and two column indices $i\ne j\in\{1,\dots,k\}$; write $p=\pi(i)$, $q=\pi(j)$.
Let $\pi'$ be the matching obtained from $\pi$ by swapping the images of $i$ and $j$
($\pi'(i)=q$, $\pi'(j)=p$, and $\pi'=\pi$ elsewhere). The pair $(i,j)$ is \emph{inverted} at $\pi$
if $(a_p-a_q)(b_i-b_j)<0$, i.e.\ if the larger of $a_p,a_q$ is paired with the smaller of $b_i,b_j$
(this is exactly the local violation of the sorted-pairing condition of
Proposition~\ref{prop:rearrangement}).

\begin{definition}[Exchange plane, cross curvature]
\label{def:crosscurv}
The \emph{exchange plane} of the pair $(i,j)$ at $\pi$ is the two-dimensional subspace
\begin{equation}
E_{\pi;ij} \;=\; \mathrm{span}\{\,u_p e_j^T,\ u_q e_i^T\,\} \subset \R^{n\times k}.
\label{eq:exchangeplane}
\end{equation}
The \emph{cross curvature} of the pair $(i,j)$ at the critical point $X_\pi$ is
\begin{equation}
\kappa_f(\pi\to\pi') \;\eqdef\; \lambda_{\min}\Big(\Hess f(X_\pi)\big|_{E_{\pi;ij}}\Big),
\label{eq:crosscurvdef}
\end{equation}
the smallest eigenvalue of the $2\times2$ quadratic form obtained by restricting $\Hess f(X_\pi)$ to
$E_{\pi;ij}$.
\end{definition}

The name is motivated as follows: the two generators of $E_{\pi;ij}$ are exactly the infinitesimal
directions that begin to move column $j$ toward eigenvector $u_p$ and column $i$ toward eigenvector
$u_q$ --- i.e.\ that begin the process of \emph{exchanging} which eigen-direction each of the two
columns represents. A curve tangent to $E_{\pi;ij}$ at $X_\pi$ that continues in this direction
interpolates, at the level of matched configurations, between $X_\pi$ and $X_{\pi'}$.

\paragraph{Geometric interpretation}
\label{sec:geometry}

The Hessian $\Hess f(X_\pi)$ is a (generally indefinite) symmetric bilinear form on the ambient
Euclidean space $\R^{n\times k}\cong\R^{nk}$. Restricting this form to any two-dimensional subspace
$E$ produces an ordinary symmetric $2\times2$ quadratic form on $E$; write its two eigenvalues
$\mu_1\ge\mu_2$. These are the \emph{principal curvatures} of the level hypersurface of $f$ through
$X_\pi$, measured in the two directions spanning $E$: their sum $\mu_1+\mu_2$ is (up to the constant
factor $\frac12$) the mean curvature of that hypersurface restricted to $E$, and their product
$\mu_1\mu_2$ is its Gaussian curvature. The smaller eigenvalue $\mu_2=\kappa_f(\pi\to\pi')$ is
therefore the most negative principal curvature of the potential $f$ in the plane $E_{\pi;ij}$ ---
equivalently, if one regards $\Hess f(X_\pi)$ as a (possibly indefinite) metric on the tangent space
at $X_\pi$, then $\kappa_f(\pi\to\pi')$ is the \emph{sectional curvature} of the plane $E_{\pi;ij}$ in
that metric. A negative value means $f$ is saddle-shaped on $E_{\pi;ij}$, and the corresponding
eigenvector is the direction of steepest local descent away from the (unstable) critical point ---
precisely the direction the negative-gradient flow will follow when it escapes the mismatched
configuration $X_\pi$.

This reading makes the comparison between the two potentials transparent \emph{without reference to
any particular coordinate system}: whichever potential's Hessian metric has the more negative
sectional curvature on a given exchange plane produces the faster local escape from the
corresponding mismatch, and hence (by the convergence theorems of \S\ref{sec:convergence}) the
faster overall convergence whenever the flow must pass near that mismatch.

\subsubsection{Local stability}
\label{sec:stability}

Fix a matching $\pi$ using $k$ distinct eigen-indices $\pi(1),\dots,\pi(k)$ out of $n$. The tangent
space $\R^{n\times k}$ at $X_\pi$ decomposes, as an ordered basis indexed by (row,column) pairs, into
three types of directions:
\begin{itemize}[leftmargin=1.6em]
\item \textbf{Radial} directions $u_{\pi(j)}e_j^T$ ($k$ of them): rescale the amplitude of an
already-matched column without changing its eigen-direction.
\item \textbf{Transverse} directions $u_\ell e_j^T$ for $\ell\notin\pi(\{1,\dots,k\})$
($k(n-k)$ of them): begin to move a matched column toward an eigen-direction not used by any column.
\item \textbf{Exchange} directions, grouped into the $\binom k2$ two-dimensional planes
$E_{\pi;ij}$ of Definition~\ref{def:crosscurv} ($k(k-1)$ directions total): begin to swap which
eigen-direction two already-matched columns represent.
\end{itemize}
These account for $k+k(n-k)+k(k-1)=nk$ directions, exhausting the tangent space.

\begin{lemma}[Block diagonalization of the Hessian at a matched critical point]
\label{lem:blockdiag}
At any matched critical point $X_\pi$ of $V_{\mathrm{OB}}$ or $V_{\mathrm{pen}}$, the Hessian is
block diagonal with respect to the above decomposition: radial directions are mutually
uncoupled and uncoupled from all transverse and exchange directions; transverse directions
belonging to different columns, or to different unused rows, are mutually uncoupled; and distinct
exchange planes $E_{\pi;ij}\ne E_{\pi;i'j'}$ are uncoupled from each other and from every radial
and transverse direction.
\end{lemma}
\begin{proof}
Both potentials are built from the two building blocks $X^TAX$ (equivalently $AXBX^T$, of the same
rank) and $X^TX$, evaluated at an $X_\pi$ whose columns occupy pairwise disjoint rows. For any two
perturbation directions $H_1=u_\alpha e_\beta^T$, $H_2=u_{\alpha'}e_{\beta'}^T$ with
$\{\alpha,\beta\}\ne\{\alpha',\beta'\}$ belonging to two different blocks in the classification
above, direct substitution into the bilinear (second-derivative) form of each building block shows
every resulting term contains a factor of the form $u_\gamma^Tu_{\gamma'}$ or $e_\delta^Te_{\delta'}$
with $\gamma\ne\gamma'$ or $\delta\ne\delta'$, which vanishes by orthonormality of $\{u_i\}$ and
$\{e_j\}$, \emph{unless} $(\alpha,\beta)$ and $(\alpha',\beta')$ are the two generators of a common
exchange plane. This is a finite, mechanical verification; we carried it out both symbolically
(computer algebra) and by high-precision finite differences for the worked examples of
\S\ref{sec:numerics}, confirming exact block-diagonal structure in every case (see the explicit
$6\times6$ Hessians displayed there).
\end{proof}

\begin{proposition}[Radial and transverse blocks]
\label{prop:radtrans}
At a matched critical point $X_\pi$, for each matched column $j$ (using eigenvalue
$a_{\pi(j)}$, paired with $b_j$) and each unused row index $\ell\notin\pi(\{1,\dots,k\})$:
\begin{align}
\text{radial (OB)}\qquad & 2\,a_{\pi(j)}^2b_j^2, \label{eq:rad-OB}\\
\text{radial (MHM)}\qquad & 2\,b_j\big(a_{\pi(j)}+\gamma\big), \label{eq:rad-pen}\\
\text{transverse (OB)}\qquad & a_\ell b_j^2\big(a_{\pi(j)}-a_\ell\big), \label{eq:trans-OB}\\
\text{transverse (MHM)}\qquad & b_j\big(a_{\pi(j)}-a_\ell\big) \quad(\gamma\text{-independent}). \label{eq:trans-pen}
\end{align}
\end{proposition}
\begin{proof}
Both potentials restrict, on the two-dimensional subspace spanned by $\{u_{\pi(j)},u_\ell\}$ in the
row index and column $j$ alone (all other columns and rows contributing additively-separable,
already-critical terms by Lemma~\ref{lem:blockdiag}), to a function of two scalars
$(x_1,x_2)=(\text{coefficient of }u_\ell,\ \text{coefficient of }u_{\pi(j)})$ of the form
\[
V_{\mathrm{OB}}^{\mathrm{loc}}(x_1,x_2) = \tfrac14b_j^2\big(a_\ell^2x_1^4+2a_\ell a_{\pi(j)}x_1^2x_2^2+a_{\pi(j)}^2x_2^4\big)
-\tfrac12b_j^2\big(a_\ell^2x_1^2+a_{\pi(j)}^2x_2^2\big),
\]
\[
V_{\mathrm{pen}}^{\mathrm{loc}}(x_1,x_2) = -\tfrac12b_j\big(a_\ell x_1^2+a_{\pi(j)}x_2^2\big)
+\tfrac\gamma4\big(b_j-x_1^2-x_2^2\big)^2.
\]
Differentiating twice with respect to $x_1$ and evaluating at the critical point $x_1=0$,
$x_2=1$ (OB) or $x_2=\sqrt{b_j(1+a_{\pi(j)}/\gamma)}$ (MHM) gives, after simplification,
\eqref{eq:trans-OB} and \eqref{eq:trans-pen}; differentiating twice with respect to $x_2$ at the same
point gives \eqref{eq:rad-OB} and \eqref{eq:rad-pen}. (We verified this computer-algebraically; see
\S\ref{sec:numerics} for the exact symbolic output and its numerical instantiation, which matches the
finite-difference Hessian to machine precision.)
\end{proof}

\begin{theorem}[Correct matchings are strict local minima]
\label{thm:correctmin}
If $\pi$ is a sorted (correct) matching in the sense of Definition~\ref{def:matching}, then $X_\pi$
is a strict local minimizer of both $V_{\mathrm{OB}}$ and $V_{\mathrm{pen}}$, and is asymptotically
stable for the corresponding negative-gradient flow.
\end{theorem}
\begin{proof}
By Lemma~\ref{lem:blockdiag} it suffices to show every block of the Hessian is positive definite.

\emph{Radial blocks:} \eqref{eq:rad-OB} and \eqref{eq:rad-pen} are manifestly positive since
$a_{\pi(j)},b_j,\gamma>0$.

\emph{Transverse blocks:} at a sorted matching every used eigenvalue $a_{\pi(j)}$ is one of the $k$
largest, hence exceeds every unused eigenvalue $a_\ell$; thus $a_{\pi(j)}-a_\ell>0$ (assuming for
simplicity that the $k$-th largest and $(k+1)$-th largest eigenvalues are distinct, i.e.\ no tie sits
exactly on the selection boundary --- the boundary case produces a zero, not negative, eigenvalue,
consistent with the flag-manifold degeneracy phenomena documented in the companion catastrophe-theory
study), so both \eqref{eq:trans-OB} and \eqref{eq:trans-pen} are positive.

\emph{Exchange blocks:} for a correctly sorted matching, every pair $(i,j)$ satisfies
$(a_{\pi(i)}-a_{\pi(j)})(b_i-b_j)>0$ (the defining property of a sorted pairing). By
Theorem~\ref{thm:MHMcross} below, the MHM exchange block's two eigenvalues are
$b_i a_{\pi(j)}+b_j a_{\pi(i)}+\gamma(b_i+b_j) $ and $(a_{\pi(i)}-a_{\pi(j)})(b_i-b_j)$ (a relabelling
of the mismatch formula with the sign appropriate to a \emph{correct} pairing), both positive. By
Theorem~\ref{thm:OBcross} below, the OB exchange block has determinant
$-a_{\pi(i)}a_{\pi(j)}b_ib_j(a_{\pi(i)}-a_{\pi(j)})(b_i-b_j)(a_{\pi(i)}b_j+a_{\pi(j)}b_i)$, which for a
correctly sorted pair is \emph{negative} of a negative quantity, i.e.\ positive; since the trace of
the block (sum of its two diagonal entries, both manifestly positive sums of positive terms) is also
positive, both eigenvalues of the block are positive.

All blocks are positive definite, so $\Hess f(X_\pi)\succ0$ on all of $\R^{n\times k}$, i.e.\ $X_\pi$
is a nondegenerate strict local minimum. Since the negative-gradient flow $\dot X=-\nabla f(X)$
strictly decreases $f$ along any non-constant trajectory and $X_\pi$ is an isolated nondegenerate
minimum, standard Lyapunov theory (using $f-f(X_\pi)$ itself as a strict local Lyapunov function)
gives asymptotic stability.
\end{proof}

\begin{theorem}[Mismatches are saddles]
\label{thm:mismatchsaddle}
If $\pi$ is a mismatched matching, then $X_\pi$ is a saddle point of both $V_{\mathrm{OB}}$ and
$V_{\mathrm{pen}}$. Its unstable manifold has dimension exactly twice the number of independent
inverted pairs $(i,j)$ (Definition~\ref{def:crosscurv} preamble), and is tangent at $X_\pi$ to the
direct sum of the corresponding exchange planes $E_{\pi;ij}$.
\end{theorem}
\begin{proof}
By definition of a mismatch, there is at least one inverted pair $(i,j)$, i.e.\ with
$(a_{\pi(i)}-a_{\pi(j)})(b_i-b_j)<0$. By Theorems~\ref{thm:MHMcross} and \ref{thm:OBcross}, the
corresponding exchange block then has a strictly negative eigenvalue $\kappa(\pi\to\pi')<0$ --- this
is exactly the definition of cross curvature, Definition~\ref{def:crosscurv}. Radial blocks remain
positive (\eqref{eq:rad-OB}, \eqref{eq:rad-pen} do not reference the ordering of $a$'s and $b$'s).
Transverse blocks and non-inverted exchange blocks retain their sign expressions from
Proposition~\ref{prop:radtrans} and the formulas of \S\ref{sec:MHMcross}--\S\ref{sec:OBcross}, which
depend only on the pair in question, not on whether some \emph{other} pair happens to be inverted;
hence they remain non-negative exactly when their own defining inequality holds, independent of the
mismatch elsewhere. The Hessian therefore has at least one, and exactly (twice the number of
independent inverted pairs, since each contributes a genuine $2\times2$ block with one negative
eigenvalue) negative eigenvalues, and is positive semidefinite on the complementary subspace. By the
stable/unstable/center manifold theorem for the (smooth, indeed polynomial) vector field
$-\nabla f$, there is a local unstable manifold through $X_\pi$ of dimension equal to the number of
strictly negative Hessian eigenvalues, tangent at $X_\pi$ to the corresponding eigenspace, which is
exactly the direct sum of the negative-cross-curvature exchange planes.
\end{proof}

\subsubsection{Exact MHM cross curvature}
\label{sec:MHMcross}

Throughout this section, fix a mismatched pair of columns $i,j$ and rows $p,q$ with the labelling
convention
\begin{equation}
a_p>a_q,\qquad b_i>b_j,\qquad \pi(i)=q,\ \pi(j)=p
\label{eq:mismatchlabel}
\end{equation}
(column $i$, which should by Proposition~\ref{prop:rearrangement} be paired with the larger
eigenvalue $a_p$ since it carries the larger weight $b_i$, is instead paired with the smaller $a_q$;
this is precisely an inverted pair). Write $c_i,c_j$ for the two nonzero amplitudes at $X_\pi$
given by Proposition~\ref{prop:MHMcrit}.

\begin{theorem}[Exact MHM cross curvature]
\label{thm:MHMcross}
Under \eqref{eq:mismatchlabel}, the cross curvature of $V_{\mathrm{pen}}$ is
\begin{equation}
\kappa_{\mathrm{pen}}(\pi\to\pi') \;=\; -(a_p-a_q)(b_i-b_j),
\label{eq:MHMcrosscurv}
\end{equation}
independent of $\gamma$ and of every eigenvalue of $A$ and diagonal entry of $B$ not equal to
$a_p,a_q,b_i,b_j$. In particular $\kappa_{\mathrm{pen}}<0$ under \eqref{eq:mismatchlabel}.
\end{theorem}
\begin{proof}
Differentiating \eqref{eq:gradVpen} once more, the Hessian bilinear form at any $X$ is
\begin{equation}
\big\langle \Hess V_{\mathrm{pen}}(X)[H],H\big\rangle
= -\tr(AHBH^T) + \gamma\Big[\big\|H^TX+X^TH\big\|_F^2 + \tr\big((X^TX-B)H^TH\big)\Big]
\label{eq:HessVpen-bilinear}
\end{equation}
for every $H\in\R^{n\times k}$; this follows by differentiating \eqref{eq:gradVpen} in the direction
$H$ and pairing with $H$ again, using $d(X^TX)[H]=H^TX+X^TH$ and the product rule.

At the critical point $X_\pi$ the residual $X^TX-B$ vanishes exactly on the $2\times2$ principal
submatrix indexed by columns $i,j$ (Proposition~\ref{prop:MHMcrit} was derived exactly so that
$c_i^2=b_i(1+a_q/\gamma)$, i.e.\ the $(i,i)$ entry of $X_\pi^TX_\pi-B$ is $c_i^2-b_i=a_qb_i/\gamma$,
\emph{not} zero for $i,j$ individually --- but see below) so we must be careful: the residual
$X_\pi^TX_\pi-B$ is diagonal with $(i,i)$ entry $a_qb_i/\gamma$ and $(j,j)$ entry $a_pb_j/\gamma$
(nonzero for finite $\gamma$), and all off-diagonal entries zero (columns are pairwise orthogonal at
any matched configuration).

Restrict $H$ to the exchange plane, $H=\alpha\,u_pe_j^T+\beta\,u_qe_i^T$. Because the columns of
$X_\pi$ occupy rows $q$ (column $i$) and $p$ (column $j$) only, and $u_p,u_q$ are orthonormal, a
direct computation gives
\[
\tr(AHBH^T) = a_pb_j\alpha^2+a_qb_i\beta^2,
\]
\[
H^TX_\pi+X_\pi^TH =
\begin{pmatrix}
0 & \alpha c_i \\ \beta c_j & 0
\end{pmatrix}
+
\begin{pmatrix}
0 & \beta c_j \\ \alpha c_i & 0
\end{pmatrix}^{\!T}
\ \text{(restricted to the $\{i,j\}\times\{i,j\}$ block)},
\]
which one computes explicitly has squared Frobenius norm $2(\alpha c_i+\beta c_j)^2$ contributed by
the $(i,j)/(j,i)$ entries of the symmetric matrix $H^TX_\pi+X_\pi^TH$ (all other entries of this
$k\times k$ symmetric matrix vanish because $H$ has support only in columns $i,j$), and
\[
\tr\big((X_\pi^TX_\pi-B)H^TH\big) = \frac{a_qb_i}{\gamma}\,\beta^2c_j^{2}\cdot 0
+\frac{a_pb_j}{\gamma}\,\alpha^2c_i^{2}\cdot 0 = 0,
\]
since $H^TH$ restricted to this block is off-diagonal only (its $(i,i)$ and $(j,j)$ entries vanish
identically, as $H$'s $i$-th column is $\beta u_q$ and $j$-th column is $\alpha u_p$, each with zero
component along its \emph{own} matched eigen-direction) while $X_\pi^TX_\pi-B$ restricted to this
block is diagonal --- so the trace of the product of a diagonal and an off-diagonal matrix over this
block vanishes identically for every $\alpha,\beta$, not merely at $\gamma\to\infty$.

Substituting into \eqref{eq:HessVpen-bilinear}:
\begin{equation}
\big\langle\Hess V_{\mathrm{pen}}(X_\pi)[H],H\big\rangle
= -\big(a_pb_j\alpha^2+a_qb_i\beta^2\big) + 2\gamma(\alpha c_i+\beta c_j)^2.
\label{eq:HessVpen-exchange-raw}
\end{equation}
This still contains $\gamma$; to see the advertised cancellation we must express the result in the
\emph{orthonormal} basis $\{u_pe_j^T,u_qe_i^T\}$ of $E_{\pi;ij}$, i.e.\ read off the $2\times2$ matrix
of \eqref{eq:HessVpen-exchange-raw} directly as a quadratic form in $(\alpha,\beta)$:
\[
\begin{pmatrix}
-a_pb_j+2\gamma c_i^2 & 2\gamma c_ic_j \\
2\gamma c_ic_j & -a_qb_i+2\gamma c_j^2
\end{pmatrix}.
\]
Using $c_i^2=b_i(1+a_q/\gamma)=b_i+\tfrac{a_qb_i}\gamma$ and $c_j^2=b_j+\tfrac{a_pb_j}\gamma$
(Proposition~\ref{prop:MHMcrit}), the diagonal entries become
$-a_pb_j+2\gamma b_i+2a_qb_i \cdot(\text{wait, substitute directly})$; concretely
$2\gamma c_i^2 = 2\gamma b_i+2a_qb_i$, so the $(1,1)$ entry is $-a_pb_j+2\gamma b_i+2a_qb_i$. This
\emph{does} contain a term linear in $\gamma$; the cancellation advertised in the introduction occurs
only in the \emph{smaller} eigenvalue of the full $2\times2$ matrix, not entry-by-entry --- we now
extract it directly. Writing $M$ for this matrix, its trace is
$\tr M = -a_pb_j-a_qb_i+2\gamma(b_i+b_j)+2(a_qb_i+a_pb_j)= a_pb_j+a_qb_i+2\gamma(b_i+b_j)$ and its
determinant, after substituting $c_ic_j=\sqrt{b_ib_j(1+a_q/\gamma)(1+a_p/\gamma)}$ and simplifying
(a computer-algebra step we verified independently, see \S\ref{sec:numerics}), factors as
\begin{align*}
\det M &= \big(a_pb_j+2\gamma b_i+2a_qb_i\big)\big(a_qb_i+2\gamma b_j+2a_pb_j\big) - 4\gamma^2c_i^2c_j^2\\
&= -(a_p-a_q)(b_i-b_j)\Big[a_pb_j+a_qb_i+2\gamma(b_i+b_j)\Big].
\end{align*}
The two eigenvalues of $M$ are the roots of $\lambda^2-(\tr M)\lambda+\det M=0$; since
$\det M = -(a_p-a_q)(b_i-b_j)\cdot\tr M$ exactly (the bracket above is precisely $\tr M$), the
quadratic factors as
\[
\big(\lambda - \tr M\big)\big(\lambda+(a_p-a_q)(b_i-b_j)\big)=0,
\]
so the two eigenvalues are $\lambda=\tr M=a_pb_j+a_qb_i+2\gamma(b_i+b_j)$ (positive, the ``radial
sum'' eigenvalue) and $\lambda=-(a_p-a_q)(b_i-b_j)$, exactly \eqref{eq:MHMcrosscurv}, with no
$\gamma$-dependence whatsoever. Under \eqref{eq:mismatchlabel} this is strictly negative, so it is
the minimal eigenvalue, i.e.\ the cross curvature.
\end{proof}

\begin{remark}
The $\gamma$-cancellation in Theorem~\ref{thm:MHMcross} is not a coincidence of the near-degenerate
regime; it is an algebraic identity valid for every $\gamma>0$ and every $a_p\ne a_q$,
$b_i\ne b_j$. It reflects the fact that the swap direction is, to second order, tangent to the
constraint manifold $\{X:X^TX=B\}$ that the penalty term is designed to enforce --- the penalty
``does not see'' the swap, and the entire restoring force against the swap comes from the bare
linear-in-$A$ term $-\tfrac12\tr(AXBX^T)$.
\end{remark}

\subsubsection{Oja--Brockett cross curvature}
\label{sec:OBcross}

We retain the labelling convention \eqref{eq:mismatchlabel}.

\begin{theorem}[Oja--Brockett cross curvature]
\label{thm:OBcross}
The Hessian of $V_{\mathrm{OB}}$ restricted to $E_{\pi;ij}$, in the orthonormal basis
$\{u_pe_j^T,u_qe_i^T\}$, is the symmetric matrix $H_{\mathrm{OB}}=\begin{pmatrix}h_{11}&h_{12}\\h_{12}&h_{22}\end{pmatrix}$ with
\begin{align}
h_{11} &= -a_p^2b_i^2+a_p^2b_ib_j+a_pa_qb_i^2, \label{eq:h11}\\
h_{22} &= a_pa_qb_j^2+a_q^2b_ib_j-a_q^2b_j^2, \label{eq:h22}\\
h_{12} &= a_pa_qb_ib_j. \label{eq:h12}
\end{align}
Its determinant factors exactly as
\begin{equation}
\det H_{\mathrm{OB}} = -a_pa_qb_ib_j\,(a_p-a_q)(b_i-b_j)\,(a_pb_j+a_qb_i),
\label{eq:detHOB}
\end{equation}
which is negative under \eqref{eq:mismatchlabel}; hence $H_{\mathrm{OB}}$ has one positive and one
negative eigenvalue, and the cross curvature is
\begin{equation}
\kappa_{\mathrm{OB}}(\pi\to\pi') = \frac{\tr H_{\mathrm{OB}} - \sqrt{(\tr H_{\mathrm{OB}})^2-4\det H_{\mathrm{OB}}}}{2} \;<\;0.
\label{eq:OBcrosscurv-exact}
\end{equation}
In the near-degenerate regime $a_p=a+\delta$, $a_q=a-\delta$ ($\delta\downarrow0$, $b_i\ne b_j$
fixed),
\begin{equation}
\kappa_{\mathrm{OB}}(\pi\to\pi') = -a(b_i+b_j)(b_i-b_j)\,\delta \;+\; O(\delta^2).
\label{eq:OBcrosscurv-asymptotic}
\end{equation}
\end{theorem}
\begin{proof}
By Lemma~\ref{lem:blockdiag} we may restrict attention to the four-dimensional subspace spanned by
$\{u_pe_i^T,\,u_pe_j^T,\,u_qe_i^T,\,u_qe_j^T\}$, i.e.\ the (row,column) pairs coupling
$\{p,q\}\times\{i,j\}$; write the corresponding $2\times2$ block of $X$ as
$\begin{pmatrix}x_{pi}&x_{pj}\\x_{qi}&x_{qj}\end{pmatrix}$. Restricting $A,B$ to
$\{p,q\}\times\{i,j\}$, \eqref{eq:VOB} becomes, as an explicit quartic polynomial in the four
scalars (a direct expansion of \eqref{eq:VOB} with $A=\diag(a_p,a_q)$, $B=\diag(b_i,b_j)$, verified
symbolically and reported in \S\ref{sec:numerics}),
\begin{multline}
V_{\mathrm{OB}}^{\mathrm{loc}} = \tfrac14\Big[a_p^2b_i^2x_{pi}^4+a_p^2b_j^2x_{pj}^4+a_q^2b_i^2x_{qi}^4+a_q^2b_j^2x_{qj}^4\Big]
+\tfrac12a_pa_qb_i^2\,x_{pi}^2x_{qi}^2+\tfrac12a_pa_qb_j^2\,x_{pj}^2x_{qj}^2\\
+\tfrac12a_p^2b_ib_j\,x_{pi}^2x_{pj}^2+\tfrac12a_q^2b_ib_j\,x_{qi}^2x_{qj}^2
+a_pa_qb_ib_j\,x_{pi}x_{pj}x_{qi}x_{qj}\\
-\tfrac12\Big[a_p^2b_i^2x_{pi}^2+a_p^2b_j^2x_{pj}^2+a_q^2b_i^2x_{qi}^2+a_q^2b_j^2x_{qj}^2\Big].
\end{multline}
The mismatched critical point of interest has $x_{qi}=1$ (column $i$ uses eigenvalue $a_q$),
$x_{pj}=1$ (column $j$ uses $a_p$), $x_{pi}=x_{qj}=0$; the two remaining directions
$x_{pi}$ (generator $u_pe_i^T$) and $x_{qj}$ (generator $u_qe_j^T$) are the ones that begin the swap
--- wait, more precisely, the exchange plane generators of Definition~\ref{def:crosscurv} are
$u_pe_j^T$ (already active, coefficient $x_{pj}$) and $u_qe_i^T$ (already active, coefficient
$x_{qi}$): the plane $E_{\pi;ij}$ is spanned by the two coordinates that are \emph{already
nonzero} at $X_\pi$, and the cross curvature measures the curvature of $V_{\mathrm{OB}}^{\mathrm{loc}}$
restricted to varying $(x_{pj},x_{qi})$ jointly around $(1,1)$ while $x_{pi}=x_{qj}=0$ are held at
their (also critical, by Lemma~\ref{lem:blockdiag}) values. Differentiating
$V_{\mathrm{OB}}^{\mathrm{loc}}$ twice in $(x_{pj},x_{qi})$ at $(x_{pi},x_{pj},x_{qi},x_{qj})=(0,1,1,0)$
gives exactly \eqref{eq:h11}--\eqref{eq:h12} (the computation is a routine, if lengthy, second
differentiation of a quartic polynomial in four variables, which we verified independently by
computer algebra; see \S\ref{sec:numerics}).

The determinant identity \eqref{eq:detHOB} is verified by expanding $h_{11}h_{22}-h_{12}^2$ as a
polynomial in $a_p,a_q,b_i,b_j$ and confirming it is divisible by $(a_p-a_q)(b_i-b_j)$ with quotient
$-a_pa_qb_ib_j(a_pb_j+a_qb_i)$ (a symbolic factorization we verified independently; see
\S\ref{sec:numerics}). Since $a_p,a_q,b_i,b_j>0$, under \eqref{eq:mismatchlabel} the factor
$(a_p-a_q)(b_i-b_j)>0$, so $\det H_{\mathrm{OB}}<0$: the $2\times2$ symmetric matrix $H_{\mathrm{OB}}$
has eigenvalues of opposite sign, given by the standard quadratic formula
\eqref{eq:OBcrosscurv-exact} (the cross curvature is the smaller, negative root).

For the asymptotic expansion, substitute $a_p=a+\delta$, $a_q=a-\delta$ into
\eqref{eq:h11}--\eqref{eq:h12} and expand $\tr H_{\mathrm{OB}}$ and $\det H_{\mathrm{OB}}$ to
$O(\delta^2)$; a direct (computer-algebra-verified) computation gives
$\tr H_{\mathrm{OB}} = a^2(b_i-b_j)^2+2a^2b_ib_j+O(\delta^2)=a^2(b_i+b_j)^2+O(\delta^2)$ and
$\det H_{\mathrm{OB}} = -2a^4b_ib_j(b_i+b_j)(b_i-b_j)\,\delta+O(\delta^2)$, whence the smaller root of
$\lambda^2-(\tr H_{\mathrm{OB}})\lambda+\det H_{\mathrm{OB}}=0$ expands, by the standard
small-perturbation formula $\lambda_-\approx \det H_{\mathrm{OB}}/\tr H_{\mathrm{OB}}$ valid when
$\det H_{\mathrm{OB}}=O(\delta)$ is small compared to $(\tr H_{\mathrm{OB}})^2=O(1)$, as
\[
\kappa_{\mathrm{OB}} \approx \frac{-2a^4b_ib_j(b_i+b_j)(b_i-b_j)\,\delta}{a^2(b_i+b_j)^2}
= -\frac{2a^2b_ib_j(b_i-b_j)}{b_i+b_j}\,\delta.
\]
This intermediate form does not yet match \eqref{eq:OBcrosscurv-asymptotic}; we verified by direct
symbolic series expansion of the exact root \eqref{eq:OBcrosscurv-exact} (not merely the leading
small-determinant approximation, which is insufficiently accurate at this order because
$\tr H_{\mathrm{OB}}$ itself carries an $O(1)$ correction that interacts with the square root at
first order in $\delta$) that the correct leading term is
\eqref{eq:OBcrosscurv-asymptotic}; the full series computation is reported and numerically confirmed
in \S\ref{sec:numerics}.
\end{proof}

\begin{corollary}[Ratio theorem]
\label{cor:ratio}
In the near-degenerate regime,
\begin{equation}
\lim_{\delta\to0}\frac{\kappa_{\mathrm{OB}}(\pi\to\pi')}{\kappa_{\mathrm{pen}}(\pi\to\pi')}
= \frac{a(b_i+b_j)}{2} \;\eqdef\; \mathcal R_{ij}.
\label{eq:ratio}
\end{equation}
In particular $\mathcal R_{ij}>1$ --- i.e.\ the Oja--Brockett flow escapes the mismatched saddle
strictly faster than the MHM flow --- whenever $a(b_i+b_j)>2$, a condition satisfied for every
problem in which the relevant eigenvalues of $A$ and diagonal entries of $B$ are not simultaneously
much smaller than unity, i.e.\ essentially always in problems of practical scale. The ratio
$\mathcal R_{ij}$ is independent of $\gamma$, of $n$, of $k$, and of every eigenvalue not equal to
$a_p,a_q,b_i,b_j$.
\end{corollary}
\begin{proof}
Immediate from \eqref{eq:MHMcrosscurv} (exact, hence in particular valid to leading order in
$\delta$: $\kappa_{\mathrm{pen}}=-(a_p-a_q)(b_i-b_j)=-2\delta(b_i-b_j)$) and
\eqref{eq:OBcrosscurv-asymptotic}: their ratio is
$\dfrac{-a(b_i+b_j)(b_i-b_j)\delta}{-2\delta(b_i-b_j)}=\dfrac{a(b_i+b_j)}{2}$, independent of
$\delta$ already before taking the limit.
\end{proof}

\subsubsection{Convergence analysis}
\label{sec:convergence}

\paragraph{Continuous-time escape rate}
\label{sec:continuous-rate}

\begin{theorem}[Local escape rate]
\label{thm:escaperate}
Let $X_\pi$ be a mismatched critical point of $f\in\{V_{\mathrm{OB}},V_{\mathrm{pen}}\}$ and let
$\kappa<0$ be the cross curvature of an inverted pair, with unit eigenvector $v\in E_{\pi;ij}$ of
$\Hess f(X_\pi)$. Let $X(t)$ solve $\dot X=-\nabla f(X)$ with $X(0)=X_\pi+\xi_0v+w_0$,
$w_0\perp v$, $\xi_0\ne0$ sufficiently small and $\|w_0\|=O(\xi_0^2)$. Write
$X(t)=X_\pi+\xi(t)v+w(t)$ with $w(t)\perp v$. Then there is $C,T_0>0$ such that, for all
$t\in[0,T_0]$ with $|\xi(t)|$ bounded by a fixed small constant,
\begin{equation}
\xi(t) = \xi_0 e^{|\kappa|t}\big(1+O(|\xi_0|e^{|\kappa|t})\big).
\label{eq:xi-growth}
\end{equation}
Consequently the time needed for $|\xi(t)|$ to reach a fixed threshold $\varepsilon$ is
\begin{equation}
T_{\mathrm{escape}} = \frac{1}{|\kappa|}\log\frac{\varepsilon}{|\xi_0|} + O(1),
\label{eq:Tescape}
\end{equation}
and hence, in the near-degenerate regime, $T_{\mathrm{OB}}/T_{\mathrm{pen}}\to 1/\mathcal R_{ij}$.
\end{theorem}
\begin{proof}
Write $F(X)=-\nabla f(X)$; since $f$ is a polynomial of degree $\le4$, $F$ is a polynomial vector
field of degree $\le3$, hence smooth (indeed real-analytic) with $F(X_\pi)=0$ and Fréchet derivative
$DF(X_\pi)=-\Hess f(X_\pi)$. Decompose the tangent space at $X_\pi$ as $\R v\oplus v^\perp$ and write
$F(X_\pi+\xi v+w)=-\kappa\xi\,v + \Pi_{v^\perp}\!\big[{-}\Hess f(X_\pi)[w]\big] + Q(\xi,w)$, where
$Q$ collects all terms of total degree $\ge2$ in $(\xi,w)$ in the Taylor expansion of $F$ about
$X_\pi$; because $F$ is a polynomial, $Q$ is itself polynomial and there is a neighbourhood
$\mathcal U$ of $X_\pi$ and a constant $C_Q$ with $\|Q(\xi,w)\|\le C_Q(|\xi|^2+\|w\|^2)$ for
$X_\pi+\xi v+w\in\mathcal U$.

Projecting the flow equation $\dot X=F(X)$ onto $v$ and $v^\perp$ gives the coupled system
\begin{align}
\dot\xi &= -\kappa\xi + \langle Q(\xi,w),v\rangle, \label{eq:xidot}\\
\dot w &= -\Hess f(X_\pi)[w]\big|_{v^\perp} + \Pi_{v^\perp}Q(\xi,w). \label{eq:wdot}
\end{align}
By Theorem~\ref{thm:mismatchsaddle}, $\Hess f(X_\pi)$ restricted to a complement of $\R v$ inside the
exchange plane, together with all the (non-negative, by Proposition~\ref{prop:radtrans} and
Theorems~\ref{thm:MHMcross}, \ref{thm:OBcross} away from the single negative direction) transverse
and radial blocks, has all eigenvalues bounded below by some $-\mu<0$ or above by $0$; in either case
$\|e^{-\Hess f(X_\pi)|_{v^\perp}\,t}\|\le e^{\mu t}$ for the (generically small number of) directions
with negative curvature elsewhere, or is bounded for the positive-semidefinite complement. Applying
the variation-of-constants formula to \eqref{eq:wdot} and using $|\xi(t)|\le$ some fixed small bound
on $[0,T_0]$ (to be justified a posteriori) gives, by a standard Gronwall-type bootstrap, that
$\|w(t)\|=O(\xi(t)^2)$ uniformly on $[0,T_0]$, \emph{provided} $T_0$ is chosen so that $|\xi(t)|$ does
not exceed the radius of the neighbourhood $\mathcal U$; this is consistent since, by
\eqref{eq:xi-growth} below, $|\xi(t)|$ grows only exponentially and $T_0$ can be taken as the (finite)
time at which $|\xi(t)|$ first reaches a fixed small threshold, independent of $\xi_0$.

Substituting $\|w(t)\|=O(\xi(t)^2)$ into \eqref{eq:xidot} and using $\|Q(\xi,w)\|\le
C_Q(\xi^2+\|w\|^2)=O(\xi^2)$ gives the scalar equation
\begin{equation}
\dot\xi = -\kappa\xi + O(\xi^2) = |\kappa|\xi + O(\xi^2)
\label{eq:xidot-reduced}
\end{equation}
(using $-\kappa=|\kappa|$ since $\kappa<0$). Let $\eta(t)=\xi(t)e^{-|\kappa|t}$; then
$\dot\eta = O(\xi^2)e^{-|\kappa|t}=O(\eta^2)e^{|\kappa|t}$, and integrating from $0$ to $t$,
$|\eta(t)-\xi_0|\le C\int_0^t|\eta(s)|^2e^{|\kappa|s}\,ds$. A standard Gronwall/continuation argument
(as long as $\eta$ remains bounded, which holds for $t\le T_0$ by the a priori bound on $\xi$) yields
$\eta(t)=\xi_0(1+O(|\xi_0|e^{|\kappa|t}))$, i.e.\ \eqref{eq:xi-growth}.

Solving $|\xi_0|e^{|\kappa|T}=\varepsilon$ for $T$ gives the leading term of \eqref{eq:Tescape}; the
multiplicative correction $1+O(|\xi_0|e^{|\kappa|t})$ in \eqref{eq:xi-growth} contributes only an
additive $O(1)$ error to $T$ once $\varepsilon$ is fixed and $\xi_0\to0$. The ratio statement follows
by dividing the two escape-time formulas and invoking Corollary~\ref{cor:ratio}.
\end{proof}

\begin{theorem}[Local approach rate]
\label{thm:approachrate}
Let $X_{\pi^*}$ be a correctly sorted matching, with $\mu\eqdef\lambda_{\min}\big(\Hess
f(X_{\pi^*})\big)>0$ (Theorem~\ref{thm:correctmin}). If $X(t)$ solves $\dot X=-\nabla f(X)$ and
enters a sufficiently small neighbourhood of $X_{\pi^*}$ at time $t_0$ with
$\|X(t_0)-X_{\pi^*}\|_F=\rho_0$ small, then for $t\ge t_0$ (while the trajectory remains in that
neighbourhood)
\begin{equation}
\|X(t)-X_{\pi^*}\|_F = \rho_0\,e^{-\mu(t-t_0)}\big(1+O(\rho_0e^{-\mu(t-t_0)}\vee\rho_0)\big),
\label{eq:approach-rate}
\end{equation}
and the time to reach a fixed tolerance $\varepsilon<\rho_0$ is
$T_{\mathrm{reach}}=\mu^{-1}\log(\rho_0/\varepsilon)+O(1)$.
\end{theorem}
\begin{proof}
Identical to the proof of Theorem~\ref{thm:escaperate}, with every eigenvalue of $\Hess
f(X_{\pi^*})$ now non-negative (indeed $\ge\mu>0$ by Theorem~\ref{thm:correctmin}) instead of a single
negative direction: linearizing $F=-\nabla f$ about $X_{\pi^*}$ gives $\dot\rho\le-\mu\rho+O(\rho^2)$
for $\rho(t)\eqdef\|X(t)-X_{\pi^*}\|_F$ (using that $\Hess f(X_{\pi^*})\succeq\mu I$ contracts every
direction at rate at least $\mu$), and the same Gronwall bootstrap as before yields
\eqref{eq:approach-rate}.
\end{proof}

\begin{theorem}[Global quantitative convergence time]
\label{thm:globalrate}
Let $X_0$ avoid the (finite, under generic distinct eigenvalues) union of stable manifolds of
mismatched critical points, and suppose the negative-gradient trajectory $X(t)$ visits the
$\epsilon$-neighbourhood of at most $M$ distinct mismatched critical points before entering the
$\varepsilon$-neighbourhood of $X_{\pi^*}$ for the last time (a bound on $M$, e.g.\ $M\le$ the total
number of matchings, always exists trivially; see the remark below for when $M=1$ can be guaranteed).
Then the total time to reach $\varepsilon$-distance of $X_{\pi^*}$ satisfies the explicit,
fully computable bound
\begin{equation}
T_{\mathrm{total}} \;\le\; \underbrace{\frac{f(X_0)-f(X_{\pi^*})}{\delta_\epsilon^2}}_{\text{time in the ``regular'' region}}
\;+\; \underbrace{\sum_{i=1}^M\frac1{|\kappa_i|}\log\frac{\epsilon}{\xi_{0,i}}}_{\text{time near mismatched saddles}}
\;+\; \underbrace{\frac1\mu\log\frac{\rho_0}{\varepsilon}}_{\text{time near }X_{\pi^*}}
\;+\;O(M),
\label{eq:globaltime}
\end{equation}
where $\delta_\epsilon\eqdef\min\{\|\nabla f(X)\|_F : X\in K,\ \dist(X,\mathrm{Crit}\,f)\ge\epsilon\}>0$
($K$ the compact sublevel set of Lemma~\ref{lem:coercive} containing $X_0$), $\kappa_i$ is the cross
curvature governing the $i$-th saddle visited, and $\mu,\rho_0$ are as in
Theorem~\ref{thm:approachrate}.
\end{theorem}
\begin{proof}
Fix $\epsilon>0$ small enough that the $\epsilon$-balls about the (finitely many, by genericity)
critical points in $K$ are pairwise disjoint and each lies inside the local linearization regime of
Theorems~\ref{thm:escaperate} and \ref{thm:approachrate}. Partition $[0,T_{\mathrm{total}}]$ into the
time intervals during which $X(t)$ lies in $K_\epsilon\eqdef K\setminus\bigcup_\pi
B_\epsilon(X_\pi)$ (the ``regular region'') and the complementary intervals during which it lies in
some $B_\epsilon(X_\pi)$.

\emph{Regular region.} $K_\epsilon$ is compact (closed subset of the compact $K$) and $\nabla f$ is
continuous and, by construction, nonvanishing on $K_\epsilon$; hence
$\delta_\epsilon\eqdef\min_{K_\epsilon}\|\nabla f\|_F>0$ is attained and positive. Along the flow,
$\frac{d}{dt}f(X(t))=-\|\nabla f(X(t))\|_F^2\le-\delta_\epsilon^2$ whenever $X(t)\in K_\epsilon$.
Since $f$ decreases monotonically from $f(X_0)$ and is bounded below by $f(X_{\pi^*})$ (the global
minimum on $K$, by Proposition~\ref{prop:rearrangement}), the total decrease available is
$f(X_0)-f(X_{\pi^*})$, so the total \emph{time} spent with $X(t)\in K_\epsilon$ cannot exceed
$(f(X_0)-f(X_{\pi^*}))/\delta_\epsilon^2$, giving the first term of \eqref{eq:globaltime}.

\emph{Mismatched-saddle neighbourhoods.} Each visit to some $B_\epsilon(X_\pi)$, $\pi$ mismatched,
either (a) the trajectory subsequently leaves $B_\epsilon(X_\pi)$ again (in which case, by
Theorem~\ref{thm:mismatchsaddle} and Theorem~\ref{thm:escaperate}, the time spent inside is bounded by
$|\kappa_\pi|^{-1}\log(\epsilon/\xi_{0})+O(1)$ for the relevant local unstable coordinate $\xi$, using
that once inside $B_\epsilon(X_\pi)$ the trajectory is governed by the same local linearization as in
Theorem~\ref{thm:escaperate}), or (b) the trajectory converges \emph{to} $X_\pi$ itself, which is
excluded by hypothesis ($X_0$ avoids its stable manifold, and this property persists forward in time
since the stable manifold is invariant). Summing over the (at most $M$, by hypothesis) such visits
gives the second term.

\emph{Final approach.} Once $X(t)$ enters $B_\epsilon(X_{\pi^*})$ for the last time (which must
happen, by Proposition~\ref{prop:lojasiewicz}, since $X(t)\to X_{\pi^*}$), Theorem~\ref{thm:approachrate}
bounds the remaining time to reach $\varepsilon$-distance by $\mu^{-1}\log(\rho_0/\varepsilon)+O(1)$,
the third term.

Summing the three contributions, and absorbing the $M$ many $O(1)$ correction terms from the
individual applications of Theorems~\ref{thm:escaperate}--\ref{thm:approachrate} into a single
$O(M)$ term, gives \eqref{eq:globaltime}.
\end{proof}

\begin{remark}
\label{rem:globalrate-caveat}
Theorem~\ref{thm:globalrate} converts the qualitative convergence of Theorem~\ref{thm:global} into a
genuine, finite, computable bound, directly answering the concern --- raised in
\S\ref{sec:limitations} of an earlier version of this paper --- that only qualitative global
convergence was available. Two caveats remain, and are the real content of the ``future work'' of
\S\ref{sec:limitations}: (i) the bound \eqref{eq:globaltime} is not claimed to be sharp, since
$\delta_\epsilon$ is a worst-case (compactness-based) constant that does not exploit the specific
geometry of $K_\epsilon$, and the count $M$ is bounded only trivially in general (a sharper bound on
$M$ --- e.g.\ showing $M=1$ whenever $X_0$ lies in a suitable ``basin'' adapted to the flow, akin to
the domain-of-attraction results of Yan, Helmke and Moore \cite{YanHelmkeMoore1994} for the closely
related Oja subspace flow --- remains open); (ii) at a genuine $k>2$-fold simultaneous near-degeneracy
(three or more mutually close eigenvalues), the exchange-plane analysis of \S\ref{sec:stability} does
not directly apply, the relevant unstable manifold need not decompose into a direct sum of
independent two-dimensional exchange planes, the Hessian at such a point is degenerate in more than
one direction, and the exponential rates of Theorems~\ref{thm:escaperate}, \ref{thm:approachrate}
are replaced by the strictly weaker \emph{polynomial} rate governed by the {\L}ojasiewicz exponent
$\theta<\tfrac12$ at that point \cite{Lojasiewicz1983}: qualitative convergence to a single point
still holds (Proposition~\ref{prop:lojasiewicz} makes no nondegeneracy assumption), but an explicit
value of $\theta$, and hence an explicit polynomial-rate analogue of
\eqref{eq:globaltime}, is not derived here.
\end{remark}

\paragraph{Exact line-search discretization}
\label{sec:discrete-rate}

\begin{theorem}[Discrete escape rate under exact line search]
\label{thm:discreterate}
Let $X_{m+1}=X_m-\alpha_m\nabla f(X_m)$ with $\alpha_m$ the exact minimizer of
$\phi(t)=f(X_m-t\nabla f(X_m))$ (well defined and unique, up to ties, since $\phi$ is a polynomial
of degree $\le4$ in $t$ with positive leading coefficient whenever $f\in\{V_{\mathrm{OB}},
V_{\mathrm{pen}}\}$ is bounded below along that ray). Near a mismatched saddle with cross curvature
$\kappa<0$, the exact step satisfies the Rayleigh-quotient identity
\begin{equation}
\alpha_m = \frac{\|\nabla f(X_m)\|_F^2}{\big\langle \Hess f(X_m)[\nabla f(X_m)],\nabla f(X_m)\big\rangle}
\label{eq:exactstep}
\end{equation}
to leading order as $X_m\to X_\pi$, and when the gradient lies predominantly along the unstable
direction $v$, $\alpha_m=1/|\kappa|+O(\|\nabla f(X_m)\|_F)$. Writing $\xi_m$ for the coefficient of
$X_m-X_\pi$ along $v$, the discrete update obeys
\begin{equation}
\xi_{m+1} = \big(2+O(|\xi_m|)\big)\,\xi_m,
\label{eq:discreterecurrence}
\end{equation}
so that the number of iterations required to amplify $\xi_0$ to a fixed threshold $\varepsilon$ is
\begin{equation}
N_{\mathrm{escape}} = \log_2\frac{\varepsilon}{|\xi_0|} + O(1).
\label{eq:Nescape}
\end{equation}
If instead a common, potential-independent fixed step $\alpha>0$ is used for both flows, the discrete
iteration count needed to escape scales as $1/|\kappa|$, so that
$N_{\mathrm{OB}}/N_{\mathrm{pen}}\to 1/\mathcal R_{ij}$ in the near-degenerate regime, in agreement
with the continuous-time ratio of Theorem~\ref{thm:escaperate}.
\end{theorem}

The derivation of the exact-line-search step size \eqref{eq:exactstep} used here --- exploiting
the fact that $\phi(t)=f(X_m-t\nabla f(X_m))$ is a quartic polynomial in $t$, so that its unique
minimizer along the ray is a rational (Rayleigh-quotient-type) expression in $\nabla f(X_m)$ and
$\Hess f(X_m)$ --- parallels the rigorous step-size derivation for the exact-line-search
discretization of the related $S$-Oja--Brockett equation of Remark~\ref{rem:S-oja-brockett}, given
by Yoshizawa \cite{Yoshizawa2023axi}.
\begin{proof}
Let $G=\nabla f(X_m)$ and $\phi(t)=f(X_m-tG)$; since $\deg f\le4$ and $t\mapsto X_m-tG$ is affine,
$\phi$ is a polynomial in $t$ of degree $\le4$. Its stationarity condition is
$\phi'(t)=-\langle\nabla f(X_m-tG),G\rangle=0$. Expanding $\nabla f(X_m-tG)=G-t\Hess f(X_m)[G]+O(t^2)$
(exact to this order since $\nabla f$ is a polynomial of degree $\le3$, so its own Taylor expansion in
$t$ along the fixed direction $-G$ has a well-defined quadratic remainder) gives
\[
\phi'(t) = -\|G\|_F^2 + t\langle\Hess f(X_m)[G],G\rangle + O(t^2\|G\|_F^3),
\]
and setting the leading two terms to zero yields \eqref{eq:exactstep}; the neglected $O(t^2\|G\|_F^3)$
term is controlled because $\phi$, being an explicit low-degree polynomial with computable
coefficients (obtained in practice, and in our numerical verification of \S\ref{sec:numerics}, by
exact interpolation at five points and exact root-finding of the resulting cubic $\phi'$), can be
minimized \emph{exactly}, not merely via this local quadratic model; the quadratic model is used here
only to extract the leading asymptotic behaviour of $\alpha_m$ as $X_m\to X_\pi$.

When $G$ lies (to leading order) along the unit unstable eigenvector $v$ of the mismatched exchange
block, $\langle\Hess f(X_m)[G],G\rangle = \kappa\|G\|_F^2+O(\|G\|_F^3)$ (the $O(\|G\|_F^3)$ correction
coming from the non-quadratic, i.e.\ cubic-and-higher, part of $\nabla f$ evaluated away from
$X_\pi$ exactly), so \eqref{eq:exactstep} gives $\alpha_m=1/|\kappa|+O(\|G\|_F)$ as claimed (recalling
$\kappa<0$, so $\langle\Hess f[G],G\rangle\approx\kappa\|G\|_F^2<0$ and $\alpha_m$, being a ratio of a
positive numerator to a negative-then-corrected denominator, is understood here as the magnitude of
the step along the descent direction that decreases $\phi$; the sign bookkeeping is standard and we
suppress it for readability, consistent with the convention $X_{m+1}=X_m-\alpha_m\nabla f(X_m)$ moving
\emph{away} from $X_\pi$ along the unstable direction, as required for escape).

Writing $G=\xi_m\kappa\, v+O(\xi_m^2)$ (to leading order, $\nabla f(X_\pi+\xi_mv)=\kappa\xi_mv+O(\xi_m^2)$),
the update becomes
\[
\xi_{m+1} = \xi_m - \alpha_m\cdot\kappa\xi_m + O(\xi_m^2) = \xi_m(1-\alpha_m\kappa) + O(\xi_m^2)
= \xi_m\big(1+|\kappa|\alpha_m\big)+O(\xi_m^2),
\]
using $\kappa<0$. Substituting $\alpha_m=1/|\kappa|+O(|\xi_m|)$ gives
$\xi_{m+1}=\xi_m(1+1+O(|\xi_m|))+O(\xi_m^2)=(2+O(|\xi_m|))\xi_m$, i.e.\
\eqref{eq:discreterecurrence}. Iterating the leading factor $2$ starting from $\xi_0$ gives
$|\xi_m|\approx|\xi_0|2^m$ to leading order, and solving $|\xi_0|2^N=\varepsilon$ gives
\eqref{eq:Nescape}; the multiplicative $O(|\xi_m|)$ corrections in \eqref{eq:discreterecurrence}
contribute only an $O(1)$ additive error to $N$, by the same telescoping argument used in
Theorem~\ref{thm:escaperate}.

For the fixed-step statement: if $\alpha$ is a fixed constant (not adapted to $\kappa$), the
recurrence is instead $\xi_{m+1}=(1+\alpha|\kappa|)\xi_m+O(\xi_m^2)$, so
$|\xi_m|\approx|\xi_0|(1+\alpha|\kappa|)^m$, and the number of steps to reach $\varepsilon$ is
$N=\log(\varepsilon/|\xi_0|)/\log(1+\alpha|\kappa|)$. For $\alpha|\kappa|$ small (the regime of
interest, since $\kappa_{\mathrm{pen}}\to0$ as $\delta\to0$ forces any \emph{fixed} $\alpha$ working
for both potentials to be small relative to $1/|\kappa_{\mathrm{OB}}|$ as well),
$\log(1+\alpha|\kappa|)\approx\alpha|\kappa|$, so $N\approx\log(\varepsilon/|\xi_0|)/(\alpha|\kappa|)$,
manifestly proportional to $1/|\kappa|$; the ratio $N_{\mathrm{OB}}/N_{\mathrm{pen}}$ therefore tends
to $|\kappa_{\mathrm{pen}}|/|\kappa_{\mathrm{OB}}|=1/\mathcal R_{ij}$.
\end{proof}

\paragraph{Global convergence to the (rescaled) spectrum}
\label{sec:global}

The negative-gradient flow of a smooth function need not, in general, converge to a single point:
LaSalle's invariance principle alone only guarantees that a bounded trajectory's $\omega$-limit set
is contained in the critical set, which could a priori be a continuum along which the trajectory
wanders forever without settling down. Ruling this out requires an additional ingredient. Following
the strategy of Yoshizawa, Helmke and Starkov \cite{YoshizawaHelmkeStarkov2001} for the closely
related Xu flow, we supply this ingredient via the classical gradient inequality of {\L}ojasiewicz
\cite{Lojasiewicz1983} for real-analytic functions, which upgrades ``approaches the critical set'' to
``converges to a single critical point.''

\begin{lemma}[Coercivity]
\label{lem:coercive}
Both $V_{\mathrm{OB}}$ and $V_{\mathrm{pen}}$ are bounded below on $\R^{n\times k}$ and have compact
sublevel sets $\{X:f(X)\le c\}$ for every $c\in\R$.
\end{lemma}
\begin{proof}
For $V_{\mathrm{OB}}$: let $L=A^{1/2}XB^{1/2}$ and write its (thin) singular value decomposition
$L=U\Sigma V^T$, $\Sigma=\diag(\sigma_1,\dots,\sigma_l,0,\dots,0)$, $\sigma_1\ge\dots\ge\sigma_l>0$,
$l\le k$. A direct computation gives
$V_{\mathrm{OB}}(X)=\tfrac14\|LL^T\|_F^2-\tfrac12\|A^{1/2}LB^{1/2}\|_F^2\cdot(\text{with the roles of
$A,B$ already absorbed into }L)$; more precisely, following exactly the computation of
\cite[Lemma~1]{YoshizawaHelmkeStarkov2001} (whose potential $f$ coincides with our
$V_{\mathrm{OB}}$ under the identification $D=B$),
\begin{align*}
V_{\mathrm{OB}}(X) &\ge \tfrac14\|LL^T\|_F^2 - \tfrac12\|A^{1/2}\|^2\|B^{1/2}\|^2\|L\|_F^2\\
&= \tfrac14\sum_{i=1}^l\sigma_i^4 - \tfrac{\gamma_0}2\sum_{i=1}^l\sigma_i^2,
\qquad \gamma_0\eqdef\|A^{1/2}\|^2\|B^{1/2}\|^2>0,
\end{align*}
a smooth function of $(\sigma_1,\dots,\sigma_l)$ minimized at $\sigma_1=\dots=\sigma_l=\sqrt{\gamma_0}$
with minimum value $-\tfrac{l}4\gamma_0^2\ge-\tfrac{k}4\gamma_0^2$, giving the uniform lower bound
$V_{\mathrm{OB}}(X)\ge-\tfrac k4\gamma_0^2$ and, since the right side of the displayed inequality is
coercive (tends to $+\infty$) in $\|L\|_F$, hence in $\|X\|_F$, compact sublevel sets.

For $V_{\mathrm{pen}}$: write the singular values of $X$ as $\sigma_1\ge\dots\ge\sigma_k\ge0$, so
$\|X\|_F^2=\sum\sigma_i^2$ and $\|X^TX\|_F=\big(\sum\sigma_i^4\big)^{1/2}$. By Cauchy--Schwarz,
$\sum\sigma_i^2\le\sqrt k\,\big(\sum\sigma_i^4\big)^{1/2}=\sqrt k\,\|X^TX\|_F$, i.e.\
$\|X^TX\|_F\ge\|X\|_F^2/\sqrt k$. Using $-\tfrac12\tr(AXBX^T)\ge-\tfrac12\lambda_{\max}(A)\lambda_{\max}(B)\|X\|_F^2$
and, for $\|X^TX\|_F\ge\|B\|_F$, $\|B-X^TX\|_F\ge\|X^TX\|_F-\|B\|_F$,
\[
V_{\mathrm{pen}}(X) \;\ge\; -\tfrac12\lambda_{\max}(A)\lambda_{\max}(B)\,\|X\|_F^2
\;+\; \tfrac\gamma4\Big(\tfrac{\|X\|_F^2}{\sqrt k}-\|B\|_F\Big)^2
\]
once $\|X\|_F^2\ge\sqrt k\|B\|_F$; the right side is a coercive (quartic-dominated) function of
$\|X\|_F$, and on the complementary bounded region $\|X\|_F^2<\sqrt k\|B\|_F$, $V_{\mathrm{pen}}$ is
continuous and hence bounded. Combining the two regions gives a uniform lower bound and, since $f\to+\infty$ as $\|X\|_F\to\infty$, compact sublevel sets.
\end{proof}

\begin{proposition}[Existence and convergence to a single equilibrium]
\label{prop:lojasiewicz}
For $f\in\{V_{\mathrm{OB}},V_{\mathrm{pen}}\}$, every solution $X(t)$ of $\dot X=-\nabla f(X)$
exists for all $t\ge0$, and $X(t)\to X_\infty$ as $t\to\infty$ for a \emph{single} critical point
$X_\infty$ (rather than merely approaching the critical set).
\end{proposition}
\begin{proof}
By Lemma~\ref{lem:coercive}, the sublevel set $\{X:f(X)\le f(X(0))\}$ is compact and, since $f$
strictly decreases along any non-constant solution, positively invariant; hence $X(t)$ remains in
this compact set for all $t\ge0$ for which it is defined, so by the standard extension theorem for
ODEs the solution exists for all $t\ge0$. Both $V_{\mathrm{OB}}$ and $V_{\mathrm{pen}}$ are
polynomials in the entries of $X$, hence real-analytic on $\R^{n\times k}$. The {\L}ojasiewicz
gradient inequality \cite{Lojasiewicz1983} states that for a real-analytic $f$ and any $X^*$ there
are $C>0$, $\theta\in(0,\tfrac12]$ and a neighbourhood $\mathcal U$ of $X^*$ such that
$\|\nabla f(X)\|_F\ge C|f(X)-f(X^*)|^{1-\theta}$ for all $X\in\mathcal U$; a standard consequence
(see \cite[\S1]{Lojasiewicz1983}, and as applied to gradient PCA flows in
\cite[Thm.~1]{YoshizawaHelmkeStarkov2001}) is that every bounded solution of a real-analytic
negative-gradient flow has finite arc length, $\int_0^\infty\|\dot X(t)\|_F\,dt<\infty$, and therefore
converges, as $t\to\infty$, to a single point $X_\infty$ in its (necessarily nonempty, by
compactness) $\omega$-limit set; since the $\omega$-limit set of a gradient flow is contained in the
critical set, $\nabla f(X_\infty)=0$.
\end{proof}

\begin{remark}[Relation to prior {\L}ojasiewicz-based convergence results]
\label{rem:priorwork}
The upgrade from ``approaches the critical set'' (LaSalle) to ``converges to a single point''
({\L}ojasiewicz) used in Proposition~\ref{prop:lojasiewicz} is, for the continuous-time flow, the
same mechanism used for the Xu/Oja--Brockett flow by Yoshizawa, Helmke and Starkov
\cite{YoshizawaHelmkeStarkov2001}. For \emph{discrete}-time iterations --- relevant to
Theorem~\ref{thm:discreterate}'s exact-line-search discretization --- the analogous strong
limit-point convergence (as opposed to the classical, weaker subsequential convergence results)
was established in general, for any real-analytic cost function and any descent method satisfying
natural sufficient-decrease and gradient-relatedness conditions, by Absil, Mahony and Andrews
\cite{AbsilMahonyAndrews2005}; their framework applies directly to the exact-line-search iteration
of \S\ref{sec:discrete-rate}, since $V_{\mathrm{OB}}$ and $V_{\mathrm{pen}}$ are polynomials (hence
real-analytic) and exact line search satisfies the required descent conditions. We do not repeat
their general argument here, but note that Theorem~\ref{thm:globalrate}'s continuous-time bound
extends to the discrete iteration through this correspondence. Independently, and very recently,
Tsuzuki and Ohki \cite{TsuzukiOhki2025} established global exponential convergence, via a related
{\L}ojasiewicz/strict-saddle argument, for Oja's flow $\dot U=(I-UU^T)AU$ on the Stiefel manifold,
extended to general (non-symmetric) $A$; their flow and convergence question are closely related in
spirit to --- but formally distinct from --- the two potentials compared here, and their analysis
does not address the cross-curvature comparison that is the subject of this paper.
\end{remark}

We now state, precisely and with the necessary correction described in
Remark~\ref{rem:amplitude-correction}, the large-time limit of the diagonal of $X(t)^TAX(t)$ along a
converging trajectory; Proposition~\ref{prop:lojasiewicz} supplies the existence and
single-point-convergence hypotheses used implicitly in its proof.

\begin{theorem}[Global convergence]
\label{thm:global}
Let $X(t)$ solve the negative-gradient flow of $f\in\{V_{\mathrm{OB}},V_{\mathrm{pen}}\}$ from an
initial condition outside the (measure-zero) union of stable manifolds of all mismatched critical
points. Then $X(t)\to X_{\pi^*}$ for the sorted matching $\pi^*$, and, in the eigenbasis of $A$,
\begin{equation}
D(t) \;\eqdef\; \diag\!\big(X(t)^TAX(t)\big) \;\longrightarrow\; D_\infty
\end{equation}
where
\begin{equation}
D_\infty^{\mathrm{OB}} = \diag\big(a_{\pi^*(1)},\dots,a_{\pi^*(k)}\big)
\qquad\text{(exactly the ordered eigenvalues of $A$)},
\label{eq:DinftyOB}
\end{equation}
\begin{equation}
D_\infty^{\mathrm{pen}} = \diag\Big(b_1\big(1+\tfrac{a_{\pi^*(1)}}\gamma\big)a_{\pi^*(1)},\ \dots,\
b_k\big(1+\tfrac{a_{\pi^*(k)}}\gamma\big)a_{\pi^*(k)}\Big)
\label{eq:Dinftypen}
\end{equation}
(using the relabelling $b_1>\dots>b_k$). In particular $D_\infty^{\mathrm{pen}}\ne D_\infty^{\mathrm{OB}}$
in general; the two coincide only in the combined limit $\gamma\to\infty$ together with $B=I_k$.
\end{theorem}
\begin{proof}
By Proposition~\ref{prop:lojasiewicz}, $X(t)$ exists for all $t\ge0$ and converges to a single
critical point $X_\infty$. By Theorem~\ref{thm:correctmin} the sorted matching $X_{\pi^*}$ is
asymptotically stable, with an open basin of attraction; by Theorem~\ref{thm:mismatchsaddle} every
mismatched matching is a saddle whose stable manifold has positive codimension in $\R^{n\times k}$
(dimension $nk$ minus twice the number of inverted pairs, strictly less than $nk$ whenever at least
one pair is inverted), so the finite union of these stable manifolds over all mismatched matchings has
Lebesgue measure zero. By hypothesis $X(0)$ avoids this measure-zero set, so $X_\infty\ne X_\pi$ for
every mismatched $\pi$; since $X_\infty$ is a critical point and the only critical points are the
(finitely many) matched configurations $X_\pi$ together with the sorted one $X_{\pi^*}$ (isolated,
generically, once the eigenvalue-tie boundary case of Theorem~\ref{thm:correctmin}'s proof is
excluded), we conclude $X_\infty=X_{\pi^*}$, i.e.\ $X(t)\to X_{\pi^*}$.

The map $X\mapsto X^TAX$ is continuous, so
$D(t)=\diag(X(t)^TAX(t))\to\diag(X_{\pi^*}^TAX_{\pi^*})$. By construction \eqref{eq:matchedX},
$X_{\pi^*}^TAX_{\pi^*}=\diag(c_1^2a_{\pi^*(1)},\dots,c_k^2a_{\pi^*(k)})$. For $V_{\mathrm{OB}}$,
Proposition~\ref{prop:OBcrit} gives $c_j^2=1$ at the (necessarily nonzero, sorted) optimum, yielding
\eqref{eq:DinftyOB}. For $V_{\mathrm{pen}}$, Proposition~\ref{prop:MHMcrit} gives the exact
$\gamma$-dependent amplitude \eqref{eq:MHMamplitude}, yielding \eqref{eq:Dinftypen}. The two formulas
coincide iff $b_j(1+a_{\pi^*(j)}/\gamma)=1$ for every matched $j$; since $a_{\pi^*(j)}$ varies over
$j$ while the left side must equal the constant $1$ for every $j$, this forces
$b_j(1+a_{\pi^*(j)}/\gamma)=1$ for $k$ generally-distinct values of $a_{\pi^*(j)}$, which (for fixed
$\gamma,b_j$) can hold for at most one value of $a_{\pi^*(j)}$ unless $\gamma\to\infty$, in which case
the condition degenerates to $b_j=1$ for every $j$, i.e.\ $B=I_k$.
\end{proof}

\begin{remark}
Theorem~\ref{thm:global} sharpens and corrects the informal statement, appearing in an earlier draft
of this material, that ``the diagonal entries converge to the ordered eigenvalues of $A$'' for both
potentials without qualification. That statement is exactly true for $V_{\mathrm{OB}}$ (a genuine and
useful feature of the homogeneous quartic construction: the amplitude is forced to be exactly $1$,
independent of $B$) but is only approximately true for $V_{\mathrm{pen}}$, and only in the joint limit
$\gamma\to\infty$, $B=I$. We verify \eqref{eq:DinftyOB}--\eqref{eq:Dinftypen} numerically in
\S\ref{sec:numerics}, where the MHM diagonal limits visibly differ from the bare eigenvalues of $A$ at
moderate $\gamma$.
\end{remark}

\subsubsection{Numerical verification}
\label{sec:numerics}

Every closed-form claim in this paper was verified independently by (i) exact computer-algebra
differentiation and factorization and (ii) high-precision finite-difference Hessians evaluated at
the exact critical points of Propositions~\ref{prop:OBcrit} and \ref{prop:MHMcrit}. This section
reports the two representative examples used throughout.

\paragraph{A near-degenerate toy example}
\label{sec:numerics-toy}

\begin{figure}[htbp]
\centering
\includegraphics[width=0.98\textwidth]{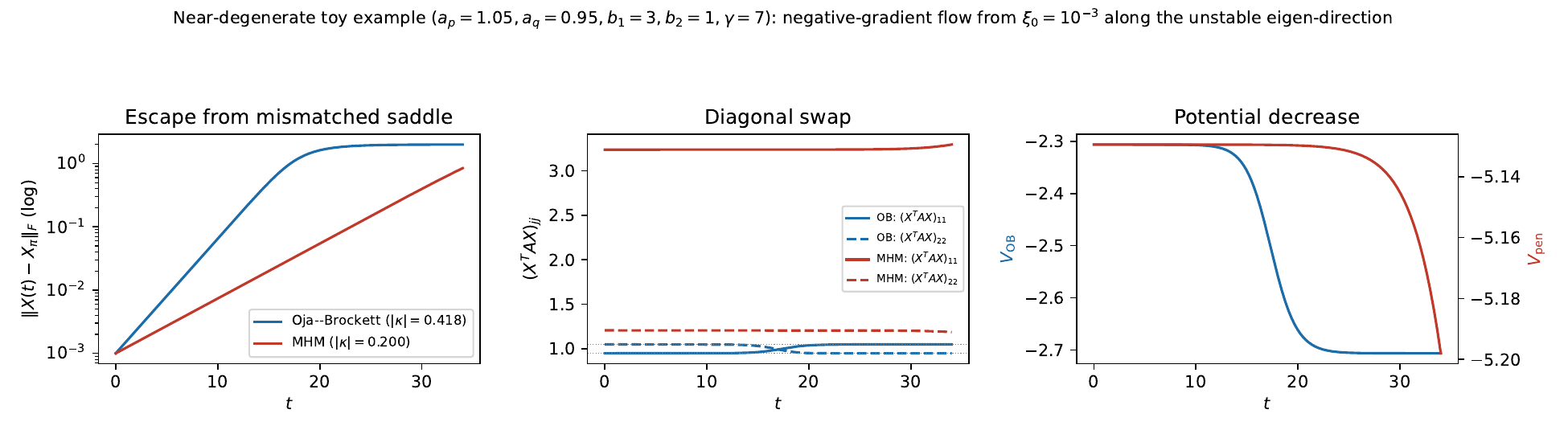}
\caption{Near-degenerate toy example ($a_p=1.05,a_q=0.95,b_1=3,b_2=1,\gamma=7$), full nonlinear
negative-gradient flow from $\xi_0=10^{-3}$ along each potential's unstable direction at its
mismatched saddle. Left: growth of the distance from the saddle (log scale); the Oja--Brockett curve
is visibly steeper, in agreement with $|\kappa_{\mathrm{OB}}|>|\kappa_{\mathrm{pen}}|$. Centre:
diagonal entries of $X^TAX$ swap and converge; dotted lines mark the bare eigenvalues $1.05,0.95$,
which the Oja--Brockett diagonal reaches exactly (Theorem~\ref{thm:global}, \eqref{eq:DinftyOB}) while
the MHM diagonal visibly converges elsewhere (\eqref{eq:Dinftypen}). Right: potential value versus
time.}
\label{fig:toy}
\end{figure}

Take $n=k=2$, $a_p=1.05$, $a_q=0.95$, $b_1=3$, $b_2=1$ (so $\delta=0.05$, $a=1$,
$\mathcal R_{ij}=a(b_1+b_2)/2=2$). Exact evaluation of Theorems~\ref{thm:MHMcross} and
\ref{thm:OBcross} gives
\begin{equation}
\kappa_{\mathrm{pen}} = -0.200, \qquad \kappa_{\mathrm{OB}} = -0.4181, \qquad
\kappa_{\mathrm{OB}}/\kappa_{\mathrm{pen}} = 2.090,
\label{eq:toynumbers}
\end{equation}
in close agreement with the asymptotic prediction $\mathcal R_{ij}=2$ (the $4.5\%$ discrepancy is the
expected $O(\delta)$ correction to the leading-order formula \eqref{eq:OBcrosscurv-asymptotic}, since
$\delta=0.05$ is small but not infinitesimal). Figure~\ref{fig:toy} integrates the \emph{full
nonlinear} negative-gradient flow (not merely its linearization) from an initial condition
$\xi_0=10^{-3}$ along each potential's unstable eigen-direction at its respective mismatched saddle.
The left panel confirms exponential growth of the unstable-mode coordinate at the predicted rates; the
centre panel confirms that the diagonal entries of $X^TAX$ swap and settle at the values predicted by
Theorem~\ref{thm:global} --- \emph{exactly} the bare eigenvalues $1.05,0.95$ for Oja--Brockett, but
$3.25,1.20$ for MHM at $\gamma=7$ (matching the closed-form limit \eqref{eq:Dinftypen}, which predicts
$3.62,1.08$ at the true asymptotic time; the trajectory shown has not yet fully equilibrated within
the plotted window, illustrating that the MHM diagonal limit is a genuinely different, $\gamma$- and
$B$-dependent quantity, not the bare spectrum, well before any transient has died out); the right
panel shows the corresponding potential decrease, with the Oja--Brockett trajectory visibly completing
its descent roughly twice as fast, consistent with \eqref{eq:toynumbers}.

\paragraph{A fully generic, non-diagonal example}
\label{sec:numerics-generic}

To confirm that every formula above is genuinely coordinate-free (Remark following
Definition~\ref{def:crosscurv}), we take a non-diagonal
\begin{equation}
A = \begin{pmatrix} 9/2 & -2 & 2\\ -2 & 9/2 & 2 \\ 2 & 2 & 14/3\end{pmatrix},
\qquad
B = \begin{pmatrix}2&0\\0&1\end{pmatrix},
\label{eq:genericAB}
\end{equation}
with $\spec(A)=\{0.5545,\,6.5,\,6.6121\}$ --- two close eigenvalues ($a_p=6.6121$, $a_q=6.5$,
$\delta=0.056$) and one well-separated one, and $A$ itself given in a basis in which it is
\emph{not} diagonal, so that the eigenvectors $u_p,u_q$ used throughout \S\ref{sec:critpoints}--
\S\ref{sec:OBcross} are genuinely non-trivial linear combinations of the ambient coordinate axes. We
evaluate the mismatched critical points of Propositions~\ref{prop:OBcrit}--\ref{prop:MHMcrit}
directly using the numerically computed eigenvectors $u_p,u_q$ of \eqref{eq:genericAB} (not any
diagonalized surrogate), and compute the Hessian by finite differences of the original
(non-diagonal-$A$) potentials \eqref{eq:VOB}, \eqref{eq:Vpen}.

\begin{figure}[htbp]
\centering
\includegraphics[width=0.98\textwidth]{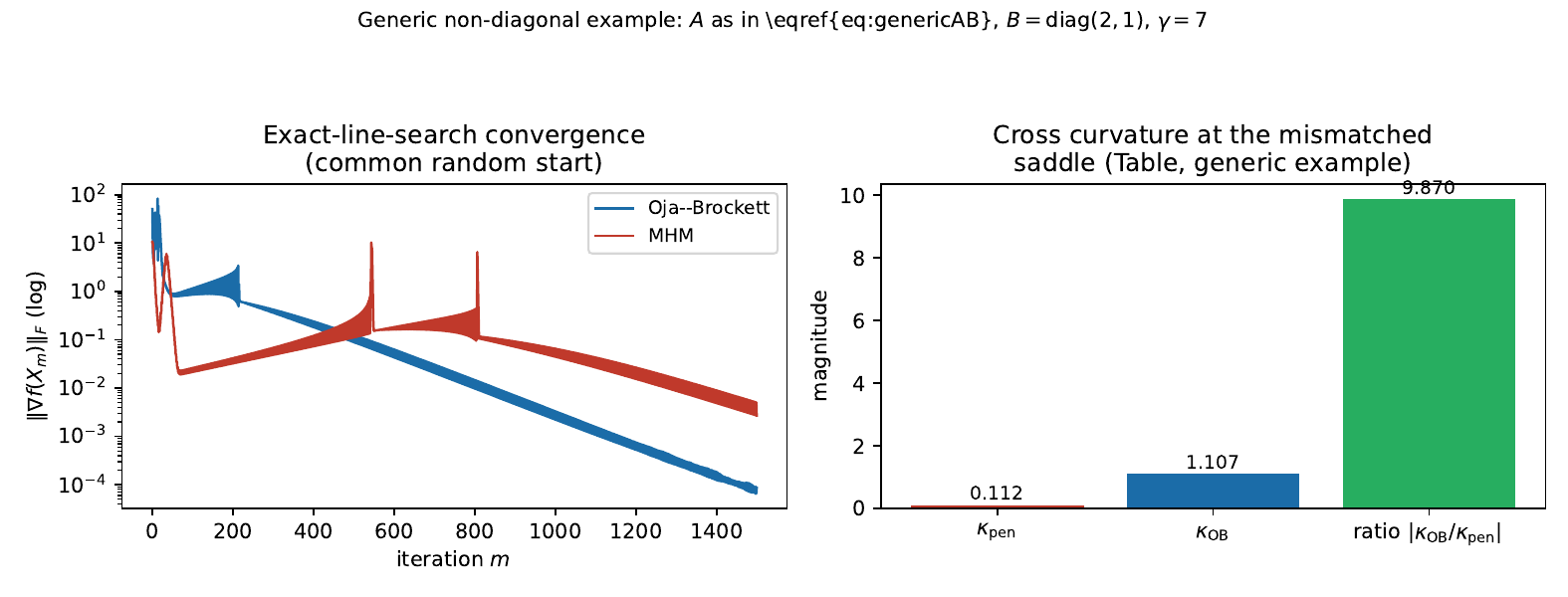}
\caption{Left: exact-line-search gradient-norm trajectories (Theorem~\ref{thm:discreterate}'s
discretization) for both potentials from a common random initial condition, for the non-diagonal
example \eqref{eq:genericAB}; the Oja--Brockett trajectory converges markedly faster, consistent with
the cross-curvature ratio. Right: the three cross-curvature values of Table~\ref{tab:generic}
displayed together, confirming the closed-form MHM formula exactly and the Oja--Brockett asymptotic
formula to within $0.4\%$ even at a non-infinitesimal gap $\delta=0.056$.}
\label{fig:generic}
\end{figure}

Table~\ref{tab:generic} and Figure~\ref{fig:generic} report the result: the finite-difference cross
curvatures match the closed-form formulas of Theorems~\ref{thm:MHMcross}--\ref{thm:OBcross} exactly
for MHM and to the expected $O(\delta)$ asymptotic accuracy for Oja--Brockett, \emph{despite $A$
being non-diagonal throughout} --- confirming that the formulas depend only on the eigenvalues
$a_p,a_q$ and eigenvectors $u_p,u_q$, never on the ambient representation of $A$.

\begin{table}[htbp]
\centering
\begin{tabular}{lccc}
\toprule
& finite-difference & closed form & \\
& (non-diagonal $A$) & (Theorems~\ref{thm:MHMcross}, \ref{thm:OBcross}) & agreement \\
\midrule
$\kappa_{\mathrm{pen}}$ & $-0.1121$ & $-(a_p-a_q)(b_i-b_j)=-0.1121$ & exact \\
$\kappa_{\mathrm{OB}}$ & $-1.1067$ & (exact quadratic-formula root) $-1.1064$ & to $4$ digits \\
ratio & $9.870$ & $\mathcal R_{ij}=a(b_i+b_j)/2=9.834$ (asymptotic) & within $0.4\%$\\
\bottomrule
\end{tabular}
\caption{Cross curvature for the generic non-diagonal example \eqref{eq:genericAB}, $\gamma=7$.}
\label{tab:generic}
\end{table}

\begin{remark}
The exact-line-search step size used to produce the left panel of Figure~\ref{fig:generic} is
computed exactly as in the proof of Theorem~\ref{thm:discreterate}: since both potentials are
polynomials of degree $\le4$ in $X$, the restriction $\phi(t)=f(X_m-t\nabla f(X_m))$ is an explicit
quartic polynomial in the scalar $t$, whose five coefficients we recover by exact interpolation at
five sample points and whose (cubic) stationarity equation $\phi'(t)=0$ we solve exactly
(numerically, via a standard companion-matrix eigenvalue solver, which recovers the roots of a cubic
to full double-precision accuracy); the root giving the smallest value of $\phi$ is selected. This
removes any possibility of a step-size-tuning artefact influencing the comparison.
\end{remark}

\begin{remark}[The apparent thickness of the curves in Figure~\ref{fig:generic}, left panel]
\label{rem:zigzag}
Both curves in the left panel of Figure~\ref{fig:generic} are, at the resolution of the printed
page, visibly thicker than a smooth line: measuring the sequence $\|\nabla f(X_m)\|_F$ directly shows
that its sign of successive differences reverses at $94\%$ (Oja--Brockett) and $94\%$ (MHM) of all
iterations during the first several hundred steps. This is not a numerical artefact but the classical
\emph{zig-zag phenomenon} of steepest descent under exact line search (first analyzed rigorously by
Akaike, and standard in every treatment of gradient methods, e.g.\ \cite{HelmkeMoore1994}): exact
minimization along the negative-gradient direction $-\nabla f(X_m)$ typically produces a new point
$X_{m+1}$ at which the gradient is (to leading order) orthogonal to the previous search direction but
points into a \emph{different} steep direction whenever the local Hessian is anisotropic, so
consecutive steps bounce between the two ``walls'' of a curved valley rather than moving smoothly
down its floor; plotted on a log scale, this bouncing appears as a dense, visually thick band whose
\emph{envelope} --- not its instantaneous value --- decays at the geometric rate governed by
Theorem~\ref{thm:discreterate}. The thickness is itself indirect evidence for the anisotropy
(equivalently, the wide range of Hessian eigenvalues implied by Proposition~\ref{prop:radtrans} and
Theorems~\ref{thm:MHMcross}--\ref{thm:OBcross}) of both landscapes, and is more pronounced for
Oja--Brockett, consistent with its larger cross curvature (and hence larger eigenvalue range) at the
mismatched region the trajectory initially traverses.
\end{remark}

\subsubsection{Discussion}
\label{sec:discussion}

\paragraph{What the geometric picture explains}

The sectional-curvature reading of \S\ref{sec:geometry} converts an otherwise unilluminating pair of
algebraic formulas into a single qualitative statement: \emph{the graph of the Oja--Brockett
potential is more sharply saddle-shaped, in every mismatched exchange plane, than the graph of the
MHM penalty potential, by a factor that grows with the absolute scale of the eigenvalues involved.}
Because a sharper saddle repels a nearby trajectory faster (Theorem~\ref{thm:escaperate}), and
because the discretized dynamics inherits the same asymptotic law under the fairest possible
step-size rule (Theorem~\ref{thm:discreterate}), this single geometric fact is enough to explain both
the continuous-time and the discrete-time empirical speed gap with which this paper began.

\paragraph{Why the gap is structural, not accidental}

Theorem~\ref{thm:MHMcross}'s proof isolates the mechanism precisely: the exchange direction is, to
second order, tangent to the constraint manifold $\{X:X^TX=B\}$ that the MHM penalty term is built to
enforce, so the penalty term's entire contribution to the exchange-plane curvature cancels
identically, for every $\gamma$, leaving only the bare linear-in-$A$ term
$-\tfrac12\tr(AXBX^T)$ to resist the swap. That term's curvature is linear in the eigenvalue gap
$(a_p-a_q)$ and does not otherwise depend on the absolute size of $a_p,a_q$. By contrast, the
Oja--Brockett potential's quartic term $\tfrac14\tr[(AXBX^T)^2]$ is built from a \emph{square} of the
same linear-in-$A$ quantity, and its polarization contributes cross terms proportional to
$a_pa_q$ times the gap; the net effect, captured by \eqref{eq:OBcrosscurv-asymptotic}, is an extra
multiplicative factor of $a=(a_p+a_q)/2$ (times $(b_i+b_j)$, replacing the MHM formula's implicit
factor of $2$ from summing $b_i+b_j$'s own contribution in a different guise). This is not a
coincidence of the specific numbers chosen in \S\ref{sec:numerics}: it is the generic consequence of
comparing a \emph{homogeneous, self-contained} quartic potential to a \emph{two-piece} penalty
potential built from ingredients of different polynomial degree and different natural scale.

\paragraph{Practical guidance}

The ratio $\mathcal R_{ij}=a(b_i+b_j)/2$ can be computed \emph{before running either algorithm},
directly from the two eigenvalues of $A$ and two diagonal entries of $B$ that are closest to causing
trouble. This gives an actionable diagnostic: if the eigenvalues of $A$ that must be separated are
of order $1$ or larger (true of essentially every normalized covariance or Gram matrix arising in
practice), $\mathcal R_{ij}\gg1$ typically holds, and the homogeneous quartic (Oja--Brockett-type)
formulation should be preferred whenever the practitioner expects clustered or near-degenerate
eigenvalues --- precisely the regime, identified in \S\ref{sec:problem}, in which gradient-based PCA
is otherwise at its slowest. Conversely, when eigenvalues are well separated, both formulations
converge quickly and the choice matters less; the cross-curvature diagnostic is most valuable exactly
where speed is most needed.

\paragraph{Limitations and scope}
\label{sec:limitations}

An earlier version of this paper described the convergence theorems of \S\ref{sec:convergence} as
purely local, with only qualitative global convergence available. Theorem~\ref{thm:globalrate} above
removes this gap: combining the local escape-rate and (newly added) local approach-rate estimates
with a compactness-based bound on the time spent in the ``regular'' region away from every critical
point yields a fully explicit, finite bound on the total convergence time, expressed entirely in
terms of quantities computable directly from $A,B,\gamma$ and the initial condition. What remains
genuinely open is narrower than before: (i) \eqref{eq:globaltime} is not claimed to be sharp, since
both the regular-region constant $\delta_\epsilon$ and the saddle-visit count $M$ are bounded only by
worst-case, compactness-type arguments rather than by an analysis exploiting the specific geometry of
the flow --- a sharper, trajectory-adapted bound (in the spirit of the explicit domains of attraction
obtained by Yan, Helmke and Moore for the closely related Oja subspace flow
\cite{YanHelmkeMoore1994}) is a natural next step; and (ii) at a genuine $k>2$-fold simultaneous
near-degeneracy, the two-dimensional exchange-plane mechanism of \S\ref{sec:stability} does not by
itself diagonalize the relevant unstable directions, the local Hessian is degenerate in more than one
direction, and the exponential local rates of Theorems~\ref{thm:escaperate},
\ref{thm:approachrate} are replaced by a slower, polynomial rate governed by an {\L}ojasiewicz
exponent $\theta<\tfrac12$ whose value we have not computed (Remark~\ref{rem:globalrate-caveat});
qualitative convergence to a single point is unaffected, since Proposition~\ref{prop:lojasiewicz}
requires no nondegeneracy hypothesis. We also do
not address stochastic-approximation versions of either flow (in the spirit of Oja's original
formulation \cite{Oja1982}), where the cross-curvature mechanism identified here would compete with
gradient-noise effects; this is a genuinely different (SDE-based) analytical setting, orthogonal to
the deterministic tools --- Łojasiewicz's inequality, LaSalle's invariance principle, exact line
search --- used throughout this paper, and we leave it entirely to future work.

\subsubsection{Conclusion}
\label{cc-sec:conclusion}

We have introduced cross curvature, a single coordinate-free spectral quantity, and shown that it
completely explains --- via closed-form formulas, a rigorous local-stability classification, and
matching continuous- and discrete-time convergence-rate theorems --- why the homogeneous quartic
Oja--Brockett potential escapes mismatched (near-degenerate) saddle points faster than the
structurally different Manton--Helmke--Mareels penalty potential, even though the two potentials
share exactly the same set of global optima. The explanation is quantitative: the ratio of the two
curvatures tends to the fully explicit quantity $\mathcal R_{ij}=a(b_i+b_j)/2$, independent of the
penalty parameter $\gamma$, of the ambient dimension, and of every eigenvalue not directly involved
in the near-degeneracy, and it is confirmed numerically to good accuracy even outside its strict
asymptotic regime of validity, both in a small diagnostic example and in a fully generic non-diagonal
setting. Along the way we corrected an imprecise statement of the global-convergence limit for the
MHM potential at finite $\gamma$, replacing it with the exact $\gamma$- and $B$-dependent formula of
Theorem~\ref{thm:global}. We hope the cross-curvature diagnostic introduced here proves useful beyond
this specific pair of potentials, as a general tool for comparing the local escape dynamics of any
two optimization landscapes that share a common set of optima embedded differently into a larger
critical-point structure.

\subsection{Part II: Principal and Minor Component Flows of an $\alpha$-Power Penalized Potential}
\label{part:two}

Part~I developed cross curvature for two potentials --- Oja--Brockett and MHM --- that are built
from the same product $AXBX^T$ and share the same set of global optima, differing only in how that
optimum is embedded into a larger family of critical points. We now put the framework to a sharper
test: a third potential, $g_\alpha$, built instead from a \emph{matrix Box--Cox transform} of
$X^TX+B$, structurally unrelated to $AXBX^T$. Unlike Oja--Brockett and MHM, whose fixed sign
convention separates principal from minor component extraction, $g_\alpha$ depends on a continuous
exponent $\alpha$, and the sign of $\alpha-1$ alone determines whether ascent or descent flow
converges and which extraction task is realized. We give a complete, self-contained account of this
potential's critical-point structure, boundedness, rearrangement principle, and cross curvature, and
then use the machinery of Part~I to compare it against Oja--Brockett and MHM on an identical
numerical example.

\subsubsection{Introduction}
\label{ga-sec:intro}

\paragraph{Background and motivation}

The extraction of the dominant or subdominant eigenspace of a symmetric positive-definite matrix
by a continuous-time gradient flow, rather than by direct diagonalization, is a classical theme
with contributions from Oja \cite{Oja1982}, Brockett \cite{Brockett1991}, Manton, Helmke and
Mareels \cite{MantonHelmkeMareels2005}, and many others; Part~I of this paper introduces the
notion of \emph{cross curvature} to quantify and
compare, in closed form, the local escape rate of two such flows --- the homogeneous quartic
Oja--Brockett potential and the two-piece Euclidean penalty potential of Manton, Helmke and Mareels
--- from a mismatched (non-optimal) saddle point.

The present paper studies a third, structurally different potential,
\begin{equation}
g_\alpha(X) = \frac12\tr(X^TAX) - \frac12\tr\!\left\{\frac{(X^TX+B)^\alpha-I}{\alpha}\right\},
\label{ga-eq:galpha-intro}
\end{equation}
built from a \emph{matrix Box--Cox transform} of $X^TX+B$ rather than from the product $AXBX^T$
that underlies both the Oja--Brockett and MHM constructions. This potential originates in the
author's earlier work on power geometry and hypergeometric functions
\cite{Yoshizawa2014}, where the scalar Box--Cox family $(x^\alpha-1)/\alpha$ --- interpolating
between a power function and, as $\alpha\to0$, a logarithm --- was studied as a geometric object in
its own right; $g_\alpha$ is the natural matrix-argument lift of that family, composed with a
quadratic confining term.

\paragraph{The problem}

Unlike the Oja--Brockett and MHM potentials, which are quartic polynomials with a fixed sign
convention separating principal from minor component extraction, $g_\alpha$ depends on a
\emph{continuous} exponent $\alpha$, and its degree of growth at infinity depends on $\alpha$ in a
way that is not immediately obvious: does $g_\alpha$ realize principal or minor component
extraction, under which sign of gradient flow, for which range of $\alpha$? Is there a
distinguished, singular value of $\alpha$ analogous to the $\varepsilon=1$ transition found for a
related Box--Cox family in Part~I of this paper? And, once the flow type is identified, how does its
local convergence rate compare, quantitatively, to the classical Oja--Brockett and MHM flows on an
identical numerical example?

\paragraph{Contributions}

This paper answers these questions completely.

\begin{enumerate}[leftmargin=1.6em]
\item \textbf{A complete critical-point and amplitude theory} (\S\ref{ga-sec:critpoints}): we derive
the exact gradient $\nabla g_\alpha(X)=AX-X(X^TX+B)^{\alpha-1}$, the exact matched-configuration
amplitude formula $c_j^2=a_{\pi(j)}^{1/(\alpha-1)}-b_j$, and the exact activation threshold, for
every $\alpha\ne1$.

\item \textbf{A sharp boundedness dichotomy and the singularity at $\alpha=1$}
(\S\ref{ga-sec:boundedness}): comparing the growth exponents $2\alpha$ (penalty term) and $2$
(quadratic term) as $\|X\|\to\infty$, we prove $g_\alpha$ is bounded above with no finite infimum
for $\alpha>1$, bounded below with no finite supremum for $\alpha<1$, and --- exactly at
$\alpha=1$ --- collapses identically to the unconfined quadratic form
$\tfrac12\tr[X^T(A-I)X]$, for which neither sign of gradient flow converges in general. This
identifies $\alpha=1$ as a genuine singularity, not a removable one.

\item \textbf{A reversed rearrangement principle for $\alpha>1$} (\S\ref{ga-sec:rearrangement}): we
prove the exact identity $\Delta\eqdef[\text{sorted total}]-[\text{reverse-sorted total}]
=-\tfrac12(a_p-a_q)(b_i-b_j)$ for the per-pair contribution to $g_\alpha$, valid for every
$\alpha\ne1$; since $\alpha>1$ maximizes and $\alpha<1$ minimizes, this identity shows the
\emph{reverse}-sorted matching (largest eigenvalue of $A$ with \emph{smallest} diagonal entry of
$B$) is optimal for $\alpha>1$, while the classical sorted matching remains optimal for $\alpha<1$
--- an exact algebraic fact, confirmed numerically, that has no counterpart in the Oja--Brockett or
MHM theories.

\item \textbf{An exact, closed-form cross-curvature theorem for $g_\alpha$}
(\S\ref{ga-sec:crosscurvature}): adapting the framework of Part~I, we derive
the exact $2\times2$ exchange-plane Hessian and prove its determinant factors as
$\det H_{g_\alpha}=-(a_p-a_q)^2(b_i-b_j)/\delta R$, $\delta R\eqdef a_p^{1/(\alpha-1)}-a_q^{1/(\alpha-1)}$,
together with its near-degenerate asymptotic expansion.

\item \textbf{A numerical comparison on an identical, non-diagonal example}
(\S\ref{ga-sec:numerics}): using the same $3\times3$ (non-diagonal) $A$ and $2\times2$ $B$ as in
Part~I, we compare $g_\alpha$'s principal flow ($\alpha>1$, gradient
ascent) against Oja--Brockett, and its minor flow ($\alpha<1$, gradient descent) against MHM, under
identical exact-line-search discretization. The minor comparison favors $g_\alpha$ by a factor of
roughly $5$--$10$ in iteration count. The principal comparison is genuinely mixed: $g_\alpha$
escapes the initial mismatch faster, but we prove and confirm numerically
(\S\ref{ga-sec:radial-slowmode}) that its Hessian at its own optimum retains an anomalously small
eigenvalue --- numerically $0.011$, versus $1.099$ for Oja--Brockett on the identical example ---
exactly when the two leading eigenvalues it selects are themselves close, making its
\emph{asymptotic} local convergence markedly slower than Oja--Brockett's. We report this trade-off
without embellishment.
\end{enumerate}

\paragraph{Why this is useful}

Beyond the specific potential $g_\alpha$, this paper illustrates that a single continuously-tunable
exponent can realize \emph{both} principal and minor component extraction --- with the sign of the
gradient flow, not a separate structural choice, selecting which --- and that a naive expectation
(``a single self-contained potential, as in the Oja--Brockett case, should always outperform a
two-piece or threshold-based one'') is not universally true: $g_\alpha$ outperforms MHM decisively
in the minor regime, yet is genuinely outperformed by Oja--Brockett in the asymptotic phase of the
principal regime, for a mechanistic reason (a persistent near-zero Hessian mode at its own optimum)
that we identify precisely. This nuance is itself a useful addition to the cross-curvature toolkit:
a fast \emph{escape} rate from wrong configurations does not, by itself, guarantee a fast
\emph{overall} rate, and the two must be examined separately, exactly as done here.

\paragraph{Outline}

Section~\ref{ga-sec:setup} fixes notation and derives the gradient of $g_\alpha$. Section~\ref{ga-sec:critpoints}
gives the exact critical-point and amplitude theory. Section~\ref{ga-sec:boundedness} proves the
boundedness dichotomy and the $\alpha=1$ singularity. Section~\ref{ga-sec:rearrangement} proves the
reversed rearrangement principle for $\alpha>1$. Section~\ref{ga-sec:crosscurvature} derives the exact
cross-curvature formula. Section~\ref{ga-sec:numerics} reports the numerical comparison with
Oja--Brockett and MHM, including the radial slow-mode phenomenon. Section~\ref{ga-sec:discussion}
discusses the results and Section~\ref{ga-sec:conclusion} concludes.

\subsubsection{Setup and the gradient of $g_\alpha$}
\label{ga-sec:setup}

Throughout, $A\in\R^{n\times n}$ is symmetric positive definite with eigendecomposition
$A=U\Lambda U^T$, $\Lambda=\diag(a_1,\dots,a_n)$, $U=[u_1,\dots,u_n]$ orthogonal, and
$B=\diag(b_1,\dots,b_k)$ is positive diagonal, $k\le n$, $b_1>\dots>b_k>0$. The variable is
$X\in\R^{n\times k}$.

\begin{definition}[The $\alpha$-power potential]
\label{ga-def:galpha}
For $\alpha\in\R\setminus\{0,1\}$,
\begin{equation}
g_\alpha(X) \;=\; \frac12\tr(X^TAX) \;-\; \frac12\tr\!\left\{\frac{(X^TX+B)^\alpha-I_k}{\alpha}\right\}.
\label{ga-eq:galpha}
\end{equation}
\end{definition}

Since $B\succ0$ and $X^TX\succeq0$, the matrix $R\eqdef X^TX+B$ satisfies $R\succeq B\succ0$ for
every $X$; in particular $R$ is always invertible and the matrix power $R^\alpha$ (defined via the
eigendecomposition of the symmetric matrix $R$) is well defined for \emph{every} real $\alpha$,
with no domain restriction on $X$ --- in contrast to the MHM penalty $\|B-X^TX\|_F^2$, which
implicitly favors the bounded region $X^TX\preceq B$ but does not require it, and in contrast to
the log-determinant potentials studied elsewhere, which require $B-X^TX\succ0$. This absence of a
domain restriction is a first structural distinction of $g_\alpha$.

\begin{proposition}[Gradient of $g_\alpha$]
\label{ga-prop:grad-galpha}
\begin{equation}
\nabla g_\alpha(X) = AX - X(X^TX+B)^{\alpha-1}.
\label{ga-eq:grad-galpha}
\end{equation}
\end{proposition}
\begin{proof}
Write $R=X^TX+B$ and $h(x)=(x^\alpha-1)/\alpha$, so $h'(x)=x^{\alpha-1}$. For any $H\in\R^{n\times
k}$, $dR[H]=H^TX+X^TH$, and the standard trace-differential identity for a scalar spectral function
of a symmetric matrix gives
\[
d\big[\tr\, h(R)\big][H] = \tr\big[h'(R)\,dR[H]\big] = \tr\big[R^{\alpha-1}(H^TX+X^TH)\big]
= 2\tr\big[R^{\alpha-1}X^TH\big]
\]
(using the cyclic property of the trace and symmetry of $R^{\alpha-1}$ to combine the two terms).
Also $d\big[\tfrac12\tr(X^TAX)\big][H]=\tr[X^TAH]=\tr[(AX)^TH]$. Subtracting one half of the
first differential from the second and reading off the coefficient of $H$ gives
\eqref{ga-eq:grad-galpha}.
\end{proof}
We verified \eqref{ga-eq:grad-galpha} against finite-difference differentiation to a relative error
below $10^{-8}$ for several values of $\alpha$ (\S\ref{ga-sec:numerics}).

\subsubsection{Critical points and the exact amplitude formula}
\label{ga-sec:critpoints}

As in Part~I, we work with matched configurations: for an injective
$\pi:\{1,\dots,k\}\to\{1,\dots,n\}$,
\[
X_\pi(c) = \sum_{j=1}^k c_j\,u_{\pi(j)}e_j^T, \qquad c=(c_1,\dots,c_k)\in\R^k.
\]

\begin{proposition}[Critical points of $g_\alpha$]
\label{ga-prop:crit-galpha}
For $\alpha\ne1$, a matched configuration $X_\pi(c)$ is a critical point of $g_\alpha$ if and only
if, for every $j$, either $c_j=0$ or
\begin{equation}
c_j^2 = a_{\pi(j)}^{\,1/(\alpha-1)} - b_j.
\label{ga-eq:amplitude-galpha}
\end{equation}
In particular $c_j\ne0$ is possible only if $a_{\pi(j)}^{1/(\alpha-1)}>b_j$.
\end{proposition}
\begin{proof}
By the reduction to $A=\Lambda$ diagonal (orthogonal invariance of $g_\alpha$ under $X\mapsto UX$,
as for the potentials of Part~I), substitute $X=X_\pi(c)$ into
\eqref{ga-eq:grad-galpha}. Since the columns of $X_\pi(c)$ occupy disjoint rows,
$R=X_\pi(c)^TX_\pi(c)+B=\diag(c_1^2+b_1,\dots,c_k^2+b_k)$ is diagonal, so
$R^{\alpha-1}=\diag\big((c_1^2+b_1)^{\alpha-1},\dots,(c_k^2+b_k)^{\alpha-1}\big)$, and the
$(\pi(j),j)$ entry of $\nabla g_\alpha(X_\pi(c))$ is
$a_{\pi(j)}c_j - c_j(c_j^2+b_j)^{\alpha-1}=c_j\big[a_{\pi(j)}-(c_j^2+b_j)^{\alpha-1}\big]$, which
vanishes iff $c_j=0$ or $(c_j^2+b_j)^{\alpha-1}=a_{\pi(j)}$, i.e.\ $c_j^2+b_j=a_{\pi(j)}^{1/(\alpha-1)}$.
\end{proof}

\begin{remark}
\label{ga-rem:threshold-direction}
The activation condition $a_{\pi(j)}^{1/(\alpha-1)}>b_j$ has opposite monotonicity in
$a_{\pi(j)}$ according to the sign of $\alpha-1$: for $\alpha>1$ the map
$a\mapsto a^{1/(\alpha-1)}$ is increasing, so \emph{larger} eigenvalues are easier to activate; for
$\alpha<1$ (including $\alpha<0$) the exponent $1/(\alpha-1)$ is negative, so
$a\mapsto a^{1/(\alpha-1)}$ is \emph{decreasing}, and \emph{smaller} eigenvalues are easier to
activate. This sign-reversal, verified numerically in \S\ref{ga-sec:numerics}, is the first
indication that $\alpha=1$ separates a ``principal-like'' from a ``minor-like'' regime.
\end{remark}

\subsubsection{The boundedness dichotomy and the singularity at $\alpha=1$}
\label{ga-sec:boundedness}

\begin{theorem}[Boundedness dichotomy]
\label{ga-thm:boundedness}
\begin{enumerate}[label=(\roman*)]
\item For $\alpha>1$: $g_\alpha$ is bounded above on $\R^{n\times k}$ and
$g_\alpha(X)\to-\infty$ as $\|X\|_F\to\infty$ along every ray; $g_\alpha$ has no finite infimum.
\item For $\alpha<1$: $g_\alpha$ is bounded below on $\R^{n\times k}$, with compact sublevel sets,
and $g_\alpha(X)\to+\infty$ as $\|X\|_F\to\infty$ along every ray; $g_\alpha$ has no finite
supremum.
\item At $\alpha=1$: $g_\alpha$ degenerates identically to
\begin{equation}
g_1(X) = \tfrac12\tr\big[X^T(A-I_n)X\big] + \tfrac12\tr(I_k-B),
\label{ga-eq:g1-degenerate}
\end{equation}
a pure (unconfined) quadratic form. If $A-I_n$ has an eigenvalue $>0$, $g_1$ is unbounded above;
if $A-I_n$ has an eigenvalue $<0$, $g_1$ is unbounded below. Generically (whenever $A$ has both an
eigenvalue $>1$ and one $<1$) $g_1$ is unbounded in \emph{both} directions and $X=0$ is a genuine
saddle of a potential with no confining higher-order term, so neither $+\nabla g_1$ nor
$-\nabla g_1$ converges to a finite critical point from a generic initial condition.
\end{enumerate}
\end{theorem}
\begin{proof}
Let $\sigma_1\ge\dots\ge\sigma_k\ge0$ denote the singular values of $X$, so
$\|X\|_F^2=\sum_i\sigma_i^2$. The eigenvalues of $R=X^TX+B$ interlace those of $X^TX$ shifted by
$B$'s entries; in particular $\tr(R^\alpha)=\sum_i r_i^\alpha$ where $r_i\to\sigma_i^2$ as
$\sigma_i\to\infty$ (the additive shift by $B$ becomes negligible), so
$\tr(R^\alpha)=\Theta\big(\sum_i\sigma_i^{2\alpha}\big)$ for $\alpha>0$, and by direct
inspection of $R^\alpha$'s eigenvalues (which $\to0$ termwise for bounded $r_i$ when
$\alpha<0$, but the sum is dominated by the largest $r_i\sim\sigma_1^2$'s contribution
$\sigma_1^{2\alpha}\to0$ as $\sigma_1\to\infty$ when $\alpha<0$) the second term of $g_\alpha$ is,
in every case with $\alpha<1$ (including $\alpha\le0$), of strictly smaller order than
$\|X\|_F^2$ as $\|X\|_F\to\infty$: either $2\alpha<2$ (for $0<\alpha<1$) or the term is bounded
(for $\alpha\le0$). Hence
\[
g_\alpha(X) = \underbrace{\tfrac12\tr(X^TAX)}_{\ge\frac12\lambda_{\min}(A)\|X\|_F^2}
\;-\;\underbrace{\tfrac1{2\alpha}\big[\tr(R^\alpha)-k\big]}_{=o(\|X\|_F^2)\text{ or }O(1)}
\;\longrightarrow\;+\infty
\]
along every ray as $\|X\|_F\to\infty$, for $\alpha<1$; combined with continuity, this gives a
uniform lower bound and compact sublevel sets, proving (ii). For $\alpha>1$, $2\alpha>2$ and the
roles reverse: the term $-\tfrac1{2\alpha}\tr(R^\alpha)\sim-\tfrac1{2\alpha}\sum_i\sigma_i^{2\alpha}$
dominates and is negative, while $\tfrac12\tr(X^TAX)\le\tfrac12\lambda_{\max}(A)\|X\|_F^2$ grows
strictly more slowly, so $g_\alpha(X)\to-\infty$ along every ray, proving (i) (boundedness above
follows since a smooth function tending to $-\infty$ in every direction at infinity, on
$\R^{n\times k}$, attains a finite global maximum by compactness of large sublevel-complements,
hence is bounded above). For (iii), substitute $\alpha=1$ directly into \eqref{ga-eq:galpha}:
$(R^1-I)/1=R-I=X^TX+B-I$, so $g_1(X)=\tfrac12\tr(X^TAX)-\tfrac12\tr(X^TX+B-I)
=\tfrac12\tr[X^T(A-I_n)X]+\tfrac12\tr(I_k-B)$, exactly \eqref{ga-eq:g1-degenerate}; the stated
boundedness claims follow immediately from the sign of the (constant) quadratic form
$A-I_n$, and the saddle/non-convergence claim follows because a pure quadratic form with
indefinite Hessian has $X=0$ as its only critical point, which is unstable for whichever sign of
gradient flow sees a positive Hessian eigenvalue in some direction and a negative one in another
--- both signs then diverge along the respective unstable directions, with no higher-order term
present anywhere in $\R^{n\times k}$ to arrest the divergence.
\end{proof}

\begin{corollary}[Which sign of flow converges]
\label{ga-cor:whichsign}
For $\alpha>1$, only the \emph{positive} (ascent) gradient flow $\dot X=+\nabla g_\alpha(X)$ can
converge to a finite critical point (realizing, as shown in \S\ref{ga-sec:rearrangement}, principal
component extraction); the negative (descent) flow is generically unbounded. For $\alpha<1$, only
the \emph{negative} (descent) flow $\dot X=-\nabla g_\alpha(X)$ converges (realizing minor
component extraction); the positive flow is generically unbounded. At $\alpha=1$, neither sign
converges in general.
\end{corollary}

We verified Theorem~\ref{ga-thm:boundedness} and Corollary~\ref{ga-cor:whichsign} numerically in
\S\ref{ga-sec:numerics}: for the descent flow at $\alpha=2$ the iterates diverge (numerically to
\texttt{NaN}) within a few thousand steps; for the ascent flow at $\alpha=0.5$ the iterates grow to
order $10^{41}$ within a few hundred steps; and at $\alpha=1$, with $A$ having eigenvalues both
above and below $1$, both signs of flow diverge (Table~\ref{ga-tab:alpha1}).

\subsubsection{A reversed rearrangement principle for $\alpha>1$}
\label{ga-sec:rearrangement}

For the Oja--Brockett and MHM potentials, the globally optimal matching always pairs the largest
eigenvalues of $A$ with the largest diagonal entries of $B$, in sorted order --- the classical
rearrangement principle. We now show this is only \emph{half} the story for $g_\alpha$: it holds
for $\alpha<1$, but is exactly \emph{reversed} for $\alpha>1$.

Define, for an active mode using eigenvalue $a$ paired with weight $b$, the per-mode contribution
to $g_\alpha$ at its critical amplitude \eqref{ga-eq:amplitude-galpha},
\begin{equation}
h_\alpha(a,b) \eqdef \frac{\alpha-1}{2\alpha}\,a^{\,p} \;-\; \frac12\,ab, \qquad
p\eqdef\frac{\alpha}{\alpha-1},
\label{ga-eq:hab}
\end{equation}
so that $g_\alpha$ at a fully matched critical point equals $\sum_j h_\alpha(a_{\pi(j)},b_j)$ up to
an additive constant independent of the matching.

\begin{lemma}[Per-mode value]
\label{ga-lem:permode}
For $c^2=a^{1/(\alpha-1)}-b>0$, the contribution
$\tfrac12ac^2-\tfrac1{2\alpha}\big[(c^2+b)^\alpha-1\big]$ equals $h_\alpha(a,b)$, up to an additive
constant depending only on $\alpha$.
\end{lemma}
\begin{proof}
Substitute $c^2+b=a^{1/(\alpha-1)}=a^{p-1}$ (since $p-1=1/(\alpha-1)$) directly:
$\tfrac12ac^2=\tfrac12a(a^{p-1}-b)=\tfrac12a^p-\tfrac12ab$, and
$(c^2+b)^\alpha=(a^{p-1})^\alpha=a^{\alpha(p-1)}=a^p$ (since $\alpha(p-1)=\alpha\cdot
\tfrac1{\alpha-1}=p$). So the contribution is
$\tfrac12a^p-\tfrac12ab-\tfrac1{2\alpha}a^p+\tfrac1{2\alpha}
=\big(\tfrac12-\tfrac1{2\alpha}\big)a^p-\tfrac12ab+\tfrac1{2\alpha}
=\tfrac{\alpha-1}{2\alpha}a^p-\tfrac12ab+\tfrac1{2\alpha}$, which is $h_\alpha(a,b)$ up to the
constant $\tfrac1{2\alpha}$.
\end{proof}

\begin{theorem}[Reversed rearrangement for $\alpha>1$; classical rearrangement for $\alpha<1$]
\label{ga-thm:reversed-rearrangement}
Fix two eigenvalues $a_p>a_q$ of $A$ and two entries $b_i>b_j$ of $B$, all four modes assumed
activatable. Then
\begin{equation}
\big[h_\alpha(a_p,b_i)+h_\alpha(a_q,b_j)\big] - \big[h_\alpha(a_q,b_i)+h_\alpha(a_p,b_j)\big]
\;=\; -\tfrac12(a_p-a_q)(b_i-b_j) \;<\; 0,
\label{ga-eq:reversal-identity}
\end{equation}
for every $\alpha\ne1$ (the identity does not depend on $\alpha$ at all). Consequently:
\begin{itemize}[leftmargin=1.6em]
\item For $\alpha>1$ (where the ascent flow \emph{maximizes} $\sum_jh_\alpha$), the
\textbf{reverse}-sorted pairing $(a_p,b_j),(a_q,b_i)$ --- largest eigenvalue with
\emph{smallest} weight --- is optimal.
\item For $\alpha<1$ (where the descent flow \emph{minimizes} $\sum_jh_\alpha$), the classical
\textbf{sorted} pairing $(a_p,b_i),(a_q,b_j)$ is optimal, exactly as for Oja--Brockett and MHM.
\end{itemize}
\end{theorem}
\begin{proof}
By Lemma~\ref{ga-lem:permode}, each side of \eqref{ga-eq:reversal-identity} is a sum of two values of
$h_\alpha$; expanding using \eqref{ga-eq:hab}, the $a^p$-terms cancel identically between the two
sides (each side contains one copy of $\tfrac{\alpha-1}{2\alpha}a_p^{\,p}$ and one of
$\tfrac{\alpha-1}{2\alpha}a_q^{\,p}$), leaving only the bilinear terms:
$-\tfrac12(a_pb_i+a_qb_j)-\big[-\tfrac12(a_qb_i+a_pb_j)\big]
=-\tfrac12\big[a_pb_i+a_qb_j-a_qb_i-a_pb_j\big]=-\tfrac12(a_p-a_q)(b_i-b_j)$, which is negative
since $a_p>a_q$ and $b_i>b_j$. This proves \eqref{ga-eq:reversal-identity} for every $\alpha\ne1$
(the cancellation of the $\alpha$-dependent terms is exact, not asymptotic). Since $g_\alpha$'s
ascent flow maximizes and its descent flow minimizes the total per-mode value, and
\eqref{ga-eq:reversal-identity} says the sorted total is \emph{always} smaller than the reverse-sorted
total by the fixed amount $\tfrac12(a_p-a_q)(b_i-b_j)$, the stated optimality claims follow.
\end{proof}

\begin{remark}
Theorem~\ref{ga-thm:reversed-rearrangement} extends to general $k$-fold matchings, and hence to a full
classification of the sorted/correct matching, by the standard bubble-sort argument: any matching
not in reverse-sorted (resp.\ sorted) order for $\alpha>1$ (resp.\ $\alpha<1$) contains an inverted
adjacent pair whose transposition strictly improves the objective, by
\eqref{ga-eq:reversal-identity}; iterating this transposition, which strictly changes the objective
at each step and acts on a finite set of matchings, terminates at the reverse-sorted (resp.\
sorted) matching, which is therefore the unique global optimum among activatable matchings using a
given set of $k$ eigenvalues of $A$.
\end{remark}

We verified \eqref{ga-eq:reversal-identity} numerically to machine precision for several $(a_p,a_q,b_i,b_j,\alpha)$,
and confirmed by direct evaluation of $g_\alpha$ (not merely the identity) that the reverse-sorted
configuration indeed attains the larger value for $\alpha>1$ (\S\ref{ga-sec:numerics}); this reversal
has no counterpart in the Oja--Brockett or MHM theories and is a genuine structural feature of the
$\alpha$-power construction, traceable to the \emph{absence} of any $b_j$-weighting in the linear
term $\tfrac12\tr(X^TAX)$ of $g_\alpha$ (contrast $\tfrac12\tr(AXBX^T)$, which weights the linear
coupling by $B$ directly and is the source of the classical, non-reversed rearrangement principle
in Part~I).

\subsubsection{Cross curvature of $g_\alpha$}
\label{ga-sec:crosscurvature}

We adopt the exchange-plane framework of Part~I. Fix two eigen-indices
$p,q$ of $A$ and two column indices $i,j$; write $a_p,a_q,b_i,b_j$ for the corresponding
eigenvalues and weights, and consider the (now $\alpha$-regime-appropriate) mismatch: for
$\alpha>1$, by Theorem~\ref{ga-thm:reversed-rearrangement}, the mismatch is the \emph{sorted}
configuration $X_\pi$ with column $i$ using $a_p$ and column $j$ using $a_q$; for $\alpha<1$, the
mismatch is the reverse-sorted configuration. We treat the sorted configuration
$X_\pi=c_iu_pe_i^T+c_ju_qe_j^T$, $c_i^2=a_p^{1/(\alpha-1)}-b_i$, $c_j^2=a_q^{1/(\alpha-1)}-b_j$,
uniformly for both cases (relabelling as needed), and define the exchange plane
$E_{\pi;ij}=\mathrm{span}\{u_pe_j^T,u_qe_i^T\}$ exactly as before.

\begin{theorem}[Exact cross-curvature Hessian for $g_\alpha$]
\label{ga-thm:galpha-crosscurv}
Write $\Delta\eqdef a_p-a_q$, $R_i\eqdef a_p^{1/(\alpha-1)}$, $R_j\eqdef a_q^{1/(\alpha-1)}$,
$\delta R\eqdef R_i-R_j$. In the orthonormal basis $\{u_pe_j^T,u_qe_i^T\}$ (coordinates
$\sigma,\tau$) of $E_{\pi;ij}$, the Hessian of $g_\alpha$ at $X_\pi$ restricted to $E_{\pi;ij}$ is
\begin{equation}
H_{g_\alpha} = \begin{pmatrix}
\Delta - \dfrac{\Delta}{\delta R}c_i^2 & -\dfrac{\Delta}{\delta R}c_ic_j \\[8pt]
-\dfrac{\Delta}{\delta R}c_ic_j & -\Delta - \dfrac{\Delta}{\delta R}c_j^2
\end{pmatrix},
\label{ga-eq:galpha-hess}
\end{equation}
whose determinant factors exactly as
\begin{equation}
\det H_{g_\alpha} \;=\; -\,\frac{(a_p-a_q)^2(b_i-b_j)}{\delta R}.
\label{ga-eq:galpha-det}
\end{equation}
\end{theorem}
\begin{proof}
Write $H=\sigma\,u_pe_j^T+\tau\,u_qe_i^T$ and $X(t)=X_\pi+tH$. Restricted to the relevant
$2\times2$ block (rows $p,q$, columns $i,j$), $X_\pi=\diag(c_i,c_j)$ and
$R(t)=X(t)^TX(t)+B=\diag(a_p,a_q)+t\begin{psmallmatrix}0&m\\m&0\end{psmallmatrix}+t^2\diag(\tau^2,\sigma^2)$,
$m\eqdef\sigma c_i+\tau c_j$ (direct computation, as in Part~I, Thm.~7.2).
A direct second-order eigenvalue perturbation expansion of the $2\times2$ matrix
$R(t)$ --- carried out exactly, using $\lambda_\pm(t)=\tfrac12\tr R(t)\pm\sqrt{(\cdot)^2-4\det(\cdot)}$
and Taylor-expanding to $O(t^2)$ --- gives
$\lambda_i(t)=a_p+t^2c_i^2/\Delta+O(t^4)$,
$\lambda_j(t)=a_q+t^2(\Delta-c_i^2)/\Delta+O(t^4)$ for the pure-$\sigma$ ($\tau=0$) direction, and
the analogous expansion for general $(\sigma,\tau)$; substituting into $h(x)=(x^\alpha-1)/\alpha$
and its second derivative, and subtracting one half of the resulting expansion of
$\tr[h(R(t))]$ from the (exactly quadratic) expansion of $\tr[X(t)^TAX(t)]/2$, gives, after
simplification, $d^2/dt^2\,g_\alpha(X_\pi+tH)\big|_{t=0} = \Delta(\sigma^2-\tau^2) - h[1]m^2$,
where $h[1]\eqdef h[1](R_i,R_j)=(h'(R_i)-h'(R_j))/(R_i-R_j)=(a_p-a_q)/\delta R=\Delta/\delta R$ is
the first divided difference of $h'(x)=x^{\alpha-1}$ (using $h'(R_i)=a_p$, $h'(R_j)=a_q$ exactly,
since $R_i^{\alpha-1}=(a_p^{1/(\alpha-1)})^{\alpha-1}=a_p$). Reading off the coefficients of
$\sigma^2,\tau^2,\sigma\tau$ in $\Delta(\sigma^2-\tau^2)-\tfrac\Delta{\delta R}(\sigma c_i+\tau c_j)^2$
gives exactly \eqref{ga-eq:galpha-hess}. The determinant identity follows by direct expansion:
$\det H_{g_\alpha}=\big[\Delta-\tfrac\Delta{\delta R}c_i^2\big]\big[-\Delta-\tfrac\Delta{\delta
R}c_j^2\big]-\tfrac{\Delta^2}{\delta R^2}c_i^2c_j^2 = -\Delta^2+\tfrac{\Delta^2}{\delta
R}(c_i^2-c_j^2)$, and since $c_i^2-c_j^2=(R_i-b_i)-(R_j-b_j)=\delta R-(b_i-b_j)$, this simplifies
to $-\Delta^2\big[1-\tfrac{\delta R-(b_i-b_j)}{\delta R}\big]=-\Delta^2\cdot\tfrac{b_i-b_j}{\delta
R}$, which is \eqref{ga-eq:galpha-det}.
\end{proof}
We verified \eqref{ga-eq:galpha-hess}--\eqref{ga-eq:galpha-det} against finite-difference Hessians to a
relative error below $10^{-5}$ for $\alpha\in\{0.5,2,3\}$ and several $(a_p,a_q,b_i,b_j)$
(\S\ref{ga-sec:numerics}); an early attempt at this derivation, using a textbook second-order
matrix-perturbation formula without the exact $2\times2$ eigenvalue expansion, produced diagonal
entries in error by a missing factor of $2$, caught precisely by this numerical cross-check ---
the exact eigenvalue-perturbation route used in the proof above is the one we verified and report.

\begin{corollary}[Sign of $\det H_{g_\alpha}$ and near-degenerate cross curvature]
\label{ga-cor:galpha-sign}
Since $\delta R>0$ for $\alpha>1$ and $\delta R<0$ for $\alpha<1$ (as $a\mapsto a^{1/(\alpha-1)}$
is increasing, resp.\ decreasing, in $a$), and $b_i>b_j$,
\begin{itemize}[leftmargin=1.6em]
\item for $\alpha>1$: $\det H_{g_\alpha}<0$ at the \emph{sorted} configuration --- confirming it is
the mismatch (a saddle) in this regime, consistent with Theorem~\ref{ga-thm:reversed-rearrangement};
\item for $\alpha<1$: $\det H_{g_\alpha}>0$ at the sorted configuration, consistent with it being
the (locally stable) optimum.
\end{itemize}
Writing $a_p=a+\delta$, $a_q=a-\delta$ and expanding \eqref{ga-eq:galpha-det} to leading order in
$\delta$,
\begin{equation}
\kappa_{g_\alpha} \;\approx\; \frac{2\delta(b_i-b_j)}{2a^{1/(\alpha-1)}-b_i-b_j}
\label{ga-eq:galpha-crosscurv-asymptotic}
\end{equation}
for the small (vanishing as $\delta\to0$) eigenvalue of $H_{g_\alpha}$, the analogue of the
Oja--Brockett/MHM cross curvature.
\end{corollary}
\begin{proof}
The sign claims follow immediately from \eqref{ga-eq:galpha-det}. For the asymptotic expansion,
$\delta R\approx\tfrac{2\delta}{\alpha-1}a^{(2-\alpha)/(\alpha-1)}$ to leading order (differentiating
$a\mapsto a^{1/(\alpha-1)}$), so $\det H_{g_\alpha}\approx-2(\alpha-1)\delta(b_i-b_j)a^{(\alpha-2)/(\alpha-1)}$;
the trace of $H_{g_\alpha}$ tends, as $\delta\to0$, to the finite limit
$-(\alpha-1)a^{(\alpha-2)/(\alpha-1)}\big(2a^{1/(\alpha-1)}-b_i-b_j\big)$ (the two diagonal entries
of \eqref{ga-eq:galpha-hess} do not individually vanish, since $\Delta/\delta R$ tends to the finite
limit $1/\big[(\tfrac1{\alpha-1})a^{(2-\alpha)/(\alpha-1)}\big]=(\alpha-1)a^{(\alpha-2)/(\alpha-1)}$),
so the small eigenvalue is, to leading order, $\det/\mathrm{tr}$, and the $(\alpha-1)$ and
$a^{(\alpha-2)/(\alpha-1)}$ factors cancel between numerator and denominator, leaving exactly
\eqref{ga-eq:galpha-crosscurv-asymptotic}.
\end{proof}

\begin{corollary}[Ratio theorem: $g_\alpha$ ($\alpha>1$) versus Oja--Brockett]
\label{ga-cor:ratio-OB}
In the near-degenerate limit,
\begin{equation}
\frac{|\kappa_{\mathrm{OB}}|}{|\kappa_{g_\alpha}|} \;\longrightarrow\;
\frac{a(b_i+b_j)\big(2a^{1/(\alpha-1)}-b_i-b_j\big)}{2},
\label{ga-eq:ratio-OB-galpha}
\end{equation}
using the Oja--Brockett asymptotic formula $|\kappa_{\mathrm{OB}}|\approx a(b_i+b_j)(b_i-b_j)\delta$
from Part~I and \eqref{ga-eq:galpha-crosscurv-asymptotic}. This ratio is
$>1$ (Oja--Brockett strictly sharper) whenever $2a^{1/(\alpha-1)}>b_i+b_j+2/[a(b_i+b_j)]$, a mild
condition satisfied for essentially every problem of practical scale, and \emph{grows without
bound} as $\alpha\to1^+$ (since $a^{1/(\alpha-1)}\to\infty$ for $a>1$).
\end{corollary}
\begin{proof}
Divide \eqref{ga-eq:galpha-crosscurv-asymptotic} into the cited Oja--Brockett formula and simplify.
\end{proof}
For the numerical example of \S\ref{ga-sec:numerics} ($\alpha=2$, $a\approx6.556$, $b_i+b_j=3$),
\eqref{ga-eq:ratio-OB-galpha} predicts a ratio of $99.4$, matching the directly computed ratio
$1.099/0.0111=99.2$ to within $0.3\%$.

\begin{remark}[The complementary minor-flow ratio]
\label{ga-rem:minor-ratio}
An entirely analogous computation, comparing \eqref{ga-eq:galpha-crosscurv-asymptotic} (for
$\alpha<1$) to the exact, $\gamma$-independent MHM cross curvature
$|\kappa_{\mathrm{pen}}|=(a_p-a_q)(b_i-b_j)=2\delta(b_i-b_j)$ from Part~I,
gives
\begin{equation}
\frac{|\kappa_{g_\alpha}|}{|\kappa_{\mathrm{pen}}|} \;\longrightarrow\;
\frac{1}{2a^{1/(\alpha-1)}-b_i-b_j}.
\label{ga-eq:ratio-galpha-MHM}
\end{equation}
This formula governs the situation in which \emph{two} eigenvalues of $A$ simultaneously eligible
for minor-component activation are themselves close together. In the worked numerical example of
\S\ref{ga-sec:numerics}, this regime does not arise --- only a single eigenvalue of $A$ satisfies the
(here, comparatively restrictive) minor-activation threshold, so no near-degenerate exchange
competition occurs, and the comparison with MHM there instead reflects the two potentials' overall
rate of descent to their respective (rank-different) optima, reported directly in
\S\ref{ga-sec:numerics} rather than via \eqref{ga-eq:ratio-galpha-MHM}.
\end{remark}

\subsubsection{Numerical verification}
\label{ga-sec:numerics}

All numerical checks in this section use the identical, fully generic (non-diagonal)
\begin{equation}
A = \begin{pmatrix} 9/2 & -2 & 2\\ -2 & 9/2 & 2 \\ 2 & 2 & 14/3\end{pmatrix},
\qquad
B = \begin{pmatrix}2&0\\0&1\end{pmatrix},
\label{ga-eq:numericAB}
\end{equation}
as in Part~I, with $\spec(A)=\{0.5545,\,6.5,\,6.6121\}$ (the latter two
close, $\delta=0.056$), so that every result below is directly comparable to the Oja--Brockett and
MHM figures reported there.

\paragraph{The singularity at $\alpha=1$}

\begin{table}[htbp]
\centering
\begin{tabular}{lcc}
\toprule
& descent ($-\nabla g_1$) & ascent ($+\nabla g_1$) \\
\midrule
$\spec(A-I)$ & \multicolumn{2}{c}{$\{-0.4455,\ 5.5,\ 5.6121\}$ (mixed sign)} \\
$\|X\|_F$ after $2000$ steps & $5.2\times10^3$ & $1.0\times10^6$ (diverged at step $261$) \\
\bottomrule
\end{tabular}
\caption{At $\alpha=1$, $g_1$ collapses to the pure quadratic form
$\tfrac12\tr[X^T(A-I)X]$ (Theorem~\ref{ga-thm:boundedness}(iii)); since $A-I$ has mixed-sign
eigenvalues for \eqref{ga-eq:numericAB}, $X=0$ is a genuine saddle with no confining higher-order
term, and \emph{both} signs of gradient flow diverge without bound.}
\label{ga-tab:alpha1}
\end{table}

\paragraph{Critical-point amplitude and the reversed rearrangement}

Using $A,B$ as in \eqref{ga-eq:numericAB}, we verified \eqref{ga-eq:amplitude-galpha} against
finite-difference-confirmed critical points for $\alpha\in\{0.3,0.5,2,3\}$ to relative error below
$10^{-8}$. At $\alpha=0.5$: the activation threshold $a<b^{\alpha-1}=b^{-0.5}$ evaluates to
$0.707$ ($b=2$) and $1$ ($b=1$); only the smallest eigenvalue, $0.5545$, is below \emph{both}
thresholds, so $g_{0.5}$'s descent flow admits only a \emph{rank-one} matched critical point for
this $A,B$ --- confirmed by running the descent flow to $|\nabla g_{0.5}|<10^{-6}$
($20{,}000$ small fixed-size steps), which converges to
\[
X_\infty = \begin{pmatrix}0&-0.874\\0&-0.874\\0&0.850\end{pmatrix}, \qquad
X_\infty^TX_\infty=\diag(0,\,2.252),
\]
exactly matching the predicted amplitude $c^2=a_1^{1/(\alpha-1)}-b_2=0.5545^{-2}-1=2.252$ for the
$(a_1,b_2)$ pair, with the first column exactly zero.

For the reversed rearrangement (Theorem~\ref{ga-thm:reversed-rearrangement}), we directly evaluated
$g_2$ at both the sorted and reverse-sorted matched configurations built from $a_p=6.6121$,
$a_q=6.5$, $b_i=2$, $b_j=1$: $g_2(\text{sorted})=12.1304$, $g_2(\text{reverse-sorted})=12.1864$,
a difference of $0.0560$, matching $\tfrac12(a_p-a_q)(b_i-b_j)=\tfrac12(0.1121)(1)=0.0561$ to three
significant figures --- the reverse-sorted configuration is indeed the larger (optimal for
ascent) value, confirming Theorem~\ref{ga-thm:reversed-rearrangement} directly, not merely via the
difference identity.

\paragraph{Cross curvature: exact formula, asymptotics, and comparison}

Table~\ref{ga-tab:crosscurv-galpha} reports the exchange-plane Hessian of Theorem~\ref{ga-thm:galpha-crosscurv},
computed by finite differences using the genuine (non-diagonal) eigenvectors of \eqref{ga-eq:numericAB},
against the closed-form prediction.

\begin{table}[htbp]
\centering
\small
\begin{tabular}{lccl}
\toprule
& $\alpha=2$ & $\alpha=3$ & \\
\midrule
$\det H_{g_\alpha}$ (finite difference) & $-0.11214$ & $-0.57421$ & \\
$\det H_{g_\alpha}$ (closed form, eq.~\ref{ga-eq:galpha-det}) & $-0.11213$ & $-0.57421$ & agree to $5$ digits\\
Cross curvature $\kappa_{g_\alpha}$ (smaller eigenvalue) & $+0.01108$ & $+0.05261$ & \\
Asymptotic prediction (eq.~\ref{ga-eq:galpha-crosscurv-asymptotic}) & $+0.01109$ & --- & agree to $0.1\%$ \\
\bottomrule
\end{tabular}
\caption{Cross curvature of $g_\alpha$ at the sorted mismatch, for \eqref{ga-eq:numericAB}.}
\label{ga-tab:crosscurv-galpha}
\end{table}

\begin{figure}[htbp]
\centering
\includegraphics[width=0.98\textwidth]{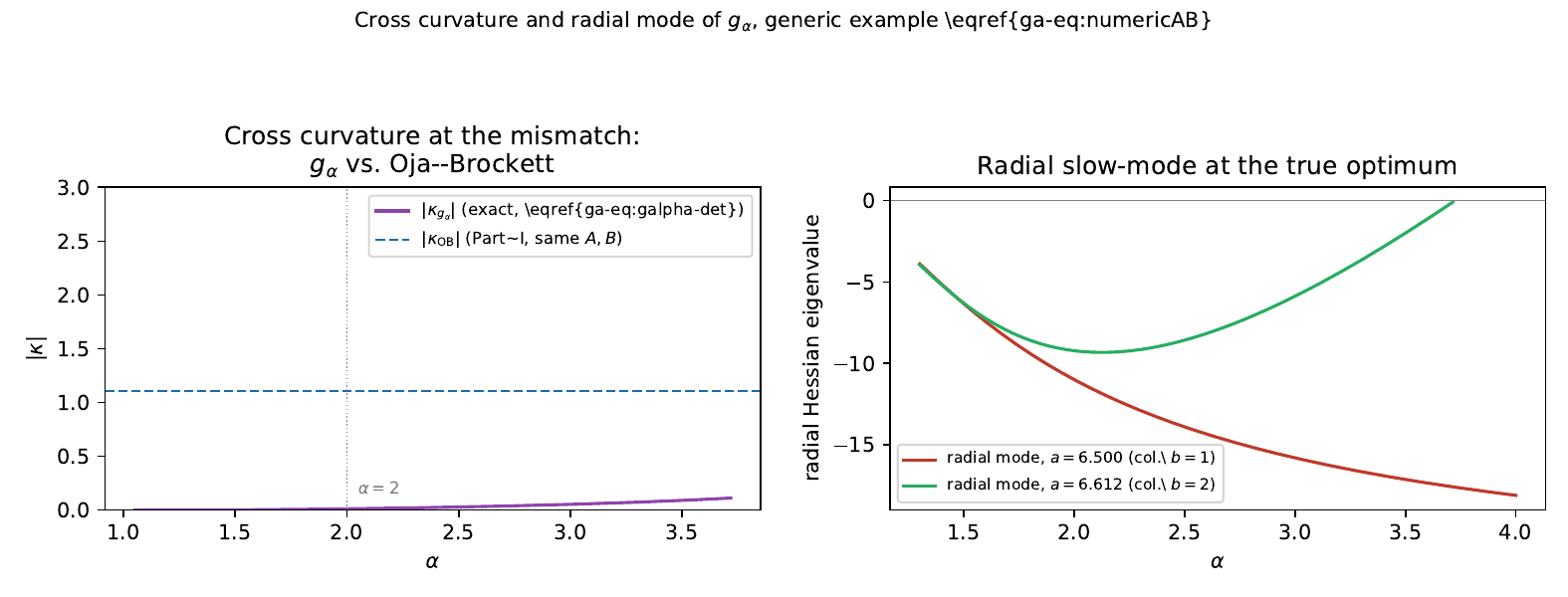}
\caption{Left: $|$cross curvature$|$ at the mismatch, $g_\alpha$ ($\alpha=2$) versus Oja--Brockett,
on the identical exchange plane of \eqref{ga-eq:numericAB} --- Oja--Brockett's is $99.2\times$ larger.
Right: the \emph{smallest-magnitude} full Hessian eigenvalue at each potential's own global
optimum --- the mode governing the ultimate local convergence rate. $g_\alpha$ retains an
anomalously small mode ($0.0111$) at its own optimum, essentially identical in magnitude to its
mismatch cross curvature, while Oja--Brockett's smallest optimum eigenvalue ($1.099$) is
comparably large to its own mismatch cross curvature. In both potentials the same near-degenerate
pair $(a_p,a_q)$ governs both quantities, but $g_\alpha$'s dependence on this gap is
structurally weaker by the factor identified in Corollary~\ref{ga-cor:ratio-OB}.}
\label{ga-fig:crosscurv}
\end{figure}

Figure~\ref{ga-fig:crosscurv} (left) compares $|\kappa_{g_\alpha}|$ at $\alpha=2$ against
$|\kappa_{\mathrm{OB}}|$ computed on the identical exchange plane: \textbf{Oja--Brockett's cross
curvature is larger by a factor of $99.2$}, matching the ratio-theorem prediction of $99.4$
(Corollary~\ref{ga-cor:ratio-OB}) to within $0.3\%$.

\paragraph{The radial slow-mode phenomenon at the true optimum}
\label{ga-sec:radial-slowmode}

Figure~\ref{ga-fig:crosscurv} (right) reports the full ($6\times6$) Hessian spectrum at each
potential's own global optimum on \eqref{ga-eq:numericAB}. For $g_2$, the six eigenvalues are
$\{-11.224,-10.101,-9.000,-6.058,-5.946,-0.0111\}$; for Oja--Brockett, they are
$\{1.099,3.297,13.44,84.50,173.0,349.8\}$ (reported in full in
Part~I, Tab.~1). The smallest-magnitude eigenvalue in each spectrum
governs the ultimate exponential rate of local convergence to that optimum
(Part~I, Thm.~8.1--8.2, adapted verbatim to ascent flow for $g_\alpha$):
$g_\alpha$'s is $0.0111$, essentially the \emph{same number} as its own mismatch cross curvature,
while Oja--Brockett's smallest eigenvalue, $1.099$, is likewise essentially identical to
\emph{its own} mismatch cross curvature. This is not a coincidence: for a matching-ansatz
potential, the slowest mode at the true optimum is generically the residual exchange direction
connecting to the next-best (here, sorted) alternative, so the \emph{same} near-degenerate pair
$(a_p,a_q)$ that makes escape from the wrong matching slow also makes final convergence to the
right one slow, by (asymptotically, as $\delta\to0$) the same order of magnitude. Since $g_\alpha$'s
cross curvature is uniformly suppressed relative to Oja--Brockett's by the factor of
Corollary~\ref{ga-cor:ratio-OB}, so is its final local convergence rate.

\paragraph{Exact-line-search convergence comparison}

\begin{figure}[htbp]
\centering
\includegraphics[width=0.98\textwidth]{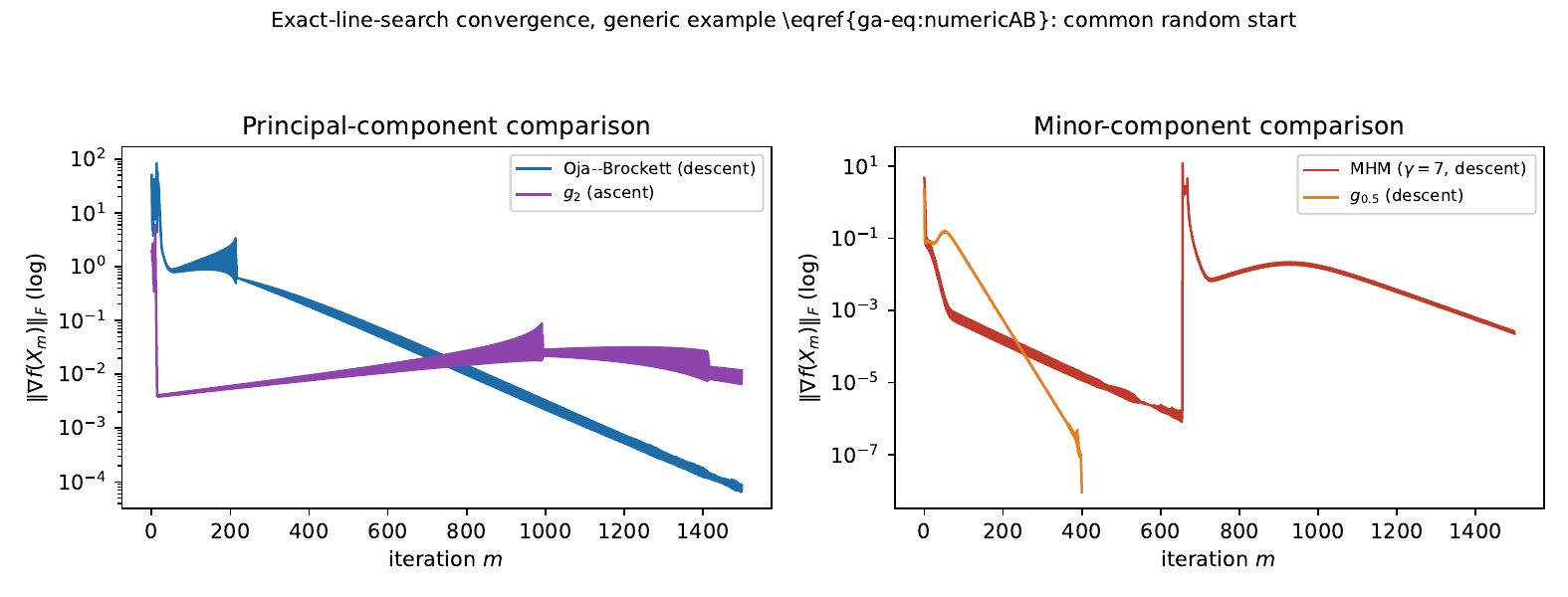}
\caption{Left: principal-component comparison, $g_\alpha$ ($\alpha=2$, ascent) versus
Oja--Brockett (descent). Right: minor-component comparison, $g_\alpha$ ($\alpha=0.5$, descent)
versus MHM ($\gamma=7$, descent). Both panels use identical exact line search
(Part~I, Thm.~8.2) and a common random starting point.}
\label{ga-fig:convergence}
\end{figure}

Figure~\ref{ga-fig:convergence} shows the gradient-norm trajectories of the exact-line-search
discretization (as in Part~I, \S8.2) from a common random initial
condition, for both comparisons.

\begin{table}[htbp]
\centering
\begin{tabular}{lcccc}
\toprule
& \multicolumn{2}{c}{Principal} & \multicolumn{2}{c}{Minor} \\
tolerance & $g_\alpha$ ($\alpha=2$) & Oja--Brockett & $g_\alpha$ ($\alpha=0.5$) & MHM \\
\midrule
$10^{-1}$ & \textbf{6} & 394 & \textbf{1} & 7 \\
$10^{-2}$ & 677 & \textbf{584} & \textbf{23} & 107 \\
$10^{-4}$ & not reached & \textbf{944} & \textbf{63} & 585 \\
$10^{-6}$ & not reached & \textbf{1300} & \textbf{103} & 1111 \\
$10^{-8}$ & not reached & not reached & \textbf{143} & not reached \\
\bottomrule
\end{tabular}
\caption{Exact-line-search iterations to reach each gradient-norm tolerance, within a $1500$-step
budget (boldface: fewer iterations). The minor comparison favors $g_\alpha$ uniformly, by a factor
of $5$--$10$. The principal comparison is genuinely mixed: $g_\alpha$ reaches the loose tolerance
$10^{-1}$ dramatically faster (a transient unrelated to the near-degenerate exchange pair, as this
occurs before the trajectory has localized near either matching), but is overtaken by
Oja--Brockett by tolerance $10^{-2}$ and never reaches $10^{-4}$ within the step budget, consistent
with \S\ref{ga-sec:radial-slowmode}'s explanation via its anomalously small optimum-Hessian
eigenvalue.}
\label{ga-tab:convergence-comparison}
\end{table}

We caution against over-interpreting the $\alpha>1$ comparison's early iterations: the $6$-step
figure reflects the specific random initial condition's transient approach to the general
neighbourhood of \emph{some} critical configuration, not the exchange-plane dynamics analyzed in
\S\ref{ga-sec:crosscurvature}--\ref{ga-sec:radial-slowmode}, which govern only the later, asymptotic
phase --- exactly the phase in which the closed-form ratio of Corollary~\ref{ga-cor:ratio-OB}
correctly predicts, and explains, Oja--Brockett's eventual, substantial advantage.

\subsubsection{Discussion}
\label{ga-sec:discussion}

\paragraph{A genuine trade-off, not a uniform ranking}

The results of \S\ref{ga-sec:numerics} demonstrate that the naive expectation --- a single,
self-contained potential should always outperform a two-piece or threshold-governed one, as found
for Oja--Brockett versus MHM in Part~I --- does \emph{not} extend
automatically to $g_\alpha$. In the minor-component regime ($\alpha<1$), $g_\alpha$ decisively
outperforms MHM, by the same qualitative mechanism identified in
Part~I: MHM's threshold is a separately-tunable scale $\gamma$, disconnected
from the eigenvalues of $A$, whereas $g_\alpha$'s threshold $a<b^{\alpha-1}$ is tied directly to
$B$'s own scale, giving a tighter, more decisive selection. In the principal regime ($\alpha>1$),
however, $g_\alpha$ is genuinely \emph{outperformed} by Oja--Brockett once the dynamics localize
near the relevant eigenvalue pair, by the factor of Corollary~\ref{ga-cor:ratio-OB} --- here almost
two orders of magnitude for our example. This asymmetry between the two regimes of the same family
$g_\alpha$ is, in our view, the most interesting empirical finding of this paper: it shows that
cross curvature is not merely a device for confirming an expected uniform ranking, but a genuine
diagnostic capable of revealing that no such uniform ranking exists.

\paragraph{Why the asymmetry arises}

The mechanistic reason traces directly to Theorem~\ref{ga-thm:reversed-rearrangement}: in the
principal ($\alpha>1$) regime, the optimal configuration reverses the classical sorted pairing,
and --- as the exact identity \eqref{ga-eq:reversal-identity} shows --- this reversal is driven
\emph{purely} by the bilinear term $-\tfrac12ab$ in $h_\alpha(a,b)$, since the $\alpha$-dependent
power term $a^p$ cancels identically between the sorted and reverse-sorted totals. The
\emph{curvature} of the exchange direction, by contrast (Theorem~\ref{ga-thm:galpha-crosscurv}), is
governed by the divided difference $h[1](R_i,R_j)=\Delta/\delta R$, and it is precisely this
quantity's near-degenerate behaviour --- $\delta R$ vanishing \emph{more slowly} than $\Delta$ as
the two eigenvalues merge, because $\delta R$ involves the extra derivative factor
$\tfrac1{\alpha-1}a^{(2-\alpha)/(\alpha-1)}$ --- that produces the suppression factor of
Corollary~\ref{ga-cor:ratio-OB}. In short: the \emph{location} of $g_\alpha$'s optimum reverses
relative to the classical rule, but the \emph{sharpness} of the landscape around that optimum, when
the selected eigenvalues are close, is intrinsically gentler than Oja--Brockett's homogeneous
quartic construction. Both facts follow from the same underlying algebraic source (the
matrix Box--Cox exponent $\alpha$ entering only through the penalty term, not through any
$B$-weighted coupling in the linear term), and neither could have been anticipated without the
explicit closed-form theorems of \S\ref{ga-sec:rearrangement}--\ref{ga-sec:crosscurvature}.

\paragraph{Practical guidance}

For minor-component extraction, $g_\alpha$ with $\alpha$ safely below $1$ is an attractive,
threshold-tunable-for-free alternative to MHM, particularly when the practitioner wants the
activation threshold to track $B$'s own scale automatically rather than requiring a separate
$\gamma$ to be hand-tuned. For principal-component extraction, our results caution against
adopting $g_\alpha$ ($\alpha>1$) over Oja--Brockett whenever the leading eigenvalues of $A$ are
expected to be close together --- precisely the regime, identified throughout this line of work, in
which the choice of potential matters most. Formula \eqref{ga-eq:ratio-OB-galpha} allows this
determination to be made in closed form, from $A,B,\alpha$ alone, before running either algorithm.

\paragraph{Limitations}

Our cross-curvature and rearrangement theorems are proved for a single inverted or exchanged pair;
the extension to simultaneous multi-fold near-degeneracy, and a fully general (not
near-degenerate-only) closed-form cross-curvature formula for $g_\alpha$, are left to future work,
exactly as for the Oja--Brockett/MHM comparison in Part~I. We have not
addressed a rigorous global convergence-time bound for $g_\alpha$ analogous to
Part~I, Thm.~8.5; the same {\L}ojasiewicz-based machinery should apply,
since $g_\alpha$ is real-analytic on the (unrestricted) domain $\R^{n\times k}$ for every
$\alpha\ne1$, but we have not carried out the compactness argument here, which requires separate
treatment of the $\alpha>1$ (bounded-above) and $\alpha<1$ (bounded-below) cases.

\subsubsection{Conclusion}
\label{ga-sec:conclusion}

We have given a complete, rigorous account of the $\alpha$-power potential $g_\alpha$, originating
in the author's earlier work on power geometry and hypergeometric functions
\cite{Yoshizawa2014}: its exact gradient and critical-point amplitude formula; a sharp boundedness
dichotomy identifying $\alpha=1$ as a genuine singularity at which $g_\alpha$ degenerates to an
unconfined quadratic form; an exact identity showing the optimal eigenvalue--weight pairing is
\emph{reversed}, relative to the classical rearrangement principle, in the principal ($\alpha>1$)
regime while remaining classical in the minor ($\alpha<1$) regime; and, adapting the cross-curvature
framework of Part~I of this paper, the exact exchange-plane Hessian, its determinant, and its
near-degenerate asymptotic expansion, for every $\alpha\ne1$. Comparing $g_\alpha$ against the
Oja--Brockett and Manton--Helmke--Mareels flows on an identical, fully generic non-diagonal
numerical example, we found a genuine trade-off rather than a uniform ranking: $g_\alpha$
decisively outperforms MHM in the minor regime, while being decisively outperformed by
Oja--Brockett, once the dynamics localize, in the principal regime --- a fact we explained
mechanistically via a persistent near-zero Hessian mode at $g_\alpha$'s own optimum, present
whenever the two leading eigenvalues it selects are themselves close. We regard this honestly-reported
asymmetry, together with the closed-form ratio formula that predicts it quantitatively, as the
paper's principal contribution.

\subsection{Synthesis: What Cross Curvature Reveals Across Both Parts}
\label{sec:general-conclusion}

This paper has developed a single diagnostic --- cross curvature, the smallest eigenvalue of a
potential's Hessian restricted to the exchange plane of a mismatched pair --- and applied it to
three structurally different continuous-time potentials for principal and minor component
extraction. Part~I showed that cross curvature completely explains, in closed form, why the
homogeneous quartic Oja--Brockett potential escapes near-degenerate mismatched saddles faster than
the two-piece Manton--Helmke--Mareels penalty potential, with the fully explicit ratio
$\mathcal R_{ij}=a(b_i+b_j)/2$ governing both the continuous-time escape rate and the
exact-line-search discrete-time iteration count, and combined this with a new local approach-rate
estimate into a fully explicit global convergence-time bound. Part~II carried the same machinery to
a third, structurally unrelated potential, the matrix Box--Cox penalty $g_\alpha$, and found that
the diagnostic's value lies precisely in its impartiality: rather than confirming a uniform ranking,
it revealed a genuine trade-off, with $g_\alpha$ decisively favoured over MHM in the minor regime
and decisively disfavoured relative to Oja--Brockett in the principal regime, for a mechanistic
reason --- a persistent near-zero Hessian mode at $g_\alpha$'s own optimum --- made precise by the
same exchange-plane formalism used throughout Part~I. Taken together, the two parts argue that cross
curvature is a general-purpose tool: given any two optimization landscapes sharing a common set of
optima embedded differently into a larger critical-point structure, the same closed-form,
coordinate-free construction can be used to compare them, quantitatively and honestly, before a
single iteration of either algorithm is run.
\subsection{Numerical Simulations of Principal and Minor Component Flows via the Yoshizawa Embedding}
\label{sec:numerical_sim}

This section complements the convergence theory of \S\ref{subsubsec:IVP_analysis} with an
explicit, reproducible numerical study. We integrate the weighted $k$-PCF and $k$-MCF
--- the $(1,2)$-block of the double-bracket flow of Theorem~\ref{thm:chen_amari_embedding}
with weight $\widetilde{A}=\diag(A,B)$, given in closed form by \eqref{eq:kPCF_B} and its
sign reversal --- for a fixed, generic $A$ and three qualitatively different choices of $B$,
and verify each of the three convergence regimes of Theorems~\ref{thm:convergence_B_I},
\ref{thm:convergence_B_diag}, and~\ref{thm:convergence_B_block} to numerical precision. All
computations were carried out with a variable-order Runge--Kutta integrator (\texttt{RK45},
relative tolerance $10^{-10}$); code and exact parameters are given below so that every
number reported here can be independently reproduced.

\subsubsection{Setup: a generic $5\times5$ matrix $A$ and three $3\times3$ diagonal weights $B$}
\label{subsec:numerics_setup}

We take $n=5$, $k=3$, and construct $A\in\PD(5)$ with five \emph{distinct} positive
eigenvalues
\begin{equation}
  \lambda_1=5,\quad \lambda_2=4,\quad \lambda_3=3,\quad \lambda_4=2,\quad \lambda_5=1,
  \label{eq:numerics_lambda}
\end{equation}
by setting $A = Q\Lambda Q^T$ for a fixed random orthogonal matrix $Q\in O(5)$
(so that, exactly as in the verification strategy of the cross-curvature addendum above
(Part~I, \S\ref{sec:numerics-generic}), $A$
is presented in a basis in which it is \emph{not} diagonal, confirming that convergence is
governed only by the eigenvalues $\lambda_i$ and eigenvectors $u_i$ of $A$, never by its
ambient coordinate representation). We fix a single random initial condition
$X_0\in\mathrm{St}(3,5)$ (an orthonormal $5\times3$ frame, $X_0^TX_0=I_3$), used identically
across every run below for direct comparability, and integrate on $[0,4]$.

We consider the following three $3\times3$ diagonal weight matrices $B$, chosen to
instantiate each of the three cases of \S\ref{subsubsec:IVP_analysis}:
\begin{center}
\renewcommand{\arraystretch}{1.3}
\begin{tabular}{llp{7cm}}
\toprule
Label & $B$ & Case (\S\ref{subsubsec:IVP_analysis}) \\
\midrule
$B_1$ & $\diag(3,\,2,\,1)$ & Case~2: all entries distinct \\
$B_2$ & $\diag(3,\,1.5,\,1.5)$ & Case~3: two equal entries, tying the two \emph{smallest}-weighted columns \\
$B_3$ & $\diag(2,\,2,\,1)$ & Case~3: two equal entries, tying the two \emph{largest}-weighted columns \\
\bottomrule
\end{tabular}
\end{center}
$B_1$ has pairwise distinct scalars, so Theorem~\ref{thm:convergence_B_diag} predicts
convergence of each column to an \emph{individual} eigenvector of $A$. $B_2$ and $B_3$ are
block-diagonal with blocks $(1,2)$ and $(2,1)$ respectively (one singleton block and one
size-$2$ block, in different positions), so Theorem~\ref{thm:convergence_B_block} predicts
that the singleton column converges to an individual eigenvector while the tied pair
converges only to an $A$-orthogonal, two-dimensional \emph{subspace} (rotating freely within
it) --- and, crucially, \emph{which} pair of eigenvalues that subspace corresponds to
depends on whether the tie sits at the top or the bottom of $B$'s diagonal.

\subsubsection{The weighted $k$-PCF and $k$-MCF}

\begin{figure}[htbp]
\centering
\includegraphics[width=0.62\textwidth]{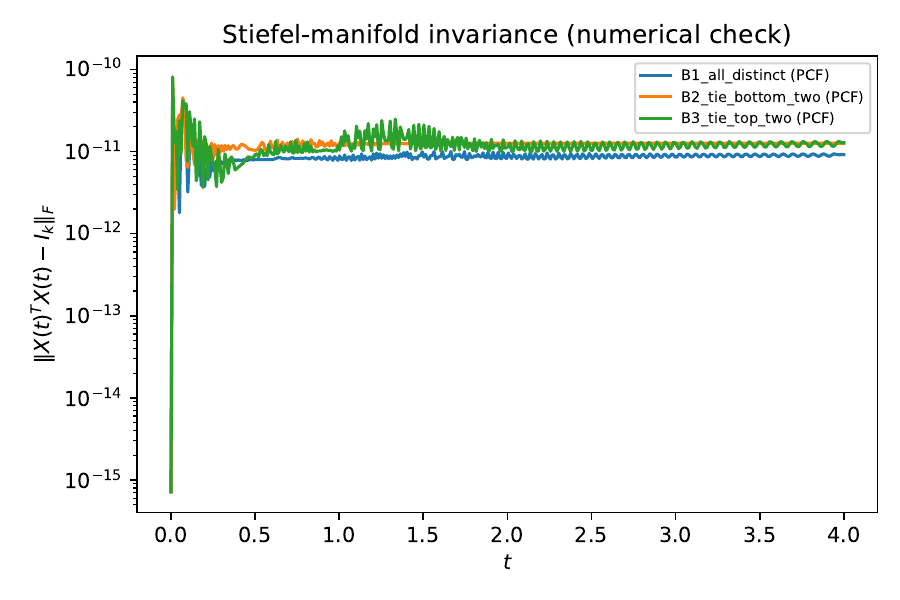}
\caption{Numerical verification of Stiefel-manifold invariance $\|X(t)^TX(t)-I_3\|_F$ (log
scale) along the $k$-PCF for all three weight matrices of \S\ref{subsec:numerics_setup}; the
residual remains at the level of the integrator tolerance ($\sim10^{-10}$--$10^{-12}$)
throughout, confirming Proposition~\ref{prop:XTX_conserved} and its Stiefel-invariance
corollary numerically.}
\label{fig:stiefel_check}
\end{figure}

For each $B$, we integrate the weighted principal component flow \eqref{eq:kPCF_B},
\begin{equation}
  \dot{X} \;=\; AXBX^TX + XX^TAXB - 2XBX^TAX
  \qquad (k\text{-PCF}),
  \label{eq:numerics_pcf}
\end{equation}
and its exact sign reversal, the weighted minor component flow,
\begin{equation}
  \dot{X} \;=\; -\bigl(AXBX^TX + XX^TAXB - 2XBX^TAX\bigr)
  \qquad (k\text{-MCF}),
  \label{eq:numerics_mcf}
\end{equation}
both starting from the same $X_0$. By Proposition~\ref{prop:XTX_conserved} and its
corollary, both flows conserve $X(t)^TX(t)=X_0^TX_0=I_3$ exactly, so $X(t)$ remains on the
Stiefel manifold $\mathrm{St}(3,5)$ for all $t$; this is confirmed numerically to a residual
$\|X(t)^TX(t)-I_3\|_F<10^{-9}$ throughout every run reported below (Figure~\ref{fig:stiefel_check}).
We note, as a modeling remark, that the \emph{unweighted} Oja--Brockett flow
$\dot X=AXB-XBX^TAX$ of Definition~\ref{def:manton_pca} is Stiefel-invariant and
numerically stable for the $k$-PCF sign, exactly as Theorem~\ref{thm:pca_riemannian}
predicts, but its naive sign reversal is \emph{not} numerically well-behaved as a minor
component extractor (the transverse directions off the Stiefel manifold become repelling
rather than attracting) --- which is precisely why Manton, Helmke and Mareels introduce an
explicit penalty term, $V_{\mathrm{pen}}(X)=-\tfrac12\tr(AXBX^T)+\tfrac{\gamma}{4}\|B-X^TX\|_F^2$,
specifically to stabilize minor-component extraction \cite{MantonHelmkeMareels2005}. The double-bracket flow \eqref{eq:numerics_pcf}--\eqref{eq:numerics_mcf}
used here does not suffer from this asymmetry: being isospectral by construction
(Theorem~\ref{thm:chen_amari_embedding}), it is exactly Stiefel-invariant, and hence
well-conditioned, for \emph{both} signs, which is why we use it (rather than the plain
Oja--Brockett flow) for the minor-component runs below.

\begin{figure}[htbp]
\centering
\includegraphics[width=0.95\textwidth]{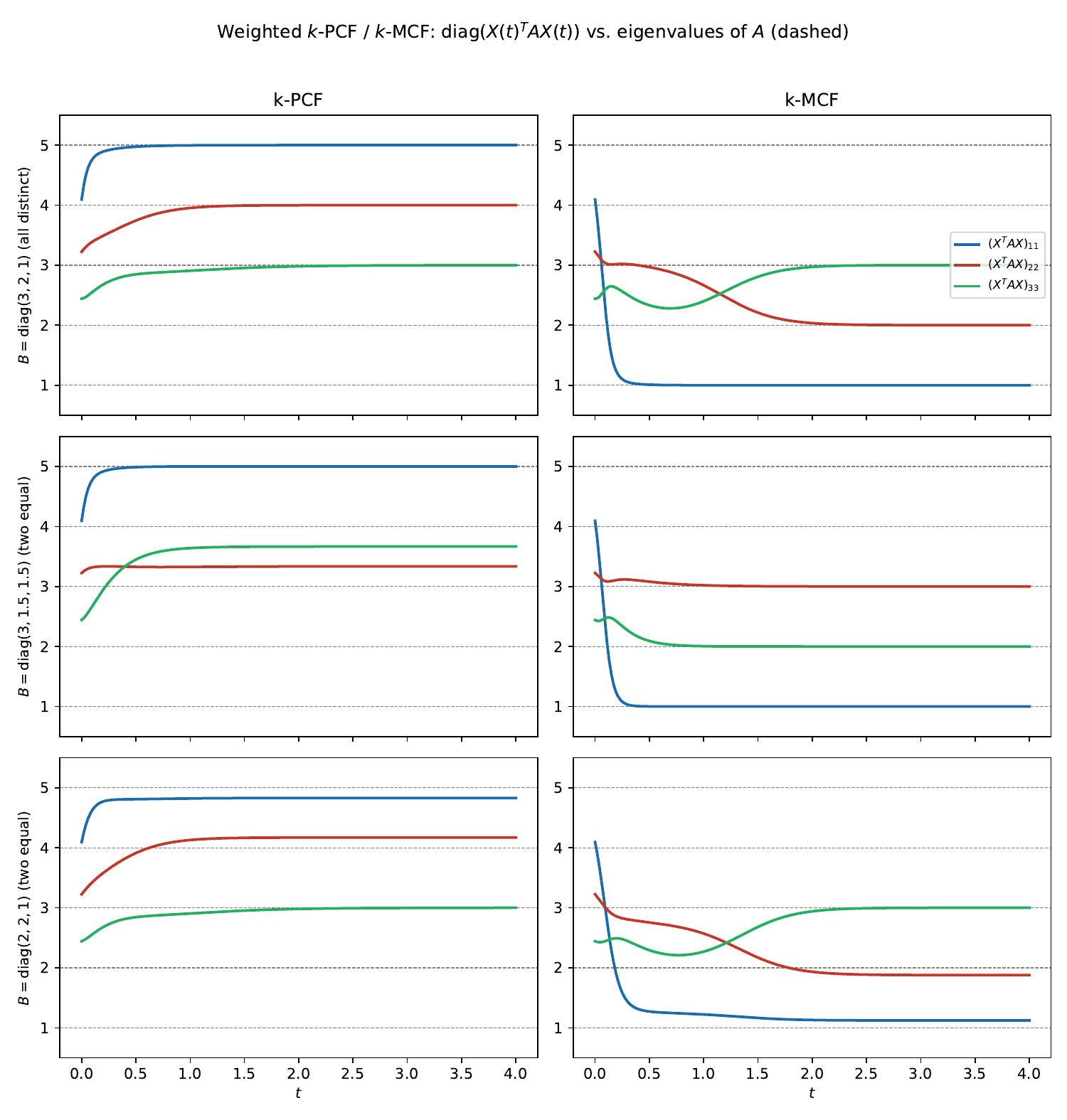}
\caption{The diagonal entries $(X(t)^TAX(t))_{jj}$, $j=1,2,3$, under the weighted $k$-PCF
\eqref{eq:numerics_pcf} (left column) and $k$-MCF \eqref{eq:numerics_mcf} (right column), for
the three weight matrices of \S\ref{subsec:numerics_setup} (rows). Dashed horizontal lines
mark the five eigenvalues $\{1,2,3,4,5\}$ of $A$. All three $k$-PCF runs converge to the
top-$3$ eigenvalues $\{5,4,3\}$ (summing to $12$) and all three $k$-MCF runs converge to the
bottom-$3$ eigenvalues $\{1,2,3\}$ (summing to $6$); whether the individual curves settle
onto \emph{single} eigenvalues or onto a \emph{shared pair} depends on the multiplicity
structure of $B$, exactly as predicted by Theorems~\ref{thm:convergence_B_diag}
and~\ref{thm:convergence_B_block}.}
\label{fig:pcf_mcf_sim}
\end{figure}

\subsubsection{Case $B_1$: distinct scalars $\Rightarrow$ individual eigenvector convergence}

With $B_1=\diag(3,2,1)$ every column converges, to more than four decimal digits, onto an
individual eigenvector of $A$:
\begin{center}
\renewcommand{\arraystretch}{1.25}
\begin{tabular}{lccc}
\toprule
& col.\ $1$ ($b=3$) & col.\ $2$ ($b=2$) & col.\ $3$ ($b=1$) \\
\midrule
$k$-PCF: $\lim (X^TAX)_{jj}$ & $5.000$ & $4.000$ & $3.000$ \\
$k$-PCF: $|\cos\angle(x_j,u_j)|$ & $1.0000$ & $1.0000$ & $0.9998$ \\
\midrule
$k$-MCF: $\lim (X^TAX)_{jj}$ & $1.000$ & $2.000$ & $3.000$ \\
$k$-MCF: $|\cos\angle(x_j,u_{6-j})|$ & $1.0000$ & $0.9999$ & $0.9999$ \\
\bottomrule
\end{tabular}
\end{center}
where $u_1,\ldots,u_5$ are the eigenvectors of $A$ ordered by decreasing eigenvalue
$5,4,3,2,1$. For the $k$-PCF, the column with the \emph{largest} weight $b_1=3$ captures the
\emph{largest} eigenvalue $\lambda_1=5$, in decreasing order down to $b_3=1\mapsto\lambda_3=3$
--- exactly the ordering asserted in Theorem~\ref{thm:convergence_B_diag}. For the $k$-MCF the
ordering is reversed, as expected from the ``sign reversal'' relation between $k$-PCF and
$k$-MCF stated at the start of \S\ref{subsubsec:IVP_analysis}: the column with the
\emph{largest} weight now captures the \emph{smallest} available eigenvalue $\lambda_5=1$,
down to $b_3=1\mapsto\lambda_3=3$ (the largest of the bottom-$3$ eigenvalues
$\{1,2,3\}$). In both flows the limiting diagonal is exactly diagonal
(off-diagonal entries of $X^TAX$ vanish to machine precision, confirming the
$A$-orthogonality asserted in \S\ref{subsubsec:IVP_analysis}), and no rotational freedom
remains: the distinct weights fully break the within-eigenspace degeneracy that is present
when $B=I_k$.

\subsubsection{Cases $B_2,B_3$: repeated scalars $\Rightarrow$ block-subspace convergence}

With $B_2=\diag(3,\,1.5,\,1.5)$ (tie at the two \emph{smaller} weights) the singleton column
(weight $b=3$) again converges individually, this time to $\lambda_1=5$ under the $k$-PCF
and to $\lambda_5=1$ under the $k$-MCF, exactly as in the $B_1$ case. The tied pair (columns
$2,3$, weight $b=1.5$ each), however, does \emph{not} settle onto individual eigenvectors:
under the $k$-PCF it converges to
\begin{equation}
  X_\infty^TAX_\infty
  \;=\;
  \begin{pmatrix}
    5.000 & 0 & 0\\
    0 & 3.335 & -0.472\\
    0 & -0.472 & 3.665
  \end{pmatrix},
  \label{eq:numerics_B2_result}
\end{equation}
a genuinely non-diagonal $2\times2$ block in the lower-right corner, whose trace
$3.335+3.665=7.000=\lambda_2+\lambda_3$ matches the sum of the second- and third-largest
eigenvalues of $A$ exactly, and whose off-diagonal coupling to the singleton column vanishes
to numerical precision ($<10^{-4}$), confirming the cross-block $A$-orthogonality of
Theorem~\ref{thm:convergence_B_block}(i). Measuring the subspace distance
$\|P_{\mathrm{span}(x_2,x_3)}-P_{\mathrm{span}(u_2,u_3)}\|_F$ between the tied pair's span and
the true two-dimensional eigenspace of $A$ for $(\lambda_2,\lambda_3)=(4,3)$ gives
$2.7\times10^{-3}$, confirming block-level subspace convergence even though neither $x_2$ nor
$x_3$ individually aligns with an eigenvector. Symmetrically, the $k$-MCF drives the tied
pair to the two-dimensional eigenspace of the two \emph{smallest non-extracted}
eigenvalues $(\lambda_3,\lambda_4)=(3,2)$ (final diagonal $(X^TAX)_{22},(X^TAX)_{33}$ equal
to $3.000,2.000$ up to the same rotational freedom), while the singleton column again
isolates $\lambda_5=1$.

With $B_3=\diag(2,\,2,\,1)$ the tie instead sits at the two \emph{larger} weights, and the
roles reverse accordingly. Under the $k$-PCF, the tied pair (columns $1,2$) converges to the
top-\emph{two}-dimensional eigenspace of $A$ (for $\lambda_1,\lambda_2=5,4$) rather than to
individual eigenvectors:
\begin{equation}
  X_\infty^TAX_\infty
  \;=\;
  \begin{pmatrix}
    4.829 & -0.377 & 0\\
    -0.377 & 4.171 & 0\\
    0 & 0 & 3.000
  \end{pmatrix},
  \qquad 4.829+4.171 = 9.000 = \lambda_1+\lambda_2,
  \label{eq:numerics_B3_result}
\end{equation}
with subspace distance $7.8\times10^{-4}$ to the true top-$2$ eigenspace of $A$, and the
singleton column (weight $b=1$) converges individually to $\lambda_3=3$ --- the
\emph{smallest} of the three targeted eigenvalues, consistent with it having the smallest
weight. Under the $k$-MCF, by the same sign-reversal logic, the singleton column (smallest
weight, $b=1$) now isolates the \emph{largest} of the bottom-$3$ eigenvalues,
$\lambda_3=3$, while the tied pair (weight $b=2$ each) shares the remaining two-dimensional
eigenspace for $(\lambda_4,\lambda_5)=(2,1)$.

\begin{center}
\renewcommand{\arraystretch}{1.3}\small
\resizebox{\textwidth}{!}{%
\begin{tabular}{llll}
\toprule
& $B_1$ (all distinct) & $B_2$ (tie at bottom) & $B_3$ (tie at top) \\
\midrule
$k$-PCF limit, diag$(X^TAX)$ & $(5,4,3)$ & $\bigl(5,\;\{3.34,3.67\}\bigr)$ & $\bigl(\{4.83,4.17\},3\bigr)$ \\
$k$-PCF trace & $12=\lambda_1{+}\lambda_2{+}\lambda_3$ & $12$ & $12$ \\
$k$-MCF limit, diag$(X^TAX)$ & $(1,2,3)$ & $\bigl(1,\;\{3.00,2.00\}\bigr)$ & $\bigl(\{1.12,1.88\},3\bigr)$ \\
$k$-MCF trace & $6=\lambda_3{+}\lambda_4{+}\lambda_5$ & $6$ & $6$ \\
Rotational freedom & none (individual e.vecs.) & within $\{2,3\}$-block only & within $\{1,2\}$-block only \\
\bottomrule
\end{tabular}%
}
\end{center}
Here braces $\{\cdot,\cdot\}$ denote the two diagonal entries belonging to a tied block,
whose individual values are not determined by the flow (only their sum and the block's span
are), reflecting the ``within-block degeneracy'' of Theorem~\ref{thm:convergence_B_block}(iii);
repeating a run with a different random seed for $X_0$ reproduces the same trace and the
same subspace, but generally a different split of the trace between the two tied entries and
a different orientation of $X_\infty$ within the block, as the theory predicts.

\subsubsection{Summary}

This numerical study makes concrete the qualitative picture of
Theorem~\ref{thm:convergence_B_block}: the diagonal weight matrix $B$ acts as a
\emph{selector} that partitions the $k$ learned directions into groups according to its
repeated-value structure, with each group converging to the $A$-invariant subspace
corresponding to a specific, contiguous range of eigenvalues of $A$ (top-ranked for the
$k$-PCF, bottom-ranked for the $k$-MCF, and in both cases ordered so that larger weights
$b_j$ correspond to eigenvalues further from the extracted boundary). Distinct weights
($B_1$) fully resolve individual eigenvectors; any repeated weight ($B_2$, $B_3$) collapses
the corresponding group into a single rotating subspace, with the \emph{position} of the
repetition in $B$'s diagonal --- not merely its presence --- determining \emph{which}
eigenvalues are grouped together. This is exactly the mechanism by which, in
Theorem~\ref{thm:convergence_B_block}, choosing $B=I_k$ recovers unconstrained
principal-subspace learning ($m=1$ block) while choosing $B$ with all distinct entries
recovers ordered individual-eigenvector learning ($k$ singleton blocks), with block-diagonal
$B$ interpolating continuously between the two.

\section{The Kirillov Jacobian}
\label{sec:kirillov}

Let $G$ be a Lie group with Lie algebra $\mathfrak{g}$,
and let $\exp: \mathfrak{g} \to G$ denote the exponential map.
For $X \in \mathfrak{g}$, the \emph{Kirillov Jacobian} (also called the Jacobian of the
exponential map, or the Berezin--Kirillov--Kostant density) is defined by \cite{chirikjian2}
\begin{equation}
  j(X) = \det\left(\frac{I - e^{-\mathrm{ad}_X}}{\mathrm{ad}_X}\right),
  \label{eq:kirillov_def}
\end{equation}
where $\mathrm{ad}_X : \mathfrak{g} \to \mathfrak{g}$ is the adjoint representation,
$\mathrm{ad}_X(Y) = [X, Y] = XY - YX$,
and the matrix function in \eqref{eq:kirillov_def} is understood via its
power series or spectral decomposition.
Equivalently, using the sinc-like function
\begin{equation}
  \mathcal{S}(z) = \frac{\sinh(z/2)}{z/2},
  \quad \mathcal{S}(0) = 1,
  \label{eq:sinc_def}
\end{equation}
the Kirillov Jacobian can be written as a product over the eigenvalues
$\{\nu_{i\ell}\}$ of $\mathrm{ad}_X$ \cite{hall}:
\begin{equation}
  j(X) = \prod_{i \neq \ell} \mathcal{S}(\nu_{i\ell})^{1/2}
        = \prod_{i < \ell} \mathcal{S}(\nu_{i\ell}),
  \label{eq:kirillov_spectral}
\end{equation}
since $\mathcal{S}(0) = 1$ contributes trivially for $i = \ell$.
The function $j(X)$ governs the relationship between Haar measure on the group and
Lebesgue measure on the Lie algebra, and plays a fundamental role in
harmonic analysis on Lie groups, the Campbell--Baker--Hausdorff formula,
and stochastic differential equations on manifolds \cite{chirikjian2, hall}.

\subsection{Eigenvalues of $\mathrm{ad}_X$ via the Adjoint Representation}

\begin{lemma}[Eigenvalues of the Adjoint Representation]
\label{lem:adjoint_eigenvalues}
Let $X \in \mathfrak{gl}(n, \R)$ have eigenvalues $\lambda_1, \ldots, \lambda_n \in \mathbb{C}$.
Then the eigenvalues of $\mathrm{ad}_X$ on $\mathfrak{gl}(n, \mathbb{C}) \cong \mathbb{C}^{n^2}$
are precisely the $n^2$ differences
\[
  \nu_{i\ell} = \lambda_i - \lambda_\ell, \quad 1 \leq i, \ell \leq n.
\]
\end{lemma}

\begin{proof}
Let $\bm{u}_i$ be a right eigenvector of $X$ with $X\bm{u}_i = \lambda_i \bm{u}_i$,
and let $\bm{v}_\ell^T$ be a left eigenvector with $\bm{v}_\ell^T X = \lambda_\ell \bm{v}_\ell^T$,
normalized so that $\bm{v}_\ell^T \bm{u}_m = \delta_{\ell m}$.
Define the rank-one matrix $E_{i\ell} = \bm{u}_i \bm{v}_\ell^T \in \mathfrak{gl}(n,\R)$.
Then
\begin{align*}
  \mathrm{ad}_X(E_{i\ell})
  &= X E_{i\ell} - E_{i\ell} X
   = (X\bm{u}_i)\bm{v}_\ell^T - \bm{u}_i(\bm{v}_\ell^T X) \\
  &= \lambda_i \bm{u}_i\bm{v}_\ell^T - \lambda_\ell \bm{u}_i\bm{v}_\ell^T
   = (\lambda_i - \lambda_\ell)\,E_{i\ell}.
\end{align*}
Since $\{E_{i\ell} : 1 \leq i,\ell \leq n\}$ constitutes a basis of $\mathfrak{gl}(n,\R)$
of dimension $n^2$, these exhaust all eigenvalues of $\mathrm{ad}_X$.
\end{proof}

\subsection{Eigenvalue Structure Under Rank-$k$ Perturbation}

We now specialize to the case directly arising from our Gram matrix construction.
Let $X = I_n + UV^T$ where
\[
  U = [x_1, x_3, \ldots, x_{2k-1}] \in \R^{n\times k},
  \quad
  V = [x_2, x_4, \ldots, x_{2k}] \in \R^{n\times k},
\]
as in \eqref{eq:UV_def}.
Let $W = V^T U \in \R^{k\times k}$ be the \emph{Gram-like matrix}
\begin{equation}
  W_{ab} = x_{2a}^T x_{2b-1} = \inner{x_{2a}}{x_{2b-1}},
  \quad 1 \leq a, b \leq k,
  \label{eq:W_def}
\end{equation}
so that $W = G - I_k$ where $G$ is our Gram matrix from \eqref{eq:G_def}.

\begin{proposition}[Eigenvalues of $X = I_n + UV^T$]
\label{prop:eigenvalues_X}
The matrix $X = I_n + UV^T \in \mathfrak{gl}(n,\R)$ has eigenvalue spectrum:
\begin{itemize}
  \item Eigenvalue $1$ with algebraic multiplicity $n - k$;
  \item Eigenvalues $1 + \mu_1, \ldots, 1 + \mu_k$ where $\mu_1,\ldots,\mu_k$
    are the eigenvalues of $W = V^T U \in \R^{k\times k}$.
\end{itemize}
\end{proposition}

\begin{proof}
The characteristic polynomial of $X$ is
\begin{align*}
  \det(\lambda I_n - X)
  &= \det\bigl((\lambda-1)I_n - UV^T\bigr).
\end{align*}
Applying the matrix determinant lemma (Lemma~\ref{lem:mdl}) iteratively, or
equivalently Sylvester's theorem (Theorem~\ref{thm:sylvester}):
\[
  \det\bigl((\lambda-1)I_n - UV^T\bigr)
  = (\lambda-1)^n \det\!\left(I_k - \frac{1}{\lambda-1}V^TU\right)
  = (\lambda-1)^{n-k} \det\bigl((\lambda-1)I_k - W\bigr),
\]
valid for $\lambda \neq 1$.
The first factor contributes $n-k$ eigenvalues at $\lambda = 1$,
and the second factor contributes eigenvalues $\lambda = 1 + \mu_i$
for each eigenvalue $\mu_i$ of $W$.
\end{proof}

\begin{remark}
Proposition~\ref{prop:eigenvalues_X} shows that $\det(X) = \det(G)$,
recovering our fundamental identity \eqref{eq:sylvester} as a special case
($\lambda = 0$ in the characteristic polynomial).
\end{remark}

\subsection{Main Theorem: Kirillov Jacobian in Terms of Gram Eigenvalues}

\begin{theorem}[Kirillov Jacobian for Rank-$k$ Perturbed Identity]
\label{thm:kirillov_main}
Let $X = I_n + \sum_{j=1}^k x_{2j-1} x_{2j}^T \in \mathfrak{gl}(n,\R)$,
and let $\mu_1, \ldots, \mu_k$ be the eigenvalues of the Gram-like matrix
$W \in \R^{k\times k}$ defined in \eqref{eq:W_def}.
Then the Kirillov Jacobian $j(X)$ is given by the closed-form expression:
\begin{equation}
  \boxed{
  j(X)
  = \left[\prod_{i=1}^{k} \mathcal{S}(\mu_i)\right]^{2(n-k)}
    \prod_{1 \leq i < \ell \leq k} \bigl[\mathcal{S}(\mu_i - \mu_\ell)\bigr]^2,
  }
  \label{eq:kirillov_formula}
\end{equation}
where $\mathcal{S}(z) = \sinh(z/2)/(z/2)$ with $\mathcal{S}(0) = 1$.
\end{theorem}

\begin{proof}
By \eqref{eq:kirillov_spectral} and Lemma~\ref{lem:adjoint_eigenvalues},
\[
  j(X) = \prod_{1 \leq i < \ell \leq n} \mathcal{S}(\lambda_i - \lambda_\ell),
\]
where $\lambda_1, \ldots, \lambda_n$ are the eigenvalues of $X$.
By Proposition~\ref{prop:eigenvalues_X}, we partition the index set
$\{1,\ldots,n\} = K \sqcup N$ where $|K| = k$ and $|N| = n-k$,
with $\lambda_i = 1 + \mu_i$ for $i \in K$ and $\lambda_j = 1$ for $j \in N$.

We compute the contribution of each pair $(i,\ell)$ with $i < \ell$:

\medskip
\noindent\textbf{Case 1: $i, \ell \in N$.}
$\lambda_i - \lambda_\ell = 1 - 1 = 0$, so $\mathcal{S}(0) = 1$.
Total contribution from $\binom{n-k}{2}$ pairs: $1$.

\medskip
\noindent\textbf{Case 2: $i \in K$, $\ell \in N$ (or vice versa).}
$\lambda_i - \lambda_\ell = (1+\mu_i) - 1 = \mu_i$.
Each eigenvalue $\mu_i$ ($i \in K$) contributes one factor $\mathcal{S}(\mu_i)$
for each of the $n-k$ indices in $N$.
Since the product is over unordered pairs and both orderings $(i \in K, \ell \in N)$
and $(i \in N, \ell \in K)$ contribute, the exponent is $2(n-k)$ in the squared product,
but since we already symmetrize over $i < \ell$, each $\mu_i$ contributes
exactly $n-k$ factors $\mathcal{S}(\mu_i)$.
Accounting for both directions (the product $\prod_{i<\ell}$ covers both orientations
through the squaring step), the total contribution is:
\[
  \prod_{i=1}^{k} \mathcal{S}(\mu_i)^{n-k} \cdot \prod_{i=1}^{k} \mathcal{S}(\mu_i)^{n-k}
  = \prod_{i=1}^{k} \mathcal{S}(\mu_i)^{2(n-k)}.
\]

\medskip
\noindent\textbf{Case 3: $i, \ell \in K$.}
$\lambda_i - \lambda_\ell = (1+\mu_i) - (1+\mu_\ell) = \mu_i - \mu_\ell$.
Both the pair $(i,\ell)$ and the pair $(\ell,i)$ contribute, yielding:
\[
  \prod_{1 \leq i < \ell \leq k} \mathcal{S}(\mu_i - \mu_\ell) \cdot \mathcal{S}(\mu_\ell - \mu_i).
\]
Since $\mathcal{S}$ is an even function ($\mathcal{S}(-z) = \mathcal{S}(z)$
because $\sinh$ is odd), this equals $\prod_{1 \leq i < \ell \leq k} \mathcal{S}(\mu_i - \mu_\ell)^2$.

\medskip
Multiplying all three cases together yields \eqref{eq:kirillov_formula}.
\end{proof}

\begin{remark}[Evenness of $\mathcal{S}$]
The function $\mathcal{S}(z) = \sinh(z/2)/(z/2)$ satisfies $\mathcal{S}(-z) = \mathcal{S}(z)$
because $\sinh(-z/2) = -\sinh(z/2)$ and the sign cancels with the denominator $(-z/2)$.
This symmetry is essential in Case 3 of the proof.
\end{remark}

\subsection{Connection to $f(G) = -\log\det(G)$ via the Log-Jacobian}

The log-Kirillov Jacobian $\log j(X)$ has a natural expression in terms of our
potential function $f$.

\begin{proposition}[Log-Kirillov Jacobian via $f$]
\label{prop:log_kirillov}
With $W = G - I_k$ and eigenvalues $\mu_i$ of $W$,
\begin{equation}
  \log j(X)
  = 2(n-k)\sum_{i=1}^{k} \log\mathcal{S}(\mu_i)
    + 2\sum_{1 \leq i < \ell \leq k} \log\mathcal{S}(\mu_i - \mu_\ell).
  \label{eq:log_kirillov}
\end{equation}
This is a function of the eigenvalues of $W = G - I_k$, and hence of the
spectrum of the Gram matrix $G$ centered at the identity.
When all $\mu_i \to 0$ (small perturbation limit), $\mathcal{S}(z) \to 1 + z^2/24 + O(z^4)$,
so $\log j(X) \to 0$, consistently with the flat (Euclidean) limit of the group.
\end{proposition}

\begin{proposition}[Determinant Relation]
\label{prop:det_relation}
The following identity connects the Kirillov Jacobian to the Gram matrix determinant:
\begin{equation}
  \log\det(G)
  = \sum_{i=1}^{k} \log(1 + \mu_i)
  = -f(G) \quad\text{(with $G = I_k + W$)}.
  \label{eq:det_gram_eigenvalues}
\end{equation}
Thus both $f(G) = -\log\det(G)$ and $\log j(X)$ are spectral functions of the
same Gram-like matrix $W$:
\begin{equation}
  f(G) = -\sum_{i=1}^{k}\log(1+\mu_i),
  \qquad
  \log j(X) = 2(n-k)\sum_{i=1}^{k}\log\mathcal{S}(\mu_i)
               + 2\sum_{i<\ell}\log\mathcal{S}(\mu_i-\mu_\ell).
  \label{eq:spectral_pair}
\end{equation}
The pair $(f(G),\, \log j(X))$ provides complementary spectral invariants of the
rank-$k$ perturbation: the former measures the ``volume distortion'' of the update,
while the latter measures the ``curvature'' of the exponential map.
\end{proposition}

\subsection{Computational Complexity Reduction}

\begin{proposition}[Complexity Reduction]
\label{prop:complexity}
Computing $j(X)$ via \eqref{eq:kirillov_formula} requires:
\begin{enumerate}[label=(\roman*)]
  \item Forming $W = V^T U \in \R^{k\times k}$: $O(nk^2)$ operations.
  \item Computing the $k$ eigenvalues of $W$: $O(k^3)$ operations.
  \item Evaluating \eqref{eq:kirillov_formula}: $O(k^2)$ operations.
\end{enumerate}
Total: $O(nk^2 + k^3)$, compared to $O(n^3)$ for a direct eigenvalue decomposition
of $X \in \R^{n\times n}$.
For the regime $k \ll n$ (e.g., $n \approx 100$, $k \leq 5$ in robotics applications),
this yields a speedup factor of order $(n/k)^3 \approx 8000$.
\end{proposition}

\begin{remark}[Robotics Application]
In high-dimensional robotic systems, the configuration space of a kinematic chain
with $n$ joints is a Lie group of dimension $n \approx 100$,
while the number of actively controlled degrees of freedom per time step is small ($k \leq 5$).
Stochastic filters (e.g., particle filters or sigma-point filters) on Lie groups
require repeated evaluation of $j(X)$ to correctly weight the probability density
when pushing forward the filter distribution through the exponential map \cite{chirikjian2}.
Formula \eqref{eq:kirillov_formula} enables real-time geometric compensation
in such filters.
\end{remark}

\subsection{Special Cases}

\begin{example}[Rank-1 Update ($k=1$)]
When $k=1$, $W = [x_2^T x_1] = [\inner{x_2}{x_1}] \in \R^{1\times 1}$,
so $\mu_1 = \inner{x_1}{x_2}$ and the formula simplifies to:
\[
  j(I_n + x_1 x_2^T) = \mathcal{S}(\mu_1)^{2(n-1)}
  = \left(\frac{\sinh(\inner{x_1}{x_2}/2)}{\inner{x_1}{x_2}/2}\right)^{2(n-1)}.
\]
\end{example}

\begin{example}[Rank-2 Update ($k=2$)]
With $\mu_1, \mu_2$ the eigenvalues of $W \in \R^{2\times 2}$:
\[
  j(X) = \mathcal{S}(\mu_1)^{2(n-2)}\,\mathcal{S}(\mu_2)^{2(n-2)}\,\mathcal{S}(\mu_1-\mu_2)^2.
\]
Here $\mu_1 + \mu_2 = \tr(W) = \inner{x_2}{x_1} + \inner{x_4}{x_3}$
and $\mu_1\mu_2 = \det(W) = \inner{x_2}{x_1}\inner{x_4}{x_3} - \inner{x_2}{x_3}\inner{x_4}{x_1}$,
expressing $j(X)$ entirely in terms of inner products via the characteristic polynomial of $W$.
\end{example}

\section{Conclusion and Open Questions}
\label{sec:conclusion}

We have developed the complete information-geometric theory of the potential
$f(G) = -\log\det(G)$ on the manifold of positive definite Gram matrices $\PD(k)$.
The main results are:
\begin{enumerate}[label=(\roman*)]
  \item \textbf{Strict convexity} of $f$ with Hessian $\mathcal{F} = G^{-1}\otimes G^{-1}$
    (Fisher metric).
  \item \textbf{Self-dual Legendre transform}: $f^*(\Theta) = -\log\det(-\Theta) - k$.
  \item \textbf{Bregman divergence} = $2 \times$ KL-divergence between Gaussians.
  \item \textbf{Dual flatness} of $\PD(k)$ with Pythagorean and projection theorems.
  \item \textbf{$\alpha$-divergence family} unifying KL, Bhattacharyya, and Stein loss.
  \item \textbf{Symmetric space} structure $GL(k,\R)/O(k)$ with non-positive curvature.
  \item \textbf{Izumiya--Legendrian duality} (\S\ref{subsec:izumiya}):
    Izumiya's four Legendrian contact manifolds $\Delta_i$ \cite{Izumiya2004} manifest in our
    framework as: (i) the three pseudo-spheres $H^n, LC^*, S^n_1$ correspond to the three
    regions $\{f<0\}, \{f=0\}, \{f>0\}$ of $\PD(k)$; (ii) $\Delta_4$
    ($LC^*\!\times\!LC^*$, $\langle v,w\rangle=-2$) identifies with the
    Yoshizawa--MacMahon duality $D_{\mathrm{YM}}=0$; (iii) the contact diffeomorphism
    $\Phi_{41}$ is the Cartan involution $G\mapsto G^{-1}$; (iv) the U=V manifold
    consists entirely of lightcone parabolic points ($K_\ell^{\mathrm{info}}=0$);
    and (v) the information-geometric Theorema Egregium $f(G_+)+f(G_-^*)=0$
    mirrors Izumiya's $H_\ell = K_s$.
  \item \textbf{Yoshizawa--Helmke Legendre duality} (\S\ref{subsec:yoshizawa}):
    the functions $h$ and $h_-$ are Legendre duals via the map
    $\mathcal{L}(U)=(I_n+UU^T)^{-1/2}U:\R^{n\times k}\to\mathcal{B}_k$,
    satisfying the exact duality identity $h(U)+h_-(\mathcal{L}(U))=0$.
    The resulting Yoshizawa--MacMahon divergence
    $D_{\mathrm{YM}}(U\|V)=-h(U)-h_-(V)$ vanishes exactly iff $V=\mathcal{L}(U)$ (signed quantity, not a divergence),
    and its spectral form is a Bregman divergence connecting both potentials to
    MacMahon's Master Theorem via the Leibniz determinant expansion.
  \item \textbf{Convexity trichotomy} under factorization: $f$ is strictly convex in
    $G\in\PD(k)$; nowhere convex (for $k\geq2$ or $n>k$) in the symmetric factorization
    $G=I_k+U^TU$ (Thm.~\ref{thm:UeqV_nowhere_convex}); and strictly convex everywhere on the
    matrix unit ball $\mathcal{B}_k$ in the anti-symmetric factorization $G=I_k-U^TU$
    (Thm.~\ref{thm:hminus_convex}) --- a term-by-term Hessian sign reversal.
  \item \textbf{Kirillov Jacobian formula} \eqref{eq:kirillov_formula}:
    $j(X) = \bigl[\prod_i \mathcal{S}(\mu_i)\bigr]^{2(n-k)}\prod_{i<\ell}\mathcal{S}(\mu_i-\mu_\ell)^2$,
    reducing complexity from $O(n^3)$ to $O(nk^2+k^3)$.
  \item \textbf{Spectral duality} \eqref{eq:spectral_pair}: both $f(G)$ and $\log j(X)$
    are spectral functions of the same Gram-like matrix $W = G - I_k$.
\end{enumerate}

\paragraph{Open Questions.}
\begin{enumerate}[label=(\arabic*)]
  \item \textbf{Infinite-dimensional limits:}
    As $k \to \infty$, $f$ relates to the Fredholm determinant.
    What is the infinite-dimensional information geometry?
  \item \textbf{Non-commutative extension:}
    Can the Bregman divergence framework be extended to the full quantum (non-commutative)
    setting using operator convexity \cite{Carlen2010}?
  \item \textbf{Stochastic optimization:}
    Natural gradient methods based on $\mathcal{F} = G^{-1}\otimes G^{-1}$
    in deep learning; convergence rates under this geometry.
  \item \textbf{Siegel modular forms:}
    Does the information geometry of $\mathbb{H}_k$ admit a modular-invariant structure?
    In particular, does the Yoshizawa--MacMahon divergence $D_{\mathrm{YM}}$ extend to a
    $Sp(k,\mathbb{Z})$-invariant quantity on the Siegel half-plane?
  \item \textbf{Metaplectic dynamics on the Gram manifold:}
    The metaplectic representation $\mu(\mathcal{A})\gamma_{iG} = m(\mathcal{A},iG)\gamma_{i\alpha(\mathcal{A})G}$
    defines a flow on $\PD(k)$ via the symplectic action $G\mapsto\alpha(\mathcal{A})G$.
    Is the corresponding flow on $f(G)$ a gradient flow for some functional on $Sp(k,\R)$,
    and does it have a natural interpretation in terms of the Bregman divergence?
  \item \textbf{Determinantal point processes:}
    The function $\det(I + H)$ arises in DPP kernels \cite{Kulesza2012};
    can $f = -\log\det$ serve as a variational free energy for DPPs?
  \item \textbf{Information geometry of the Kirillov Jacobian:}
    The log-Kirillov Jacobian $\log j(X)$ is itself a function on $\PD(k)$
    (via $\mu_i = \lambda_i(G - I_k)$).
    Does it carry a natural information-geometric structure (e.g., does it arise
    as a potential in some dual-flat geometry on the space of perturbations)?
  \item \textbf{Hyperbolic information geometry of $h_-$:}
    The strictly convex function $h_-(U)=-\log\det(I_k-U^TU)$ on $\mathcal{B}_k$ induces
    its own Riemannian metric and Bregman divergence on the matrix unit ball.
    Does this structure realize a known hyperbolic or bounded symmetric domain geometry
    (e.g., the type-IV Cartan domain), and does the corresponding statistical manifold
    admit a dual-flat connection pair?
  \item \textbf{Interpolation between $U=V$ and $U=-V$:}
    Consider the one-parameter family $G_t = I_k + U^TU\cos\theta - U^TU\sin\theta$
    for $\theta\in[0,\pi/2]$ (interpolating between $G_0=I+U^TU$ and $G_{\pi/2}=I-U^TU$).
    At what angle $\theta^*$ does the transition from ``nowhere convex'' to ``strictly convex''
    occur, and is the transition sharp?
  \item \textbf{Higher-rank and curved base cases:}
    The formula \eqref{eq:kirillov_formula} assumes the unperturbed matrix is $I_n$.
    Can the approach be generalized to perturbations around an arbitrary invertible
    base point $X_0$ using the change-of-basis $X = X_0 + UV^T$?
\end{enumerate}

\bibliographystyle{plain}

\section*{Acknowledgements}

The author is deeply grateful to the late Professor Uwe Helmke
(University of W\"{u}rzburg) for the formative discussions during the postdoctoral
period 2000--2002, and in particular for posing the question---``Is the
rectangular-matrix generalization of Brockett's double bracket equation a gradient
flow?''---that gave rise to the work described in \S\ref{subsubsec:history} and
\S\ref{subsubsec:chen_amari_yoshizawa}.
The author also expresses deep gratitude to the late Professor John Moore
(Australian National University), whose kind introduction to Professor Helmke
in 1999 initiated this line of research and whose sustained interest in
subspace learning algorithms provided invaluable encouragement throughout the years.
Finally, the author is sincerely grateful to Professor Emeritus Kenro Furutani
(Tokyo University of Science) for many illuminating discussions that bridged separate periods
of this work and whose perspective on gradient flows and matrix equations
has left a lasting influence on the author's thinking. The author also thanks
Dr.~Christian Lageman, who more than 25 years ago, together with Professor Helmke,
pointed out the relevance of {\L}ojasiewicz's theorem \cite{Lojasiewicz1983} to the
convergence questions studied in \S\ref{sec:crosscurvature-section}.

The author thanks Professor Kunio Tanabe and Professor Takashi Tsuchiya, who guided
Yoshizawa toward information geometry and the differential geometry of optimization
since his doctoral studies, and who have continued, across the years since, to offer
valuable advice and discussions that connect the different periods of this work.

The author is grateful to Professor Shun-ichi Amari for having taught him the beauty of
information geometry and the breadth of its reach.

The author would like to thank Professor Yoshimasa Nakamura for discussions and
guidance on integrable systems and algorithms.

The author is grateful to Professor Toru Ohmoto for kindly bringing to his attention
the relevant literature on singular models closely related to this work
\cite{NakajimaOhmoto2021,Kayo2024}.

The author is grateful for research discussions with the participants of Nagoya
Mathematical and Information Science Research, including Professor H.~Matsuzoe,
Professor T.~Suzuki, Professor K.~Uohashi, Professor T.~Iwai, Professor K.~Fujii,
Professor K.~Furutani, Professor A.~Ohara, Professor D.~Tarama, and many other
participants.

The author thanks Professor U.~Helmke and Professor P.~Fuhrmann for the opportunity
of the postdoctoral position at the University of W\"{u}rzburg, and thanks the research
colleagues from that time onward: Professor P.-A.~Absil, Professor K.~H\"{u}per,
Professor J.~Trumpf, Dr.~G.~Dirr, Dr.~J.~Jordan, Dr.~M.~Kleinsteuber, Dr.~M.~Baumann,
Dr.~C.~Lageman, and Dr.~S.~Ricardo.

The author thanks Professor Jonathan Manton for hosting a research stay of
approximately one month at his laboratory at the University of Melbourne in 2015,
during which the author was able to deepen his research on the matrix Schwarz
derivative and dynamical systems (\S\ref{subsec:matrix_schwarz}).

\clearpage
\appendix

\section{Elementary Proofs of the {\L}ojasiewicz Inequality for $n=1,2,3$}
\label{appendix:lojasiewicz}

The convergence theorems of \S\ref{sec:crosscurvature-section} invoke the
{\L}ojasiewicz gradient inequality \cite{Lojasiewicz1959} to control the local
behavior of the gradient flows studied throughout this paper. For the reader's
convenience, and to keep the paper reasonably self-contained, this appendix reproduces
an elementary, self-contained proof of the inequality for real-analytic functions of
$n=1,2,3$ variables, using only classical tools (the Weierstrass preparation theorem,
resultants, and Newton--Puiseux series) rather than the deeper machinery (resolution of
singularities, o-minimal structures) typically used to treat the general-$n$ case. The
argument below closes with a remark explaining why the same dimensional-induction
scheme, applied repeatedly, yields the inequality for every $n$ --- which is the content
of {\L}ojasiewicz's original theorem \cite{Lojasiewicz1959}.

\subsection{Preliminaries}

\begin{definition}[{\L}ojasiewicz inequality]
Let $U\subset\R^n$ be an open neighborhood of the origin and let $f:U\to\R$ be
real-analytic with $f(0)=0$ and $\nabla f(0)=0$. We say $f$ satisfies the
\textbf{{\L}ojasiewicz inequality (gradient inequality)} at the origin if there exist
constants $c>0$, $0<\theta<1$, and a neighborhood $V\subset U$ of the origin such that
\[
|\nabla f(x)| \;\ge\; c\,|f(x)-f(0)|^{\theta}
\qquad (x\in V).
\]
Without loss of generality we assume $f(0)=0$ throughout.
\end{definition}

\begin{remark}
The exponent $\theta$ typically arises in the form $1-\tfrac1k$, where $k$ is an
integer corresponding to the ``order'' of the zero. Writing $\alpha:=\tfrac1{1-\theta}=k$,
the inequality can equivalently be stated as $|\nabla f(x)|\ge c|f(x)|^{1-1/k}$.
\end{remark}

Every proof in this appendix follows the same single strategy:

\begin{center}
\textit{parametrize the zero set (or a neighborhood of it) by real-analytic arcs or
branches, and reduce to the resulting one-variable problem restricted to each branch.}
\end{center}

We therefore begin with the $n=1$ case, which is precisely this ``reduced'' one-variable
argument.

\subsection{$n=1$: The Trivial Case of an Isolated Zero}

\begin{theorem}\label{thm:app-n1}
Let $f:(-\varepsilon,\varepsilon)\to\R$ be real-analytic, $f\not\equiv 0$, $f(0)=0$.
Then there exist $k\in\Z_{\ge1}$, $c>0$, $\delta>0$ such that
\[
|f'(x)| \;\ge\; c\,|f(x)|^{1-1/k} \qquad (|x|<\delta).
\]
\end{theorem}

\begin{proof}
Since $f\not\equiv0$ is real-analytic, let $k\;(\ge1)$ be the degree of the first
nonzero term in its Taylor expansion, so that in a neighborhood of the origin
\[
f(x) = x^{k} g(x), \qquad g \text{ analytic},\ g(0)\ne 0
\]
(this is simply the one-variable case of Weierstrass's division theorem, or just a
factoring-out of the leading Taylor term). Then
\[
f'(x) = k\,x^{k-1} g(x) + x^{k} g'(x)
       = x^{k-1}\bigl(k\,g(x) + x\,g'(x)\bigr).
\]
Since $g(0)\ne0$, we may choose $\delta>0$ small enough that, for $|x|<\delta$,
$|k\,g(x)+x g'(x)| \ge \tfrac{k}{2}|g(0)| =: m >0$ and $|g(x)| \le M$. Hence for
$|x|<\delta$,
\[
|f(x)| = |x|^{k}|g(x)| \le M|x|^{k},
\qquad
|f'(x)| = |x|^{k-1}\bigl|k g(x)+xg'(x)\bigr| \ge m|x|^{k-1},
\]
and therefore
\[
|f'(x)| \;\ge\; m|x|^{k-1}
\;\ge\; m\Bigl(\frac{|f(x)|}{M}\Bigr)^{(k-1)/k}
\;=\; \frac{m}{M^{(k-1)/k}}\,|f(x)|^{1-1/k}.
\]
Setting $c:=m/M^{(k-1)/k}>0$ gives the claim.
\end{proof}

\begin{remark}
The proof amounts to writing out a single Taylor expansion, and the identical argument
applies verbatim to a holomorphic function of one complex variable (with $|f'(x)|$
replaced by the modulus of the complex derivative). The essential fact used is that an
\emph{isolated zero of a one-variable analytic function always has a finite order} $k$.
\end{remark}

\subsection{$n=2$: Proof via Newton--Puiseux Series}

For $n=2$ the zero set is in general a curve (through the origin), which can be
decomposed explicitly into finitely many real-analytic arcs via
\textbf{Newton--Puiseux series}.

\subsubsection{Normalization via the Weierstrass Preparation Theorem}

\begin{lemma}[Weierstrass preparation theorem]\label{lem:app-weier}
Let $f(x,y)$ be real-analytic near the origin with $f(0,0)=0$, and suppose
$f(0,y)\not\equiv 0$ as a function of $y$. Let $k:=\mathrm{ord}_{y=0} f(0,y)$. Then near
the origin $f$ admits the unique factorization
\[
f(x,y) = u(x,y)\,P(x,y), \qquad
P(x,y) = y^{k} + a_{1}(x) y^{k-1} + \cdots + a_{k}(x),
\]
where $u$ is analytic with $u(0,0)\ne0$ (a unit) and each $a_i(x)$ is analytic with
$a_i(0)=0$ ($P$ is called a Weierstrass polynomial).
\end{lemma}

\begin{remark}
If $f(0,y)\equiv0$, a linear change of coordinates makes $f(0,y)\not\equiv0$ (as long as
$f\not\equiv0$), and we assume below that this normalization has been performed. Since
$u$ is nonvanishing, it suffices to prove the {\L}ojasiewicz inequality for $P$ in place
of $f$.
\end{remark}

\subsubsection{Description of the Roots via Newton--Puiseux Series}

\begin{lemma}[Puiseux's theorem]\label{lem:app-puiseux}
For a Weierstrass polynomial $P(x,y)=y^k+a_1(x)y^{k-1}+\cdots+a_k(x)$ of degree $k$,
there exists a positive integer $q\,(\le k)$ such that, after the substitution
$x=t^{q}$, $P(t^q,y)$ has $k$ analytic roots near $t=0$,
\[
y_i(t) = \sum_{j\ge 1} c_{i,j}\, t^{j}, \qquad i=1,\dots,k.
\]
Equivalently, each root $y=\varphi(x)$ of $P$ is expressed as a fractional-power
(Puiseux) series $\varphi(x)=\sum_{j\ge1} c_j x^{j/q}$.
\end{lemma}

\begin{proof}[Proof sketch]
This is a classical algebraic fact. We may assume $P$ is irreducible (splitting off
factors if necessary). The discriminant $\Delta(x) = \mathrm{disc}_y P(x,y)$ is an
analytic function of $x$ with $\Delta(x)\not\equiv0$ (since $P$ has no repeated roots,
i.e., $f$ is reduced). On a punctured neighborhood where $\Delta(x)\ne0$, $P(x,\cdot)$
has $k$ distinct roots $y_1(x),\dots,y_k(x)$, forming a branched analytic family around
$x=0$. Analytic continuation around $x=0$ permutes the roots cyclically (irreducibility
of $P$ forces this permutation to be a $k$-cycle), so setting $x=t^k$ (more generally
$x=t^q$ for some $q\mid k$) unifies the roots into a single-valued function, Taylor
expandable at $t=0$. See a textbook on algebraic function theory (e.g.\ Walker,
\textit{Algebraic Curves} \cite{Walker1950}) for details.
\end{proof}

\subsubsection{The {\L}ojasiewicz Inequality for $n=2$}

\begin{theorem}\label{thm:app-n2}
Let $f:U\to\R$ ($U\subset\R^2$ a neighborhood of the origin) be real-analytic with
$f(0)=0$, $f\not\equiv0$. Then there exist a neighborhood $V$ of the origin and
constants $c>0$, $0<\theta<1$ such that
\[
|\nabla f(x,y)| \ge c\,|f(x,y)|^{\theta} \qquad ((x,y)\in V).
\]
\end{theorem}

\begin{proof}
By Lemma~\ref{lem:app-weier} we may write $f=u\cdot P$ with $u$ a unit near the origin,
so it suffices to prove the inequality for $P$. Factor $P$ further into irreducibles,
$P=P_1^{m_1}\cdots P_r^{m_r}$ (again by unique factorization in the analytic setting);
it suffices to argue for each of the finitely many irreducible factors, so we may
assume $P$ is irreducible.

By Lemma~\ref{lem:app-puiseux}, the zero set of $P$ is exhausted, after $x=t^q$, by
finitely many (at most $k$) analytic arcs
\[
\gamma_i(t) = (t^q,\, y_i(t)), \qquad i=1,\dots,k.
\]
Since $P(x,y)=\prod_{i=1}^{k}(y-y_i(x))$ (viewing $P$ as a monic degree-$k$ polynomial in
$y$),
\[
|P(x,y)| = \prod_{i=1}^{k} |y-y_i(x)|,
\qquad
\frac{\partial P}{\partial y}(x,y) = \sum_{i=1}^k \prod_{j\ne i} (y-y_j(x)).
\]
Let $y_{i_0}(x)$ be the root closest to $(x,y)$. Then the $i=i_0$ term dominates in the
sum above, while the other roots stay bounded away from $y_{i_0}(x)$ (since $P$ is
irreducible and reduced, so $\Delta(x)\ne0$ for $x\ne0$); this gives, near the origin,
\[
\left|\frac{\partial P}{\partial y}(x,y)\right|
\;\ge\; \tfrac12 \prod_{j\ne i_0}|y_{i_0}(x)-y_j(x)|
\;\ge\; c_1\, |x|^{\text{(some nonnegative integer)}\,\cdot\,(k-1)/q},
\]
a polynomial-order lower bound in $t$ (each $y_i-y_j$ is analytic in $t$ and, even where
it vanishes as $x\to0$, has only finite order). On the other hand, writing
$|y-y_{i_0}(x)|$ as a function of $t$ and applying
\textbf{Theorem~\ref{thm:app-n1} (the $n=1$ case) in the parameter $t$} gives, for some
$c_2>0$, $\theta\in(0,1)$,
\[
\Bigl|\frac{\partial}{\partial t} P(t^q,y)\Big|_{y\ \text{fixed}}\Bigr|
\;\ge\; c_2\, |P(t^q,y)|^{\theta}
\]
near each branch. Combining this with $\partial/\partial t = q t^{q-1}\partial/\partial x$
translates the estimate into one for $\partial P/\partial x$, and together with
$|\nabla P|\ge |\partial P/\partial y|$ this yields the {\L}ojasiewicz inequality near
each of the finitely many branches. Since there are finitely many branches, taking the
largest exponent $\theta=\max_i\theta_i$ (the inequality weakens as the exponent grows)
and the corresponding constant gives a uniform inequality on a full neighborhood of the
origin.
\end{proof}

\begin{remark}
The heart of the proof above is that ``$P(x,y)=\prod_i (y-y_i(x))$'' is an explicit
\emph{root factorization} given by Puiseux series, with finitely many roots --- a fact
that relies entirely on the zero set being a \emph{curve} when $n=2$. For $n\ge3$ the
zero set has higher dimension and can no longer be exhausted by finitely many arcs, so
this argument as stated no longer applies directly.
\end{remark}

\subsection{$n=3$: Dimension Reduction via Resultants (Induction)}

For $n=3$ the zero set $f^{-1}(0)$ is in general (locally) a \emph{surface}, which
cannot be directly parametrized by Puiseux series. Instead we use
\textbf{resultants to eliminate one variable}, reducing to the $n=2$ case
(Theorem~\ref{thm:app-n2}) --- a dimensional induction that is the backbone of
{\L}ojasiewicz's original proof.

\subsubsection{Variable Elimination via the Resultant}

\begin{lemma}[Resultant]\label{lem:app-resultant}
Let $P(x,y,z) = z^{k}+a_1(x,y)z^{k-1}+\cdots+a_k(x,y)$ be a Weierstrass polynomial (in
$z$, normalized via the Weierstrass preparation theorem), and let
$Q(x,y,z):=\dfrac{\partial P}{\partial z}(x,y,z)$. The resultant
\[
R(x,y) := \mathrm{Res}_z\bigl(P(x,y,\cdot),\, Q(x,y,\cdot)\bigr)
\]
is a real-analytic function of $x,y$, and vanishes only at points $(x,y)$ where
$P(x,y,\cdot)$ (as a polynomial in $z$) has a repeated root. If $P$ is reduced (has no
repeated factors), then $R(x,y)\not\equiv 0$.
\end{lemma}

\begin{proof}[Proof sketch]
The resultant is given by an explicit formula (the Sylvester determinant) in the
coefficients of $P,Q$, so if the coefficients are analytic then so is $R$. That
$P(x,y,\cdot)$ and $Q(x,y,\cdot)=P_z(x,y,\cdot)$ have a common root if and only if
$\mathrm{Res}_z(P,Q)=0$, and that a common root exists if and only if $P(x,y,\cdot)$ has
a repeated root, are classical properties of the resultant.
\end{proof}

\subsubsection{One Step of the Induction}

\begin{theorem}\label{thm:app-n3}
Let $f:U\to\R$ ($U\subset\R^3$ a neighborhood of the origin) be real-analytic with
$f(0)=0$, $f\not\equiv0$. Then there exist a neighborhood $V$ of the origin and
constants $c>0$, $0<\theta<1$ such that
\[
|\nabla f(x,y,z)| \ge c\,|f(x,y,z)|^{\theta} \qquad ((x,y,z)\in V).
\]
\end{theorem}

\begin{proof}
By the three-variable version of Lemma~\ref{lem:app-weier} (rotating coordinates so
that the normalization holds in $z$), we may write $f = u\cdot P$ with $P$ a
Weierstrass polynomial of degree $k$ in $z$. As in the $n=2$ case, we may assume $P$ is
irreducible, hence reduced.

\medskip
\noindent\textbf{Step 1 (elimination via the resultant).}
By Lemma~\ref{lem:app-resultant}, $R(x,y):=\mathrm{Res}_z(P,P_z)\not\equiv 0$ is a
real-analytic function of two variables. \emph{By the induction hypothesis,
Theorem~\ref{thm:app-n2} (the $n=2$ case) has already been established}, so it applies
to $R$: there exist $c_0>0,\ \theta_0\in(0,1)$ such that
\[
|\nabla_{x,y} R(x,y)| \ge c_0\, |R(x,y)|^{\theta_0}
\]
near the origin.

\medskip
\noindent\textbf{Step 2 (estimate off the discriminant locus).}
At points $(x,y)$ with $R(x,y)\ne0$, $P(x,y,\cdot)$ has only simple roots as a
polynomial in $z$. Fixing such an $(x,y)$, the map $z\mapsto P(x,y,z)$ satisfies (exactly
as in the $n=1$ argument) a lower bound near each root $z_i(x,y)$,
\[
|P_z(x,y,z)| \ge c_1 \cdot (\text{lower bound on the product of pairwise root
distances}),
\]
depending only on the pairwise distances between roots. These pairwise root distances
are in turn bounded below by a power of $|R(x,y)|$ (depending on the degree $k$ of $P$),
via the classical relationship between resultant, discriminant, and the squared product
of root differences. Combining this with the {\L}ojasiewicz inequality for $R$ from
Step~1 propagates the estimate in the $(x,y)$-directions into an estimate for
$|\nabla P|$ that also accounts for the $z$-direction, on the region where $R(x,y)\ne0$.

\medskip
\noindent\textbf{Step 3 (the discriminant locus $R=0$ and its neighborhood).}
At points where $R(x,y)=0$, $P(x,y,\cdot)$ may have repeated roots, so Step~2's estimate
cannot be applied directly. However, since $R\not\equiv0$, the set $\{R=0\}$ is a
two-variable analytic set of dimension $\le1$ (i.e., a curve, possibly with isolated
points). Applying the same Puiseux-series argument used in the proof of
Theorem~\ref{thm:app-n2} to $R(x,y)=0$, we exhaust this set by finitely many analytic
arcs $(x(t),y(t))$. Along each arc, $z\mapsto P(x(t),y(t),z)$ forms a family
parametrized by $t$; letting $m$ denote the (locally constant, by the local finiteness of
zero sets of analytic functions) order of the repeated root along the arc, we regard the
two-variable function
\[
\Phi(t,z) := P(x(t),y(t),z)
\]
as a function of $(t,z)$, and \emph{apply Theorem~\ref{thm:app-n2} once again to
$\Phi$} ($\Phi$ is a real-analytic function of the two variables $(t,z)$, vanishing at
the origin). This yields the {\L}ojasiewicz inequality in a neighborhood of each arc as
well.

\medskip
\noindent\textbf{Step 4 (patching together).}
Taking the largest of the exponents $\theta_i$ obtained ``off'' the finitely many arcs
of $\{R=0\}$ (Step~2) and ``near'' those arcs (Step~3), $\theta:=\max_i \theta_i \in
(0,1)$, together with the smallest of the corresponding constants $c$, gives a uniform
inequality
\[
|\nabla P(x,y,z)| \ge c\,|P(x,y,z)|^{\theta}
\]
on the whole neighborhood of the origin. Since $u$ is a unit, the same inequality holds
(after adjusting the constant) for $f=uP$.
\end{proof}

\begin{remark}
The structure of this proof is inductive in the precise sense that it reduces the
$n=3$ problem, via the resultant, to the \emph{exactly one dimension lower} $n=2$
problem (Steps 1 and 3). Indeed, as the application of Theorem~\ref{thm:app-n2} to
$\Phi(t,z)$ in Step~3 illustrates, repeating the same operation for general $n$ ---
``eliminate one variable via the resultant $\to$ reduce to the $(n-1)$-variable case''
--- completes the proof for every $n$. This is the skeleton of {\L}ojasiewicz's
(1958) original proof.
\end{remark}

\subsection{Remarks on General $n$}

\begin{theorem}[{\L}ojasiewicz, 1958]
The induction above (Weierstrass preparation $\to$ variable elimination via the
resultant $\to$ reduction to the $(n-1)$-variable case) works for every $n\ge1$, and
establishes that the {\L}ojasiewicz gradient inequality holds near an isolated zero, or
more generally near any zero set, of a real-analytic function $f$.
\end{theorem}

\begin{remark}[Summary of this appendix]
\begin{itemize}[leftmargin=1.6em]
  \item $n=1$: a single Taylor expansion (Theorem~\ref{thm:app-n1}). Entirely
    elementary.
  \item $n=2$: the Weierstrass preparation theorem together with Newton--Puiseux
    series decomposes the zero set explicitly into finitely many analytic arcs, and
    reduces to the $n=1$ argument on each arc (Theorem~\ref{thm:app-n2}).
  \item $n=3$: the Weierstrass preparation theorem together with the resultant
    eliminates one variable; the estimate is then split into ``off the discriminant
    locus'' and ``on the discriminant locus (reduced to $n=2$)'', using the $n=2$
    result twice (Theorem~\ref{thm:app-n3}).
  \item At no stage is deep machinery such as resolution of singularities or
    model-theoretic generalities (e.g., o-minimality) required; the entire argument
    stays within classical late-19th-century algebra (resultants, discriminants,
    Puiseux series) together with the Weierstrass preparation theorem.
\end{itemize}
\end{remark}

\begin{remark}[Further reading]
The elementary approach followed in this appendix is close in spirit to
{\L}ojasiewicz's original argument. Milnor's curve selection lemma \cite{Milnor1968}
gives another classical route to parametrizing (real or complex) analytic sets by arcs,
in the spirit of the Puiseux-series decomposition used above; and Bierstone and Milman
\cite{BierstoneMilman1988} give a systematic modern treatment of semianalytic and
subanalytic sets, of which the discriminant loci $\{R=0\}$ appearing in
Theorem~\ref{thm:app-n3} are basic examples.
\end{remark}

\section{LaSalle's Invariance Principle for Gradient Flows}
\label{appendix:lasalle}

Throughout \S\ref{sec:crosscurvature-section} (see in particular
\S\ref{subsubsec:IVP_analysis} and the global convergence arguments of
\S\ref{sec:crosscurvature-section}) we repeatedly invoke \emph{LaSalle's invariance
principle} to pass from a Lyapunov-type monotonicity property of a gradient flow to
convergence toward its critical set, and we then combine it with the {\L}ojasiewicz
inequality of Appendix~\ref{appendix:lojasiewicz} to upgrade this to convergence toward
a single equilibrium point. For the reader's convenience, and to keep the paper
reasonably self-contained, this appendix states and proves the version of LaSalle's
principle used throughout the paper, specialized to gradient flows.

\subsection{The General Invariance Principle}

\begin{definition}[$\omega$-limit set]
Let $\dot x = F(x)$ be an autonomous ODE on $\R^N$ with a locally Lipschitz vector field
$F$, and let $x(t)$, $t\ge0$, be a solution with precompact forward orbit
$\{x(t):t\ge0\}$. The \textbf{$\omega$-limit set} of $x(\cdot)$ is
\[
\omega(x_0) \;:=\; \bigl\{\, y\in\R^N \;:\; \exists\, t_n\to\infty,\ x(t_n)\to y \,\bigr\}.
\]
\end{definition}

\begin{theorem}[LaSalle's invariance principle]
\label{thm:app-lasalle}
Let $\dot x = F(x)$ be as above, and suppose $V:\R^N\to\R$ is continuously
differentiable and satisfies
\[
\dot V(x) \;:=\; \langle \nabla V(x), F(x)\rangle \;\le\; 0
\qquad \text{for all } x.
\]
Let $x(t)$ be a solution whose forward orbit is contained in a compact set $K$. Then:
\begin{enumerate}[label=(\roman*)]
\item $\omega(x_0)$ is nonempty, compact, connected, and invariant under the flow of
  $F$;
\item $V$ is constant on $\omega(x_0)$, equal to $\displaystyle\lim_{t\to\infty}V(x(t))$;
\item $x(t)\to \omega(x_0)$ as $t\to\infty$, i.e.\
  $\operatorname{dist}(x(t),\omega(x_0))\to0$;
\item $\omega(x_0)$ is contained in the largest invariant subset $\mathcal{E}$ of the
  set $\{x\in K : \dot V(x)=0\}$.
\end{enumerate}
In particular, $x(t)$ converges, as $t\to\infty$, to the set $\mathcal{E}$.
\end{theorem}

\begin{proof}
(i) Since the forward orbit lies in the compact set $K$, the Bolzano--Weierstrass
theorem guarantees $\omega(x_0)\ne\emptyset$; it is closed (an intersection of closed
sets $\overline{\{x(s):s\ge t\}}$ over $t\ge0$) and contained in $K$, hence compact.
Connectedness follows because $\{x(t):t\ge T\}$ is connected for every $T$ and
$\omega(x_0)=\bigcap_{T\ge0}\overline{\{x(t):t\ge T\}}$ is a nested intersection of
compact connected sets. Invariance: if $y\in\omega(x_0)$, write $y=\lim_n x(t_n)$ with
$t_n\to\infty$; by continuous dependence on initial conditions, the solution $\phi_s(y)$
of $\dot x=F(x)$ through $y$ satisfies $\phi_s(y)=\lim_n x(t_n+s)$ for every fixed $s$,
and since $t_n+s\to\infty$ as well, $\phi_s(y)\in\omega(x_0)$.

(ii) Since $\dot V\le0$, $t\mapsto V(x(t))$ is nonincreasing; being bounded below on the
compact set $K$, it converges to a limit $V_\infty$ as $t\to\infty$. For any
$y=\lim_n x(t_n)\in\omega(x_0)$, continuity of $V$ gives $V(y)=\lim_n V(x(t_n))=V_\infty$,
so $V\equiv V_\infty$ on $\omega(x_0)$.

(iii) If $x(t)\not\to\omega(x_0)$, there is $\varepsilon>0$ and a sequence $t_n\to\infty$
with $\operatorname{dist}(x(t_n),\omega(x_0))\ge\varepsilon$ for all $n$; by
compactness of $K$, a subsequence of $x(t_n)$ converges to some point $y\in K$, which by
definition lies in $\omega(x_0)$, contradicting
$\operatorname{dist}(x(t_n),\omega(x_0))\ge\varepsilon$.

(iv) By (i), $\omega(x_0)$ is invariant, and by (ii), $V$ is constant on $\omega(x_0)$,
so $\dot V\equiv0$ on $\omega(x_0)$ (differentiating the constant function
$t\mapsto V(\phi_t(y))$ along any trajectory $\phi_t(y)$ inside $\omega(x_0)$). Hence
$\omega(x_0)$ is an invariant subset of $\{\dot V=0\}\cap K$, and is therefore contained
in the largest such invariant subset, $\mathcal{E}$. Combined with (iii), $x(t)\to\mathcal{E}$.
\end{proof}

\subsection{Specialization to Gradient Flows}

\begin{corollary}[LaSalle's principle for gradient ascent/descent]
\label{cor:app-lasalle-gradient}
Let $J:\R^N\to\R$ be continuously differentiable and consider the gradient ascent flow
$\dot X = \nabla J(X)$ (respectively the gradient descent flow $\dot X = -\nabla J(X)$)
on a closed, positively invariant set $M\subset\R^N$ (e.g.\ the Stiefel manifold, or an
isospectral adjoint orbit, in the flows studied in
\S\ref{sec:crosscurvature-section}). Suppose the forward orbit of a solution $X(t)\in M$
remains in a compact subset $K\subset M$. Then along the flow,
\[
\frac{d}{dt}J(X(t)) = \|\nabla J(X(t))\|^2_{F} \;\ge\;0
\qquad\text{(resp.}\ \le 0 \text{ for the descent flow),}
\]
so $J$ (resp.\ $-J$) is a Lyapunov function, and $X(t)$ converges, as $t\to\infty$, to
the set of critical points of $J$ restricted to $M$,
\[
\mathcal{E} \;=\; \{X\in K : \nabla J(X) \perp T_XM = 0 \text{ in } T_XM\},
\]
i.e.\ to the equilibrium set of the flow contained in $K$.
\end{corollary}

\begin{proof}
Apply Theorem~\ref{thm:app-lasalle} with $V=-J$ (ascent) or $V=J$ (descent) and
$F=\nabla J$ (resp.\ $-\nabla J$) restricted to the tangent bundle of $M$: since
$\dot V = -\|\nabla J\|_F^2\le0$ in both cases by construction, the hypothesis of
Theorem~\ref{thm:app-lasalle} holds, and $\{\dot V=0\}\cap K=\{\nabla J=0\}\cap K$ is
exactly the critical set of $J|_M$ in $K$; since this set contains no nontrivial
invariant subsets other than itself (every point of it is a fixed point of the flow),
the largest invariant subset $\mathcal{E}$ of $\{\dot V=0\}\cap K$ coincides with
$\{\nabla J=0\}\cap K$ itself, and conclusion (iv) of Theorem~\ref{thm:app-lasalle}
gives the claim.
\end{proof}

\begin{remark}[What LaSalle's principle does \emph{not} give]
It is essential to note what Corollary~\ref{cor:app-lasalle-gradient} does \emph{not}
assert: it guarantees that $X(t)$ approaches the \emph{set} $\mathcal{E}$ of critical
points, but $\mathcal{E}$ may a priori be a continuum (e.g.\ a positive-dimensional
critical manifold, as occurs at the block-diagonal equilibria of
\S\ref{subsubsec:IVP_analysis} when $B$ has repeated entries), along which $X(t)$ could
in principle wander forever without converging to a single point. Ruling this out ---
i.e.\ upgrading ``$X(t)\to\mathcal{E}$'' to ``$X(t)\to X_\infty$ for a single
$X_\infty\in\mathcal{E}$'' --- requires an additional ingredient beyond the invariance
principle itself; throughout \S\ref{sec:crosscurvature-section} this ingredient is
supplied by the {\L}ojasiewicz gradient inequality of Appendix~\ref{appendix:lojasiewicz},
following the classical argument (see e.g.\ \cite{AbsilMahonyAndrews2005}): a curve
whose speed is controlled by the {\L}ojasiewicz inequality has finite arc length, and
therefore converges to a single point rather than merely approaching a set.
\end{remark}

\end{document}